\pdfoutput=1
\documentclass[a4paper,10pt,english,leqno]{book}

\usepackage[utf8]{inputenc}
\usepackage[T1]{fontenc}
\usepackage[english]{babel}

\usepackage{amsmath}
\usepackage{amsthm}
\usepackage{amssymb}
\usepackage{amscd}
\usepackage{mathrsfs}
\usepackage[all]{xy}

\usepackage{XCharter} 
\usepackage[xcharter]{newtxmath} 
\usepackage{microtype}
\usepackage{textcomp}

\usepackage[x11names]{xcolor}
\definecolor{MainColor}{RGB}{24, 60, 112} 
\definecolor{DefColor}{RGB}{34, 139, 34}  
\definecolor{PropColor}{RGB}{128, 24, 34} 

\definecolor{ThmBg}{RGB}{245, 248, 252}   
\definecolor{DefBg}{RGB}{245, 252, 245}   
\definecolor{PropBg}{RGB}{253, 248, 248}  

\usepackage[headheight=21pt]{geometry}
\usepackage{lscape}

\usepackage{fancyhdr}
\renewcommand{\headrulewidth}{0.5pt}
\renewcommand{\headrule}{\hbox to\headwidth{\color{MainColor}\leaders\hrule height \headrulewidth\hfill}}

\usepackage{titlesec}
\titleformat{\chapter}[display]
  {\normalfont\Huge\bfseries\color{MainColor}}
  {\filleft\MakeUppercase{\chaptertitlename} \Huge\thechapter}
  {2ex}
  {\titlerule\vspace{1ex}\filleft}
  [\vspace{1ex}\titlerule]
\titleformat{\section}
  {\normalfont\Large\bfseries\color{MainColor}}{\thesection}{1em}{}
\titleformat{\subsection}
  {\normalfont\large\bfseries\color{MainColor}}{\thesubsection}{1em}{}

\usepackage{thmtools}
\usepackage[framemethod=TikZ]{mdframed}

\mdfsetup{skipabove=1.5em,skipbelow=1.5em}

\declaretheoremstyle[
    headfont=\bfseries\color{MainColor},
    bodyfont=\itshape,
    mdframed={
        linewidth=2pt,
        rightline=false, topline=false, bottomline=false,
        linecolor=MainColor,
        backgroundcolor=ThmBg,
        nobreak=false
    }
]{thmstyle}

\declaretheoremstyle[
    headfont=\bfseries\color{PropColor},
    bodyfont=\itshape,
    mdframed={
        linewidth=2pt,
        rightline=false, topline=false, bottomline=false,
        linecolor=PropColor,
        backgroundcolor=PropBg,
        nobreak=false
    }
]{propstyle}

\declaretheoremstyle[
    headfont=\bfseries\color{DefColor},
    bodyfont=\normalfont,
    mdframed={
        linewidth=2pt,
        rightline=false, topline=false, bottomline=false,
        linecolor=DefColor,
        backgroundcolor=DefBg,
        nobreak=false
    }
]{defstyle}

\declaretheoremstyle[
    headfont=\itshape\color{black!70},
    bodyfont=\normalfont,
    notefont=\itshape,
    notebraces={(}{)},
    spaceabove=1em,
    spacebelow=1em
]{remstyle}

\declaretheorem[style=thmstyle, name=Theorem, numbered=no]{thm}
\declaretheorem[style=thmstyle, name=Lemma, numbered=no]{lemme}
\declaretheorem[style=thmstyle, name=Variant, numbered=no]{variante}

\declaretheorem[style=propstyle, name=Proposition, numbered=no]{prop}
\declaretheorem[style=propstyle, name=Corollary, numbered=no]{cor}

\declaretheorem[style=defstyle, name=Definition, numbered=no]{defi}
\declaretheorem[style=defstyle, name=Scholium, numbered=no]{scolie}

\declaretheorem[style=remstyle, name=Remark, numbered=no]{rmq}
\declaretheorem[style=remstyle, name=Remarks, numbered=no]{rmqs}
\declaretheorem[style=remstyle, name=Example, numbered=no]{exemple}
\declaretheorem[style=remstyle, name=Examples, numbered=no]{exemples}
\declaretheorem[style=remstyle, name=Notation, numbered=no]{notation}

\usepackage{tocloft}
\usepackage{makeidx}
\usepackage{index}

\disableindex{default}
\newindex{not}{ndx}{ndd}{List of notation}
\newindex{ter}{tdx}{tdd}{Subject Index}
\makeindex

\usepackage[colorlinks=true, linkcolor=MainColor, citecolor=DefColor, urlcolor=MainColor]{hyperref}
\hypersetup{hypertexnames=false}

\usepackage{verbatim}

\usepackage{emptypage} 

\newcommand{\frM}{\mathfrak{M}}  \newcommand{\frO}{\mathfrak{O}}
\newcommand{\frP}{\mathfrak{P}}

 \newcommand{\frZ}{\mathfrak{Z}}

\newcommand{\fra}{\mathfrak{a}}  
  
\newcommand{\frg}{\mathfrak{g}}

\newcommand{\frs}{\mathfrak{s}} \newcommand{\frt}{\mathfrak{t}} 
  
 \newcommand{\frz}{\mathfrak{z}}

\newcommand{\bbA}{\mathbb{A}}  \newcommand{\bbC}{\mathbb{C}}
  \newcommand{\bbF}{\mathbb{F}}
\newcommand{\bbG}{\mathbb{G}}  
  
 \newcommand{\bbN}{\mathbb{N}} 
\newcommand{\bbP}{\mathbb{P}}  \newcommand{\bbR}{\mathbb{R}}
 \newcommand{\bbT}{\mathbb{T}} \newcommand{\bbU}{\mathbb{U}}
  
 \newcommand{\bbZ}{\mathbb{Z}}

\newcommand{\caA}{\mathcal{A}} \newcommand{\caB}{\mathcal{B}} 
 \newcommand{\caE}{\mathcal{E}} \newcommand{\caF}{\mathcal{F}}
\newcommand{\caG}{\mathcal{G}} \newcommand{\caH}{\mathcal{H}} \newcommand{\caI}{\mathcal{I}}
  \newcommand{\caL}{\mathcal{L}}
\newcommand{\caM}{\mathcal{M}} \newcommand{\caN}{\mathcal{N}} \newcommand{\caO}{\mathcal{O}}
\newcommand{\caP}{\mathcal{P}} \newcommand{\caQ}{\mathcal{Q}} \newcommand{\caR}{\mathcal{R}}
\newcommand{\caS}{\mathcal{S}} \newcommand{\caT}{\mathcal{T}} \newcommand{\caU}{\mathcal{U}}
\newcommand{\caV}{\mathcal{V}} \newcommand{\caW}{\mathcal{W}} \newcommand{\caX}{\mathcal{X}}
 \newcommand{\caZ}{\mathcal{Z}}

  \newcommand{\scrC}{\mathscr{C}}
\newcommand{\scrD}{\mathscr{D}} \newcommand{\scrE}{\mathscr{E}} 
  \newcommand{\scrI}{\mathscr{I}}
\newcommand{\scrJ}{\mathscr{J}}

\newcommand{\una}{\mathbf{1}_A}

\DeclareMathOperator{\coker}{coker}
\DeclareMathOperator{\Dim}{\Dim}

\DeclareMathOperator{\End}{End}
\DeclareMathOperator{\Hom}{Hom}
\DeclareMathOperator{\im}{Im}

\DeclareMathOperator{\Ind}{Ind}
\DeclareMathOperator{\supp}{Supp}
\DeclareMathOperator{\tr}{Tr}
\DeclareMathOperator{\Id}{Id}
\DeclareMathOperator{\Int}{Int}
\DeclareMathOperator{\ind}{ind}
\DeclareMathOperator{\res}{Res}
\DeclareMathOperator{\vol}{Vol}

\newcommand{\spcheck}{\check{}\, }
\newcommand{\bil}[2]{\langle  #1,#2 \rangle }
\newcommand{\bilo}{\langle \,  .\, , . \, \rangle }

\renewcommand{\thechapter}{\Roman{chapter}}
\renewcommand{\thesection}{\thechapter.\arabic{section}}
\renewcommand{\thesubsection}{\thesection.\arabic{subsection}}

\makeatletter\AtBeginDocument{\renewcommand{\thethmt@dummyctr}{\ifnum\value{subsection}=0 \thesection\else\thesubsection\fi}}\makeatother
\begin{document}
\binoppenalty=10000 \relpenalty=10000 \emergencystretch=1.5em

\numberwithin{equation}{subsection}
\renewcommand{\theequation}{\thesubsection.\arabic{equation}}
\setcounter{tocdepth}{1}

\begin{titlepage}
    \thispagestyle{empty}

    \begin{tikzpicture}[remember picture, overlay]
        \fill[MainColor] (current page.north west) rectangle ([yshift=-3.5cm]current page.north east);

        \fill[MainColor] (current page.south west) rectangle ([yshift=2.5cm]current page.south east);

        \begin{scope}[shift={(current page.center)}, opacity=0.05, MainColor, line width=3pt]
            \draw (0,0) circle (4cm);
            \draw (2,0) circle (4cm);
            \draw (-2,0) circle (4cm);
            \draw (0,2) circle (4cm);
            \draw (0,-2) circle (4cm);
        \end{scope}
    \end{tikzpicture}

    \vspace*{3.5cm}

    \begin{flushleft}
        {\Huge \bfseries \textcolor{MainColor}{Representations of \\[0.3em] $p$-adic Reductive Groups}}\\[1cm]

        \textcolor{MainColor}{\rule{12cm}{3pt}}\\[2cm]

        {\LARGE \bfseries David Renard}\\[0.5cm]
        {\large \itshape École Polytechnique}\\[2cm]
    \end{flushleft}

    \vfill

    \vspace*{2cm}
\end{titlepage}

\cleardoublepage
\thispagestyle{empty}
\null
\vspace*{0.3\textheight}
\begin{flushright}
\Large \itshape To my father, \\
André Renard \\
(1948--2026)
\end{flushright}
\vfill
\cleardoublepage

\tableofcontents

\mainmatter

\cleardoublepage 
\phantomsection  
\addcontentsline{toc}{chapter}{Preface 2026}

\chapter*{Preface 2026}

Since the publication of the first edition, numerous errors in the text have been pointed out to me:
simple typographical errors, embarrassing spelling mistakes, ill-advised copy-pastes, and even some
mathematical problems.

I would like to warmly thank Ammar Yasir Kiliç, Eyal Kaplan, Yuta Nakayama, and Professor Naoki Imai who
pointed out many of these errors to me. I have made available on my webpage

\url{https://perso.pages.math.cnrs.fr/users/david.renard/recherche.html}

\noindent a continuously corrected version, as well as an erratum detailing the changes made.

Given the progress made by AI in mathematical writing, it seemed appropriate to use it to carry out a
complete proofreading of the text and detect remaining errors (specifically Google AI Studio Gemini 3.1
pro-preview). The verdict is damning: the list of typos found is impressively long. It seems justified to
make this new version accessible to as many people as possible, by having it translated into English at
the same time. A second, more systematic proofreading was then carried out with Claude (Anthropic); it
made it possible to correct further typos, but also some mathematical errors, and to complete several
proofs. The writing and any remaining errors are my own. Some proofs, or parts of proofs, were
nevertheless written by AI: the injectivity of the map (\ref{1P}) (Gemini), which seems entirely
convincing to me, unlike the proof I had written trying to follow that of \cite{Ca}, which was
unfortunately completely wrong; and (Claude) the complete proof of Langlands' combinatorial lemma
\ref{comblangl}, the end of the proof of the proposition of \ref{finitude}, the first proof of
(\ref{JHNL}) and some additions to that of the geometric lemma \ref{geomlemma}, the volume estimates and
the final step of the proof of Casselman's criterion \ref{CasCrit}, as well as some arguments in
\ref{langtrip}.

It does not seem useful to me to publish this new version as a book with an official publisher. The
proofreading work done at the time by the reviewers was excellent, and the value added by the publisher
was absolutely zero. In fact, a (small and unimportant) part of the numerous typos is due to the
processing of my \LaTeX{} file by the latter, and running it through a spell checker, which my \LaTeX{}
editor at the time did not offer, could have proved useful. It is highly likely that a new publisher
would not do any better.

I was deeply impressed by the efficiency of the help provided by the AI. I am curious to see the
evolution of mathematical text production in the near future.

\cleardoublepage 
\phantomsection  
\addcontentsline{toc}{chapter}{Preface  to the 2010 Edition}

\chapter*{Preface to the 2010 edition}
The goal of this book is to present, in a more or less complete form, the foundations of the
representation theory of $p$-adic reductive groups. This subject blends a geometric perspective, which
can only be developed through the careful study of explicit examples, and the use of very general results
and principles, falling under what is often called ``abstract nonsense'', i.e., category theory
\footnote{The most important concept in this regard is that of adjunction, which is embodied in
representation theory under the name of ``Frobenius reciprocity''.}. One of the main promoters of the use
of these techniques is J. Bernstein, to whom many fundamental results are due, in particular those around
what is called the \emph{Bernstein center}.

The initial idea was to present it as sketched in the more or less official notes of J. Bernstein
\cite{Be1}, \cite{Be2}, but which has not, to our knowledge, been completely presented in a publishable
form. There are older texts where one can find the foundations of the theory and the proofs of its main
results, due to Bernstein, Casselman, Harish-Chandra, Howe, Jacquet, Zelevinskii (\cite{BeZe1},
\cite{BeZe2}, \cite{Ca}, \cite{Sil}). Then there is the very clear and complete write-up by P. Deligne
\cite{Del} of Bernstein's ideas on the ``center''. However, the texts \cite{Be1} and \cite{Be2},
subsequent to Deligne's, propose a significantly different approach to these results, which seems to us
even more elegant and of greater scope. On the other hand, specialists in the field have sometimes made
their notes on the subject available on the web, giving or completing the proofs of \cite{Be1},
\cite{Be2}, for example \cite{DeB}, \cite{Ro}. But none of these texts covers the entirety of Bernstein's
theory. The need for a book therefore became apparent. Conversations with colleagues confirmed this
impression, which motivated me to undertake this task. My background being in real rather than $p$-adic
groups, I was not the best placed for this work, which is undoubtedly reflected in the imperfections of
this text. But as many know, the best way to learn a subject is to teach it to students, or to write a
book on it. I nevertheless hope that others, novices or specialists, will benefit from it.

Let us now give an overview of the book's content, as well as indications on how to read it. Each chapter
is preceded by an introduction describing its content and is followed by notes where we have tried to
attribute the results and indicated the sources used for writing the proofs. We will therefore be brief
in this general introduction.

The theory of smooth representations of $p$-adic reductive groups is presented in Chapters VI and VII,
which constitute the core of the book. The main tools of this study are the functors of parabolic
induction and restriction, the notion of supercuspidal representation, and their finiteness properties
arising from the fine structure of $p$-adic reductive groups. The category of smooth representations of a
$p$-adic reductive group is studied in Chapter VI, where, among other results, the Bernstein
decomposition theorem and the determination of the center of this category are established. Chapter VII
deals with square-integrable representations, tempered representations, intertwining operators, and
culminates with the Langlands classification theorem.

The preceding chapters place this study in a more general framework: a $p$-adic reductive group is in
particular a totally disconnected topological group. The representation theory of this class of groups is
developed in Chapters III and IV. Chapter II presents general results on the topology of totally
disconnected spaces (functions, distributions, sheaves) and totally disconnected groups (Haar measure,
convolution, Hecke algebra). The language of category theory is used extensively: let $G$ be a totally
disconnected topological group. The main object of study in Chapter III is the category $\caM(G)$ of
smooth (i.e., continuous) representations of $G$ in $\bbC$-vector spaces. We show that this category is
equivalent to that of non-degenerate modules over the Hecke algebra of $G$, which is an algebra with
idempotents. The notions of an algebra with idempotents and of non-degenerate modules for such an algebra
are the subject of Chapter I of the book. Many results on the category $\caM(G)$ then follow immediately
from those obtained in the purely algebraic framework of Chapter I. In particular, the forgetful and base
change functors from module theory yield functors between categories of smooth representations of totally
disconnected groups. Chapter IV studies specific classes of representations: compact, unitary, and
square-integrable modulo the center. Compact representations behave like representations of a compact
group, and their properties are absolutely essential for the rest of the theory.

To go further, one must restrict the class of groups studied to groups possessing more structure.
$p$-adic reductive groups have a very rich structure, which comes on the one hand from the general theory
of reductive algebraic groups, and on the other hand from Bruhat-Tits theory. These structural results
are simply recalled in Chapter V. However, we have tried to present them carefully, introducing the
necessary notation, ordering them, and giving precise references for the unproven results.

This book does not necessarily have to be read linearly. The reader already having some knowledge of the
subject can start reading directly at Chapter VI and refer to the preceding chapters as needed.

The first four chapters require very few mathematical prerequisites, but some familiarity with the
language of category theory is necessary. The rudiments of this theory are given in Appendix A. Chapters
VI and VII are only accessible to readers with a good knowledge of reductive groups, and Chapter V can at
most serve as a guide for the student who will have to refer to the given bibliography. In particular,
the study of examples is indispensable at this stage, and the reader will not find any in this book.

\section*{Acknowledgments}

I would like to warmly thank J. Bernstein for giving his consent to the writing of this book presenting
many (but not all) of his ideas on the representation theory of $p$-adic reductive groups. I would also
like to thank the specialists in the subject who kindly answered my questions during my learning process:
A. Badulescu, G. Chenevier, G. Muic, P. Garrett, and more particularly S. DeBacker, A. Roche, and J-F.
Dat who also agreed to let me reuse, practically word for word, certain proofs from their notes or
articles. Y. Laszlo helped me, with patience, to understand some notions of algebraic geometry, which are
necessary for reading (and therefore writing) this book. An anonymous referee noticed numerous errors,
sometimes embarrassing, and suggested some simplifications of proofs.

Of course, the remaining errors and imperfections are my sole responsibility. Finally, S. Aicardi's
absolute knowledge of the arcana of computer science and the mysteries of \LaTeX{} was invaluable.

 \def\indexter#1{\index[ter]{#1}}

\chapter{Algebras with idempotents} \label{algaidem}

The representation theory of a finite group $G$ (say, with values in complex vector spaces) reduces to
the study of the category $\caM(\bbC[G])$ of unital modules over the group algebra $\bbC[G]$. Multiple
constructions and results then follow. For more general groups, an important idea is to define a ``good''
category of representations to study, i.e., rich enough to contain the representations we are likely to
be interested in, and for which we can nevertheless show that it is equivalent to a category of modules
over a certain $\bbC$-algebra $\caH(G)$, so that the constructions made in the case of finite groups
generalize. It is possible to define such a category when $G$ is the group of rational points of a
connected reductive algebraic group defined over a local field $\bbF$. In the case where $\bbF$ is a
$p$-adic field, we will see in the following chapters the definition of the category of smooth
representations of the group $G$, and the equivalence of this category with a category of modules over a
certain $\bbC$-algebra $\caH(G)$, the \emph{Hecke algebra} of $G$.

Similarly, if $\bbF$ is Archimedean, the category of Harish-Chandra modules is equivalent to a certain
category of modules over the Hecke algebra of $G$ (see \cite{KV}). In both cases, these Hecke algebras
are not unital. This could a priori ruin our strategy, because the constructions for finite groups
crucially use the fact that one considers the category of unital modules over the unital algebra
$\bbC[G]$. It turns out that the Hecke algebras in question are equipped with a structure close to that
given by an identity element. We use the terminology ``algebra with idempotents'' for this structure,
whereas \cite{KV} uses that of an algebra with approximate identity, which perhaps conveys more the idea
of a generalization of the unital case. In this chapter, we study algebras with idempotents, and their
categories of non-degenerate modules (a notion that generalizes that of unital modules over a unital
ring), as well as certain functors between these categories.

\noindent {\bf Notation and Conventions.} We will use the following notation throughout the text:

\medskip

\noindent $\bbZ-\mathbf{mod}$\index[not]{Zmod@$\bbZ-\mathbf{mod}$} for the category of abelian groups,

\noindent $k-\mathbf{Vect}$ \index[not]{kVect@$k-\mathbf{Vect}$} for the category of vector spaces over a
field $k$.

\medskip

If $A$ and $B$ are rings:

\noindent $A-\mathbf{mod}$ \index[not]{Amod@$A-\mathbf{mod}$ } for the category of left $A$-modules,

\noindent $\mathbf{mod}-A$ \index[not]{modA@$\mathbf{mod}-A$} for the category of right $A$-modules,

\noindent $A-\mathbf{mod}-B$ \index[not]{AmodB@$A-\mathbf{mod}-B$ } for the category of $A-B$-bimodules.

\medskip

Recall that $M$ is an $A-B$-bimodule if $M$ is equipped with a left $A$-module structure and a right
$B$-module structure, and the actions of $A$ and $B$ on $M$ commute.

If $V$ is a vector space, we denote by $V^*$ its algebraic dual, and by $\bilo$ the canonical bilinear
form on $V^*\times V$.

Rings are not assumed to be unital unless explicitly stated otherwise. If $A$ is a unital ring
\indexter{unital!(ring)}, we denote by $\una$ \index[not]{$\una$} its identity element. If $A$ is a
unital ring, modules over $A$ are not assumed to be unital a priori. If they are (i.e., if $\una$ acts as
the identity), we specify this by saying that the $A$-module $M$ is \emph{unital}.
\indexter{unital!(module)}

In general, if $e$ is an idempotent of $A$, i.e., such that $e^2=e$ (for example, if $A$ is unital,
$e=\una$), and if $M$ is an $A$-module, then $M$, as a $\bbZ$-module, decomposes as:
\[
M=e\cdot M \oplus (1-e)\cdot M,
\]
where $e\cdot M=\{ e\cdot m \mid m \in M \}$ and $(1-e)\cdot M=\{ m-e\cdot m \mid m \in M \}$. The
idempotent $e$ acts as the identity on $e\cdot M$ and annihilates $(1-e)\cdot M$. The ring $eAe$ is
unital (with identity $e$), and $e\cdot M$ is a unital $eAe$-module.

If $A$ is a unital ring, we denote by \index[not]{M(A)@$\caM(A)$}$\caM(A)$ the full subcategory of
$A-\mathbf{mod}$ whose objects are the unital $A$-modules.

\begin{rmq}
The notation $(1-e)\cdot M$ is merely a notational convenience, and one should be careful that the $1$
here does not denote the identity element of the ring $A$.
\end{rmq}

\section{Rings with idempotents}\label{sec_idem}

\subsection{Order on idempotents}\label{idemgen}

Let $A$ be a ring and let $\scrI=\mathrm{Idem}(A)$ \index[not]{Idem(A)@$\mathrm{Idem}(A)$} denote the set
of idempotent elements of $A$, equipped with the partial order:
\[
e\leq f \quad \text{ if } \quad eAe \subset fAf.
\]

\begin{rmqs}
If $e$ is an idempotent of $A$, then $eAe$ is a subalgebra of $A$ admitting $e$ as its identity.
Furthermore:
\begin{align*}
 eae=a &\Leftrightarrow a \in eAe \quad (a\in A),\\
 eAe&=eA \cap Ae.
\end{align*}
\end{rmqs}

\begin{lemme}
Let $e, f \in \scrI$, and $a\in A$.

\noindent $(i)$ The following conditions are equivalent:

a) $e\leq f$ (i.e., $eAe \subset fAf$).

b) $e \in fAf$.

c) $e=fef$.

Moreover, we then have $e=ef=fe$.

\noindent $(ii)$ If $e\leq f$ and if $ea=a$, then $fa=fea=ea=a$.
\end{lemme}

\begin{proof}
It is clear that $c) \Rightarrow b) \Rightarrow a)$. Suppose that $a)$ holds. Then we can write
$e \in eAe$ as $faf$ for some $a\in A$. We then have $fef=ffaff=faf=e$. Similarly, $fe=e$ and $ef=e$.
Point $(ii)$ follows immediately.
\end{proof}

\subsection{Rings with idempotents: definitions} \label{idem}

\begin{defi}
Let $A$ be a ring (not a priori unital). We say that $A$ is a ring with idempotents if for any finite set
of elements $\{a_i\}$ of $A$, there exists an idempotent $e$ in $A$ ($e^2=e$) such that $a_i=ea_ie$ for
all $i$.
\end{defi}

\begin{lemme}
Let $A$ be a ring with idempotents. The order $\leq$ on $\scrI=\mathrm{Idem}(A)$ is directed, i.e., any
finite family of elements of $\scrI$ admits an upper bound.
\end{lemme}

\begin{proof}
This follows easily from the definitions.
\end{proof}

We can reformulate the above definition by saying that a ring with idempotents is a ring $A$ such that:
\[
A=\varinjlim_{e \in \scrI} \; eAe =
\bigcup_{e\in \scrI} eAe
\]
The inductive limit is formed with respect to the inclusion morphisms $eAe \hookrightarrow fAf$,
$e\leq f$ in $\scrI$ (see \ref{colimites} for the notion of inductive limit).

\begin{defi}
A module $M$ over the ring with idempotents $A$ is said to be \index[ter]{non-degenerate!(module)}
non-degenerate if for every element $m \in M$, there exists an idempotent $e$ of $A$ such that
$e\cdot m=m$. Equivalently, $M$ is non-degenerate if $\displaystyle M=\varinjlim_{e \in \scrI} e\cdot M$.
\end{defi}

\begin{rmqs}
$ $
\begin{itemize}
\item[1.] A unital ring is a ring with idempotents, and its unital modules are exactly the non-degenerate
modules.

\item[2.] We will sometimes need the notion of a \indexter{directed system of idempotents} directed
system of idempotents. In a ring with idempotents, a subset $\scrJ$ of $\scrI$ will be called directed if
any finite family of elements of $\scrI=\mathrm{Idem}(A)$ admits an upper bound in $\scrJ$. For a module
$M$ to be non-degenerate, it suffices that every element of $M$ be fixed by an idempotent in a directed
system $\scrJ$. For such a $\scrJ$, and for any non-degenerate module $M$, we have (see \ref{colimites}
for the definition of inductive limits),
\[
A=\varinjlim_{e \in \scrJ} \; eAe \quad \text{ and }
\quad M=\varinjlim_{e \in \scrJ} \; e\cdot M .
\]

\item[3.] More generally, if $\scrJ$ is a subset of $\mathrm{Idem}(A)$ such that any two idempotents in
$\scrJ$ admit an upper bound in $\scrJ$, one can generalize the notion of a non-degenerate module,
relative to this system $\scrJ$.
\end{itemize}
\end{rmqs}

For any module $M$ over $A$, we denote by \index[not]{MA_@$M_A$}$M_A$ the largest non-degenerate
submodule of $M$. We have
\[
M_A=\varinjlim_{e \in \scrI} \; e\cdot M =
\bigcup_{e \in \scrI} e\cdot M.
\]

We denote by \index[not]{M(A)@$\caM(A)$}$\caM(A)$ the category of non-degenerate (left) modules over the
ring with idempotents $A$. It is a full subcategory of the category $A-\mathbf{mod}$ of all modules over
$A$. One easily verifies that this subcategory is stable under passing to submodules and quotients: it is
therefore an abelian category (cf. \ref{abcat}). In the case where $A$ is a unital ring, $\caM(A)$ is
then the category of unital $A$-modules, which is consistent with the notation introduced previously. We
will also sometimes need to consider right modules. We denote by
\index[not]{M(A)_d@$\caM(A)_r$}$\caM(A)_r$ the category of non-degenerate right $A$-modules.

\begin{prop}
The functor $A-\mathbf{mod}\rightarrow \caM(A)$, $M\mapsto M_A$ is exact.
\end{prop}

\begin{proof}
First, note that if $f\colon M\rightarrow N$ is a morphism in $A-\mathbf{mod}$, it sends $M_A$ into
$N_A$. This defines the functor on morphisms by restriction to the non-degenerate part. If
\[
L \stackrel{f}{\rightarrow} M
\stackrel{g}{\rightarrow} N
\]
is an exact sequence in $A-\mathbf{mod}$, it is clear that $\im(f_{|L_A})\subset \ker(g_{|M_A})$.
Conversely, if $m \in \ker g \cap M_A$, there exists $l\in L$ such that $f(l)=m$. Let $e$ be an
idempotent of $A$ fixing $m$. We then have $f(l)=m=e\cdot m=e\cdot f(l)=f(e\cdot l)$, and since
$e\cdot l\in L_A$, we indeed have $m \in \im(f_{|L_A})$. This shows that $\im(f_{|L_A})=\ker(g_{|M_A})$.
\end{proof}

\subsection{Topology of $A$. Completion} \label{completeA}
\index[ter]{completion!of a ring with idempotents}

We refer the reader to \cite{BourTG} or \cite{Dug} for the generalities on topological spaces used in
this section, such as the notions of a family of pseudometrics, uniform structure, completion, etc.

Let $A$ be a ring with idempotents and let $e \in \scrI=\mathrm{Idem} (A)$. As a $\bbZ$-module, we have a
decomposition
\begin{align*}
 & A=A(1-e) \oplus Ae, \\
\text{ where } \quad &A(1-e)=\{ a-ae, \, a\in A \}
=\{ a\in A\, \vert \, ae=0 \} \\
\text{ and } \quad &Ae=\{ ae, \, a\in A \}
=\{ a\in A\, \vert \, ae=a \}.
\end{align*}

There exists a natural topology on $A$: a neighborhood basis of $0$ is given by the subsets $A(1-e)$,
$e \in \scrI$. Indeed, if $e\leq f$, we have $e=fef=fe$, hence $(a-af)e=0$ for all $a\in A$, and
therefore:
\[
a-af=(a-af)- (a-af)e \in A(1-e),
\]
which shows that $A(1-f) \subset A(1-e)$. Since any finite set of idempotents admits an upper bound, we
see that any finite intersection of subsets of the form $A(1-e)$ contains a subset of this form. We may
therefore choose the subsets $A(1-e)$, $e\in \scrI$, as a neighborhood basis of $0$. A neighborhood basis
of the element $b\in A$ is obtained by the translation $a\mapsto b+a$. This equips $A$ with a uniform
structure. An equivalent way to view this is to introduce on $A$ the family of pseudometrics $d_e$,
$e \in \scrI$, where
\[
d_e\colon A\times A \rightarrow \{0,1\},
\]
\[
d_e(a,b)=0 \text{ if } a-b\in A(1-e), \;
d_e(a,b)=1 \text{ otherwise. }
\]
It is immediate to verify that the $d_e$ satisfy the axioms of pseudometrics, that the family of
pseudometrics $(d_e)_{e \in \scrI}$ is directed (i.e., if $e\leq f$, then $d_e \leq d_f$), and that the
topology defined above is the one defined by this family of pseudometrics.

We will now construct the completion \index[not]{Abar@$\overline{ A}$} $\overline{ A }$ of the
topological space $A$ equipped with this uniform structure. The situation studied here is close to the
one where $A$ is a unital ring, $J$ is an ideal of $A$, and $A$ is equipped with the left-translation
invariant topology where a neighborhood basis of $0$ is given by $\{J^n\}_{n\in \bbN}$. We then know that
$\overline{A}$ is given by the projective limit of the $A/J^n$. Similarly, in our case, if $e$ and $f$
are two idempotents of $A$ such that $e\leq f$, let us introduce the transition morphisms
$A/ A(1-f)\rightarrow A/A(1-e)$ which will be used to form the projective limit (cf. \ref{limites}).
Since $A=A(1-e) \oplus Ae$, we have $A/A(1-e)\simeq Ae$. We therefore seek to define morphisms
$Af\rightarrow Ae$ when $e\leq f$. It is natural to take:
\[
Af \rightarrow Ae,\quad b\mapsto be.
\]
Let us now form the projective limit $\displaystyle \overline{A}=\varprojlim Ae$ for these transition
morphisms. Explicitly, $\overline{A}$ is the set of families $\{\bar a(e)\}_{e\in \scrI}$, such that for
any $e$, $f$ in $\scrI$ with $e\leq f$, we have $\bar a(e)\in Ae$ and $\bar a(f)e=\bar a(e)$. The space
$\overline{A}$ is equipped with the left-translation invariant topology such that a neighborhood basis of
$0$ is given by the set of subsets of the form
\[
\overline{A}(1-e):=\{ \bar a \in \overline{A}\, \mid \,
\bar a(e)=0\},
\]
with $e\in \scrI$. It is clear that $A$ embeds canonically and continuously into $\overline{A}$ via
$a\mapsto \bar a$, $\bar a(e)=ae$. We equip $\overline{A}$ with a ring structure extending that of $A$.
If $\bar a, \bar b \in \overline{A}$, $e\in \scrI$, we set
\[
(\bar a\bar b)(e)= \bar a(f)\bar b(e)
\]
where $f\in \scrI$ is such that $f\bar b(e)=\bar b(e)$. One easily verifies that this is well-defined
(i.e., does not depend on the choices made). The ring $\overline{A}$ is then unital, with identity
$\mathbf{1}_{\overline{A}}$, where $\mathbf{1}_{\overline{A}}(e)=e$ for all $e$ in $\scrI$.

It is also clear that $A$ is dense in $\overline{A}$. Indeed, for any $\bar a \in \overline{A}$, for any
$e \in \scrI$, for any $f \geq e$, we have
\[
(\bar a- \bar a(f))e= \bar a(e)-\bar a(f)e=
\bar a(e)-\bar a(e)=0,
\]
and therefore $\bar a- \bar a(f) \in \overline{A}(1-e)$.

We have not shown that $\overline{A}$ is a complete uniform space, the completion of $A$, but we will not
need this in what follows, and we therefore leave it to the reader. On the other hand, the fact that the
action of $A$ on any non-degenerate $A$-module $M$ extends to an action of $\overline{A}$ is absolutely
crucial. If $e\cdot m=m$, $e\in \scrI$, $m\in M$, we set $\bar a\cdot m=\bar a(e) \cdot m$,
$\bar a \in \overline{A}$. This does not depend on the choice of $e$ and equips $M$ with a unital
$\overline{A}$-module structure. In other words, if $M$ is equipped with the discrete topology, the
uniformly continuous map $A \rightarrow M$, $a \mapsto a\cdot m$ extends by continuity to $\overline{A}$.

\subsection{Completion of a non-degenerate module} \label{compmod}
\index[ter]{completion!of a non-degenerate module}

Let us generalize the constructions of the previous section, by defining the completion
\index[not]{Mbar@$\overline{M}$}$\overline{M}$ of a non-degenerate module $M$ over the ring with
idempotents $A$. The proofs of the following assertions easily adapt from those of the previous section.
We equip $M$ with the neighborhood basis of $0$ given by the subsets of the form:
\[
(1-e)\cdot M:= \{ m-e\cdot m, \, m \in M \}
\]
where $e$ runs through $\scrI=\mathrm{Idem}(A)$. A neighborhood basis of any element $m\in M$ is obtained
by translation $m' \mapsto m+m'$ in $M$. We obtain a uniform structure on $M$. This topology on $M$ is
given by the family of pseudometrics $d_e$, $e \in \scrI= \mathrm{Idem}(A)$,
\[
d_e\colon M \times M \rightarrow \{0,1\},
\]
\[
d_e(m,m')=0 \text{ if } m-m'\in (1-e)\cdot M,\;
d_e(m,m')=1 \text{ otherwise.}
\]

The completion $\overline{M}$ of $M$ for this uniform structure is then the projective limit
\[
\overline{M} = \varprojlim_{e\in \scrI} \; e\cdot M
\]
for the transition morphisms
\[
f\cdot M \rightarrow e \cdot M~, \quad
m \mapsto e\cdot m, \quad (e\leq f).
\]
One can view $\overline{M}$ as the set of families $\{ \bar m(e)\}_{e\in \scrI}$ such that for any $e,f$
in $\scrI$ with $e \leq f$, we have $\bar m(e) \in e\cdot M$ and $e\cdot \bar m (f)=\bar m(e)$. Moreover,
$M$ embeds into $\overline{M}$ via
\[
m \mapsto \bar m, \quad \bar m(e)=e\cdot m, \quad
(e\in \scrI).
\]

\begin{prop}
With the above notation, the morphism
\[
\Hom_A(A,M) \rightarrow \overline{M}, \quad
\phi \mapsto \{\phi(e)\}_{e \in \mathrm{Idem}(A)}
\]
is an isomorphism.
\end{prop}

\begin{proof}
Let $\phi \in \Hom_A(A,M)$ and let $a \in A$. Let $e$ be an idempotent of $A$ such that $ae=a$. We then
have $\phi(a)=\phi(ae)=a\phi(e)$, so $\phi$ is determined by the family $\{ \phi(e) \}_{e\in \scrI}$ and
it is clear that if $e \leq f$, then $e\cdot \phi(f)=\phi(e)$. Conversely, if we have a family
$\{ \bar m(e) \}_{e\in \scrI}$ of elements of $M$, we can define $\phi \in \Hom_A(A,M)$ by
$\phi(a)=a \cdot \bar m(e)$, $e\in \scrI$, $ae=a$ provided that $\{ \bar m(e) \}_{e\in \scrI}$ is an
element of $\overline{M}$.
\end{proof}

\medskip

The completion $\overline{M}$ is equipped with an $A$-module structure, apparent in the above
description, since $\Hom_A(A,M)$ is naturally a left $A$-module:
\[
(a\cdot \phi) (b)=\phi (ba), \quad (a,b \in A), \;
\phi \in \Hom_A(A,M)
\]
This structure extends that of $M$, but the module $\overline{M}$ is no longer necessarily
non-degenerate.

\begin{rmq}
In the previous section, these constructions were made by considering $A$ as a right module over itself.
Since the constructions of this section obviously adapt to right modules, we see that
$\overline{A} \simeq \End_{A^\circ}(A)$ (the notation $\Hom_{A^\circ}$ is there to remind us that we are
considering right module structures).
\end{rmq}

\subsection{Smoothness/completion adjunction} \label{lisscomp}

Let $A$ be a ring with idempotents and $\scrI= \mathrm{Idem}(A)$. We saw in section \ref{idem} that the
``smoothness'' functor
\[
A-\mathbf{mod}\rightarrow \caM(A), \quad
M \mapsto M_A
\]
is an exact functor. In the previous section, we defined a ``completion'' functor
\[
\caM(A)\rightarrow A-\mathbf{mod}, \quad
N \mapsto \overline{N}.
\]
This functor being defined in an obvious way on morphisms by continuous extension, i.e., if
$f \in \Hom_A(M,N)$ is a morphism in $\caM(A)$, it extends to a morphism
\[
\bar f \colon \overline{M} \rightarrow \overline{N}, \quad \bar f(\bar m)(e)=f(\bar m(e)), \quad (e \in \scrI).
\]

\begin{prop}
For any module $M \in \caM(A)$, we have $(\overline{M})_A= M$. The completion functor is faithful. It is
the right adjoint of the smoothness functor.
\end{prop}

\begin{proof}
Let us explicitly give the form of the natural inclusion of $M$ into $\overline{M}$ when $\overline{M}$
is identified with $\Hom_A(A,M)$ by Proposition \ref{compmod}:
\begin{equation} \label{HAAM}
M \hookrightarrow \Hom_A(A,M), \quad m \mapsto \phi_m,
\end{equation}
where $\phi_m(a)=a\cdot m$, for all $a \in A$. Recall that the $A$-module structure on $\Hom_A(A,M)$ is
given by:
\[
(a\cdot \phi)(a')=\phi(a'a), \quad (a,a' \in A),
(\phi \in \Hom_A(A,M)).
\]
Suppose that $\phi \in \Hom_A(A,M)$ is fixed by an idempotent $e$ of $A$. Then for all $a\in A$,
$a\cdot \phi(e)=\phi(ae)=(e\cdot \phi) (a)=\phi(a)$, and therefore $\phi$ is in the image of
(\ref{HAAM}). We deduce that $M$ identifies with the non-degenerate part of $\overline{M}=\Hom_A(A,M)$.
We note the canonical isomorphism
\begin{equation}\label{MbarA}
M\simeq \Hom_A(A,M)_A \simeq \overline{M}_A .
\end{equation}

If $f \in \Hom_A(M,N)$ is a morphism in $\caM(A)$ such that $\bar f=0$, then $f=0$. This shows that the
completion functor is faithful. We will now establish the adjunction of the two functors. For
generalities concerning adjoint functors, we refer to Appendix \ref{adjfonct}. We must exhibit a natural
isomorphism
\[
\Hom_A(N ,\overline{M}) \simeq \Hom_A(N_A,M)
\]
for any $M \in \caM(A)$ and any $N \in A-\mathbf{mod}$. Since the image of a non-degenerate module is
non-degenerate, we have
\[
\Hom_A(N_A,\overline{M})=
\Hom_A(N_A,(\overline{M})_A)=\Hom_A(N_A,M).
\]
On the other hand, any morphism $\phi \in \Hom_A(N,\overline{M})$ defines by restriction a morphism
$\phi_A$ in $\Hom_A(N_A,\overline{M})=\Hom_A(N_A,M)$. If $\phi_A=0$, we have
$\phi(e\cdot n)=0=e \cdot \phi(n)$ for all $n \in N$ and all $e\in \scrI$, and therefore $\phi$ is zero,
because $M$ is non-degenerate. This restriction map is therefore injective. It is surjective, because if
$\phi_A\in \Hom_A(N_A,M)$, we define $\phi$ in $\Hom_A(N,\overline{M})$ by $\phi(n)=\bar m$ where
$\bar m(e)= \phi_A(e \cdot n)$ for all $e \in \scrI$. One easily verifies that the restriction of $\phi$
to $N_A$ is $\phi_A$. The naturality of these constructions is clear.
\end{proof}

\begin{rmq}
If $M$ is a non-degenerate $A$-module, and $N$ a submodule, then $\overline{N}$ is the closure of $N$ in
$\overline{M}$, i.e.,
\begin{equation}\label{compsousmod}
\overline{N}=\{ \bar n \in \overline{M} \, \mid \,
\bar n(e) \in e\cdot N, \, (e \in \scrI) \}
=\{ \bar n \in \overline{M} \, \mid \,
A\cdot \bar n \subset N \}.
\end{equation}
\end{rmq}

\subsection{Other descriptions of $\overline{A}$} \label{autredescrA}

Since $A$ is naturally, by left multiplication, a non-degenerate left $A$-module, it is also a unital
left $\overline{A}$-module, and any element $\bar a\in \overline{A}$ defines an endomorphism of abelian
groups of $A$ which commutes with right multiplications. One can therefore identify (cf. remark
\ref{compmod}) $\overline{A}$ with the subgroup $\End_{A^\circ}(A)$ of $\End_\bbZ(A)$.

Let us equip $\End_\bbZ(A)$ with the topology of pointwise convergence: it is a topological abelian group
(for addition) where a neighborhood basis of $0$ is given by the family of open sets
\[
\caO(U_1,\ldots,U_n,a_1,\ldots,a_n):=
\{ \beta \in \End_\bbZ(A) \, \mid \, \beta(a_i)\in U_i \}
\]
where the $a_i$ are elements of $A$, and the $U_i$ are open neighborhoods of $0$ in $A$. One can take the
$U_i$ of the form $U_i=A(1-e_i)$ with $e_i \in \scrI$. The induced topology on $\overline{A}$ then
coincides with the one defined previously. Indeed
\[
\caO(U_1,\ldots,U_n,a_1,\ldots,a_n)\cap \overline{A}=
\{ \bar a \in \overline{A} \, \mid \, \bar a a_i \in U_i \},
\]
so if the $U_i$ are of the form $U_i=A(1-e_i)$ with $e_i \in \scrI$, and if we take $f \in \scrI$ such
that $f a_i = a_i$ for all $i$, then for all $\bar a \in \overline{A}(1-f)$, we have
$\bar a a_i=\bar a f a_i = 0$, which shows that
$\overline{A}(1-f) \subset \caO(U_1,\ldots,U_n,a_1,\ldots,a_n)\cap \overline{A}$. Conversely,
$\overline{A}(1-f)=\caO( A(1-f),f)\cap \overline{A}$.

It is clear that $\overline{A}$ is contained in the closure, for this topology of pointwise convergence,
of the set of left multiplications by elements of $A$. Conversely, this closure is of course itself
contained in the commutant of right multiplications, i.e., in $\End_{A^\circ}(A)=\overline{A}$.

We have thus obtained the following lemma.

\begin{lemme}
The ring $\overline{A}$ identifies with the closure in $\End_\bbZ(A)$, for the topology of pointwise
convergence, of the set of left multiplications by elements of $A$.
\end{lemme}

One can view $\overline{A}$ in yet another way. Let
\[
\varpi:\caM(A)\rightarrow \bbZ-\mathbf{mod}
\]
be the forgetful functor.

\begin{prop}
The completion $\overline{A}$ of $A$ identifies with $\Hom(\varpi ,\varpi)$, the ring of endomorphisms of
the forgetful functor $\varpi$, via the map which associates to a natural transformation $\theta$ of the
forgetful functor to itself the element $\theta_A \in \End_\bbZ(A)$.
\end{prop}

\begin{proof}
It may be useful to make the right-hand side of this equality explicit. It is a $\Hom$ in the category of
additive functors from $\caM(A)$ to $\bbZ-\mathbf{mod}$, and it therefore naturally admits a ring
structure. Let $\theta$ be a natural transformation from $\varpi$ to itself. Let $f$ in $\Hom_A(A,A)$ be
a morphism given by a right multiplication in $A$. From the definition of a natural transformation, we
see that $\theta_A$ commutes with $f$, and therefore, according to the previous lemma, this gives us an
element $\bar a$ of $\overline{A}$. Conversely, if $\bar a \in \overline{A}$, and if $M \in \caM(A)$, we
take for $\theta_M$ the action of the element $\bar a$ on $M$. One easily verifies that this defines a
natural transformation from $\varpi$ to itself. The two constructions are inverse to each other.
\end{proof}

\subsection{Center of $\caM(A)$}
\index[ter]{center!of $\caM(A)$} \label{idemmod}

We refer the reader to \ref{centrecat} for generalities on the center of a category. Let $A$ be a ring
with idempotents. We propose to make explicit the center of the category $\caM(A)$, i.e., the set of
natural transformations of the identity functor of the category to itself. We consider $A$ as a left
module over itself.

\begin{lemme}
Let $\phi\in \overline{A}=\End_{A^\circ}(A)$ be a morphism also commuting with the action by left
multiplication of $A$ on itself, i.e.,
\[
\phi(abc)=a\phi(b)c, \quad (a,b,c\in A).
\]
Then for any non-degenerate $A$-module $M$, there exists a morphism
$\tilde \phi_M \colon M \rightarrow M$ such that:
\[
\tilde \phi_M(a\cdot m)= \phi(a)\cdot m, \quad (a\in A),\, (m\in M).
\]
Moreover, $\phi \mapsto \tilde \phi_M$ is natural in $M$, i.e., if $M$ and $N$ are non-degenerate
$A$-modules, then for any morphism $f \in \Hom_A(M,N)$, we have
$\tilde \phi_N \circ f=f \circ \tilde \phi_M$.
\end{lemme}

\begin{proof}
Since $M$ is non-degenerate, for any $m\in M$, let us choose an idempotent $e$ of $A$ such that
$e\cdot m=m$. We then set $\tilde \phi_M(m)=\phi(e)\cdot m$. This does not depend on the choice of $e$,
because if $e$ and $f$ satisfy $e\cdot m=f \cdot m=m$, we can, according to Lemma \ref{idem}, find an
upper bound $h$ for $e$ and $f$. Since $e=he$ and $f=hf$, we then have:
\[
\phi(e)\cdot m=\phi(he)\cdot m=(\phi(h)e)\cdot m=
\phi(h)\cdot m= \phi(f)\cdot m.
\]
Let $e$ be an idempotent fixing $a\cdot m$, and let $f$ be an idempotent fixing $a$, and let us take $h$
as an upper bound for $e$ and $f$. We then have:
\begin{align*}
\tilde \phi_M(a\cdot m)&= \phi(e)\cdot (a\cdot m)=
\phi(he)\cdot(a \cdot m)= (\phi(h)e)\cdot(a \cdot m)\\
&=\phi(h)\cdot(a \cdot m)
= \phi(ha)\cdot m=\phi(hfa)\cdot m=\phi(fa)\cdot m\\
&= \phi(a)\cdot m.
\end{align*}
The rest of the lemma is easily proved.
\end{proof}

In other words, any $\phi\in \overline{A}$ commuting with the action by left multiplication of $A$ on
itself defines a natural transformation of the functor $\caI d_{\caM(A)}$ to itself. Conversely, such a
natural transformation $\theta$ determines a morphism
$\phi=\theta_A \in \overline{A}=\End_{A^\circ } (A)$ commuting with left multiplications. These two
constructions are inverse to each other.

\begin{thm}
The center of $\caM(A)$ identifies with the commutant of $A$ in $\overline{A}$, or equivalently with the
center of $\overline{A}$. It is the projective limit over $\scrI=\mathrm{Idem}(A)$ of the centers of the
rings $eAe$, $e \in \scrI$.
\end{thm}

\begin{proof}
The center of $\caM(A)$ identifies with the subgroup of $\overline{A}$ of elements commuting with the
elements of $A$, or by continuity, with the center of $\overline{A}$. We verify the last assertion. The
transition morphisms
\[
fAf \rightarrow eAe, \quad a \mapsto eae, \quad
(e,f \in \scrI, \; e\leq f),
\]
induce transition morphisms between the centers of these rings. Indeed, if $z$ is in the center of $fAf$,
it commutes with $e$, and we have for all $x \in eAe$,
\[
[eze,x]=ezex-xeze=ezx-xze=zex-xez=zx-xz=0.
\]
This shows that the image of $z$ by the transition morphism is central in $eAe$. This defines the
projective system whose limit we want to identify.

If $\bar a \in \overline{A}$ is central, then it is clear that for all $e \in \scrI$,
$e \bar a e=\bar a(e)$ is in the center of $eAe$. For the converse, we first show that
$e\bar a =\bar a e$ for all $e \in \scrI$. Let us take $f \geq e$. We have
\[
( e\bar a -\bar a e)f =e\bar a f - \bar a e f =
e(f \bar a f)- (\bar a f) e = e \bar a (f) -
\bar a (f) e =0
\]
because $\bar a (f) = f \bar a f$ is central in $fAf$, and in particular commutes with $e$. This being
established, for $\bar a$ to commute with $eAe$, it is necessary and sufficient that
$\bar a e=e \bar a e=\bar a (e)$ be in the center of $eAe$, and we conclude by the fact that $A$ is the
union of the $eAe$.
\end{proof}

\medskip

We denote by \index[not]{Z(A)@$\frZ(A)$}$\frZ(A)$ the center of the category $\caM(A)$. It is clear from
the above that $\frZ(A)$ is equipped with a ring structure.

\begin{rmq}
When $A$ is unital, the center $\frZ(A)$ identifies with the center of the ring $A$, a well-known result.
\end{rmq}

\subsection{The center as an algebra of functions on irreducibles} \label{centF}

Let $A$ be a $\bbC$-algebra with idempotents, of finite or countable dimension. Let $\mathbf{Irr}(A)$
\index[not]{Irr(A)@$\mathbf{Irr}(A)$} denote the set of isomorphism classes of simple modules in
$\caM(A)$. Schur's lemma \ref{ASchur} shows that any element $z$ of the center $\frZ(A)$ acts by a scalar
operator on simple modules. For any simple $A$-module $M$, and any element $z \in \frZ(A)$, let $z_M$
denote this scalar. It obviously depends only on the class of $M$ in $\mathbf{Irr}(A)$. Let
$\mathrm{Funct}(\mathbf{Irr}(A))$ denote the algebra of functions on $\mathbf{Irr}(A)$ with values in
$\bbC$. We have highlighted an algebra morphism:
\begin{equation} \label{centrefonctions}
\frZ(A) \rightarrow \mathrm{Funct}(\mathbf{Irr}(A)),
\quad z \mapsto [M \mapsto z_M].
\end{equation}

In good situations, this morphism is injective. Let us introduce the following separation property:

({\bf Sep}): for any non-zero $a \in A$, there exists a non-degenerate simple module $M$ and $m \in M$
such that $a\cdot m\neq 0$.

Equivalently, any non-degenerate $A$-module $M$ is given by a morphism \footnote{It is explained in the
next section how a module in $\caM(A)$ is naturally equipped with a $\bbC$-vector space structure which
makes the action of $A$ linear.} $ \pi_M \colon  A \rightarrow \End_\bbC(M)$, and the property ({\bf
Sep}) says that the algebra morphism:
\begin{equation} \label{plongEnd}
A \rightarrow \prod_{M \in \mathbf{Irr}(A)} \End_\bbC (M)
\end{equation}
is injective.

The extension of (\ref{plongEnd}) to $\overline{A}$ remains injective and the same is therefore true of
its restriction to $\frZ(A)$.

In the case where $A$ satisfies the property ({\bf Sep}), any element of $\frZ(A)$ is therefore
characterized by the function on $\mathbf{Irr}(A)$ that it defines. As we will see, this applies to the
case where $A$ is the Hecke algebra of a totally disconnected topological group $G$. When $G$ is a
$p$-adic reductive group, Bernstein's description of the center is then essentially the determination of
the image of (\ref{centrefonctions}).

\section{Forgetful functors and their adjoints}\label{AoubliB}

\subsection{Vector space structure on non-degenerate modules} \label{AMODUn}
Let $A$ be a unital $\bbC$-algebra, with identity element $\una$. If $M$ is a unital $A$-module, it is
canonically equipped with a $\bbC$-vector space structure by:
\begin{equation*}
 \lambda m:= (\lambda \una)\cdot m ,\qquad (m\in M), \, (\lambda \in \bbC).
\end{equation*}
Moreover, we have: $(\forall m\in M), \, (\forall \lambda \in \bbC), (\forall a \in A)$,
\begin{equation}\label{lam}
 \lambda(a\cdot m)= (\lambda a )\cdot m=a \cdot(\lambda m).
\end{equation}

Suppose that $A$ is a $\bbC$-algebra with idempotents. It is still possible to equip any non-degenerate
$A$-module $M$ with a $\bbC$-vector space structure by: $(\forall m\in M)$, $(\forall \lambda \in \bbC)$,
$(\forall e \in \mathrm{Idem}(A) \mid e\cdot m=m )$,
\begin{equation*}
 \lambda m:= (\lambda e)\cdot m.
\end{equation*}
One easily verifies that this does not depend on the choice of $e$, and that property (\ref{lam}) is
still satisfied.

Suppose that $X$ and $Y$ are two non-degenerate $A$-modules. From the above, they are also two
$\bbC$-vector spaces. It is then immediate that any morphism in $\Hom_{A}(X,Y)$ is also $\bbC$-linear. In
particular, $\Hom_{A}(X,Y)$ is a $\bbC$-vector space: for any $\phi$ in $\Hom_{A}(X,Y)$ and for any
$\lambda$ in $\bbC$,
\begin{equation*}
\lambda \phi \colon  x \mapsto  \lambda(\phi(x))=\phi(\lambda x).
\end{equation*}
Similarly, if $X$ is a right $A$-module, and $Y$ a left $A$-module, both non-degenerate, the tensor
product $X \otimes_A Y$ inherits a vector space structure and $X \otimes_A Y$ is a quotient of
$X \otimes_\bbC Y$.

If $A$ is a non-unital $\bbC$-algebra, there is no canonical $\bbC$-vector space structure on
$\Hom_{A}(X,Y)$ or on $X \otimes_A Y$, even when $X$ and $Y$ are vector spaces. This is extremely
inconvenient, and we remedy it by introducing the unital algebra
\[
\widetilde A= A \oplus \bbC\cdot \una,
\]
where we decide that $\una$ is the identity of $\widetilde A$ (this determines a unique algebra structure
on $\widetilde A$). If $X$ and $Y$ are two $A$-modules and also $\bbC$-vector spaces, they naturally
become unital $\widetilde A$-modules and we have
\begin{align*}
&\Hom_{\widetilde A}(X,Y)= \Hom_{A}(X,Y) \cap \Hom_{\bbC}(X,Y), \\
&X\otimes_{\widetilde A} Y=  X \otimes_A Y/\mathrm{Span}(\{ \lambda x\otimes y- x\otimes \lambda y, \; x\in X, y \in Y, \lambda \in \bbC\}  ).
\end{align*}
where $\mathrm{Span}(P)$ denotes the subspace spanned by a subset $P$ of a vector space. If $A$ is an
algebra with idempotents, from the above, if $X$ and $Y$ are two non-degenerate $A$-modules,
\begin{equation} \label{adndeg}
 \Hom_{\widetilde A}(X,Y)= \Hom_{A}(X,Y), \quad  X \otimes_{\widetilde A} Y=  X \otimes_A Y
\end{equation}
for the canonical $\bbC$-vector space structures on $X$ and $Y$. This is true in particular if $A$, as
well as the modules $X$ and $Y$, are unital.

\subsection{Associativity formulas} \label{KVC20}

When $A$ is a unital algebra, and $X$ is a module in $\caM(A)$, we can define the functors:
\[
\bullet \otimes_A X: \caM(A)_r \longrightarrow \mathbf{\bbC -Vect}, \qquad Y \mapsto Y \otimes_A X,
\]
\[
\Hom_A(\, \bullet \, ,X): \caM(A) \longrightarrow  \mathbf{\bbC -Vect}, \qquad Y \mapsto \Hom_A(Y, X)
\]
and
\[
\Hom_A(X,\, \bullet \, ): \caM(A) \longrightarrow  \mathbf{\bbC -Vect}, \qquad Y \mapsto \Hom_A(X,Y)
\]
Their properties are well known (see for example \cite{Faisc}). By (\ref{adndeg}), these results
generalize to the case where $A$ is an algebra with idempotents.

\begin{prop}
Let $A$ be a $\bbC$-algebra with idempotents, and $X$ a non-degenerate $A$-module.

The functor $ \bullet \otimes_A X$ is covariant and right exact, the functor $\Hom_A( \, \bullet \, ,X)$
is contravariant and left exact, the functor $\Hom_A(X, \, \bullet \, )$ is covariant and left exact.
\end{prop}

We refer the reader to \cite{Knapp88}, pages 218-219 for a proof.

\medskip

Let $R,S$ be two $\bbC$-algebras. We will write $X^R$, ${}^S Y^R$, $Z^R$ to say that $X$ is a right
$R$-module, $Y$ an $S-R$-bimodule (left over $S$ and right over $R$), etc. The $R$-modules naturally
become unital $\tilde R$-modules. The same applies to $S$.

\begin{thm}
Let $R,S$ be two $\bbC$-algebras and $X^R$, ${}^R Y^S$, $Z^S$ be modules. Then we have a natural
isomorphism:
\begin{align}\label{AD1}
\Hom_{\widetilde S}(X\otimes_{\widetilde R}Y,Z)&\simeq  \Hom_{\widetilde R}(X,\Hom_{\widetilde S}(Y,Z))\\
\phi &\mapsto \psi, \quad \psi(x)(y)=\phi(x \otimes y).
\end{align}

Similarly, for any modules ${}^R X$, ${}^S Y^R$ and ${}^S Z$, we have a natural isomorphism:
\begin{equation} \label{AD2}
\Hom_{\widetilde S}(Y\otimes_{\widetilde R}X,Z)\simeq \Hom_{\widetilde R}(X,\Hom_{\widetilde S}(Y,Z)).
\end{equation}

Finally, for any modules $X^R$, ${}^R Y^S$ and ${}^S Z$, we have a natural isomorphism:
\begin{equation}  \label{AD3}
(X\otimes_{\widetilde R}Y)\otimes_{\widetilde S} Z \simeq X\otimes_{\widetilde R}(Y\otimes_{\widetilde S} Z).
\end{equation}
\end{thm}

A proof can be found in \cite{Knapp88}, pages 97-98.

\begin{cor}
Let $R$, $S$ be $\bbC$-algebras with idempotents. For any non-degenerate modules $X^R$, ${}^RY^S$ and
$Z^S$, we have a natural isomorphism:
\begin{equation}\label{ADId1}
\Hom_{S}(X\otimes_{R}Y,Z) \simeq\Hom_{R}(X,\Hom_{S}(Y,Z)_R).
\end{equation}
For any non-degenerate modules ${}^R X$, ${}^S Y^R$ and ${}^S Z$, we have a natural isomorphism:
\begin{equation} \label{ADId2}
\Hom_{S}(Y\otimes_{R}X,Z) \simeq \Hom_{R}(X,\Hom_{ S}(Y,Z)_R).
\end{equation}

Finally, for any non-degenerate modules $X^R$, ${}^R Y^S$ and ${}^S Z$, we have a natural isomorphism:
\begin{equation}  \label{ADId3}
(X\otimes_{R}Y)\otimes_{S} Z\simeq  X\otimes_{R}(Y\otimes_{S} Z).
\end{equation}
\end{cor}

\begin{proof}
This follows immediately from the theorem, from (\ref{adndeg}) and, for example for the second formula,
from the fact that
\[
\Hom_{R}(X,\Hom_{S}(Y,Z)_R)= \Hom_{R}(X,\Hom_{S}(Y,Z))
\]
because the image of a non-degenerate module is non-degenerate.
\end{proof}

\begin{rmq}
As a consequence of the theorem and its corollary, we obtain the well-known adjunction of the $\Hom$ and
$\bigotimes$ functors. But the theorem gives the naturality of the isomorphisms in the three variables
$X$, $Y$ and $Z$.
\end{rmq}

\subsection{Forgetful and induction functors} \label{Oublietadjoints}
\indexter{functor!forgetful} \indexter{functor!induction} Let $A$ and $B$ be unital $\bbC$-algebras, and
$\phi \colon A \rightarrow B$ a morphism of unital algebras. The ``forgetful functor''
\[
\phi_* \colon \caM(B) \rightarrow \caM(A)
\]
which transforms a unital $B$-module $N$ into a unital $A$-module by
\[
a\cdot n:= \phi(a)\cdot n,\quad (a\in A), (n \in N),
\]
admits left and right adjoints, respectively:
\[
{}^*\phi,\phi^* \colon   \caM(A)  \rightarrow  \caM(B)
\]
which we will describe. Note that by $\phi$, $B$ is equipped with an $A$-bimodule structure. The left
adjoint ${}^*\phi$ is then given by the tensor product
\[
{}^*\phi \colon  M \mapsto B\otimes_A M
\]
and the right adjoint $\phi^*$ by the $\Hom$ functor
\[
\phi^* \colon  M \mapsto \Hom_A(B,M).
\]

Note also that we adopt the conventions of algebraic geometry: $\phi \mapsto \phi_*$ is contravariant,
but despite this we put the star as a subscript rather than a superscript. This is because from a
geometric point of view we consider $\phi$ as a morphism between the varieties $\mathrm{Spec}(B)$ and
$\mathrm{Spec}(A)$, which restores covariance.

On the other hand, in the construction of ${}^*\phi$ (resp. $\phi^*$), we only use the right $A$-module
structure on $B$ induced by $\phi$ and the natural left module structure of $B$ on itself (resp. the left
$A$-module structure on $B$ induced by $\phi$ and the natural right module structure of $B$ on itself)
and the fact that the left and right actions commute. It is therefore interesting to describe the
forgetful functor $\phi_*$ in terms of these a priori weaker structures. Now, for any unital $B$-module
$N$, we have the canonical isomorphisms
\[
N \simeq  \Hom_B(B,N)  \simeq B\otimes_B N,
\]
and the construction of $\Hom_B(B,N)$ (resp. $B\otimes_B N$) as an $A$-module only involves the right
$A$-module and left $B$-module structures (resp. left $A$-module and right $B$-module structures) of $B$
and the fact that the left and right actions commute. The adjunctions between $\phi_*$ and
${}^*\phi,\phi^*$ are then consequences of the second identity of Theorem \ref{KVC20}.

When $A$ and $B$ are no longer assumed to be unital, the constructions made above no longer work. In the
situation that interests us, where $A$ and $B$ are $\bbC$-algebras with idempotents, we use the corollary
of Theorem \ref{KVC20} to restore the situation. The context is the following: let $A$ and $B$ be rings
with idempotents. We assume that $B$ is a non-degenerate $A-B$-bimodule (if the bimodule structure comes
from a morphism $\phi:A\rightarrow B$, it is not necessarily non-degenerate. Care must be taken, in this
situation, to ensure that this hypothesis is satisfied). We define the forgetful functor:
\begin{equation}\label{Foubli}
\index[not]{FAB@$\caF_A^B$}\caF_A^B \colon \caM(B) \rightarrow \caM(A), \quad N \mapsto B\otimes_B N,
\end{equation}
where $B$ is considered as an $A-B$-bimodule. Note that $B\otimes_B N$ is indeed a non-degenerate left
$A$-module, since $B$ is.

If $B$ is a non-degenerate $B-A$-bimodule, we define similarly the \emph{pseudo-forgetful functor}:
\index[ter]{pseudo-forgetful functor}
\begin{equation}\label{pseudooubli}
{} \index[not]{FABa@$\spcheck\caF_A^B$}\spcheck\caF_A^B \colon  \caM(B) \rightarrow \caM(A),\quad  N \mapsto \Hom_B(B,N)_A.
\end{equation}
Note that $\Hom_B(B,N)$ is the completion of the $B$-module $N$, in the sense of Section \ref{compmod}.
The $A$-module structure on $\Hom_B(B,N)$ is given by:
\[
(a\cdot \phi)(b')=\phi(b'a), \quad (a \in A),\,  (b' \in B), (\phi \in \Hom_B(B,N)).
\]

\begin{rmq}
The terminology pseudo-forgetful functor is taken from \cite{KV}, and serves to convey the idea that it
is a functor resembling a forgetful functor. It should be noted that this has no relation with
pseudo-functors.
\end{rmq}

\begin{lemme}
Let $B$ be a $\bbC$-algebra with idempotents and $N$ a $B$-module. Then $B\otimes_B N$ identifies with
$N_B$, the non-degenerate part of $N$ via
\begin{equation}\label{BBNN}
 B\otimes_B N \rightarrow  N, \quad b \otimes n \mapsto b\cdot n.
\end{equation}
If $N$ is non-degenerate, we then have $B\otimes_B N \simeq N$. Still in the case where $N$ is
non-degenerate, the natural morphism
\begin{equation}\label{HBBN}
N \rightarrow  \Hom_B(B,N),\quad n \mapsto \phi_n
\end{equation}
\[
\phi_n(b)=b\cdot n, \; (b\in B), \, (n \in N).
\]
is injective, and its image is the non-degenerate part $\Hom_B(B,N)_B$ of the $B$-module \allowbreak
$\Hom_B(B,N)$.
\end{lemme}

\begin{proof}
Since $B$ is non-degenerate as a module over itself, it is clear that $B\otimes_B N$ is non-degenerate,
and its image under any morphism is non-degenerate. For any $n$ in $N_B$, there exists an idempotent $e$
in $B$ such that $e\cdot n=n$, and therefore the map (\ref{BBNN}) is surjective onto the non-degenerate
part $N_B$. Now suppose that $\sum_{i=1}^r b_i\cdot n_i=0$, $b_i \in B$, $n_i\in N$. Let $e$ be an
idempotent of $B$ such that $eb_i=b_i$ for all $i=1,\ldots ,r$. Then
\[
\sum_{i=1}^r b_i\otimes  n_i=  \sum_{i=1}^r eb_i\otimes  n_i=e\otimes \sum_{i=1}^r b_i\cdot n_i=0,
\]
which shows that the map (\ref{BBNN}) is injective. The last assertion has already been proved in the
proof of Proposition \ref{lisscomp} (see (\ref{MbarA})).
\end{proof}

\begin{rmq}
The pseudo-forgetful functor is obtained by taking the non-degenerate part of $\Hom_B(B,N)$ as an
$A$-module. When the non-degenerate parts of $\Hom_B(B,N)=\overline{N}$ as an $A$ or $B$ module coincide,
we therefore have ${}\spcheck \caF_A^B(N)\simeq N$.
\end{rmq}

Let us now construct the adjoints (respectively right and left) of $\caF_A^B$ and ${}\spcheck \caF_A^B$.
We will call these functors ``induction functors''.
\begin{defi}
If $B$ is a non-degenerate $A-B$-bimodule, we set
\[
I_A^ B: \caM(A) \rightarrow \caM(B),\quad M \mapsto \Hom_A(B,M)_B.
\]
and if $B$ is a non-degenerate $B-A$-bimodule, we set
\[
P_A^ B: \caM(A) \rightarrow \caM(B),\quad M \mapsto B \otimes_A M.
\]
\end{defi}
\index[not]{IAB@$I_A^B$}\index[not]{PAB@$P_A^B$} We then have:
\begin{thm}
The functor $I_A^ B$ is the right adjoint of the forgetful functor $\caF_A^B$, and the functor $P_A^ B$
is the left adjoint of the pseudo-forgetful functor ${}\spcheck  \caF_A^B$. We therefore have, for any
$N \in \caM(A)$ and $M \in \caM(B)$, natural isomorphisms
\[
\Hom_B(M, I_A^ B(N))\simeq \Hom_A( \caF_A^B(  M),N),
\]
\[
\Hom_B(P_A^ B(N),M)\simeq \Hom_A(N, {}\spcheck \caF_A^B( M)).
\]
\end{thm}

\begin{proof}
This follows immediately from the definitions and the corollary of Theorem \ref{KVC20}.
\end{proof}

\begin{rmq}
We will need to further clarify the situation in the particular case where $B$ is equipped with an
$A$-bimodule structure coming from a morphism $\phi \colon A \rightarrow B$. We always assume $B$ is
non-degenerate for this $A$-bimodule structure. As we noted above, this is not automatic. The forgetful
functor is then equal to $\phi_*$ by (\ref{BBNN}). The $A$-bimodule $B$ being assumed non-degenerate,
this means that for any finite family $(b_i)$ in $B$, there exists an idempotent $e$ in $A$ such that for
all $i$, $\phi(e)b_i=b_i\phi(e)=b_i$. In other words
\[
\{ \phi(e),\, e \in \mathrm{Idem}(A)  \}
\]
is a directed system of idempotents for $B$ (cf. Remark \ref{idem}). In particular, the non-degenerate
part of a $B$-module coincides with its non-degenerate part for the $A$-module structure given by the
forgetful functor. Let $N$ be a non-degenerate $B$-module, and $\Hom_B(B,N)_A$ its image under the
pseudo-forgetful functor. From what we have just said, the non-degenerate part of $\Hom_B(B,N)$ for the
$A$-module and $B$-module structures coincide, and therefore, according to the previous remark,
\[
\Hom_B(B,N)_A=\Hom_B(B,N)_B\simeq N.
\]
We see that the pseudo-forgetful functor is in this case naturally isomorphic to the forgetful functor.
\end{rmq}

\subsection{Exactness properties of forgetful and induction functors}

The notation is the same as in the previous section, where we defined the functors $\caF_A^B$,
${}\spcheck \caF_A^B$, $I_A^B$ and $P_A^B$. We now give the exactness properties of these functors, and
some consequences thereof.

\begin{thm}
The functor $\caF_A^B$ is exact, the functors $I_A^B$ and ${}\spcheck \caF_A^B$ are left exact. The
functor $P_A^B$ is right exact.
\end{thm}
\begin{proof}
This results from the right exactness of the $\otimes$ functors and the left exactness of the $\Hom$
functors (Proposition \ref{KVC20}) and the fact that the functor
\[
A-\mathbf{mod}\rightarrow \caM(A), \quad M \mapsto M_A
\]
is exact (Proposition \ref{idem}). For the first assertion, one must also verify the left exactness of
the forgetful functor, but this follows trivially from (\ref{BBNN}).
\end{proof}

\begin{cor}
The functor $I_A^B$ preserves the injectivity of modules.
\end{cor}

\begin{proof}
This follows from the fact that $I_A^ B$ is the right adjoint of an exact functor. We recall the
argument, which is standard. Let $X$ be an injective module in $\caM(A)$. We need to prove that the
functor
\[
\caM(B)\rightarrow \bbZ-\mathbf{mod},\quad Y \mapsto \Hom_B(Y, I_A^B (X) )
\]
is exact. But this functor, by adjunction, is isomorphic to the functor
\[
\caM(B)\rightarrow \bbZ-\mathbf{mod}, \quad Y \mapsto \Hom_A(\caF_A^B (Y),  X ).
\]
We see that the latter is exact, being the composition of the forgetful functor $\caF_A^B$ which is exact
according to the theorem, and the functor $\Hom_A(\, \bullet \, ,  X )$ which is exact by hypothesis.
\end{proof}

\section{The functor $j_e \colon M \mapsto e\cdot M$}

\subsection{Exactness}\label{exactje}
Let $A$ be a $\bbC$-algebra with idempotents, and $e$ an idempotent of $A$. As we noted above, any
non-degenerate $A$-module $M$ decomposes into
\begin{equation}\label{emem}
M=e\cdot M \oplus (1-e)\cdot M
\end{equation}
and $e\cdot M $ is a unital $eAe$-module. It is easy to see that this defines a functor
\[
j_e \colon  \caM(A) \rightarrow \caM(eAe), \quad M \mapsto e\cdot M.
\]
\begin{prop}
The functor $j_e$ defined above is exact.
\end{prop}
\begin{proof}
It is clear that any morphism of $A$-modules preserves the decompositions (\ref{emem}). Let
$M_1 \stackrel{\psi}{\rightarrow}  M_2 \stackrel {\phi}{\rightarrow}  M_3$ be an exact sequence in
$\caM(A)$. We need to check that if $e\cdot m_2 \in \ker \phi$, there exists $e\cdot m_1 \in e\cdot M_1$
such that $\psi(e\cdot m_1)=e\cdot m_2$. Let us take $m_1 \in M_1$ such that $\psi(m_1)=e\cdot m_2$. We
then have $e\cdot(e\cdot m_2)=e\cdot m_2=e\cdot \psi(m_1)=\psi(e\cdot m_1)$.
\end{proof}

\subsection{Simple modules}\label{35}
Let $A$ be a $\bbC$-algebra with idempotents and $e \in \mathrm{Idem}(A)$. Let us denote:
\medskip

\noindent \index[not]{M(Ae)@$\caM(A,e)$}$\caM(A,e)$ the full subcategory of $\caM(A)$ of modules $M$ such
that $M=Ae\cdot M$.

\noindent \index[not]{Irr(A)@$\mathbf{Irr}(A)$}$\mathbf{Irr}(A)$ the set of isomorphism classes of
non-degenerate simple $A$-modules,

\noindent \index[not]{IrreAe@$\mathbf{Irr}(eAe)$}$\mathbf{Irr}(eAe)$ the set of isomorphism classes of
unital simple $eAe$-modules,

\noindent \index[not]{IrrAe@$\mathbf{Irr}(A,e)$}$\mathbf{Irr}(A,e)$ the subset of $\mathbf{Irr}(A)$ of
isomorphism classes of modules $M$ such that $e \cdot M\neq 0$.

\medskip

The elements of $\caM(A,e)$ are therefore the modules $M$ in $\caM(A)$ generated by $e\cdot M$, and
$\mathbf{Irr}(A,e)$ is the set of isomorphism classes of irreducible objects of $\caM(A,e)$.

\begin{lemme}
Consider the induction functor:
\[
i \colon  \caM(eAe)\rightarrow \caM(A) \quad Z \mapsto A\otimes_{eAe}Z.
\]
Then $j_e \circ i$ is naturally isomorphic to the identity of $\caM(eAe)$, i.e., for any
$Z \in \caM(eAe)$,
\begin{align}\label{riz}
 Z  &\simeq j_e \circ i(Z), \\
\nonumber z &\mapsto e\otimes z
\end{align}
these isomorphisms being natural in $Z$.
\end{lemme}

Note that we are not quite in the conditions of Section \ref{Oublietadjoints} where the functor $P$ is
defined similarly, because $A$ is not a unital $eAe$-module, as $e$ does not act as the identity of $A$.
Nevertheless, if $Z$ is a unital $eAe$-module, it is in particular a $\bbC$-vector space. The algebra $A$
is as well, and the two induced $\bbC$-vector space structures on the tensor product $A\otimes_{eAe}Z$
coincide.

\begin{proof}
By definition $j_e\circ i(Z)=e( A\otimes_{eAe}Z)=e A\otimes_{eAe}Z$, and we have an injection
$z \mapsto e\otimes z$ from $Z$ into $e A\otimes_{eAe}Z$. Indeed, we need to show that if $z\in Z$ is
non-zero, $e\otimes z$ does not vanish in $A\otimes_{eAe}Z$. By the definition of the tensor product, it
suffices to find a balanced bilinear map (for the $\bbZ$-module structure)
$B \colon A\times Z\rightarrow W$ which does not vanish on $(e,z)$ because its $\bbZ$-linear lift to the
tensor product will not vanish on $e\otimes z$. Let us take $W=\Hom_\bbC(A,Z)$ and
\[
B \colon A\times Z \longrightarrow   \Hom_\bbC(A,Z) , \quad (a,z)\mapsto [\phi_{a,z}  \colon b\mapsto ebae \cdot z].
\]
It is clear that $B$ is $\bbZ$-bilinear. It is balanced because
\[
B(aece,z) =[b\mapsto eb(aece)e\cdot z  =   ebaece\cdot z  ]=B(a, ece\cdot z).
\]
We have $B(e,z)= [\phi_{e,z} \colon b\mapsto ebee\cdot z=ebe\cdot z]  $, and
$\phi_{e,z}(e)=e\cdot z=z\neq 0$ so $\phi_{e,z}\neq 0$.

We will now show that $z \mapsto e\otimes z$ from $Z$ into $e A\otimes_{eAe}Z$ is also surjective: if
$ea\otimes z \in e A\otimes_{eAe}Z$,
\[
ea\otimes z = ea\otimes ez= eae\otimes z  = e\otimes eaez     .
\]
The naturality of these isomorphisms is routinely verified.
\end{proof}

\medskip

We immediately deduce that $A\cdot(e\cdot i(Z))=i(Z)$, i.e., the functor $i$ takes values in $\caM(A,e)$.

The functor $j_e$ may annihilate certain modules of $\caM(A)$, and therefore we cannot hope to
reconstruct all non-degenerate $A$-modules from modules in $\caM(eAe)$. Nevertheless, we will obtain all
the irreducibles in $\caM(A,e)$.

\begin{defi}
Let $M$ be a non-degenerate $A$-module and $e \in \mathrm{Idem}(A)$. We define the non-degenerate
$A$-module:
\[
M_e:= M/F(eA,M), \quad \text{ where } \quad F(eA,M)=\{ m\in M\mid eA\cdot m=0  \}.
\]
\end{defi}

\begin{lemme}
Let $M$ be a non-degenerate $A$-module. We then have $(M_e)_e=M_e$.
\end{lemme}
\begin{proof}
Let $\bar m\in M_e$ such that $eA\cdot \bar m=0$. Let $m \in M$ be a lift of $\bar m$ in $M$. We have
$eA \cdot m \subset F(eA,M)$ and therefore $eAeA\cdot m=0$. In particular $eeA \cdot m=eA\cdot m=0$,
hence $\bar m=0$. It follows that $(M_e)_e=M_e$.
\end{proof}

\begin{prop}
The map $M \mapsto e\cdot M$ realizes a bijection between $\mathbf{Irr}(A,e)$ and $\mathbf{Irr}(eAe)$,
whose inverse is given by $W \mapsto  (A\otimes_{eAe}W)_e$.
\end{prop}

\begin{proof}
We will now show that if $M$ is a non-degenerate simple $A$-module such that $e\cdot M\neq 0$, then
$e\cdot M$ is a unital simple $eAe$-module. Let $e\cdot m\in e\cdot M$ be non-zero. Since $M$ is simple,
we have
\[
(eAe)e\cdot m=(eA)e\cdot m=e(Ae\cdot m)=e\cdot M.
\]
Therefore $e\cdot M$ does not contain any proper $eAe$-submodule, so it is a simple module. This shows
that $M \mapsto e\cdot M$ defines a map from $\mathbf{Irr}(A,e)$ to $\mathbf{Irr}(eAe)$.

Let us now show that if $W$ is a unital simple $eAe$-module, then $(A\otimes_{eAe}W)_e$ is a simple
$A$-module. The canonical projection $(A\otimes_{eAe}W) \rightarrow (A\otimes_{eAe}W)_e$ is transformed
by the functor $j_e$ into a morphism of unital $eAe$-modules
\begin{equation}\label{eAWe}
e\cdot(A\otimes_{eAe}W) \rightarrow e\cdot (A\otimes_{eAe}W)_e.
\end{equation}
This morphism is surjective, since $j_e$ is exact (Proposition \ref{exactje}). It is also injective
because an element of the kernel is an invariant under $e$ of $F(eA,A \otimes_{eAe} W)$, hence zero by
definition.

This shows that the morphism (\ref{eAWe}) is an isomorphism. Composing it with (\ref{riz}), we obtain an
isomorphism of $eAe$-modules $W\simeq   e\cdot (A\otimes_{eAe}W)_e$.

Let $\bar w$ be a non-zero element of $(A\otimes_{eAe}W)_e$. We want to show that
$A\cdot \bar w= (A\otimes_{eAe}W)_e$. Indeed, this implies that $(A\otimes_{eAe}W)_e$ has no proper
submodules, and therefore that it is a simple module. Suppose $e\cdot \bar w\neq 0$. Since
$W\simeq  e\cdot (A\otimes_{eAe}W)_e$ and $W$ is a simple $eAe$-module,
\[
eAe\cdot (e\cdot \bar w)=  e\cdot (A\otimes_{eAe}W)_e.
\]
We have seen that $A\otimes_{eAe}W$ is generated by $e\cdot (A\otimes_{eAe}W)$, and consequently
$(A\otimes_{eAe}W)_e$ is generated by $e\cdot (A\otimes_{eAe}W)_e$. We therefore have
$A\cdot \bar w= Ae\cdot \bar w=(A\otimes_{eAe}W)_e$, which shows the assertion in this case.

It remains to treat the case where $e\cdot \bar w= 0$. Let us use the previous lemma for this, which
asserts that $(A\otimes_{eAe}W)_e=((A\otimes_{eAe}W)_e)_e$, and therefore that $0\neq eA\cdot \bar w$.
Let $\bar w'$ be a non-zero element of the form $\bar w'=a\cdot \bar w$ such that $e\cdot \bar w'\neq 0$.
The preceding shows that $A\cdot \bar w'= (A\otimes_{eAe}W)_e$ and we deduce that
$A\cdot \bar w=(A\otimes_{eAe}W)_e$. This completes the proof that $W \mapsto  (A \otimes_{eAe}W)_e$
defines a map from $\mathbf{Irr}(eAe)$ to $\mathbf{Irr}(A,e)$.

We now need to see that the two maps thus defined are inverse to each other. We have already shown that
if $W$ is in $\mathbf{Irr}(eAe)$, $W\simeq e\cdot (A\otimes_{eAe}W)_e$. Let $M$ be in
$\mathbf{Irr}(A,e)$. Let
\[
\phi\colon A\otimes_{eAe}e\cdot M \rightarrow M, \quad a\otimes e\cdot m \mapsto a e\cdot m.
\]

This morphism is surjective by hypothesis on $M$. Moreover $eA\cdot (\ker \phi)=0$. Indeed, suppose that
$$\phi \left(\sum_i a_i\otimes e\cdot m_i\right)=\sum_i a_i e\cdot m_i= 0.$$ Let $ea \in eA$. We have:
\[
ea\cdot \left( \sum_i a_i\otimes e\cdot m_i \right)= \sum_i eaa_ie\otimes e\cdot m_i
=e\otimes ea\left(\sum_i  a_i  e\cdot m_i\right)=0.
\]
This shows that $\ker \phi \subset F(eA, A\otimes_{eAe}e\cdot M )$. Since $M$ is a simple $A$-module, and
$\phi$ is surjective, $\ker \phi$ is a maximal proper submodule and therefore
$\ker \phi = F(eA,A\otimes_{eAe}e\cdot M)$ because
$A\otimes_{eAe}e\cdot M \neq  F(eA,A\otimes_{eAe}e\cdot M)$.

The morphism $\phi$ therefore factors through
\[
\bar \phi \colon  ( A\otimes_{eAe}e\cdot M)_e \rightarrow M,\quad a \otimes e\cdot m \mapsto ae\cdot m.
\]
Moreover, $\bar \phi$ is injective and remains surjective, so it is an isomorphism. This completes the
proof of the proposition.
\end{proof}

\medskip

In certain cases, the functor $j_e$ is in fact an equivalence of categories between $\caM(A,e)$ and
$\caM(eAe)$. The following result gives a criterion for this.

\begin{thm}
The functor $j_e$ is an equivalence of categories between $\caM(A,e)$ and $\caM(eAe)$ if and only if the
subcategory $\caM(A,e)$ of $\caM(A)$ is stable under passing to subquotients (in particular, it is then
an abelian category). The inverse functor is given by the functor $i=P_{eAe}^A$.
\end{thm}

\begin{proof}
To show that the condition is necessary, it suffices to verify that if $V$ is in $\caM(A,e)$, and if $W$
is an $A$-submodule of $V$, then $W$ and $V/W$ are still in $\caM(A,e)$. It is clear that if
$Ae\cdot V=V$, then $Ae\cdot (V/W)=V/W$. It remains to show that $W \in \caM(A,e)$. On the other hand, by
replacing $V$ with $V/Ae\cdot W$, we reduce to the case where $Ae\cdot W=0$. We then need to show that
$W=0$. Consider the canonical projection $p \colon V \rightarrow V/W$, and its image under the functor
$j_e$,
\[
j_e(p) \colon  e \cdot V \rightarrow e \cdot (V/W).
\]
Since $j_e$ is exact, $j_e(p)$ is surjective, but we also see that it is injective. Indeed, if
$e\cdot v=w \in W$, we have $e\cdot v=w=e\cdot w=0$ since $e\cdot W=0$. If we assume that $j_e$ is an
equivalence of categories, the projection $p$ must be an isomorphism, and then $W=0$.

We will now show that the condition is sufficient. It is clear that for any $V \in \caM(A,e)$, the
submodule $F(eA,V)$ is trivial. The proof of the previous proposition then suggests that the inverse of
$j_e$ is given by the induction functor $i$. We have already shown in (\ref{riz}) that $j_e \circ i$ is
naturally isomorphic to the identity of $\caM(eAe)$. If $V$ is an object of $\caM(A,e)$, we have a
canonical surjection
\[
\psi=\psi_V \colon (i\circ j_e)(V)\rightarrow V, \quad a \otimes e\cdot v \mapsto ae\cdot v
\]
Let $W$ be the kernel of $\psi$. We have $e\cdot W=W \cap e \cdot V=W \cap (j_e\circ i\circ j_e)(V)=0$
because $\psi$ is injective on $(j_e\circ i\circ j_e)(V)$. But by hypothesis, $W=Ae\cdot W$, hence $W=0$,
which shows that $\psi$ is an isomorphism. We deduce that the $\psi_V$ give a natural isomorphism between
$i\circ j_e$ and the identity of $\caM(A,e)$.
\end{proof}

\section{Duality}
\index[ter]{duality}
\subsection{The functor $M\mapsto \widetilde M$}\label{dualite}
Let $A$ be a $\bbC$-algebra with idempotents, and $M$ a non-degenerate (left) $A$-module. The dual of $M$
is
\[
M^*=\Hom_\bbC(M,\bbC)
\]
which is naturally equipped with a right $A$-module structure. Let $\widetilde M$ denote the
non-degenerate part of this $A$-module. This defines a contravariant functor:
\[
\caM(A) \rightarrow \caM(A)_d, \quad M \mapsto \widetilde M.
\]

Similarly, if $M$ is a non-degenerate right $A$-module, its dual $\widetilde M$ is a non-degenerate left
$A$-module, and we have a (contravariant) functor
\[
\caM(A)_r \rightarrow \caM(A), \quad M \mapsto \index[not]{Mt@$\widetilde M$}\widetilde M.
\]

In practice, the algebra $A$ is often equipped with an anti-involution via which one can identify right
and left modules, and duality will then be a functor from $\caM(A)$ to $\caM(A)$.

\begin{lemme}\label{dualitelem}
Let $M$ be a non-degenerate $A$-module, and let $e$ be an idempotent of $A$. We then have
\[
\widetilde M \cdot e \simeq  (e\cdot M)^*.
\]
The duality functor $M \mapsto \widetilde M$ is exact.
\end{lemme}
\begin{proof}
The decomposition $M=e\cdot M \oplus (1-e)\cdot M$ allows us to identify $(e\cdot M)^*$ with the space of
linear forms on $M$ vanishing on $(1-e)\cdot M$. But these are clearly the linear forms fixed by the
(right) action of $e$, i.e., $(e\cdot M)^*\simeq   M^* \cdot e \simeq \widetilde M \cdot e $.

Let us now show the exactness of the duality functor. If
\[
0 \rightarrow M_1 \rightarrow M_2 \rightarrow M_3 \rightarrow 0
\]
is an exact sequence in $\caM(A)$, we want to show the exactness of the sequence
\[
0 \rightarrow \widetilde M_3 \rightarrow \widetilde M_2 \rightarrow \widetilde M_1 \rightarrow 0.
\]
Since these modules are non-degenerate, it suffices to show that for any idempotent $e$ of $A$, the
sequence
\[
0 \rightarrow  \widetilde M_3 \cdot e \rightarrow \widetilde M_2 \cdot e \rightarrow   \widetilde M_1\cdot e  \rightarrow 0
\]
is exact. From the above, this can be rewritten in the form
\[
0 \rightarrow (e\cdot M_3)^* \rightarrow (e\cdot  M_2)^* \rightarrow (e\cdot M_1)^* \rightarrow 0.
\]
But the functors $j_e \colon M\mapsto e\cdot M$ and $M \mapsto M^*$ are exact. This completes the proof
of the lemma.
\end{proof}

\begin{cor}\label{dualitecor}
We have a canonical injection $M\hookrightarrow \widetilde {\widetilde M}$.
\end{cor}

\begin{proof}
Any vector $m\in M$ is fixed by a certain idempotent $e$ of $A$. Consequently $m \in e\cdot M$ and $m$
defines an element of $((e\cdot M)^*)^*$. Therefore, since
$ \widetilde M\cdot e= M^*\cdot e=(e\cdot M)^*$, $m$ defines an element of the dual of
$\widetilde M\cdot e$. The preceding equality applied to $\widetilde M$ tells us that
$m\in  (\widetilde M\cdot e)^*= e\cdot (\widetilde M^*)=e\cdot(\widetilde {\widetilde M})$. The morphism
$M \rightarrow \widetilde {\widetilde M}$ is injective because for any non-zero $m \in M$, we can find
$\lambda \in M^*$ such that $\lambda(m)\neq 0$, and if $e\cdot m=m$, then $(\lambda\cdot e)(m)\neq 0$,
$\lambda \cdot e \in \widetilde M$.
\end{proof}

\subsection{Duality and forgetful functor}\label{duaetoub}

The goal of this section is to show that $M \mapsto \widetilde M$ intertwines the forgetful functor and
the pseudo-forgetful functor, as well as their respective adjoints $I_A^B$ and $P_A^B$.

\begin{thm}
Let $A$ and $B$ be algebras with idempotents, and suppose that $B$ is equipped with a non-degenerate
$A$-bimodule structure. We assume that the right actions of $A$ and $B$ on $B$ commute with the left
actions, so that $B$ is also an $A-B$ and a $B-A$-bimodule. The forgetful functor $\caF_A^B$ and the
pseudo-forgetful functor ${}\spcheck  \caF_A^B$, and their respective right and left adjoints $I_A^ B$
and $P_A^ B$ are then defined as in the previous section (with the obvious adaptations when considering
right modules rather than left modules). We then have, for any non-degenerate left $B$-module $N$,
\[
{}\spcheck  \caF_A^B(\widetilde N)=   \caF_A^B(N)^\sim ,
\]
and for any non-degenerate left $A$-module $M$,
\[
I_A^ B(\widetilde M)=P_A^B(M)^\sim.
\]
\end{thm}
The equalities of functors stated in this theorem are obviously to be interpreted as simply the existence
of natural isomorphisms (which are moreover explicit) between them. This is a very widespread abuse of
language and notation for convenience.

\begin{proof}
We have
\begin{align*}
{}\spcheck \caF_A^B(\widetilde N)&= \Hom_B(B, \Hom_\bbC(N,\bbC)_B)_A
                  = \Hom_B(B, \Hom_\bbC(N,\bbC) )_A \simeq  \Hom_\bbC(B   \otimes_B N,\bbC)_A\\
&                  \simeq   \Hom_\bbC(\caF_A^B(N), \bbC)_A \simeq (\caF_A^B(N))^\sim
\end{align*}
We first use the fact that the image of a non-degenerate module is non-degenerate, then the first
identity of Theorem \ref{KVC20}, where $N$ is viewed as a left $B$-module and a right $\bbC$-module, and
$\bbC$ as the trivial right $\bbC$-module.

The second assertion of the theorem is shown in the same way using Theorem \ref{KVC20}. It can also be
deduced from the first by adjunction. We use the following identity for this.

\begin{lemme}
For any right $A$-module $V$ and for any left $A$-module $W$, both non-degenerate
\begin{align}\label{VWWV}
\Hom_A(V,\widetilde W)&\longrightarrow \Hom_A(W,\widetilde V)\\
\nonumber \phi &\mapsto \psi, \quad \psi(w)(v)=\phi(v)(w), \quad (\forall v \in V, w \in W),
\end{align}
is a natural isomorphism in $V$ and $W$.
\end{lemme}
\begin{proof}
We have
\begin{align*}
 \Hom_A(V,\widetilde W)&= \Hom_A(V, \Hom_\bbC(W, \bbC)_A)= \Hom_A(V,\Hom_\bbC(W, \bbC))\simeq \Hom_\bbC(V\otimes_A W, \bbC).
\end{align*}
This is obtained from the first identity of Theorem \ref{KVC20}, by considering $W$ as an
$A-\bbC$-bimodule. Similarly, by considering $V$ as a $\bbC-A$-bimodule, the second identity of Theorem
\ref{KVC20} gives us
\[
\Hom_A(W,\widetilde V)\simeq \Hom_\bbC(V\otimes_A W,\bbC).
\]
All these identities are functorial in $V$ and $W$.
\end{proof}

\medskip

We now complete the proof of the theorem. We will now show that we can deduce from the lemma that
$I_A^ B(\widetilde M)=P_A^B(M)^\sim $. We have, for any $M \in \caM(A)$, $N \in \caM(B)$,
\begin{align*}
\Hom_A(M, {}\spcheck  \caF_A^B(\widetilde N))&\simeq   \Hom_B(P_A^ B(M),\widetilde N) \simeq  \Hom_B(N, P_A^ B(M)^\sim\, ).
\end{align*}
This is obtained by using the adjunction between ${}\spcheck  \caF_A^B$ and $P_A^ B$ and the identity of
the lemma. On the other hand
\begin{align*}
\Hom_A(M, {}\spcheck  \caF_A^B(\widetilde N))&\simeq \Hom_A(M, \caF_A^B(N)^\sim\, )
\simeq  \Hom_A(\caF_A^B(N), \widetilde  M ) \simeq  \Hom_B(N,I_A^ B(\widetilde M)).
\end{align*}
Hence
\[
\Hom_B(N, P_A^ B(M)^\sim\, )\simeq  \Hom_B(N,I_A^ B(\widetilde M)),
\]
which, by the uniqueness of the adjoint up to isomorphism, shows the desired identity.
\end{proof}

\subsection{Duality and completion}\label{dualetcomp}
Let $A$ be a $\bbC$-algebra with idempotents and let $M$ be a non-degenerate $A$-module. We defined
$\widetilde M$ as the non-degenerate part of $M^*$. It is easy to see directly that when we complete
$\widetilde M$, we recover $M^*$. We will provide a proof of this fact which follows formally from what
has already been established. For any right $A$-module $N$, we indeed have:
\begin{align*}
\Hom_A(N,\overline{(M^*)_A}) &\simeq \Hom_A(N_A,(M^*)_A)
\simeq \Hom_A(N_A,\Hom_\bbC(M,\bbC)_A)\\
&\simeq  \Hom_A(N_A,\Hom_\bbC(M,\bbC))
\simeq \Hom_\bbC(N_A \otimes_A M,\bbC).
\end{align*}
The first equality comes from the adjunction \ref{lisscomp}, the last from (\ref{ADId2}), with $R=A$ and
$S=\bbC$. On the other hand:
\begin{align*}
\Hom_A(N,M^*)&=  \Hom_A(N,\Hom_\bbC(M,\bbC))
= \Hom_{\widetilde A}(N,\Hom_\bbC(M,\bbC))
\simeq \Hom_\bbC(N \otimes_{\widetilde A} M,\bbC)\\
& \simeq \Hom_\bbC(N \otimes_{A} M,\bbC)
\end{align*}
We used (\ref{AD2}) and the related notation for the second equality.

We show that the canonical morphism $\iota \colon N_A \otimes_A M \to N \otimes_A M$ is an isomorphism.
Since $M$ is non-degenerate, for any $m \in M$, there exists an idempotent $e \in A$ such that
$e \cdot m = m$. Define a map $\phi \colon N \times M \to N_A \otimes_A M$ by
$\phi(n, m) = (n \cdot e) \otimes m$, where $e$ is such an idempotent. We first verify that this map is
well-defined, i.e., it does not depend on the choice of $e$. If $f \in A$ is another idempotent such that
$f \cdot m = m$, choose an idempotent $h \in A$ such that $he=e$, $eh=e$, $hf=f$ and $fh=f$ (which is
always possible for a finite set of elements in a ring with idempotents). Then, in the tensor product
$N_A \otimes_A M$, we have:
\[
(n \cdot h) \otimes m = (n \cdot h) \otimes (e \cdot m) = (n \cdot he) \otimes m = (n \cdot e) \otimes m.
\]
By symmetry, $(n \cdot f) \otimes m = (n \cdot h) \otimes m$, which proves that $\phi$ is well-defined.
We next verify that $\phi$ is $A$-balanced. For $n \in N$, $a \in A$ and $m \in M$, choose an idempotent
$h \in A$ such that $h \cdot m = m$, $a h = a$ and $h a = a$. Then $h$ fixes $m$ and satisfies
$h \cdot (a \cdot m) = (ha) \cdot m = a \cdot m$, so $h$ also fixes $a \cdot m$. We then have:
\[
\phi(n \cdot a, m) = (n \cdot a h) \otimes m = (n \cdot a) \otimes m,
\]
and
\[
\phi(n, a \cdot m) = (n \cdot h) \otimes (a \cdot m) = (n \cdot h a) \otimes m = (n \cdot a) \otimes m.
\]
The bilinear map $\phi$ therefore induces a morphism of abelian groups
$\Phi \colon N \otimes_A M \to N_A \otimes_A M$. For any $n \in N_A$ and $m \in M$, choose an idempotent
$e$ fixing both $n$ (on the right) and $m$ (on the left). We have
$\Phi(\iota(n \otimes m)) = \Phi(n \otimes m) = (n \cdot e) \otimes m = n \otimes m$. Thus,
$\Phi \circ \iota = \mathrm{Id}_{N_A \otimes_A M}$, which proves that $\iota$ is injective. Since
$\iota \circ \Phi$ is trivially the identity on the generators of $N \otimes_A M$, $\iota$ is an
isomorphism.

We obtain for any right $A$-module $N$ a (natural) isomorphism:
\[
\Hom_A(N,\overline{(M^*)_A}) \simeq  \Hom_{A}(N,M^*)
\]
We deduce from \ref{egalitecat} that
\[
\overline{\widetilde{M}}=\overline{(M^*)_A}\simeq M^*.
\]

This gives us a realization of $\overline{M}$ which is convenient for computations. We embed $M$ into a
module $\widetilde L$ (for example, we can use the embedding of $M$ into $\widetilde{\widetilde M}$). We
then have, according to (\ref{compsousmod}),
\begin{equation}\label{dualetcompeq}
\overline{M} =\{ l \in L^*\, \mid \,
A\cdot  l  \subset M \}.
\end{equation}

\section{Some classes of modules}

\subsection{Admissible modules}
\indexter{admissible!(module)} \label{modadm} As in the previous sections, $A$ is a $\bbC$-algebra with
idempotents.

\begin{defi}
A non-degenerate module $M$ in $\caM(A)$ is said to be admissible if for any idempotent $e$ of $A$, the
vector space $e\cdot M$ is finite-dimensional.
\end{defi}

\begin{prop}\label{modadmprop}
A module $M$ is admissible if and only if the inclusion $M \hookrightarrow \widetilde {\widetilde M}$ is
an isomorphism.
\end{prop}

\begin{proof}
Suppose $M$ is admissible. Then for any idempotent $e$ of $A$, $e\cdot M$ is finite-dimensional. The
injection of $e\cdot M$ into $e\cdot (\widetilde{\widetilde M})$ exhibited in the proof of Corollary
\ref{dualitecor} is then an isomorphism. Since every vector of $\widetilde{\widetilde M}$ is fixed by an
idempotent of $A$, we obtain $M \simeq \widetilde{\widetilde  M}$. Conversely, suppose that
$M\simeq  \widetilde{\widetilde M}$. Let $e$ be an idempotent of $A$. We therefore have, according to the
isomorphisms proved in \ref{dualite},
\[
e\cdot M= e\cdot \widetilde{\widetilde M}
=e\cdot ((\widetilde M)^*)=(\widetilde M\cdot e)^*
=((e\cdot M)^*)^*.
\]
Thus $e\cdot M$ is finite-dimensional.
\end{proof}

\begin{cor}\label{modadmcor}
Let $M$ be an admissible non-degenerate $A$-module. Then $\widetilde M$ is admissible. Moreover, $M$ is
simple if and only if $\widetilde M$ is simple.
\end{cor}

\begin{proof}
The first point is a corollary of the proof of the above proposition, rather than of the result itself.
For the second point, it suffices to note that if $M$ admits a proper submodule $M_1$, then its
orthogonal in $\widetilde M$ is also a proper submodule of $\widetilde M$, and therefore $\widetilde M$
is reducible. The equality $M=\widetilde{\widetilde M}$ shows the converse.
\end{proof}

\subsection{Projective and injective modules}

Let $A$ be a $\bbC$-algebra with idempotents. Recall that a module $X \in \caM(A)$ is said to be
projective (resp. injective) if the functor $\Hom_A(X, \, \bullet \,  )$ (resp.
$\Hom_A(\, \bullet \, , X )$) is exact. \index[ter]{projective!(module)} \index[ter]{injective!(module)}

\begin{thm}
The category $\caM(A)$ has enough projectives (i.e., every module in $\caM(A)$ is a quotient of a
projective module).
\end{thm}

\begin{proof}
Let $e \in A$ be an idempotent, and consider the module $P_e=Ae$. The module $P_e$ is projective since
for any $A$-module $M$ in $\caM(A)$, $\Hom_A(P_e,M)=e\cdot M=j_e(M)$ and the functor $j_e$ is exact (\ref
{exactje}). On the other hand, any direct sum of projective modules is again projective. Indeed, for any
family $(P_i)_{i\in I}$ of modules $$\Hom_A(\oplus_{i\in I}P_i,M)= \prod_{i\in I} \Hom_A (P_i,M). $$ If
$M$ is in $\caM(A)$, and if $m\in M$, by definition, there exists an idempotent $e$ of $A$ such that
$m=e\cdot m$. Thus $m$ is in the image of $P_e \rightarrow M$, $a\cdot e \mapsto a \cdot m$. For each
$m\in M$, we choose such an $e$. By taking the sum over all $m\in M$ of the corresponding $P_e$, we see
that $M$ is a quotient of a projective module.
\end{proof}

\medskip

By using duality in $\caM(A)$, we deduce the existence of enough injectives from the existence of enough
projectives. For this, let us first show the following

\begin{lemme}
Let $X$ be a projective module in $\caM(A)_r$. Then $\widetilde X$ is injective in $\caM(A)$.
\end{lemme}

\begin{proof}
We need to show the exactness of the functor $\Hom_A(\, \bullet\, ,\widetilde X )$. However, according to
(\ref{VWWV}), this functor is isomorphic to the functor $\Hom_A(X,\,  \widetilde \bullet \,  )$ which is
the composition of the duality functor, exact according to Lemma \ref{dualitelem}, and the functor
$\Hom_A(X,\,  \bullet \, )$, exact by hypothesis.
\end{proof}

\begin{cor}
The category $\caM(A)$ has enough injectives (i.e., every module in $\caM(A)$ is a submodule of an
injective module).
\end{cor}

\begin{proof}
Let $M$ be in $\caM(A)$ and $P$ a projective module of $\caM(A)_r$ such that $\widetilde M$ is a quotient
of $P$. Then $\widetilde{\widetilde M}$ is a submodule of the injective module $\widetilde P$. Note that
we are dealing here with right modules, but the proof of the theorem above is also valid for the category
$\caM(A)_r$. Since $M$ embeds into $\widetilde{\widetilde M}$ (Corollary \ref{dualitecor}), it embeds
into $\widetilde P$.
\end{proof}

\subsection{Free modules}\label{modlibres}
Let $A$ be a $\bbC$-algebra with idempotents. We say that $M \in \caM(A)$ is a free module if it is
isomorphic to a direct sum of modules isomorphic to $A$.

Every module $M$ is a quotient of a free module. The proof is similar to that of the fact that $\caM(A)$
has enough projectives. Indeed, let $M\in \caM(A)$, and let $m\in M$, fixed by a certain idempotent $e$
of $A$. Then $m$ is in the image of the morphism $\phi \colon A \rightarrow M$, defined by
$\phi(a)=a\cdot m$. We thus construct a surjective morphism from $\bigoplus_{m\in M}A$ to $M$.

Unlike the case of unital rings, $A$ is not necessarily projective in $\caM(A)$, and therefore a direct
factor of a free module is not necessarily projective. On the other hand, every projective module is a
direct factor of a free module, the proof remaining the same as in the case of unital rings (cf.
\cite{Lang}).

\section{Properties of $\caM(A)$}
Let $A$ be a $\bbC$-algebra with idempotents. Since the category $\caM(A)$ is a full subcategory of
$A-\mathbf{mod}$ stable under passing to submodules and quotients, it is an abelian category. We refer
the reader to Appendix \ref{catab} for generalities on abelian categories and the terminology used.

\begin{prop}
The category $\caM(A)$ is complete and cocomplete and satisfies the axioms $\mathbf{AB5}$ and
$\mathbf{AB6}$ (cf. \ref{catab}). Inductive limits in $\caM(A)$ are exact.
\end{prop}

\begin{proof}
The first assertion means that arbitrary limits and colimits exist in $\caM(A)$, and to do this, it
suffices to verify the existence of products and direct sums in $\caM(A)$. We use the existence of
arbitrary products and direct sums in $A-\mathbf{mod}$. If $(M_i)_{i\in I}$ is a family of objects in
$\caM(A)$, then $\prod_i M_i$ and $\bigoplus_i M_i$ exist in $A-\mathbf{mod}$. It is immediate that
$\bigoplus_i M_i$ is non-degenerate, and is therefore the direct sum of the $M_i$ in $\caM(A)$. On the
other hand, the product $\prod_i M_i$ is not necessarily non-degenerate, but it is clear that its
non-degenerate part $(\prod_i M_i)_A$ satisfies the universal property of the product in $\caM(A)$.

Axiom $\mathbf{AB5}$ asserts that for any increasing directed system $M_i$ of subobjects of an object
$M$, and for any subobject $N$ of $M$, we have
\[
N\cap \left( \sum_i M_i  \right)=
\sum_i (N\cap M_i).
\]
Since all objects are modules, hence sets, we can use standard set-theoretic reasoning to show the above
equality, by showing inclusion in both directions.

Axiom $\mathbf{AB6}$ asserts that every object of $\caM(A)$ is an increasing sum of finitely generated
objects. Let us point out that the numbering of the $\mathbf{AB}$ axioms adopted in Appendix \ref{catab}
does not coincide with Grothendieck's: our axiom $\mathbf{AB6}$ expresses that the category is locally
finitely generated. Here again, this property is obvious in the case of a category of modules stable
under passing to submodules.

Finally, the last assertion means that if
\[
((M_i)_{i\in \caI}, (f_{ij})_{i\leq j}),\quad ((N_i)_{i\in \caI}, (g_{ij})_{i\leq j}),\quad ((P_i)_{i\in \caI}, (h_{ij})_{i\leq j})
\]
form inductive systems in $\caM(A)$ such that for all $i\in \caI$ we have an exact sequence
\[
0 \rightarrow M_i \rightarrow N_i\rightarrow P_i
\rightarrow 0
\]
then the sequence
\[
0 \rightarrow \varinjlim_i M_i \rightarrow
\varinjlim_i N_i\rightarrow \varinjlim_i  P_i
\rightarrow 0
\]
is exact. We can either deduce this from the same statement in $A-\mathbf{mod}$, proved for example in
\cite{Bal2}, Chap. II, \S 6, Prop. 3, by noting that all inductive limits are in fact in $\caM(A)$, or by
reproducing the proof.
\end{proof}

\section{Notes on Chapter I}

Generalities on rings with idempotents can be found in several places in the literature. I used
\cite{Del} for the section on the completion of the ring $A$ and on the center of $\caM(A)$. I followed
\cite{KV} for the definition of the induction functors $I_A^B$ and $P_A^B$, as well as for the
terminology ``pseudo-forgetful functor''. The proof of the fact that duality intertwines the functors
$I_A^B$ and $P_A^B$ is also taken from \cite{KV}. Section \ref{35} is inspired by \cite{DeB}, the results
found therein being due to J. Bernstein. The notion of the completion of a non-degenerate module is
introduced in \cite{Be2}.

\chapter[T.d. spaces and groups]{Totally disconnected spaces and groups} \label{chaptd}

In this chapter, we develop the study of locally compact totally disconnected (t.d. for short)
topological spaces (and more particularly groups), i.e., such that each point of the space admits a
neighborhood basis of compact open sets. Indeed, $p$-adic reductive groups are equipped with such a
topology, coming from that of the fields over which they are defined. The appropriate notion of a
``smooth'' function on such spaces is that of a locally constant function. A smooth function with compact
support is then a linear combination of characteristic functions of disjoint compact open sets. This fact
gives functional analysis on these spaces a purely algebraic character, which makes it much simpler than
in the case of differentiable manifolds, for example. The same is true for the theory of sheaves (of
abelian groups) on such a space. We also study continuous actions of t.d. topological groups on t.d.
topological spaces. Here again, the theory is free from the pathologies encountered in the general case,
or even in the differential framework. A fundamental tool in the study of locally compact topological
groups is the Haar measure. A Haar measure is a measure invariant under the action by left translation of
the group on itself. The proof of the existence and uniqueness up to a scalar factor of such a measure on
locally compact totally disconnected topological groups is quite elementary. The action by right
translation of the group on the Haar measure defines a ``modular'' function (or character) on the group.
The calculation of this function can be carried out in terms of indices of compact open subgroups.

Another fundamental notion for the rest of the theory is that of the convolution of distributions. It is
this convolution that will allow us to define the algebras playing a role in the representation theory of
t.d. topological groups. The space of compactly supported distributions on such a group $G$ becomes an
algebra once equipped with the convolution product, whose identity element is the Dirac distribution at
the identity element of the group. The Haar measures on the compact open subgroups of $G$ provide a
directed system of idempotents for this algebra. The Hecke algebra of $G$ is the subalgebra consisting of
compactly supported distributions fixed by left and right convolution by one of these idempotents. We
obtain an algebra with idempotents (without identity), denoted $\caH(G)$. The theory of Chapter
\ref{algaidem} applies to this algebra. We explicitly describe the completion and the center of the
category of non-degenerate $\caH(G)$-modules in terms of so-called essentially compactly supported
distributions. The importance of these results comes from the fact that the category of non-degenerate
$\caH(G)$-modules is equivalent to the category of smooth representations of $G$, which will be the
object of study in the next chapter.

\section{Topology of t.d. spaces}

In this section, we give some general results on the topology of locally compact totally disconnected
spaces.

If $X$ is a set and $Y$ a subset of $X$, we denote by \index[not]{ZZchiY@$\chi_Y$}$\chi_Y$ the
characteristic function of $Y$ in $X$. If $X$ is a topological space and $Y$ is a subset of $X$, we
denote by $\overline{Y}$ the closure of $Y$ in $X$.

\subsection{T.d. spaces}\label{recouv}

A locally compact totally disconnected topological space (t.d. space) is a Hausdorff topological space
such that each point admits a neighborhood basis of compact open sets.

A locally compact totally disconnected topological group is a topological group whose underlying
topological space is t.d. To do this, it suffices that the identity element $e$ admits a neighborhood
basis consisting of compact open subgroups. This condition is in fact also necessary (van Dantzig's
theorem). For convenience, one can take this a priori stronger condition as the definition of t.d.
topological groups.

\begin{lemme}
$(i)$ Let $X$ be a t.d. space and $Y$ a locally closed subset of $X$. Then $Y$ is a t.d. space for the
induced topology.

$(ii)$ If $K$ is a compact subset of $X$ and $\bigcup_\alpha U_\alpha$ is an open cover of $K$, then
there exists a finite subordinate cover of $K$ by disjoint compact open sets $V_i$ (i.e., for all $i$,
there exists $\alpha$ such that $V_i \subset U_\alpha$).
\end{lemme}

\begin{proof}
Recall that $Y \subset X$ is locally closed if for every point $y\in Y$, there exists a neighborhood $U$
of $y$ in $X$ such that $U \cap Y$ is closed in $U$. Here, we can assume $U$ is a compact open set, and
$(i)$ then follows quite easily from the definitions. Let us move on to $(ii)$. For each point $x$ in
some $U_\alpha$, let us choose a compact open neighborhood $V_{x,\alpha} \subset U_\alpha$. We then have
$K \subset \bigcup_{\alpha,x} V_{x,\alpha}$, and by the compactness of $K$, we can extract a finite
subcover $K \subset \bigcup_{i} V_{x_i,\alpha_i}$. Since the $V_{x_i,\alpha_i}$ are both open and closed,
the same is true for their pairwise intersections and symmetric differences. From this, we can easily
construct a finite cover of $K$ satisfying the desired properties in the statement.
\end{proof}

\begin{prop}
A t.d. space satisfies the \index[ter]{Baire category theorem}Baire category theorem, i.e., a countable
intersection of dense open sets is dense.
\end{prop}
\begin{proof}
This is true for any locally compact space (cf. \cite{Dug}, Theorem XI. 10.1).
\end{proof}

\subsection{Functions and distributions}
\index[ter]{distribution} \label{FetD}

Let $X$ be a t.d. space. We consider the following functional spaces on $X$:

\noindent \index[not]{CX@$\scrC^\infty(X)$}$\scrC^\infty(X)$, the space of locally constant
complex-valued functions,

\noindent \index[not]{DX@$\scrD(X)$}$\scrD(X)$, the subspace of $\scrC^\infty(X)$ of compactly supported
functions,

\noindent \index[not]{DX'@$\scrD'(X)$}$\scrD'(X)$, the (algebraic) dual of $\scrD(X)$, the space of
distributions (note that these spaces are not equipped with topologies).

\medskip

Recall that if $f$ is a complex-valued function on a t.d. space $X$, its
support\index[ter]{support!of a function}, denoted $\supp f$, is defined as the complement of the union
of the open sets $U$ of $X$ where the restriction of $f$ is zero.

\begin{lemme}\label{FetDlem}
Any locally constant compactly supported function on a t.d. space $X$ can be written as a finite sum of
scalar multiples of characteristic functions of disjoint compact open sets.
\end{lemme}

\begin{proof}
Let $f \in \scrD(X)$. For each point $x \in \supp f$, let us choose a compact open neighborhood $\caU_x$
where $f$ is constant. Since $\supp f$ is compact, we can extract a finite subcover from the cover by the
$\caU_x$, and according to $(ii)$ of Lemma \ref{recouv}, we can assume these open sets are disjoint. The
result is then clear.
\end{proof}

\begin{prop}
Let $X$ be a t.d. space and $U$ an open subset of $X$. Let $Z=X \setminus U$. We then have an exact
sequence
\begin{equation}\label{uxz}
 0 \rightarrow\scrD(U) \stackrel{i_U} {\longrightarrow}\scrD(X)  \stackrel{p_Z} {\longrightarrow}\scrD(Z)
\rightarrow 0.
\end{equation}
\end{prop}

\begin{proof}
The injection $i_U$ is given by extending functions on $U$ by $0$ outside of $U$. The map $p_Z$ is the
restriction of functions from $X$ to $Z$. Since the subset $Z$ is closed, we indeed obtain locally
constant and compactly supported functions in this way. To show that the restriction from $X$ to $Z$ is
surjective, we must show how to extend a function $f$ from $\scrD(Z)$ to the whole of $X$. Since the
support of $f$ (let us denote it $Z'$) is compact, and since $f$ is locally constant, we can cover $Z'$
by open sets of $X$ such that $f$ is constant on the intersection of $Z'$ and one of these open sets.
According to $(ii)$ of Lemma \ref{recouv}, there exists a finite cover of $Z'$ by disjoint compact open
sets of $X$ having this same property. We extend $f$ to a constant function on these compact open sets.
Outside these compact open sets, we extend $f$ by $0$. The exactness of the sequence (\ref{uxz}) is now
easy to verify: it is clear that $p_Z\circ i_U=0$, and if $p_Z(f)=0$, since $f$ is locally constant, $f$
is zero in a neighborhood of $Z$. Its support is therefore in $U$.
\end{proof}

\medskip

By duality, we obtain:

\begin{cor}\label{FetDcor}
The sequence:
\begin{equation}\label{zxu}
 0 \rightarrow\scrD'(Z) \stackrel{p_Z^*} { \longrightarrow}\scrD'(X)
  \stackrel{i_U^*} {\longrightarrow} \scrD'(U)\rightarrow 0
\end{equation}
is exact.
\end{cor}

If $T$ is a distribution on $X$, we define its support\index[ter]{support!of a distribution}, denoted
$\supp T$: it is the complement of the union of the open sets $U$ of $X$ such that $i_U^*(T)=0$. We will
say that $f$ (resp. $T$) is supported in a subset $Y$ of $X$ if $\supp f\subset Y$ (resp.
$\supp T\subset Y$). If $Z$ is a closed subset of $X$, then the exact sequence (\ref{zxu}) shows that any
distribution $T$ on $X$ supported in $Z$ can be written uniquely as $p_Z^*(T_0)$, $T_0 \in \scrD'(Z)$. We
will therefore often identify distributions on $X$ supported in $Z$ with distributions on $Z$.

We can now define another functional space on $X$:

\noindent \index[not]{EX'@$\scrE'(X)$}$\scrE'(X)$, the space of
\index[ter]{distribution!compactly supported}compactly supported distributions on $X$.
\medskip

Let us define a duality between the spaces $\scrC^\infty(X)$ and $\scrE'(X)$: if $T \in \scrE'(X)$, its
support $F$ is compact. Consider the exact sequence (\ref{zxu}) with $Z=F$. Then there exists a unique
distribution $T_0$ on $F$ such that $p_F^*(T_0)=T$. Since $F$ is compact, for any function
$f \in \scrC^\infty(X)$, $p_F(f)$ is in $\scrD(F)$ and we set:
\[
\langle T,f \rangle= \langle T_0,p_F(f) \rangle.
\]

When convenient, we will also use the following notation for the duality between $\scrD(X)$ and
$\scrD'(X)$ (resp. between $\scrC^\infty(X)$ and $\scrE'(X)$):
\[
\langle T,f \rangle=\int_X f\, dT= \int_X f(x)\, dT(x), \quad (f\in \scrD(X)), \; (T \in \scrD'(X)).
\]

We can easily extend the above results to functions with values in a $\bbC$-vector space $V$. We denote
by \index[not]{C(X,V)@$\scrC^\infty(X,V)$} $\scrC^\infty(X,V)$ (respectively $\scrD(X,V)$)
\index[not]{D(X,V)@$\scrD(X,V)$} the set of locally constant functions on $X$ with values in $V$
(respectively, locally constant with compact support). We have a canonical injection
\begin{align}\label{XVXV}
 \scrC^\infty(X)\otimes V &\rightarrow  \scrC^\infty(X,V)\\
  \nonumber           f \otimes v &\mapsto [x \mapsto f(x)v].
\end{align}
And similarly
\begin{align}\label{XVXV2}
 \scrD(X)\otimes V &\rightarrow  \scrD(X,V)\\
  \nonumber           f \otimes v &\mapsto [x \mapsto f(x)v].
\end{align}
But any function in $\scrD(X,V)$ is a finite sum of constant functions on disjoint compact open sets of
$X$ (Lemma \ref{FetDlem}), so it is clear that (\ref{XVXV2}) is surjective, and it is therefore an
isomorphism.

\subsection{Tensor product of distributions}\label{tens}

Let $X$ and $Y$ be t.d. spaces and let us equip $X\times Y$ with the product topology, which makes it a
t.d. space. Recall that according to Lemma \ref{FetDlem}, any function in $\scrD(X\times Y)$ is a linear
combination of characteristic functions of disjoint compact open sets. Since the subsets of the form
$U\times V$, where $U$ (resp. $V$) is a compact open set of $X$ (resp. of $Y$), form a basis of the
topology of $X \times Y$, we can, according to Lemma \ref{recouv}, $(ii)$, write any compact open set of
$X\times Y$ as a finite disjoint union of compact open sets of the form $U\times V$. Any function in
$\scrD(X\times Y)$ is therefore a linear combination of characteristic functions of disjoint compact open
sets of the form $U\times V$. This shows that $\scrD(X\times Y) \simeq \scrD(X)\otimes\scrD(Y) $, the
isomorphism being given by
\[
(f\otimes g) (x,y)= f(x) \, g(y),\quad x\in X,\,  y\in Y,\, f\in \scrD(X),\, g\in \scrD(Y).
\]

We can then define the tensor product of distributions: if $T_1 \in \scrD'(X)$, $T_2 \in \scrD'(Y)$, then
$T_1 \otimes T_2$ is defined on $\scrD(X\times Y)$ by:
\[
\langle T_1 \otimes T_2, f \otimes g \rangle =\langle T_1, f \rangle \langle  T_2,  g \rangle.
\]
On the other hand, the argument used above shows that Fubini's formulas are valid: for any
$f \in \scrD(X\times Y)$,
\begin{align*}
 \int_{X\times Y} f(x,y)\, d(T_1 \otimes T_2)(x,y)&=\int_X  \left( \int_Y f(x,y) \, dT_2(y) \right)  dT_1(x)= \int_Y  \left(   \int_X f(x,y) \, dT_1(x) \right)  dT_2(y).
\end{align*}
Indeed, the integrals reduce to finite sums.

\begin{prop}
We have
\[
\mathrm{Supp}(T_1\otimes T_2)=\mathrm{Supp}(T_1)\times\mathrm{Supp}(T_2).
\]
In particular, the tensor product of two compactly supported distributions is compactly supported.
\end{prop}
\begin{proof}
Let $\Omega$ be an open set of $X\times Y$ and, up to shrinking it, let us assume it is of the form
$U\times V$ where $U$ is an open set of $X$ and $V$ is an open set of $Y$. Let $T_{1,U}$ (resp.
$T_{2,V}$) denote the restriction of $T_1$ to $U$ (resp. $T_2$ to $V$). It is clear that the restriction
of $T_1\otimes T_2$ to $\Omega$ is equal to $T_{1,U}\otimes T_{2,V}$. The first assertion is deduced from
this by a sequence of tautologies. The second follows immediately.
\end{proof}

\section{Sheaves on a t.d. topological space}\label{faiscetoptd}

\subsection{Generalities on sheaf theory}\label{genFaisc}

Let $X$ be a topological space and let $\caT$ \index[not]{T@$\caT$} be the topology of $X$, i.e., the set
of open subsets of $X$. A \emph{basis for the topology} of $X$ is a subset $\caT_0 \subset \caT$
\index[not]{T_0@$\caT_0$} such that for any $x \in X$ and for any open set $U$ containing $x$, there
exists an open set $U_0 \in \caT_0$ containing $x$ such that
\[
U_0 \subset U.
\]
Of course, a basis determines the topology. We will always assume that $\caT_0$ is stable under finite
intersections (this is the case for $\caT$, and for $\caT_c$ when $X$ is t.d.).

We will introduce the notions of presheaf and sheaf (of abelian groups) on $X$, with basis $\caT_0$. This
slightly generalizes the definitions given in most treatises on sheaf theory (for example \cite{Faisc}),
which only consider the case $\caT_0=\caT$. This greater generality is in fact only apparent, since we
will see that any sheaf with basis $\caT_0$ naturally extends to a sheaf with basis $\caT$ (this is not
true for presheaves), but this flexibility in the definitions is convenient when $X$ is a t.d. space, for
which the choice $\caT_0=\caT_c$, the set of compact open subsets of $X$, is particularly well-suited.

\begin{defi}
A presheaf $\caF$ (of abelian groups) on $X$, with basis $\caT_0$, consists, for each open set $U$ in
$\caT_0$, of an abelian group $\caF(U)$ (if $\emptyset \in \caT_0$, then we require that
$\caF(\emptyset)=\{0\}$), and for any pair $U \subset V$ of open sets in $\caT_0$, of a restriction
morphism
\[
\rho_{V,U}\colon \caF(V) \rightarrow \caF(U)
\]
satisfying $\rho_{U,U}=\Id_{\caF(U)}$ and when $U\subset V\subset W$ are elements of $\caT_0$,
$\rho_{V,U}\circ \rho_{W,V}=\rho_{W,U}$. When we assume that the $\caF(U)$ are equipped with additional
structures (rings, vector spaces over a field $k$, etc.), and that the restriction morphisms preserve
these structures, we speak of presheaves of rings, of $k$-vector spaces, etc.
\end{defi}

An equivalent and slightly more sophisticated way to define presheaves is the following: let us again
denote by $\caT_0$ the category whose objects are the elements of $\caT_0$, and whose morphisms are given
by:
\[
\Hom_{\caT_0}(U,V)=\emptyset \quad \text{ if } U \text{ is not included in } V,
\]
\[
\Hom_{\caT_0}(U,V)\quad  \text{is the inclusion of } U \text { into } V  \text{ if } U \text{ is
included in } V.
\]

A presheaf of abelian groups (resp. of rings, of $k$-vector spaces, etc.) is then a (contravariant)
functor from $\caT_0$ to the category of abelian groups (resp. of rings, of $k$-vector spaces, etc.).

The presheaves (of abelian groups) on $X$, with basis $\caT_0$ (we will now simply say presheaf on $X$),
form a category, denoted $p\caS h (X)$ \index[not]{pShX@$p\caS h (X)$} whose morphisms are defined as
follows. Let $\caF_1$ and $\caF_2$ be two presheaves on $X$ with basis $\caT_0$. A morphism $\phi$
between $\caF_1$ and $\caF_2$ is the data, for each $U \in \caT_0$, of a group morphism
\[
\phi_U\colon \caF_1(U) \rightarrow \caF_2(U).
\]
It is moreover required that if $U \subset V$ are open sets in $\caT_0$,
\begin{equation}\label{morfai}
\phi_U \circ \rho_{V,U}= \rho_{V,U}  \circ \phi_V.
\end{equation}

\begin{defi}
A presheaf $\caF$ is a \emph{sheaf} if the following condition is satisfied:

\medskip

\noindent {\bf (F)}: let $V \in \caT_0$ and $\{V_i\}_{i\in I}$ be a cover of $V$ by open sets in
$\caT_0$. Suppose that for each $i \in I$ an element $f_i$ of $\caF(V_i)$ is given such that for any
$i,j$ in $I$ and for any open set $V_{ij} \subset  V_i \cap V_j$ in $\caT_0$, we have
\[
\rho_{V_i,V_{ij}}(f_i)=\rho_{V_j,V_{ij}}(f_j).
\]
Then there exists a unique $f \in \caF(V)$ such that for all $j \in I$, $\rho_{V,V_{j}}(f)=f_j$.
\end{defi}

The sheaves on $X$ form a category (the morphisms are the presheaf morphisms), denoted $\caS h (X)$
\index[not]{ShX@$\caS h (X)$}. Let us now define the stalk $\caF_x$ \index[not]{Fx@$\caF_x$} of a
presheaf $\caF$ on $X$. The set of open sets $U$ in $\caT_0$ containing $x$, denoted $\caV_{\caT_0}(x)$
\index[not]{V_T0x@$\caV_{\caT_0}(x)$}, is equipped with the partial order:
\[
V \leq U  \text{ if and only if }  U \subset V.
\]
This set is \emph{directed}, i.e., if $U$ and $V$ are open sets in $\caT_0$ containing $x$, there exists
an open set $W$ in $\caT_0$ such that $W \geq V$, $W \geq U$ (this simply follows from the fact that
$\caT_0$ is a basis for the topology). The morphisms $\rho_{V,U}$, $U\subset V$, $U,V \in \caT_0$, form
an inductive system of abelian groups (see Examples \ref{colimites}). We set
\[
\caF_x= \varinjlim_{U \in \caV_{\caT_0}(x)} \caF(U).
\]
This inductive limit can be described explicitly: let us form the disjoint union
$\coprod_{U \in \caV_{\caT_0}(x)} \caF(U)$ which we equip with the following equivalence relation:
$s\in \caF(U)$ is equivalent to $t \in \caF(V)$ if there exists $W \in  \caV_{\caT_0}(x)$,
$W \subset U\cap V$ such that $\rho_{U,W}(s)=\rho_{V,W}(t)$. The inductive limit $\caF_x$ is then the set
of equivalence classes, equipped with the group structure and the canonical morphisms
$\rho_{U,x}: \caF(U)\rightarrow \caF_x$.

We will now define a functor
\[
p\caS h(X) \longrightarrow \caS h (X), \qquad \caF \mapsto \caF^+.
\]
Let $x \in X$ and $U \in \caT_0$, with $x \in U$ and let
\[
\rho_{U,x}\colon \caF(U) \rightarrow \caF_x
\]
be the canonical projection. Let us set
\[
\widehat{\caF} =  \coprod_{x \in X} \caF_x.
\]
This is the disjoint union of the $\caF_x$. The space $\widehat{\caF}$ \index[not]{Fhat@$\widehat{\caF}$}
is called the \emph{étale space} \indexter{etale space@étale space} over the presheaf $\caF$.

For any open set $U$ of $X$ (not necessarily in $\caT_0$), a \emph{section} \indexter{section} of
$\widehat{\caF}$ over $U$ is a map
\[
s\colon U \rightarrow  \widehat{\caF}
\]
such that for all $x \in U$, $s(x) \in \caF_x$ and such that there exists a cover of $U$ by open sets
$(U_i)_{i \in I}$ in $\caT_0$, and elements $f_i \in \caF(U_i)$ satisfying:
\begin{equation} \label{section}
s(x)=\rho_{U_i,x}(f_i), \; (\forall x \in U_i\cap U).
\end{equation}

Let $\widehat{\caF}$ be the étale space for the presheaf $\caF$. For any $U \in \caT_0$, let us define
$\caF^+(U)$ as the set of sections of $\widehat{\caF}$ over $U$. If $U \subset V$ are open sets in
$\caT_0$, let us define the restriction morphism $\rho^+_{V,U}$ by $\rho^+_{V,U}(s)=s_{|U}$,
$s\in \caF^+(V)$. We must ensure that this is indeed a section of $\widehat{\caF}$ over $U$, in the sense
we defined above. By definition, there exists a cover of $V$ by open sets $(V_i)_{i \in I}$ in $\caT_0$,
and elements $f_i \in \caF(V_i)$ satisfying
\[
s(x)=\rho_{V_i,x}(f_i), \; (\forall x \in V_i\cap V).
\]
Let us take a cover $(U_j)_{j \in J}$ of $U$ such that for all $j\in J$, $U_j \in \caT_0$ and
$U_j \subset V_i$ for some $i \in I$. Let us then set $g_j=\rho_{V_i,U_j}(f_i)$. We have for all
$x \in U_j$,
\[
\rho_{U_j,x}(g_j)=\rho_{U_j,x}( \rho_{V_i,U_j}(f_i))=\rho_{V_i,x}(f_i)=s(x)=s_{|U}(x),
\]
which shows that $ s_{|U}$ is indeed a section. It is clear that $\caF^+$ is a presheaf of abelian
groups, the group structure on the $\caF^+(U)$ being clear.

Let us now show that $\caF^+$ is indeed a sheaf by verifying property {\bf (F)}. Let $V \in \caT_0$ and
$\{V_i\}_{i\in I}$ be a cover of $V$ by open sets in $\caT_0$. Suppose that for each $i \in I$ an element
$s_i$ of $\caF^+(V_i)$ is given such that for any $i,j$ in $I$ and for any open set
$V_{ij} \subset V_i \cap V_j$ in $\caT_0$, we have
\[
\rho^+_{V_i,V_{ij}}(s_i)=\rho^+_{V_j,V_{ij}}(s_j).
\]
For any $x \in V$, let us choose $V_i$ containing $x$ and set $s(x)=s_i(x)$. It is clear that $s(x)$ does
not depend on the choice of $V_i$ containing $x$, because if $x \in V_i \cap V_j$, we take an open set
$V_{ij} \in \caT_0$ such that $x \in V_{ij}\subset V_i \cap V_j $, and we then have:
\[
s_i(x)=(s_i)_{|V_{ij}}(x)= \rho^+_{V_i,V_{ij}}(s_i)(x)=\rho^+_{V_j,V_{ij}}(s_j)(x)=(s_j)_{|V_{ij}}(x)=s_j(x).
\]

Let us verify that $s$ is a section of $\widehat{\caF}$ over $V$. Each $s_i$, $i \in I$, is a section,
and therefore there exists a cover $(U^i_j)_{j \in J_i}$ of $V_i$ by open sets in $\caT_0$ and elements
$f^i_j \in \caF(U_j^i)$ such that $s_i(x)=\rho_{U_j^i,x}(f_j^i)$. We then have
$s(x)=\rho_{U_j^i,x}(f_j^i)$ for all $x \in U_j^i$, which shows that $s$ is a section. The existence
property in the definition of sheaves is therefore established. It remains to show uniqueness. Suppose
that two sections $s$ and $s'$ in $\caF^+(V)$ satisfy
\[
\rho_{V,V_i}(s)=\rho_{V,V_i}(s')=s_i, \quad (\forall i \in I).
\]
We immediately have, for all $x \in V$, by choosing a $V_i$ containing $x$, $s(x)=s'(x)=s_i(x)$, and
therefore $s=s'$.

\bigskip

Let $\phi\colon \caF \rightarrow \caG$ be a morphism of presheaves on $X$. Let us now define
$\phi^+\colon \caF^+ \rightarrow \caG^+$. Let $s \in \caF^+(U)$. By definition, there exists a cover
$(U_i)_{i \in I}$ of $U$ by open sets in $\caT_0$ and elements $f_i \in \caF(U_i)$ such that
$s(x)=\rho_{U_i,x}(f_i)$. Let us set, for all $i \in I$, $g_i=\phi_{U_i}(f_i)$ and
$t(x)=\rho_{U_i,x}(g_i)$. Standard verifications show that this is independent of the choices made and
that condition (\ref{morfai}) is satisfied for the $\rho^+_{V,U}$ and the $\phi_U^+$. We also leave to
the reader the verification of the fact that
\[
\caF \mapsto \caF^+, \qquad [\caF \stackrel{\phi}{\rightarrow} \caG] \mapsto [\caF^+ \stackrel{\phi^+}{\longrightarrow} \caG^+]
\]
is a functor.

\begin{prop}
The functor $\caF \mapsto \caF^+$ from $p \caS h(X)$ to $\caS h (X)$ is the left adjoint of the forgetful
functor from $\caS h (X)$ to $p \caS h(X)$, i.e., for any $\caF\in p \caS h(X)$ and for any
$\caG \in \caS h (X)$, we have a natural isomorphism:
\[
\Hom_{p \caS h(X)}(\caF,\caG)\simeq  \Hom_{\caS h(X)}(\caF^+,\caG)
\]
\end{prop}

\begin{proof}
For any presheaf $\caF$ on $X$, let us define a morphism
\[
\iota^\caF \colon  \caF \rightarrow \caF^+.
\]
Let $U \in \caT_0$ and $f \in \caF(U)$. Let us set
\[
\iota^\caF_U  (f)=s_f,
\]
where $s_f\colon U \rightarrow  \coprod_{x \in U} \caF_x  $ is the section given by
$s_f(x)=\rho_{U,x}(f)$. This is indeed a section. This can be seen by choosing $U$ itself as a cover of
$U$. A straightforward verification shows that for any presheaf morphism
$\caF \stackrel{\phi}{\rightarrow} \caG$, the following diagram commutes
\[
\xymatrix{
\caF  \ar[r]^{\phi}\ar[d]_{\iota^\caF}
&  \caG \ar[d]_{\iota^\caG}  \\
\caF^+\ar[r]_{ \phi^+}&  \caG^+    \\ }
\]

Let $\psi: \caF \rightarrow \caG$ be a morphism of presheaves on $X$, where $\caG$ is a sheaf. We will
define $\Psi: \caF^+ \rightarrow \caG$ factoring $\psi$ through $\iota^\caF$. Let $s \in \caF^+(U)$, and
the $U_i$, $f_i$ as in (\ref{section}). Let us set $f'_i= \psi_{U_i}(f_i) \in \caG(U_i)$. We then have
according to (\ref{morfai}), for any $i$ and $j$ in $I$:
\[
\rho_{U_i,U_i\cap U_j}(f'_i)=\rho_{U_i,U_i\cap U_j}(\psi_{U_i}(f_i))= \psi_{U_i\cap U_j}(\rho_{U_i,U_i\cap U_j}(f_i)),
\]
and similarly
\[
\rho_{U_j,U_i\cap U_j}(f'_j)=\rho_{U_j,U_i\cap U_j}(\psi_{U_j}(f_j))= \psi_{U_i\cap U_j}(\rho_{U_j,U_i\cap U_j}(f_j)).
\]

We will now show that
\begin{equation} \label{SE2}
\rho_{U_i,U_i\cap U_j}(f'_i)= \rho_{U_j,U_i\cap U_j}(f'_j).
\end{equation}
We have for all $x \in U_i\cap U_j$:
\[
s(x)=\rho_{U_i,x}(f_i)=\rho_{U_j,x}(f_j) \quad \text{in } \caF_x.
\]
The presheaf morphism $\psi \colon \caF \to \caG$ induces a morphism on stalks
$\psi_x \colon \caF_x \to \caG_x$ such that $\psi_x(\rho_{U,x}(f)) = \rho_{U,x}(\psi_U(f))$ for any
$f \in \caF(U)$. Applying $\psi_x$ to the equality above, we obtain:
\[
\rho_{U_i,x}(\psi_{U_i}(f_i)) = \rho_{U_j,x}(\psi_{U_j}(f_j)) \quad \text{in } \caG_x,
\]
which means $\rho_{U_i,x}(f'_i) = \rho_{U_j,x}(f'_j)$. We can rewrite this as:
\[
\rho_{U_i\cap U_j,x}(\rho_{U_i,U_i\cap U_j}(f'_i)) = \rho_{U_i\cap U_j,x}(\rho_{U_j,U_i\cap U_j}(f'_j)),
\]
hence
\begin{equation} \label{SE3}
\rho_{U_i\cap U_j,x}(\rho_{U_i,U_i\cap U_j}(f'_i)-\rho_{U_j,U_i\cap U_j}(f'_j))=0 \quad \text{in } \caG_x.
\end{equation}

Now we have the following result:
\begin{lemme}
Let $\caF$ be a sheaf on $X$, $U\in \caT_0$ and $f\in \caF(U)$ such that $\rho_{U,x}(f)=0 \in \caF_x$ for
all $x \in U$. Then $f=0$.
\end{lemme}
\begin{proof}
By definition of the inductive limit, for any $x \in U$, there exists an open set $U_x \in \caT_0$
containing $x$ and contained in $U$ such that $\rho_{U,U_x}(f)=0$. The open sets $(U_x)_{x \in U}$ form a
cover of $U$. Property {\bf (F)} of sheaves then shows that $f=0$.
\end{proof}

We now complete the proof of the proposition. Applying the lemma to the sheaf $\caG$ and the section
$\rho_{U_i,U_i\cap U_j}(f'_i)-\rho_{U_j,U_i\cap U_j}(f'_j) \in \caG(U_i\cap U_j)$, equation (\ref{SE3})
implies (\ref{SE2}). Property {\bf (F)} for the sheaf $\caG$ then implies that there exists a unique
$f' \in \caG(U)$ such that $\rho_{U,U_i}(f')=f'_i$. Let us then set $\Psi_U(s)=f'$.

Let $f \in \caF(U)$, $ \iota_U^\caF(f)=s_f$ and $f'=\Psi_U(s_f)$ obtained as above. We then have:
\[
\Psi_U(\iota_U^\caF(f))=\Psi_U(s_f)=f'
\]
and $\rho_{U,U_i}(f')=f'_i=\psi_{U_i}(f_i)$ for all $i \in I$. This shows that
$\Psi\circ \iota^\caF= \psi$.

Conversely, if $\Psi$ is in $\Hom_{\caS h(X)}(\caF^+,\caG)$, we define $\phi$ in
$ \Hom_{p \caS h(X)}(\caF,\caG)$ by
\[
\phi_U(f)=\Psi_U(s_f), \quad (f \in \caF(U)).
\]
We then immediately see that the operations $\phi \mapsto \Psi$ and $\Psi \mapsto \phi$ described above
are inverse to each other, using the fact that $\Psi$ is uniquely determined by the data of the
$\Psi(s_f)$, $f \in \caF(U)$, $U \in \caT_0$.
\end{proof}

\begin{rmq}
The sheaf $\caF^+$ and the morphism $\iota^\caF$ are solutions to the universal problem (on the right)
posed by $\caF$ and the forgetful functor from sheaves to presheaves (see \ref{propuniv}). There is
therefore uniqueness of $(\caF^+, \iota^\caF)$ up to a unique isomorphism. In particular, if $\caF$ is a
sheaf, $\caF$ is naturally isomorphic to $\caF^+$ (the isomorphism is $ \iota^\caF$). A sheaf is
therefore identified with the sheaf of sections of its étale space $\widehat{\caF}$. When $\caF$ is a
sheaf, we will no longer distinguish between $\caF$ and $\caF^+$ in what follows.
\end{rmq}

Let us end this section by noting that the category of sheaves on $X$ with basis $\caT_0$ is equivalent
to the category of sheaves with basis $\caT$. One obviously passes from a sheaf with basis $\caT$ to a
sheaf with basis $\caT_0$ by restriction. In the other direction, if $\caF$ is a sheaf with basis
$\caT_0$ and if $U$ is any open set of $X$, we define $\caF^+(U)$ as the set of sections
\[
s\colon U \rightarrow \coprod_{x \in U} \caF_x
\]
exactly as in (\ref{section}). This defines a sheaf with basis $\caT$ extending the sheaf $\caF^+$ with
basis $\caT_0$ constructed in the proposition above (the proof that it is a sheaf is the same as for
$\caT_0$). Let us denote it (temporarily) $\caF^+_\caT$. Since $\caF$ is a sheaf, according to the remark
above, $\caF$ is naturally isomorphic to $\caF^+$. This shows that the functor which consists of starting
from a sheaf $\caF$ with basis $\caT_0$, constructing the sheaf $\caF^+_\caT$ and restricting it to the
basis $\caT_0$ is naturally isomorphic to the identity functor of the category of sheaves on $X$ with
basis $\caT_0$. In the other direction, let us start from a sheaf $\caG$ on $X$ with basis $\caT$. It is
naturally isomorphic to the sheaf $\caG^+$. Let us restrict $\caG$ to the basis $\caT_0$ to obtain
$\caG_{\caT_0}$, then let us form $(\caG_{\caT_0})^+$. It is clear that this sheaf is equal to $\caG^+$.
This shows that we obtain a functor isomorphic to the identity functor of the category of sheaves on $X$
with basis $\caT$ and finishes establishing the equivalence of categories.

\subsection{Support of a section}
\indexter{support!of a section}

We retain the notation of the previous section, in particular $\caF$ is a sheaf on $X$ with basis
$\caT_0$.
\begin{defi}
Suppose that $s\colon U \rightarrow \widehat{\caF}$ is a section. Its support, denoted
$\mathrm{Supp}(s)$, is the set of $x \in U$ such that $s(x) \neq 0$.
\end{defi}

\begin{prop}
The support of a section $s\colon U \rightarrow \widehat{\caF}$ of the sheaf $\caF$ is closed in $U$.
\end{prop}

\begin{proof}
Suppose that $x \notin \mathrm{Supp}(s)$. By definition of a section, $x$ has an open neighborhood
$W\in \caT_0$ such that there exists an element $f \in \caF(W)$ with $s(y)=\rho_{W,y}(f)$ for all
$y \in W$. Since by hypothesis $\rho_{W,x}(f)=s(x)=0$, by definition of the inductive limit, there exists
a neighborhood $V \in \caT_0$ of $x$ contained in $W$ such that $\rho_{W,V}(f)=0$. Then for all
$y \in V$,
\[
s(y)=\rho_{V,y}\circ \rho_{W,V}(f)=0.
\]
This shows that $V$ does not intersect the support of $s$. The complement of this support is therefore
open.
\end{proof}

\medskip

By the identification of $\caF$ and $\caF^+$ we can speak of the support of an element of $\caF(U)$. If
$U$ is any open set, we denote by $\caF_c(U)$ the set of elements of $\caF(U)$ with compact support.

\subsection{Sheaf of modules}

Let $\caR$ be a sheaf of rings on $X$, with basis $\caT_0$ (we say that $(X,\caR)$ is a
\index[ter]{ringed space}\emph{ringed space}). A sheaf of $\caR$-modules on $X$ is a sheaf $\caM$ on $X$
such that for any open set $U \in \caT_0$, $\caM(U)$ is a (left) module over $\caR(U)$. If $U \subset V$
are open sets in $\caT_0$, then the restriction map $\rho_{V,U}: \caR(V) \rightarrow \caR(U)$ is a ring
morphism by which any $\caR(U)$-module becomes an $\caR(V)$-module. We then require that the restriction
morphisms $\rho_{V,U}: \caM(V)\rightarrow \caM(U)$ be morphisms of $\caR(V)$-modules.

\subsection{Characterization of sheaves on a t.d. space}

Let $X$ be a t.d. space. The set $\caT_c$ of compact open subsets of $X$ is a basis for the topology of
$X$. We will consider in what follows the sheaves (of abelian groups) with basis $\caT_c$. The very
particular nature of the locally compact totally disconnected topology gives them a very simple
structure, as shown by the following proposition.

\begin{prop}
Let $X$ be a t.d. space and let $\caF$ be a sheaf (of abelian groups) with basis $\caT_c$. Then there
exists an abelian group $\caF_c(X)$ and a sheaf $\caF_1$ isomorphic to $\caF$ such that

\begin{itemize}
\item[1.] for all $U \in \caT_c$, $\caF_1(U)\subset \caF_c(X)$,

\item[2.] if $U\subset V$ are open sets in $\caT_c$, then $\caF_1(U)\subset \caF_1(V)$,

\item[3.] if $U\subset V$ are open sets in $\caT_c$, then for all $f \in \caF_1(U)$, $\rho_{V,U}(f)=f$
while $\rho_{V,V\setminus U}(f)=0$,

\item[4.] $\caF_c (X)=\bigcup_{U \in \caT_c} \caF_1(U)$.
\end{itemize}
\end{prop}

\begin{proof}
Since $\caF$ is a sheaf, $\caF$ is isomorphic to $\caF^+$, the sheaf of sections of the étale space
$\widehat{\caF}$. Let $\caF_c(X)$ be the space of compactly supported sections in $\caF^+(X)$. If $U$ is
a compact open set of $X$, we embed $\caF^+(U)$ into $\caF_c(X)$ by extending the sections
$s \in \caF^+(U)$ by $0$ outside of $U$. Since $X \setminus U$ is open, it is clear that this indeed
defines a section in $\caF^+(X)$, and it is supported in $U$, hence compact. The other assertions of the
proposition are then easily verified.
\end{proof}

\begin{exemple}
Let $\scrC^\infty_X$ be the sheaf whose space of sections over an open set $U$ is the space
$\scrC^\infty(U)$ of locally constant functions on $U$. The abelian group given by the proposition is in
this case $\scrD(X)$, the space of locally constant compactly supported functions on $X$. Note that the
$\scrC^\infty(U)$ are $\bbC$-algebras (for pointwise multiplication of functions), in particular, they
are rings, and therefore $(X,\scrC^\infty_X)$ is a ringed space. Note also that the algebra $\scrD(X)$ is
generally not unital, unless $X$ is compact, in which case $\scrC^\infty(X)=\scrD(X)$, whose identity is
$\chi_X$. On the other hand, it is an algebra with idempotents (see \ref{sec_idem}), since each compact
open set $U$ of $X$ gives an idempotent $\chi_U$ and these $\chi_U$ form a directed system of idempotents
(see Remark \ref{idem}, 2).
\end{exemple}

\subsection{Sheaves of $\bbC$-vector spaces}\label{FCev}

Let $X$ be a t.d. space and let $\scrC^\infty_X-\caM od (X)$
\index[not]{CXmod@$\scrC^\infty_X-\caM od (X)$} be the category of sheaves of $\scrC^\infty_X$-modules.

\begin{prop}
Let $X$ be a t.d. space and $\caM$ a sheaf of $\bbC$-vector spaces on $X$ with basis $\caT_c$. Then
$\caM$ is naturally a sheaf of $\scrC^\infty_X$-modules.
\end{prop}

\begin{proof}
We identify $\caM$ and the sheaf of sections of the associated étale space $\widehat{\caM}$. Let
$U \in \caT_c$, $s: U \rightarrow \widehat{\caM}$ be a section and $f \in \scrC^\infty(U)$. Then
$x \mapsto f(x)s(x)$ is defined, and since $f$ is locally constant, it is clear that it is a section.
Thus $\caM$ becomes a $\scrC^\infty_X$-module.
\end{proof}

\begin{rmq}
It is obvious that conversely, a sheaf of $\scrC^\infty_X$-modules is naturally a sheaf of $\bbC$-vector
spaces on $X$. Indeed, $\bbC$ naturally embeds into all the $\scrC^\infty(U)$ by associating to a scalar
in $\bbC$ the constant function equal to this scalar on $U$. Thus the notions of a sheaf of $\bbC$-vector
spaces on $X$ and of a sheaf of $\scrC^\infty_X$-modules are equivalent (i.e., the respective categories
are equivalent). Moreover, by virtue of the last paragraph of \ref{genFaisc}, this remains true for
sheaves with basis $\caT$.
\end{rmq}

Our goal is now to prove the following theorem. Recall that $\scrD(X)$ is an algebra with idempotents,
with the $\{\chi_U\}$, $U$ compact open, as a system of idempotents. The algebra $\scrC^\infty(X)$, on
the other hand, is a unital algebra, with identity $\chi_X$. The notion of a non-degenerate module
defined in the previous chapter for algebras with idempotents easily extends to the case of an algebra
equipped with a directed system of idempotents (cf. the third remark of \ref{idem}) and $\scrC^\infty(X)$
also admits the $\{\chi_U\}$, $U$ compact open, as a directed system of idempotents. In what follows, we
consider non-degenerate modules over $\scrC^\infty(X)$ relative to this directed system of idempotents.

\begin{thm}\label{faisceqcat}
For any $\caM \in \scrC^\infty_X-\caM od (X)$, $\caM_c(X)$ is a non-degenerate $\scrD(X)$-module and
$\caM \mapsto \caM_c(X)$ realizes an equivalence of categories between $\scrC^\infty_X-\caM od (X)$ and
$\caM(\scrD(X))$.
\end{thm}

Note first that we can replace $\scrD(X)$ in the statement by $\scrC^\infty(X)$.

\begin{lemme}
Let $M$ be a non-degenerate $\scrD(X)$-module. Then the action of $\scrD(X)$ on $M$ extends uniquely to
an action of $\scrC^\infty(X)$ which makes $M$ a non-degenerate $\scrC^\infty(X)$-module.
\end{lemme}

\begin{proof}
Let $m \in M$ and let $U$ be a compact open set such that $\chi_U \cdot m=m$. Let us define, for
$f \in \scrC^\infty(X)$,
\[
f\cdot m=(f\chi_U)\cdot m.
\]
One easily verifies that the right-hand side is independent of the choice of $U$ and that this definition
makes $M$ a non-degenerate $\scrC^\infty(X)$-module, for the directed system $\{\chi_U\}$, $U$ compact
open, and does so uniquely.
\end{proof}

\medskip

\begin{proof}[Proof of the theorem]
Let $s \in \caM_c(X)$, and let $U$ be a compact open set of $X$ such that $\supp s \subset U$. The
section $s$ is then fixed by the idempotent $\chi_U$. This shows that $\caM_c(X)$ is a non-degenerate
$\scrD(X)$-module.

Let us define the inverse functor: for any non-degenerate $\scrD(X)$-module $M$, and for any compact open
set $U$ of $X$, let us define $\caM(U)=\chi_U \cdot M$. If $U \subset V$ are compact open sets, the
restriction morphism
\[
\rho_{V,U}: \caM(V)=\chi_V \cdot M \rightarrow  \caM(U)=\chi_U \cdot M
\]
is given by $m \mapsto \chi_U\cdot m$. We will verify that $\caM$ is indeed a sheaf.

First, we have a simple characterization of the modules $\caM(U)$, $U$ compact open:
\[
\caM(U)=\{ m\in M\mid \chi_U\cdot m=m  \}.
\]
The fact that $\caM$ is a presheaf is simply the translation of the fact that if $U \subset V$ are two
compact open sets, then $\chi_U\chi_V=\chi_U$. Let us now verify property {\bf (F)} of sheaves. Let
$V \in \caT_c$ and $\{V_i\}_{i\in I}$ be a cover of $V$ by open sets in $\caT_c$. Suppose that for each
$i \in I$ an element $m_i$ of $\caM(V_i)$ is given such that for any $i,j$ in $I$ and for any open set
$V_{ij} \subset V_i \cap V_j$ in $\caT_c$, we have
\[
\rho_{V_i,V_{ij}}(m_i)= \chi_{V_{ij}}\cdot m_i=\rho_{V_j,V_{ij}}(m_j)= \chi_{V_{ij}}\cdot m_j.
\]

We want to show that there exists a unique $m \in \mathcal{M}(V)$ such that for all $j \in I$,
\[
\rho_{V,V_{j}}(m)= \chi_{V_j}\cdot m= m_j.
\]
Let us fix $k$ and $l$ in $I$, and show that we can replace $V_k$ and $V_l$ in the cover by the single
compact open set $V_0=V_k \cup V_l$. More precisely, we will see that there exists a unique element
$m_0 \in \mathcal{M}(V_0)$ such that
\begin{equation}\label{V0Vk1}
\rho_{V_0,V_k}(m_0)=m_k, \quad \rho_{V_0,V_l}(m_0)=m_l.
\end{equation}
Moreover, for all $j \in I$,
\begin{equation}\label{V0Vk2}
\rho_{V_0,V_0\cap V_j}(m_0)=\rho_{V_j,V_0\cap V_j}(m_j).
\end{equation}

Let us set
\[
m_{kl}= \chi_{V_k\cap V_l}\cdot m_k= \chi_{V_k \cap V_l}\cdot m_l.
\]
Since $\chi_{V_k}\chi_{V_l}=\chi_{V_k\cap V_l}$ and $\chi_{V_k}\cdot m_k=m_k$,
$\chi_{ V_l}\cdot m_l=m_l$, we obtain
\[
m_{kl}=  \chi_{V_l}\cdot m_k= \chi_{V_k}\cdot m_l.
\]
Let us set $V_0=V_k \cup V_l$. Since $\chi_{V_0}=\chi_{V_k}+\chi_{V_l}-\chi_{V_k\cap V_l}$, we then have
\[
m_{kl}=  \chi_{V_k \cap V_l}\cdot m_{kl}=  \chi_{V_l}\cdot m_{kl}= \chi_{V_k}\cdot m_{kl}=\chi_{V_0}\cdot m_{kl}.
\]
Let us define $m_0=m_k+m_l-m_{kl}$. We have $\chi_{V_k} \cdot m_0=m_k$, $\chi_{V_l}\cdot m_0=m_l$ and
$\chi_{V_k\cap V_l}\cdot m_0 =m_{kl}$, hence $\chi_{V_0}\cdot m_0=m_0$. Thus $m_0 \in \mathcal{M}(V_0)$,
and
\[
\rho_{V_0,V_k}(m_0)=m_k, \quad \rho_{V_0,V_l}(m_0)=m_l,
\]
which establishes (\ref{V0Vk1}). For all $j \in I$, we have
\[
\chi_{V_j}\cdot m_0= \chi_{V_j}\cdot (m_k+m_l-m_{kl})=( \chi_{V_k}+  \chi_{V_l}-\chi_{V_k\cap V_l})\cdot m_j=\chi_{V_0}\cdot m_j,
\]
which establishes (\ref{V0Vk2}).

For uniqueness, note that if $m_0'$ satisfies (\ref{V0Vk1}), then
\[
m_0'=\chi_{V_0}\cdot m_0'= (\chi_{V_k}+  \chi_{V_l}-\chi_{V_k\cap V_l})\cdot m_0'=m_k+m_l-m_{kl}=m_0.
\]

Having established this pairwise gluing, we proceed by induction. For any finite subset $J \subset I$, we
can construct a unique element $m_J \in \mathcal{M}(\bigcup_{j \in J} V_j)$ such that
$\chi_{V_j} \cdot m_J = m_j$ for all $j \in J$. Since $V$ is compact, we can extract a finite subcover
$V_{i_1}, \ldots, V_{i_n}$ from the cover $(V_i)_{i \in I}$. Let $J = \{i_1, \ldots, i_n\}$. We thus
obtain an element $m = m_J \in \mathcal{M}(V)$ such that $\chi_{V_j} \cdot m = m_j$ for all $j \in J$.

It remains to verify that $\chi_{V_i} \cdot m = m_i$ for all $i \in I$, including those not in the finite
subcover $J$. Let $i \in I$. For any $j \in J$, the compatibility hypothesis gives
$\chi_{V_i \cap V_j} \cdot m_i = \chi_{V_i \cap V_j} \cdot m_j$. On the other hand, since
$\chi_{V_j} \cdot m = m_j$, we have
$\chi_{V_i \cap V_j} \cdot m = \chi_{V_i} \cdot (\chi_{V_j} \cdot m) = \chi_{V_i} \cdot m_j = \chi_{V_i \cap V_j} \cdot m_j$.
Therefore, $\chi_{V_i \cap V_j} \cdot (\chi_{V_i} \cdot m - m_i) = 0$ for all $j \in J$.

Since $(V_i \cap V_j)_{j \in J}$ is a finite cover of $V_i$, the idempotent $\chi_{V_i}$ can be expressed
as a polynomial in the idempotents $\chi_{V_i \cap V_j}$ without constant term. Because each term in this
polynomial expansion contains at least one factor of $\chi_{V_i \cap V_j}$, it follows that
$\chi_{V_i} \cdot (\chi_{V_i} \cdot m - m_i) = 0$. Since $m_i \in \mathcal{M}(V_i)$, we have
$\chi_{V_i} \cdot m_i = m_i$, which simplifies the equation to $\chi_{V_i} \cdot m = m_i$. This proves
the existence of the desired global section $m$.

Finally, the uniqueness of $m$ on $V$ follows from the same polynomial identity: if $m'$ is another
section such that $\chi_{V_i} \cdot m' = m_i$ for all $i \in I$, then $\chi_{V_j} \cdot (m - m') = 0$ for
all $j \in J$. The polynomial expansion of $\chi_V$ in terms of $\chi_{V_j}$ then implies
$\chi_V \cdot (m - m') = 0$, hence $m = m'$. This completes the proof that $\mathcal{M}$ is a sheaf.

It is clear that the two functors are inverse to each other: if we start from a non-degenerate
$\mathscr{D}(X)$-module $M$, we form the corresponding sheaf $\mathcal{M}$ as above. We then have
\[
\mathcal{M}_c(X)= \bigcup_{U\in \mathcal{T}_c} \mathcal{M}(U)= \bigcup_{U\in \mathcal{T}_c} \chi_U\cdot M= M,
\]
because $M$ is non-degenerate. Conversely, if we start from a sheaf of $\mathscr{C}^\infty_X$-modules
$\mathcal{M}$, and we set $M=\mathcal{M}_c(X)$, we see that for any compact open set $U$,
$\chi_U\cdot M= \chi_U \cdot \mathcal{M}(U)=\mathcal{M}(U) $. This shows that the sheaf constructed as
above from $M$ is indeed the initial sheaf $\mathcal{M}$.
\end{proof}
\medskip

Let $M$ be a non-degenerate $\scrD(X)$-module and let $\caM$ be the corresponding sheaf of
$\scrC^\infty_X$-modules given by the theorem. The stalk $\caM_x$ at any point $x \in X$ is the inductive
limit of the $\caM(U)=\chi_U \cdot M $ over the system of compact open sets $U$ of $X$ containing $x$.
More explicitly, let us equip $M$ with the equivalence relation defined by $ m \equiv m'$ if there exists
a compact open neighborhood $U$ of $ x $ such that $ \chi_{U}\cdot m=\chi_{U}\cdot m'$. The stalk
$\caM_x$ is then the set of equivalence classes.

We can also directly define the stalk $\caM_x$ as $M/M(x)$ where
\[
M(x)=\{ m\in M \mid \exists f \in \scrD(X), \, f(x)\neq 0 \text { and } f\cdot m=0  \}.
\]
Indeed, if $m \in M$, its image in $M/M(x)$ depends only on its equivalence class for the above relation:
if for some $m' \in M$ and some compact open neighborhood $U$ of $x$,
$ \chi_{U}\cdot m=\chi_{U}\cdot m'$, we have $\chi_U\cdot (m-m')=0$, which shows that $m-m'$ is in $M(x)$
(take $f=\chi_U$). This allows us to define a map from $\caM_x$ to $M/M(x)$. We then use $(ii)$ of Lemma
\ref{recouv} to show that this map is injective: if $m\in M(x)$ because $f \cdot m=0$ for some function
$f$ satisfying the required properties, we write it in the form $f=\sum_{j=0}^{r} \alpha_j \, \chi_{U_j}$
where the $U_j$ are disjoint compact open sets of $X$, $U_0$ contains $x$ and $\alpha_0$ is non-zero.
Since for all $j\neq 0$, $\chi_{U_0}\chi_{U_j}\cdot m=0$, we see that $\chi_{U_j}\cdot m \equiv 0$. On
the other hand
\[
f\cdot m = \sum_{j=0}^{r} \alpha_j \, \chi_{U_j} \cdot m=0
\]
and therefore $\alpha_0 \, \chi_{U_0}\cdot m \equiv 0$, hence $ m \equiv 0$, which shows the injectivity
of the map. The surjectivity being clear, we indeed have $\caM_x \simeq M/M(x)$.

\begin{exemple}
If $M=\scrD(X)$ then $\caM$ is the sheaf $\scrC^\infty_X$ of locally constant functions on $X$.
\end{exemple}

\begin{exemple}
If $q \colon X \rightarrow Y$ is a continuous map between two t.d. spaces, any non-degenerate
$\scrD(X)$-module $M$ becomes a non-degenerate $\scrD(Y)$-module via
\[
g\cdot m= (g \circ q) \chi_U\cdot   m,\quad (g \in \scrD(Y)), \, (m\in M),
\]
where $\chi_U$ is an idempotent of $\scrD(X)$ which fixes $m$. One immediately verifies that this
definition is independent of the choice of $\chi_U$ and indeed gives a non-degenerate $\scrD(Y)$-module
structure.

In fact $q$ defines an algebra morphism, still denoted $q$:
\[
q \colon \scrD(Y) \rightarrow \scrC^\infty(X), \quad g \mapsto g \circ q.
\]
This equips $\scrD(X)$ with a non-degenerate $\scrD(Y)$-bimodule structure as follows. Note first that
the algebras considered are commutative, and there is no need to distinguish between right and left
modules. Next, if $f \in \scrD(X)$, let $U$ be a compact open set of $X$ such that $\chi_U \cdot f=f$ and
let $U_Y$ be a compact open set of $Y$ containing the compact set $q(U)$. We then have
$(\chi_{U_Y}\circ q)\chi_U=\chi_U(\chi_{U_Y}\circ q)=\chi_U$, hence
\[
\chi_{U_Y} \cdot f= (\chi_{U_Y}\circ q)\cdot f=(\chi_{U_Y}\circ q)\chi_U\cdot  f= \chi_U \cdot f=f.
\]
We can therefore apply the results of the last paragraph of \ref{Oublietadjoints}. By the equivalence of
categories established in the theorem, the forgetful and pseudo-forgetful functors coincide respectively
with functors usually denoted $q_{!}$ and $q_*$ in sheaf theory (see \cite{Faisc}).

Recall that the forgetful functor admits a right adjoint
\[
I_Y^X \colon  \caM(\scrD(Y)) \rightarrow \caM(\scrD(X)),\quad M \mapsto \Hom_{\scrD(Y)}(\scrD(X), M)_{\scrD(X)},
\]
and the pseudo-forgetful functor a left adjoint
\[
P_Y^X \colon  \caM(\scrD(Y)) \rightarrow \caM(\scrD(X)),\quad  M \mapsto \scrD(X) \otimes_{\scrD(Y)} M.
\]
Again by the equivalence of categories of the theorem, we obtain functors
\[
I_Y^X, P_Y^X \colon   \scrC^\infty_Y-\caM od  (Y) \rightarrow  \scrC^\infty_X-\caM od  (X).
\]
In sheaf theory, the traditional notations for $I_Y^X$ and $ P_Y^X$ are respectively $q^{!}$ and
$q^{-1}$.
\end{exemple}

\medskip

If $\caF$ is a sheaf on the t.d. space $X$ and $Y$ is a locally closed subset of $X$, then we define a
sheaf $\caF_{|Y}$ on $Y$, which we call the restriction of $\caF$ to $Y$. The stalk of $\caF_{|Y}$ at a
point $y\in Y$ is equal to the stalk $\caF_y$, and the sections over an open set $U$ of $Y$ are those
which coincide in a neighborhood of each point of $U$ with the restriction of a section of $\caF$. When
we identify $\caF$ with the sheaf of sections of $\widehat{\caF}$, $\caF_{|Y}$ is identified with the
sheaf of sections over the open sets of $Y$.

\begin{defi}
Let $\caF \in \scrC^\infty_X-\caM od (X)$. The space \index[not]{F_c(X)*@$\caF_c(X)^*$}$\caF_c(X)^*$
(this is the algebraic dual of the vector space $\caF_c(X)$) is called the space of distributions on
$\caF$. It is clear that $\caF_c(X)^*$ is a $\scrC^\infty(X)$-module, the action being defined by
\[
\bil{f\cdot T}{\phi}= \bil{T}{f\phi}, \quad (T \in \caF_c(X)^*), \, (\phi \in \caF_c(X)), \, ( f \in \scrC^\infty(X)).
\]
\end{defi}

Let $U$ be an open set of $X$, $Z$ a closed set of $X$ and $\caF \in \scrC^\infty_X-\caM od (X)$. As in
(\ref{uxz}), we define two morphisms:
\[
i_U  \colon \caF_c(U) \rightarrow  \caF_c(X)  \quad    \text{ and }  \quad p_Z\colon \caF_c(X) \rightarrow  (\caF_{|Z})_c(Z).
\]
For simplicity, we write $\caF_c(Z)=(\caF_{|Z})_c(Z)$. We then have a result which generalizes
(\ref{uxz}) and (\ref{zxu}), and whose proof is similar:
\begin{prop}
With the above notation, and $Z=X \setminus U$, we have an exact sequence
\[
0\longrightarrow \caF_c(U) \stackrel{i_U}{\longrightarrow}  \caF_c(X) \stackrel {p_Z}  {\longrightarrow}   \caF_c(Z)   \longrightarrow 0.
\]
And dually, an exact sequence:
\[
0\longrightarrow \caF_c(Z)^* \stackrel{p_Z^*}  {\longrightarrow}  \caF_c(X)^* \stackrel{i_U^*} {\longrightarrow}  \caF_c(U)^* \longrightarrow 0.
\]
\end{prop}

\begin{rmq}
The above proposition asserts that any compactly supported section $s\in \caF_c(U)$ is the restriction of
a section in $\caF_c(X)$. We can then describe the stalk at a point $x\in X$ in the following way: if
$\phi \in \caF_c(X)$, let $\phi_{|U}$ denote its restriction to an open set $U$.
\begin{align*}
 \caF_x&=\caF_c(X)/ \{\phi \in \caF_c(X) \mid \exists  U \text{ open neighborhood of } x,\,  \phi_{|U} =0 \}\\
&=\caF_c(X)/ \langle  \, f\cdot \phi, \, \phi   \in \caF_c(X),\, f\in \scrD(X), f(x)=0 \rangle.
\end{align*}
To obtain this last equality, we note that
\[
\langle  \, f\cdot \phi, \, \phi   \in \caF_c(X),\, f\in \scrD(X), f(x)=0 \rangle
\]
is included in
\[
\{\phi \in \caF_c(X) \mid \exists  U \text{ open neighborhood of } x,\;   \phi_{|U } =0 \}.
\]
Let $\phi$ be in this latter set and $Z$ its support. Let $U$ be a compact open set of $X$ containing
$Z$, but not $x$ (such an open set is easy to construct using Lemma \ref{recouv}). Let $f$ be the
characteristic function of $U$. We then have $f\cdot \phi=\phi$ and $f(x)=0$, which shows the inclusion
in the other direction.
\end{rmq}

We end this section by returning to the notion of isomorphism of sheaves. Let
$\gamma \colon X\rightarrow Y$ be a homeomorphism between t.d. topological spaces. This induces an
algebra isomorphism $\scrD(Y)\simeq \scrD(X)$ by which we identify $\scrD(X)$ and $\scrD(Y)$-modules.

Let $X$ and $Y$ be two t.d. topological spaces and let $\caF$, $\caE$ be two sheaves in
$\scrC^\infty_X-\caM od (X)$ and $\scrC^\infty_Y-\caM od (Y)$ respectively. An isomorphism between $\caF$
and $\caE$ is therefore the data of a homeomorphism $\gamma \colon X\rightarrow Y$ and an isomorphism of
$\scrD(X)$-modules between $\caF_c(X)$ and $\caE_c(Y)$. It is clear that this induces an isomorphism of
$\bbC$-vector spaces $\caF_x \rightarrow \caE_{\gamma(x)}$ for all $x \in X$.

If $\gamma  \colon  \caF \rightarrow \caE$ is such an isomorphism, we still denote by
\[
\gamma \colon  \caE_c(Y)^* \rightarrow \caF_c(X)^*
\]
its transpose.

\section{T.d. topological groups}

The identity element of a group $G$ will be denoted by \index[not]{000@$\mathbf{1_G}$} $\mathbf{1}_G$.

\subsection{Groups and group actions}
If $X$ is a set, $G$ a group acting on $X$, and $E$ a vector space of functions on $X$, then $G$ acts
linearly on $E$ by:
\[
g\cdot f(x)=f(g^{-1}\cdot x), \quad (x\in X), \,  (g\in G),\, (f\in E).
\]

If $E'$ is a functional space on $X$ in duality with a space of functions $E$, then $G$ acts on $E'$ by:
\[
\langle g\cdot T, \, f \rangle =  \langle  T,\,  g^{-1}\cdot f \rangle, \quad  (g\in G),\,
(f\in E), \, (T\in E').
\]

In particular, if $G$ is a group, its automorphism group $\mathrm{Aut} \, G$ is a group acting on $G$. We
can then apply the preceding conventions. If $g_0$ is an element of $G$, we denote by $\Int(g_0)$ the
inner automorphism of $G$, $g\mapsto g_0gg_0^{-1}$.

The group $G$ acts on itself by left and right translation. Let us denote by $l$ and $r$ respectively
these actions:
\[
l(g)\cdot g'=gg',\quad r(g)\cdot g'=g'g^{-1},\quad  (g,\ g' \in G).
\]
These actions induce actions of $G$ on all functional spaces on $G$, which we will still denote by $l$
and $r$.

We denote by $f\mapsto \check f$ the linear map from a functional space on $G$ to itself induced by the
anti-automorphism $g\mapsto g^{-1}$ of $G$.

\subsection{A finiteness result}\label{fini}
In this paragraph, $G$ is a t.d. group. We will often use, without further comment, the following result:
\begin{lemme}
$(i)$ Let $K\subset K_1$ be compact open subgroups of $G$. Then $[K_1:K]$ is finite.

$(ii)$ If $K$ is a compact open subgroup of $G$, then $G/K$ equipped with the quotient topology is
discrete. If moreover $G$ is $\sigma$-compact, then $G/K$ is countable.
\end{lemme}
\begin{proof}
$(i)$ The group $K_1$ is the disjoint union of right cosets modulo $K$, which forms a cover of $K_1$ by
disjoint open sets. By compactness, we can extract a finite subcover. But this subcover cannot be
strictly smaller. Therefore, the set of right cosets is finite.

$(ii)$ Since left multiplication by $g\in G$ is a homeomorphism, $gK$ is open in $G$; hence each point of
$G/K$ is open for the quotient topology, and $G/K$ is discrete. Suppose that $G$ is a countable union of
compact subsets $A_n$, $n\in \bbN$. From the cover of $A_n$ by the open sets $gK$, $g\in G$, we extract a
finite subcover, and thus $G$ is covered by a countable family $(g_i K)_i$, which shows that $G/K$ is at
most countable.
\end{proof}

\subsection{Actions of t.d. groups}\label{actiontd}

When considering an action of a t.d. group $G$ on a t.d. space $X$, it is always a continuous left
action.

\begin{prop}\label{actiontdprop}
Suppose that $G$ is a $\sigma$-compact t.d. group acting on a t.d. space $X$ with a finite number of
orbits. Then there exists an open orbit $X_0$, and for any point $x_0 \in X_0$, the map
$G\rightarrow X_0$, $g\mapsto g\cdot x_0$ is open. It follows that all orbits are locally closed in $X$.
\end{prop}

\begin{proof}
Let $K$ be a compact open subgroup of $G$. Since we assumed the latter to be $\sigma$-compact, we have a
countable system of representatives $(g_i)_{i\in I}$ for the left cosets of $K$ in $G$. Let us also
choose a system of representatives $x_0,\ldots, x_m$ for the orbits of $G$ in $X$. We then have
\[
X=\bigcup_{i,j}  (g_iK)\cdot x_j.
\]
Thus $X$ is a countable union of compact subsets. Since $X$ satisfies the Baire category theorem
(Proposition \ref{recouv}), one of the subsets $ (g_iK)\cdot x_j$ is a neighborhood of one of its points.
Let $(g_ik)\cdot x_j $ be such a point. Then the subset $$K\cdot x_j=(g_ik)^{-1}(g_iK)\cdot x_j$$ is a
neighborhood of $x_j= (g_ik)^{-1}(g_ik)\cdot x_j$. By renumbering the orbits, if necessary, we may assume
that $j=0$. Let $K'$ be a compact open subgroup of $G$ contained in $K$. Since $K'$ has finite index in
$K$, there exists a finite system of representatives $\{k_1,\ldots, k_{r}\}$ for the left cosets of $K'$
in $K$. We can then repeat the same reasoning as above: one of the subsets $k_iK'\cdot x_0$ is a
neighborhood of one of its points $k_ik'\cdot x_0$ and $K'\cdot x_0 = (k_ik')^{-1}k_iK'\cdot x_0$ is a
neighborhood of $x_0$. This shows that the map $g\mapsto g\cdot x_0$ from $G$ to $X$ is open. For the
last assertion, we apply the above to the action of $G$ on the closed $G$-invariant subset
$X\setminus X_0$, which has strictly fewer orbits; by induction on the number of orbits, each orbit is
open in a closed subset of $X$, hence locally closed; this completes the proof of the proposition.
\end{proof}

\medskip

Let $G$ be a t.d. group acting on a t.d. space $X$. Let $G\backslash X$ denote the set of orbits, and
$p: X \rightarrow G\backslash X$ the natural projection. We equip $G\backslash X$ with the quotient
topology ($U \in G\backslash X$ is open if and only if $p^{-1}(U)$ is open in $X$). It is clear that $p$
is then continuous, but $G\backslash X$ is a topological space which is not necessarily Hausdorff.

\begin{lemme}\label{actiontdlem}
The map $p$ is open. If $Z$ is a closed (resp. locally closed) $G$-invariant subset of $X$, then $p(Z)$
is closed (resp. locally closed) in $G\backslash X$.
\end{lemme}

\begin{proof}
Let $U$ be an open subset of $X$. We have $p^{-1}(p(U))=\bigcup_{g\in G}g\cdot U $, and this is an open
subset of $X$. Thus $p(U)$ is open in $G\backslash X$, which shows that the map $p$ is open. If $Z$ is a
closed $G$-invariant subset of $X$, its complement $U$ is an open $G$-invariant subset of $X$, and since
$p$ is surjective, the complement of $p(Z)$ in $G\backslash X$ is the image of $U$ under $p$: it is an
open subset of $G\backslash X$. If now $Z$ is locally closed, $Z=U\cap F$, where $U$ is open and $F$ is
closed in $X$. We then have $Z=U\cap \overline{ Z}$, and since $\overline{ Z}$ is $G$-invariant,
$p(Z)=p(U)\cap p(\overline{ Z})$ (indeed, $Z=(G\cdot U)\cap \overline{Z}$, since $Z$ and $\overline{Z}$
are $G$-invariant, and $p(G\cdot U)=p(U)$). This shows that $p(Z)$ is locally closed in $G\backslash X$.
\end{proof}

Suppose that the topology of $G\backslash X$ is Hausdorff. Then, since $p$ is open and the image of a
compact set under $p$ is compact, it is clear that the space $G\backslash X$ is a totally disconnected
space.

\begin{cor}\label{actiontdcor}
Let $G$ be a $\sigma$-compact t.d. group, acting transitively on a t.d. space $X$. Let $x_0$ be in $X$,
and let $H$ be its stabilizer in $G$. Then the natural map
\[
G/H \rightarrow X, \quad gH\mapsto g\cdot x_0
\]
is a homeomorphism.
\end{cor}

\subsection{Quotient spaces}\label{topquo}
We apply the above results to the case of a quotient group:
\begin{lemme}
Let $N$ be a closed subgroup of $G$. Then the space $G/N$, equipped with the quotient topology, is a t.d.
space. If $N$ is normal, the group $M=G/N$ is a t.d. group.
\end{lemme}
\begin{proof}
The only thing left to verify is that the space $G/N$ is Hausdorff. We need to show that the diagonal
\[
\Delta=\{(xN,xN)\in G/N \times G/N \}
\]
is closed in $G/N \times G/N$. Consider the action of $N\times N$ on $G\times G$ by right multiplication,
so that the set of orbits $(G\times G)/(N\times N)$ identifies canonically with $G/N \times G/N$. Let
$p_2$ be the natural projection of $G\times G$ onto $G/N \times G/N$. By Lemma \ref{actiontdlem}, it
suffices to verify that $p_2^{-1}(\Delta)$ is closed. But
\[
p_2^{-1}(\Delta)=\{ (g_1,g_2)\in G\times G \, \mid \,  g_1^{-1}g_2\in N    \}
\]
is closed because $N$ is a closed subgroup.
\end{proof}

\subsection{Functional spaces on $G$}\label{espfonc}
In addition to the functional spaces $\scrC^\infty(G)$, $\scrD(G)$, $\scrD'(G)$, $\scrE'(G)$ already
defined above, and on which $G$ acts by $l$ and $r$, we introduce for any compact open subgroup $K$ of
$G$ the following spaces:

\[
\index[not]{CGKl@$\scrC^\infty(G,K,l)$}\scrC^\infty(G,K,l)= \{   f\in \scrC^\infty(G) \mid
l(k)\cdot f=f, \,  k\in K \},
\]
\[
\index[not]{CGKr@$\scrC^\infty(G,K,r)$}\scrC^\infty(G,K,r)= \{   f\in \scrC^\infty(G) \mid
r(k)\cdot f=f, \, k\in K \},
\]
\[
\index[not]{CGK@$\scrC^\infty(G,K)$}\scrC^\infty(G,K)= \scrC^\infty(G,K,l) \cap \scrC^\infty(G,K,r)
\]

We define analogously:
\[
\scrD(G,K,l), \quad \scrD(G,K,r)\quad  \text{ and } \quad \scrD(G,K),
\]
\[
\scrD'(G,K,l), \quad \scrD'(G,K,r) \quad\text{ and }\quad \scrD'(G,K),
\]
\[
\scrE'(G,K,l), \quad \scrE'(G,K,r) \quad\text{ and }\quad \scrE'(G,K).
\]

We also set:

\[
\scrC^\infty_{unif,l}(G)=\bigcup_{K} \scrC^\infty(G,K,l),\quad
\scrC^\infty_{unif,r}(G)=\bigcup_{K} \scrC^\infty(G,K,r),
\]
\[
\text{ and } \quad \scrC^\infty_{unif}(G)=\bigcup_{K} \scrC^\infty(G,K)
\]
where $K$ runs through the set of compact open subgroups of $G$.

\medskip

Similarly we define:
\[
\scrD_{unif,l}(G), \quad \scrD_{unif,r}(G) \text{ and } \scrD_{unif}(G),
\]
\[
\scrD'_{unif,l}(G), \quad \scrD'_{unif,r}(G) \text{ and } \scrD'_{unif}(G),
\]
\[
\scrE'_{unif,l}(G),  \quad\scrE'_{unif,r}(G) \text{ and } \scrE'_{unif}(G).
\]

\enlargethispage{2\baselineskip}
\begin{lemme} \label{KgK}
$(i)$ Let $f$ be a function in $\scrD(G)$. Then there exists a compact open subgroup $K$ of $G$, a finite
set of representatives $\{g_j \}$ for the double cosets modulo $K$ in $G$, and scalars $a_j$ such that
\[
f=\sum_j a_j \, \chi_{ K g_j K}
\]
where $\chi_{X}$ denotes the characteristic function of a subset $X$ of $G$.

$(ii)$ Let $f \in \scrD(G,K,l)$ (resp. $\scrD(G,K,r)$). Then the support of $f$ is a finite union of
right (resp. left) cosets of $K$ in $G$. Let $g_1,\ldots,g_r$ be representatives of these cosets. Then we
can write $f$ as:
\[
f=\sum_i c_i \; \chi_{Kg_i},\qquad \text{resp. }\quad   f=\sum_i c_i\;  \chi_{g_iK},
\]
where only finitely many $c_i$ are non-zero.
\end{lemme}

\begin{proof}
Let $F=\mathrm{Supp}(f)$, a compact subset of $G$. Since $f$ is locally constant, for any $g\in F$, there
exists a compact open subgroup $K_g$ such that $f$ is constant on $K_g g K_g$. From the cover
$F \subset \bigcup_{g\in F} K_g g K_g$, we can extract a finite subcover (by compactness of $F$), say
$F \subset \bigcup_{i} K_{g_i}g_i K_{g_i}$. Let $K= \bigcap_i K_{g_i}$: this is a compact open subgroup
of $G$. Since for all $i$, $[K_{g_i}:K ]$ is finite (Lemma \ref{fini}), we can find a finite set
$\{g_j\}$ of representatives for the double cosets modulo $K$ such that $F \subset \bigcup_{j} K g_j K$,
and therefore we can write $f$ as:
\[
f=\sum_j a_j \, \chi_{ K g_j K},
\]
where the $a_j$ are scalars. The proof of $(ii)$ is similar.
\end{proof}

\begin{cor}
We have:
\[
\scrD(G) \subset \scrC^\infty_{unif}(G) \subset \scrC^\infty(G).
\]
In particular $\scrD_{unif,l}(G)=\scrD_{unif,r}(G)=\scrD_{unif}(G)=\scrD(G)$.
\end{cor}

\subsection{Haar measure}

\index[ter]{Haar measure}

We show here the existence of a (left) Haar measure on a t.d. group. The result remains valid for any
locally compact topological group, but the proof for t.d. groups given here is simpler than that of the
general case.

\begin{prop}
Let $G$ be a t.d. group. There exists, up to a scalar multiplicative factor, a unique distribution
\index[not]{ZZmu@$\mu_G$}$\mu_G$ on $G$ invariant under left translation. We can choose $\mu_G$ to be
positive (i.e., such that $\langle \mu_G, f \rangle > 0$ for any non-identically zero function
$f\in \scrD(G)$ with non-negative values). A right-invariant measure $\nu_G$ can be defined similarly.
\end{prop}
\begin{proof}
Let $K_\alpha$ be a fundamental system of neighborhoods of the identity element in $G$ consisting of
compact open subgroups, and suppose moreover that one of them, say $K_{\alpha_0}$, contains all the
others. For all $\alpha$, let us set $\scrD_\alpha=\scrD(G,K_\alpha,r)$. It is clear that if
$K_\alpha \subset K_\beta$, then $\scrD_\beta \subset \scrD_\alpha $ and that
$\scrD(G)=\bigcup_{\alpha} \scrD_\alpha$. On the other hand, each $\scrD_\alpha$ is stable under the
action $l$ of $G$. To define a left-invariant distribution on $G$, it suffices to define in a compatible
way for each $\alpha$ a left-invariant element $\mu_\alpha$ of the dual of $\scrD_\alpha$, i.e., if
$\scrD_\beta \subset \scrD_\alpha $, then ${\mu_\alpha}_{|\scrD_\beta}=\mu_\beta$. Since $\scrD_\alpha$
is generated by the left translates of the characteristic function of $K_\alpha$, this element of the
dual is completely determined by its value on $\chi_{K_ \alpha}$. It follows that $\mu_G$ is indeed
unique up to multiplication by a scalar factor. To prove existence, we must normalize the $\mu_\alpha$ in
a compatible way, and for this, we use the compact open subgroup $K_{\alpha_0}$ above. We choose
$\mu_{\alpha_0}$ by setting $\langle \mu_{\alpha_0},\chi_{K_{\alpha_0}} \rangle =1$, and we then
normalize the $\mu_\alpha$ by setting
\[
\langle \mu_{\alpha},\chi_{K_{\alpha}} \rangle = [ K_{\alpha_0}:K_\alpha]^{-1}.
\]
It is easy to verify that this indeed defines a compatible family $\{\mu_\alpha\}$.
\end{proof}

\medskip

The Haar measure $\mu_G$ is an element of $\scrD'(G)$. For any function $f \in \scrC^\infty(G)$, let us
define the distribution
\[
f \mu_G:  \phi \in \scrD(G) \mapsto \int_G  \phi(g)  f(g) \,  d\mu_G(g).
\]
This allows us to embed (canonically, up to a scalar multiple) $\scrC^\infty(G)$ into $\scrD'(G)$. If $f$
is moreover compactly supported, the same is true for the distribution $f \mu_G$, yielding an embedding
of $\scrD(G)$ into $\scrE'(G)$.

\subsection{Modular function}\label{foncmod}
\index[ter]{modular!(function or character)} Let $G$ be a t.d. group and let $\mu_G$ be a positive left
Haar measure on $G$. The uniqueness up to a multiplicative factor of the Haar measure allows us to define
the modular function $\delta_G$, using the formula
\[
r(g)\cdot \mu_G = \index[not]{ZZdeltaG@$\delta_G$}\delta_G(g)\,  \mu_G.
\]
Indeed, $r(g)\cdot \mu_G$ is a positive left Haar measure and in particular $\delta_G$ takes values in
$\bbR^\times_+$. It is easy to verify that $\delta_G$ is a character of $G$. In practice, the formula to
use is therefore the following: for any function $f \in \scrD(G)$, for all $h,h'\in G$,
\[
\int_G f(h'gh^{-1})\; d\mu_G(g)= \int_G f(gh^{-1})\; d\mu_G(g)= \delta_G(h) \;  \int_G f(g)\;
d\mu_G(g).
\]

\begin{prop}
$(i)$ The restriction of $\delta_G$ to any compact subgroup of $G$ is identically equal to $1$.

$(ii)$ The distribution $\delta_G^{-1} \mu_G$ is right-invariant. If $T \mapsto \check T$ is the
involution of $\scrD'(G)$ induced by the homeomorphism $g\mapsto g^{-1}$, then
$\check \mu_G=\delta_G^{-1} \mu_G$.
\end{prop}
\begin{proof}
The image of a compact subgroup $K$ of $G$ under $\delta_G$ is a compact subgroup of $\bbR^\times_+$:
there is only the trivial subgroup. This shows $(i)$. For $(ii)$, we first verify that we have
$r(g)\cdot \delta_G=\delta_G(g) \delta_G$ for all $g\in G$. We then have:
\[
r(g)\cdot(\delta_G^{-1} \mu_G)= (r(g)\cdot\delta_G^{-1}) \, (r(g)\cdot\mu_G)=
\delta_G(g)^{-1}\delta_G^{-1}\delta_G(g)\mu_G=\delta_G^{-1}\mu_G
\]
and therefore $\delta_G^{-1}\mu_G$ is right-invariant. Since it is clear that $\check \mu_G$ is as well,
these two distributions are proportional, and by evaluating them on the characteristic function of a
compact open subgroup, we see that they are equal.
\end{proof}

\begin{defi}
A group $G$ such that $\delta_G$ is identically equal to $1$ is said to be unimodular.
\index[ter]{unimodular}
\end{defi}
\begin{exemples}
$ $
\begin{itemize}
\item[1.] An abelian group is unimodular.
\item[2.] A compact group is unimodular (because the only compact subgroup of $\bbR^\times_+$ is
$\{1\}$).
\item[3.] If $G$ is the union of its compact subgroups, then $G$ is unimodular.
\end{itemize}
\end{exemples}

\medskip

We can more generally define the modulus of an automorphism $\sigma$ of $G$ by the formula
$\sigma\cdot \mu_G=\delta_G(\sigma)\; \mu_G $. More explicitly, since for any function $f \in \scrD(G)$,
\[
\bil{\sigma\cdot \mu_G}{f}= \bil{ \mu_G}{ \sigma^{-1}\cdot f}=
\int_G  (\sigma^{-1}\cdot f)(g) \; d\mu_G(g)= \int_G  f(\sigma (g)) \; d\mu_G(g),
\]
we obtain the formula
\begin{equation}\label{fmod}
 \int_G  f(\sigma(g)) \; d\mu_G(g)  = \delta_G(\sigma)   \int_G f(g) \; d\mu_G(g).
\end{equation}
In particular, if $\sigma=\Int (h)$, we obtain $\delta_G(\sigma)=\delta_G(h)$.

\subsection{Calculation of the modular function}
We will give a formula for $\delta_G$ in terms of indices of compact open subgroups. As in the previous
section, $G$ is a t.d. group and $\sigma$ is an automorphism of $G$. We fix a left Haar measure $\mu_G$
on $G$. Let $K$ be a compact open subgroup of $G$ and let us set $K_1=K\cap \sigma^{-1}(K)$. Since
\[
\mu_G(K)=\int_G \chi_K(g)\, d\mu_G(g),
\]
equation (\ref{fmod}) gives us:
\begin{align*}
\delta_G(\sigma)\mu_G(K)&=\delta_G(\sigma) \int_G \chi_K(g)\, d\mu_G(g)
= \int_G  \chi_K(\sigma(g))\, d\mu_G(g) = \int_G \chi_{\sigma^{-1}(K)}(g)\, d\mu_G(g)\\
  &=\mu_G(\sigma^{-1}(K)).
\end{align*}
Now,
\[
\mu_G(K)= [K:K_1 ]  \mu_G(K_1), \quad   \mu_G(\sigma^{-1}(K))= [\sigma^{-1}(K):K_1 ] \mu_G(K_1),
\]
hence,
\begin{equation}\label{fmodK}
\delta_G(\sigma)= \frac{[\sigma^{-1}(K):K_1 ]}{[K:K_1 ]}.
\end{equation}

\subsection{Invariant measure on a homogeneous space}\label{invmesquo}
\index[ter]{invariant measure}

Let $H$ be a closed subgroup of a t.d. group $G$. There does not always exist a measure on the t.d.
topological space $G/H$ (cf. Lemma \ref{topquo}) invariant under the action of $G$ by left translations.
Let us describe what happens in general, but for reasons that will appear later, let us consider instead
the space $H\backslash G$ equipped with the action of $G$ by right translations. Let us set
\index[not]{ZZdeltaHG@$\delta_{H\backslash G}$}$\delta_{H\backslash  G}=\delta_G/\delta_H$: this is a
character of $H$. Consider the set of locally constant functions $f$ on $G$ satisfying:

$a)$ $f(hg)= \delta_{H\backslash G} (h)\, f(g)$, $g\in G$, $h\in H$.

$b)$ $f$ is compactly supported modulo $H$, i.e., there exists a compact subset $K_f$ of $G$ such that
$\supp (f)\subset H\cdot K_f$.

We denote by $\scrD(G,H, \delta_{H\backslash G})$ the space of these functions. It is clear that this
space is invariant under the right action of $G$. If $\delta_{H\backslash G}$ is identically equal to the
constant function $1$, then $\scrD(G,H, \delta_{H\backslash G})$ is isomorphic to
$\scrD(H \backslash G)$. \index[not]{DGHdelta@$\scrD(G,H, \delta_{H\backslash G})$}

\begin{thm}
There exists, up to a multiplicative factor, a unique distribution $\nu_{H \backslash G}$ on \allowbreak
$\scrD(G,H, \delta_{H\backslash G})$ invariant under the right action of $G$. We can choose it to be
positive. \index[not]{ZZnuHG@$\nu_{H \backslash G}$}
\end{thm}

We denote:
\[
\langle\nu_{H \backslash G},f \rangle = \int_{H \backslash G}f(g)\; d\nu_{H \backslash G}(g),
\quad f\in \scrD(G,H, \delta_{H\backslash G}).
\]

\begin{proof}
Let $\mu_G$ and $\mu_H$ be left Haar measures on $G$ and $H$ respectively. Let us define an operator
$p \colon \scrD(G) \rightarrow \scrC^\infty(G)$ by
\[
p(f)(g)=\int_H  f(hg) \delta_G(h)^{-1}\,  d \mu_H(h).
\]
It is clear that $p$ commutes with the right action of $G$ on $\scrD(G)$ and $\scrC^\infty(G)$. On the
other hand, for all $h_0 \in H$,
\begin{align*}
   p(f)(h_0g)&=\int_H  f(hh_0g) \delta_G(h)^{-1} \, d\mu_H(h) = \delta_G(h_0)\int_H  f(hh_0g) \delta_G(hh_0)^{-1}   \, d\mu_H(h)\\
&=\delta_G(h_0) \delta_H(h_0)^{-1}    \int_H  f(hg) \delta_G(h)^{-1} \, d\mu_H(h)
=\delta_{H\backslash G}(h_0) \int_H  f(hg) \delta_G(h)^{-1}   \, d\mu_H(h)\\
&=\delta_{H\backslash G}(h_0) \,  p(f)(g)
\end{align*}
and therefore the image of $p$ is in $\scrD(G,H, \delta_{H\backslash G})$, the condition on the support
being obviously satisfied.

Let us now show that $p$ is surjective. For any $g$ in $G$ and any compact open subgroup $K$ of $G$, let
$\scrD(G)_g^K$ (respectively $\scrD(G,H, \delta_{H\backslash G})_g^K$) denote the subspace of $\scrD(G)$
(resp. $\scrD(G,H, \delta_{H\backslash G})$) of functions supported in $HgK$ and fixed under the action
by right translation of $K$. It is clear that
\[
p(\scrD(G)_g^K) \subset \scrD(G,H, \delta_{H\backslash G})_g^K.
\]
On the other hand, any function in $\scrD(G,H, \delta_{H\backslash G})_g^K$ is determined by its value at
$g$, and therefore this space is $1$-dimensional.

Let $f \in \scrD(G,H, \delta_{H\backslash G})$ and $K_f$ a compact subset of $G$ such that
$\supp (f) \subset H\cdot K_f$. For each $x\in K_f$, let $K_x$ be a compact open subgroup of $G$ such
that $f$ is constant on $xK_x$. From the cover $K_f \subset \bigcup_{x\in K_f} xK_x$ we extract by
compactness a finite cover $K_f \subset \bigcup_{i} x_iK_{x_i}$, and we set $K=\bigcap_i K_{x_i}$. Since
each $[K_{x_i}: K]$ is finite, by considering representatives of the left cosets of $K$ in the $K_{x_i}$,
we find a finite cover of $K_f$ of the form $K_f \subset \bigcup_{j} x_jK$. From this, we deduce that $f$
can be written as a finite sum $f=\sum_j f_j$ where $\supp (f_j)=Hx_jK$, the $f_j$ being right-invariant
under $K$ (and $f_j$ is determined by its value at $x_j$). Any function in
$\scrD(G,H,\delta_{H\backslash G})$ can be written as a finite sum of functions in some
$\scrD(G,H,\delta_{H\backslash G})_g^K$. Consequently, it suffices to show the surjectivity of
$p \colon \scrD(G)_g^K \rightarrow \scrD(G,H, \delta_{H\backslash G})_g^K$. Since the target space is
one-dimensional, it suffices to show that $p (\scrD(G)_g^K)$ is non-zero. But it is clear that $p(f)$ is
non-zero as soon as $f$ takes non-negative values and is not identically zero. This completes the proof
of the surjectivity of $p$.

Let $f \in\scrD(G)_g^K $. Since $f$ is compactly supported, there exists a compact subset $F$ of $H$ such
that $f$ is supported in $FgK$. For each $y \in F$, $ygK$ is a neighborhood of $yg$ in $G$, and therefore
there exists a compact open subgroup $K_y$ of $G$ such that $K_y yg \subset ygK$. From the cover
$F \subset \bigcup_{y\in F} K_y y$ we extract a finite subcover $F \subset \bigcup_{i} K_{y_i} y_i$, and
we have $FgK \subset \bigcup_{i} K_{y_i} y_igK \subset \bigcup_{i} y_igK$. Since $f$ is invariant under
the right action of $K$, we see that it is determined by its values at the points $y_i$. The space
$\scrD(G)_g^K$ is therefore identified with the space of finitely supported functions on the space
$HgK/K$, a space on which $H$ acts transitively by left translation. Since $HgK/K \subset G/K$ and the
latter is discrete (Lemma \ref{fini}, $(ii)$), the same is true for $HgK/K$. The space $\scrD(G)_g^K$ is
therefore generated by the $l(h)\cdot \chi_{gK}$, $h\in H$. Two linear forms $T_i$, $i=1,2$, on
$\scrD(G)_g^K$ such that $T_i(l(h)\cdot f)= \delta_G(h)^{-1} T_i(f)$ are therefore proportional. For any
$g\in G$, let $\delta_g$ denote the Dirac distribution at $g$. Since
$\scrD(G,H, \delta_{H\backslash G})_g^K$ is $1$-dimensional, its dual is as well, and it is generated by
$\delta_g$. It is easy to see that the distributions $T_1$ and $T_2$ defined by
$T_1(f):= \langle \delta_g,p(f)\rangle$ and $T_2(f):= \langle \nu_G,f \rangle$ both satisfy the above
property, and are therefore proportional. We deduce that if $f\in \scrD(G)$, $p(f)=0$ as soon as
$\langle \nu_G,f \rangle=0$. We have thus shown that
$\scrD(G,H, \delta_{H\backslash G}) \simeq \scrD(G)/\ker p$ and that any distribution on $G$ invariant
under the right action of $G$ vanishes on $\ker p$. The right Haar measure $\nu_G$ on $G$ therefore
induces an invariant linear form $\nu_{H\backslash G}$ on $\scrD(G,H, \delta_{H\backslash G})$, unique up
to a multiplicative factor. Since we can take $\nu_G$ to be positive, the same is true for
$\nu_{H\backslash G}$.
\end{proof}

\begin{rmq}
The space $\scrD(G,H,\delta_{H\backslash G})$ is a $\scrD(H\backslash G)$-module (by pointwise
multiplication of functions) and therefore defines a sheaf on $H\backslash G$ by the equivalence of
categories of Theorem \ref{faisceqcat}.
\end{rmq}

\subsection{Convolution}\label{lract}

Let $G$ be a t.d. group and let $T_1$ and $T_2$ be two distributions in $\scrE'(G)$. We can then define
their convolution product by the formula:
\[
\int_G f(g)\, d( T_1 *T_2)(g) := \int_{G \times G} f(gh)\, d(T_1 \otimes T_2)(g,h), \quad
(f \in \scrD(G)).
\]
This is well-defined, because $T_1 \otimes T_2$ is compactly supported according to Proposition
\ref{tens}. If $T \in \scrD'(G)$ and $f \in \scrD(G)$, or if $T \in \scrE'(G)$ and
$f \in \scrC^\infty(G)$, we also define $T*f$ and $f*T$ in $\scrC^\infty(G)$ by:
\[
(T*f)(g_0):= \int_G f(g^{-1}g_0)\, dT(g)=\bil{T}{l(g_0)\cdot \check f}
\]
\[
(f*T)(g_0):= \int_G  f(g_0g^{-1}) \, dT(g) = \bil{T}{r(g_0^{-1})\cdot \check f}.
\]
If $T$ and $f$ are compactly supported, then $T*f$ and $f*T$ are in $\scrD(G)$.

\begin{prop}
For all distributions $T_1,T_2,T_3$ in $\scrE'(G)$ and any function $f$ in $\scrD(G)$, we have

$(i)$ $\mathrm{Supp}(T_1*T_2) \subset \mathrm{Supp}(T_1)\cdot \mathrm{Supp}(T_2)$. In particular
$T_1*T_2$ is compactly supported.

$(ii)$ $(T_1*T_2)*T_3=T_1*(T_2*T_3)$.

$(iii)$ $ (T_1*T_2)*f=T_1*(T_2*f), \qquad (T_1*(f*T_2))=(T_1*f)*T_2,$
$ \quad   (f*T_1)*T_2=f*(T_1*T_2) .$

Let $\mu_G$ and $\nu_G$ be respectively left and right Haar measures. For any $T \in \scrE'(G)$ and any
$f \in \scrD(G)$, we have,

$(iv)$ $(T*f\mu_G)= (T*f)\mu_G, \quad (f\nu_G*T)=(f*T)\nu_G $.

For any $T \in \scrD'(G)$ and any $f \in \scrD(G)$, we have,

$(v)$ $\bil{T}{f} = (T*\check f)(\mathbf{1}_G)= (\check f *T) (\mathbf{1}_G)$.
\end{prop}

\begin{proof}
Point $(i)$ follows easily from Proposition \ref{tens}, and point $(ii)$ from Fubini's formulas. The
definitions of $T*f$ and $f*T$ are made so that $(iv)$ is satisfied. It is easy to deduce $(iii)$ from
$(ii)$ and $(iv)$, or directly from Fubini's formulas. The last point is immediate.
\end{proof}

\medskip

The space $\scrE'(G)$ is stable under the convolution product and this product is associative and
bilinear for the vector space structure of $\scrE'(G)$. Let $\delta_g$ denote the Dirac distribution at
the point $g\in G$. It is obviously an element of $\scrE'(G)$. It is clear using Fubini's formulas that
for any $T \in \scrE'(G)$, we have
\[
\delta_{\mathbf{1}_G}*T=T*\delta_{\mathbf{1}_G}=T.
\]
In other words, $\delta_{\mathbf{1}_G}$ is the identity element for the convolution product on
$\scrE'(G)$. To summarize, the space $\scrE'(G)$ of compactly supported distributions on $G$ equipped
with the convolution product is an associative algebra, with identity element $\delta_{\mathbf{1}_G}$.
\medskip

The group $G$ acts on the convolution algebra $\scrE'(G)$ by $l$ and $r$.

\begin{lemme}
Let $T \in \scrE'(G)$ and let $g,\, g'\in G$. We then have
\[
\delta_g *T =l(g)\cdot T,\quad T*\delta_g=r(g^{-1})\cdot T,
\quad \delta_g*\delta_{g'}=\delta_{gg'}.
\]
\end{lemme}
\begin{proof}
This is the result of straightforward calculations.
\end{proof}

\subsection{Some formulas}\label{Gamma}
If $\Gamma$ is a compact subgroup of $G$, then the Haar measure on $\Gamma$, normalized by
$\mu_\Gamma(\Gamma)=1$, defines an element $e_\Gamma$ of $\scrE'(G)$:
\[
\langle e_\Gamma,f \rangle=\int_\Gamma f(\gamma) \, d\mu_\Gamma(\gamma), \qquad
(f \in \scrC^\infty(G)).
\]

\begin{lemme}\label{Gammalem}
Let $\Gamma_1$ and $\Gamma_2$ be compact subgroups of $G$, and let $g\in G$. Then there exists a unique
compactly supported distribution $T$ on $G$ such that $\supp T \subset \Gamma_1 g \Gamma_2$, invariant
under the action of $\Gamma_1$ by $l$ and of $\Gamma_2$ by $r$, and normalized by $\int_G dT=1$.
\end{lemme}
\begin{proof}
We define an action $\rho$ of $\Gamma_1 \times \Gamma_2$ on $ \Gamma_1 g \Gamma_2$ by:
\[
\rho(\gamma_1,\gamma_2)\cdot x = \gamma_1x\gamma_2^{-1},\quad (x\in \Gamma_1 g \Gamma_2), \,
(\gamma_1\in \Gamma_1),\, (\gamma_2\in \Gamma_2).
\]
Since this action is transitive, we deduce from Theorem \ref{invmesquo} and the fact that
$\Gamma_1 \times \Gamma_2$ is unimodular, that there exists a unique
$(\Gamma_1 \times \Gamma_2)$-invariant distribution $T_0$ on $\Gamma_1 g \Gamma_2$ such that
$\int_{\Gamma_1 g \Gamma_2} dT_0=1$. According to Corollary \ref{FetDcor}, this defines a distribution on
$G$ having the desired properties. One easily verifies that this distribution is none other than
$e_{\Gamma_1}*\delta_g*e_{\Gamma_2}$.
\end{proof}

\begin{prop}\label{Gammaprop}
Let $\Gamma$, $\Gamma_1$ and $\Gamma_2$ be compact subgroups of $G$, such that
$\Gamma=\Gamma_1\cdot \Gamma_2$ and let $g$ be an element of $G$. We have:

$(i)$ $e_\Gamma=e_{\Gamma_1}  *e_{\Gamma_2} =e_{\Gamma_2}  *e_{\Gamma_1} $.

$(ii)$ $\delta_g * e_\Gamma * \delta_{g^{-1}}= e_{g\Gamma g^{-1}}$.

In particular $e_\Gamma * e_\Gamma=e_\Gamma$ and if $\Gamma_1 \subset \Gamma_2$,
$e_{\Gamma_2}=e_{\Gamma_1}  *e_{\Gamma_2} =e_{\Gamma_2}  *e_{\Gamma_1}$.
\end{prop}
\begin{proof}
This follows without difficulty from the previous lemma.
\end{proof}

\subsection{The Hecke algebra}    \label{HGK}
We will slightly change the notation introduced in paragraph \ref{espfonc}.

\begin{lemme}
Let us fix a left Haar measure $\mu_G$ on $G$. Then multiplication by $\mu_G$ induces an isomorphism
$$\scrD(G)\simeq  \scrE'_{unif,l}(G).$$ Moreover, any distribution $T$ in $\scrE'_{unif,l}(G)$ is also in
$\scrE'_{unif,r}(G)$, hence in $\scrE'_{unif}(G)$.
\end{lemme}
We therefore have $\scrE'_{unif,l}(G)=\scrE'_{unif,r}(G)=\scrE'_{unif}(G)$ and we will henceforth denote
by \index[not]{HG@$\caH(G)$}$\caH(G)$ this space, generally called the \index[ter]{Hecke algebra} Hecke
algebra of $G$.

\begin{proof}
Let $T$ be a distribution in $\scrE'_{unif,l}(G)$, and let $K$ be a compact open subgroup such that $T$
is invariant under the action of $K$ by $l$. According to Lemma \ref{Gammalem}, the restriction of $T$ to
any subset of the form $Kg$ of $G$ is proportional to the restriction of $\mu_G$. We deduce that $T$ is
of the form $T=f\, \mu_G$ for a function invariant under the action of $K$ by $l$. Such a function is
locally constant. Moreover, clearly, $\supp T= \supp f$ and therefore since $T$ is compactly supported,
$f$ is in $\scrD(G)$. Conversely, it is clear that if $f\in \scrD(G)$, then
$f\, \mu_G \in \scrE'_{unif,l}(G)$. We deduce the last assertion.
\end{proof}

\medskip

If $\Gamma$ is a compact open subgroup of $G$, then $e_\Gamma$ is in $\caH(G)$. On the other hand, for
any $g \in G$, $\delta_g$ is in $\scrE'(G)$ but not in $\caH(G)$.

\begin{rmqs}
$ $
\begin{itemize}
\item[1.] It is clear that $\caH(G)$ is stable under the convolution product. But since
$\delta_{\mathbf{1_G}}$ is not in $\caH(G)$, the algebra $(\caH(G),*)$ is not unital.

\item[2.] The subalgebra $\caH(G)$ is in fact a two-sided ideal of $(\scrE'(G),*)$. Indeed,
\[
T*f\mu_G=(T*f)\mu_G,\quad    f\nu_G*T=(f*T)\nu_G,\quad (T \in \scrE'(G)), \,  (f \in \scrD(G)).
\]

\item[3.] The isomorphism $\scrD(G) \simeq \caH(G)$ induces by transport of structure a convolution
product on $\scrD(G)$:
\[
(f_1\mu_G)*  (f_2\mu_G)=(f_1*f_2)\mu_G, \quad  (f_1,f_2\in \scrD(G)).
\]
The explicit formula for this product is:
\begin{equation}\label{fconvg}
(f_1*f_2)(g_0)=\int_G f_1(g)f_2(g^{-1}g_0)\, d\mu_G(g), \quad (f_1,f_2\in \scrD(G)).
\end{equation}
\end{itemize}
\end{rmqs}

\begin{prop}
The algebra $\caH(G)$ is an algebra with idempotents. The family $(e_K)_K$, where $K$ runs through the
set of compact open subgroups of $G$, is a directed system of idempotents. It is equipped with the
anti-involution $T \mapsto \check T$ induced by $g\mapsto g^{-1}$. The idempotents $e_K$ satisfy
$\check e_K=e_K$.
\end{prop}
\begin{proof}
We have seen that any compact open subgroup $K$ of $G$ defines an idempotent $e_K$ in $\caH(G)$. For any
distribution $T$ in $\caH(G)$, we can find a compact open subgroup $K$ of $G$ small enough such that $T$
is invariant under the action of $K$ by left and right translation. We immediately deduce that
$e_K* T*e_K=T$. If we have a finite family of distributions $\{T_i\}$ in $\caH(G)$, we find for each of
them a compact open subgroup $K_i$ of $G$ such that $e_{K_i}* T_i * e_{K_i}=T_i$. It then suffices to
take $K =\bigcap_i K_i$, which is indeed a compact open subgroup of $G$, and we have $e_K*T_i*e_{K}=T_i$.
It is clear that $T \mapsto \check T$ preserves the support and the invariance by left and right
translation and therefore preserves $\caH(G)$. It is obvious that $\check e_K=e_K$ for any compact open
subgroup $K$.
\end{proof}

\medskip

Note that $e_{K_1}\leq e_{K_2}$ for the order on $\mathrm{Idem}(\caH(G))$ (cf. paragraph \ref{idemgen})
if and only if $K_2 \subset K_1$.

Let $K$ be a compact open subgroup of $G$. Let us set
\index[not]{HGK@$\caH(G,K)$}$$\caH(G,K)=e_K*\caH(G)*e_K.$$ It is a unital subalgebra of $\caH(G)$, with
identity element $e_K$.

\medskip

\begin{rmqs}
$ $
\begin{itemize}
\item[1.] The space $\caH(G,K)$ is none other than the space $\scrE'(G,K)$ introduced in paragraph
\ref{espfonc}.

\item[2.] A basis for $\caH(G,K)$ is given by the distributions $$a_{g,K}:=e_K*\delta_g*e_K, $$ where $g$
runs through a system of representatives for the double cosets $K\backslash G/K$.
\end{itemize}
\end{rmqs}

\begin{proof}
Let $T=e_K*T*e_K$ be an element of $\caH(G,K)$, and let $k \in K$. We have:
\[
l(k)\cdot T=l(k)\cdot(e_K*T*e_K)=\delta_k*e_K*T*e_K=e_K*T*e_K
\]
according to Lemma \ref{lract} and Proposition \ref{Gammaprop}. Thus the action $l$ of $K$ on $\caH(G,K)$
is trivial. Similarly for the action $r$. This shows that $\caH(G,K)\subset  \scrE'(G,K)$. Conversely, if
$T \in \scrE'(G,K)$, we use Fubini to then show that $e_K*T=T$ and $T*e_K=T$, and therefore
$T \in \caH(G,K)$.

It is clear that $a_{g,K}=a_{g',K}$ for all $g'$ in the double coset $KgK$. The isomorphism
$\caH(G)\simeq \scrD(G)$ given by the choice of a Haar measure induces an isomorphism
$\caH(G,K)\simeq \scrD(G,K)$, and according to Lemma \ref{espfonc}, any function in $\scrD(G,K)$ is a
linear combination of the $\chi_{KgK}$. It is easy to see that $\chi_{KgK}$ corresponds up to a scalar
multiple to $a_{g,K}$ under the isomorphism $\caH(G,K)\simeq \scrD(G,K)$, and that the $\{\chi_{KgK}\}_g$
where $g$ runs through a system of representatives for the double cosets $K\backslash G/K$ is a basis for
$\scrD(G,K)$. Remark 2 follows from this.
\end{proof}

\subsection[Description of the completion and the center]{Completion of $\caH(G)$ and center of $\caM(\caH(G))$} \label{essencomp}
In this paragraph, we will explicitly describe in terms of distributions on $G$ the completion
\index[ter]{completion!of $\caH(G)$}$\overline{\caH(G)}$ of the algebra with idempotents $\caH(G)$, in
the sense of Section \ref{completeA}. This will also allow us to describe the center
\index[ter]{center!of $\caM(\caH(G))$} of the category $\caM(\caH(G))$ in terms of invariant
distributions on $G$.

\begin{defi}
For any distribution $T\in \scrD'(G)$, for any element $h$ of $\caH(G)$, and for any function $f$ in
$\scrD(G)$, let us set
\[
\bil{T \bullet h}{f}:=\bil{T} {f*\check h} \quad \text{ and }   \quad \bil{h \bullet T}{f}:=\bil{T}{\check h *f}.
\]
It is clear that $T\bullet h$ and $h \bullet T$ define elements of $\scrD'(G)$.
\end{defi}

\begin{rmq}
With the notation of the definition, if $T$ is compactly supported, $T*h$ is defined, compactly
supported, and using Proposition \ref{lract}, $(v)$,
\begin{align*}
 \bil{T * h}{f}&=((T*h)*\check f)(\mathbf{1}_G)=(T*(h*\check  f))(\mathbf{1}_G)\\
&= (T*(f* \check h )\, \check{}\, )(\mathbf{1}_G) =  \bil{T}{f*\check   h} =\bil{T\bullet h}{f}.
\end{align*}
This shows that $T\bullet h=T*h$ in this case. Similarly, we have $h\bullet T=h*T$. This allows us to
abandon the notation $T\bullet h$ and replace it with $T*h$, even when $T$ is not compactly supported.
Note, however, that in this case $T*h$ does not denote the convolution of two distributions, and one must
return to the definition to perform calculations. As an exercise, we leave it to the reader to verify
that
\[
(T*h)\, \check{} = \check h*\check T, \quad (T \in \scrD'(G)), \; (h\in \caH(G)),
\]
and that
\[
T*(h_1*h_2)=(T*h_1)*h_2   , \quad (T \in \scrD'(G)), \; (h_1,\, h_2 \in \caH(G)).
\]
The above calculation also shows that for $T \in \scrD'(G)$ and $h \in \caH(G)$, an equivalent definition
of $T*h$ is:
\begin{equation}\label{defequiv}
\bil{T*h}{f}= (T*(h*\check f)) (\mathbf{1}_G).
\end{equation}
\end{rmq}

\begin{defi}
We say that a \index[ter]{distribution!essentially compact}distribution $T$ in $ \scrD'(G)$ is
essentially compact on the right (resp. on the left), if for any $h\in \caH(G)$, $T*h$ (resp. $h*T$) is a
compactly supported distribution. We denote by $\scrD'(G)_{r.e.c.}$ (resp. $\scrD'(G)_{l.e.c.}$) the
space of essentially right-compact (resp. left-compact) distributions on $G$. We say that a distribution
$T$ in $ \scrD'(G)$ is essentially compact if it is essentially compact on the right and on the left. We
denote
\[
\scrD'(G)_{e.c.}= \scrD'(G)_{r.e.c.} \cap \scrD'(G)_{l.e.c.}.
\]
The space $\scrD'(G)_{e.c.}$ contains $\scrE'(G)$.
\end{defi}

It is clear that if $T$ is essentially compact on the right (resp. on the left), then $T*h$ (resp. $h*T$)
is in $\caH(G)$ for all $h\in \caH(G)$. Indeed, by definition, $T*h$ (resp. $h*T$) is compactly
supported. Moreover, there exists an idempotent $e_K$ of $\caH(G)$ such that $h*e_K=h$ (resp. $e_K*h=h$),
hence $T*h \in \scrE'(G,K,r)\subset \caH(G)$ (resp. $h*T \in \scrE'(G,K,l)\subset \caH(G)$).

\begin{lemme}
The space $\scrD'(G)_{r.e.c.}$ is equipped with a $\bbC$-algebra structure which extends that of
$\scrE'(G)$. The same is true for $\scrD'(G)_{l.e.c.}$ and $\scrD'(G)_{e.c.}$.
\end{lemme}

\begin{proof}
We must define the product $T_1*T_2$ of two distributions $T_1,T_2$ in $\scrD'(G)_{r.e.c.}$, then verify
that this product is a $\bbC$-algebra product. Let us set, for any $f$ in $\scrD(G)$
\[
\bil {T_1*T_2}{f}:= \bil{T_1*(T_2*e_K)}{f}
\]
where $e_K$ is an idempotent of $\caH(G)$ such that $f*e_K=f$. This is well-defined because
$e_K \in \caH(G)$, so $(T_2*e_K)\in \caH(G)$ since $T_2$ is essentially compact, and
$T_1*(T_2*e_K)\in \caH(G)$, since $T_1$ is essentially compact. Standard verifications show that this
definition is independent of the choice of the idempotent $e_K$ of $\caH(G)$. In the case where
$T_2=h \in \caH(G)$, the definitions of the products
\[
T_1*T_2 \quad \text{ and } \quad T_1*h
\]
coincide, as is easily verified using (\ref{defequiv}).

To verify the associativity of this product, let us first show by the above remark that for $T_1$, $T_2$
essentially right-compact and $h\in \caH(G)$, we have
\[
(T_1*T_2)*h=T_1*(T_2*h)   .
\]
It immediately follows that $ T_1*T_2$ is essentially right-compact and the associativity of the
convolution product of essentially right-compact distributions is then quickly established. Its
bilinearity is obvious. The algebra $\scrD'(G)_{r.e.c.}$ admits $\delta_{\mathbf{1}_G}$ as identity
element. A similar argument applies to $\scrD'(G)_{l.e.c.}$ by exchanging left and right. For
$\scrD'(G)_{e.c.}$, we have two definitions of the product (a right version, and a left one). Let us show
that they coincide. Let $T_1,T_2 \in \scrD'(G)_{e.c.}$ and $f \in \scrD(G)$. Let us choose the idempotent
$e_K$ such that
\[
e_K*f=f*e_K=f.
\]
We then have
\begin{align*}
\bil{T_1*(T_2*e_K)}{f}&=((T_1*(T_2*e_K))*\check f)(\mathbf{1}_G)=(T_1*(T_2*(e_K*\check f)))(\mathbf{1}_G)\\
&=(T_1*(T_2*(f*e_K)\, \check{}\, ))(\mathbf{1}_G)=(T_1*(T_2*\check f))(\mathbf{1}_G).
\end{align*}
Similarly for the left version, we show that
\[
\bil{(e_K*T_1)*T_2}{f}=((\check f*T_1)*T_2) (\mathbf{1}_G).
\]
We must therefore show that
\[
(T_1*(T_2*\check f))(\mathbf{1}_G)= ((\check f*T_1)*T_2)(\mathbf{1}_G).
\]
We know that if $h_1,h_2 \in \caH(G)$, then
\[
(h_1*(h_2*\check f))(\mathbf{1}_G)= ((\check f*h_1)*h_2) (\mathbf{1}_G).
\]
To be able to use this, we choose an idempotent $e_K$ such that $ e_K*f=f*e_K=f$,
$e_K*(T_2*\check f  )=(T_2*\check f  )*e_K=T_2*\check f$ and
$e_K*(\check f*T_1  )=(\check f*T_1  )*e_K=\check f *T_1$. We then have:
\begin{align*}
&(T_1*(T_2*\check f))(\mathbf{1}_G)= (T_1*(e_K*(T_2*\check f)))(\mathbf{1}_G)= ((T_1*e_K)*(e_K*(T_2*\check f)))(\mathbf{1}_G)\\
&=  ((T_1*e_K)*(e_K*(T_2*(e_K*\check f))))(\mathbf{1}_G)= ((T_1*e_K)*(e_K*T_2))*(e_K*\check f)(\mathbf{1}_G)\\
&=(e_K*\check f)*((T_1*e_K)*(e_K*T_2))(\mathbf{1}_G)=\check f *((T_1*e_K)*(e_K*T_2))(\mathbf{1}_G)\\
&=((\check f *T_1)*e_K)*(e_K*T_2)(\mathbf{1}_G)=((\check f *T_1)*T_2)(\mathbf{1}_G).
\end{align*}
\end{proof}

We therefore have the following inclusions:
\[
\caH(G)\subset \scrE'(G) \subset \scrD'(G)_{e.c.} \subset \scrD'(G)_{r.e.c.}
\]

\begin{thm}
The $\bbC$-algebra $\scrD'(G)_{r.e.c.}$ is identified with the completion $\overline {\caH(G)}$ of the
algebra with idempotents $\caH(G)$.
\end{thm}

\begin{proof}
Recall that an element $\bar a$ of $\overline {\caH(G)}$ is a map
\[
\bar a \colon \mathrm{Idem}(\caH(G)) \rightarrow \caH(G),
\]
such that $\bar a(e) \in \caH(G)*e$ for all $e \in \mathrm{Idem}(\caH(G)) $, and $\bar a(f)*e=\bar a (e)$
for any pair of idempotents $(e,f)$ such that $e\leq f$. To such an element $\bar a$ of
$\overline {\caH(G)}$, let us associate the distribution $A$ on $G$ defined by
\[
\bil{A}{f}= \bil{\bar a (e)}{f}
\]
for any $f \in \scrD(G)$ and any idempotent $e$ of $\caH(G)$ such that $f*e=f$ and $\check e=e$. We
verify that this is independent of the choice of the idempotent $e$. We will now show that $A$ is
essentially right-compact. For any $h\in \caH(G)$, for any $f \in \scrD(G)$,
\begin{align*}
 \bil{A*h}{f}&=\bil{A}{f*\check h}
 = \bil{\bar a (e)}{f*\check h} \quad \text{ where } \check h  *e= \check h\,  \text{ and } \check e=e \\
 &= \bil{\bar a (e)* h}{f}.
\end{align*}
This shows that $A*h=\bar a(e)*h$ where $e$ is such that $\check h  *e= \check h$ and $\check e=e$ (hence
$e*h=h$). In particular, $A*h$ is in $\caH(G)$ and therefore $A$ is essentially right-compact.

Conversely, let us define a map
\[
\scrD'(G)_{r.e.c.} \rightarrow \overline{\caH(G)}, \quad T \mapsto \overline{T}, \quad
\overline{T}(e)=T*e.
\]
It is easy to verify that this map is the inverse of the map $\bar a \mapsto A$ defined previously, and
that these are indeed $\bbC$-algebra morphisms. The theorem follows from this.
\end{proof}

\begin{cor}
The center of the category $\caM(\caH(G))$ is identified with the subspace of $\scrD'(G)_{e.c.}$ of
distributions invariant under conjugation.
\end{cor}

\begin{proof}
According to the results of Section \ref{idemmod}, the center of $\caM(\caH(G))$ is identified with the
center of $\scrD'(G)_{r.e.c.}$. It is also the commutant of $\caH(G)$ in $\scrD'(G)_{r.e.c.}$. A
distribution $T$ is therefore in the center if $T*h=h*T$, for all $h \in \caH(G)$. Note then that
$T \in \scrD'(G)_{e.c.}$. In other words, the center of $\scrD'(G)_{r.e.c.}$ is also the center of
$\scrD'(G)_{e.c.}$.

Let $T$ be in the center of $\scrD'(G)_{e.c.}$. We then have
\begin{equation}\label{centre1}
  \bil{T}{f*h}= \bil{T}{h*f}, \quad (f\in \scrD(G)), \; (h\in \caH(G)).
\end{equation}

On the other hand, a distribution $T$ invariant under conjugation is a distribution satisfying
$r(g)\cdot l(g)\cdot T=T$ for all $g \in G$, which can be rewritten
\[
\bil{l(g)\cdot T}{f}= \bil{r(g^{-1})\cdot T}{f}, \quad (f\in \scrD(G)),\, (g\in G),
\]
i.e.,
\[
\bil{ T}{l(g^{-1})\cdot f}= \bil{T}{r(g)\cdot f}, \quad (f\in  \scrD(G)), \, (g\in G),
\]
or again, substituting $g$ for $g^{-1}$,
\[
\bil{ T}{\delta_g * f}= \bil{T}{ f *\delta_g}, \quad (f\in \scrD(G)),\, (g\in G),
\]
that is
\begin{equation}\label{centre2}
\bil{ T}{e_K*\delta_g *e_K* f}= \bil{T}{ f*e_K *\delta_g*e_K}, \quad (f\in \scrD(G)), \, (g\in G),
\end{equation}
where $e_K$ is an idempotent of $\caH(G)$ such that
\[
e_K*f=f*e_K=f,\quad  e_K*(\delta_g*f)=\delta_g*f, \quad (f*\delta_g)*e_K=f*\delta_g.
\]

It is then clear that (\ref{centre1}) implies (\ref{centre2}) by taking $h=e_K*\delta_g*e_K$. The
converse is immediate by the formulas defining $h*T$ and $T*h$.
\end{proof}

\section{Notes on Chapter II}

Most of the results in this chapter are in \cite{BeZe1}. The section on sheaves on t.d. topological
spaces is largely inspired by the treatment given in \cite{Bump}. The description of the center and the
completion of $\caH(G)$ in terms of distributions is in \cite{Del}, but the notes from a talk by A. Moy
and remarks by M. Tadi\'c were very useful to me.

\chapter[Representations of t.d. groups]{Representations of totally disconnected groups}\label{chapreptd}

Having defined locally compact totally disconnected topological groups in the previous chapter, we now
tackle the representation theory of these groups. We only consider representations in complex vector
spaces, for which it turns out that the good category to consider is that of smooth representations,
i.e., those which are continuous when the representation space is equipped with the discrete topology.
This category is then equivalent to that of non-degenerate modules over the Hecke algebra of the group
defined in Chapter II. This allows us to use the constructions and results of Chapter I. We deduce from
it, for example, directly the properties of the functor associating to a representation the space of its
vectors fixed by a compact open subgroup, or those of the duality functor (contragredient
representation). The core of the chapter is devoted to the induction and forgetful functors between
categories of representations of two t.d. groups defined by a morphism between these two groups. Here
again, we particularize the results of Chapter I. The most important case is that where the morphism is
an inclusion of a closed subgroup. We recover in this case the usual notion of induced representation,
``Frobenius reciprocity'' being the fact that induction is the right adjoint of the forgetful functor. We
also identify in this framework compact induction and the functor $P$ of Chapter I, up to a modular
character, compact induction being therefore the left adjoint of the pseudo-forgetful functor. This is
particularly useful when the pseudo-forgetful functor is naturally isomorphic to the forgetful functor,
for example when the subgroup is open. Another important type of functor between categories of
representations is given by the construction of the space of ``coinvariants''. In a particular case which
will be used in Chapter VI for the construction of the parabolic induction of reductive groups, the
Jacquet functor is a functor of type $P$, hence left adjoint to the pseudo-forgetful functor (which is
isomorphic to the forgetful functor in this case).

Another fundamental tool of representation theory is Schur's lemma. The version we need here is proved in
the appendix. Irreducible representations admit a central character. More generally, the center of the
category of smooth representations acts on irreducible representations by scalars. We also show that
there are enough irreducible smooth representations in the sense that irreducible representations
``separate the points'' of the Hecke algebra of the group.

To summarize, our object of study in this chapter is not so much the category of smooth representations
of a totally disconnected group as the functors that can be defined between such categories. These
functors are particular cases of the induction and forgetful functors of Chapter I, via the equivalences
of categories with non-degenerate modules over Hecke algebras.

\section{Smooth representations}
\index[ter]{representation!smooth}
\subsection{Generalities}
Let $G$ be a group. A representation $(\pi,V)$ of $G$ is the data of a complex vector space $V$ and a
group morphism $\pi$ from $G$ to $\mathrm{GL}(V)$. This defines a linear action of $G$ on $V$. A
representation $(\pi,V)$ is said to be unitary if $V$ is equipped with a Hermitian inner product
(antilinear with respect to the first variable) $\bilo_V$ satisfying
\[
\bil{\pi(g)\cdot v}{w}_V=\bil{v}{\pi(g^{-1})\cdot w}_V, \quad (g\in G), (v,w \in V).
\]

If $W$ is a subspace of $V$ stable under the action of $G$ given by $\pi$, we denote by $\pi_W$ the
morphism from $G$ to $\mathrm{GL}(W)$ induced by $\pi$. The representation $(\pi_W,W)$ (often simply
denoted $(\pi,W)$) is a subrepresentation of $(\pi,V)$. Similarly, $\pi$ induces a morphism $\pi_{V/W}$
from $G$ to $\mathrm{GL}(V/W)$. The representation $(\pi_{V/W},V/W)$ (often simply denoted $(\pi,V/W)$)
is a quotient representation of $(\pi,V)$. More generally, if $W_1 \subset W_2$ are two subspaces of $V$
stable under the action of $G$, the morphism $\pi_{W_2/W_1}$ from $G$ to $\mathrm{GL}(W_2/W_1)$ induced
by $\pi$ defines a representation $(\pi_{W_2/W_1}, W_2/W_1)$ (or simply $(\pi, W_2/W_1)$) which is said
to be a subquotient of $(\pi,V)$.

A representation $(\pi,V)$ of $G$ with $V$ of dimension $1$ is called a character. \indexter{character}
We can in this case identify $V$ with $\bbC$, and we then sometimes denote this character by $\bbC_\pi$.
If it takes values in the complex numbers of modulus $1$, it is said to be unitary.
\indexter{unitary!(character)}

\begin{defi}
Let $G$ be a t.d. group. A representation $(\pi,V)$ of $G$ in a complex vector space $V$ is said to be
smooth\index[ter]{smooth}, if for all $v$ in $V$, there exists an open subgroup of $G$ which fixes $v$.
\end{defi}

\begin{rmq}
Let us equip $V$ with the discrete topology. If we assume that for all $v$ in $V$, the map:
\[
\phi_v \colon G \rightarrow V,\quad g\mapsto \pi(g)\cdot v
\]
is continuous, then the inverse image of the open set $\{ v\}$ of $V$ is an open set, which shows that
$\pi$ is smooth. Conversely, if $(\pi,V)$ is smooth, the inverse image of any open set of $V$, which is
the union (composed of open sets) of its points, is a union of open sets, hence an open set.
\end{rmq}

If $(\pi,V)$ is any representation of $G$, let $V_0$ denote the set of vectors of $V$ fixed by an open
subgroup of $G$. It is clear that $V_0$ is a subspace of $V$ stable under the action of $G$. The
representation $(\pi,V_0)$ is then trivially smooth. It is called the smooth part of $(\pi,V)$.

\subsection{Examples}

The group $G$ acts on $\scrC^\infty(G)$ by the left and right regular representations:
\[
(\forall f \in \scrC^\infty(G)), \, (\forall  g\in G),\,  (\forall  x\in G),\quad
l(g)\cdot f\, (x)=f(g^{-1}x), \quad r(g)\cdot f\,  (x)= f(xg).
\]

\begin{prop}
The representations $l$ and $r$ preserve the subspaces $\scrC^\infty_{unif}(G)$ and $\scrD(G)$. The
subrepresentations $(l,\scrC^\infty_{unif}(G))$, $(r,\scrC^\infty_{unif}(G))$, $(l,\scrD(G))$, and
$(r,\scrD(G))$ are smooth.
\end{prop}

\begin{proof}
Let $f \in \scrC^\infty_{unif}(G)$ and let $K$ be such that $f \in \scrC^\infty(G,K)$. We then have, for
all $k,k' \in K$, $x,g \in G$,
\[
(l(g)\cdot f)(gkg^{-1}xk')=f(kg^{-1}xk')=f(g^{-1}x)=(l(g)\cdot f)(x).
\]

Thus $l(g)\cdot f$ is in $\scrC^\infty(G,K \cap gKg^{-1}) \subset \scrC^\infty_{unif}(G)$. The proof is
similar for $(r,\scrC^\infty_{unif}(G))$, $(l,\scrD(G))$, $(r,\scrD(G))$. If $f \in \scrC^\infty(G,K)$,
$l(K)$ and $r(K)$ fix $f$, by definition, so $(l,\scrC^\infty_{unif}(G))$, $(r,\scrC^\infty_{unif}(G))$
are smooth. The same is true for $(l,\scrD(G))$ and $(r,\scrD(G))$ which are subrepresentations.
\end{proof}

\subsection{The category $\caM(G)$} \label{Kinv}

Let $(\pi,V)$ and $(\tau,W)$ be two smooth representations of the group $G$. An intertwining operator
\index[ter]{intertwining operator} $f$ between $(\pi,V)$ and $(\tau,W)$ is an element of $\Hom_\bbC(V,W)$
satisfying
\[
f(\pi(g)\cdot v)=\tau(g)\cdot f(v),\quad (v\in V),\, (g\in G).
\]
Let $\Hom_G((\pi,V),(\tau,W))$ (or for short simply $\Hom_G(V,W)$, or even $\Hom_G(\pi,\tau)$) denote the
space of intertwining operators between $(\pi,V)$ and $(\tau,W)$.

Let $\caM(G)$\index[not]{M(G)@$\caM(G)$} denote the category of smooth representations of $G$, the
morphisms being the intertwining operators. It is an abelian category. We do not prove it now, although a
direct verification is elementary, because this follows from the equivalence of categories of Section
\ref{eqcat}. Let $\mathbf{Irr}(G)$\index[not]{Irr(G)@$\mathbf{Irr}(G)$} denote the set of isomorphism
classes of irreducible smooth representations of $G$.

Given a smooth representation $(\pi,V)$ of $G$, and an element $T$ of $\scrE'(G)$, we define on $V$ the
operator $\pi(T)$ as follows. Let
\[
\phi_v\colon G\rightarrow V,\quad g\mapsto \pi(g)\cdot v.
\]
It is a locally constant function with values in $V$, an element of $\scrC^\infty(G , V)$ (\ref{XVXV}).
The duality between $\scrC^\infty(G,V)$ and $\scrE'(G)$ (see \ref{FetD})
\[
\langle\,  .\, ,\, .\,  \rangle\colon \scrC^\infty(G, V) \times\scrE'(G) \rightarrow V
\]
allows us to define
\[
\pi(T)\cdot v=\langle T ,\phi_v \rangle = \int_G \pi(g)\cdot v \; dT(g).
\]
One easily verifies that:
\[
\pi(T_1*T_2)=\pi(T_1)\circ \pi(T_2),\quad  \pi(\delta_g)=\pi(g),\quad (T_1,T_2 \in  \scrE'(G)), \, (g\in G).
\]
This defines an $\scrE'(G)$-module structure (and by restriction, an $\caH(G)$-module structure) on $V$.

\begin{exemple}
Consider the smooth representation $(l,\scrD(G))$. Then for all $T \in \scrE'(G)$, for all
$f\in \scrD(G)$,
\[
l(T)\cdot f= T*f.
\]
\end{exemple}

\begin{rmq}
Let us fix a Haar measure $\mu_G$ on $G$, and identify $\scrD(G)$ with $\caH(G)$ via the isomorphism
$f\mapsto f \mu_G$. By transport of structure, if $(\pi,V)$ is a smooth representation of $G$, we obtain
a representation of $\scrD(G)$, given explicitly by the formula:
\[
\pi(f)= \int_G f(g)\pi(g) \, d\mu_G(g), \quad (f \in \scrD(G)).
\]
\end{rmq}

\subsection{An equivalence of categories} \label{eqcat}
\label{replongfin} Any smooth representation $(\pi,V)$ of $G$ induces an $\caH(G)$-module structure on
$V$.
\begin{thm}
If $(\pi,V)$ is a smooth representation of $G$, the $\caH(G)$-module $V$ is non-degenerate. Conversely,
any non-degenerate $\caH(G)$-module $V$ defines a smooth representation of $G$. These two operations are
inverse to each other and induce an equivalence of categories between $\caM(G)$ and $\caM(\caH(G)).$
\end{thm}

\begin{proof}
Let $v \in V$, and let $K$ be a compact open subgroup of $G$ which fixes $v$. Then $\pi(e_K)\cdot v =v$
and therefore $v$ is fixed by an idempotent of $\caH(G)$. This shows that the $\caH(G)$-module $V$ is
non-degenerate. Moreover, any morphism in $\caM(G)$ is obviously a morphism in $\caM(\caH(G))$.
Conversely, if $V$ is a non-degenerate $\caH(G)$-module, for all $v \in V$, there exists a compact open
subgroup $K$ of $G$ such that $e_K \cdot v =v$. If $T \in \scrE'(G)$, we define
$T\cdot v:= (T* e_K) \cdot v$. This does not depend on the choice of $K$, because this construction is a
particular case of the construction of Lemma \ref{idemmod}: indeed any distribution $T \in \scrE'(G)$
defines an operator $h \mapsto T*h$ in $\caH(G)$ which commutes with right convolution in $\caH(G)$.
Lemma \ref{idemmod} then allows us to functorially define $\pi(T): V\rightarrow V$. We can then set
$\pi(g)\cdot v:=\pi(\delta_g)\cdot v$, $v\in V$. It is also clear that any morphism in $\caM(\caH(G))$
induces a morphism in $\caM(G)$. The two functors thus defined establish an equivalence of categories.
\end{proof}

\medskip

\begin{rmq}
By virtue of this equivalence of categories, we will denote the action of an element $T$ of $\caH(G)$ on
a vector $v$ of a smooth representation $(\pi,V)$ of $G$, either by $\pi(T)\cdot v$, or simply by
$T\cdot v$. This latter notation will be preferred when the underlying representation $\pi$ is clearly
identified by the context.
\end{rmq}

In Appendix \ref{JoHo}, we recall the notion of objects of finite length \index[ter]{finite length} and
of composition series\index[ter]{composition series} (also called Jordan-Hölder
series\index[ter]{Jordan-Hölder}) in certain abelian categories. We can apply this to modules in
$\caM(\caH(G))$ and by the equivalence of categories between $\caM(\caH(G))$ and $\caM(G)$, we transpose
the terminology of \ref{JoHo} to smooth representations of $G$.

Thus any non-zero smooth representation of $G$ admits an irreducible subquotient; if it is finitely
generated, it admits an irreducible quotient; and if it is of finite length, it admits an irreducible
subrepresentation.

\subsection{Fixed points of $K$}\label{piK}

Let $K$ be a compact open subgroup of $G$, and let $e_K$ be the idempotent of $\caH(G)$ given by the
normalized Haar measure on $K$. Any smooth representation $(\pi,V)$ of $G$ is a non-degenerate
$\caH(G)$-module, according to the equivalence of categories established in the previous paragraph. We
defined in Section \ref{exactje} a functor $j_{e_K}$, denoted more simply by $j_K$, from $\caM(\caH(G))$
to $\caM(e_K*\caH(G)*e_K)$. Recall that $e_K*\caH(G)*e_K=\caH(G,K)$. We give some properties of this
functor.

\begin{prop}
$(i)$ The module $j_K(V)=e_K\cdot V$ is the subspace $V^K$ of $V$ of vectors fixed under the action of
$K$ in $V$.

$(ii)$ The kernel of $\pi(e_K)$ in $V$ is the subspace $V(K)$ of $V$ generated by vectors of the form
$\pi(e_K)\cdot v-v$, $v\in V$.

$(iii)$ We have $V=V^K \oplus V(K)$.

$(iv)$ The functor $j_K$ is exact.
\end{prop}
\begin{proof}
Point $(i)$ is immediate. It is clear that
\[
e_K\cdot(e_K\cdot v-v)=0,\quad  (v\in V),
\]
so $V(K) \subset \ker e_K$. On the other hand, since $\pi(e_K)$ is a projector,
\[
V=\ker e_K \oplus \im e_K.
\]
Since $v= e_K \cdot v - (e_K \cdot v-v)$, we have $V=e_K \cdot V + V(K)$. We deduce that
$V(K) = \ker \pi(e_K)$ and we obtain $(ii)$ and $(iii)$. The exactness of $j_K$ was proved in Section
\ref{exactje}.
\end{proof}

\medskip

We restate, in the more specific framework of this section, the results of Proposition \ref{35}.

\begin{thm}
$(i)$ The representation $(\pi,V)$ is irreducible or zero if and only if for any compact open subgroup
$K$ of $G$, the $\caH(G,K)$-module $V^K$ is simple or zero.

$(ii)$ Let $(\pi_i,V_i)$, $i=1,2$, be two irreducible smooth representations of $G$, and let $K$ be a
compact open subgroup of $G$ such that $V_i^K \neq 0$. Then $\pi_{1}$ and $\pi_{2}$ are equivalent if and
only if $V_1^K$ and $V_2^K$ are isomorphic $\caH(G,K)$-modules.

$(iii)$ For each simple module $W$ of $\caH(G,K)$, there exists an irreducible representation $\pi$ of
$G$ such that $V^K=W$.
\end{thm}

\subsection{Contragredient representation}\label{contrag}

Let $(\pi,V)$ be a smooth representation of $G$. The dual space $V^*$ is equipped with the representation
$\pi^*$ of $G$ defined by:
\[
\pi^*(g)\cdot\lambda(v)=\lambda(\pi(g^{-1})\cdot v), \quad (\lambda \in V^*), (v\in V).
\]
This representation is not necessarily smooth. Let $\widetilde V$ be the smooth part of $V^*$. It is
clear that $\widetilde V$ is stable under the action of $G$. We will denote by
$(\tilde \pi,\widetilde V)$ the restriction of $\pi^*$ to $\widetilde V$ \index[not]{Vt@$\widetilde V$}:
it is a smooth representation of $G$, which is called the contragredient
representation\index[ter]{contragredient} of $(\pi,V)$.

In Section \ref{dualite}, we constructed a duality functor from $\caM(\caH(G))$ to $\caM(\caH(G))_r$.
Since the algebra $\caH(G)$ is equipped with the anti-involution $T \mapsto \check T$ induced by
$g \mapsto g^{-1}$, we can identify right and left modules, and therefore view the duality functor as an
endofunctor of $\caM(\caH(G))$. It is clear that by the equivalence of categories of Theorem \ref{eqcat},
the contragredient representation $(\tilde \pi,\widetilde V)$ of $(\pi,V)$ corresponds to the
non-degenerate $\caH(G)$-module $\widetilde V$ defined in Section \ref{dualite}.

\begin{lemme}
Let $T \in \scrE'(G)$. For all $\lambda \in \widetilde V$, for all $v\in V$, we have:
\[
(\tilde \pi (\check T)\cdot \lambda) (v)= \lambda ( \pi (T)\cdot v).
\]
In particular,
\[
(\tilde \pi (e_K)\cdot \lambda) (v)= \lambda ( \pi (e_K)\cdot v)
\]
for any compact subgroup $K$. If $\lambda$ is moreover $K$-invariant, we then have:
\[
\lambda(v)=\lambda(\pi(e_K)\cdot v).
\]
Therefore $(V^*)^K\simeq (V^K)^*$, and if $K$ is open, $\widetilde V^K=(V^*)^K\simeq (V^K)^*$.
\end{lemme}
\begin{proof}
Recall that the homeomorphism $g\mapsto g^{-1}$ of $G$ induces the endomorphism $T\mapsto \check T$ of
$\scrE'(G)$. We have:
\begin{align*}
(\tilde \pi (\check T)\cdot \lambda) (v)&=\int_G (\tilde \pi (g)\cdot \lambda) (v)\;  d\check T(g)=
 \int_G  \lambda (\pi(g^{-1})\cdot v)\; d\check T(g)=  \int_G  \lambda (\pi(g)\cdot v)\, dT(g)\\
& =\lambda ( \pi (T)\cdot v).
\end{align*}
Since $K$ is a unimodular group, we have $\check e_K=e_K$, hence the second assertion of the lemma. The
rest is a reformulation of Lemma \ref{dualitelem}.
\end{proof}

\medskip

Let us also reformulate Corollary \ref{dualitecor}:
\begin{cor}
We have a canonical injection $V\hookrightarrow \widetilde {\widetilde V}$.
\end{cor}

\subsection{Admissible representations} \label{Admis}
The admissibility condition of a representation is a technical finiteness criterion allowing one to
obtain many important results. We will show later that in the case of a $p$-adic reductive group, all
irreducible smooth representations are admissible\index[ter]{admissible!(representation)}.

\begin{defi}
A smooth representation $(\pi,V)$ of the t.d. group $G$ is said to be admissible, if for any open
subgroup $N$ of $G$, the subspace of fixed vectors $V^N$ is finite-dimensional.
\end{defi}

Let us translate this: since the compact open subgroups form a neighborhood basis of the identity in $G$,
the representation $(\pi,V)$ is admissible if and only if $V^K$ is finite-dimensional for any compact
open subgroup $K$ of $G$. In terms of an $\caH(G)$-module, this means that $e_K \cdot V$ is
finite-dimensional for any $K$ as above. But for any idempotent $e$ of $\caH(G)$, there exists such a $K$
such that $e \leq e_K$. We therefore have $e \cdot V=e_Ke\cdot V \subset e_K \cdot V$, and we see that
the representation $(\pi,V)$ is admissible if and only if the non-degenerate $\caH(G)$-module $V$ is
admissible in the sense of \ref{modadm}.

\begin{prop}
A smooth representation $(\pi,V)$ of $G$ is admissible if and only if
$V \simeq \widetilde{\widetilde V}$.
\end{prop}
\begin{proof}
This is the translation of Proposition \ref{modadmprop}. Corollary \ref{modadmcor} can also be stated in
terms of representations rather than $\caH(G)$-modules.
\end{proof}

\subsection{Schur's lemma} \label{schur}
\index[ter]{Schur's lemma}

\begin{prop}
Let $G$ be a $\sigma$-compact t.d. group, and let $(\pi,V)$ be an irreducible smooth representation of
$G$. Then the space $\Hom_G(V,V)$ is the set of scalar operators. If $(\tau,W)$ is another irreducible
smooth representation of $G$ not equivalent to $(\pi,V)$, then $\Hom_G(V,W)=0$.
\end{prop}
\begin{proof}
Let us first show that $V$ is of at most countable dimension: let us fix a non-zero $v\in V$ and let $K$
be a compact open subgroup of $G$ fixing $v$. Then by irreducibility, $V$ is generated by the vectors
$\pi(g)\cdot v$, where $g$ runs through a system of representatives of $G/K$. But such a system of
representatives is countable (Lemma \ref{fini}).

The proposition then follows from Theorem \ref{ASchur} and the equivalence of categories \ref{eqcat}:
\[
\Hom_G(V,V) \simeq \End_{\caH(G)} (V).
\]
The second point is obvious, since a non-zero intertwining operator between $V$ and $W$ is necessarily an
isomorphism, which is excluded by hypothesis.
\end{proof}

\begin{rmq}
If in the proposition, we assume $(\pi,V)$ is admissible, we can replace the hypothesis that $G$ is
$\sigma$-compact by the fact that the identity element of $G$ admits a countable neighborhood basis.
Indeed, $V$ is of countable dimension, since $V=\bigcup_K V^K$, where $K$ runs through a (countable)
family of compact open subgroups of $G$ forming a neighborhood basis of the identity, and each $V^K$ is
finite-dimensional.
\end{rmq}

\subsection{An application of Schur's lemma}\label{bilcan}
The following result is an application of Schur's lemma that we will use later. We therefore still assume
$G$ is $\sigma$-compact.
\begin{prop}
Let $(\pi,V)$ be an irreducible admissible representation of $G$. If
$B\colon V \times \widetilde V \rightarrow \bbC$ is a $G$-invariant bilinear form, then there exists a
constant $c_B$ in $\bbC$ such that
\[
B(v,\lambda)=c_B \; \lambda(v).
\]
In particular, if $B$ is non-zero, it is non-degenerate.
\end{prop}
\begin{proof}
Since $(\pi,V)$ is admissible and irreducible, recall that the same is true for
$(\tilde \pi,\widetilde V)$ and that $\widetilde{\widetilde V}\simeq V$. If $B$ is degenerate, there
exists a non-zero element $\lambda \in \widetilde V$ such that $B(v,\lambda)=0$ for all $v\in V$.
Therefore the kernel of the map
\[
\alpha \colon \widetilde V\rightarrow \widetilde V,\quad \lambda\mapsto \left(v\mapsto B(v,\lambda) \right)
\]
is non-trivial and $G$-stable. Since $\widetilde V$ is irreducible, we have $\ker \alpha=\widetilde V$
and therefore $B=0$. One argues in the same way with the left kernel
$\{v\in V\,|\,B(v,\lambda)=0, \ \forall \lambda\in\widetilde V\}$, which is a $G$-stable subspace of $V$,
using the irreducibility of $V$.

If $B$ is non-zero, consider two isomorphisms between $V$ and $\widetilde{\widetilde V}$: the first is
the canonical isomorphism
\[
\beta_1: v \mapsto \left(\lambda \mapsto \lambda(v)\right),
\]
and the second is given by $\beta_2: v \mapsto \left(\lambda \mapsto B(v,\lambda)\right)$. According to
Schur's lemma, $\beta_1$ and $\beta_2$ are proportional.
\end{proof}

\subsection{Central character}\label{caraccentral}
\index[ter]{character!central}

Suppose that $G$ is $\sigma$-compact and let $Z$ be a subgroup contained in the center $Z(G)$. Let
$(\pi,V)$ be an irreducible smooth representation of $G$. According to Schur's lemma, for any element
$z\in Z$, $\pi(z)$ acts as a scalar in $V$. Let us denote:
\[
\pi(z)=\chi_\pi(z) \; \Id_V, \quad (z\in Z).
\]
It is clear that $\chi_\pi$ is a character of $Z$. When $Z=Z(G)$, it is called the central character of
the representation $(\pi,V)$.

If $(\pi,V)$ is a smooth representation of $G$, and if there exists a character $\chi$ of $Z(G)$ such
that:
\[
\pi(z)=\chi(z) \; \Id_V, \quad (z\in Z(G)),
\]
we say that $(\pi,V)$ admits the central character $\chi$ (it is of course uniquely determined).

\begin{lemme}
If $(\pi_1,V_1)$ and $(\pi_2,V_2)$ are two smooth representations of $G$, admitting respectively the
central characters $\chi_1$ and $\chi_2$, which are assumed to be distinct, then there are no non-trivial
intertwining operators between $V_1$ and $V_2$.
\end{lemme}
\begin{proof}
Let $z\in Z(G)$ such that $\chi_1(z)\neq \chi_2(z)$ and $A\in \Hom_G(V_1,V_2)$. Then for all $v \in V_1$,
\[
\chi_1(z) A(v)=A(\chi_1(z)v)=A(\pi_1(z)\cdot v)=\pi_2(z)\cdot A(v)=\chi_2(z)A(v),
\]
and therefore $A(v)=0$.
\end{proof}

\enlargethispage{\baselineskip} \needspace{144pt}
\begin{prop}
If $(\pi,V)$ is a smooth representation of $G$ that is admissible or of finite length, then, for any
character $\chi$ of $Z$, let us set
\[
V_\chi\index[not]{Vchi@$V_\chi$}=\{ v \in V\mid \exists m\in \bbN,\, \forall z\in Z, \, (\pi(z)-\chi(z)\Id)^m\cdot v=0   \},
\]
and call this space the characteristic subspace of $V$ for the character $\chi$. Then
$V=\bigoplus_{\chi} V_\chi$. If $(\pi,V)$ is of finite length, only finitely many factors are non-zero.
\end{prop}

We denote by $\mathrm{Exp}(Z,V)$\index[not]{Exp@$\mathrm{Exp}(Z,V)$} the set of characters $\chi$ of $Z$
such that $V_\chi \neq \{0\}$.

\begin{proof}
If $(\pi,V)$ is of finite length, we proceed by induction on the length, the case where $(\pi,V)$ is
irreducible being treated by Schur's lemma. If $V$ is not irreducible, we consider an irreducible
subrepresentation $V_0$, and we use the induction hypothesis to decompose
$\overline{V}=V/V_0= \bigoplus_\chi \overline{V}_\chi$. Let $\chi_0$ be the character by which $Z$ acts
in the irreducible representation $V_0$. Let $\bar v \in \overline{V}_\chi$, and $v\in V$ a lift of
$\overline{v}$. We then have, for some $m\in \bbN$,
\[
(\pi(z)-\chi(z)\Id)^m\cdot v \in V_0, \quad (z \in Z),
\]
and therefore
\[
(\pi(z)-\chi_0(z)\Id)(\pi(z)-\chi(z)\Id)^m\cdot v=0,  \quad (z \in Z).
\]
We conclude by classical arguments that
\[
V=\bigoplus_{\chi} V_\chi \quad \text{where} \quad \chi \in \mathrm{Exp}(Z,\overline{V})\cup \{\chi_0\} .
\]
In the case where $(\pi,V)$ is admissible, for any compact open subgroup $K$ of $G$, we decompose $V^K$,
which is finite-dimensional, into characteristic subspaces $V^K_\chi$ for the action of the elements of
$Z$ (since the actions of $Z$ and $K$ commute, $V^K$ is stable under the action of $Z$). The
decompositions obtained are compatible with the inclusions $V^{K_1} \subset V^{K_2}$ when
$K_2 \subset K_1$, and we then set
\[
V_\chi=  \bigcup_K V^K_\chi.
\]
The decomposition $V=\bigoplus_\chi V_\chi$ easily follows from this.
\end{proof}

\begin{rmq}
We defined in \ref{idemmod} the center of the category $\caM(\caH(G))$ as being the algebra of
endomorphisms of the identity functor. Let us denote it $\frZ(G)$. It is clear that an element of
$Z=Z(G)$ defines, by its action on the modules of $ \caM(G)\simeq \caM(\caH(G))$, such an endomorphism,
and therefore an element of $\frZ(G)$. But, for any element $z$ of $\frZ(G)$, Schur's lemma asserts that
the action of $z$ on a simple module $V$ is a scalar operator. As above, this therefore defines for any
simple module $V$ a character of $\frZ(G)$, which is still called the central character of the module $V$
and which extends the character of $Z(G)$ described above.
\end{rmq}

\subsection{Completeness of irreducible representations}\label{separlemme}

In this paragraph, we show that if $G$ is a $\sigma$-compact t.d. group, then $G$ has enough irreducible
smooth representations, in the following sense.

\begin{thm}
Let $T \in \caH(G)$, $T\neq 0$. Then there exists an irreducible smooth representation $\pi$ of $G$ such
that $\pi(T)\neq 0$.
\end{thm}

\begin{proof}
For a sufficiently small compact open subgroup $K$, we have $T=e_K*T*e_K\neq 0$, so we can assume that
$T \in \caH(G,K)$. Write $T=f \mu_G$, $f\in \scrD(G,K)$, and set $f^*(g)=\overline{ f(g^{-1})}$,
$g \in G$. We verify directly using (\ref{fconvg}) that:
\[
(f^* * f)(\mathbf{1}_G)=\int_G |f(g^{-1})|^2  \, d\mu_G(g),
\]
and therefore, since $f\neq 0$, we find $(f^* * f)(\mathbf{1}_G)\neq 0$. Set $h= f^* * f$. We have
$h^*=h$ and $h\neq 0$, hence by the same reasoning $h^{*2}= h *h\neq 0$, $h^{*4}= h^{*2} *h^{*2}\neq 0$,
etc. Set $S=h \mu_G$. We then have $S*S=  (h *h)\,  \mu_G \neq 0$, $S*S*S*S=  h^{*4}   \, \mu_G\neq 0$,
etc. Thus $S$ is not a nilpotent element. On the other hand, $S$ lies in the algebra $\caH(G,K)$, which
is of countable dimension because $G/K$, and therefore \textit{a fortiori} $K\backslash G/K$, is discrete
and countable (Lemma \ref{fini}). We can conclude using Lemma \ref{nilp} that there exists a simple
$\caH(G,K)$-module $W$ such that $S\cdot W\neq 0$, and then, by Theorem \ref{piK}, $(iii)$, that there
exists an irreducible smooth representation $(\pi,V)$ of $G$ such that $V^K\simeq W$. We then have
$\pi(S)\neq 0$, hence $\pi(T)\neq 0$.
\end{proof}

\begin{rmq}
The algebra $\caH(G)$ therefore satisfies property ({\bf Sep}) of \ref{centF}. Thus, besides the
description of $\frZ(G)$ obtained in \ref{essencomp}, one can try to determine it by characterizing its
image under the morphism (\ref{centrefonctions}).
\end{rmq}

\subsection{Characters of admissible representations}

Let $(\pi,V)$ be an admissible representation of $G$, let $T\in \caH(G)$ and let $K$ be a compact open
subgroup of $G$ such that $T=e_K*T*e_K$. The operator $\pi(T)$ then has its image in $e_K\cdot V=V^K$,
which is finite-dimensional. We can therefore define the trace of $\pi(T)$, denoted by $\tr \pi(T)$ or
$\Theta_\pi(T)$. We can view $\tr \pi$ as a linear form on $\caH(G)$. If we fix a left Haar measure
$\mu_G$, for any $f\in \scrD(G)$, we set $\Theta_\pi(f)=\tr \pi(f):= \tr \pi(f \mu_G)$. We thus define a
distribution on $G$, which is called the \emph{character} \index[ter]{character} of $\pi$.

\begin{lemme}
Let $g_0\in G$ and $\Int(g_0)$ be the inner automorphism of $G$ given by conjugation by $g_0$. We then
have, for any admissible representation $(\pi,V)$ of $G$:
\[
\Int(g_0)\cdot( \tr \pi)=\delta_G(g_0)\,  \tr \pi.
\]
In particular if $G$ is unimodular, $\tr \pi$ is a distribution invariant under conjugation.
\end{lemme}

\begin{proof}
We have for all $T \in \caH(G)$:
\begin{align*}
(\Int(g_0)\cdot  \tr \pi)(T)&=\tr \pi(\Int(g_0^{-1})(T))=
\tr \pi (\delta_{g_0^{-1}}*T* \delta_{g_0})
= \tr(\pi(g_0^{-1})\pi(T)\pi(g_0))\\
&=\tr \pi(T)
\end{align*}
Therefore if $f \in \scrD(G)$, we have:
\begin{align*}
 (\Int(g_0)\cdot\tr \pi) (f)&=  \tr \pi(\Int(g_0^{-1}) (f))=
 \tr \pi(\Int(g_0^{-1}) (f) \, \mu_G)
 =\delta_G(g_0)\,  \tr \pi(\Int(g_0^{-1})\cdot   (f \mu_G))\\
 &=\delta_G(g_0)\,  \tr \pi( f \mu_G)=\delta_G(g_0)\,  \tr \pi(f).
\end{align*}
\end{proof}

\begin{exemple}
Let $(\pi,V)$ be an admissible representation of $G$. The character of the contragredient representation
$(\tilde \pi,\widetilde V)$ is given by
\[
\Theta_{\tilde \pi}(T)=\Theta_\pi(\check T ), \quad (T\in \caH(G)), \quad
\Theta_{\tilde \pi}(f)=\Theta_\pi(\check f  \delta_G^{-1} ), \quad
(f \in \scrD(G)).
\]
Recall that $\check f$ is defined by $\check f(g)=f(g^{-1})$.
\end{exemple}

\subsection{Linear independence of characters}

Let $G$ be a t.d. group and $(\pi_1,V_1)$, \ldots , $(\pi_n,V_n)$ be irreducible admissible
representations of $G$, pairwise non-equivalent. We then have:

\begin{prop}
The characters $\tr \pi_1$, \ldots, $\tr \pi_n$ are linearly independent.
\end{prop}

\begin{proof}
Since the $\pi_i$ are smooth, we can find a compact open subgroup $K$ such that for all $i$,
$V_i^{K}\neq \{0\}$. According to Theorem \ref{piK}, the $\caH(G,K)$-modules $V_i^K$ are simple,
non-isomorphic and finite-dimensional. Linear independence then follows from \cite{BourAlgVIII}, Chap.
VIII, \S 20, no.~6, Prop. 6 (the traces are non-zero since $\tr(e_K)=\dim V_i^K$).
\end{proof}

\begin{cor}
Two irreducible admissible representations $\pi_1$ and $\pi_2$ of $G$ are equivalent if and only if their
characters are equal.
\end{cor}

\subsection{Tensor products of representations}\label{tensprod}

For all $i=1,2,\ldots ,r$, suppose that $G_i$ is a t.d. group. The product group $\prod_i G_i$ is then a
t.d. topological group for the product topology. Let $(\pi_i,V_i) \in \caM(G_i)$, for $i=1,2,\ldots ,r$.
We then define the tensor product:
\[
\boxtimes_i (\pi_i,V_i)=(\pi_1\otimes \ldots \otimes \pi_r,
V_1\otimes \ldots \otimes V_r)
\]
as being the representation of $\prod_i G_i $ in $V_1\otimes \ldots \otimes V_r$ given by:
\[
(\pi_1\boxtimes \ldots \boxtimes \pi_r)(g_1,\ldots ,g_r)=
\pi_1(g_1)\otimes \ldots \otimes \pi_r(g_r),\quad (g_i \in G_i, \, i=1, \ldots,r).
\]
The tensor products here are tensor products of vector spaces over $\bbC$. However, to lighten the
notation, we simply write $\otimes$ rather than $\otimes_\bbC$.

If $G_1=G_2=   \ldots = G_r=G$, we define a representation $\pi$ of $G$ in $\bigotimes_i V_i$ by
$\pi(g)=(\otimes_i \pi)(\Delta(g))$ where
\[
\Delta \colon G \rightarrow \prod_i G, \quad g \mapsto (g,\ldots, g).
\]
We denote this representation of $G$ by $\bigotimes_i (\pi_i,V_i)$.

\begin{prop}
With the above notation, suppose that all the representations $(\pi_i,V_i)$ are admissible. Then
$\boxtimes_i (\pi_i,V_i)$ is admissible. Furthermore,
$(\boxtimes_i (\pi_i,V_i))^\sim= \boxtimes_i (\tilde \pi_i,\widetilde   V_i)$. If all the representations
$(\pi_i,V_i)$ are irreducible, then the same is true for $\boxtimes_i  (\pi_i,V_i)$. Finally, any
irreducible admissible smooth representation $(\pi,V)$ of $\prod_i G_i$ is of the form
$ \boxtimes_i (\pi_i,V_i)$ where the $(\pi_i,V_i)$ are irreducible admissible and uniquely determined up
to isomorphism by $(\pi,V)$.
\end{prop}

\begin{proof}
It suffices to treat the case $r=2$. Let $K_i \subset G_i$, $i=1,2$, be compact open subgroups and set
$K = K_1 \times K_2 \subset G=G_1 \times G_2$. One easily verifies that
$\caH(G,K)=\caH(G_1,K_1)\otimes \caH(G_2,K_2)$ and that $(V_1\otimes V_2)^K=V_1^{K_1}\otimes V_2^{K_2}$,
so that $(V_1\otimes V_2)^K$ is finite-dimensional. The natural map:
\[
\widetilde V_1 ^{K_1}\otimes \widetilde V_2^{K_2}
=(V_1^{K_1})^* \otimes (V_2^{K_2})^* \rightarrow  ((V_1\otimes V_2)^K)^*
=\widetilde{(V_1\otimes V_2)}^K
\]
is an isomorphism. This proves the first two assertions. The rest then follows easily from Theorem
\ref{piK} and the following result, a proof of which can be found in \cite{Bump}, Prop. 3.4.1.
\end{proof}

\needspace{160pt}
\begin{scolie}
Let $A_1$ and $A_2$ be unital algebras and let $(\tau_1,V_1)$, $(\tau_2,V_2)$ be finite-dimensional
irreducible representations of $A_1$ and $A_2$ respectively. Then the representation
$$(\tau_1 \boxtimes \tau_2,V_1 \otimes V_2)$$ of $A_1\otimes A_2$ is irreducible. Conversely, any
finite-dimensional irreducible representation of $A_1\otimes A_2$ can be written in this form, uniquely
up to isomorphism.
\end{scolie}

If $(\pi,V)$ is a smooth representation of $G$ and if $(\omega,\bbC)$ is a (smooth) character of $G$,
then we will simplify the notation and write $(\pi\omega,V)$ for the representation
$(\pi\otimes \omega,V\otimes \bbC)$ (of course $V\otimes \bbC$ is canonically identified with $V$).

\begin{lemme}
Let $(\pi,V)$ be an admissible smooth representation of $G$ and $(\omega,\bbC)$ a (smooth) character of
$G$. Then $(\pi\omega,V)$ is admissible.
\end{lemme}

\begin{proof}
Let $K$ be a compact open subgroup of $G$ on which $\omega$ is trivial. Let $K'$ be another compact open
subgroup of $G$. Any vector $v$ of $V$ is fixed by $K'\cap K$ under $\pi$ if and only if it is fixed for
$\pi\omega$.
\end{proof}

\subsection{The functors $ \otimes$ and $\Hom$ }\label{OHom}

In the previous section, we defined the tensor product of representations of $\caM(G)$. If we fix a
representation $(\pi,V)$ of $\caM(G)$, we can then define a functor
\[
\bullet \otimes V \colon \caM(G)\rightarrow \caM(G),\quad (\rho,W)\mapsto (\rho\otimes \pi,W\otimes  V).
\]
We will define a $\Hom$ functor in a similar way. If $(\rho,W)$ and $(\pi,V)$ are two representations of
$G$ (not necessarily smooth), $\Hom(W,V)=\Hom_\bbC(W,V)$ is equipped with a representation $\pi^\rho$ of
$G$ given by: for all $w\in W$, for all $g\in G$ and all $\phi \in \Hom(W,V)$,
\[
(\pi^\rho(g)\cdot \phi)(w)=\pi(g)\cdot \phi(\rho(g^{-1})\cdot w).
\]
Even if $\rho$ and $\pi$ are smooth, the representation thus obtained is not necessarily so. To be able
to define a functor from the category $\caM(G)$ to itself, we must therefore take its smooth part
$\Hom(W,V)_0$. We still denote by $\pi^\rho$ the restriction of $\pi^\rho$ to the smooth part. Thus we
obtain a functor:
\[
\Hom(W,\, \bullet\, )_0 \colon \caM(G)\rightarrow \caM(G), \quad  (\pi,V)\mapsto (\pi^\rho,\Hom(W,V)_0)
\]

\medskip

We can also of course define the contravariant functor
\[
\Hom(\, \bullet\,  ,V)_0 \colon \caM(G)\rightarrow \caM(G), \quad  (\rho,W)\mapsto (\pi^\rho,\Hom(W,V)_0).
\]

\begin{exemple}
If we take for $(\pi,V)$ the trivial representation of $G$, we obtain the functor
$\rho\mapsto \tilde \rho$ which associates to any representation its contragredient.
\end{exemple}

\medskip

The two functors we have just defined are linked by an adjunction relation, as shown by the following
result.

\needspace{150pt}
\begin{prop}
Let $(\pi_i,V_i)$, $i=1,2,3$, be smooth representations of $G$. We then have a natural isomorphism (in
all variables)
\begin{align*}
 \Hom(V_1\otimes V_2,V_3)_0 &\simeq \Hom(V_1,\Hom(V_2,V_3)_0)_0\\
\phi& \mapsto \Phi, \quad \Phi(v_1)(v_2)=\phi(v_1\otimes v_2).
\end{align*}
\end{prop}

\begin{proof}
In the category of complex vector spaces, we have a natural isomorphism:
\[
\Hom(V_1\otimes V_2,V_3)\simeq  \Hom(V_1,\Hom(V_2,V_3))
\]
given by $\phi \mapsto \Phi$, $\Phi(v_1)(v_2)=\phi(v_1 \otimes v_2)$ (this is a particular case of
\ref{AD1}). We will now show that it intertwines the actions of $G$. For all $g \in G$, $g\cdot \phi$ is
the linear map
\[
(g\cdot \phi) \colon v_1\otimes v_2 \mapsto \pi_3(g)\cdot (\phi(\pi_1(g^{-1})\cdot v_1 \otimes \pi_2(g^{-1})\cdot v_2)),
\]
and $g\cdot \Phi$ is the linear map
\[
g\cdot \Phi:  V_1\rightarrow \Hom(V_2,V_3), \quad  v_1  \mapsto  g\cdot (\Phi(\pi_1(g^{-1})\cdot v_1))
\]
where
\begin{align*}
 g\cdot (\Phi(\pi_1(g^{-1})\cdot v_1))(v_2)&=   \pi_3(g)\cdot \left(    (\Phi(\pi_1(g^{-1})\cdot v_1))(\pi_2(g^{-1})\cdot v_2)\right)\\
& =   \pi_3(g)\cdot (  \phi( \pi_1(g^{-1})\cdot v_1 \otimes \pi_2(g^{-1})\cdot v_2)).
\end{align*}
Therefore $\phi \mapsto \Phi$ is $G$-equivariant. By taking the smooth part on each side, we obtain
\[
\Hom(V_1\otimes V_2,V_3)_0\simeq  \Hom(V_1,\Hom(V_2,V_3))_0.
\]
To be able to conclude, we use the following result, which we will use again later.
\end{proof}

\begin{lemme}
Let $(\pi_i,V_i)$, $i=1,2$ be representations of $G$ (not necessarily smooth). Then,
\[
\Hom((V_1)_0,V_2)_0= \Hom((V_1)_0,(V_2)_0)_0.
\]
\end{lemme}

\begin{proof}
Let $\phi \in \Hom((V_1)_0,V_2)_0$, $v_1 \in (V_1)_0$. Let us choose a compact open subgroup $K$ of $G$
small enough to simultaneously fix $\phi$ and $v_1$. Then for all $k\in K$,
\begin{align*}
& \pi_2(k^{-1})\cdot \phi(v_1)= \pi_2(k^{-1})\cdot
((k\cdot \phi)(\pi_1(k)\cdot v_1))\\
=&   \pi_2(k^{-1})\cdot (\pi_2(k)\cdot (\phi  (\pi_1(k^{-1})\cdot (\pi_1(k)\cdot v_1))))
=\phi(v_1).
\end{align*}
This shows that $\phi(v_1)$ is also fixed by $K$.
\end{proof}

By taking the $G$-invariant morphisms on each side of the equality established by the proposition, we
obtain

\begin{cor}\label{OHomcor}
$\Hom_G(V_1\otimes V_2,V_3)\simeq \Hom_G(V_1,\Hom(V_2,V_3)_0).$
\end{cor}

\begin{rmqs}
$ $
\begin{itemize}
\item[1.] The corollary shows that in particular the functor $\bullet \otimes V_2$ is the left adjoint of
the functor $\Hom(V_2,\, \bullet \,  )_0$. But the proposition also shows naturality in the variable
$V_2$.

\item[2.] The above adjunction implies that the functor $\bullet  \otimes V$ preserves colimits, hence in
particular is right exact and preserves inductive limits, while the functor $\Hom(V,\, \bullet\, )_0$
preserves limits, hence in particular is left exact and preserves projective limits.

\item[3.] If $(\omega,\bbC_\omega)$ is a 1-dimensional representation of $G$, then
$\bullet  \otimes \bbC_\omega$ realizes an equivalence of categories with inverse
$\bullet  \otimes \bbC_{\omega^{-1}}$. Moreover these two functors are adjoints, in both possible
directions, i.e., for any smooth representations $(\pi,V)$ and $(\rho,W)$ of $G$, we have a natural
isomorphism
\[
\Hom_G(V,W\otimes \bbC_\omega)\simeq  \Hom_G(V\otimes \bbC_{\omega^{-1}} ,W),
\]
the other adjunction being obtained by exchanging the roles of $\omega$ and $\omega^{-1}$. This shows in
this case that the functor $\bullet \otimes \bbC_{\omega}$ is a right and left adjoint and therefore that
it preserves limits and colimits. In particular it is exact.
\end{itemize}
\end{rmqs}

Let $(\pi_i,V_i)$, $i=1,2$ be smooth representations of $G$. Let us now look at the space of coinvariants
of $V_1\otimes V_2$ for the action of $G$: by definition, it is the quotient of $V_1\otimes V_2$ by the
subspace generated by the vectors of the form
$\pi_1(g)\cdot v_1\otimes \pi_2(g)\cdot v_2-v_1\otimes v_2$, $v_1 \in V_1$, $v_2 \in V_2$, $g \in G$. Let
us denote it $(V_1\otimes V_2)_G$. By the equivalence of categories \ref{eqcat}, it is immediate that
$(V_1\otimes V_2)_G$ is identified with $V_1\otimes_{\caH(G)}V_2$, where the left $\caH(G)$-module
structure on $V_1$ given by $\pi_1$ is transformed into a right module structure by the anti-involution
$T\mapsto \check T$ of $\caH(G)$. Let $(\pi_i,V_i)$, $i=1,2,3$ be smooth representations of $G$. From the
natural isomorphism (this is a particular case of \ref{AD3})
\[
(V_1 \otimes V_2)\otimes V_3 \simeq  V_1 \otimes (V_2\otimes V_3),
\]
we deduce, by taking the coinvariants for the action of $G$ on each side,
\begin{align}\label{V1V2V3}
(V_1 \otimes V_2)\otimes_{\caH(G)} V_3 \simeq  V_1 \otimes _{\caH(G)} (V_2\otimes V_3).
\end{align}
This formula should be compared with that of the corollary above.

\begin{thm}
Let $(\pi,V)$ be a smooth representation of $G$.

$(i)$ The functor $\bullet \otimes V$ is exact, and preserves projective objects in $\caM(G)$.

$(ii)$ The (contravariant) functor $\Hom(\, \bullet\, ,V)_0$ is exact, and sends projective objects to
injective objects.

$(iii)$ The functor $\Hom(V,\, \bullet \, )_0$ is exact, and preserves injective objects.
\end{thm}

\begin{proof}
The proof of these three points is similar, so we will only prove $(iii)$. First, since the category of
complex vector spaces is semisimple, all objects in it are projective, injective and flat. The functor
$\Hom_\bbC(V,\, \bullet\, )$ from the category of complex vector spaces to itself is therefore exact. Let
\begin{equation} \label{suiexK}
   0 \xrightarrow{} W_1 \xrightarrow {\phi}
W_2  \xrightarrow {\psi}W_3 \xrightarrow {}   0
\end{equation}
be an exact sequence in $\caM(G)$. We obtain an exact sequence of complex vector spaces
\[
0 \longrightarrow \Hom_\bbC(V,W_1) {\longrightarrow} \Hom_\bbC(V,W_2)
{\longrightarrow }\Hom_\bbC(V,W_3)  \longrightarrow    0,
\]
these spaces being equipped with (non-smooth) representations of $G$. It is easy to see that taking the
smooth part preserves left exactness. We therefore have an exact sequence of smooth representations of
$G$:
\[
\{0\} \longrightarrow \Hom_\bbC(V,W_1)_0 {\longrightarrow} \Hom_\bbC(V,W_2)_0
{\longrightarrow}\Hom_\bbC(V,W_3)_0.
\]
We now establish the surjectivity of the last arrow. Let $f\in  \Hom_\bbC(V,W_3)_0$: there exists a
compact open subgroup $K$ of $G$ such that $k\cdot f=f$ for all $k\in K$, i.e.,
\[
k\cdot f(k^{-1}\cdot v)=f(v), \; (\forall k\in K, \, \forall v\in V).
\]
Consider the exact sequence (\ref{suiexK}) as an exact sequence in the category $\caM(K)$ of smooth
representations of the group $K$. This category is semisimple (we prove this below), and the sequence is
therefore split: there exists an intertwining operator (for the actions of $K$)
$s \colon W_3\rightarrow W_2$ such that $\psi\circ s=\Id_{W_3}$. Set $g=s\circ f  \in \Hom_\bbC(V, W_2)$.
We have
\[
(\forall k\in K, \, \forall v\in V), \quad (k\cdot g)(v) =
k \cdot(( s\circ f)(k^{-1}\cdot v))= s(k\cdot f(k^{-1}\cdot v))= s\circ f(v)=g(v).
\]
Thus $g\in  \Hom_\bbC(V, W_2)_0$ and satisfies $\psi\circ g=f$, which shows the assertion. We indeed have
an exact sequence of smooth representations of $G$
\[
\{0\} \longrightarrow \Hom_\bbC(V,W_1)_0 {\longrightarrow} \Hom_\bbC(V,W_2)_0
{\longrightarrow }\Hom_\bbC(V,W_3)_0\longrightarrow \{0\},
\]
which establishes the exactness of the functor $\Hom(V,\, \bullet\, )_0$.

Let us now show that this functor preserves injective objects. Suppose $W$ is injective in $\caM(G)$. We
need to show that $\Hom(V,W)_0$ is injective, i.e., that the functor
\[
\bullet\mapsto \Hom_G(\, \bullet\, , \Hom(V,W)_0)
\]
is exact. But according to the corollary above, we have a natural isomorphism
\[
\Hom_G(\, \bullet\,  ,\Hom(V,W)_0)\simeq \Hom_G(\bullet \otimes V,W) .
\]
But the functor $ \bullet \mapsto  \Hom_G(\, \bullet \, \otimes V,W) $ is exact as a composition of two
exact functors.

It remains to establish that $\caM(K)$ is a semisimple category. Let $(\pi,V)\in \caM(K)$, it suffices to
show that $V$ is a sum of simple submodules, according to Lemma \ref{AVII}. Let $v\in V$, and let $K_1$
be a normal compact open subgroup of $K$ fixing $v$. The subrepresentation $V_v$ generated by $v$ is
finite-dimensional because $K/K_1$ is finite, and completely reducible as a representation of the finite
group $K/K_1$, and therefore as a representation of $K$. Thus $v$ is contained in an irreducible
subrepresentation of $K$.
\end{proof}

\section{Induction and restriction}

\subsection{Forgetful functor and adjoints}\label{FOA}

We will apply the results of Section \ref{Oublietadjoints} to define functors between categories
$\caM(H)$ and $\caM(G)$, where $H$ and $G$ are t.d. topological groups. Let
$\phi \colon  H \rightarrow G$ be a continuous group morphism. We then have a forgetful functor
$\caF_\phi \colon \caM(G)\rightarrow \caM(H)$ defined as follows: if $(\pi,V)$ is a smooth representation
of $G$, then $\caF_\phi(V)=V$ and $\caF_\phi(\pi)$ is the representation of $H$ given by
\begin{equation}
\caF_\phi(\pi)(h)\cdot v =\pi(\phi(h))\cdot v, \quad (v\in V), \;
(h\in H). \label{oubl}
\end{equation}
We indeed obtain a smooth representation $V$ of $H$ because for all $v\in V$, $\mathrm{Stab}_G(v)$ is
open in $G$, $V$ being a smooth representation of $G$, and $\phi$ being continuous,
$\mathrm{Stab}_H (v)= \phi^{-1}(\mathrm{Stab}_G (v))$ is open in $H$.

Since we have the equivalences of categories $\caM(G) \simeq \caM(\caH(G))$ and
$\caM(H) \simeq \caM(\caH(H))$, we deduce by transport of structure a forgetful functor
\[
\caF_\phi  \colon \caM(\caH(G))\rightarrow \caM(\caH(H)).
\]
We must now identify it with a forgetful functor of the type defined in \ref{Oublietadjoints}, and to do
this, we must equip $\caH(G)$ with a non-degenerate left $\caH(H)$-module structure. We will go slightly
further by equipping $\caH(G)$ with a non-degenerate $\caH(H)$-bimodule structure such that the left
(resp. right) action of $\caH(H)$ on $\caH(G)$ commutes with the action by right (resp. left)
multiplication of $\caH(G)$ on itself. In such a situation, we defined in \ref{Oublietadjoints} not only
a forgetful functor, but also a pseudo-forgetful functor, as well as their respective right and left
adjoints.

Recall that $(\caH(G),l)$ and $(\caH(G),r)$ are two smooth representations of $G$, which via $\caF_\phi$
define two smooth representations of $H$, and therefore two non-degenerate left $\caH(H)$-module
structures on $\caH(G)$. We use the anti-involution induced by $h\mapsto h^{-1}$ on $\caH(H)$ to change
the left module structure given by $(\caH(G),r)$ into a right module structure. This therefore gives us a
non-degenerate $\caH(H)$-bimodule structure on $\caH(G)$. We can therefore, as in \ref{Oublietadjoints},
define the functors:
\[
\caF_H^G:=  \caF_{\caM(H)}^{\caM(G)}, \quad {}\spcheck \caF_H^G:=
{}\spcheck\caF_{\caM(H)}^{\caM(G)},\quad   I_H^G:=
I_{\caM(H)}^{\caM(G)},\quad   P_H^G:= P_{\caM(H)}^{\caM(G)}.
\]
It is clear that we then have $\caF_\phi= \caF_H^G$. On the other hand, recall that $ I_H^G$ is the right
adjoint of $\caF_H^G$ and that $P_H^G$ is the left adjoint of $ {}\spcheck \caF_H^G$.

\medskip

\begin{exemple}
Let $\sigma\colon H \rightarrow H_1$ be an isomorphism of t.d. groups. Let $(\pi,V)$ be a representation
of $H_1$. Applying the forgetful functor $\caF_\sigma$, we obtain a representation $\caF_\sigma (\pi,V)$
of $H$. We rather denote this representation by $\caF_\sigma (\pi,V)=({}^\sigma \pi,V)$. Its space is
$V$, and the action of $H$ is given by
\[
{}^{\sigma} \pi(h)\cdot v= \pi(\sigma(h))\cdot v.
\]
(The notation is mnemonic: $\sigma$ pulls $\pi$ to the left.)

This applies in particular when $\sigma$ is an automorphism of $H$. Thus, if $H$ is a subgroup of the
group $G$, and if $\sigma$ is the automorphism of $H$ induced by $\Int (g^{-1})$, for some $g\in G$
normalizing $H$, we denote by $(\pi^g,V)$ the representation $({}^\sigma \pi,V)$ of $H$
($\sigma=\Int (g^{-1}) \colon gHg^{-1} \rightarrow H$ and $g$ pushes $\pi$ from $\caM(H)$ to
$\caM(gHg^{-1})$, hence the notation $\pi^g$). In practice, the formula to use is therefore:
\begin{equation}\label{pig}
\pi^{g}(ghg^{-1})=\pi(h), \qquad (h\in H, g \in G).
\end{equation}
\end{exemple}

\subsection{Induction from a closed subgroup}\label{Ind}

Let us now place ourselves in the case where the morphism $\phi \colon H \rightarrow G$ is the inclusion
of a closed subgroup. In this case, a classical construction, the induction of representations from $H$
to $G$, provides a right adjoint to the forgetful functor. By uniqueness of the adjoint, this functor is
of course naturally isomorphic to $I_H^G$. A left adjoint to the pseudo-forgetful functor is constructed
in a similar way, by imposing a compactness condition on the support. Let us start with the definition of
induced representations; we will then show the adjunction properties.

\begin{defi}
Let $G$ be a t.d. group and $H$ a closed subgroup of $G$. Let $(\rho,E)$ be a smooth representation of
$H$. We define the induced representation
\[
\Ind_H^G (\rho,E)=(\Ind_H^G \rho, \Ind_H^G E)
\]
as follows. The space $\Ind_H^G E$ consists of functions $f \colon G\rightarrow E$ satisfying:

\medskip

a) $f(hg)=\rho(h)\cdot f(g),\quad (g\in G),\; (h\in H).$

b) There exists a compact open subgroup $K$ of $G$ (depending on $f$), such that $f(gk)=f(g)$, $g\in G$,
$k\in K$.

\medskip

This space is equipped with the representation $\Ind_H^G \rho$ defined by
\[
(\Ind_H^G \rho) (g)\cdot f =r(g)\cdot f
\]

The condition b) guarantees that $(\Ind_H^G \rho, \, \Ind_H^G E)$ is smooth.

We also define $\ind_H^G E$ as the subspace of $\Ind_H^G E$ of functions satisfying moreover the
following condition on their support:

\medskip

c) there exists a compact subset $F$ of $G$ such that $\supp f \subset H.F$

\medskip

It is clear that $\ind_H^G E$ is stable under the action $\Ind_H^G \rho$ of $G$. We denote by
$\ind_H^G \rho$ the restriction of $\Ind_H^G \rho$ to $\ind_H^G E$ and
$\ind_H^G (\rho,E)=(\ind_H^G \rho, \ind_H^G E)$.
\end{defi}

\begin{lemme}\label{Indlem}
Let $K$ be a compact open subgroup of $G$, and let $\Omega$ be a system of representatives for the double
cosets $H\backslash G/K$. For all $g\in \Omega$, set $K_g= H \cap gKg^{-1}$. The restriction of functions
from $\Ind_H^G E$ (resp. $\ind_H^G E$) to $\Omega$ defines an isomorphism of $(\Ind_H^G E)^K$ (resp.
$(\ind_H^G E)^K$) with the space $\caF(\Omega,E)$ (resp. $\caF_c(\Omega,E)$) of functions
$f\colon \Omega \rightarrow E$ such that $f(g)\in E ^{K_g}$ for all $g\in \Omega$ (resp. with the space
of such functions having finite support).
\end{lemme}

\begin{proof}
Since any function $f$ in $(\Ind_H^G E)^K$ satisfies:
\begin{align}\label{indO}
 \qquad f(hgk)=\rho(h)\cdot f(g),\quad (h\in H,\, g\in G,\, k\in
K),
\end{align}
it is clear that the restriction of $f$ to $\Omega$ determines $f$. On the other hand if
$h=gkg^{-1} \in K_g$, we have
\[
\rho(h)\cdot f(g)=f(hg)=f(gkg^{-1}g)=f(gk)=f(g),
\]
and therefore $f(g) \in E^{K_g}$. Conversely, any function $f$ on $\Omega$ with values in $E$ satisfying
$f(g)\in E^{K_g}$, can be lifted by formula (\ref{indO}) to a function in $(\Ind_H^G E)^K$. The lifting
is independent of the choices made, because if $hgk=h_1gk_1$, with $h,h_1\in H$ and $k,k_1\in K$, we have
$h_1^{-1}h=gk_1k^{-1}g^{-1} \in H \cap gKg^{-1}=K_g$, and therefore
\[
\rho(h)\cdot f(g)=  \rho(h_1h_1^{-1}h)\cdot
f(g)=\rho(h_1)\rho(h_1^{-1}h)\cdot f(g)= \rho(h_1)\cdot f(g).
\]
Since the space $G/K$ is discrete (cf. Lemma \ref{fini}), the same is true for $\Omega$, and a compactly
supported function on a discrete space is a finitely supported function.
\end{proof}

\begin{prop}
Induction from $H$ to $G$:
\[
(\rho,E)\mapsto \Ind_H^G  (\rho,  E), \quad (\text{resp. } \quad    \ind_H^G  (\rho,  E))
\]
defines an exact functor from $\caM(H)$ to $\caM(G)$.
\end{prop}

Let us describe, for once, the effect of this functor on morphisms: if $(\rho_1,E_1)$ and $(\rho_2,E_2)$
are two smooth representations of $H$, and if $\phi$ is an intertwining operator between $\rho_1$ and
$\rho_2$, the intertwining operator $\Ind_H^G(\phi)$ between $\Ind_H^G (\rho_1,E_1)$ and
$\Ind_H^G (\rho_2,E_2)$ is given by
\[
\Ind_H^G(\phi) (f)=\phi\circ f, \quad (f\in \Ind_H^G (E_1)).
\]

\begin{proof}
We have seen that the induced representation of a smooth representation is a smooth representation.
Suppose we are given an exact sequence
\[
0 \rightarrow E_1  \rightarrow E_2 \rightarrow E_3\rightarrow 0
\]
in $\caM(H)$. We need to show that the sequence
\[
0 \rightarrow  \Ind_H^G E_1  \rightarrow \Ind_H^G E_2 \rightarrow
\Ind_H^G E_3 \rightarrow 0
\]
is exact in $\caM(G)$. It suffices to show that for any compact open subgroup $K$ of $G$, the sequence
\[
0 \rightarrow  (\Ind_H^G E_1)^K  \rightarrow (\Ind_H^G E_2)^K \rightarrow
(\Ind_H^G E_3)^K \rightarrow 0
\]
is exact. According to the previous lemma $(\Ind_H^G E_i)^K$ is isomorphic to the space
\[
\caF(\Omega,E_i)\simeq \prod_{g\in \Omega} E_i^{K_g}.
\]
But the functor $E \mapsto E^{K_g}$ is exact (Proposition \ref{piK}). The proof is the same for
$\ind_H^G$.
\end{proof}

\subsection{Induction and admissibility} \label{IndAdm}

We keep the notation of the previous paragraph.

\begin{lemme}
Suppose that $H\backslash G$ is compact. Let $(\rho,E)$ be an admissible smooth representation of $H$.
Then
\[
\Ind_H^G  (\rho,  E)= \ind_H^G  (\rho,  E)
\]
is admissible.
\end{lemme}

\begin{proof}
This follows from the considerations of the previous paragraph, noting that here $\Omega=H\backslash G/K$
is a finite set.
\end{proof}

\subsection{Intertwining}\label{entrel}

Let $G$ be a t.d. group, and let $H$ be a closed subgroup of $G$.

\begin{lemme}
Suppose that $g\in G$ normalizes $H$. Then, with the notation of Example \ref{FOA},
$\Ind_H^G (\rho^g, E)$ is equivalent to $\Ind_H^G (\rho, E)$. Similarly $\ind_H^G (\rho^g, E)$ is
equivalent to $\ind_H^G (\rho, E)$.
\end{lemme}

\begin{proof}
Let us define the intertwining operator which realizes the equivalence of these representations: if
$f \in \Ind_H^G (\rho, E)$ (or $\ind_H^G (\rho, E)$),
\[
A\cdot f (g')= f(g^{-1}g').
\]
We then have:
\begin{align*}
 (A\cdot f) (hg')&= f(g^{-1}hg')= f((g^{-1}hg)g^{-1}g')=
  \rho(g^{-1} h g )\cdot (A \cdot f)(g')\\
&=\rho^g(h)\cdot A\cdot f(g'),
\end{align*}
and therefore $ A\cdot f$ is indeed in $\Ind_H^G (\rho^g, E)$ (or $\ind_H^G (\rho^g, E)$). It is clear
that $A$ is invertible (exchange the roles of $(\rho, E)$ and $(\rho^g, E)$, and replace $g$ by $g^{-1}$
to get the inverse).
\end{proof}

\subsection{Frobenius reciprocity} \label{Frob}
\index[ter]{Frobenius reciprocity}

Let $G$ be a t.d. group and let $H$ be a closed subgroup of $G$. Frobenius reciprocity is here the
assertion that the functor $\Ind_H^G$ is the right adjoint of the forgetful functor. We adopt here a more
classical terminology for this forgetful functor, denoting it by $\res_H^G$.

\begin{thm}
The functor $\Ind_H^G$ is the right adjoint of the functor $\res_H^G$. For all $(\pi,V)$ in $\caM(G)$ and
all $(\tau,E)$ in $\caM(H)$, we therefore have a natural isomorphism:
\[
\Hom_G(V,\Ind_H^G E) \simeq  \Hom_H(\res_H^G V,E).
\]
\end{thm}

\begin{proof}
Let $(\pi,V)$ in $\caM(G)$ and $(\tau,E)$ in $\caM(H)$. Let
\[
\alpha \colon \res_H^G(\Ind_H^G(\tau,E))\rightarrow (\tau,E),\quad \beta \colon (\pi,V) \rightarrow\Ind_H^G(\res_H^G (\pi,V))
\]
be the morphisms defined by
\[
\alpha \colon  f \mapsto f(\mathbf{1}_G), \; (f\in \Ind_H^G(\tau,E)), \quad \beta\colon v\mapsto (f_v \colon g \mapsto \pi(g)\cdot v)
\]
We immediately verify that $\alpha$ is a morphism in $\caM(H)$ and $\beta$ is a morphism in $\caM(G)$. We
check that they are adjunction morphisms, in the sense of \ref{adjfonct}. We first look at the
composition
\[
\Ind_H^G(\tau,E) \longrightarrow  \Ind_H^G \res_H^G  \Ind_H^G(\tau,E)
\rightarrow  \Ind_H^G(\tau,E),
\]
where the first arrow is the morphism $\beta$ for the space $\Ind_H^G(E)$. It sends $f \in \Ind_H^G(E)$
to
\[
f_f \colon  g \mapsto \Ind_H^G(\tau)(g)\cdot f =r(g)\cdot f.
\]
The second arrow is obtained by applying the functor $\Ind_H^G$ to the $H$-morphism $\alpha$. The effect
of the induction functor on morphisms is obtained by composition. When we evaluate at $r(g)\cdot f$, we
obtain $\alpha(r(g)\cdot f)=(r(g)\cdot f)(\mathbf{1}_G)=f(g)$. The composition of these two arrows is
therefore the identity of $\Ind_H^G(\tau,E)$. Similarly for the other composition
\[
\res_H^G(\pi,V) \rightarrow \res_H^G   \Ind_H^G \res_H^G (\pi,V) \rightarrow
\res_H^G (\pi,V),
\]
the first arrow is obtained by applying the functor $\res_H^G$ to the $G$-morphism $\beta$. It sends
$v\in V$ to the function $f_v\colon h \mapsto \pi(h)\cdot v$. The second arrow is the morphism $\alpha$
for the representation $\res_H^G(\pi,V)$. It sends the function $f_v$ to
$f_v(\mathbf{1}_G)=\pi( \mathbf{1}_G)\cdot v=v$. This composition is therefore the identity of
$\res_H^G (\pi,V)$. This completes the verification of the properties of the adjunction morphisms.
\end{proof}

\begin{rmq}
Let us identify the categories $\caM(G)$ and $\caM(\caH(G))$ (cf. Section \ref{eqcat}). By uniqueness of
the adjoint, the functor $I_H^G$ defined in \ref{FOA} and $\Ind_H^G$ are isomorphic and therefore for any
smooth representation $(\rho, W)$ of $H$, we have a natural isomorphism
\begin{equation} \label{IsoIInd}
\Ind_H^G(\rho,W) \simeq \Hom_{\caH(H)}(\caH(G),W)_{\caH(G)}.
\end{equation}
\end{rmq}

\begin{proof}
It is instructive to give the explicit form of this isomorphism.

Let $\phi \in \Hom_{\caH(H)}(\caH(G),W)_{\caH(G)}$. Let us define a function $f_\phi$ on $G$ by
\begin{equation}\label{IInd1}
f_\phi(g) = \phi(\delta_g).
\end{equation}
Since the distribution $\delta_g$ is not in $\caH(G)$, this definition calls for some explanation. Let us
choose an idempotent $e_K$ of $\caH(G)$, for some compact open subgroup $K$ of $G$, such that $\phi$ is
fixed by $e_K$. The distribution $\delta_g*e_K$ is then in $\caH(G)$, and we set
$ \phi(\delta_g)= \phi(\delta_g*e_K)$. We can immediately verify that this does not depend on the choice
of $K$. Moreover, for all $k \in K$, for all $g \in G$, we have
\[
f_\phi(gk)=  \phi(\delta_{gk})=  \phi(\delta_{g}*\delta_k*e_K)=
\phi(\delta_{g}*e_K)=f_\phi(g),
\]
which shows that $f_\phi$ is fixed by $r(K)$. We will show how $f_\phi$ transforms under left translation
by $h\in H$: for all $h \in H$, for all $g \in G$, we have
\begin{equation}\label{fphihg}
 f_\phi(hg)=   \phi(\delta_{hg})=\phi(\delta_{hg}*e_K)=\phi(\delta_h*\delta_g*e_K).
\end{equation}
On the other hand $ \rho(h)\cdot f_\phi(g)= (\delta_h *f)\cdot f_\phi(g)$ where $f$ is an idempotent of
$\caH(H)$ which fixes $\delta_g*e_K$. The idempotent $f$ then also fixes $f_\phi(g)= \phi(\delta_g*e_K)$.
We then have
\[
\rho(h)\cdot f_\phi(g)= (\delta_h *f)\cdot f_\phi(g)=
(\delta_h *f)\cdot \phi(\delta_g *e_K)= \phi(   (\delta_h *f)\cdot
\delta_g *e_K)),
\]
where $\delta_h *f \in \caH(H)$ acts on $\delta_g *e_K$ via the $\caH(H)$-module structure of $\caH(G)$.
More explicitly
\begin{align} \label{dhfdg}
 (\delta_h *f)\cdot (\delta_g *e_K)&= \int_H l(h')\cdot (\delta_g *e_K)\, d(\delta_h *f)(h')\\
  \nonumber &= \int_H \int_H  l(h_1h_2) \cdot (\delta_g *e_K)\, d(\delta_h)(h_1) \;  df(h_2)\\
\nonumber &  =    l(h)\cdot \left(  \int_H l(h_2) \cdot (\delta_g *e_K)  \, df(h_2) \right)\\
    \nonumber& =    l(h)\cdot (f\cdot (\delta_g *e_K))=   l(h)\cdot
    (\delta_g *e_K)=\delta_h*\delta_g*e_K.
\end{align}
We obtain by comparing (\ref{fphihg}) and (\ref{dhfdg})
\[
\rho(h)\cdot f_\phi(g)=  \phi(\delta_h*\delta_g*e_K)= f_\phi(hg).
\]
This shows that $f_\phi$ is in $\Ind_H^G W$.

Conversely, if $f \in \Ind_H^G W$, let us define $\phi_f \in \Hom_\bbC(\caH(G),W)$ by
\begin{equation}\label{IInd2}
  \phi_f(T)= (T\cdot f)(\mathbf{1}_G),  \quad (T\in \caH(G)),
\end{equation}
where $T\cdot f$ denotes the result of the action of $T\in \caH(G)$ on the element $f$ of the
$\caH(G)$-module $\Ind_H^G W$. If $f$ is fixed by $r(K)$, for some compact open subgroup $K$ of $G$, then
for all $T \in \caH(G)$,
\[
(e_K \cdot \phi_f)(T)=\phi_f(T*e_K)=((T*e_K)\cdot f)(\mathbf{1}_G)=
(T\cdot(e_K\cdot f))(\mathbf{1}_G)=\phi_f(T).
\]
If $S \in \caH(H)$, we have
\begin{align*}
 \phi_f(S\cdot T) &= ((S\cdot T)\cdot f) (\mathbf{1}_G)
=  \left(     \int_G  r(g) \cdot f \; d(S\cdot T)(g)      \right)  (\mathbf{1_G})\\
&= \left(    \int_G  r(g) \cdot f \; d \left(  \int_H l(h)\cdot T \; dS(h)  \right)
(g) \right) (\mathbf{1}_G) \\
& = \left(   \int_G   \int_H  l(h^{-1})\cdot( r(g) \cdot f) \; dT(g)\, dS(h)
\right)(\mathbf{1}_G) \\
&=    \int_H (l(h^{-1})\cdot (T\cdot f))(\mathbf{1}_G)\, dS(h)\\
&=  \int_H  (T\cdot f)(h) \, dS(h)   =  \int_H  \rho(h)\cdot ((T\cdot f) (\mathbf{1}_G) )\, dS(h)  \\
&=      S \cdot ((T\cdot f)(\mathbf{1}_G) )= S\cdot  \phi_f(T).
\end{align*}
This shows that $\phi_f$ is in $\Hom_{\caH(H)}(\caH(G),W)_{\caH(G)}$. One also easily verifies that
$f \mapsto \phi_f$ and $\phi \mapsto f_\phi$ are inverse to each other. The verification of the
naturality in $(\rho,W)$ of the isomorphisms $\Ind_H^G (W) \simeq \Hom_{\caH(H)}(\caH(G),W)_{\caH(G)}$
is, as often in this kind of thing, tedious but elementary. This shows that we indeed have a natural
isomorphism (\ref{IsoIInd}), realized by (\ref{IInd1}) and (\ref{IInd2}).
\end{proof}

\subsection{Compact induction and pseudo-forgetful functor}\label{indcomp}
\index[ter]{induction!compact}

In what follows, $H$ is a closed subgroup of the t.d. group $G$. The notation is the same as in Section
\ref{Ind}. The pseudo-forgetful functor $\check{}\, \caF_H^G$ has as its left adjoint the functor
$P_H^G$. We now want to identify this functor in terms of compact induction. Let then $(\rho,W)$ be a
smooth representation of $H$. Recall that $P_H^G(W)=\caH(G)\otimes_{\caH(H)} W$, that multiplication by
$\mu_G$ gives an algebra isomorphism $\scrD(G)\simeq \caH(G)$, and that similarly multiplication by
$\mu_H$ gives an algebra isomorphism $\scrD(H)\simeq \caH(H)$. By transport of structure, $\scrD(G)$ is
therefore equipped with a right $\scrD(H)$-module structure, which we will now identify explicitly. Let
$T=f\mu_G$, and $S=a\mu_H$ respectively in $\caH(G)$ and $\caH(H)$. We have by definition
\[
T\cdot S= \int_H (r(h)^{-1}\cdot T) \; dS(h)=\int_H (r(h^{-1})\cdot
(f\mu_G)) \; dS(h).
\]
Now, for any test function $\phi \in \scrC^\infty(G)$,
\begin{align*}
&\bil{r(h)^{-1}\cdot (f\mu_G) }{\phi}=\bil{(f\mu_G) }{r(h)\cdot\phi}
= \int_G f(g)\phi(gh)\, d\mu_G(g)\\
 =& \int_G f(gh^{-1})\phi(g) \delta_G(h)^{-1}  \, d\mu_G(g) =
 \bil{  \delta_G(h^{-1}) (r(h)^{-1}\cdot f)\mu_G }{\phi}.
\end{align*}
Therefore $r(h)^{-1}\cdot (f\mu_G)=  \delta_G(h)^{-1} (r(h)^{-1}\cdot f)\mu_G$. We deduce that
\begin{align*}
  T\cdot S&=\int_H  \delta_G(h)^{-1} (r(h)^{-1}\cdot f)\mu_G \,
dS(h)
= \left( \int_H  \delta_G(h)^{-1} (r(h)^{-1}\cdot f) a(h) \, d\mu_H(h)  \right)
\mu_G.
\end{align*}
The right action of $a \in \scrD(H)$ on $f \in \scrD(G)$ is therefore written
\begin{equation}\label{fpara}
f\cdot a =  \int_H \delta_G(h)^{-1} a(h)    (r(h)^{-1}\cdot f)  \, d\mu_H(h).
\end{equation}
With the right $\scrD(H)$-module structure thus defined on $\scrD(G)$, we then have
\[
P_H^G(W) \simeq \scrD(G)\otimes_{\scrD(H)} W.
\]
Let $f \in \scrD(G)$, and $v\in W$ and let us define a function $p(f\otimes v)$ from $G$ to $W$ by
\[
p(f\otimes v)(g)= \int_H \check f(hg) ((\delta_G\rho)(h^{-1})\cdot v)
\, d\mu_H(h).
\]

Let us first show that this is well-defined for any element $f \otimes v$ of
$\scrD(G)\otimes_{\scrD(H)} W$. We need to see that if $a \in \scrD(H)$, we indeed have
\[
p(f\cdot a\otimes v)=p(f\otimes (a\mu_H) \cdot v).
\]
Now
\begin{align*}
& p(f\cdot a \otimes v)(g)= \int_H \check {(f\cdot  a)}(hg)
((\delta_G \rho)(h^{-1}) \cdot v) \, d\mu_H(h) \\
=& \int_H \left[\int_H \delta_G(h')^{-1} a(h')    (r(h')^{-1}\cdot  f)((hg)^{-1})
\, d\mu_H(h')\right]
 ((\delta_G \rho)(h^{-1})\cdot v) \, d\mu_H(h) \\
=& \int_H \int_H \delta_G(h'h)^{-1} \check f (h'hg) a(h')  (\rho(h^{-1})\cdot v)
\, \, d\mu_H(h') d\mu_H(h) \\
=& \int_H \delta_G(h)^{-1} \check f (hg)  \rho(h)^{-1} \cdot
\left(\int_H a(h')  (\rho(h')\cdot v) \,  d\mu_H(h') \right) d\mu_H(h) \\
=& p(f \otimes (a\mu_H) \cdot v)(g).
\end{align*}

A similar calculation shows that if $f$ is fixed by a certain compact open subgroup $K$ of $G$ for the
action by left translation, then $p(f \otimes v)$ is also fixed by $K$ (for the action by right
translation).

We now calculate how $p(f \otimes v)$ transforms under left translation by an element of $H$.
\begin{align*}
 & p(f \otimes v)(h_0g)=  \int_H \check f(hh_0g) ((\delta_G \rho)(h^{-1})\cdot v)
 \, d\mu_H(h)\\
 &= \delta_G(h_0) \rho(h_0)\cdot    \int_H \check f(hh_0g)
 ((\delta_G \rho)((hh_0)^{-1})\cdot v) \,    d\mu_H(h)\\
&= \delta_H(h_0)^{-1} \delta_{G}(h_0) \rho(h_0) \cdot \left( \int_H \check
f(hg)((\delta_G \rho)(h^{-1})\cdot v) \,  d\mu_H(h) \right)\\
&= (\delta_{H\backslash G} \rho)(h_0) \cdot \left( \int_H \check
f(hg)((\delta_G \rho)(h^{-1})\cdot v) \,   d\mu_H(h)\right)\\
&=  (\delta_{H\backslash G} \rho)(h_0)\;   p(f \otimes v)(g).
\end{align*}

On the other hand, since $f$ is compactly supported, $ p(f \otimes v)$ is compactly supported modulo $H$.
This shows that $p(f \otimes v)$ is a function in $\ind_H^G (W\otimes  \delta_{H\backslash G})$. By
linearity, we have thus defined a linear map:
\begin{equation}\label{IsoPind}
 p\colon P_H^G (W) \rightarrow \ind_H^G  (W\otimes  \delta_{H\backslash G}).
\end{equation}
A calculation similar to the one performed above shows that $p$ intertwines the actions of $G$ on
$P_H^G (W)$ and $\ind_H^G  (W\otimes  \delta_{H\backslash G})$.

\begin{thm}
The morphism $p\colon P_H^G (W) \rightarrow \ind_H^G (W\otimes  \delta_{H\backslash G})$ is a (natural)
isomorphism in the category $\caM(G)$.
\end{thm}

\begin{proof}
Let us fix a compact open subgroup $K$ of $G$. It suffices to show that $p$ induces an isomorphism
\[
(\scrD(G)\otimes_{\scrD(H)} W)^K \simeq ( \ind_H^G (W\otimes
\delta_{H\backslash G}))^K.
\]
We saw in Lemma \ref{Indlem} that this latter space is identified with the space of finitely supported
functions on a system of representatives $\Omega$ of $H\backslash G/K$, satisfying $f(g) \in W^{K_g}$ for
all $g$ in $\Omega$, i.e.,
\begin{equation}\label{indcompiso1}
( \ind_H^G (W\otimes  \delta_{H\backslash G}))^K \simeq \bigoplus_{g\in
  \Omega}  W^{K_g}.
\end{equation}
Let us now study the space
\[
(\scrD(G)\otimes_{\scrD(H)} W)^K= (e_K*\scrD(G)) \otimes_{\scrD(H)} W.
\]
Note that through the equivalences of categories
\[
\caM(H)\simeq \caM(\caH(H)) \simeq \caM(\scrD(H)),
\]
$\scrD(G)\otimes_{\scrD(H)} W$ is identified with $\scrD(G)\otimes_{H} W$ where this latter space is by
definition the quotient of $\scrD(G)\otimes W$ by the space generated by the vectors of the form
$r(h^{-1})\cdot f\otimes w -f \otimes \rho(h)\cdot w$. On the other hand, for all $g \in \Omega$, $KgH$
is an open set of $G$, so we have a decomposition of $G$ into disjoint open sets
\[
G=\coprod_{g^{-1}\in \Omega} KgH.
\]
This induces a decomposition of $e_K*\scrD(G)$ into
\[
e_K*\scrD(G)=\bigoplus_{g^{-1}\in \Omega} \scrD(G)^K_g
\]
where $\scrD(G)^K_g$ is the space of functions $f$ in $\scrD(G)$ invariant under left translation by $K$
and supported in $KgH$. Note that since $\Omega$ is a system of representatives in $G$ for the double
cosets $H\backslash G/K$, the set $\{g^{-1}\}_{g\in \Omega}$ is a system of representatives for
$K\backslash G/H$. Since each $\scrD(G)^K_g$ is stable under the action of $H$ by right translation,
\[
(e_K*\scrD(G))\otimes_{H} W = \bigoplus_{g^{-1}\in \Omega} \scrD(G)^K_g \otimes_{H} W.
\]
The attentive reader will have long since noticed the resemblance between the proof we are completing and
that of Theorem \ref{invmesquo}. In particular, the definition of the function $p$ is similar in both
cases, and we also introduce (with a slightly different convention, but we can restore things by passing
to the inverse in the group $G$) the spaces $\scrD(G)^K_g$. We had seen that $\scrD(G)^K_g$ is isomorphic
to the space of finitely supported functions on the discrete space $K\backslash KgH$, a space equipped
with a transitive action of $H$ by right translation. We conclude that
\[
\scrD(G)^K_g \otimes_{H} W  \simeq \bbC \otimes _{\mathrm{Stab}_H (K\backslash Kg)}W
\simeq W^{K_{g^{-1}}} .
\]
and therefore
\begin{equation} \label{indcompiso2}
 P_H^G(W)^K\simeq   \bigoplus_{g^{-1}\in \Omega} W^{K_{g^{-1}}}=
 \bigoplus_{g \in \Omega} W^{K_g}.
\end{equation}
The isomorphism $ P_H^G(W)^K \simeq  (\ind_H^G (W\otimes \delta_{H\backslash G}))^K$ is deduced from
(\ref{indcompiso1}) and (\ref{indcompiso2}): indeed, $p$ respects these decompositions along the double
cosets. If $f\in\scrD(G)^K_g$, then $p(f\otimes v)$ is supported in $Hg^{-1}K$, and for $f=\chi_{Kg}$,
the function $p(f\otimes v)$ takes at $g^{-1}$ a non-zero multiple of $\rho(e_{K_{g^{-1}}})\cdot v$. Thus
$p$ induces, factor by factor, the above isomorphisms.
\end{proof}

\medskip

\begin{rmq}
We end this section by studying the case of induction from an open subgroup $H$ of $G$ (the topology
being totally disconnected, $H$ is also closed in $G$). In this case, we have an obvious inclusion
$\caH(H) \hookrightarrow \caH(G)$, where we identify $\caH(H)$ with the space of distributions in
$\caH(G)$ supported in $H$. The non-degenerate $\caH(H)$-bimodule structure of $\caH(G)$ described
previously and using the representations $l$ and $r$ of $H$ on $\caH(G)$ can be described here more
simply as the bimodule structure given by the algebra inclusion of $\caH(H)$ into $\caH(G)$. In this
framework, we noted at the end of Section \ref{Oublietadjoints} that the forgetful functor and the
pseudo-forgetful functor are isomorphic. What interests us here is the fact that consequently the
forgetful functor $\res_H^G$ admits a left adjoint, in addition to the right adjoint $\Ind_H^G$, and that
this is given by $\ind_H^G$ (the factor $\delta_{H\backslash G}$ is trivial in this case, since $\mu_H$
is, up to the choice of a multiplicative factor, the restriction of $\mu_G$ to $H$). We therefore have,
for any representation $(\pi,V)$ of $\caM(G)$, and any representation $(\rho,W)$ of $\caM(H)$, a natural
isomorphism
\begin{equation}\label{indres}
 \Hom_G(\ind_H^G(\rho,W),(\pi,V))\simeq
  \Hom_H((\rho,W),\res_H^G(\pi,V)).
\end{equation}
\end{rmq}

\subsection{Induction and duality} \label{IndandDual}

We can combine the results of the previous section and Theorem \ref{duaetoub} to obtain the following
statement:

\begin{thm}\label{inddualthm}
Let $H$ be a closed subgroup of the t.d. group $G$, and let $(\rho,W)$ be a smooth representation of $H$.
We then have a natural isomorphism
\[
\Ind_H^G(\widetilde W)\simeq \ind_H^G(W\otimes  \delta_{H\backslash G})^\sim.
\]
The duality between $\Ind_H^G(\widetilde W)$ and $\ind_H^G(W\otimes  \delta_{H\backslash  G})$ is given
explicitly by the formula
\[
\bil{f}{h}=\int_{H\backslash G}  h(g)(f(g))\, d\nu_{H\backslash G}, \qquad
(f \in \ind_H^G(W\otimes  \delta_{H\backslash  G})),  \; (h\in  \Ind_H^G(\widetilde W)).
\]
\end{thm}

\begin{proof}
It remains to prove the duality formula. One easily verifies that the integral is well-defined, and it
gives a $G$-equivariant morphism from $\Ind_H^G(\widetilde W)$ to
$\ind_H^G(W\otimes  \delta_{H\backslash  G})^\sim$. We check that this morphism is indeed an isomorphism
by Lemmas \ref{Indlem} and \ref{dualitelem}.
\end{proof}

\subsection{Induction in stages}

Let $J\subset H$ be two closed subgroups of $G$. Since induction from $H$ to $G$ is the right adjoint of
restriction from $G$ to $H$, and since trivially:
\[
\res_J^G=\res_J^H \circ \res_H^G,
\]
we obtain the

\begin{prop}
We have an isomorphism of functors:
\[
\Ind_J^G \simeq \Ind_H^ G \circ \Ind_J^H.
\]
\end{prop}

\begin{proof}
This follows from the uniqueness of the adjoint (\ref{adjfonct}).
\end{proof}

\subsection{Coinvariants }\label{FoncJac1}
\index[ter]{coinvariants}

Let $H$ be a t.d. group and $(\tau,E)$ a representation of $H$. We denote by
\[
E(H,\tau)
\]
the subspace of $E$ generated by the vectors $(\tau(h)\cdot v-v)$, $v\in E$, $h\in H$. The coinvariant
space is by definition the quotient:
\[
E_{H,\tau}= E/E(H,\tau).
\]
The dual $(E_{H,\tau})^*$ of $E_{H,\tau}$ is identified with the subspace of $E^*$ of linear forms
vanishing on $E(H,\tau)$, i.e., such that $\tau^{*}(h)\cdot (\lambda)=\lambda$, for all $h\in H$.

\begin{exemple}
The theorem of existence and uniqueness of the Haar measure on a t.d. group $G$ can be reformulated by
saying that
\[
\dim (\scrD(G)_{G,l})=1.
\]
\end{exemple}

\begin{rmqs}
$ $
\begin{itemize}
\item[1.] If $H$ is a closed subgroup of a t.d. group $G$, and if $(\pi,V)$ is a smooth representation of
$G$, we simplify the notation for $V(H,\res_H^G(\pi))$ (resp. $V_{H,\res_H^G(\pi)}$) to $V(H,\pi)$ (resp.
$V_{H,\pi}$), or even simply $V(H)$ (resp. $V_H$) if the representation $\pi$ is clearly identified by
the context. In this case, any subgroup $N$ of $G$ normalizing $H$ acts on $V(H,\pi)$ and therefore on
$V_{H,\pi}$. We will denote by $\pi_{H}$ the action of $N$ on $V_{H,\pi}$.

\item[2.] If $H$ and $N$ are closed subgroups of a t.d. group $G$, and if $N$ normalizes $H$, then
\[
(\pi,V)\mapsto (\pi_{H}, V_{H,\pi})
\]
defines a functor $j_H$ from $\caM(G)$ to $\caM(N)$. Indeed, any intertwining operator $A$ between two
smooth representations $(\pi,V)$ and $(\rho,W)$ of $G$ induces an intertwining operator between
$(\pi_{H}, V_{H,\pi})$ and $(\rho_{H}, W_{H,\rho})$.
\end{itemize}
\end{rmqs}

\begin{lemme}
$(i)$ Let $H$ be a t.d. group and suppose that $H_1$ and $H_2$ are closed subgroups of $H$ such that
$H=H_1H_2$, where $H_1$ normalizes $H_2$. Let $(\tau,E)$ be a representation of $H$. We then have
\[
(E_{H_2,\tau})_{H_1, \tau_{H_2}}= E_{H,\tau},
\]
where $\tau_{H_2}$ denotes the representation of $H_1$ on $E_{H_2,\tau}$ obtained from $\tau$ by passing
to the quotient (cf. Remark 1 above).

$(ii)$ Let $H,H_1,H_2$ be as above, $H$ being a closed subgroup of a t.d. group $G$, and let $N,N_2$ be
closed subgroups of $G$ such that $N\subset N_2$, $N$ normalizes $H$ and $H_1$, $N_2$ normalizes $H_2$.
Let $J$ be the subgroup of $G$ generated by $H_1$ and $N_2$. Then the functor
\[
j_H \colon \caM(G)\rightarrow \caM(N)
\]
is isomorphic to the composition of the functors
\[
j_{H_2} \colon  \caM(G) \rightarrow \caM(J) \text { and }  j_{H_1} \colon  \caM(J)\rightarrow \caM(N).
\]
\end{lemme}

\begin{proof}
For all $h_1 \in H_1$, for all $h_2\in H_2$ and for all $v\in E$, we have:
\[
\tau(h_1h_2)\cdot v-v= (\tau(h_1)\cdot (\tau(h_2)\cdot v)-(\tau(h_2)\cdot v))
+ ( \tau(h_2)\cdot v -v).
\]
This shows that $E(H,\tau)=E(H_1,\tau)+E(H_2,\tau)$. Point $(i)$ follows immediately, and $(ii)$ is a
consequence of $(i)$.
\end{proof}

\begin{defi}
Let $X$ be a set and $(X_i)_{i\in I}$ a family of subsets of $X$. We say that $X$ is the increasing
directed union of the $X_i$ if $X=\bigcup_{i\in I}X_i$ and for any finite subset $J\subset I$, there
exists $k\in I$ such that $\bigcup_{i\in J}X_i\subset X_k$.
\end{defi}

\begin{rmq}
Suppose that $H$ is a t.d. group and that $H$ is the increasing directed union of its compact subgroups.
It is then clear that for any compact subset $C$ of $H$, there exists a compact subgroup $K$ of $H$ which
contains $C$. (Cover $C$ by a finite family of open sets of the form $g_jK_j$, $j\in J$, where $g_j\in H$
and $K_j$ is a compact open subgroup of $H$ and take $K$ containing all the $g_j$ and all the $K_j$).
\end{rmq}

\needspace{156pt}
\begin{prop}
Suppose that $H$ is the increasing directed union of its compact subgroups. Let $(\tau,E)$ be a smooth
representation of $H$. Then $E(H,\tau)$ is the space of vectors $v\in E$ such that there exists a compact
subgroup $K$ of $H$ satisfying
\[
\tau(e_K)\cdot v=0.
\]
\end{prop}

\begin{proof}
The proposition follows from Proposition \ref{piK}, which treats the case where $H$ is itself compact.
\end{proof}

\begin{cor}
With the hypotheses of the above proposition, if $E'$ is a subrepresentation of $E$, then
$E'(H,\tau)=E(H,\tau)\cap E'$ and if
\[
0 \rightarrow  E' \rightarrow E \rightarrow E'' \rightarrow 0
\]
is an exact sequence in $\caM(H)$, then
\[
0 \rightarrow  E'_{H,\tau'} \rightarrow E_{H,\tau} \rightarrow
E''_{H,\tau''} \rightarrow 0
\]
is exact. In particular, the functor $j_H$ defined in Remark 2 above is exact as soon as $H$ is the
increasing directed union of its compact subgroups.
\end{cor}

\begin{proof}
The first assertion is a consequence of the above proposition. The only delicate point to verify then is
that $E'_{H,\tau'} \rightarrow E_{H,\tau}$ is an injection, which follows from the first point.
\end{proof}

\subsection{A particular case of the forgetful functor} \label{Jacmod}

Suppose now that $P$ is a t.d. group and that $N$ is a closed normal subgroup of $P$. Let us set $M=P/N$.
It is a t.d. group (Lemma \ref{topquo}). Let $\phi: P\rightarrow M$ denote the natural projection. It
induces a forgetful functor $\caF_P^M \colon \caM(M) \rightarrow \caM(P)$ and a pseudo-forgetful functor
$\check{}\, \caF_P^M$ (\ref{FOA}). On the other hand, we define
$\phi^*  \colon \scrD(M)\rightarrow \scrC^\infty(P), \; \psi\mapsto \psi\circ \phi$, and by duality we
define $\phi_*  \colon  \scrE'(P)\rightarrow \scrD'(M)$. We therefore have
$\bil{\phi_* (T)}{\psi}= \bil{T} {\psi\circ \phi}$ for $\psi\in \scrD(M)$.

We show that in this context, these two functors are naturally isomorphic. Let $K$ be a compact open
subgroup of $P$, and let us denote $K_M=\phi(K)$. It is a compact open subgroup of $M$, and the family of
$K_M$, as $K$ runs through the family of compact open subgroups of $P$, forms a neighborhood basis of the
identity of $M$, by definition of the quotient topology. Consequently, the family of idempotents
$e_{K_M}$ forms a directed system of idempotents for $\caH(M)$.

The Hecke algebra $\caH(M)$ is equipped with a non-degenerate $\caH(P)$-bimodule structure. We will show
that the action (say on the right, but the proof would be similar for the left action) of an element
$a_{t,K}= e_K*\delta_t*e_K$ of $\caH(P)$, where $t$ is an element of $P$, on a distribution $T$ of
$\caH(M)$ is given by
\begin{equation}\label{HPHM}
T\cdot a_{t,K}=T*a_{\phi(t),K_M}
\end{equation}
where $a_{\phi(t),K_M}=e_{K_M}*\delta_{\phi(t)}*e_{K_M}$. Indeed, for any test function $\psi$,
\begin{align*}
\bil{T\cdot a_{t,K}}{\psi}&= \int_M \psi(m) \, d(T\cdot a_{t,K})(m)
 = \int_M \psi(m)  d\left( \int_P  r(\phi(p^{-1}))\cdot T  \; d(a_{t,K})(p)
 \right)(m)\\
&=  \int_M  \int_P \psi(m\phi(p) ) \, dT(m) \,  d(a_{t,K}) (p)
=  \int_M  \int_{M}   \psi(mm') \, dT(m)  d(\phi_*(a_{t,K}))(m')\\
&=  \int_M  \int_{M}   \psi(mm') \, dT(m)  d(a_{\phi(t),K_M})(m')
=\bil{T*a_{\phi(t),K_M}}{\psi}.
\end{align*}
The equality $\phi_*(a_{t,K})= a_{\phi(t),K_M} $ follows from the characterization of $a_{\phi(t),K_M}$
as the unique distribution supported in $K_M\phi(t)K_M$, invariant under the action of $K_M$ by left and
right translation, and normalized (Lemma \ref{Gammalem}). Indeed, we verify the right translation
invariance of $\phi_*(a_{t,K})$ under $K_M$. Let $k\in K_M$ and let $\tilde k \in K$ such that
$\phi(\tilde k)=k$. Then, for any test function $\psi \in \scrD(M)$,
\begin{align*}
\bil{r(k)\cdot \phi_*(a_{t,K}) }{\psi}&= \bil{\phi_*(a_{t,K}) }{r(k^{-1})\cdot
\psi}= \bil{a_{t,K}}{ (r(k^{-1})\cdot  \psi)\circ  \phi}
=  \bil{a_{t,K}}{ r(\tilde k^{-1})\cdot(  \psi\circ  \phi)}\\
&= \bil{ r(\tilde k)\cdot a_{t,K}}{  \psi\circ  \phi} = \bil{ a_{t,K}}{  \psi\circ  \phi}\\
&= \bil{\phi_*(a_{t,K}) }{\psi}.
\end{align*}
Left translation invariance and normalization are verified in the same way.

A particular case of (\ref{HPHM}) is
\[
T\cdot e_K= T*e_{K_M}.
\]
Therefore, for any finite family of elements $T_i$ of $\caH(M)$, there exists a compact open subgroup $K$
such that $T_i= T_i * e_{K_M}= T_i\cdot e_K$. From this, it follows that for any module $V$ in
$\caM(\caH(M))$,
\begin{align*}
  \check{}\, \caF_P^M(V)&=\Hom_{\caH(M)}(\caH(M), V)_{\caH(P)}
=\Hom_{\caH(M)}(\caH(M), V)_{\caH(M)} \simeq V \simeq  \caF_P^M(V).
\end{align*}
All these isomorphisms are of course natural in $\caM(P)$.

The pseudo-forgetful functor $\check{}\, \caF_P^M$ admits a left adjoint $P_P^M$. Now let us describe it
in the language of representations. If $(\pi,V)$ is a smooth representation of $P$, according to Remark 2
of the previous paragraph, $(\pi_N,V_N)$ is a smooth representation of $P$, whose kernel contains $N$. We
can therefore view $(\pi_N,V_N)$ as a representation of $M$, and we will do so without changing notation.
This defines a functor $j_N$ from $\caM(P)$ to $\caM(M)$.

\begin{prop}
The functors $j_N$ and $P_P^M$ are naturally isomorphic. In particular $j_N$ is the left adjoint of the
forgetful functor.
\end{prop}

\begin{proof}
Let $(\pi,V)$ be a smooth representation of $P$. For all $v$ in $V$, let us choose a compact open
subgroup $K$ of $P$ fixing $v$, and define
\begin{equation}\label{VHMHPV}
\Phi \colon  V\rightarrow \caH(M)\otimes_{\caH(P)} V, \quad v \mapsto e_{K_M}\otimes v.
\end{equation}
First of all, note that $e_{K_M}\otimes v$ is in fact independent of the choice of $K$. Indeed, if
$K' \subset K$ is another compact open subgroup, we have
\[
e_{K'_M} \otimes v= e_{K'_M}\otimes e_K\cdot v= e_{K'_M}\cdot e_{K}\otimes v=  e_{K'_M}*e_{K_M}\otimes v= e_{K_M} \otimes v.
\]
This shows in particular that the map (\ref{VHMHPV}) is linear. Let us now show that $V(N)$ is in
$\ker \Phi$. Let $v\in  V$. We have for all $n \in N$,
\[
\Phi(\pi(n)\cdot v-v)= e_{K_M}\otimes (\pi(n)\cdot v-v),
\]
where $K$ is chosen such that $e_K\cdot v=v$ and $e_K\cdot (\pi(n)\cdot v)= \pi(n)\cdot v$. Then
\begin{align*}
\Phi(\pi(n)\cdot v-v)&=e_{K_M}\otimes e_K*\delta_n*e_K\cdot v- e_{K_M}\otimes  v
=e_{K_M}\otimes a_{n,K}\cdot v- e_{K_M}\otimes  v\\
 &=e_{K_M}\cdot  a_{n,K}\otimes  v - e_{K_M}\otimes  v=e_{K_M} *a_{\phi(n),K_M}  \otimes  v - e_{K_M}\otimes  v\\
& =e_{K_M} *e_{K_M}\otimes  v - e_{K_M}\otimes  v\\
&=0.
\end{align*}
We used the fact that $\phi(n)=\mathbf{1}_M$ and therefore $a_{\phi(n),K_M}=e_{K_M}$. By passing to the
quotient, $\Phi$ induces a linear map
\[
\overline{ \Phi}\colon V_N \rightarrow \caH(M)\otimes_{\caH(P)}V.
\]
This linear map is an intertwining operator for the actions of $M$ on $V_N$ and
$\caH(M)\otimes_{\caH(P)}V$. Indeed, for all $m\in M$ and for all $\bar v\in V_N$, if we choose $p$ in
$P$ such that $\phi(p)=m$, $v\in V$ which lifts $\bar v$, and $K$ a compact open subgroup of $P$ such
that $e_K\cdot v=v$ and $e_K\cdot(\pi(p)\cdot v)=\pi(p)\cdot v$, then
\begin{align*}
\overline{ \Phi} (\pi_N(m)\cdot \bar v)&= \overline{ \Phi }(\overline{\pi(p)\cdot  v})
=e_{K_M}\otimes \pi(p)\cdot v
= e_{K_M}\otimes e_K*\delta_p*e_K\cdot v= e_{K_M}\otimes a_{p,K}\cdot v\\
&=e_{K_M}\cdot a_{p,K}\otimes v= e_{K_M}*a_{m,K_M}\otimes v
=   e_{K_M}*\delta_m*e_{K_M}\otimes v\\
& = e_{K_M}\cdot(l(m)\cdot \overline{\Phi}(\bar v)).
\end{align*}
By replacing $K$ by a smaller compact open subgroup, we can assume that $e_{K_M}$ fixes
$l(m)\cdot \overline{ \Phi} (\bar v)$.

We must now construct an inverse for $\overline{ \Phi}$. Since $V_N$ is a smooth representation of $M$,
it is a non-degenerate $\caH(M)$-module. We then define
\[
\Psi\colon \caH(M)\otimes V \rightarrow V_N ,\quad T\otimes v
\mapsto T\cdot \bar v.
\]
We will show that $\Psi$ induces a map
\[
\Psi\colon \caH(M)\otimes_{\caH(P)} V \rightarrow V_N .
\]
To do this, it suffices to verify that for all $T \in \caH(M)$, for all $S\in \caH(P)$ and for all
$v\in V$,
\[
\Psi (T\cdot S \otimes v- T\otimes S\cdot v)=0.
\]
Now $\Psi (T\cdot S \otimes v- T\otimes S\cdot v)=(T\cdot S)\cdot \bar v- T\cdot (S\cdot \bar v)$. For a
sufficiently small compact open subgroup $K$ of $P$, the distribution $S$ is in $\caH(P,K)$ and is
therefore a linear combination of distributions of the form $a_{p,K}$. It therefore suffices to verify
that
\[
(T\cdot a_{p,K})\cdot \bar v- T\cdot (a_{p,K}\cdot \bar v)=0.
\]
Since $T\cdot a_{p,K}= T * a_{\phi(p),K_M}$, it suffices to verify that
$a_{p,K}\cdot \bar v= a_{\phi(p),K_M}\cdot \bar v$. But this follows from the definitions. We complete
the proof by verifying that $\Psi$ is the inverse of $\overline{ \Phi}$, which follows from a
straightforward calculation.
\end{proof}

\begin{cor}
Suppose that $P$ is a closed subgroup of the t.d. group $G$ and that $N$, $M$ are as above. Then the
functor $\Ind_P^G \circ  \caF_P^M$ from $\caM(M)$ to $\caM(G)$ is the right adjoint of the functor
$j_N\circ\res_P^G \colon V \mapsto V_N$.
\end{cor}

\begin{proof}
This is clear by composition of adjoint functors.
\end{proof}

\subsection{Mackey isomorphisms}\label{Mackey}

The Mackey isomorphism\index[ter]{isomorphism!Mackey} establishes the commutation relations between the
$\Hom$ (resp. $\otimes$) functor of Section \ref{OHom} and the induction functor $I_H^G$ (resp. $P_H^G$)
of Section \ref{FOA}. Let us state the result.

\begin{thm}
Let $\phi \colon H \rightarrow G$ be a continuous morphism of t.d. groups by means of which we construct
the induction functors $I_H^G$ and $P_H^G$ (\ref{FOA}). Then, for any representation $(\rho,W)$ in
$\caM(H)$ and any representation $(\pi,V)$ in $\caM(G)$, we have a natural isomorphism
\[
\Hom(V,I_H^G(W))_0 \simeq I_H^G(\Hom( \caF_H^G(V),W)_0).
\]
Moreover, if the pseudo-forgetful functor $\check{}\, \caF_H^G$ is isomorphic to the forgetful functor
$\caF_H^G$, then
\[
P_H^G(W)\otimes V \simeq    P_H^G(W \otimes \caF_H^G(V) ).
\]
\end{thm}

\begin{proof}
For any representation $(\sigma,E)$ in $\caM(G)$, we have
\begin{align}
\label{x1} \Hom_G(E, \Hom(V,I_H^G(W))_0)&\simeq \Hom_G(E\otimes V,
I_H^G(W))\\
 \label{x2}&\simeq  \Hom_H(\caF_H^G(E\otimes V), W)\\
 \label{x3}&\simeq  \Hom_H(\caF_H^G(E)\otimes \caF_H^G(V), W)\\
\label{x4}&\simeq  \Hom_H(\caF_H^G(E), \Hom(\caF_H^G(V), W)_0)\\
\label{x5}&\simeq  \Hom_G(E, I_H^G(\Hom(\caF_H^G(V), W)_0)).
\end{align}
Indeed (\ref{x1}) is the application of Corollary \ref{OHomcor}, (\ref{x2}) is the adjunction of the
functors $\caF_H^G$ and $I_H^G$, (\ref{x3}) is trivial, (\ref{x4}) is again Corollary \ref{OHomcor}, and
(\ref{x5}) is again the adjunction of $\caF_H^G$ and $I_H^G$.

A general argument from category theory (cf. \ref{egalitecat}) then allows us to conclude that we indeed
have a natural isomorphism
\[
\Hom(V,I_H^G(W))_0 \simeq I_H^G(\Hom( \caF_H^G(V),W)_0).
\]

Let us now move on to the second Mackey isomorphism, and try to copy the argument used for the first.
\begin{align}
\label{y1} \Hom_G(V\otimes P_H^G(W),E)&\simeq \Hom_G(P_H^G(W),\Hom(V,E)_0)\\
 \label{y2}&\simeq  \Hom_H(W,\check{}\, \caF_H^G(\Hom (V,E)_0) )\\
 \label{y3}&\simeq  \Hom_H(W,\Hom (\caF_H^G(V), \check{}\, \caF_H^G( E))_0) \\
\label{y4}&\simeq  \Hom_H(\caF_H^G(V)\otimes W, \check{}\, \caF_H^G( E) )\\
\label{y5}&\simeq  \Hom_G(P_H^G( \caF_H^G(V)\otimes W), E).
\end{align}
We used in this calculation Corollary \ref{OHomcor} and the adjunction between the functors $P_H^G$ and
$ \check{}\, \caF_H^G$. The only unjustified step is equality (\ref{y3}), which unlike (\ref{x3}) above,
is not trivial in general, but obviously becomes so if $\check{}\, \caF_H^G$ is isomorphic to $\caF_H^G$.
\end{proof}
\section{Representations in sections of sheaves}

\subsection{Action on a sheaf}\label{actionfaisceau}

\begin{defi}
Let $X$ be a t.d. space, $G$ a t.d. group, and $\caF \in \scrC^\infty_X-\caM od(X)$ a sheaf of
$\bbC$-vector spaces on $X$. An action of $G$ on $\caF$ consists, for all $g\in G$, of an isomorphism
$\gamma_g$ from $\caF$ to itself (over a homeomorphism of $X$, still denoted $\gamma_g$; see \ref{FCev}),
such that $\gamma_g\circ \gamma_{g'} =\gamma_{gg'}$ for any pair $(g,g')$ of elements of $G$ and
$\gamma_e=\Id_\caF$ if $e$ is the identity element of $G$. We further assume that the underlying map
$\gamma \colon G\times X \rightarrow X$ is continuous. The space of compactly supported global sections
$\caF_c(X)$ is then by definition a representation space for $G$. By the equivalence of categories of
\ref{FCev}, we see that the data of a sheaf equipped with an action of $G$ is equivalent to the data of a
non-degenerate $\scrD(X)$-module, equipped with a linear action of $G$ compatible with that of
$\scrD(X)$. If the representation of $G$ in $\caF_c(X)$ is smooth, we say that $\gamma$ is smooth. A
distribution $T \in \caF_c(X)^*$ is said to be $G$-invariant if $\gamma_g\cdot T=T$ for all $g $ in $G$.
\end{defi}

In such a situation, let us examine the space of coinvariants $(\caF_c(X))_G$. Suppose there exist a t.d.
space $Y$ and a continuous map $q\colon X \rightarrow Y$, such that $q(\gamma_g(x))=q(x)$, for all
$x\in X$ and all $g\in G$. Then $\caF_c(X)$ is naturally equipped with a $\scrD(Y)$-module structure, and
$\caF_c(X)(G)$ is a $\scrD(Y)$-submodule of $\caF_c(X)$. It follows that $\caF_c(X)_G$ is a
$\scrD(Y)$-module, which implies that it corresponds via the equivalence of categories of Theorem
\ref{FCev} to a sheaf on $Y$, which we will denote by $\caF'$. By definition $\caF_c(X)_G$ and
$\caF'_c(Y)$ are isomorphic as $\scrD(Y)$-modules.

\begin{prop}
The stalk $\caF_y'$ of $\caF'$ at a point $y\in Y$ is naturally isomorphic to $\caF_c(q^{-1}(\{y\}))_G$.
\end{prop}

\begin{proof}
Let $Z=q^{-1}(\{y\})$. It is a closed subset of $X$, and the restriction
$p_Z \colon \caF_c(X)\rightarrow \caF_c(Z)$ is surjective (second proposition of \ref{FCev}). The kernel
of $p_Z$ is the subspace $L\subset \caF_c(X)$ generated by the sections of the form $f\cdot \phi$, where
$f\in \scrD(Y)$, $f(y)=0$, and $\phi\in \caF_c(X)$. Indeed, it is clear that $L \subset \ker p_Z$.
Conversely, if $\phi \in \ker p_Z$, then $q(\supp \phi)$ is a compact set in $Y$ which does not contain
$y$. We can cover $q(\supp \phi)$ by a compact open set $U$ which does not contain $y$ (this is a
consequence of Lemma \ref{recouv}). If $f$ is the characteristic function of $U$, we have
$f\in \scrD(Y)$, $f(y)=0$ and $f\cdot \phi=\phi$, hence $\phi \in L$. We therefore have
$\caF_c(Z)\simeq \caF_c(X)/L$, hence $\caF_c(Z)_G\simeq \caF_c(X)/L'$ where $L'$ is the subspace of
$\caF_c(X)$ generated by $L$ and $\caF_c(X)(G)$. Consequently
\begin{align*}
 \caF_c(Z)_G &=\caF_c(X)/(\caF_c(X)(G)+L)\\
&=\left( \caF_c(X)/\caF_c(X)(G) \right) /\left( L/(L\cap  \caF_c(X)(G)) \right)
\end{align*}
Now, $  \caF_c(X)/\caF_c(X)(G)=\caF'_c(Y)$ and the description of the stalk of the sheaf $\caF'$ at $y$
made in Remark \ref{FCev} is
\[
\caF'_y=  \caF'_c(Y)/L''
\]
where $L''$ is the space generated by the $f\cdot \psi$, $\psi \in \caF'_c(Y)$ and $f \in \scrD(Y)$ such
that $f(y)=0$. We then see that $L''\simeq L/(L\cap  \caF_c(X)(G))$ which completes the proof.
\end{proof}

\begin{cor}
With the above hypotheses, if there are no $G$-invariant distributions for the restrictions of the sheaf
$\caF$ to the fibers $Z=q^{-1}(\{y\})$, then there are no $G$-invariant distributions in $\caF_c(X)^*$.
\end{cor}

\begin{proof}
By hypothesis, for every $y\in Y$, the space $\caF_c(q^{-1}(\{y\}))_G$ has no non-zero linear form, hence
is zero. By the proposition, all the stalks of the sheaf $\caF'$ vanish, so $\caF'=0$, and
$\caF_c(X)_G\simeq\caF'_c(Y)=0$. Now a $G$-invariant distribution in $\caF_c(X)^*$ vanishes on
$\caF_c(X)(G)$, hence defines a linear form on $\caF_c(X)_G=0$: it is zero.
\end{proof}

\subsection{Induction and sheaves}\label{indetfaisc}

Let us fix a continuous action $\gamma_0$ of a t.d. group $G$ on a t.d. space $X$. Let us form the
category $\scrC^\infty_{X,G}-\caM od$ whose objects are the sheaves of $\bbC$-vector spaces on $X$
equipped with a smooth action $\gamma$ of $G$ compatible with $\gamma_0$, and whose morphisms are the
$G$-equivariant sheaf morphisms. For example, if $X=\{x\}$ is a singleton, then
$\scrC^\infty_{X,G}-\caM od $ is simply $\caM(G)$. The correspondence $\caF \mapsto \caF_c(X)$ (which
induces an equivalence of categories between $\scrC^\infty_{X}-\caM od$ and $\caM(\scrD(X))$ according to
Theorem \ref{FCev}) gives by restriction to the subcategory $\scrC^\infty_{X,G}-\caM od $ a functor
\[
\Gamma_c \colon  \scrC^\infty_{X,G}-\caM od\longrightarrow  \caM(G).
\]
If $Q$ is a closed subgroup of $G$, and $Z$ a locally closed, $Q$-invariant subset of $X$, the functor
\[
\caF \mapsto \caF_{|Z}
\]
induces a functor
\[
\caR_{Z,Q}^{X,G}  \colon   \index[not]{RZQ@$\caR_{Z,Q}^{X,G}:$}
\scrC^\infty_{X,G}-\caM od \longrightarrow \scrC^\infty_{Z,Q}-\caM od.
\]
In particular, if $Z=\{x\}$ is a point of $X$, then $\caR_{Z,Q}^{X,G} \caF$ is the smooth representation
of $Q$ in the stalk $\caF_x$ of $\caF$.

We now place ourselves in the case where the action of $G$ on $X$ is transitive. If we assume moreover
that $G$ is $\sigma$-compact, $X$ is then homeomorphic to the quotient space $ G/H$, where $H$ is the
stabilizer of any point of $X$ (Corollary \ref{actiontdcor}). For reasons of compatibility with our
definition of induction, it is more appropriate to work with $H\backslash  G$ rather than $G/H$ (the
action of $G$ on $H\backslash G$ is given by $(g_0,Hg)\mapsto Hgg_0^{-1}$). In this case, we define an
``induction'' functor
\[
\caI_H^G:\index[not]{IHG@$\caI_H^G:$} \caM(H) \longrightarrow
\scrC^\infty_{X,G}-\caM od,
\]
as follows. Let $(\tau, E)$ in $\caM(H)$. The space of the induced representation $\ind_H^G E$ is
naturally a non-degenerate $\scrD(H\backslash G)$-module, by pointwise multiplication of functions: if
$\phi \in \scrD(H\backslash G)$, and $f \in \ind_H^G E$, $\phi f$ is in $\ind_H^G E$. According to the
equivalence of categories of Theorem \ref{faisceqcat}, $\ind_H^G E$ is the space of compactly supported
sections of a certain sheaf $\caF^{\tau}$ on $H\backslash G$. The functor $\caI_H^G$ sends $(\tau,E)$ to
$\caF^{\tau}$. Moreover it is clear that the representation of $H$ in the stalk of the sheaf
$\caF^{\tau}$ over the point $o=H\mathbf{1}_G$ of $H\backslash G$ is isomorphic to $(\tau,E)$.

Conversely, if $\caF$ is a sheaf on $H\backslash G$ equipped with a smooth action of $G$ compatible with
the action of $G$ on $H\backslash G$, let $(\tau,E)$ denote the representation of $H$ in the stalk of the
sheaf $\caF$ over the point $o$ of $H\backslash G$. Then the representation of $G$ in
$\caF_c(H\backslash G)$ is isomorphic to $\ind_H^G (\tau,E)$ and the smooth part of the representation of
$G$ in $\caF(H\backslash G)$ is isomorphic to $\Ind_H^G (\tau,E)$. More explicitly, the first of these
isomorphisms is given by
\[
\alpha \colon \caF_c(H\backslash G)\rightarrow \ind_H^G (\tau,E), \quad  \phi\mapsto  [ g \mapsto (\gamma_g\cdot \phi)(o)],
\]
with inverse
\[
\beta \colon  f \mapsto [ Hg^{-1}= \gamma_g (o)\mapsto \gamma_g(f(g^{-1}))],
\]
and the second by analogous formulas.

One easily verifies that these isomorphisms are natural. We have therefore obtained the

\enlargethispage{2\baselineskip}
\begin{thm}
The functors
\[
\caI_H^G \colon  \caM(H) \longrightarrow \scrC^\infty_{H\backslash G,G}-\caM od
\quad \text { and } \quad \caR_{H\mathbf{1}_G,H}^{H\backslash G,G} \colon
\scrC^\infty_{H\backslash G,G}-\caM od  \longrightarrow \caM(H)
\]
are inverse to each other and define an equivalence of categories between $\caM(H)$ and
$\scrC^\infty_{H\backslash G,G}-\caM od$.
\end{thm}

\section{Notes on Chapter III}
Most of this chapter is taken from \cite{BeZe1}. We have added some elements inspired by the philosophy
of \cite{KV}, which consists of viewing the induction and Jacquet functors as particular cases of the
forgetful ($\caF,  {}\spcheck\caF$) or induction ($I,P$) functors defined in the framework of algebras
with idempotents. Beyond the conceptual gain, this allows one to quickly determine the adjunctions
between composites of such functors. The part on representations in spaces of sections of equivariant
sheaves is taken from \cite{BeZe2}.

\chapter[Compact, unitary representations...]{Compact, square-integrable, and unitary representations}

We will now study particular classes of smooth representations of a t.d. group. We begin by recalling the
representation theory of compact groups, a theory which, in the framework of t.d. groups, is particularly
elementary. The main result is that the category of smooth representations of a compact t.d. group is
semisimple. This is not the case in general if the group $G$ is not compact, but it is then possible to
define a class of representations having good properties, by a support condition on their matrix
coefficients: the compact representations. Any compact representation is semisimple, and if we fix an
isomorphism class $\tau$ of irreducible compact representations, any smooth representation of $G$
decomposes into a direct sum where the first factor is a direct sum of representations in the class
$\tau$, and the second is a smooth representation of which no irreducible subquotient is in $\tau$. This
is interpreted in terms of a decomposition of categories. We give a simple finiteness criterion for the
category $\caM(G)$ to decompose into a ``compact part'', which is a product, over the isomorphism classes
$\tau$ of compact representations, of the categories $\caM(G)_\tau$ whose elements are the direct sums of
representations in the class $\tau$, and a non-compact part, whose elements are the representations
having no compact subquotient. This finiteness criterion will be verified for $p$-adic reductive groups,
which is the basis of the Bernstein decomposition theorem. We then study unitary representations and we
introduce a twisted version of the induction functors which preserves the unitarity of representations.
Admissible unitary representations are semisimple. Finally we look at smooth representations whose
central character is unitary and whose matrix coefficients are square-integrable functions (modulo the
center). Such a representation is admissible and unitary and the coefficients of such representations
satisfy the Schur orthogonality formulas.

\section{Compact representations}

\subsection{Representations of compact groups}\label{rgc}
Let $K$ be a compact t.d. group.

\needspace{142pt}
\begin{thm}
$(i)$ Any irreducible smooth representation of $K$ is finite-dimensional.

$(ii)$ For any finite-dimensional smooth representation $(\pi, V)$ of $K$, there exists a normal compact
open subgroup $N$ in $K$ such that $N$ acts trivially on $V$.

$(iii)$ Any smooth representation of $K$ is semisimple.

$(iv)$ Any smooth representation of $K$ is unitary.
\end{thm}

\begin{proof}
$(i)$ Let $(\pi, V)$ be an irreducible smooth representation of $K$, and let $v\in V$ be non-zero. Since
$\pi$ is smooth, there exists a compact open subgroup $K_v$ of $K$ which fixes $v$. Since $K/K_v$ is
finite (cf. Lemma \ref{fini}), the vector space $W$ generated by $K\cdot v$ is finite-dimensional. Since
it is $K$-stable and non-zero, by the irreducibility of $\pi$, $W=V$.

$(ii)$ Let $(\pi, V)$ be a finite-dimensional smooth representation of $K$. Let $N=\ker \pi$. It is a
normal subgroup of $K$. Let $\{v_i\}_{i=1,\ldots,n}$ be a basis of $V$. Then
\[
N =\bigcap_{i=1}^n \mathrm{Stab}_K(v_i).
\]
Since $\pi$ is smooth, each $\mathrm{Stab}_K(v_i)$ is an open subgroup of $K$, so $N$ is open. Since $K$
is compact, any open subgroup of $K$, being closed, is compact.

$(iii)$ Let us first show that $V$ is a sum of irreducible subrepresentations. Let $v \in V$ and let $W$
be the subspace generated by $K\cdot v$. As in the proof of $(i)$, we see that $W$ is finite-dimensional,
and therefore according to $(ii)$, there exists a normal compact open subgroup $N$ in $K$ such that $N$
acts trivially on $W$. Since $K/N$ is finite, we can according to the representation theory of finite
groups (see for example \cite{Se}) decompose $W$ into a finite direct sum of irreducible representations
of $K$. This shows the assertion. Lemma \ref{lemmesemsimP} then allows us to conclude.

$(iv)$ Let $(\pi, V)$ be a smooth representation of $K$. Let us choose a positive definite Hermitian
inner product $\bilo$ on $V$. Let us set:
\[
\langle v,\, w\rangle_0= \int_K  \langle \pi(k)\cdot v ,\,
\pi(k)\cdot w \rangle  \; d\mu_K(k),\qquad (v,w \in V)
\]
where $\mu_K$ is a Haar measure on $K$. Then $\bilo_0$ is a $K$-invariant positive definite Hermitian
inner product on $V$.
\end{proof}

\medskip

Let $\widehat K $ denote the set of equivalence classes of irreducible smooth representations of $K$.

\begin{cor}
The category $\caM(K)$ is semisimple. More precisely, any smooth representation $(\pi, V)$ of $K$
decomposes canonically into:
\[
V=\bigoplus_{(\sigma,V_\sigma)\in \widehat K} V(\sigma),
\]
where $V(\sigma)$ is the image of the canonical morphism:
\[
\Hom_K(V_\sigma,V)\otimes  V_\sigma, \quad \phi\otimes v \mapsto \phi(v).
\]
The subrepresentation $V(\sigma)$ is a direct sum of irreducible representations isomorphic to
$(\sigma,V_\sigma)$.
\end{cor}

\begin{proof}
Let $V=\bigoplus_{i\in I} V_i$ be a decomposition of $V$ into irreducible subrepresentations. If
$V_i \simeq V_\sigma$, there exists an injective intertwining operator
$\phi_i\colon V_\sigma \rightarrow V$ whose image is $V_i$. Therefore $V_i \subset V(\sigma)$ and we
deduce that $V =\sum V(\sigma)$. Let $I_\sigma$ be the set of $i$ such that $V_i \simeq V_\sigma$. We
then have $\bigoplus_{i\in I_\sigma} V_i \subset V(\sigma)$. For any morphism
$\phi \in \Hom_K(V_\sigma,V)$, the composition of $\phi$ with a projection $p_j\colon V \rightarrow V_j$
is zero if $j \notin I_\sigma$, and therefore $V(\sigma)\subset \bigoplus_{i\in I_\sigma} V_i$.
\end{proof}

\medskip

We call $V(\sigma)$ \index[not]{V_s@$V(\sigma)$} the isotypic component \indexter{isotypic component} of
type $\sigma$ of $V$. The dimension of $\Hom_K(V_\sigma,V)$, which we will denote by $m(\sigma)$, is
called the multiplicity of $\sigma$ in $V$. When this dimension is finite, we have:
\[
m(\sigma)= \dim V(\sigma)/\dim V_\sigma
\]

\subsection{An admissibility criterion}

From the above, we derive a necessary and sufficient condition for a smooth representation of a t.d.
group to be admissible.

\begin{prop}
Let $(\pi, V)$ be a smooth representation of a t.d. group $G$. Then $(\pi, V)$ is admissible if and only
if for any compact open subgroup $K$ of $G$, any class $\sigma \in \widehat K$ has finite multiplicity in
$\res_K^G (\pi,V)$.
\end{prop}

\begin{proof}
Suppose that for a certain compact open subgroup $K$ of $G$, and a certain class $\sigma \in \widehat K$,
the multiplicity $m(\sigma)$ is not finite. Let $N$ be a normal compact open subgroup of $K$ such that
$N \subset \ker \sigma$. Then $\dim V^N$ is not finite and therefore $(\pi, V)$ is not admissible.
Conversely if $K$ is a compact open subgroup of $G$, the trivial representation of $K$ being irreducible,
$\dim V^K$ is finite.
\end{proof}

\subsection{Compact representations}\label{comprep}
\indexter{representation!compact}

Let $G$ be a t.d. group. If $G$ is not compact, the category $\caM(G)$ is generally not semisimple.
Nevertheless, there exists a class of smooth representations of $G$ which behave like the representations
of compact groups. In particular they are semisimple.

Throughout this section, we assume $G$ is unimodular.

\begin{defi}
Let $(\pi, V)$ be a smooth representation of $G$. We say that $(\pi, V)$ is compact if for all $v\in V$,
and for all compact open subgroups $K$ of $G$, the function:
\[
f_{K,v}  \colon G \rightarrow V ,\quad g \mapsto \pi(e_K)\pi(g^{-1})\cdot v
\]
is compactly supported.
\end{defi}

It is clear that if $(\pi, V)$ is compact, any subquotient of $(\pi, V)$ is also compact.

Let $(\pi, V)$ be a smooth representation of $G$ and let $v \in V$, $\lambda \in \widetilde V$. Then the
matrix coefficient $\phi_{v,\lambda}$ is the locally constant function:
\[
\phi_{v,\lambda} \colon G \rightarrow \bbC,\quad g \mapsto \lambda(\pi(g)\cdot v).
\]

\begin{thm}
A smooth representation of $G$ is compact if and only if all its matrix coefficients are compactly
supported.
\end{thm}

\begin{proof}
Suppose $(\pi,V)$ is compact. Let $v\in V$, $\lambda \in \widetilde V$ and $K$ a compact open subgroup of
$G$ such that $\lambda \in \widetilde V^K$. Then it is clear that the support of $\phi_{v,\lambda}$ is in
the support of $\check f_{K,v}$ (there is a $\pi(g^{-1})$ in the definition of the function $f_{K,v}$ and
a $\pi(g)$ in that of $\phi_{v,\lambda}$, which explains the composition with the involution
$\check f_{K,v}(g)=f_{K,v}(g^{-1})$). For the converse, we will find a finite number of
$\lambda_i \in \widetilde V^K$ such that
$\supp f_{K,v} \subset \bigcup_i  \supp \check \phi_{v,\lambda_i}$. The image of $f_{K,v}$ is in $V^K$.
Let $E_v$ denote the subspace of $V^K$ generated by this image. Extract a basis $(v_i)_{i\in I}$ of $E_v$
from the generating system $(\pi(e_K)\pi(g)\cdot v)_{g \in G}$ and choose a linear form $\lambda_0$ on
$V^K$ taking the value $1$ on the $v_i$. Extend $\lambda_0$ to a linear form $\tilde \lambda_0$ on $V$,
vanishing on $V(K)$ (recall that $V=V^K\oplus V(K)$ according to Proposition \ref{piK}). Then
$e_K \cdot \tilde \lambda_0= \tilde \lambda_0$ and
\[
g \mapsto  \tilde \lambda_0(\pi(g)\cdot v)=\lambda_0(\pi(e_K)\pi(g)\cdot v)
\]
is compactly supported, invariant under the action by left translation of $K$, and therefore its support
is covered by a finite number of subsets of the form $K\cdot g_j$. This shows that the subspace $E_v$ is
finite-dimensional, say $\ell$. Then choose $\ell$ elements $\lambda_1,\ldots , \lambda_\ell$ of
$\widetilde V^K$ separating the points of $E_v$. Then
$\supp f_{K,v} \subset \bigcup_{i=1}^\ell  \supp \check \phi_{v,\lambda_i}$.
\end{proof}

\begin{prop}
Any finitely generated compact representation is admissible.
\end{prop}

\begin{proof}
Let $(\pi,V)$ be such a representation. Suppose that $V$ is generated by $v_1,\ldots,v_l$. If $K$ is a
compact open subgroup of $G$, then $V^K=\pi(e_K)\cdot V$ is generated by the vectors of the form
$\pi(e_K)\pi(g)\cdot v_i$, $g\in G$. It therefore suffices to see that the space generated is
finite-dimensional for fixed $i$. Now, we can repeat the argument given in the proof of the theorem
above, consisting of constructing a linear form $\lambda$ fixed by $K$, taking the value $1$ on a maximal
linearly independent family $\pi(e_K)\pi(g_j)\cdot v_{i}$ and $0$ on $V(K)$. The compactness of the
support of the matrix coefficient $\phi_{v_{i}, \lambda}$ shows that this space is finite-dimensional.
\end{proof}

\medskip

Compact representations of a group $G$ do not always exist. A necessary condition is given in the
following lemma.

\begin{lemme}
Let $G$ be a $\sigma$-compact t.d. group, and suppose that there exists a compact representation
$(\pi,V)$ of $G$. Then the center of $G$ is compact.
\end{lemme}

\begin{proof}
By taking an irreducible subquotient, we can assume that $(\pi,V)$ is irreducible. Then $(\pi,V)$ admits
a central character $\chi_\pi$. Let $v \in V$ be non-zero and $K$ a compact open subgroup of $G$ fixing
$v$. Then for all $z \in Z(G)$,
\[
f_{K,v}(z)=\pi(e_K)\pi(z^{-1})\cdot v= \chi_\pi(z^{-1})v.
\]
The support of $ f_{K,v}$ contains $Z(G)$, which is therefore compact.
\end{proof}

\subsection{Decomposition of representations}\label{repcomp2}

Assume that we are under the hypotheses of the previous paragraph and that $G$ is $\sigma$-compact. The
main result on compact representations is the following:

\begin{thm}
Let $(\tau,W)$ be an irreducible compact representation of $G$. Then any representation $(\pi,V)$ of
$\caM(G)$ decomposes into:
\[
V=V(\tau)\oplus V(\tau)^\perp
\]
where $V(\tau)$ is a direct sum of representations isomorphic to $(\tau,W)$, and none of the irreducible
subquotients of $V(\tau)^\perp$ is isomorphic to $(\tau,W)$.
\end{thm}

The proof requires intermediate results and will span the next few sections. Recall that we assumed $G$
is unimodular. Let us fix a Haar measure $\mu_G$ (hence bi-invariant). Let us identify $\scrD(G)$ and
$\caH(G)$ via the choice of this Haar measure. The spaces $\scrD(G)$ and $\End(W)$ are equipped with a
representation structure of $G\times G$:
\[
(g_1,g_2)\cdot f(x)= f(g_1^{-1}xg_2) \quad (g_1,g_2,x\in G), \, (f \in \scrD(G)),
\]
\[
((g_1,g_2)\cdot \lambda)(v)= g_1\cdot\lambda(g_2^{-1}\cdot v),\quad
(g_1,g_2\in G), \, (\lambda \in \End(W)), (v\in W).
\]
Moreover, the morphism
\[
\tau \colon  \scrD(G)\rightarrow \End(W),\quad f\mapsto \tau(f)
\]
intertwines these actions of $G\times G$. Since $\scrD(G)$ is a smooth representation, it follows that
the image of $\tau$ takes values in the smooth part $\End(W)_0$ of $\End(W)$. We will see that this
morphism is surjective, and construct a section.

Let us first show that as a representation of $G \times G$:
\begin{equation}\label{WWEW}
W\otimes \widetilde W \simeq \mathrm{End} (W)_0
\end{equation}
To do this, let us define the morphism
$\alpha \colon  W\otimes \widetilde W \rightarrow \mathrm{End}(W)_0 $ by:
\begin{equation}\label{AWLV}
  \alpha(w\otimes \lambda)(v)= \lambda(v) w.
\end{equation}
Injectivity and $G\times G$-equivariance are clear.

Since irreducible compact representations are admissible, we can show surjectivity by establishing that
for any compact open subgroup $K$ of $G$, the induced map:
\begin{equation}\label{AKW}
\alpha^K \colon   (W\otimes \widetilde W)^{K\times K}  \longrightarrow  \mathrm{End} (W)^{K\times K}
\end{equation}
is surjective by a dimension argument. The space $\mathrm{End} (W)^{K\times K}$ is naturally a subspace
of $\mathrm{End} (W^{K})$, via the restriction map. Indeed, $W$ decomposes into $W=W^K\oplus W(K)$
(Proposition \ref{piK}) and if $A \in \mathrm{End}(W)^{K \times  K}$, we have $A=\tau(e_K)A\tau(e_K)$,
and the restriction of $A$ to $W(K)=\ker \tau(e_K)$ is zero. Therefore if the restriction of $A$ to $W^K$
is zero, $A$ is too. We then have:
\[
\dim (\mathrm{End} (W)^{K\times K})\leq (\dim  W^{K})^2 =   \dim (W\otimes \widetilde W)^{K\times K}.
\]
which completes the proof of the surjectivity of (\ref{AKW}) and therefore of (\ref{AWLV}). Note that the
proof in fact establishes (\ref{WWEW}) for any admissible representation $W$.

Let $\phi \colon  W\otimes \widetilde W \rightarrow \scrD(G)$ be the linear map
$v\otimes \lambda \mapsto \check{\phi}_{v,\lambda}$ and let $\psi$ be the unique linear map from
$\End(W)_0$ to $\scrD(G)$ such that $\psi\circ \alpha=\phi$. It is clear that $\psi$ is
$G\times G$-equivariant. We now show that, up to a scalar coefficient, it is the sought-after section of
the morphism $\tau$. Since $\tau\circ \psi$ is an intertwining operator for the representation of
$G\times G$ in $\End(W)_0$, and since this is isomorphic to $W\otimes \widetilde W$, hence irreducible
according to Proposition \ref{tensprod}, Schur's lemma tells us that there exists a scalar $\kappa(\tau)$
such that $\tau\circ \psi=\kappa(\tau)\Id$.

\begin{prop}
$(i)$ If $(\rho,E)$ is an irreducible smooth representation of $G$ not equivalent to $(\tau,W)$, then for
all $f$ in the image of $\psi$, $\rho(f)=0$.

$(ii)$ The scalar $\kappa(\tau)$ is non-zero.
\end{prop}

\begin{proof}
$(i)$ Let $v \in E$, and consider the morphism of $G$-modules
\[
W\otimes \widetilde W \rightarrow E, \quad w\otimes \lambda \mapsto
\rho (\phi( w\otimes \lambda))\cdot v .
\]
As a representation of $G$ ($G$ acting only on the first factor), $W\otimes \widetilde W$ is a direct sum
of irreducible representations equivalent to $(\tau,W)$. When $(\rho,E)$ is not equivalent to $(\tau,W)$,
the image of this morphism is necessarily zero, and therefore $\rho(f)\cdot v=0$ for all $v\in E$, for
all $f$ in the image of $\phi$ (and therefore of $\psi$). This proves $(i)$.

Let $f$ be a non-zero element in the image of $\psi$. We want to show that $\tau(f)$ is non-zero. To see
this, it suffices to invoke the completeness theorem (Theorem \ref{separlemme}) and $(i)$.
\end{proof}

We set $d(\tau)=\kappa(\tau)^{-1}$. \index[not]{d_t@$d(\tau)$} We call this scalar the formal degree of
$\tau$.

\begin{rmqs}
$ $
\begin{itemize}
\item[1.] The scalar $d(\tau)$ depends on the choice of the Haar measure.

\item[2.] If $G$ is compact, and we normalize the Haar measure by $\int_G d\mu_G=1$, then the formal
degree of an irreducible smooth representation of $G$ is its dimension.

\item[3.] We will more generally define the formal degree of any irreducible square-integrable smooth
representation of $G$ in Section \ref{orthschur}.
\end{itemize}
\end{rmqs}

We continue the proof of the main theorem of this section in the following sections.

\subsection{Projection onto a compact representation}\label{cncomp}

Let $(\tau,W)$ be an irreducible compact representation of $G$. Let $K$ be a compact open subgroup of
$G$. The results of the previous paragraph allow us to obtain the:

\begin{thm}
$(i)$ There exists a unique distribution $e_{K,\tau}$ \index[not]{e_Kt@$e_{K,\tau}$} in $\caH(G)$ such
that $\tau(e_{K,\tau})=\tau(e_K)$ and $\rho(e_{K,\tau})=0$ if $(\rho,E)$ is an irreducible smooth
representation of $G$ not equivalent to $(\tau,W)$.

$(ii)$ If $K'$ is a compact open subgroup of $G$ contained in $K$, we have:
\[
e_{K',\tau}* e_{K,\tau}=e_{K}*e_{K',\tau}=e_{K',\tau}*e_K=e_{K,\tau}.
\]
In particular $e_{K,\tau}$ is an idempotent.

$(iii)$ If $g\in G$, then $\delta_g*e_{K,\tau}*\delta_{g^{-1}}=e_{gKg^{-1},\tau}$.
\end{thm}

\begin{proof}
The uniqueness of such a distribution is a direct consequence of the completeness theorem
\ref{separlemme}. Let us set:
\[
e_{K,\tau}= d(\tau)\psi(\tau(e_K)) \in \caH(G).
\]
Then,
\[
\tau( e_{K,\tau})= d(\tau) (\tau\circ \psi\circ \tau)(e_K)=\tau(e_K)
\]
Similarly since $\rho\circ \psi=0$:
\[
\rho( e_{K,\tau})= d(\tau) (\rho\circ \psi \circ \tau)(e_K)=0.
\]
To show the first assertion of $(ii)$, by the completeness theorem \ref{separlemme}, it suffices to
verify that evaluating either of these distributions at $\rho$, for any irreducible smooth representation
$(\rho,E)$ of $G$, yields the same result. It is clear that according to $(i)$, we always obtain $0$ if
$(\rho,E)$ is not equivalent to $(\tau,W)$, and always $\tau(e_K)$ if $(\rho,E)$ is equivalent to
$(\tau,W)$ (because $e_K*e_{K'}=e_{K'}*e_K=e_K$). The last assertion is proved similarly.
\end{proof}

Let $(\pi,V)$ be a smooth representation of $G$, and let $v\in V$. If $v\in V^K$, we have
$\pi(e_{K',\tau})\cdot v=\pi(e_{K,\tau})\cdot v$ for all $K' \subset K$. Let us denote this vector simply
by $\pi(e_{\tau})\cdot v$. This defines an operator $\pi(e_\tau)$ on $V$. We will sometimes lighten the
notation by simply writing $e_{\tau}\cdot v$ and $e_\tau$ when the context clearly indicates which
representation $(\pi,V)$ we are working with. Note that despite what this notation suggests, we cannot
define an element $e_\tau$ in $\caH(G)$ inducing the operators $\pi(e_\tau)$. We will see later that
$e_\tau$ in fact defines an element of the center of the category $\caM(G)$.

\begin{prop}
$(i)$ $\pi(e_\tau)$ is a projector of $V$ which commutes with the action of $G$.

$(ii)$ If $(\pi',V')$ is another smooth representation of $G$, and if $A\in \Hom_G(V,V')$, then
$A\circ \pi(e_\tau)=  \pi'(e_\tau)\circ A$.

$(iii)$ $V = \im \pi(e_\tau)\oplus \ker \pi(e_\tau)$ is a decomposition of $V$ into subrepresentations,
and $\im \pi(e_\tau)$ is a direct sum of irreducible subrepresentations equivalent to $(\tau,W)$. None of
the irreducible subquotients of $\ker \pi(e_\tau)$ is equivalent to $(\tau,W)$.
\end{prop}

\begin{proof}
$(i)$ The fact that $\pi(e_\tau)$ is a projector follows immediately from $(ii)$ of the previous theorem.
We have, for all $K$ small enough:
\begin{align*}
 \pi(g)\pi(e_\tau)\cdot v&=  \pi(g)\pi(e_{K,\tau})\cdot v=
\pi(g)\pi(e_{K,\tau})\pi(g^{-1})\pi(g)\cdot v\\
& = \pi(e_{gKg^{-1},\tau})\pi(g)\cdot v =\pi(e_\tau)\pi(g)\cdot v
\end{align*}
and therefore $\pi(e_\tau)$ commutes with the action of $G$. Point $(ii)$ is proved in the same way. The
decomposition of $(iii)$ is immediate, since $\pi(e_\tau)$ is a $G$-equivariant projector.

We will show that $\im \pi(e_\tau)$ is a direct sum of irreducible subrepresentations equivalent to
$(\tau,W)$. According to Lemma \ref{lemmesemsimP}, it suffices to prove that $\im \pi(e_\tau)$ is
generated by subrepresentations equivalent to $(\tau,W)$. Let $w=\pi(e_\tau)\cdot v$ be in the image of
$\pi(e_\tau)$. Let $K$ be a compact open subgroup of $G$ such that $v \in V^K$. Then
$w=\pi(e_\tau)\cdot v=\pi(e_{K,\tau})\cdot v$. In the proof of Theorem \ref{repcomp2}, we introduced the
$G\times G$-morphism:
\[
\phi \colon W\otimes \widetilde W \rightarrow \scrD(G)\simeq \caH(G).
\]
Since $W\otimes \widetilde W $ is irreducible according to Proposition \ref{tensprod}, and since $\phi$
is non-zero, $\phi$ is injective. By definition, $e_{K,\tau}$ is in the image of $\phi$, and therefore
$\caH(G)*e_{K,\tau}$ embeds into $\im \phi$. Since $W\otimes \widetilde W $ is, as a representation of
$G$, a direct sum of irreducible representations equivalent to $W$, the same is true for
$\caH(G)*e_{K,\tau}$. The image of the $G$-morphism
\[
\caH(G)*e_{K,\tau}\rightarrow V,\quad     h* e_{K,\tau}
\mapsto \pi( h* e_{K,\tau})\cdot v
\]
is therefore a direct sum of irreducible representations equivalent to $W$. This completes the proof of
the fact that $\im \pi(e_\tau)$ is a direct sum of representations isomorphic to $(\tau,W)$.

The functor:
\begin{equation} \label{projV}
\caM(G) \rightarrow \caM(G), \quad  V \mapsto   \pi( e_{\tau})\cdot V
\end{equation}
is exact. Indeed, it is clearly left exact. Suppose that
\[
\beta \colon  (\pi_1,V_1) \rightarrow (\pi_2,V_2)
\]
is a surjective morphism of $G$-modules. For all $w\in V_2$, there exists $v\in V_1$ such that
$\beta(v)=\pi_2(e_\tau)\cdot w$, and therefore
\[
\beta(\pi_1(e_\tau)\cdot v)=\pi_2(e_\tau)\cdot\beta(v)= \pi_2(e_\tau) \circ \pi_2(e_\tau)\cdot w= \pi_2(e_\tau)\cdot w     .
\]
This shows the assertion.

\medskip

We can now show the last assertion of the proposition. Since $\pi(e_\tau)$ annihilates any subquotient of
$\ker \pi(e_\tau)$, by the exactness of the functor above and the fact that $\tau(e_\tau)$ is the
identity, no subquotient of $\ker \pi(e_\tau)$ is isomorphic to $(\tau,W)$.
\end{proof}

\bigskip

Theorem \ref{repcomp2} is now completely proved.

\begin{rmqs}
$ $
\begin{itemize}
\item[1.] The subspace $\ker \pi(e_\tau)$ is the only $G$-invariant complement of $\im  \pi(e_\tau)$.

\item[2.] If $(\tau',W')$ is another irreducible compact representation of $G$, not equivalent to
$(\tau,W)$, then:
\[
\pi(e_\tau)\circ \pi(e_{\tau'}) = \pi( e_{\tau'})\circ \pi(e_\tau)=0.
\]

\item[3.] It is easy to verify that $V\mapsto \pi(e_\tau)\cdot V$ defines a natural transformation of the
identity functor of the category $\caM(G)$ to itself. We can therefore interpret $e_\tau$ as an element
of the center of $\caM(G)$.
\end{itemize}
\end{rmqs}

\subsection{Semisimplicity of compact representations} \label{conscomp}

Let us now deduce some consequences from the results obtained above.

\begin{cor}
Any compact representation is semisimple.
\end{cor}

\begin{proof}
According to Lemma \ref{lemmesemsimP}, it suffices to prove that a compact representation is generated by
its irreducible subrepresentations. Let $(\rho,W)$ be a compact representation of $G$, and let $W^f$ be
the sum of all irreducible subrepresentations of $W$. We must therefore show that $W/W^f$ is zero.
Suppose otherwise, and let $\pi$ be an irreducible subquotient of $W/W^f$. It is a compact
representation, since any subquotient of a compact representation is a compact representation. We then
have $e_\pi\cdot (W/W^f)\neq 0$ and therefore $e_\pi\cdot W$ is not included in $W^f$. But this leads to
a contradiction, because by $(iii)$ of Proposition \ref{cncomp}, $e_\pi\cdot W\subset W^f$.
\end{proof}

Let us introduce some notation:

\begin{notation}
Let $(\tau,W)$ be an irreducible compact representation of the t.d. group $G$. Let $\caM(G)_\tau$
\index[not]{M(G)_t@$\caM(G)_\tau$} denote the full subcategory of $\caM(G)$ consisting of representations
which are a direct sum of irreducible representations equivalent to $(\tau,W)$ and
$[\caM(G)\setminus\tau]$\index[not]{M(G)_t1@$[\caM(G)\setminus\tau]$} the full subcategory of $\caM(G)$
consisting of representations of which no irreducible subquotient is isomorphic to $(\tau,W)$. More
generally, given non-isomorphic irreducible compact representations $(\tau_i,W_i)$, $i=1,\ldots, r$, let
$[\caM(G)\setminus\tau_1,\ldots,\tau_r]$ denote the full subcategory of $\caM(G)$ consisting of
representations of which no irreducible subquotient is isomorphic to one of the $(\tau_i,W_i)$.

Let $\caM(G)_{c}$ \index[not]{M(G)c@$\caM(G)_{c}$} denote the full subcategory of $\caM(G)$ consisting of
representations of which all irreducible subquotients are compact representations, and $\caM(G)_{nc}$
\index[not]{M(G)_nc@$\caM(G)_{nc}$} the full subcategory of $\caM(G)$ consisting of representations of
which no irreducible subquotient is a compact representation. Let $\mathbf{Irr}(G)_c$
\index[not]{Irr(G)_c@$\mathbf{Irr}(G)_{c}$} denote the set of isomorphism classes of irreducible compact
representations.
\end{notation}

The category $\caM(G)_{c}$ is therefore semisimple, and any simple module in it is both projective and
injective.

\begin{prop}
Any irreducible compact representation $(\tau,W)$ of $G$ is projective and injective in $\caM(G)$.
\end{prop}

\begin{proof}
Let $(\pi,V)$ be in $\caM(G)$. We have:
\[
\Hom_G(W,V)=\Hom_G(W,\pi(e_\tau)\cdot V),
\]
\[
\Hom_G(V,W)=\Hom_G(\pi(e_\tau)\cdot V,W).
\]
This reduces us to showing that $W$ is projective and injective in $\caM(G)_\tau$, but as we noted above,
this is indeed the case.
\end{proof}

\subsection{Decomposition of $\caM(G)$} \label{decMG}

Let us reinterpret the results of the previous paragraph in terms of category decompositions
(\ref{deccat}). Let $(\tau,W)$ be an irreducible compact representation of the t.d. group $G$. The
category $\caM(G)$ decomposes into:
\[
\caM(G)= \caM(G)_\tau \times  [\caM(G)\setminus\tau].
\]

More generally, given non-isomorphic irreducible compact representations $(\tau_i,W_i)$, $i=1,\ldots, r$,
we have a decomposition:
\[
\caM(G)= \prod_i \caM(G)_{\tau_i} \times  [\caM(G)\setminus\tau_1,\ldots,\tau_r].
\]

On the other hand $\caM(G)_{c}$ decomposes into:
\[
\caM(G)_{c}=\prod_{\tau \in \mathbf{Irr}(G)_c} \caM(G)_\tau.
\]

It is natural to ask whether the category $\caM(G)$ decomposes as
\[
\caM(G)= \caM(G)_c \times \caM(G)_{nc},
\]
where $\caM(G)_{nc}$ denotes the category of smooth representations of $G$ of which no subquotient is a
compact representation. We give a criterion for this.

\begin{prop}
The category $\caM(G)$ splits into a product of categories
\[
\caM(G)= \caM(G)_c \times \caM(G)_{nc},
\]
if and only if the following condition is satisfied:

{\bf(KF)}: For any compact open subgroup $K$ of $G$, the number of isomorphism classes of irreducible
compact representations $(\pi,V)$ of $G$ such that $V^K\neq 0$ is finite.
\end{prop}

\begin{proof}
Suppose condition {\bf(KF)} is satisfied. We can then define for each smooth representation $(\rho,W)$ of
$G$, an operator
\[
\rho(e_c)= \sum_{\tau \in \mathbf{Irr}(G)_c}  \rho(e_\tau)
\]
Indeed, any vector $w \in W$ is fixed by a certain compact open subgroup $K$ of $G$, and when we apply
$\rho(e_c)$ to $w$, by hypothesis only a finite number of $\rho(e_\tau)\cdot w$ are non-zero.

Since $\rho(e_{\pi_1}) \rho(e_{\pi_2})=0$ if $(\pi_1,V_1)$ and $(\pi_2,V_2)$ are in $\mathbf{Irr}(G)_c$
and not equivalent, we see that $\rho(e_c)$ is a projector. Therefore
\[
W= \rho(e_c)\cdot W\oplus (1- \rho(e_c))\cdot W
\]
as $G$-modules. Since $ \rho(e_c)\rho(e_\pi) =\rho(e_\pi) \rho(e_c)=\rho(e_\pi)$ for any representation
$(\pi,V)\in \mathbf{Irr}(G)_c$, we have
\[
\rho(e_c)\cdot W=\bigoplus_{(\pi,V)\in \mathbf{Irr}(G)_c }\rho(e_\pi)\cdot W.
\]
It follows that $\rho(e_c)\cdot W$ is the largest submodule of $W$ whose irreducible subquotients are
compact representations. On the other hand, since
\[
(1-\rho(e_c))\rho(e_\pi) =\rho(e_\pi)(1- \rho(e_c))= 0
\]
for any representation $(\pi,V)\in \mathbf{Irr}(G)_c$, according to Proposition \ref{cncomp}, $(iii)$,
$(1- \rho(e_c))\cdot W$ admits no compact subquotient.

Finally, we show the uniqueness of this decomposition. Suppose there exists a second one of the same form
\[
W=W'_c\oplus W'_{nc}.
\]
Since $W'_c$ is in $\caM(G)_c$, $\rho(e_c)$ acts as the identity on $W'_c$, and therefore
$W'_c\subset W_c$. On the other hand, $\rho(e_c)\cdot W'_{nc}=0$, otherwise $W'_{nc}$ would contain a
compact representation. Therefore $W'_{nc}\subset W_{nc}$, which shows the uniqueness of the
decomposition.

Conversely, let $K$ be a compact open subgroup of $G$ and set $(\sigma,E)=\ind_K^G 1$ (the compact
induction from $K$ to $G$ of the trivial representation of $K$). By hypothesis, $(\sigma,E)$ decomposes
into a compact component and a non-compact component
\[
E=E_c\oplus E_{nc}.
\]
Let $(\pi,V)$ be an irreducible compact representation of $G$ such that $V^K \neq 0$. By Frobenius
reciprocity for compact induction (\ref{indres})
\[
\Hom_G( \ind_K^G 1, V) \simeq    \Hom_K( 1, \res_K^G V) \simeq V^K
\]
and $\pi$ is therefore a quotient of $\ind_K^G 1$. Since $\pi$ is projective (Proposition
\ref{conscomp}), $\pi$ is in fact realized as a subrepresentation of $E=\ind_K^G 1$, and it is by
definition in the compact part $E_c$. We deduce that the irreducible compact representations $(\pi,V)$
such that $V^K \neq 0$ appear as subrepresentations of $E_c$. But $E_c$ is finitely generated because it
is a quotient of
\[
(\sigma,E)=\ind_K^G 1 = P_K^G (1)= \caH(G)\otimes_{\caH(K)} \bbC
\]
(see \ref{indres}) which is finitely generated (and even cyclic, generated by $e_K\otimes 1$). A finitely
generated semisimple representation is of finite length, which completes the proof.
\end{proof}

\section{Unitary representations}

\subsection{Hermitian representations} \label{herm}

If $V$ is a vector space over $\bbC$, let $V_\bbR$ denote the underlying real vector space and
$\overline{ V}$ the complex vector space obtained from $V$ by redefining scalar multiplication: the
scalar multiplication of a vector $v$ of $V$ by a complex number $\lambda\in \bbC$ is defined as
$\bar \lambda v$. In particular, $\overline{ V}_\bbR=V_\bbR$. If $(\pi,V)$ is a smooth representation of
$G$, we define in an obvious way a representation $( \bar\pi,\overline{ V})$ by
$\bar \pi(g)\cdot v= \pi(g)\cdot v$, $g\in G$, $v\in V$. We also set
$(\pi^h,V^h)=(\bar{\tilde \pi},\bar{\widetilde V})$. We note that we also have
$(\pi^h,V^h)=(\tilde{\bar \pi},\widetilde{\overline{ V}})$.

\begin{defi}
We say that the representation $(\pi,V)$ is Hermitian \indexter{representation!Hermitian} if there exists
a non-degenerate Hermitian form $\bilo$ on $V$ such that for all $v, w \in V$, for all $g \in G$
\[
\bil{v}{w}=  \bil{\pi(g)\cdot v}{\pi(g)\cdot w}.
\]
\end{defi}

\begin{rmqs}
1. If $(\pi,V)$ is irreducible and admissible, and if $(\pi^h,V^h)$ is isomorphic to $(\pi,V)$, then it
is Hermitian. Indeed, the canonical duality between $V$ and $V^h$ gives the existence of a non-degenerate
$G$-invariant sesquilinear form $B$ on $V$. By Proposition \ref{bilcan}, the form
$(v,w)\mapsto \overline{B(w,v)}$ is equal to $cB$ for some $c\in\bbC$, necessarily of modulus $1$; after
multiplying $B$ by a square root of $c$, we obtain a Hermitian form. Conversely, if $(\pi,V)$ is
admissible, the Hermitian inner product being non-degenerate, it induces an injection of $\overline{ V}$
into $\widetilde V$. By passing to the dual, we obtain a surjection of $V=\widetilde{{\widetilde V}}$
onto $(\overline{ V})^\sim$ (note that this is where we use the hypothesis that $V$ is admissible). We
deduce by conjugation a surjection of $\overline{ V}$ onto $\widetilde V$, which is equal to the
injection $\overline{ V} \hookrightarrow \widetilde V$ obtained previously. All these morphisms being
$G$-equivariant, we obtain $\overline{ V}\simeq \widetilde V$ as a representation of $G$, hence
$(\pi^h,V^h)\simeq (\pi,V)$.

2. If $(\pi,V)$ is irreducible admissible Hermitian, any other non-degenerate $G$-invariant Hermitian
form on $V$ is equal to $\bilo$ up to multiplication by a non-zero scalar, according to Proposition
\ref{bilcan}.
\end{rmqs}

Let $H$ be a closed subgroup of the t.d. group $G$ and let $(\tau,E)$ be a smooth representation of $H$.
We assume that $H\backslash G$ is compact, so that $\Ind_H^G(\tau,E)=\ind_H^G(\tau,E)$ (see
\ref{IndAdm}). We would like the induction functor to preserve Hermitian representations, but this is not
the case, as the following calculation shows: according to Theorem \ref{inddualthm},
\[
\ind_H^G(\tilde \tau)=\ind_H^G(\tau\otimes \delta_{H\backslash G}) ^\sim.
\]
Since it is clear that for any smooth representation $\rho$ of $H$,
$\overline{\ind_H^G \rho}=\ind_H^G\bar \rho$, we obtain
\[
\ind_H^G(\tau^h)=\ind_H^G(\tau\otimes \delta_{H\backslash G})^h.
\]

To obtain an induction that preserves Hermitian representations, we must normalize it by the factor
$\delta_{H\backslash G}^{1/2}$. Since $\delta_{H\backslash G}$ takes values in $\bbR^\times_+$, this is
defined unambiguously. We therefore set:
\begin{equation}\label{iHG}
 i_H^G\, \tau = \ind_H^G (\tau\otimes
  \delta_{H\backslash G}^{1/2}).
\end{equation}
Noting that $\tilde \delta_{H\backslash G}=\delta_{H\backslash G}^{-1}$, the above calculation
immediately shows that with this new definition of induction we now have
\begin{equation}\label{UiHG}
 (i_H^G\, \tau)^\sim=i_H^G(\tilde \tau),\quad
 (i_H^G\, \tau)^h=i_H^G(\tau^h).
\end{equation}

\subsection{Unitary representations}\label{repsunit}
\index[ter]{unitary!(representation)}

The unitary representations of a locally compact topological group $G$ are the continuous representations
of $G$ in a Hilbert space which preserve the Hermitian inner product. In general, for a t.d. group such
representations are not smooth. In a way, smoothness is incompatible with the completeness of the
representation space. We wish here to stay within the framework of smooth representations, which pushes
us to slightly redefine the notion of unitary representation. The results cited in \cite{BeZe1}, 4.21
show that the nuance is ultimately not very important.

\begin{defi}
Let $(\pi,V)$ be a smooth representation of the group $G$. It is said to be unitary
\indexter{representation!unitary} if the complex vector space $V$ is equipped with a $G$-invariant
(positive definite) Hermitian inner product $\bilo_V$, i.e.,
\[
\langle \pi(g)\cdot u ,\,  \pi(g)\cdot v  \rangle_V= \langle u ,\,  v  \rangle_V,
\quad (u,v\in V),\; (g\in G).
\]
\end{defi}

It follows easily from Proposition \ref{bilcan} that if $(\pi,V)$ is irreducible admissible unitary, a
$G$-invariant positive definite Hermitian inner product on $V$ is equal to $\bilo_V$ up to multiplication
by a strictly positive real number.

The admissibility hypothesis will allow us to recover the usual properties of unitary representations in
Hilbert spaces.

\begin{prop}
Let $(\pi,V)$ be an admissible unitary smooth representation of the t.d. group $G$ for the Hermitian
inner product $\langle .\, ,\,. \rangle_V$. If $V_1$ is a $G$-stable subspace of $V$ then
\[
V_1^\perp:=\{ v\in V \mid \forall v_1 \in V_1,  \langle v,\,v_1 \rangle_V =0\}
\]
is also a $G$-stable subspace of $V$ and
\[
V=V_1 \oplus  V_1^\perp
\]
\end{prop}

\begin{proof}
The only non-trivial part of the proposition is to show that if $v \in V$, then $v \in V_1 +  V_1^\perp$.
We can find a compact open subgroup $K$ of $G$ such that $v\in V^K$. Note that $V_1^K= V^K \cap V_1$, and
let $W$ be the orthogonal of $V_1^K$ in $V^K$ for the restriction $\bilo_{V^K}$ of the Hermitian inner
product to $V^K$. Since $V^K$ is finite-dimensional because $(\pi,V)$ is admissible, we have:
\[
V^K=V_1^K \oplus W.
\]
It remains to show that $W$ is included in $V_1^\perp$. Suppose that this is not the case. Then there
exists $w \in W$, $v_1 \in V_1$ such that $\langle  w,\,   v_1  \rangle_V \neq 0$. Now:
\[
\langle  w ,\,  v_1  \rangle_V =\langle \pi(e_K)\cdot w ,\,  v_1  \rangle_V=
\langle w ,\, \pi(e_K) \cdot v_1  \rangle_V =\langle w ,\, \pi(e_K) \cdot v_1
\rangle_{V^K},
\]
and this cannot be non-zero according to the definition of $W$.
\end{proof}

\begin{cor}
Suppose that $G$ admits a countable neighborhood basis of the identity. Let $(\pi,V)$ be an admissible
unitary representation of $G$.

$(i)$ $\dim \End_G (V)=1$ if and only if $(\pi,V)$ is irreducible.

$(ii)$ $(\pi,V)$ is completely reducible (semisimple).
\end{cor}

\begin{proof}
The first point is clear: if $(\pi,V)$ is irreducible, then $\dim \End_G (V)=1$ by Schur's lemma (see
\ref{schur}). Conversely, if $(\pi,V)$ is not irreducible, there exists a subrepresentation $V_1$ of $V$
and from the above $V=V_1\oplus V_1^\perp$ is a decomposition of $V$ into subrepresentations. The
projections onto $V_1$ and $V_1^\perp$ generate a two-dimensional subspace of $\End_G (V)$.

The second point follows from Lemma \ref{lemmesemsimP}.
\end{proof}

\subsection{Induction of unitary representations}\label{unitaires}

Let $H$ be a closed subgroup of $G$ such that $H\backslash G$ is compact. If $(\tau,E)$ is a unitary
representation of $H$, we know from the previous section that $i_H^G\, \tau$ is Hermitian. Let us
explicitly describe a $G$-invariant sesquilinear form on $i_H^G\, \tau$ inducing this Hermitian
structure. We will then show that it is positive definite, and therefore that $i_H^G\, \tau$ is unitary.

We set, for any pair $(f_1,f_2)$ of functions in $i_H^G\, E$,
\begin{equation}\label{prodscal}
 \bil{f_1}{f_2}=\int_{H\backslash G}
\bil{f_1(g_x)}{f_2(g_x)}_{E} \; d\nu_{H\backslash G}(x),
\end{equation}
where $g_x$ is a representative in $G$ of the class $x\in H\backslash G$ (more precisely, the integral
denotes the invariant linear form of \ref{invmesquo} applied to the function (\ref{gf1f2})). We must show
that this is well-defined, i.e., that the function
\begin{equation}\label{gf1f2}
g \mapsto \bil{f_1(g)}{f_2(g)}_{E}
\end{equation}
is indeed in $\scrD(G,H,  \delta_{H\backslash G})$ (see \ref{invmesquo}).

Now,
\[
\bil{f_1(hg)}{f_2(hg)}_E= \bil{\delta_{H\backslash G}^{\frac{1}{2}}(h)\tau(h)
\cdot  f_1(g)}{\delta_{H\backslash G}^{\frac{1}{2}}(h)\tau(h)\cdot f_2(g)}_{E}
\]
and since $\bilo_{E}$ is $H$-invariant, and $\delta_{H\backslash G}$ takes values in $\bbR$, we obtain
\[
\bil{f_1(hg)}{f_2(hg)}_E=\delta_{H\backslash G}(h)\bil{f_1(g)}{f_2(g)}_E.
\]
The smoothness and support properties being clear, we see that (\ref{gf1f2}) is in
$\scrD(G,H,  \delta_{H\backslash G})$. We immediately verify that $\bilo$ defines a sesquilinear form on
$i_H^G\, E$ and that this form is positive. On the other hand, if $\bil{f}{f}=0$, then since the integral
reduces to a finite sum, we see that this implies $f=0$, and therefore $\bilo$ is indeed a Hermitian
inner product.

\section{Square-integrable representations} \label{carreint}

\subsection{The space $L^2(G,\chi,dg^*)_0$}\label{L2chi}

We assume that $G$ is unimodular. Since $Z(G)$ is central, conjugation by an element of $G$ is trivial on
it, and the modulus of $G$ factors through that of $G/Z(G)$; the quotient group $G/Z(G)$ is therefore
also unimodular. Let $dg^*$ denote a bi-invariant Haar measure on $G/Z(G)$.

\begin{defi}
Let $\chi$ be a unitary character of $Z(G)$ and let $L^2(G,\chi,dg^*)_0$ be the space of functions $f$ in
$\scrC^\infty(G)$ such that

$a)$ $f(gz)=\chi(z)f(g), \quad (g\in G), \, (z\in Z(G))$

$b)$ $\int_{G/Z(G)}|f(g)|^2 \; dg^* < \infty$.

\medskip

Then $L^2(G,\chi,dg^*)_0$ is a $G$-stable subspace of $\scrC^\infty(G)$ for the right regular
representation, and
\[
\langle f_1,f_2  \rangle_{L^2,\chi}=  \int_{G/Z(G)} \overline{f_1(g)} f_2(g)
\; dg^*
\]
defines on $(L^2(G,\chi,dg^*)_0,r)$ a positive definite invariant Hermitian inner product. The
representation $(L^2(G,\chi,dg^*)_0,r)$ of $G$ is therefore unitary.
\end{defi}

\begin{rmq}
The space $L^2(G,\chi,dg^*)_0$ is the space of smooth vectors in its Hilbert space completion
$L^2(G,\chi,dg^*)$.
\end{rmq}

\subsection{Admissibility and unitarity}\label{unit2int}

We retain the hypotheses of the previous paragraph.

\begin{defi}
Let $(\pi,V)$ be a smooth representation of $G$. We say that $(\pi,V)$ is square-integrable modulo the
center, or that $(\pi,V)$ belongs to the discrete series if
\indexter{representation!square-integrable}\indexter{discrete series}

$a)$ The center $Z(G)$ of $G$ acts on $V$ by a unitary character $\chi$.

$b)$ Any matrix coefficient $\phi_{v,\lambda}$ of $(\pi,V)$ is in $L^2(G,\chi,dg^*)_0$.

If there exists a character $\omega$ of $G$ such that $(\pi\omega,V)$ is square-integrable modulo the
center, we say that $(\pi,V)$ is essentially square-integrable modulo the center.
\end{defi}

\begin{rmq}
Let $(\pi,V)$ be an irreducible smooth representation of $G$, admitting a unitary central character.
Then, for $(\pi,V)$ to be square-integrable modulo the center, it suffices that one of the non-zero
matrix coefficients of $\pi$ be square-integrable modulo the center. This is of course false if we do not
assume $(\pi,V)$ is irreducible. Since we do not use this result, we leave the (easy) proof to the
reader.
\end{rmq}

\begin{prop}
Let $(\pi,V)$ be an irreducible smooth representation, essentially square-integrable modulo the center of
$G$. Then $(\pi,V)$ is admissible.
\end{prop}

\begin{proof}
According to Lemma \ref{tensprod}, we can assume that $(\pi,V)$ is a smooth representation
square-integrable modulo the center. We retain the idea of the proof of Theorem \ref{comprep}. Suppose
that $(\pi,V)$ is not admissible. Then there exists a compact open subgroup $K$ of $G$, a non-zero vector
$v$ in $V^K$ and a sequence $(g_n)_{n\in \bbN}$ of elements of $G$ such that
$\{\pi(e_K)\pi(g_n)\cdot v\}$ is linearly independent in $V^K$ (indeed, $V^K$ is generated by the
$\{\pi(e_K)\pi(g)\cdot v\}$ according to \ref{piK}). Note then that the double cosets $Kg_nZ(G)$ are
distinct in $K\backslash G/Z(G)$. Let us choose $\lambda \in \widetilde V^K$ such that
$\lambda(\pi(e_K)\pi(g_n)\cdot v)=1$ for all $n$. Let $\overline{ K}$ be the image of $K$ in $G/Z(G)$, so
that $\overline{ K}$ is a compact open subgroup of $G/Z(G)$. We then have
\begin{align*}
\int_{G/Z(G)}|\phi_{v,\lambda}(g)|^2 \; dg^*
&=\int_{G/Z(G)}|\lambda(\pi(g)\cdot v)|^2 \; dg^* \\
&=\int_{G/Z(G)}|(\pi(e_K)\cdot \lambda)(\pi(g)\cdot v)|^2 \; dg^*\\
&= \int_{G/Z(G)}|\lambda(\pi(e_K)\pi(g)\cdot v)|^2 \; dg^*\\
&=\sum_{\bar g \in \overline{ K}\backslash (G/Z(G))} |  \lambda( \pi(e_K)\pi(g)\cdot v)|^2 \; \vol_{dg^*}(\overline{ K})\\
&\geq \sum_{i=0}^\infty |
  \lambda(\pi(e_K)\pi(g_i)\cdot v)|^2 \; \vol_{dg^*}(\overline{ K})\geq
  \sum_{i=0}^\infty \vol_{dg^*}(\overline{ K}) =\infty
\end{align*}
and we reach a contradiction.
\end{proof}

We will now show that an admissible representation square-integrable modulo the center is unitary.

\begin{lemme}\label{unit2intlem}
Let $(\pi,V)$ be a smooth representation of $G$, admissible and square-integrable modulo the center. Let
$\chi$ be its (unitary) central character. Then $(\pi,V)$ is unitary and semisimple.
\end{lemme}

\begin{proof}
Let $W \subset V$ be a finitely generated subrepresentation, generated by vectors $\{w_1,\ldots ,w_n\}$.
Let $\mathrm{Ann}(W)$ denote the annihilator of $W$ in $\widetilde V$, so that
$\widetilde W \simeq \widetilde V/\mathrm{Ann}(W)$. Let us set, for
$\lambda_1,\lambda_2 \in \widetilde W $,
\[
(\lambda_1,\lambda_2)=  \sum_{i=1, \ldots,n} \int_{G/Z(G)}
\overline{{\lambda_1} (\pi(g)\cdot w_i)} \; {\lambda_2}(\pi(g)\cdot  w_i) \; dg^*.
\]
This defines a $G$-invariant positive definite Hermitian form on $\widetilde W$ (which depends on the
choice of generators $\{w_1,\ldots ,w_n\}$ and is therefore not canonical). We deduce that $\widetilde W$
is unitary, and therefore that $W$ is unitary (indeed, $W$ being admissible,
$W\simeq\widetilde{\widetilde W}$, and the contragredient of an admissible unitary representation $U$ is
isomorphic to $\overline U$, which is unitary). According to Corollary \ref{repsunit}, $W$ is then
semisimple. Since $V$ is the union of its finitely generated subrepresentations, Lemma \ref{lemmesemsimP}
tells us that $V$ is semisimple. Let us then write $V=\bigoplus_i V_i$, with irreducible $V_i$; each
$V_i$, being finitely generated, is unitary, and the orthogonal sum of the Hermitian inner products of
the $V_i$ makes $V$ a unitary representation.
\end{proof}

\subsection{Schur orthogonality}\label{orthschur}
\indexter{Schur orthogonality}

Recall that for any function $f$ on a group $G$, the function $\check f$ is defined by
$\check f(g)=f(g^{-1})$.

\begin{lemme}
$(i)$ Let $(\pi,V)$ be an irreducible representation of $G$, essentially square-integrable modulo the
center. Then there exists a strictly positive real number $d(\pi)$, called the formal degree of $\pi$,
such that for all $v_1,v_2\in V$, for all $\lambda_1,\lambda_2\in \widetilde V$, we have:
\[
\int_{G/Z(G)}\phi_{v_1,\lambda_1}(g)\check\phi_{v_2,\lambda_2}(g)\; dg^*=
\frac{\lambda_1(v_2)\lambda_2(v_1)}{d(\pi)}.
\]

$(ii)$ Let $(\pi_1,V_1)$, $(\pi_2,V_2)$ be two irreducible representations of $G$, essentially
square-integrable modulo the center, non-equivalent, and with the same central character. Then for all
$v_1 \in V_1$, $v_2\in V_2$, $\lambda_1 \in \widetilde V_1$, $\lambda_2 \in \widetilde V_2$,
\[
\int_{G/Z(G)}\phi_{v_1,\lambda_1}(g)\check
\phi_{v_2,\lambda_2}(g)\; dg^*=0.
\]
\end{lemme}

\begin{proof}
Let $\omega$ be a smooth character of $G$ such that $\pi\omega$ is square-integrable modulo the center.
It is easy to see that the integral we are trying to calculate is the same if we replace $\pi$ by
$\pi\omega$, the contributions due to $\omega$ canceling out. We may therefore assume $\pi$ is
square-integrable modulo the center. According to the lemma of the previous paragraph, we can equip $V$
with an invariant Hermitian inner product $\bilo_V$. Consider the representation $V\otimes \widetilde V$
of $G\times G$: according to Proposition \ref{tensprod}, it is irreducible, admissible, and its
contragredient is naturally isomorphic to $\widetilde V\otimes V $.

Consider the bilinear form $B$ on $(V\otimes \widetilde V)\times(\widetilde V\otimes V)$ defined by:
\[
B(v_1\otimes \lambda_1, \lambda_2 \otimes v_2)= \int_{G/Z(G)}\phi_{v_1,\lambda_1}(g)\check
\phi_{v_2,\lambda_2}(g)\; dg^*.
\]
It is easy to see that this bilinear form is $G\times G$-invariant because $G/Z(G)$ is unimodular.
According to Proposition \ref{bilcan}, to show that this form is non-degenerate, it suffices to show that
it is non-zero. Let us choose $v\in V$ non-zero, and let $\lambda_v\in \widetilde V$ be the linear form
on $V$ defined by $\lambda_v(w)=\langle v,w \rangle_V$, $w\in V$. Since
\begin{align*}
 \lambda_v(\pi(g^{-1})\cdot v)&= \langle v, \pi(g^{-1})\cdot v
\rangle_V= \langle \pi(g)\cdot v,v \rangle_V=  \overline{\langle v,
  \pi(g)\cdot v \rangle }_V\\
& = \overline{ \lambda_v(\pi(g)\cdot v)},
\end{align*}
we have:
\begin{align*}
B(v\otimes \lambda_v, \lambda_v \otimes v)&= \int_{G/Z(G)}
  \lambda_v(\pi(g)\cdot v) \lambda_v(\pi(g^{-1})\cdot v)\; dg^*\\
&=\int_{G/Z(G)} |\lambda_v(\pi(g)\cdot v)|^2 \; dg^* >0
\end{align*}
According to Proposition \ref{bilcan}, there exists a non-zero complex constant $c$ such that the
bilinear form $B$ on $(V\otimes \widetilde V)\times(\widetilde V\otimes V)$ is equal to $c$ times the
canonical bilinear form on
$(V\otimes \widetilde V)\times(\widetilde V\otimes V)\simeq(V\otimes \widetilde V)\times \widetilde {(V\otimes \widetilde V)} $.
This is given by
\[
(v_1 \otimes \lambda_1,\lambda_2\otimes v_2)_0=\lambda_1(v_2)\lambda_2(v_1).
\]
Since
\[
(v\otimes \lambda_v, \lambda_v \otimes v)_0= \lambda_v(v)^2= \langle
v,v \rangle^2>0,
\]
the constant $c$ is real, strictly positive. We therefore set $d(\pi)=c^{-1}$.

Let us prove $(ii)$. We fix $\lambda_1 \in \widetilde  V_1$ and $v_2 \in V_2$. We define an intertwining
operator between $V_1$ and $(\widetilde V_2)^\sim \simeq V_2$ by:
\[
v_1 \mapsto \left( \lambda_2 \mapsto \int_{G/Z(G)}\phi_{v_1,\lambda_1}(g)\check
\phi_{v_2,\lambda_2}(g)\; dg^* \right).
\]
Since $V_1$ and $V_2$ are non-equivalent, this intertwining operator is zero.
\end{proof}

\begin{cor}
Let $(\pi_i,V_i)$, $i=1,\ldots,r$, be irreducible representations of $G$, essentially square-integrable
modulo the center, pairwise non-equivalent, and with the same central character. For each $i$, let us
choose $v_i \in V_i$ and $\lambda_i \in \widetilde V_i$ non-zero. Then the matrix coefficients
$\phi_{v_i,\lambda_i}$ are linearly independent. The same conclusion holds if we consider irreducible
compact representations.
\end{cor}

\begin{proof}
Suppose that $\sum_i c_i \phi_{v_i,\lambda_i}=0$. Let us fix $j$ between $1$ and $r$, and choose
$\lambda'_j \in \widetilde V_j$ such that $\lambda'_j(v_j)=1$ and $v'_j \in V_j$ such that
$\lambda_j(v'_j)=1$. Then according to $(ii)$ of the lemma, we have
\begin{align*}
0&=\int_{G/Z(G)}  \left(\sum_i c_i \phi_{v_i,\lambda_i}\right) \check
\phi_{v'_j,\lambda'_j}\;  dg^* = c_j \; \int_{G/Z(G)} \phi_{v_j,\lambda_j}        \check
\phi_{v'_j,\lambda'_j}\;  dg^*\\
&= c_j \, d(\pi_j)^{-1}.
\end{align*}
Therefore $c_j=0$. If the representations are compact, the center of $G$ is compact (Lemma
\ref{comprep}). There is then no need to quotient by the center and we adapt the previous proofs.
\end{proof}

\begin{rmq}
More generally, if $A$ is a subgroup of $Z(G)$ such that $Z(G)/A$ is compact, we can replace the
integrals over $G/Z(G)$ by integrals over $G/A$ without affecting their convergence, and the results
obtained above are still valid.
\end{rmq}

\section{Notes on Chapter IV}

The results of this chapter are found in \cite{BeZe1}. In writing this chapter, I also used the notes
\cite{DeB} and \cite{Ro}. For example, the condition $(\mathbf{KF})$ comes from \cite{Ro}. The proof of
Lemma \ref{unit2intlem} is taken from \cite{Wald}.

\chapter[Structure of $p$-adic groups]{Structure of $p$-adic reductive groups}\label{chapstruct}

This chapter is a catalog of results, not all proved here, but referenced as much as possible. We
introduce our main object of study, $p$-adic reductive groups. Our exposition, hitherto roughly
``self-contained'', ceases to be so, for it is out of the question to treat the theory of algebraic
groups, or even reductive algebraic groups, in this book, especially since this is very well done in
numerous works (\cite{Bo}, \cite{Hum}, \cite{Spr}, \cite{BoTi}). Similarly, we do not feel courageous (or
competent) enough to present Bruhat-Tits theory (\cite{BrTi}, \cite{BrTi2}, \cite{Tits}). The class of
groups we study is that of the groups of rational points of a connected reductive algebraic group defined
over a $p$-adic field. Bruhat-Tits theory constructs the building of such groups, and the action of the
group on its building. We cannot recommend strongly enough that the reader complete this brief
presentation by studying two or three concrete examples of classical groups, starting with the group
$GL(n)$. These examples are developed in the literature we indicate.

\section{Non-Archimedean local fields}

We refer the reader to \cite{CaFr} for more details concerning non-Archimedean local fields, as well as
for the proofs of the results stated here.

Let $\bbF$ \index[not]{F@$\bbF$} be a non-Archimedean local field
\index[ter]{non-Archimedean local field} equipped with its discrete valuation \index[ter]{valuation}
$v_\bbF$ \index[not]{v_F@$v_\bbF$} taking values in $\bbZ \cup \{ \infty \} $, normalized so that the
image of $v_\bbF$ is $\bbZ \cup \{ \infty\}$. Let $\frO_\bbF$ \index[not]{O_P@$\frO_\bbF$} denote the
ring of integers of $\bbF$,
\[
\frO_\bbF=\{ x\in \bbF \mid \, v_\bbF(x)\geq 0 \}
\]
and $\frP_\bbF$ \index[not]{P_F@$\frP_\bbF$} its unique maximal ideal ($\frO_\bbF$ is a local ring),
\[
\frP_\bbF= \{ x\in \bbF \mid \, v_\bbF(x) > 0 \} .
\]
Let $\frO_\bbF^\times$ denote the group of units of $\frO_\bbF$. It is the set of elements of $\frO_\bbF$
of valuation zero. Let us fix a uniformizer \index[ter]{uniformizer} $\varpi$
\index[not]{ZZpivar@$\varpi$} of $\bbF$, i.e., an element satisfying $v_\bbF(\varpi)=1$. We can then
write any $x \in \bbF^\times$ uniquely in the form $x=u\varpi^n$ where $n\in \bbZ$ and
$u \in \frO_\bbF^\times$.

Let $k$ be the residue field $\frO_\bbF/\frP_\bbF$. It is a finite field, hence isomorphic to $\bbF_q$,
for some power $q$ of a prime number $p$. The normalized absolute value on $\bbF$ is given by
\[
|x|_\bbF= q^{-v_\bbF(x)}, \quad (x \in \bbF)
\]
\index[not]{abs@$\vert \, .\, \vert _\bbF$} This absolute value equips $\bbF$ with a metric space
structure, making it a totally disconnected topological group (for addition). A neighborhood basis of $0$
is given by the compact open subgroups $\frP_\bbF^n$.

The additive group $(\bbF,+)$ is equipped with a Haar measure $dx$, which can be normalized for example
by setting the volume of the compact set $\frO_\bbF$ equal to 1. Absolute value and Haar measure are
related by the formula:
\[
\int_\bbF f(ax) \, dx =  |a|^{-1}_\bbF\,  \int_\bbF f(x) \, dx, \quad (a
\in \bbF ^\times), (f \in \scrC_c^\infty(\bbF)).
\]

The algebraic closure of $\bbF$ is denoted by $\overline{\bbF}$ and its separable closure by $\bbF^s$.
The one-variable multiplicative group $\mathbf{G}_m$ is the algebraic group defined over $\bbF$ whose
group of points over $\bbF$ is identified with $\bbF^\times$.

Any finite-dimensional vector space over $\bbF$ is equipped with the analytic topology. This is the case
in particular for the algebra $M_n(\bbF)$ of square matrices of size $n$ with coefficients in $\bbF$. The
group $GL_n(\bbF)$ of invertible matrices in $M_n(\bbF)$ is equipped with the induced topology, making it
a t.d. group. A neighborhood basis of the identity is given by the compact open subgroups
\[
K_n= \Id_n + \frP_\bbF^n \; M_n(\frO_\bbF).
\]

\section{$p$-adic reductive groups}

\subsection{Linear algebraic groups}
We recall some elements of the theory of linear algebraic groups. We refer the reader to \cite{Spr} for
the notions and proofs of the statements recalled here. We fix a non-Archimedean local field $\bbF$ as in
the previous section.

Let $\bbG$ \index[not]{G@$\bbG$} be a connected linear algebraic group defined over $\bbF$ and let $G$
denote the group of its points over $\bbF$. We identify $\bbG$ and the group of its points over
$\overline{\bbF}$. We denote respectively by $R(\bbG)$ \index[not]{R(G)@$R(\bbG)$} and $R_u(\bbG)$
\index[not]{R_u(G)@$R_u(\bbG)$} the radical and the unipotent radical \index[ter]{radical}
\index[ter]{unipotent radical} of $\bbG$. The group $R(\bbG)$ is a connected solvable normal algebraic
subgroup of $\bbG$, and maximal for this property. It is unique. Similarly $R_u(\bbG)$ is the connected
unipotent normal algebraic subgroup of $\bbG$ and maximal for this property. If $R_u(\bbG)$ is trivial,
the group $\bbG$ is said to be reductive. If $R(\bbG)$ is trivial, the group $\bbG$ is said to be
semisimple. Let $\scrD\bbG$ \index[not]{D(G)@$\scrD\bbG$} be the derived group \index[ter]{derived group}
of $\bbG$. It is always defined over $\bbF$. If $\bbG$ is reductive, then $R(\bbG)$ is central in $\bbG$,
$\scrD\bbG$ is semisimple, $\bbG$ is generated by $ R(\bbG)$ and $\scrD\bbG$, and
$R(\bbG) \cap \scrD\bbG$ is finite. Let $Z(\bbG)$ \index[not]{Z(G)@$Z(\bbG)$} denote the center of
$\bbG$. If $\bbG$ is reductive, $R(\bbG)$ is the identity component of $Z(\bbG)$ and the subgroups
$R(\bbG)$ and $Z(\bbG)$ are defined over $\bbF$. We denote respectively by $R(G)$ and $Z(G)$ their points
over $\bbF$. \index[not]{R(G)@$R(G)$} \index[not]{Z(G)@$Z(G)$}

If $\bbG$ is reductive, a maximal torus in $\bbG$ is its own centralizer. There exist maximal tori of
$\bbG$ defined over $\bbF$ and if $\bbF$ is of characteristic zero, the set of conjugacy classes under
$G$ of maximal tori of $\bbG$ defined over $\bbF$ is finite. All tori of $\bbG$ defined and split over
$\bbF$ and maximal for this property are conjugate under $G$.

A split component \index[ter]{split component} of $\bbG$ is a torus defined and split over $\bbF$,
contained in the radical $R(\bbG)$ and maximal for this property. If $\bbG$ is reductive, a split
component of $\bbG$ is then a maximal split torus in the center of $\bbG$. Since two such tori must be
conjugate by $R(G)$ which is central, we deduce the uniqueness of the split component for a reductive
group.

The group $\bbG$ can be algebraically embedded into a group $GL_N(\overline{\bbF})$ for some $N$, the
embedding being defined over $\bbF$. The group $G$ therefore embeds into $GL_N(\bbF)$. Equipped with the
induced topology, $G$ becomes a totally disconnected topological group, for which the notions and results
of Chapters \ref{chaptd} and \ref{chapreptd} therefore apply. This topology is independent of the choice
of the embedding into a $GL_N(\bbF)$. Moreover, if $\bbG$ is reductive, $G$ is unimodular.

\subsection{Rational characters}\label{carrat}
\index[ter]{character!rational} Let $\bbG$ be a connected reductive algebraic group defined over $\bbF$
and let $X^*(\bbG)$ \index[not]{X(G)@$X^*(\bbG)$} be the group of its algebraic characters (i.e., the
group of algebraic morphisms defined over $\overline{\bbF}$ from $\bbG$ to the multiplicative group
$ \mathbf{G}_m$). We denote by $X^*(G)$ \index[not]{X(G)1@$X^*(G)$} the subgroup of $X^*(\bbG)$ of
characters defined over $\bbF$, and we call it the group of rational characters of $\bbG$ (or of $G$).
Since the restriction of such an algebraic character to $\scrD\bbG$ is trivial because $ \mathbf{G}_m$ is
commutative, $X^*(\bbG)$ is identified with the group of algebraic characters of $\bbG/\scrD\bbG$. Since
$\bbG=R(\bbG)\scrD\bbG$, we have
\[
\bbG/\scrD\bbG \simeq R(\bbG)/(R(\bbG)\cap \scrD\bbG).
\]
Now $R(\bbG)$ is a torus, and the group of its algebraic characters $X^*(R(\bbG))$ is a finitely
generated, torsion-free abelian group (a lattice, hence isomorphic to $\bbZ^r$ for some $r\in \bbN$). The
group $R(\bbG)\cap \scrD\bbG$ being finite, $X^*(\bbG)$ is a full-rank sublattice of $X^*(R(\bbG))$. The
group of rational characters $X^*(G)$ of $G$ is a sublattice, possibly trivial, of the lattice
$X^*(\bbG)$.

Let $X_*(\bbG)$ \index[not]{X(G)2@$X_*(\bbG)$} be the dual group,
$X_*(\bbG):= \Hom_\bbZ(X^*(\bbG),\bbZ)$. Let $\langle \cdot, \cdot \rangle$ denote the natural duality
between $X^*(\bbG)$ and $X_*(\bbG)$ and note that this duality induces an isomorphism
$X^*(\bbG)\simeq \Hom_\bbZ(X_*(\bbG),\bbZ)$.

We also denote by $X_*(G)$ \index[not]{X(G)3@$X_*(G)$} the lattice $\Hom_\bbZ(X^*(G),\bbZ)$. Similarly
$X^*(G)\simeq \Hom_\bbZ(X_*(G),\bbZ)$.
\medskip

\begin{rmqs}
{ \sl Case of tori. }

1. If $\bbT$ is an algebraic torus defined over $\bbF$, we can moreover identify $X_*(\bbT)$ with the
group of morphisms from $\mathbf{G}_m$ to $\bbT$ so that if $t \in \mathbf{G}_m$, $\phi \in X_*(\bbT)$
and $\chi \in X^*(\bbT)$,
\[
\chi(\phi(t))= t^{\langle \chi , \phi  \rangle}.
\]

The elements of $X^*(\bbT)$ form a basis of the vector space of polynomial functions on $\bbT$. The group
algebra $\overline{\bbF}[X^*(\bbT)]$ is the algebra of polynomial functions on $\bbT$, and $\bbT$ is
identified with the spectrum of this algebra.

2. A torus $\bbT$ defined over $\bbF$ is said to be anisotropic if $X^*(T)=\{0\}$ and split if
$X^*(T)=X^*(\bbT)$. Any torus $\bbT$ defined over $\bbF$ admits a unique maximal split subtorus $\bbA$,
and a unique maximal anisotropic subtorus $\bbT_{an}$ and $\bbT$ is the almost direct product
\index[ter]{almost direct product} of $\bbA$ and $\bbT_{an}$ (this means that they generate $\bbT$ and
that their intersection is finite).

\end{rmqs}

\subsection{Unramified characters}\label{carnonram}
Let $\bbG$ be a connected reductive algebraic group defined over $\bbF$. Let us define a morphism
\footnote{The reader should note that the group $X_*(G)$ being a lattice, its operation is written
additively, whereas the operation of $G$ is written multiplicatively.}
\begin{equation}
H_G: G \rightarrow X_*(G)
\end{equation}
\label{v} \index[not]{H_G@$H_G$} by the formula:
\[
\langle \, \chi, H_G(g) \,    \rangle= v_\bbF(\chi(g)), \quad g\in G,
\chi \in X^*(G).
\]
Let ${}^0 G$ \index[not]{G0@${}^0 G$} be the kernel of this morphism, and $\Lambda(G)$
\index[not]{ZZLambda(G)@$\Lambda(G)$} its image. We therefore have an exact sequence:
\[
1 \rightarrow {}^0 G \longrightarrow G \stackrel{H_G}{\longrightarrow}
\Lambda(G)\rightarrow 1
\]
It is immediate that an equivalent definition of ${}^0 G$ is given by
\[
{}^0 G:=\bigcap_{\chi \in  {X}^*(G)} \ker |\chi|_\bbF .
\]

\begin{lemme}
Let $\chi \in X^*(G)$. Then $|\chi|_\bbF\colon G \rightarrow \bbR^\times_+$ is a smooth character of $G$
whose kernel contains all compact subgroups of $G$.
\end{lemme}
\begin{proof}
The kernel of $|\chi|_\bbF$ is the inverse image under $\chi$ of the compact open set $\frO_\bbF^\times$
of $\bbF^\times$, so it is an open set of $G$. Consequently, $|\chi|_\bbF$ is a smooth character of $G$,
taking values in $\bbR^\times_+$. In particular its kernel contains compact open subgroups as small as we
want.

For any compact subgroup $K$ of $G$, $|\chi|_\bbF(K)$ is a compact subgroup of $\bbR^\times_+$, and the
only such subgroup is $\{1\}$.
\end{proof}

\begin{prop}
The group ${}^0G$ is an open, closed and normal unimodular subgroup of $G$. Any compact subgroup of $G$
is contained in ${}^0 G$. The derived group $\scrD G$ is contained in ${}^0 G$.
\end{prop}
\begin{proof}
According to the lemma, any compact subgroup of $G$ is in ${}^0 G$. On the other hand, $\ker |\chi|_\bbF$
is a normal subgroup of $G$, and the intersection of normal subgroups is again normal, so ${}^0 G$ is
normal in $G$. Since the characters $|\chi|_\bbF$ are smooth, the $\ker |\chi|_\bbF$ are closed in $G$.
It follows that ${}^0 G$ is closed in $G$. Since ${}^0 G$ contains compact open subgroups, it must be
open in $G$. It is therefore unimodular because $G$ is. The last assertion follows from the fact that any
character is trivial on the derived group.
\end{proof}

\begin{defi}
An unramified character \index[ter]{character!unramified} of $G$ is a group morphism from $G$ to
$\bbC^\times$ trivial on ${}^0G$. We denote by $\caX(G)$ the set of unramified characters of $G$.
\end{defi}

\subsection{Variety structure of $\caX(G)$} \label{varXG}
An unramified character of $G$ is of the form $u \circ H_G$ where $u$ is a morphism from $\Lambda(G)$ to
$\bbC^\times$. The unramified characters of $G$ are thus identified with the characters of $\Lambda(G)$,
i.e., with the elements of the group $\Hom_\bbZ(\Lambda(G),\bbC^\times)$. Note that $\Lambda(G)$, a
subgroup of $X_*(G)$, is a finitely generated free abelian group. Let us introduce the group algebra
$\bbC[\Lambda(G)]$. We have a \index[not]{CLambda(G)@$\bbC[\Lambda(G)]$} natural isomorphism:
\[
\Hom_\bbZ(\Lambda(G),\bbC^\times)\simeq \Hom_{alg}(\bbC[\Lambda(G)], \bbC),
\]
and the above objects are identified with the maximal spectrum of $\bbC[\Lambda(G)]$, a complex algebraic
torus denoted by $\bbU$. We denote by $t \mapsto \omega_t$ the isomorphism between the complex torus
$\bbU$ and the group of unramified characters of $G$. This equips $\caX(G)$ with a complex torus
structure. The algebra of polynomial functions on $\caX(G)$ is explicitly identified with
$\bbC[\Lambda(G)]$ in the following way: if $f=\sum_{\bar g \in \Lambda(G)} a_{\bar g} \bar g$ is an
element of $\bbC[\Lambda(G)]$ and if $\chi \in \caX(G)$, then
$f(\chi)=\sum_{\bar g \in \Lambda(G)} a_{\bar g} \chi(g)$, where the $g\in G$ lift the
$\bar g \in \Lambda(G)$. If we identify $\bbC[\Lambda(G)]$ with the space of finitely supported functions
on the group $\Lambda(G)$, the product becomes the convolution of these functions.

For all $g \in G$, let $\mathbf{ev}_g$ \index[not]{evg@$\mathbf{ev}_g$} denote the evaluation morphism at
$g$:
\begin{equation*}
 \mathbf{ev}_g \colon  \caX(G)\rightarrow \bbC^\times, \quad
\chi \mapsto \chi(g).
\end{equation*}
It is clear that $\mathbf{ev}_g$ is a polynomial function on $\caX(G)$, which is given in the above
identification by the class $\bar g$ of $g$ in $G/{}^0G$, viewed as an element of
$\Lambda(G) \subset  \bbC[\Lambda(G)]$. We denote by $\chi_{un}$ \index[not]{ZZchi_un@$\chi_{un}$} the
morphism
\[
\chi_{un}\colon G \rightarrow \bbC[\Lambda(G)]^\times, \quad g \mapsto
\mathbf{ev}_g.
\]
Any point $\chi \in \caX(G)$ determines an algebra morphism $\Psi_\chi$ \index[not]{ZZPsichi@$\Psi_\chi$}
from $\bbC[\Lambda(G)]$ to $\bbC$, and we have the relation
\begin{equation}\label{Psichi}
\Psi_\chi \circ \chi_{un}= \chi.
\end{equation}

The group $G$ acts on $G/{}^0G$ by left translation, so it acts on $\Lambda(G)$ and on
$\bbC[\Lambda(G)]$. A simple calculation shows that for all $g \in G$, for all $f \in \bbC[\Lambda(G)]$,
and for all $\chi \in \caX(G)$,
\begin{equation} \label{Xun}
 (g\cdot f)(\chi)=\chi(g)f(\chi)=
(\chi_{un}(g)f)(\chi).
\end{equation}
The action of $G$ on $\bbC[\Lambda(G)]$ is therefore given by $\chi_{un}$. It factors through $G/{}^0G$
and we recover the natural action of $\Lambda(G)=G/{}^0G$ on its group algebra by left translation.

The variety $\caX(G)$ is a complex torus, hence a group, which acts on itself by left translation. The
group structure is also of course the one defined by multiplication of characters. The group $\caX(G)$
therefore acts on the space of its polynomial functions $F$ by
\[
(\chi\cdot f)(\psi)=f(\chi^{-1}\psi),\quad  (\chi,\psi \in \caX(G)),
\, (f \in F).
\]

\subsection{Case of tori}\label{castordep}

Let us specify the constructions of the previous paragraph in the case where $\bbG=\bbT$ is a torus over
$\bbF$.

First recall that if $\bbT$ is an algebraic torus defined over $\bbF$, if $\bbA$ is its split component,
and $\bbT_{an}$ \index[not]{Tan@$\bbT_{an}$} its anisotropic component, then $\bbT$ is the almost direct
product $\bbT=\bbT_{an} \bbA$ (cf. Remarks \ref{carrat}). The natural restriction map from $X^*(\bbT)$ to
$X^*(\bbA)$ is injective, $X^*(A)=X^*(\bbA)$ and $X^*(T_{an})=\{0\}$, and $X^*(T)$, viewed as a
sublattice of $X^*(A)$, is of finite index (it is the sublattice of characters whose restriction to
$\bbA \cap \bbT_{an}$ is trivial).

Let $\bbA$ be a split torus over $\bbF$. We then have:
\[
A\simeq \Hom_\bbZ(X^*(A),\bbF^\times)\simeq  \Hom_\bbZ(X^*(A),\bbZ)
\otimes_\bbZ \bbF^\times = X_*(A)\otimes_\bbZ \bbF^\times,
\]
the isomorphisms being natural. Let $\varpi\in \bbF^\times$ be a uniformizer of $\bbF$. The set of its
powers forms a subgroup $\bbZ_\varpi$ \index[not]{Zpi@$\bbZ_\varpi$} of $\bbF^\times$, and since by
definition $v_\bbF(\varpi)=1$, $v_\bbF$ realizes an isomorphism between $\bbZ_\varpi$ and $\bbZ$. Let us
set $$C_A=\Hom_\bbZ(X^*(A),\bbZ_\varpi)\simeq \Hom_\bbZ(X^*(A),\bbZ)=X_*(A). $$ \index[not]{C_A@$C_A$}
This identifies $C_A$ with a subgroup of $A$, $H_A$ giving the isomorphism between $C_A$ and $X_*(A)$. In
particular, $H_A$ is surjective, so $\Lambda(A)=X_*(A)$ and $C_A$ is a section of $H_A$.

Note also that if $\bbA$ is a split torus, then $A$ is isomorphic to a product of copies of $\bbF^\times$
and the group ${}^0A$ is the maximal compact subgroup of $A$. Indeed, we trivially have
${}^0\bbF^\times = \frO_\bbF^\times$, the group of units of the field $\bbF$.

If $\bbA_1$ and $\bbA_2$ are two split tori defined over $\bbF$, with $\bbA_1 \subset \bbA_2$, then
$\Lambda(A_1) \subset \Lambda (A_2)$. Indeed, from the above, we easily see that
${}^0A_1=A_1 \cap {}^0A_2$. We can also identify $C_{A_1}$ with a subgroup of $C_{A_2}$.

\subsection{Unramified characters (continued)} \label{Lambda}

We return to the hypotheses and notation of \ref{carnonram}. Let $A_G$ be a split component of $G$. The
constructions made in \ref{carnonram} and \ref{castordep} provide exact sequences
\[
1 \rightarrow {}^0 G \rightarrow G \rightarrow \Lambda(G)\rightarrow 1
\]
and
\[
1 \rightarrow {}^0 A_G \rightarrow A_G \rightarrow \Lambda(A_G)\simeq
X_*(A_G) \rightarrow 1.
\]

\begin{lemme}
We have ${}^0 A_G = {}^0 G \cap A_G$. The lattice $\Lambda(A_G)=X_*(A_G)$ embeds into
$\Lambda(G)\subset X_*(G)$. Moreover the index of $X_*(A_G)$ in $\Lambda(G)$ (and even in $X_*(G)$) is
finite (they are lattices of the same rank).
\end{lemme}

\begin{proof}
Consider the restriction map from $X^*(G)$ to $X^*(A_G)$. We claim that it is an injection. Indeed, we
saw in \ref{carrat} that $X^*(G)=X^*(G')$ where
\[
\bbG'=\bbG/\scrD\bbG=R(\bbG)/(R(\bbG)\cap \scrD \bbG).
\]
Let us set $\bbA'_G= \bbA_G/(\bbA_G \cap \scrD \bbG)$. It is clear that $\bbA'_G$ is a split component of
$\bbG'$. Since $\bbG'$ is a torus, the lattice $X^*(G)=X^*(G')$ is of finite index in $X^*(A'_G)$ and
therefore of finite index in $X^*(A_G)$. Let us take $g \in {}^0G \cap A_G$ and $\chi \in X^*(A_G)$. Then
for some $m\in \bbN^*$, $\chi^m$ is the restriction to $A_G$ of some character in $X^*(G)$. We therefore
have
\[
|\chi^m(g)|_\bbF= |\chi(g)|^m_\bbF=1.
\]
But $ |\chi(g)|_\bbF \in \bbR^\times_+$ so $ |\chi(g)|_\bbF =1$. We deduce that $g \in  {}^0 A_G$. The
reverse inclusion ${}^0 A_G\subset  {}^0G \cap A_G$ is clear.

We now have a well-defined and injective morphism:
\[
X_*(A_G)= \Lambda(A_G)\simeq A_G/{}^0A_G \rightarrow   G/{}^0G
=\Lambda(G)\subset X_*(G).
\]
Since $ X_*(A_G)$ has the same rank as $ X_*(G)$, we deduce the last assertion.
\end{proof}

\begin{prop}
The quotient $G/{}^0GA_G$ is finite, and ${}^0G \cap  Z(G)$ is compact. If $Z(G)$ is compact, then
${}^0G=G$.
\end{prop}
\begin{proof}
We have
\[
(G/{}^0G)/(A_G/{}^0A_G)\simeq  G/{}^0G A_G.
\]
The finiteness of the cardinality of the left-hand side of this equality having been established in the
lemma, we deduce that $G/{}^0G A_G$ is finite. If $Z(G)$ is compact, we have $A_G=\{1\}$, and ${}^0G$ is
then of finite index in $G$. Suppose $G \neq {}^0G$, and let $g \in G \setminus {}^0G$. Then there exists
a rational character $\chi \in X^*(G)$ such that $|\chi(g)|_\bbF \neq 1$. But some power of $g$ is in
${}^0G$, say $g^m$. We then have $|\chi(g^m)|_\bbF=1= |\chi(g)|^m_\bbF $ and we obtain a contradiction.
Therefore we have $G={}^0G$ in this case.

The group $R(G)$ is of finite index in $Z(G)$ ($R(\bbG)$ is the identity component of $Z(\bbG)$) and
$R(G)_{an}  A_G$ is cocompact in $R(G)$ ($R(\bbG)$ is the almost direct product of $R(\bbG)_{an}$ and
$\bbA_G$), hence in $Z(G)$, where $R(G)_{an}$ is compact. Therefore $R(G)_{an}\subset {}^0G$ and since
$A_G\cap {}^0G={}^0A_G$, we have $R(G)_{an}A_G\cap {}^0G=  R(G)_{an}  {}^0A_G$, which is compact.
Consequently ${}^0G \cap  Z(G)$ is compact.
\end{proof}

\subsection{Inertial classes of irreducible representations}\label{finXGpi}

The group of unramified characters $\caX(G)$ acts on the set of equivalence classes of irreducible smooth
representations of $G$, $\mathbf{Irr}(G)$ by
\[
(\omega,\pi)\mapsto \pi\otimes \omega,\quad (\omega \in \caX(G)),
\; (\pi \in \mathbf{Irr}(G)).
\]
The stabilizer of $\pi \in \mathbf{Irr}(G)$ is denoted by $\caX(G)(\pi)$,
\index[not]{XGpi@$\caX(G)(\pi)$} i.e.,
\[
\caX(G)(\pi)=\{  \omega \in \caX(G)\mid \pi\otimes \omega \simeq \pi \}.
\]

\begin{lemme}
Let $\pi \in \mathbf{Irr}(G)$. Then $\caX(G)(\pi)$ is finite.
\end{lemme}

\begin{proof}
We retain the notation of the previous section. In particular we recall that $\caX(G)$ is identified with
the set of characters of the lattice $\Lambda(G)$, that the latter possesses a finite-index sublattice
$\Lambda(A_G)$, and that we can lift the latter to a subgroup $C_{A_G}$ contained in $A_G$ (hence in
particular central in $G$). The group $\caX(G)$ acts on $\caX(A_G)$ by multiplication, via the natural
restriction of characters $\caX(G) \rightarrow \caX(A_G)$. Moreover, for any $\pi \in \mathbf{Irr}(G)$,
Schur's lemma gives us a character of the central group $C_{A_G}$ by which this group acts in the space
of $\pi$. This defines a map
\[
\mathbf{Irr}(G) \rightarrow  \Hom_\bbZ(C_{A_G},\bbC^\times) \simeq
\caX(A_G)
\]
which is easily seen to be equivariant for the actions of $\caX(G)$ on $ \mathbf{Irr}(G)$ and $\caX(A_G)$
described above. It therefore suffices to show that the action of $\caX(G)$ on $\caX(A_G)$ admits finite
stabilizers. This stems from the fact that $\Lambda(A_G)$ is of finite index in $\Lambda(G)$, and
therefore that the number of characters of $\Lambda(G)$ which are trivial on $\Lambda(A_G)$ is finite.
\end{proof}

\begin{cor}
The orbits of $\caX(G)$ in $\mathbf{Irr}(G)$ are isomorphic (but non-canonically) to complex algebraic
tori.
\end{cor}

\begin{proof}
We know that $\caX(G)$ is equipped with a complex algebraic torus structure (cf. Section \ref{varXG}),
and the quotient of such a complex torus by a finite subgroup is again a complex torus. The
non-canonicity comes from the fact that one must choose a base point.
\end{proof}

\begin{notation}
The orbits of the action of $\caX(G)$ in $\mathbf{Irr}(G)$ are called inertial classes. We denote by
$[\pi]$ \index[not]{ZZpi1@$[\pi]$} the inertial class of a representation $\pi$ and $[\mathbf{Irr}(G)]$
\index[not]{Irr(G)1@$[\mathbf{Irr}(G)]$} the set of inertial classes. \index[ter]{inertial class}
\end{notation}

Let $F$ denote the algebra of polynomial functions on the variety $\caX(G)$. We saw in \ref{varXG} that
$F \simeq \bbC[\Lambda(G)]$. Let $D$ be an inertial class in $\mathbf{Irr}(G)$ and $(\pi,V)$ a
representation in this inertial class and $\caX(G)(\pi)$ the stabilizer of $\pi$. Then $D$ is identified
with the homogeneous space $\caX(G)/\caX(G)(\pi)$. The algebra of polynomial functions on $D$ is then
$F^{\caX(G)(\pi)}$, the algebra of invariants of $F$ under the action of the finite group $\caX(G)(\pi)$.
It is convenient to say that $f$ is a polynomial function on $D$ if there exists a polynomial function
$\tilde f$ on $\caX(G)$ such that $$f(\pi\otimes \chi)= \tilde f(\chi).$$ Similarly, if $U$ is a Zariski
open set of $D$, we say that $f$ is a rational function on $U$ if there exist $f_1$ and $f_2$ in $F$ such
that $ f_2(\chi)f(\pi\otimes \chi)=f_1(\chi)$ and $f_2(\chi)\neq 0$ for all $\chi \in \caX(G)$ such that
$\pi \otimes \chi \in U$.

\section{Parabolic subgroups}\label{sgp}

We refer the reader to \cite{BoTi} or \cite{Spr} for the proofs of the results stated in this section.
Let $\bbG$ be a connected reductive algebraic group defined over $\bbF$ and let us set $G=\bbG(\bbF)$. An
algebraic subgroup (as groups defined over $\overline{\bbF}$) $\bbP$ of $\bbG$ is a parabolic subgroup if
the variety $\bbP \backslash \bbG$ is projective. A parabolic subgroup \index[ter]{parabolic subgroup}
$P$ of $G$ is the group of $\bbF$-points of a parabolic subgroup $\bbP$ defined over $\bbF$.
\index[ter]{parabolic subgroup} We will say that ``$P$ is a parabolic subgroup of $G$'' to express the
fact that $P$ is the group of $\bbF$-points of a parabolic subgroup $\bbP$ of $\bbG$ defined over $\bbF$.
We will do the same for other subgroups of $\bbG$ defined over $\bbF$.

The radical and the unipotent radical of a parabolic subgroup $\bbP$ of $\bbG$ defined over $\bbF$ are
defined over $\bbF$.

\subsection{}
Let $P$ be a parabolic subgroup of $G$ and $N$ its unipotent radical. Then there exists a (non-unique)
reductive subgroup $M$ of $G$ in $P$ normalizing $N$ such that $P=M \rtimes N$. We then say that $M$ is a
Levi factor \index[ter]{Levi factor} of $P$ and a Levi subgroup \index[ter]{Levi subgroup} of $G$. Levi
subgroups are obtained in the following way. Let $A$ be a split component of $P$ and set $M=Z(G,A)$. Then
$M$ is a Levi factor of $P$. Moreover $A$ is the split component of $M$. If $P$ is given, the choice of a
split component of $P$ is therefore equivalent to the choice of a Levi factor. The unipotent radical $N$
acts transitively by conjugation on the set of split components of $P$, and therefore on the set of Levi
factors of $P$. The Levi subgroups $M$ of $G$ are the centralizers of split tori $A$ in $G$ such that $A$
is the split component of $M$. We will often use an expression like ``$P=MN$ is a parabolic subgroup of
$G$'' to summarize the situation where $P$ is a parabolic subgroup of $G$, $N$ is its unipotent radical,
and $M$ a Levi factor. If $M$ is a Levi subgroup of $G$, we denote by $\caP(M)$
\index[not]{P(M)@$\caP(M)$} the set of parabolic subgroups of $G$ admitting $M$ as a Levi factor. It is a
finite set. Let $P \in \caP(M)$. We then denote by $\overline{ P}=M\overline{N}$ the opposite parabolic
subgroup to $P$, i.e., the unique parabolic subgroup $Q \in \caP(M)$ such that $P\cap Q=M$, and $A_M$
\index[not]{A_M@$A_M$} (or simply $A$) the split component of $M$.

\subsection{}
The minimal parabolic subgroups of $G$ are those whose split components are maximal split tori of $G$.
Two split tori of $G$, maximal for this property, are conjugate in $G$. Two minimal parabolic subgroups
of $G$ are therefore conjugate in $G$.

\subsection{}
Let $P=MN$ be a parabolic subgroup of $G$, and $A$ the split component of $M$. We set:
\[
\fra =X_*(A)\otimes_\bbZ\bbR,\quad
\fra^*=X^*(A)\otimes_\bbZ \bbR.
\]
\index[not]{a@$\fra$}\index[not]{a*@$\fra^*$} The perfect duality between $X_*(A)$ and $X^*(A)$ extends
to a vector space duality between $\fra$ and $\fra^*$, which is consistent with the notation.

\subsection{}
Let $P=MN$ be a parabolic subgroup of $G$, and $A$ the split component of $M$. Let $\Sigma'(A)$
\index[not]{ZZSigma'(A)@$\Sigma'(A)$} denote the set of roots of $A$ in $\mathrm{Lie}(G)$. Then
$\Sigma'(A)$ is identified with a subset of $\fra^*$. We denote by $\Sigma(A)$
\index[not]{ZZSigma(A)@$\Sigma(A)$} the set of reduced roots of $\Sigma'(A)$, i.e., the elements $\alpha$
of $\Sigma'(A)$ such that $\alpha/n\notin \Sigma'(A)$ if $n\geq 2$.

If $P=MN \in \caP(M)$, we denote by $\Sigma'(P)$ \index[not]{ZZSigma(P)'@$\Sigma'(P)$} the system of
positive roots in $\Sigma'(A)$ relative to $P$ (the roots of $A$ in $N$) and $\Sigma(P)$
\index[not]{ZZSigma(P)@$\Sigma(P)$} the system of reduced positive roots.

\subsection{} \label{actfid}
Let $P=MN$ be a parabolic subgroup of $G$, and $A$ the split component of $M$. The Weyl group
\begin{equation}\label{WA}
W(A)=N_G(A)/Z_G(A)=N_G(A)/M=N_G(M)/M
\end{equation}
\index[not]{W(A)@$W(A)$} acts on $X^*(A)$ and therefore on $X_*(A)\simeq \Lambda(A)$. Consequently,
$W(A)$ also acts on $\fra^*$ and $\fra$. The action of $W(A)$ on $X^*(A)$ or on
$X_*(A) \simeq \Lambda(A)$ is faithful.

\subsection{}\label{splitcenter}
Let $P=MN$ be a parabolic subgroup of $G$, with split component $A_M$. Let $A_G$ be the split component
of $G$. We have
\[
A_G \subset A_M \subset M \subset G,
\]
which induces morphisms given by the restriction of characters:
\[
\xymatrix{
X^*(G)   \ar[r] & X^*(M)  \ar[r] &  X^*(A_M)\ar[r]&     X^*(A_G)
}.
\]
Let us tensor by $\bbR$. We obtain, according to Lemma \ref{Lambda}
\[
\xymatrix{
\fra_M^*=X^*(M) \otimes_\bbZ \bbR  \ar[r] &  X^*(A_M)\ar[d] \otimes_\bbZ \bbR =\fra_M^*  \\
\fra_G^*=X^*(G) \otimes_\bbZ \bbR  \ar[u]  \ar[r]&  X^*(A_G)\otimes_\bbZ \bbR =\fra_G^*\\
},
\]
the horizontal arrows being isomorphisms. Let $(\fra_M^G)^*$ denote the kernel of the right vertical
arrow. We then have
\[
\fra_M^*=\fra_G^* \oplus  (\fra_M^G)^*.
\]
Similarly, dually, we obtain
\[
\fra_M=\fra_G \oplus  \fra_M^G.
\]

\subsection{}
We continue with the same notation. We will assume $\fra_M$ and $\fra_M^*$ are always equipped with inner
products (denoted $(.,.)$, the duality between $\fra_M$ and $\fra_M^*$ being denoted $\bilo$) invariant
under the action of $W(A_M)$.

When $P$ is a minimal parabolic subgroup, $\Sigma(A_M)$ is a root system in $\fra_M^*$ equipped with an
inner product as above (\cite{BoTi}, Cor. 5.8) and we can thereby identify $W(A_M)$ with the Weyl group
of the root system $\Sigma(A_M)$ (\cite{BoTi}, 5.3). Moreover, if $\alpha \in\Sigma(A_M)$, we denote by
$\check \alpha$ its coroot, i.e., the element of $\fra_M$ satisfying
\[
\langle \check \alpha,\lambda \rangle = 2\frac{(\alpha,\lambda)}{(\alpha,\alpha)} , \qquad (\lambda \in \fra^*)    .
\]

The set of coroots, denoted $\Sigma^\vee(A_M)$, is a root system in $\fra_M$. The roots (resp. the
coroots) generate the subspace $(\fra_M^G)^*$ (resp. $\fra_M^G$) defined in \ref{splitcenter}.

\subsection{}
Let us fix a maximal split torus $A_\emptyset$ \index[not]{A_0@$A_\emptyset$} of $G$ and let
$M_\emptyset$ \index[not]{M_0@$M_\emptyset$} denote its centralizer in $G$. It is a Levi factor of a
minimal parabolic subgroup in $G$. We fix a minimal parabolic subgroup
$P_\emptyset=M_\emptyset N_\emptyset$,\index[not]{P_0@$P_\emptyset$} with Levi factor $M_\emptyset$. If
$P$ is a parabolic subgroup of $G$, we say that $P$ is semi-standard if $A_\emptyset\subset P$ and
standard if $P_\emptyset\subset P$. In both cases, $P$ has a unique Levi factor $M$ containing
$A_\emptyset$. Since the split component of $M$ is clearly contained in $A_\emptyset$, we have
$M_\emptyset \subset M$. We then say that $M$ is a semi-standard Levi subgroup
\index[ter]{Levi subgroup!semi-standard} of $G$. We say that such an $M$ is standard if it is a Levi
factor of a standard parabolic subgroup of $G$. The split torus $A_\emptyset$ is again maximal for this
property in any standard Levi subgroup $M$, and we can thus define the standard or semi-standard Levi
subgroups and parabolic subgroups
\index[ter]{parabolic subgroup!standard}\index[ter]{parabolic subgroup!semi-standard} of $M$ with
reference to this same $A_\emptyset$.

\subsection{} \label{conjstand}
Let us fix a minimal parabolic subgroup $P_\emptyset=M_\emptyset N_\emptyset$ of $G$, with split
component $A_\emptyset$. Then any parabolic subgroup of $G$ is conjugate to a standard parabolic
subgroup.

\subsection{}\label{decsemistandard}
Let $P=MN$ and $Q=LU$ be two semi-standard parabolic subgroups of $G$. Then $M\cap Q$ is a semi-standard
parabolic subgroup of $M$ with Levi factor $M\cap L$ and unipotent radical $M\cap U$. Moreover $N\cap Q$
decomposes into $N\cap Q=(N\cap L)(N\cap U)$.

\subsection{}\label{semst}
Let $P=MN$ be a semi-standard parabolic subgroup of $G$. There exists a bijection between the set of
standard (resp. semi-standard) parabolic subgroups of $M$ and the set of standard (resp. semi-standard)
parabolic subgroups $Q=LU$ of $G$ contained in $P$, whose inverse is obtained as follows: to $Q=LU$, we
associate the parabolic subgroup $M\cap Q= L( M\cap U)$. The standard minimal parabolic subgroup of $M$
is therefore $M_\emptyset(M\cap N_\emptyset)$.

\subsection{} \label{psi}
We fix a minimal parabolic subgroup $P_\emptyset=M_\emptyset N_\emptyset$, with split component
$A_\emptyset$. We set
\[
\Sigma'_\emptyset=\Sigma'(A_\emptyset), \;  \Sigma_{\emptyset}=
\Sigma(A_\emptyset), \;  \Sigma_\emptyset^+=\Sigma(P_\emptyset),
\]
and we denote by $\Delta_\emptyset=\Delta(P_\emptyset)$ the set of simple roots in
$\Sigma_{\emptyset}^+$. \index[not]{ZZSigma'_0@$\Sigma'_\emptyset$}
\index[not]{ZZSigma_0@$\Sigma_\emptyset$} \index[not]{ZZDelta_0@$\Delta_\emptyset$} Similarly, we use the
notation ${\Sigma'_\emptyset}^\vee$, $\Sigma_{\emptyset}^\vee$, $(\Sigma_\emptyset^\vee)^+$ for the
coroots. \index[not]{ZZSigma'c_0@${\Sigma'_\emptyset}^\vee$}
\index[not]{ZZSigma_0c@$\Sigma_\emptyset^\vee$} \index[not]{ZZSigma_0cp@$(\Sigma_\emptyset^\vee)^+$}
\index[not]{ZZDelta_0c@$\Delta_\emptyset^\vee $}

The set of simple roots $\Delta_\emptyset$ is a basis of $(\fra_\emptyset^G)^*$ and the set of simple
coroots $\Delta_\emptyset^\vee$ is a basis of $\fra_\emptyset^G$. Let
$\widehat\Delta_\emptyset= \{ \varpi_\alpha , \, \alpha \in \Delta_\emptyset\}$ be the basis of
$(\fra_\emptyset^G)^*$ dual to $\Delta_\emptyset^\vee$, whose elements are called the fundamental
weights, \index[ter]{fundamental weights} \index[not]{ZZDelta_0h@$\widehat \Delta_\emptyset$} and
$\widehat \Delta_\emptyset^\vee=  \{  \varpi_\alpha^\vee, \;  \alpha \in \Delta_\emptyset\}$ the basis of
$\fra_\emptyset^G$ dual to $\Delta_\emptyset$, whose elements are called the fundamental coweights.
\index[not]{ZZDelta_0hc@$\widehat \Delta_\emptyset^\vee$}

If $P=MN$ is a standard parabolic subgroup, we denote by $\Delta_\emptyset^M$ the analogue of
$\Delta_\emptyset$ when we replace $G$ by $M$ in the definition.
\index[not]{ZZDelta_0M@$\Delta_\emptyset^M$}

\subsection{}  \label{traceP}
Let $P=MN\subset Q=LU$ be two standard parabolic subgroups of $G$. We then have
$A_G \subset A_L\subset A_M \subset A_\emptyset$ and
\[
\fra_G \subset \fra_L \subset \fra_M \subset \fra_\emptyset
\]
where $\fra_G$ (resp. $\fra_M$, resp. $\fra_L$) is the orthogonal of the roots
$\alpha \in  \Delta_\emptyset$ (resp. $\Delta_\emptyset^M$, resp. $\Delta_\emptyset^L$). Similarly
\[
\fra_G^* \subset \fra_L^* \subset \fra_M^*\subset \fra_\emptyset^*
\]
are described in the same way by replacing the roots with the coroots. Let $\fra_\emptyset^G$ (resp.
$(\fra_\emptyset^G)^*$) denote the subspace of $\fra_\emptyset$ (resp. of $\fra_\emptyset^*$) generated
by the $\check \alpha \in \Delta_\emptyset^\vee$ (resp. the $\alpha \in \Delta_\emptyset$). We then have
the decompositions (see \ref{splitcenter})
\[
\fra_\emptyset = \fra_G \oplus \fra_\emptyset^G, \quad   \fra_\emptyset^*=
\fra_G^* \oplus (\fra_\emptyset^G)^*.
\]

More generally, let $(\fra_M^L)^*$ be the subspace of $\fra_M^*$ generated by the restrictions
$\alpha_{|\fra_M}$ of the roots $\alpha \in \Delta_\emptyset^L\setminus \Delta_\emptyset^M$.
\index[not]{a1(ML)@$(\fra_M^L)^*$} We then have the decomposition
\[
\fra_M^*= \fra_L^* \oplus (\fra_M^{L})^*.
\]
In particular
\[
\fra_M^*= \fra_G^* \oplus (\fra_M^{G})^*,
\]
where $(\fra_M^G)^*$ is generated by the restrictions $\alpha_{|\fra_M}$ of the roots
$\alpha \in \Delta_\emptyset\setminus \Delta_\emptyset^M$.

We obtain similar decompositions of $\fra_M$ by exchanging the role of roots and coroots. These
decompositions are orthogonal for the fixed inner product.

We denote by $p_M^L$ \index[not]{pML@$p_M^L$} the projections respectively of $\fra_M$ and $\fra_M^*$
onto $\fra_M^{L}$ and $(\fra_M^{L})^*$ in the above decompositions.

\subsection{}
Let $P=MN$ be a parabolic subgroup of $G$, and $A_M$ the split component of $M$. Let us choose a minimal
parabolic subgroup $P_\emptyset=M_\emptyset N_\emptyset$, with split component $A_\emptyset$, contained
in $P$ with $A_M$ contained in $A_\emptyset$ (this is possible according to \ref{conjstand}). The
inclusion $A_M \subset A_\emptyset$ induces an inclusion of $\fra_M$ into $\fra_\emptyset$. Let us define
\[
\Delta(P):=\{ \alpha_{|\fra_M}, \; \alpha \in \Delta_\emptyset\setminus
\Delta_\emptyset^M  \}.
\]
\[
\widehat \Delta(P):=\{ {\varpi_\alpha}_{|\fra_M}, \; \alpha \in
\Delta_\emptyset\setminus \Delta_\emptyset^M  \}.
\]
We define dually, by replacing roots and fundamental weights with coroots and fundamental coweights,
$\Delta^\vee(P)$ and $\widehat \Delta^\vee(P)$. If $\beta= \alpha_{|\fra_M} \in  \Delta(P)$,
$\alpha \in  \Delta_\emptyset\setminus \Delta_\emptyset^M $, we denote by $\check \beta$ the projection
onto $\fra_M$ of the coroot $\check \alpha$.

\index[not]{ZZDeltaP@$\Delta(P) $}\index[not]{ZZDeltaPh@$\widehat\Delta(P) $}
\index[not]{ZZDeltaPv@$\Delta^\vee(P) $}\index[not]{ZZDeltaPhv@$\widehat \Delta^\vee(P) $}

Note that the use of this notation is not standard because $\Delta(P) $ is not a set of simple roots of a
root system. In particular, if $\beta \in \Delta(P)$, $\check \beta$ is not a coroot.

Moreover, $(\fra_M^G)^*$ is generated by $\Delta(P)$, because the roots of $\Delta(P)$ are the
restrictions to $\fra_M$ of the roots in $\Delta_ \emptyset\setminus  \Delta_ \emptyset^M$, and
$\fra_M^G$ is generated by the projections of the coroots in
$\Delta_ \emptyset^\vee \setminus (\Delta_ \emptyset^\vee)^M$.

\subsection{}  \label{APG}
Let us continue with the notation of the previous paragraph. We set
\[
{}^+[\fra^*]_P^G= \{   \chi =  \sum_{\alpha \in \Delta(P)}
c_\alpha \; \alpha, \quad  c_\alpha>0 \}
\]
\[
{}^+\overline{[\fra^*]}_P^G= \{   \chi =  \sum_{\alpha \in \Delta(P)}
c_\alpha \; \alpha, \quad  c_\alpha \geq 0 \}
\]
\[
{}_P^G [\fra^*]^+ = \{ \chi
\in \fra^* \mid \bil{\chi}{ \check \alpha}> 0, \; \alpha \in \Delta(P)  \}
\]
\[
{}_P^G \overline { [\fra^*]}^+ = \{ \chi
\in \fra^* \mid \bil{\chi}{ \check \alpha} \geq  0,  \; \alpha \in
\Delta(P) \}.
\]
\index[not]{a1+PG@${}^+[\fra^*]_P^G$}\index[not]{a1+PG@${}^+\overline{[\fra^*]}_P^G$}
\index[not]{a1PG+@$ {}_P^G [\fra^*]^+$} \index[not]{a1PG+@$ {}_P^G \overline { [\fra^*]}^+$} We define
similarly the subsets ${}^+[\fra]_P^G$, ${}^+\overline{[ \fra]}_P^G$, ${}_P^G [\fra]^+$,
${}_P^G \overline{[ \fra]}^+$ by exchanging the role of roots and coroots. We also denote
\[
{}^-[\fra]_P^G= -({}^+[\fra]_P^G),\;    {}^-\overline{[ \fra]}_P^G=
-({}^+\overline{[ \fra]}_P^G), \,  {}_P^G [\fra]^-=-({}_P^G [\fra]^+), \;
{}_P^G \overline{[ \fra]}^-=-({}_P^G \overline{[ \fra]}^+).
\]

Note that an element $\chi \in  {}_P^G [\fra^*]^+$ is an element which can be written
\[
\chi =  \sum_{\alpha \in \Delta(P)}     c_\alpha \; \varpi_\alpha \; + \mu,
\quad  c_\alpha>0,\;  \mu \in \fra_G^*,
\]
(similarly for ${}_P^G\overline{[\fra^*]}^+$ with the $c_\alpha \geq 0$).

\subsection{} \label{plusam}
Let $P=MN$ be a standard parabolic subgroup of $G$. The subset $ {}_P^G { [\fra_M^*]}^+$ is not an open
convex cone (in the strict sense) in general, since if $x \in   {}_P^G { [\fra_M^*]}^+$, then
$x+\fra_G^* \subset  {}_P^G { [\fra_M^*]}^+$. On the other hand, the intersection
\[
{}_P^G { [(\fra_M^G)^*]}^+  = {}_P^G { [\fra_M^*]}^+  \cap  (\fra_M^{G})^* =
\{   \chi =  \sum_{\alpha \in \Delta(P)}     c_\alpha \; \varpi_\alpha, \quad
c_\alpha>0 \}
\]
is an open convex cone and its closure
\[
{}_P^G \overline { [(\fra_M^G)^*]}^+ =  {}_P^G \overline { [\fra_M^*]}^+
\cap  (\fra_M^{G})^*  = \{   \chi =  \sum_{\alpha \in \Delta(P)}
c_\alpha \; \varpi_\alpha, \quad  c_\alpha\geq 0 \}
\]
is a closed convex cone.

We have
\[
{}_P^G  [\fra_M^*]^+=  {}_P^G { [(\fra_M^G)^*]}^+  \oplus \fra_G^*,
\]
\[
{}_P^G \overline { [\fra_M^*]}^+=   {}_P^G \overline { [(\fra_M^G)^*]}^+
\oplus \fra_G^*.
\]

The cone $ {}_P^G \overline { [(\fra_M^G)^*]}^+$ admits a cellular decomposition into cones of the same
kind coming from the standard parabolic subgroups $Q$ containing $P$:
\[
{}_P^G \overline { [(\fra_M^G)^*]}^+ = \coprod_{P \subset Q=LU}
{}_Q^G { [(\fra_L^G)^*]}^+.
\]

We deduce a decomposition
\[
{}_P^G \overline { [\fra_M^*]}^+ = \coprod_{P \subset Q=LU}
{}_Q^G { [\fra_L^*]}^+  .
\]

\subsection{} \label{fait24}
We note the following fact:

\begin{lemme}
Let $P=MN$ be a standard parabolic subgroup of $G$ and let $\mu$ be an element of
$ {}^G_P[(\fra_M^G)^*]^+$. Then the set of roots $\gamma$ in $\Sigma_\emptyset$ such that
$\bil{\mu}{\check \gamma}>0$ is the set of roots of $A_\emptyset$ in $N$.
\end{lemme}

\begin{proof}
By definition, $  {}^G_P[(\fra_M^G)^*]^+$ is the set of elements $\mu \in (\fra_\emptyset^G)^*$ such that
$(\mu, \alpha) = 0$ for all $\alpha \in \Delta_\emptyset^M$, and $(\mu, \alpha) > 0$ for all
$\alpha \in \Delta_\emptyset \setminus \Delta_\emptyset^M$. Any root $\gamma \in \Sigma_\emptyset$ can be
written as a linear combination $\gamma=\sum_{\alpha \in \Delta_\emptyset} n_\alpha \, \alpha$, where the
coefficients $n_\alpha \in \bbZ$ are either all non-negative or all non-positive. Using the
$W(A_\emptyset)$-invariant inner product on $\fra_\emptyset^*$, we have the relation
$\bil{\mu}{\check \gamma} = \frac{2(\mu, \gamma)}{(\gamma, \gamma)}$. By linearity of the inner product
in the second variable, we obtain:
\[
(\mu, \gamma) = \sum_{\alpha \in \Delta_\emptyset} n_\alpha (\mu, \alpha)
= \sum_{\alpha \in \Delta_\emptyset \setminus \Delta_\emptyset^M} n_\alpha (\mu, \alpha).
\]
Since $(\mu, \alpha) > 0$ for all $\alpha \in \Delta_\emptyset \setminus \Delta_\emptyset^M$, the inner
product $(\mu, \gamma)$ (and therefore $\bil{\mu}{\check \gamma}$) is strictly positive if and only if
all $n_\alpha \geq 0$ and there exists at least one
$\alpha \in \Delta_\emptyset \setminus \Delta_\emptyset^M$ such that $n_\alpha > 0$. This is exactly the
condition for $\gamma$ to be a root of $A_\emptyset$ in the unipotent radical $N$.
\end{proof}

\subsection{Langlands' combinatorial lemma} \label{comblangl}
\index[ter]{Langlands' combinatorial lemma}

Let $P=MN$ be a standard parabolic subgroup of $G$ and $\mu \in \fra_M^*$. Langlands' combinatorial lemma
asserts that there exists a unique standard parabolic subgroup $Q$ containing $P$ such that $\mu$
decomposes into
\begin{equation} \label{LangComb}
\mu=\mu_G+\mu^+ +\mu^-
\end{equation}
with $\mu_G \in \fra_G^*$, $\mu^+ \in   {}_Q^G  { [(\fra_L^G)^*]}^+$ and
$\mu^- \in - p_M^L( {}^+\overline{[\fra_M^*]}_P^G )$. This decomposition is unique. Moreover $\mu_G$ is
the orthogonal projection of $\mu$ onto $\fra_G^*$, $\mu^+$ is the orthogonal projection of $\mu$ onto
the cone $ {}_P^G \overline { [(\fra_M^G)^*]}^+$ and $\mu^-$ is the orthogonal projection of $\mu$ onto
the cone ${}^-\overline{[\fra_M^*]}_P^G$. The elements $\mu_G$, $\mu^+$ and $\mu^-$ are pairwise
orthogonal. Note that if we consider $\mu$ as an element of $\fra_\emptyset^*$, and apply the above
result to $(\mu,P_\emptyset)$, we find the same decomposition.

We give a proof of this lemma, adapted from \cite{Carm}, at the end of the chapter.

\subsection{}
Let $P=MN$ be a parabolic subgroup of $G$. We give another interpretation of the space $\fra^*$. Since
$$X_*(A)\subset \Lambda(M) \subset X_*(M)$$ are lattices of the same rank, we have
\[
\fra=X_* (A) \otimes_\bbZ \bbR=\Lambda(M)\otimes_\bbZ \bbR=X_*(M)\otimes_\bbZ \bbR.
\]
Therefore
\begin{align*}
\fra^*&= \Hom_\bbR(\Lambda(M)\otimes_\bbZ \bbR,\bbR)\\
&= \Hom_\bbZ(\Lambda(M),\bbR)\\
&=\{ \phi \in   \Hom_\bbZ(M,\bbR) \mid \phi_{|{}^0M}\equiv 0   \}\\
&= \Hom_{\bbZ, smooth}(M,\bbR)\\
\end{align*}
where $ \Hom_{\bbZ, smooth}(M,\bbR)$ is the group of smooth morphisms from $M$ to $\bbR$. Care must be
taken that the operation of $\Lambda(M)$ is written additively, whereas that of $M$ is written
multiplicatively. It remains to justify the last equality, by showing that any smooth morphism from $M$
to $\bbR$ is trivial on ${}^0M$. We know that the group $R(M)\scrD M$ is cocompact in $M$ and that
$\scrD M \subset {}^0M$. The group $(R(M)\cap {}^0M)\scrD M$ is therefore cocompact in ${}^0M$. Moreover,
we saw in the proof of Proposition \ref{Lambda} that $R(M)\cap {}^0M$ is compact, and therefore any
smooth morphism from this group with values in $\bbR$ is necessarily trivial (see Lemma \ref{carnonram}).
Of course, any morphism from $M$ to the abelian group $\bbR$ is trivial on commutators. Any smooth
morphism $\phi$ from $M$ to $\bbR$ is therefore trivial on a finite-index subgroup of ${}^0M$ (its kernel
being open, the subgroup $(R(M)\cap {}^0M)\scrD M\cdot\ker\phi$ of ${}^0M$ is open and cocompact, hence
of finite index). Since $\bbR$ contains no element of finite order other than $0$, we see that such a
morphism is trivial on ${}^0M$.

For any character $\chi$ of $M$, $\ln |\chi|  \in \fra^*$.

\subsection{} \label{ReChi}
Let us resume the notation of the previous paragraph. Let us set
\[
\fra_\bbC= X_*(A)\otimes_\bbZ \bbC \text{ and }
\fra_\bbC^*= X^*(A)\otimes_\bbZ \bbC.
\]
We have as above
\[
\fra_\bbC=X_* (A) \otimes_\bbZ \bbC=\Lambda(M)\otimes_\bbZ \bbC=X_*(M)\otimes_\bbZ \bbC
\]
and
\begin{align*}
\fra_\bbC^*&= \Hom_\bbC(\Lambda(M)\otimes_\bbZ \bbC,\bbC)\\
&= \Hom_\bbZ(\Lambda(M),\bbC).
\end{align*}

Recall that $q$ denotes the cardinality of the residue field $\frO_\bbF/\frP_\bbF$. The group morphism
\[
\bbC \rightarrow \bbC^\times, \quad s \mapsto q^{-s}
\]
is surjective, with kernel $\frac{2i\pi}{\ln q} \bbZ$. It induces by composition a morphism
\[
\fra_\bbC^*= \Hom_\bbZ(\Lambda(M),\bbC) \rightarrow
\Hom_\bbZ(\Lambda(M),\bbC^\times ) \simeq \caX(M) .
\]
The group $\Lambda(M)$ being a lattice, we easily see that this morphism is surjective, by lifting to
$\bbC$ the values taken on a basis of this lattice. The kernel of this morphism is
\[
\Hom_\bbZ(\Lambda(M), \frac{2i\pi}{\ln q} \bbZ ) =
\frac{2i\pi}{\ln q}   \Hom_\bbZ(\Lambda(M), \bbZ )=\frac{2i\pi}{\ln q}\caR
\]
where $\caR= \Hom_\bbZ(\Lambda(M), \bbZ ) $ is a sublattice of $\fra^*$.

We have thus obtained a surjective morphism:
\[
\fra_\bbC^*   =X^*(M)\otimes_\bbZ \bbC \rightarrow \caX(M)
\]
which is easily verified to be given by
\[
\chi \otimes s \mapsto [\; g \mapsto |\chi(g)|_\bbF^s \;  ]
\]
and whose kernel is the lattice $\frac{2i\pi}{\ln q}\caR$. We thus recover the algebraic torus structure
of $\caX(M)$. If $\lambda \in \fra_\bbC^*$, we denote by $e^\lambda$ the corresponding unramified
character of $M$.

If $\chi \in \caX(M)$, let us lift it to an element $\lambda \in \fra_\bbC^*$. The real part
$\Re (\lambda) \in \fra^*$ is independent of the choice of $\lambda$. We denote it by $\Re(\chi)$.
\index[not]{R(chi)@$\Re(\chi)$}

Any element $\chi \in  \Hom_\bbZ(M,\bbC^\times)$ decomposes into
\[
\chi=\frac{\chi}{|\chi|} \, |\chi|.
\]
The character $\frac{\chi}{|\chi|} $ is unitary, $|\chi|$ is unramified, and we still denote by
$\Re(\chi)\in\fra^*$ the element such that $|\chi(m)|=q^{-\langle\Re(\chi),H_M(m)\rangle}$, $m\in M$,
i.e., $\ln|\chi|=-\ln q\,\Re(\chi)$ in $\Hom_\bbZ(M,\bbR)=\fra^*$; this is compatible with the previous
definition.

We denote by $\im(\caX(M))$ \index[not]{ImX(M)@$\im(\caX(M))$} the group of unitary unramified characters
of $M$. It is the set of $\chi \in \caX(M)$ such that $ \Re(\chi)=0$. It is also the image of $i\fra^*$
under the surjective morphism $ \fra_\bbC^*  \rightarrow \caX(M)$. It is a real submanifold of the
complex torus $\caX(M)$, and more precisely, it is a compact torus, isomorphic to
$i\fra^*/\frac{2i\pi}{\ln q}\caR$.

\subsection{}  \label{bof}
Let $P=MN$ be a parabolic subgroup of $G$. Recall the map $H_M \colon M \rightarrow X_*(M)\subset \fra_M$
defined in \ref{carnonram}. We denote by $M^+$ (resp. $M^{++}$) the inverse image
\index[not]{M1@$M^+$ }\index[not]{M2@$M^{++}$} under $H_M$ of $ {}_P^G \overline{[ \fra]}^+\cap X_*(M)$
(resp. $ {}_P^G [\fra]^+ \cap X_*(M)$). More generally, if $H$ is a subset of $M$, we denote
$H^+=H\cap M^+$ and $H^{++}= H \cap M^{++}$. It follows directly from the definitions that
\begin{align*}
 A_M^+&=\{a \in A_M,\, v_\bbF(\alpha(a))\geq 0, \; (\alpha \in \Delta(P))\}\\
&=\{a \in A_M,\, |\alpha(a)|_\bbF \leq 1, \; (\alpha \in \Delta(P))\}\\
& \\
A_M^{++}&=\{a \in A_M,\, v_\bbF(\alpha(a))> 0, \; (\alpha \in \Delta(P))\}\\
&=\{a \in A_M,\,   |\alpha(a)|_\bbF < 1 , \; (\alpha \in \Delta(P))\}.
\end{align*}
\index[not]{A_M1@$A_M^+$}\index[not]{A_M2@$A_M^{++}$} Note that the choice of $P$ does not appear in the
notation, although the subsets $A_M^+$ and $A_M^{++}$ depend on it. We will mostly use this notation with
standard parabolic subgroups, which will make this ambiguity harmless. For example, if $P$ is standard,
the decomposition \ref{plusam} induces, for the subgroups $C_{A_M}\subset A_M$ and $C_{A_L}\subset A_L$
introduced in \ref{FM0}, a decomposition
\begin{equation}\label{coneA}
C_{A_M}^+=    \coprod_ {Q} C_{A_L}^{++},
\end{equation}
the sum being over the standard parabolic subgroups $Q=LU$ containing $P$, where
$C_{A_M}^+=C_{A_M}\cap A_M^+$ and $C_{A_L}^{++}=C_{A_L}\cap A_L^{++}$; here we use that $C_{A_L}$
identifies with $X_*(A_L)=X_*(A_M)\cap\fra_L$. (Such a decomposition does not hold for $A_M^+$ itself:
for $G=SL_2$, $-1\in A_\emptyset^+$ belongs to no $A_L^{++}$.)

\subsection{}
Let us set, for all $\epsilon >0$,
\begin{equation}\label{epsilon}
 A^{+}(\epsilon)=\{a \in A_M,\,       |\alpha(a)|_\bbF <\epsilon, \, (\alpha \in
\Delta(P))\},
\end{equation}
\index[not]{A1E@$A^+(\epsilon)$} and for any subset $\Omega$ of $A$,
$\Omega^+(\epsilon)=\Omega \cap A^{+}(\epsilon)$.

\subsection{}  \label{FM0}
Let $P=M N$ be a parabolic subgroup of $G$, with split component $A$. The lattice $\Lambda(M)=M/{}^0M$
does not in general admit a natural lift to a subgroup of $M$, but as we saw in \ref{castordep}, the
sublattice $\Lambda(A)=X_*(A)$ does admit a lift, denoted $C_{A}$, to a subgroup of $A$. Since $X_*(A)$
is a finite-index sublattice of $\Lambda(M)$, there exists a finite set $F_M$ of elements of $M$ such
that $\widetilde C_A:= C_AF_M$ is \index[not]{CAt@$\widetilde C_A$} a set-theoretic lift of $\Lambda(M) $
in $M$. We can take the elements of $F_M$ in $M^{+}$. Indeed, we can always replace $f \in F_M$ by $t^mf$
with $t\in C_A^{++}$, and we will then have $t^mf \in M^{+}$ for $m$ large enough. Let
$C_A^+= C_A \cap M^+$, $\widetilde C_A^+:= \widetilde C_A \cap M^+$,
$C_A^+(\epsilon)=C_A \cap  A^{+}(\epsilon)$. In general, we cannot arrange for
$\widetilde C_A^+= C_A^+F_M$. On the other hand, there exists a finite set $F_M^+\subset\widetilde C_A^+$
such that $\widetilde C_A^+= C_A^+F_M^+$, the writing of an element of $\widetilde C_A^+$ in this form
being no longer unique. Indeed, let $t_\alpha\in C_A^+$, $\alpha\in\Delta(P)$, be elements constructed as
in \ref{SSM}, such that $\langle\alpha,H_M(t_\alpha)\rangle=N_0>0$ and
$\langle\beta,H_M(t_\alpha)\rangle=0$ for $\beta\neq\alpha$. Any $m\in\widetilde C_A^+$ can be written
$m=\big(\prod_\alpha t_\alpha^{k_\alpha}\big)m'$, where $k_\alpha$ is the integer part of
$\langle\alpha,H_M(m)\rangle/N_0$; then $m'\in\widetilde C_A^+$ and
$0\leq\langle\alpha,H_M(m')\rangle<N_0$ for all $\alpha\in\Delta(P)$, which leaves only finitely many
possibilities for $m'$ modulo $C_{A_G}\subset C_A^+$. \index[not]{FM+@$F_M^+$} When
$P= P_\emptyset=M_\emptyset N_\emptyset$, we denote $\Lambda(M_\emptyset)=\Lambda_\emptyset$,
$C_{A_\emptyset }=C_\emptyset$, $F_\emptyset=F_{M_\emptyset}$, etc.\index[not]{FM@$F_M$}
\index[not]{F0@$F_\emptyset$}

\subsection{} \label{SSM}
Since $A_G\subset A_\emptyset$, according to a remark made in \ref{castordep}, $C_{A_G}$ is a subgroup of
$C_\emptyset$. Since $A_G \subset A_\emptyset^+$, we have $C_{A_G}\subset C_\emptyset^{+}$. Let us
introduce on $C_\emptyset^{+}$ the equivalence relation:
\begin{equation}
\label {eqC}
a_1\sim a_2   \text{ if }  a_1 a_2^{-1} \in C_{A_G}.
\end{equation}
Set $S:=C_\emptyset^+/C_{A_G}=C_\emptyset^+/\sim$ and $d=|\Delta_\emptyset|$. The simple roots are
linearly independent and vanish on $\fra_G$; hence there exist an integer $N_0>0$ and, for every
$\alpha\in\Delta_\emptyset$, an element $t_\alpha\in C_\emptyset$ such that
$|\alpha(t_\alpha)|_\bbF=q^{-N_0}$ and $|\beta(t_\alpha)|_\bbF=1$ for $\beta\in\Delta_\emptyset$,
$\beta\neq\alpha$: via the isomorphism $C_\emptyset\simeq X_*(A_\emptyset)$ of \ref{castordep}, one has
to solve in $X_*(A_\emptyset)$ the system $\langle\beta,x\rangle=N_0\delta_{\alpha\beta}$,
$\beta\in\Delta_\emptyset$, which has a rational solution. The $t_\alpha$ belong to $C_\emptyset^+$, and
their images generate a submonoid $S_0$ of $S$ isomorphic to $\bbN^d$. Let $R$ be the set of classes of
the $z\in C_\emptyset^+$ such that $|\alpha(z)|_\bbF>q^{-N_0}$ for all $\alpha\in\Delta_\emptyset$; it is
a finite set, since the image of $X_*(A_\emptyset)$ in $\fra_\emptyset/\fra_G$ is a lattice and
$X_*(A_G)$ has finite index in $X_*(A_\emptyset)\cap\fra_G$. We have
\[
S=\bigcup_{r\in R}rS_0:
\]
if $z\in C_\emptyset^+$ and if $n_\alpha$ denotes the integral part of $v_\bbF(\alpha(z))/N_0$, the
element $z\prod_\alpha t_\alpha^{-n_\alpha}$ belongs to $C_\emptyset^+$ and its class lies in $R$. We
deduce from \ref{coneA} a decomposition
\[
S=   \coprod_ {Q} S_L^{++},
\]
the sum being over the standard parabolic subgroups $Q=LU$ and $S_L^{++}= C_{A_L}^{++}/C_{A_G}$.

\subsection{}\label{pasmieux}
Let $P=MN$ be a parabolic subgroup of $G$, with split component $A$. We assume $P \neq G$, so that $A_G$
is strictly included in $A$. According to Remark \ref{castordep}, we can view $C_{A_G}$ as a subgroup of
$C_A$. Since $A_G$ is the direct product of ${}^0 A_G$ and $C_{A_G}$, and $A_G \cap {}^0G ={}^0A_G$
(Lemma \ref{Lambda}) we have:
\[
C_{A_G} \cap {}^0G=   C_{A_G}\cap A_G \cap {}^0G=  C_{A_G}\cap {}^0A_G=\{1\}.
\]
We deduce that
\[
C_A \cap {}^0GA_G= C_{A_G}(C_A \cap {}^0G),
\]
and therefore that the inclusion $C_A \hookrightarrow G$ induces an inclusion
\[
C_A/ C_{A_G}(C_A \cap {}^0G)  \hookrightarrow G/{}^0GA_G.
\]
Since the right-hand side is finite according to Proposition \ref{Lambda}, the left-hand side is too.

Let us take an element $t$ in $C_A^{++}$. Its powers $t^n$, $n \in \bbN^*$ are in $C_A^{++}$ and one of
them is in $C_{A_G}(C_A \cap {}^0G)$. This shows that $C_A^{++}\cap {}^0G$ is non-empty.

\section{Weyl groups}\label{Weylgroups}

We fix a minimal parabolic subgroup $P_\emptyset=M_\emptyset N_\emptyset$ of $G$ and we use the notation
of \ref{psi}. For all $\gamma \in \Sigma_\emptyset$, we denote by $U_\gamma$ the unipotent subgroup of
$G$ normalized by $A_\emptyset$ such that the adjoint action of $A_\emptyset$ on the Lie algebra of
$U_\gamma$ admits only the weights $\gamma$ and $2\gamma$ and is maximal for this property.
\index[not]{UG@$U_\gamma$}

\subsection{Bruhat decomposition}

The Weyl group
\[
W(A_\emptyset)=N_G(A_\emptyset)/Z_G(A_\emptyset)=N_G(A_\emptyset)/M_\emptyset
\]
is denoted by $W_G$.\index[not]{W_G@$W_G$} We denote by $S_\emptyset$ the set of reflections $s_\alpha$
in $W_G$ associated to the simple roots $\alpha \in \Delta_\emptyset$.

\needspace{147pt}
\begin{prop}
[\cite{BoTi}] The quadruple $(G,P_\emptyset,N_G(A_\emptyset),S_\emptyset)$ is a Tits system. In
particular, we have the following decomposition of $G$ (Bruhat decomposition)
\[
G=\coprod_{w\in W_G} P_\emptyset w  P_\emptyset.
\]
\end{prop}

\subsection{}\label{Botits}
Let $T$ be a subset of $\Sigma_\emptyset$. We say that $T$ is closed if
$$(T+T)\cap \Sigma_\emptyset \subset T.$$ We say that $T$ is convex if $T$ is the intersection of
$\Sigma_\emptyset$ with a closed convex cone of $\fra^*$. A closed subset $T$ of $\Sigma_\emptyset$ is
said to be symmetric if $T=-T$. In this case, $T$ is a root system, and we denote by $W_T \subset W_G$
its Weyl group. A closed subset $T$ of $\Sigma_\emptyset$ is said to be unipotent if
$T \subset w\cdot \Sigma_\emptyset^+$ for some $w \in W_G$.

For any closed subset $T$ of $\Sigma_\emptyset$, we denote by $G(T)$ the algebraic subgroup of $G$
generated by $M_\emptyset$ and the $U_\gamma$, $\gamma \in T$. If $T$ is unipotent, we denote by $U(T)$
the subgroup of $G$ generated by the $U_\gamma$, $\gamma \in T$.

The following result is proved in \cite{BoTi}, 3.22.
\begin{prop}
Let $S$, $T$ be closed subsets of $\Sigma_\emptyset$.

$(i)$ If $S$ and $T$ are convex, then $G(S) \cap G(T)= G(S\cap T)$.

$(ii)$ If $T$ is unipotent, then $G(S) \cap U(T)= U(S\cap T)$.
\end{prop}

We call a closed subset $\caP$ of $\Sigma_\emptyset$ parabolic if $\Sigma_{\emptyset}^+ \subset \caP$. In
this case, we set $\caM=\caP \cap (-\caP)$, $\caN=\caP \setminus \caM$, and we say that
$(\caP,\caM,\caN)$ is a parabolic triplet \index[ter]{parabolic triplet} in $\Sigma_\emptyset$.

\begin{prop}[\cite{BoTi}, 5.12-5.18]
Let $\Theta$ be a subset of $\Delta_\emptyset$, and let $\caP$, $\caM$ be the closed subsets of
$\Sigma_\emptyset$ generated by $\Sigma_\emptyset^+ \cup (-\Theta)$ and $(-\Theta) \cup \Theta$
respectively. Then $(\caP,\caM,\caN=\caP \setminus \caM )$ is a parabolic triplet in $\Sigma_\emptyset$
and $P=MN$ is a standard parabolic subgroup of $G$, where $P=G(\caP)$, $M=G(\caM)$ and $N=U(\caN)$. Any
standard parabolic subgroup $P=MN$ of $G$ is obtained in this way. Moreover, the subset $\Theta$ is
uniquely determined by the data of $(\caP,\caM,\caN)$ or of $P=MN$. In particular, the subsets $\caP$,
$\caM$ and $\caN$ are convex.
\end{prop}

\subsection{}\label{PGQ1}
For any standard Levi subgroup $M$ of $G$, with split component $A$, let us set
\begin{equation}\label{WM}
W_M=N_M(A_\emptyset)/Z_M(A_\emptyset).
\end{equation}
Since
\[
Z_G(A_\emptyset)=M_\emptyset \subset M,
\]
we obtain a canonical injection $W_M \hookrightarrow W_G$. If $w \in W_G$, and if $g$ is a lift of $w$ in
$G$, we set, for any standard Levi subgroup of $G$, $w\cdot M=gMg^{-1}$. Since $M_\emptyset \subset M$,
it is clear that this is independent of the chosen lift. Let us set, for two standard Levi subgroups $M$
and $L$ of $G$,
\begin{equation}\label{WML}
W(L,M)=\{w \in W_G \mid w\cdot L=M  \}.
\end{equation}
\index[not]{WLM@$W(L,M)$} We of course have $W_M \cdot   W(L,M)\cdot W_L=W(L,M)$. Note that
\[
W(M,M)=(N_G(M)\cap N_G(A_\emptyset))/M_\emptyset,
\]
and that $M_\emptyset \subset M $. Therefore we have a well-defined map
\[
W(M,M) \rightarrow W(A)=N_G(A)/Z_G(A)= N_G(M)/M
\]
whose kernel is $W_M$. It follows that we have an injection
\begin{equation*}
W(M,M)/W_M \rightarrow W(A) ,
\end{equation*}
and we immediately verify that it is surjective because the maximal split tori are conjugate in $M$. Thus
we obtain
\begin{equation}\label{WMWA}
W(M,M)/W_M \simeq W(A) .
\end{equation}

\subsection{} \label{WoM}
Let $M$ be a standard Levi subgroup of $G$. Let us set
\[
W(\, *\, ,M)= \bigcup_L W(L,M),\quad W(M, \, *\, )= \bigcup_L W(M,L),
\]
and
\[
l(M)=|W_M\backslash W(\, *\, ,M)|=|W(M,\, *\, )/W_M|
\]
where the sum is over the set of standard Levi subgroups of $G$. If $M$ is a non-standard Levi subgroup,
we set $l(M)=l(M')$ where $M'$ is a standard Levi subgroup of $G$ conjugate to $M$.

\subsection{}\label{WMWGWM}
We note the following consequence of the Bruhat decomposition. We refer the reader to \cite{BoTi},
5.15-5.20 for a proof.

\begin{lemme}
Let $P=MN$ and $Q=LU$ be two semi-standard parabolic subgroups of $G$, where $M$ and $L$ are standard
Levi subgroups. Then the set of orbits of $Q$ in $X=P\backslash G$ is in bijection with the set of double
cosets $W_M\backslash W_G/W_L$.
\end{lemme}

We then equip $W_M\backslash W_G/W_L$ with the Bruhat order
\[
\bar{w_1} \leq_{PQ} \bar{w_2}  \text{ if }   P\bar{w_1} Q \subset
\overline { P\bar{w_2} Q}.
\]
\index[not]{01@$\leq_{PQ}$} In this situation, we then choose a total order $\preceq_{PQ}$ finer than the
Bruhat order on $ W_M\backslash W_G/W_L$.\index[not]{02@$\preceq_{PQ}$}

\begin{cor}
Let $P=MN$ and $Q=LU$ be two parabolic subgroups of $G$. Then the set of orbits of $Q$ in
$X=P\backslash G$ is finite.
\end{cor}
\begin{proof}
This is immediately deduced from the lemma above and the fact that any parabolic subgroup of $G$ is
conjugate to a standard parabolic subgroup.
\end{proof}

\subsection{}\label{WeylgroupsS}
Let $P$ and $Q$ be standard parabolic subgroups of $G$. We will describe an explicit system of
representatives of the double cosets $W_M\backslash W_G/W_L$ possessing useful properties. Let us set:
\begin{equation}\label{Wml}
 W^{L,M}=\{w \in W_G \mid w \cdot (L \cap P_\emptyset)
\subset P_\emptyset, \,  w^{-1}\cdot  (M \cap P_\emptyset) \subset P_\emptyset   \}.
\end{equation}
\index[not]{WLM1@$W^{L,M}$}

We then have
\begin{lemme}

$(i)$ In each double coset $W_M\backslash W_G /W_L$, there exists a unique element of $W^{L,M}$. It is
the element whose length is minimal in this double coset.

$(ii)$ If $w \in  W^{L,M}$, then $M \cap  w\cdot Q$ is a standard parabolic subgroup of $M$, with Levi
factor $M \cap w\cdot L$ and unipotent radical $M \cap  w\cdot U$ and $w^{-1}\cdot P \cap L$ is a
standard parabolic subgroup of $L$, with Levi factor $w^{-1}\cdot M \cap L$ and unipotent radical
$w^{-1}\cdot N \cap L$. We also have the decompositions
\[
P \cap  w\cdot Q=  (P \cap  w\cdot L)(P \cap  w\cdot U), \quad  Q
\cap w^{-1} \cdot P=  (Q\cap w^{-1} \cdot M) (Q\cap w^{-1} \cdot N),
\]
\[
N \cap  w\cdot Q=  (N \cap  w\cdot L)(N \cap  w\cdot U), \quad  U
\cap w^{-1} \cdot P=  (U\cap w^{-1} \cdot M) (U\cap w^{-1} \cdot N).
\]

$(iii)$ $ W^{L,M}\cap W(L,M)  \simeq W_M\backslash  W(L,M)/W_L\simeq W(L,M)/W_L.$
\end{lemme}

\begin{proof}
Let us start with $(ii)$. Let $\caP,\caM,\caN,\caQ,\caL,\caU$ be the subsets of $\Sigma_\emptyset$
corresponding via the bijection of the second proposition of \ref{Botits} respectively to the subgroups
$P,M,N,Q,L,U$ of $G$. Let us set $\caM^ +=\caM \cap \Sigma_\emptyset^+$ and
$\caL^ +=\caL \cap \Sigma_\emptyset^+$. We clearly have $W_M=W_\caM$ and $W_L=W_\caL$. Moreover,
\[
W^{L,M}=\{  w \in W_G \mid w\cdot \caL^+ \subset \Sigma_\emptyset^+, \,
w^{-1}\cdot \caM^+ \subset \Sigma_\emptyset^+\}.
\]
Let $w \in W^{L,M}$. Then $\caM \cap w \cdot \caQ$ contains $\caM^+$. It is therefore a parabolic subset
of $\caM$, and
\[
(\caM \cap w  \cdot \caQ,\caM \cap w \cdot \caL, \caM \cap w \cdot \caU)
\]
is a parabolic triplet of $\caM$. We deduce from the second proposition of \ref{Botits} that
\[
M\cap w  \cdot Q=(M\cap w  \cdot L)(M\cap w  \cdot U)
\]
is a standard parabolic subgroup of $M$. We show in the same way that
\[
L\cap   w^{-1} \cdot  P=(L\cap w^{-1} \cdot  M)(L\cap w^{-1} \cdot N)
\]
is a standard parabolic subgroup of $L$. The other decompositions are obtained without difficulty.

$(i)$ We use the following results, a proof of which can be found in \cite{St}, Append. I, II:

- the length $l(w)$ of an element $w \in W_G$ is equal to the number of roots $\gamma$ in
$\Sigma_\emptyset^+$ such that $w\cdot \gamma \in - \Sigma_\emptyset^+$.

- If $w \in W_\caL$, then
$w\cdot (\Sigma_\emptyset^+ \setminus \caL^+)\subset \Sigma_\emptyset^+ \setminus \caL^+$.

- If $w \in W_G$ and $\gamma \in \Delta_\emptyset$, then $l(ws_\gamma)>l(w)$ if and only if
$w \cdot \gamma \in   \Sigma_\emptyset^+ $.

It follows that
\begin{multline}
  W^{L,M}=\{  w \in W_G \mid \forall \gamma \in
\Delta_\emptyset  \cap \caL, \;    l(ws_\gamma)>l(w),\\
   \forall  \gamma \in  \Delta_\emptyset \cap \caM,\;
   l(w^{-1}s_\gamma)>l(w^{-1}) \} .
\end{multline}
The desired result then follows from \cite{B}, Chap. IV, \S 1, exercise 3. See also \cite{Carter}, Prop.
2.7.3.

$(iii)$ is obvious.
\end{proof}

\subsection{} \label{PWQsemi}
We now assume that $P$ and $Q$ are arbitrary parabolic subgroups of $G$. From the previous paragraph, we
deduce the existence of a system of representatives $\caW^{Q,P}$ \index[not]{WQP@$\caW^{Q,P}$} in $G$ of
the double cosets $P\backslash G/Q$ satisfying the following properties:

If $z \in  \caW^{Q,P}$, then $M \cap  z\cdot Q$ is a parabolic subgroup of $M$, with Levi factor
$M \cap z\cdot L$ and unipotent radical $M \cap  z\cdot U$, $z^{-1}\cdot P \cap L$ is a parabolic
subgroup of $L$, with Levi factor $z^{-1}\cdot M \cap L$ and unipotent radical $z^{-1}\cdot N \cap L$.
Moreover
\[
P \cap  z\cdot Q=  (P \cap  z\cdot L)(P \cap  z\cdot U), \quad  Q
\cap z^{-1} \cdot P=  (Q\cap z^{-1} \cdot M) (Q\cap z^{-1} \cdot N),
\]
\[
N \cap  z\cdot Q=  (N \cap  z\cdot L)(N \cap  z\cdot U), \quad  U
\cap z^{-1} \cdot P=  (U\cap z^{-1} \cdot M) (U\cap z^{-1} \cdot N).
\]

Indeed, the previous paragraph treats the case where $P$ and $Q$ are standard and where we take for
$\caW^{Q,P}$ a system of representatives in $G$ of the elements of $W^{L,M}$. But $P$ and $Q$ are
conjugate respectively to standard parabolic subgroups $P_1=M_1N_1$ and $Q_1=L_1U_1$, say by elements $g$
and $h$ respectively:
\[
P_1=g\cdot P, \quad Q_1=h \cdot Q.
\]
We also assume that
\[
M_1=g\cdot M,\quad L_1=h\cdot L.
\]
From the decomposition
\[
G=\coprod_{w \in \caW^{Q_1,P_1}} P_1 w Q_1
\]
of the previous paragraph, we deduce
\begin{small}
\[
G=\coprod_{w \in \caW^{Q_1,P_1}} P_1 w Q_1= \coprod_{w \in \caW^{Q_1,P_1}}
(g\cdot P)w (h\cdot Q)=\coprod_{w \in \caW^{Q_1,P_1}}
g(P(g^{-1}wh)Q)h^{-1}.
\]
\end{small}
We therefore have
\[
G= \coprod_{w \in \caW^{Q_1,P_1}} P(g^{-1}wh)Q.
\]
We may therefore set $\caW^{Q,P}=g^{-1}\caW^{Q_1,P_1}h  $. Let us now show that for all
$z=   g^{-1}wh  \in  \caW^{Q,P}  $, $M \cap  z\cdot Q$ is indeed a parabolic subgroup of $M$, with Levi
factor $M \cap z\cdot L$ and unipotent radical $M \cap  z\cdot U$. We have
\[
M \cap  z\cdot Q=(g^{-1}\cdot M_1)\cap (g^{-1}wh)\cdot (h^{-1}\cdot
Q_1)  = g^{-1}\cdot (M_1 \cap w \cdot Q_1),
\]
and therefore $M \cap  z\cdot Q$ is a parabolic subgroup of $M$. Moreover
\[
M \cap  z\cdot Q= g^{-1}\cdot ((M_1 \cap w \cdot L_1)(M_1\cap w\cdot U_1))=
(M \cap z\cdot L)(M \cap  z\cdot U).
\]
The remaining assertions are proved in the same way.

\subsection{Maximal Levi subgroups} \label{levmax}
\index[ter]{maximal!(Levi subgroup)} We say that a Levi subgroup is maximal if it is a Levi factor of a
proper maximal parabolic subgroup of $G$. The standard maximal Levi subgroups of $G$ correspond by the
second proposition of \ref{Botits} to the subsets of $\Delta_\emptyset$ of the form
\begin{equation*}\label{Thetaalpha}
 \Theta=\Delta_\emptyset \setminus \{\alpha\},
\end{equation*}
$\alpha$ being a simple root. It is not difficult to see using the considerations of Section
\ref{WeylgroupsS} that a Levi subgroup $M$ is maximal if and only if $l(M)=2$, where $l(M)$ is defined in
\ref{WoM}.

Suppose $M$ is standard given by the subset $\Delta_\emptyset \setminus \{\alpha\}$. Let $w$ be a
representative of the non-trivial element of $W(M,*)/W_M$ in $G$. Then if $P=MN$ is the standard
parabolic subgroup containing $M$, and if $P'$ is the standard parabolic subgroup containing
$M'=w\cdot M$, we have
\[
w \cdot P =\overline{P'}.
\]

\section{Compact subgroups}

The results of this section are proved in \cite{BrTi}.

\subsection{Maximal compact subgroups}\label{Cartan}

We fix a minimal parabolic subgroup $P_\emptyset=M_\emptyset N_\emptyset$ of $G$ as in the previous
section.
\begin{thm}
There exists a maximal compact subgroup $K_0$ of $G$, open, such that

$(1)$ $G=K_0 P_\emptyset=  P_\emptyset K_0 $.

$(2)$ For all $s \in W_G$, there exists a representative of $s$ in $K_0$.

$(3)$ Let $P=MN$ be a standard parabolic subgroup. Then $$K_0 \cap P = (K_0\cap M)(K_0\cap N).$$

$(4)$ (Cartan decomposition)\index[ter]{decomposition!Cartan}
\[
G=\coprod_{a \in \widetilde C_\emptyset^+ }  K_0 a  K_0=
\bigcup_{a \in C_\emptyset^+ ,f \in F_\emptyset^+}  K_0 af  K_0.
\]
(See \ref{FM0} for the notation.)

$(5)$ For any standard parabolic subgroup $P=MN$, the maximal compact subgroup $K_0 \cap M$ of $M$
satisfies properties $(1)$ to $(4)$ for the minimal parabolic subgroup $P_\emptyset \cap M$ of $M$.
\end{thm}

Such a maximal compact subgroup $K_0$ is said to be adapted to $A_\emptyset$.

\begin{cor}
If $P$ is any parabolic subgroup of $G$, it is conjugate under $K_0$ to a standard parabolic subgroup.
Moreover $G/P$ is compact.
\end{cor}

\subsection{Iwahori decomposition of compact open subgroups}\label{KN}
\index[ter]{decomposition!Iwahori} We fix $P_\emptyset$, $A_\emptyset$, $K_0$ as in the previous
paragraph.

\begin{thm}[\cite{Del}, 2.1]
There exists a neighborhood basis of the identity in $G$ consisting of compact open subgroups $K$ such
that:

$(i)$ $K$ is a normal subgroup of $K_0$.

$(ii)$ For any standard parabolic subgroup $P=MN$, we have a decomposition
\[
K=K_{\overline{N}}K_M K_N, \quad K_{\overline{N}}=K\cap \overline{N},\quad
K_N=K\cap N,\quad K_M=K\cap M
\]
called the Iwahori decomposition of $K$ with respect to $P$. Moreover, $C_A$ normalizes $K_M$, and for
all $m\in C_A^+$, we have $mK_Nm^{-1}\subset K_N$ and $m^{-1}K_{\overline{N}}m\subset K_{\overline{N}}$.

$(iii)$ Let $m \in   C_A^{++}$. Then $N=\bigcup_n m^{-n} K_N m^n$, and
$\overline{N}=\bigcup_n m^n  K_{\overline{N}} m^{-n}$. The unipotent radical of a parabolic subgroup is
therefore the union of its compact open subgroups. More generally, for any compact subset $\Gamma$ of $N$
(resp. $\overline{ \Gamma}$ of $\overline{N}$), there exists $\epsilon >0$ such that for all
$m \in C_A^+(\epsilon)$, $\Gamma\subset m^{-1}K_Nm$ (resp.
$\overline{ \Gamma }\subset  m K_{\overline{N}} m^{-1}$).

$(iv)$ With the notation of $(iii)$, the set of $N_i:= m^i  K_N m^{-i}$, $i\in \bbN$, forms a
neighborhood basis of the identity in $N$. Similarly the set of
$\overline{N}_i:= m^{-i}  K_{\overline{N}} m^i$, $i\in \bbN$, forms a neighborhood basis of the identity
in $\overline{N}$. More generally, for any compact open subgroup $\Gamma$ of $N$ (resp.
$\overline{ \Gamma}$ of $\overline{N}$), there exists $\epsilon >0$ such that for all
$m \in C_A^+(\epsilon)$, $mK_Nm^{-1} \subset \Gamma$ (resp.
$m^{-1}K_{\overline{N}}m \subset \overline{  \Gamma}$).
\end{thm}

\subsection{Decomposition of the Hecke algebra $\caH(G,K)$} \label{HCDH}
We fix $P_\emptyset$, $A_\emptyset$, $K_0$ as in the previous paragraph. Recall that if $K$ is a compact
open subgroup of $G$, and if $g$ is an element of $G$, we defined in the second remark of \ref{HGK} the
distribution $a_{g,K}:=e_K*\delta_g*e_K$ in $\caH(G,K)$. When there is no risk of confusion, for example
if $K$ is fixed, we will denote this distribution simply by $a_g$.

Let us fix a compact open subgroup $K$ of $G$, contained and normal in $K_0$ and admitting an Iwahori
decomposition with respect to the standard parabolic subgroups. Let us choose a system of representatives
$\{ x_1,\ldots, x_r \}$ for the cosets of $K$ in $K_0$.

Let us set $\caH_0:= \caH(K_0,K)$. It is clear that $\caH_0$ is a finite-dimensional subalgebra of
$\caH(G,K)$ of which a basis is $\{a_{x_i}\}_{i=1,\ldots,r}$.

\begin{lemme}
For all $g \in G$, we have:
\[
a_ga_{x_i}=a_{gx_i}, \quad a_{x_i}a_g=a_{x_ig}.
\]
\end{lemme}

\begin{proof}
Since $K$ is normal in $K_0$, $e_K*\delta_{x_i}=\delta_{x_i}*e_K$, hence
\begin{align*}
a_ga_{x_i}&=(e_K*\delta_g*e_K)*(e_K*\delta_{x_i}*e_K)=
e_K*\delta_g*e_K*\delta_{x_i}*e_K \\
 &=e_K*\delta_g*\delta_{x_i}*e_K=a_{gx_i}.
\end{align*}
Similarly, $a_{x_i}a_g=a_{x_ig}$.
\end{proof}

\begin{lemme}
Let $z_1,z_2 \in   C_\emptyset^+$. Then $a_{z_1}a_{z_2}=a_{z_1z_2}$.
\end{lemme}
\begin{proof}
This is not completely trivial because the elements of $C_\emptyset^+ $ do not normalize $K$. We will
therefore use the Iwahori decomposition of $K$ with respect to $P_\emptyset$ which we write simply as
$K=K_{\overline{  N}}K_MK_N$ to lighten the notation. Therefore, we have the following equality in
$\scrE'(G)$ (cf. Proposition \ref{Gammaprop})
\[
e_K=e_{K_{\overline{N}}}*e_{K_M}*e_{K_N} =e_{K_{N}}*e_{K_M}*e_{K_{\overline {N}}}.
\]
The second equality being obtained by taking the inverse. We moreover have
\[
z_1K_Nz_1^{-1}\subset K_N\subset K, \quad z_2^{-1}K_{\overline{N}}z_2\subset
K_{\bar N}\subset K,\quad  z_1K_Mz_1^{-1}= K_M\subset K
\]
hence:
\[
e_{K}*e_{  z_1K_Nz_1^{-1}}=e_{K},\quad   e_{z_2^{-1}K_{\bar N}z_2}*e_K=e_K,
\quad   e_K*e_{K_M}=e_K.
\]
Finally
\begin{align*}
a_{z_1}a_{z_2}&= (e_K*\delta_{z_1}*e_K)*(e_K*\delta_{z_2}*e_K)=
e_K*\delta_{z_1}*e_K*\delta_{z_2}*e_K\\
&= e_K*\delta_{z_1}*e_{K_{N}}*e_{K_M}*e_{K_{\overline{N}}}*\delta_{z_2}*e_K\\
&=  e_K*e_{z_1 K_N z_1^{-1} } * \delta_{z_1} * e_{K_M} * \delta_{z_2}
* e_{ z_2^{-1}K_{\overline{N}}z_2 } * e_K\\
&=e_K*\delta_{z_1}*e_{K_M}* \delta_{z_2}*e_K
=e_K*e_{K_M}*\delta_{z_1}* \delta_{z_2}*e_K\\
&=e_K*\delta_{z_1z_2}*e_K=a_{z_1z_2}.
\end{align*}
\end{proof}

Let $\caA$ be the subalgebra of $\caH(G,K)$ generated by the $a_z$, $z\in C_\emptyset^+$. Since
$C_\emptyset ^+$ is an abelian monoid, we deduce:
\begin{cor}
The subalgebra $\caA$ is abelian: for all $z_1,z_2 \in  C_\emptyset^+$, $a_{z_1}a_{z_2}=a_{z_2}a_{z_1}$.
\end{cor}

Let $\caM$ denote the sub-$\caA$-module of $\caH(G,K)$ generated by the elements $a_{zf}$ with
$z\in  C_\emptyset^+$ and $f\in F_\emptyset^{+}$.

\begin{lemme}
The $\caA$-module $\caM$ is finitely generated.
\end{lemme}

\begin{proof}
Let $z_1,z_2\in C_\emptyset^+$ and $f\in F_\emptyset^{+}$
\begin{align*}
a_{z_1}a_{z_2f}&= (e_K*\delta_{z_1}*e_K)*(e_K*\delta_{z_2f}*e_K)=
e_K*\delta_{z_1}*e_K*\delta_{z_2f}*e_K\\
&= e_K*\delta_{z_1}*e_{K_{N}}*e_{K_M}*e_{K_{\overline{N}}}*\delta_{z_2f}*e_K=
e_K*\delta_{z_1}*e_{K_{\overline{N}}}*\delta_{z_2f}*e_K\\
&= e_K*\delta_{z_1z_2}*\delta_{z_2^{-1}}*e_{K_{\overline{N}}}*\delta_{z_2}*
\delta_f*e_K= e_K*\delta_{z_1z_2}*e_{z_2^{-1}K_{\overline{N}}z_2}* \delta_f*e_K
\end{align*}
On the other hand
\begin{align*}
& e_K* \delta_{z_1z_2}* e_{z_2^{-1}z_1^{-1}K_{\overline{N}} z_1z_2} *  \delta_f
* e_K=e_K * e_{K_{\overline{N}}}*\delta_{z_1z_2}* \delta_f*e_K=
e_K*\delta_{z_1z_2f}*e_K=a_{z_1z_2f}
\end{align*}
We conclude that if $  z_2^{-1}K_{\overline{N}}z_2fK=z_2^{-1}z_1^{-1}K_{\overline{N}} z_1z_2fK$, then
$a_{z_1}a_{z_2f}=a_{z_1z_2f}$. Note that the groups $z_2^{-1}K_{\overline{N}}z_2$ and
$z_2^{-1}z_1^{-1}K_{\overline{N}} z_1z_2$ are subgroups of $K$. For any subgroup $\Gamma$ of $K$, let us
set
\[
\Vert \Gamma \Vert =\sum_{f\in F_\emptyset^{+}}  \lvert  \Gamma f K  /K\rvert
\]
where $\vert  \Gamma f K  /K\vert$ denotes the number of orbits of the action of $K$ on $\Gamma f K$. It
is clear that if $\Gamma'\subset \Gamma$ then $\Vert \Gamma' \Vert \leq \Vert \Gamma \Vert $, and if
$\Vert \Gamma' \Vert = \Vert \Gamma \Vert $, then $ \Gamma' f K = \Gamma f K $ for all
$f\in F_\emptyset^{+}$. For all $z\in C_\emptyset^+$, let us set
\[
\nu(z)= \Vert    z^{-1}K_{\overline{N}}z  \Vert
\]
Recall that in \ref{SSM} we set $S=C_\emptyset^+/C_{A_G}$. Since $C_{A_G}$ is central in $G$, it is clear
that $\nu(z)$ depends only on the class of $z$ in $S$. We also saw that $S$ is the union of the
translates $rS_0$, $r\in R$, where $R$ is finite and $S_0\simeq\bbN^d$. Fix $r\in R$; by transport of
structure, $\underline n\mapsto \nu(r\prod_\alpha t_\alpha^{n_\alpha})$ defines a function, still denoted
$\nu$, on $\bbN^d$. We equip $\bbN^d$ with the partial order defined by
$\underline n_1\leq \underline n_2$ if $\underline n_2-\underline n_1\in \bbN^d$. If the elements $z_1$
and $z_2$ of $C_\emptyset ^+$ correspond to elements $\underline n_1$ and $\underline n_2$ of $\bbN^d$
with $\underline n_1\leq \underline n_2$, it is then clear that
$z_2^{-1}K_{\overline{N}}z_2 \subset z_1^{-1}K_{\overline{N}}z_1$ and therefore
$\nu(\underline n_2)\leq \nu(\underline n_1)$. We say that $\underline n \in \bbN^d$ is critical for
$\nu$ if for all $\underline n' \in \bbN^d$ satisfying $\underline n'< \underline n $ we have
$\nu(\underline n')>\nu(\underline n)$. Since the $\caA$-module $\caM$ is generated by the $a_{zf}$,
$z\in C_\emptyset^+$, $f\in F_\emptyset^{+}$, from the above we see that the $\caA$-module $\caM$ is
generated by the $a_{zf}$ where $z$ now runs through a set of representatives of the classes in
$S=C_\emptyset^+/C_{A_G}$ corresponding, for some $r\in R$, to an element $\underline n$ critical for the
function $\nu$ attached to $r$. It therefore remains to show that the set of critical points for $\nu$ in
$\bbN^d$ is finite.

This reduces to a purely combinatorial problem. For a critical $\underline n\in \bbN^d$, we have
$0\leq \nu(\underline n) < \nu(0)$ and by induction on $\nu(0)$, we show that the number of critical
points in $\underline n+\bbN^d$ is finite. The complement of $\underline n+\bbN^d$ in $\bbN^d$ is covered
by a finite number of submonoids of rank less than or equal to $d-1$ (isomorphic to $\bbN^r$ for
$r\leq d-1$). By induction on $d$, we show that there is a finite number of critical points in these
submonoids. Finally, this shows that there are indeed only finitely many critical points in $\bbN^d$.
\end{proof}

Let $\{D_j\}_{j=1,\ldots,q}$ be a finite generating system of the $\caA$-module $\caM$. Recall that we
set $\caH_0=\caH(K_0,K)$.
\begin{thm}
The Hecke algebra $\caH(G,K)$ decomposes into:
\[
\caH(G,K)=\sum_{j=1}^q\caH_0\caA D_j \caH_0.
\]
\end{thm}
\begin{proof}
Let us start from the Cartan decomposition
\[
G=K_0 C_\emptyset^+  F_\emptyset^{+} K_0.
\]
Since $K_0=\bigcup_{i=1}^r Kx_i= \bigcup_{i=1}^r x_iK$, we deduce
\[
G=\bigcup_{i,j=1}^r  Kx_i   C_\emptyset^+  F_\emptyset^{+}   x_j K.
\]
A generating family of $\caH(G,K)$ is therefore given by the
\[
\{ a_{x_i  z f  x_j} \}, i,j=1,\ldots, r, f \in F_\emptyset^{+} ,
z\in   C_\emptyset^+  .
\]
Now $a_{x_ifzx_j}=a_{x_i} a_{zf}  a_{x_j}$ and $a_{zf} \in  \caM =  \sum_{j=1}^q \caA D_j $ according to
the lemmas above.
\end{proof}

\begin{rmq}
By considering the right $\caA$-module structure of $\caM$ and adapting the proofs above, we also obtain
a decomposition of $\caH(G,K)$ in the form
\[
\caH(G,K)=\sum_{j=1}^{q'}\caH_0  E_j \caA  \caH_0,
\]
where the $E_j$ are elements of $\caM$.
\end{rmq}

\subsection{Calculations of modular functions}\label{calcfonctmod}

We have fixed a left Haar measure $\mu_G$ on $G$.
\begin{prop}
The group $G$ is unimodular.
\end{prop}
\begin{proof}
We need to show that the modular function $\delta_G$ defined in \ref{foncmod} is identically equal to
$1$. According to the Cartan decomposition of $G$, every element of $G$ can be written $k_1zfk_2$, with
$k_1,k_2\in K_0$, $z\in C_\emptyset\subset A_\emptyset$ and $f\in F_\emptyset$. Since $\delta_G$ is a
character taking values in $\bbR^\times_+$, it is trivial on any compact group; moreover,
$X_*(A_\emptyset)$ being of finite index in $\Lambda_\emptyset$, some power $f^m$ of $f$ lies in
${}^0M_\emptyset A_\emptyset$, with ${}^0M_\emptyset$ compact, so that $\delta_G(f)$ is determined by the
values of $\delta_G$ on $A_\emptyset$. It therefore suffices to show that $\delta_G(a)=1$ for all
$a\in A_\emptyset$. Let $K$ be a compact open subgroup of $G$ admitting an Iwahori decomposition with
respect to $P_\emptyset$. Let us also choose a compact open subgroup $K_1$ contained in
$K \cap a^{-1} Ka$, also admitting an Iwahori decomposition with respect to $P_\emptyset$. Let us
simplify the notation by setting $P=P_\emptyset=MN$. Let us write these Iwahori decompositions:
\[
K=K_{\overline{N}}K_M K_N,\quad  K_1=K_{1,\overline{N}}K_{1,M} K_{1,N}.
\]

A formula for $ \delta_{G}(a)$ is given by equation (\ref{fmodK}):
\[
\delta_{G}(a)=\frac{[a^{-1}Ka:K_1]}{[K:K_1]}
\]
Now,
\begin{align*}[K:K_1]
&=[K_{\overline{N}}K_M K_N:K_{1,\overline{N}}K_{1,M} K_{1,N} ]\\
&=[K_{\overline{N}} :K_{1,\overline{N}}][K_M: K_{1,M} ][K_N: K_{1,N}],
\end{align*}
and similarly
\begin{align*}
 [a^{-1}Ka:K_1]&=[a^{-1}K_{\overline{N}}K_M K_Na:K_{1, \overline{N} }K_{1,M} K_{1,N} ]\\
&=[a^{-1}K_{\overline{N}}a :K_{1,\overline{N}}][a^{-1}K_Ma: K_{1,M} ][a^{-1}K_Na: K_{1,N}].
\end{align*}
But $a^{-1}K_Ma=K_M$, hence finally
\[
\delta_{G}(a)= \frac{[a^{-1}K_{\overline{N}}a:K_{1,\overline{N}}]}{[K_{ \overline{N} }:K_{1,\overline{N}}]}
\frac{[a^{-1}K_{ N}a:K_{1,N}]}{[K_{N}:K_{1,N}]}.
\]
Recall that $\overline{N}$ is generated by the unipotent subgroups $U_{-\gamma}$,
$\gamma \in \Sigma_\emptyset^+$. If we denote by $m_\gamma$ the multiplicity of the root $\gamma$ in
$\Sigma'_\emptyset$ and $m_{2\gamma}$ the multiplicity of the root $2\gamma$, we therefore have
\[
\frac{[a^{-1}K_{\overline{N}}a:K_{1,\overline{N}}]}{[K_{\overline{N}}:K_{1, \overline{N} }]}=
\prod_{\gamma \in \Sigma_\emptyset^+ }|\gamma(a)^{m_\gamma+2m_{2\gamma}}|_\bbF.
\]
and similarly
\[
\frac{[a^{-1}K_{ N}a:K_{1,N}]}{[K_{N}:K_{1, N}]}=\prod_{\gamma \in
\Sigma_\emptyset^+ }|\gamma(a)^{-(m_\gamma+2m_{2\gamma})}|_\bbF.
\]
This shows that $\delta_G(a)=1$.
\end{proof}

\bigskip

We will now give another expression for $\delta_{P\backslash G}(a)$. Let $\mu_P$ be a left Haar measure
on $P$ and $\delta_P$ the corresponding modular function.

\begin{lemme}
With the above notation, for all $m \in M$,
\[
\delta_{P\backslash G}(m)=\frac{\delta_{G}(m)}{\delta_{P}(m)}=
\delta_P(m)^{-1}=\delta_N(m)^{-1}=\delta_{\overline{N}}(m),
\]
and if $a\in A$, we moreover have
\[
\delta_{P\backslash G}(a) = \prod_{\gamma \in \Sigma(P) }
|\gamma(a)^{m_\gamma+2m_{2\gamma}}|_\bbF.
\]
\end{lemme}

\begin{proof}
We have by definition, using the notation introduced above,
\[
\int_{P} \chi_{K_MK_N}(apa^{-1}) \, d\mu_P(p)   =    \delta_P(a) \int_{P}
\chi_{K_MK_N}(p) \, d\mu_P(p),
\]
i.e., $\mu_P(a^{-1}K_MK_Na)=\delta_P(a)\mu_P(K_MK_N)$.

Since $K_MK_N$ is a compact subgroup of $P$, equation (\ref{fmodK}) gives here
\[
\delta_{P}(a)=\frac {[a^{-1}K_MK_Na: K_{1,M} K_{1,N}] }{ [K_MK_N: K_{1,M}K_{1,N}]}
\]
Moreover we have
\begin{align*}
  [K_MK_N: K_{1,M}K_{1,N}]&=[K_M:K_{1,M}][K_N: K_{1,N}], \\
[a^{-1}K_MK_Na: K_{1,M}K_{1,N}]&=[a^{-1}K_Ma:K_{1,M}][a^{-1}K_Na: K_{1,N}].
\end{align*}
Since $a^{-1}K_Ma=K_M$, we obtain
\[
\delta_{P}(a)= \frac {[a^{-1}K_Na:  K_{1,N}] }{ [K_N:K_{1,N}]}
\]
As above, $N$ is generated by the unipotent subgroups $U_\gamma$, $\gamma \in \Sigma(P)$. A calculation
similar to the one done above shows that we then have
\[
\delta_{P}(a) = \prod_{\gamma \in \Sigma(P) }
|\gamma(a)^{-(m_\gamma+2m_{2\gamma})}|_\bbF.
\]
Note that this calculation shows in passing that
\[
\delta_{P}(a) =\delta_N(a)=\delta_{\overline{N}}(a)^{-1}.
\]
for all $a \in A$. By the same reasoning as above using the Cartan decomposition (but this time for the
group $M$ rather than $G$), we deduce the equalities for $m \in M$.
\end{proof}

\subsection{Proof of Langlands' combinatorial lemma}

Langlands' combinatorial lemma is based on the properties of projection onto a cone in a Euclidean space.
Let $V$ be such a space, finite-dimensional, where we denote by $\bilo$ the inner product. Let us fix an
open cone $C$ of $V$, and denote by $\check{}\, C$ its dual
\[
\check{}\, C=\{ v \in V \mid \forall w \in C, \bil{v}{w}>0  \} .
\]
Let $\overline{C}$ and $\overline{\check{}\, C}$ denote the respective closures of $C$ and
$\check{}\, C$.

\begin{prop}
For all $v\in V$, there exists a unique element $v_0=p_{\overline{C}}(v)$ of $\overline{C}$, called the
projection of $v$ onto $\overline{C}$, satisfying

$(i)$ $\forall w \in\overline{C}$, $||v-v_0||\leq ||v-w||$.

Moreover, the element $v_0$ is characterized by one of the following properties, equivalent to $(i)$,

$(ii)$ $\forall w \in\overline{C}$, $\bil{v-v_0}{w-v_0}\leq 0$,

$(iii)$ the vector $v_0-v$ is in $\overline{\check{}\, C}$ and $\bil{v-v_0}{v_0}= 0$.
\end{prop}

The existence of the projection is a standard result for any closed convex subset of a Hilbert space. The
equivalence of conditions $(i)$, $(ii)$ and $(iii)$ for cones is elementary and left to the reader.

\bigskip

Let us now fix a basis $\{\alpha_1,\ldots ,\alpha_n\}$ of $V$, and denote by
$\{\omega_1,\ldots ,\omega_n\}$ the dual basis. We henceforth assume that $C$ is the open cone generated
by the vectors $\omega_1,\ldots ,\omega_n$, i.e.,
\[
C=\{ v=\sum_{i=1}^n t_i \, \omega_i, \; t_i>0   \}.
\]
We then have
\[
\check{}\, C=\{ w=\sum_{i=1}^n s_i \, \alpha_i, \; s_i>0   \}.
\]

\begin{cor}
For all $v \in V$, there exists a unique subset $F(v)$ of $\{1,\ldots,n\}$ such that
\[
v= \sum_{i\notin F(v)} t_i \, \omega_i-\sum_{j \in F(v)} s_j \, \alpha_j,
\]
with $t_i>0$, $s_j\geq 0$. We then have
\[
v_0=p_{\overline{C}}(v)=  \sum_{i\notin F(v)} t_i \, \omega_i.
\]
\end{cor}

\begin{proof}
Let $v \in V$. Its projection $v_0$ onto the cone $\overline{C}$ is written
\[
v_0=\sum_{i=1}^n t_i \, \omega_i
\]
where the $t_i$ are non-negative. Let $F(v)=F(v_0)$ be the set of $i$ such that $t_i$ is zero, so that
\[
v_0=\sum_{i \notin F(v)} t_i \, \omega_i, \quad t_i>0.
\]
We decompose $v-v_0$ in the basis $\{\alpha_1,\ldots ,\alpha_n\}$:
\[
v_0-v= \sum_{j=1}^n s_j \, \alpha_j.
\]
According to $(iii)$ of the proposition, $v_0-v$ is in $\overline{\check{}\, C}$, which translates to
$s_j\geq 0$, and on the other hand $v_0-v$ is orthogonal to $v_0$, which gives $s_j =0$ if
$j\notin F(v)$.

Conversely, if $v=\sum_{i\notin F(v)} t_i \, \omega_i - \sum_{j \in F(v)} s_j \, \alpha_j$, we see that
$v_0=p_{\overline{C}}(v)=  \sum_{i\notin F(v)} t_i \, \omega_i$, which shows the uniqueness of $F(v)$.
\end{proof}

\bigskip

We retain the notation of $\ref{plusam}$. The vector space $V$ is now $(\fra_M^G)^*$, of which a basis is
$\Delta(P)=\{ \alpha_1,\ldots ,\alpha_l\}$, the $\alpha_i$ being the restrictions to $\fra_M$ of the
roots in $\Delta_\emptyset \setminus \Delta_\emptyset^M$. Let $\{\omega_1,\ldots ,\omega_l\}$ denote the
dual basis. The cone $C$ above is therefore the cone denoted ${}_P^G[(\fra_M^G)^*]^+$ in \ref{plusam},
the cone $\check{}\, C$ being the one denoted ${}^+[(\fra_M^G)^*]_P^G$.

Let $\mu \in (\fra_M^G)^*$. From the above, there exists a subset $F(\mu)$ of $\{1,\ldots ,l\}$ such that
we can write $\mu$ as
\[
\mu= \sum_{i\notin F(\mu)} t_i \, \omega_i-\sum_{j \in F(\mu)} s_j \, \alpha_j,
\]
with $t_i>0$ and $s_j\geq 0$. Now such a subset $F(\mu)$ of $\{1,\ldots ,l \}$ determines a unique
standard parabolic subgroup $Q=LU$ of $G$ containing $P$ such that $\Delta(Q)$ is the restriction to
$\fra_L$ of the roots $\alpha_i$ with $i$ in $\{1,\ldots ,l\}\setminus F(\mu)$.

We set $\mu^+= \sum_{i\notin F(\mu)} t_i \, \omega_i$ and $\mu^-= -\sum_{j \in F(\mu)} s_j \, \alpha_j$.
It is clear that $\mu^+ \in {}_Q^G[(\fra_L^G)^*]^+$ and
$\mu^- \in - p_M^L({}^+\overline{[\fra_M^*]}_P^G)$. Moreover, by the above corollary, $\mu^+$ is the
orthogonal projection of $\mu-\mu_G$ onto the cone $\overline C={}_P^G\overline{[(\fra_M^G)^*]}^+$. Let
us show that $\mu^-$ is the orthogonal projection of $\mu-\mu_G$ onto the cone
$-\overline{\check{}\,C}={}^-\overline{[\fra_M^*]}_P^G$. For $w=-\sum_j u_j\alpha_j$, with $u_j\geq0$, we
have $\bil{\mu^+}{w}=-\sum_{i\notin F(\mu)}t_iu_i\leq0$, since $\bil{\omega_i}{\alpha_j}=\delta_{ij}$,
and $\bil{\mu^+}{\mu^-}=0$ for the same reason. Hence
\[
\bil{(\mu-\mu_G)-\mu^-}{w-\mu^-}=\bil{\mu^+}{w}-\bil{\mu^+}{\mu^-}\leq0,
\]
which is the characterization $(ii)$ of the projection, applied to the cone $-\overline{\check{}\,C}$.
Finally, since $\mu_G$ is orthogonal to $(\fra_M^G)^*$, these projections are also those of $\mu$.

Let us finally show that we obtain the same decomposition by viewing $\mu$ as an element of
$\fra_\emptyset^*$ and applying the above to $(\mu,P_\emptyset)$. Let $p$ denote the orthogonal
projection of $\fra_\emptyset^*$ onto $\fra_M^*$; the restriction to $\fra_M$ of an element of
$\fra_\emptyset^*$ identifies with its projection by $p$. Let us first show that, for
$\alpha\in\Delta_\emptyset\setminus\Delta_\emptyset^M$, $p(\alpha)-\alpha$ is a linear combination with
non-negative coefficients of elements of $\Delta_\emptyset^M$. Set $\beta=\alpha-p(\alpha)$; it is an
element of the subspace spanned by $\Delta_\emptyset^M$, and since $p(\alpha)$ is orthogonal to this
subspace, we have $\bil{\beta}{\gamma}=\bil{\alpha}{\gamma}\leq 0$ for all $\gamma\in\Delta_\emptyset^M$,
since two distinct simple roots form an obtuse or right angle (\cite{B}, Chap. VI, \S 1, no.~5, Remark 3
following Theorem 2). Write $\beta=\beta_+-\beta_-$, where $\beta_+$ and $\beta_-$ are linear
combinations with non-negative coefficients of elements of $\Delta_\emptyset^M$, with disjoint supports.
Then $\bil{\beta}{\beta_+}\leq0$ and $\bil{\beta_-}{\beta_+}\leq 0$, hence
$\bil{\beta_+}{\beta_+}=\bil{\beta}{\beta_+}+\bil{\beta_-}{\beta_+}\leq 0$, and $\beta_+=0$. Thus
$p(\alpha)=\alpha+\beta_-$ is a linear combination with non-negative coefficients of elements of
$\Delta_\emptyset$, which lie in $\Delta^L_\emptyset$ if $\alpha\in\Delta^L_\emptyset$.

Let us return to the decomposition $\mu=\mu_G+\mu^++\mu^-$ of $\mu\in\fra_M^*$ relative to $P$. We have
$\mu^-=-\sum_{j\in F(\mu)}s_j\alpha_j$, where the $\alpha_j$, $j\in F(\mu)$, are the restrictions to
$\fra_M$ of the roots of $\Delta^L_\emptyset\setminus\Delta^M_\emptyset$. By the above, $\mu^-$ is the
opposite of a linear combination with non-negative coefficients of elements of $\Delta^L_\emptyset$,
hence $\mu^-\in-p^L_{M_\emptyset}({}^+\overline{[\fra_\emptyset^*]}_{P_\emptyset}^G)$. Since
$\mu^+\in{}_Q^G[(\fra_L^G)^*]^+$ does not depend on $P$, the decomposition $\mu=\mu_G+\mu^++\mu^-$
satisfies the conditions of the lemma for $(\mu,P_\emptyset)$; by uniqueness, it is the one obtained for
$(\mu,P_\emptyset)$.
\qed

\chapter[Representations of $p$-adic groups]{Representations of $p$-adic reductive groups}\label{VI}

This chapter constitutes the core of the book. It presents the theory of the ``Bernstein center'', which
gives an explicit description of the center of the category $\caM(G)$ of smooth representations of a
$p$-adic reductive group $G$. Defined abstractly, the center of an abelian category $\caA$ is the set,
equipped with a ring structure, of natural transformations from the identity functor of the category
$\caA$ to itself. In the case where $\caA$ is the category of non-degenerate modules over a
$\bbC$-algebra with idempotents $A$, we gave in Chapter I a description of this center as the center of
the completed ring $\overline{A}$. When $\caA=\caM(G)$, the category of smooth representations of a
totally disconnected group $G$, the equivalence of categories $\caM(G)\simeq \caM(\caH(G))$ established
in Chapter III allows us to give a description of the center of this category as a convolution algebra of
essentially compact invariant distributions. This description is of a geometric nature. Bernstein's
theory gives another one, of a ``spectral'' nature. To clarify this concept, it may be useful to recall
what happens for finite groups. If $G$ is a finite group, the category of representations of $G$ is
equivalent to the category of unital $\bbC[G]$-modules, $\bbC[G]$ being the algebra (over $\bbC$) of the
group $G$, which is viewed as the algebra of complex-valued functions on $G$, equipped with the
convolution product. In this case, we know that the center of the category is naturally identified with
the center of the ring $\bbC[G]$, i.e., with the algebra $\bbC[G]^G$ of constant functions on the
conjugacy classes of $G$. Recall there are two natural bases for $\bbC[G]^G$: one coming from geometry,
consisting of the characteristic functions of the conjugacy classes of $G$, and the other, spectral,
consisting of the characters of the irreducible representations. We have seen the analogue of the first
of these bases for t.d. groups. The theory of the Bernstein center gives the analogue for $p$-adic
reductive groups of the second.

The structure of $p$-adic reductive groups recalled in the previous chapter highlights the fact that such
a group possesses remarkable subgroups, the Levi subgroups, which are again $p$-adic reductive groups.
The idea underlying the entire representation theory of these groups since its origin is that one can
study the representations of the group $G$ via the parabolic induction functors, defined by the general
constructions of Chapter III, between categories of smooth representations of Levi subgroups and the
category of smooth representations of $G$.

We therefore begin by defining these functors, denoted $i_P^G$, where $P=MN$ is a parabolic subgroup of
$G$, as a composition of functors studied in Chapter III. We immediately deduce the existence of a left
adjoint, denoted $r_P^G$, given explicitly (parabolic restriction functor, also called Jacquet functor),
and the existence of a right adjoint. The determination of this right adjoint as being the Jacquet
functor $r_{\bar P}^G$ is a deep result, called Bernstein's second adjunction theorem, and will be
established much later in the chapter. These functors are normalized so that $i_P^G$ preserves the
unitarity of representations. Moreover these functors are exact, and the induction functors are
additionally faithful.

Having defined the induction functors, the notion of supercuspidal representation naturally emerges.
Heuristically, these representations are those that do not come from Levi subgroups via the induction
functors. More precisely, a representation is supercuspidal if its image under all Jacquet functors
$r_P^G$ is zero, where $P$ runs through the set of proper parabolic subgroups of $G$. An irreducible
smooth representation of $G$ is either supercuspidal, or it appears as a subrepresentation of an induced
representation from a Levi subgroup. Thus a good part of the representation theory of reductive groups
boils down to two problems: the study of supercuspidal representations and the study of parabolic
induction functors.

A first result, Harish-Chandra's theorem, characterizes supercuspidal representations by a support
condition on their matrix coefficients: their support is compact modulo the center $Z(G)$ of $G$. We thus
almost recover the definition of compact representations from Chapter IV, with the nuance introduced by
the center of the group. It is therefore important to understand the consequences of this nuance, and to
see which results from the theory of compact representations can be preserved. First of all, like compact
representations, irreducible supercuspidal representations are admissible. We deduce the admissibility of
all irreducible smooth representations from the fact that parabolic induction functors preserve
admissibility. This result is a fairly immediate generalization of the case of compact representations,
but it can be crucially improved by exploiting the structure of reductive groups via a uniform
admissibility theorem: if $G$ is a $p$-adic reductive group, and $K$ a compact open subgroup of $G$,
there exists a positive constant $c$ such that the dimension of $V^K$, for any irreducible smooth
representation $V$ of $G$, is bounded by $c$. This fact forms the basis of all subsequent developments in
Bernstein's theory. It leads to a first decomposition of the category $\caM(G)$, starting from the
category decompositions obtained in Chapter IV using compact representations. We must take into account
the role played by the center, when it is non-compact, and introduce an equivalence relation on
supercuspidal representations, whose classes are called inertial classes: two representations are in the
same class if they differ by an unramified character of $G$. Since the set of unramified characters of
$G$ is equipped with a complex algebraic variety structure (it is a torus), and since an inertial class
$D$ is isomorphic (non canonically, one must choose a base point) to a quotient of this torus by a finite
subgroup, $D$ also admits a complex algebraic variety structure, and it is again a torus. Having fixed an
inertial class $D$ of irreducible supercuspidal representations of $G$, any smooth representation
$(\pi,V)$ can be written as a direct sum of a representation all of whose irreducible subquotients are in
$D$, and a representation none of whose irreducible subquotients are in $D$. But the uniform
admissibility theorem, through the finiteness criterion ({\bf KF}) of Chapter IV, implies a much finer
decomposition of $\caM(G)$: any smooth representation $(\pi,V)$ can be written as a direct sum of a
supercuspidal representation, itself a direct sum over the inertial classes $D$ of representations all of
whose irreducible subquotients are in $D$, and a representation none of whose irreducible subquotients
are supercuspidal. This constitutes a first step in Bernstein's decomposition theorem, which describes
$\caM(G)$ as a product of indecomposable categories. It remains to decompose the non-cuspidal part, which
we will do by a careful study of the parabolic induction functors. Before doing so, let us completely
describe the supercuspidal factors of the decomposition. Let us fix an inertial class $D$ of irreducible
supercuspidal representations of $G$, and let $\caM(G)_D$ denote the full subcategory of smooth
representations all of whose irreducible subquotients are in $D$. This category admits a small (i.e.,
finitely generated) progenerator, which we exhibit, and which by general arguments of category theory,
yields an equivalence of categories between $\caM(G)_D$ and the category of right modules over the
endomorphism ring of this progenerator. The structure of this ring is well understood: it is a
$\bbC$-algebra of Laurent polynomials, slightly twisted, i.e., it is commutative only up to
multiplication by scalar factors. Furthermore, the center of this algebra is indeed an algebra of Laurent
polynomials. It is, in fact, the algebra of polynomial functions on the variety $D$. An element $z$ of
the center of $\caM(G)_D$ is identified with a function on $D$, and the evaluation of this function at a
point $\pi$ of $D$ is the scalar by which, according to Schur's lemma, $z$ acts in the space of the
representation $\pi$. This completes the analysis of the supercuspidal part of $\caM(G)$. Note
incidentally that it follows from the description of the center of $\caM(G)_D$ as an algebra of Laurent
polynomials that it is a Noetherian algebra.

The key result allowing us to analyze the induced part of $\caM(G)$ is the geometric
restriction-induction lemma, which analyzes the composition of the functors $r_Q^G \circ i_P^G$ when
$P=MN$ and $Q=LU$ are two parabolic subgroups of $G$. This lemma, based on the structure of the orbits of
$Q$ in the complete variety $P\backslash G$, gives a filtration of $r_Q^G \circ i_P^G$ whose associated
graded module is a direct sum of functors of the form $i_{Q'}^L \circ w \circ r_{P'}^M$, where $P'=M'N'$
is a parabolic subgroup of $M$, $w$ denotes both an element of $G$ and the inner automorphism it defines,
and $Q'=L'U'$ is a parabolic subgroup of $L$ conjugate to $P'$ such that $w\cdot L'=M'$. Using this
result, it is possible to understand the structure of induced representations of irreducible
supercuspidal representations of a Levi subgroup. First, we deduce that an induced representation admits
no supercuspidal subquotient and is of finite length.

A cuspidal datum is a pair $(M,\rho)$ where $M$ is a Levi subgroup of $G$ and $\rho$ an irreducible
supercuspidal representation of $M$. Two cuspidal data are associated if they are conjugate under $G$. If
$(M_1,\rho_1)$ and $(M_2,\rho_2)$ are two cuspidal data, and if $P_1$, $P_2$ are parabolic subgroups of
$G$ with Levi factors $M_1$ and $M_2$ respectively, then the induced representations $i_{P_1}^G\rho_1$
and $i_{P_2}^G\rho_2$ admit a common irreducible subquotient if and only if the cuspidal data
$(M_1,\rho_1)$ and $(M_2,\rho_2)$ are associated. The Jordan-Hölder series of $i_{P_1}^G\rho_1$ and
$i_{P_2}^G\rho_2$ are then equivalent.

We continue the study of the functors $i_P^G$, $r_P^G$, by establishing that they preserve certain
classes of representations. It is elementary to show that the functors $i_P^G$ preserve admissibility,
and that the functors $r_P^G$ preserve finite generation. Jacquet's lemma, in a first version applying
only to admissible representations, shows that the functors $r_P^G$ also preserve admissibility. It is a
simple consequence of the geometric lemma that the functors $i_P^G$ preserve representations of finite
length. Howe's theorem establishes that a representation is of finite length if and only if it is
finitely generated and admissible. It follows that the functors $r_P^G$ preserve representations of
finite length.

We then continue the study of the induced part of the category $\caM(G)$. For this, we must equip the set
of equivalence classes (under conjugation) of cuspidal data $\Omega(G)$ of $G$ with a variety structure.
This is inherited from the fact that the inertial classes of irreducible supercuspidal representations of
the Levi subgroups of $G$ are equipped with such a structure; we have even seen that they are complex
tori. We must take conjugation into account, and we then see that $\Omega(G)$ is a disjoint union
(infinite) of connected components, which are themselves complex algebraic varieties, more precisely
quotients of tori by the action of a finite group. We denote by $\caB(G)$ the set of connected components
of $\Omega(G)$. To each irreducible smooth representation of $G$, we can associate, by the properties of
induced representations stated above, an element of $\Omega(G)$ called the cuspidal support (or also, by
analogy with the theory of real groups, infinitesimal character) of the representation, depending only on
its isomorphism class. We can now state Bernstein's decomposition theorem. For any element $\frs$ of
$\caB(G)$ (i.e., a connected component of $\Omega(G)$), let $\caM(G)_\frs$ be the full subcategory of
representations of $\caM(G)$ all of whose irreducible subquotients have a cuspidal support in $\frs$.
Then $\caM(G)$ is the product over $\caB(G)$ of the categories $\caM(G)_\frs$. The proof of this theorem
relies on the properties of the induction and restriction functors previously established. As a
corollary, we deduce a Noetherian property of the category $\caM(G)$, which implies that the functors
$i_P^G$ preserve finitely generated representations.

Following the ideas that guided the analysis of the cuspidal part, we would now like to analyze each of
the components $\caM(G)_\frs$, by describing these categories as categories of modules over a unital
algebra, and thereby obtain a description of their centers. To do this, as for the cuspidal components,
we must exhibit a small progenerator of these categories. A natural approach is to construct it via
parabolic induction from a small progenerator of a cuspidal component $D$ of a Levi subgroup $M$ such
that the pair $(M,D)$ determines the connected component $\frs$. We have seen that the induction functors
preserve finite generation, which shows that the object thus obtained is small. A technical result shows
that this object is independent of the choices made for $M$ and $D$, and it follows quite easily that it
is indeed a generator of the category. It remains to show that it is also a projective object. Since
$i_P^G$ admits a right adjoint, a general result tells us that it suffices for this adjoint to be an
exact functor. We stated above that the right adjoint of $i_P^G$ is $r_{\bar P}^G$, and therefore is
exact. The proof of this result, Bernstein's second adjunction theorem, is delicate. It relies on two
intermediate results, which will be used again later. The first is the strong version (without
admissibility hypothesis) of Jacquet's lemma, which asserts that if $K$ is a compact open subgroup of $G$
admitting an Iwahori decomposition $K=K_{\overline{N}}K_MK_N$ with respect to the standard parabolic
subgroups, then for any smooth representation $(\pi,V)$ of $G$, the natural map from $V^K$ to its Jacquet
module $r_P^G(V)$ realizes a surjection onto $r_P^G(V)^{K_M}$. Moreover, we obtain a canonical section of
this surjection. From this, we deduce the second result, which is the existence of a non-degenerate
$M$-equivariant duality \footnote{Duality of Casselman, who established it in the case where $\pi$ is
admissible. It is necessary to free oneself from this hypothesis to obtain the second adjunction
theorem.} between the spaces $r_P^G (V)$ and $r_{\bar P}^G(\widetilde V)$, for any smooth representation
$(\pi,V)$ of $G$, and any parabolic subgroup $P=MN$ of $G$. Casselman's duality is equivalent to the
second adjunction theorem, and we therefore have explicit progenerators of the categories $\caM(G)_\frs$.
From there, we obtain, as for the cuspidal components, a description of the center of the category
$\caM(G)_\frs$ as an algebra of polynomial functions on the variety $\frs$, the fact that the categories
$\caM(G)_\frs$ are indecomposable, and the description of the center of the entire category $\caM(G)$ as
an algebra of polynomial functions on $\Omega(G)$.

Throughout this chapter, $G$ denotes a $p$-adic reductive group in the sense of the previous chapter,
whose notation we follow.

\section{The functors $i_P^G$ and $r_P^G$}

\subsection{Adjunction properties}\label{frobPG}
Let $P=MN$ be a parabolic subgroup of $G$. Let us specialize the results of paragraphs \ref{FoncJac1} and
\ref{Jacmod} to the present situation. Since $M$ normalizes $N$, we have according to Remark
\ref{FoncJac1}, a functor
\[
j_N \colon \caM(G) \rightarrow \caM(M),\quad (\pi,V) \mapsto (\pi_N,V_N).
\]

More exactly, this functor is the composition of the forgetful functor
$\caF^G_P \colon  \caM(G)\rightarrow \caM(P)$ and the functor $j_N \colon  \caM(P) \rightarrow \caM(M)$.
The fact that we have two different mathematical objects denoted in the same way should not cause any
confusion here.

In the opposite direction, if $(\tau,E)$ is a smooth representation of $M$, we extend it trivially to
$P$, i.e., $\tau(mn)=\tau(m)$ for all $m\in M$ and all $n\in N$ (if $\phi \colon  P\rightarrow M$ is
given by the composition of the canonical projection $P\rightarrow P/N$ and the isomorphism
$P/N\simeq M$, this representation of $P$ is $\caF^M_P(E)$ where $\caF_P^M$ is the forgetful functor
associated to $\phi$, but note that according to \ref{Jacmod}, it is also $\check{}\, \caF^M_P(E)$). By
then inducing from $P$ to $G$, we obtain a representation of $G$. We will simplify the notation by simply
writing $\Ind_P^G (\tau,E)$ instead of $\Ind_P^G(\caF_P^M(\tau,E))$. This defines a functor from
$\caM(M)$ to $\caM(G)$. We give the first properties of these functors.

\begin{thm}
The functor $j_N \colon  (\pi,V) \mapsto (\pi_N,V_N)$ from $\caM(G)$ to $\caM(M)$ is the left adjoint of
the induction functor $\Ind_P^G$. For any smooth representation $(\tau,E)$ of $M$ and any smooth
representation $(\pi,V)$ of $G$, we have a natural isomorphism
\[
\Hom_M(  (\pi_N,V_N), (\tau,E))\simeq  \Hom_G((\pi,V),\Ind_P^G (\tau,E)).
\]
The functors $\Ind_P^G$ and $j_N$ are exact.
\end{thm}

\begin{proof}
As we noted above, $j_N \colon \caM(G) \rightarrow \caM(M)$ is the composition of the forgetful functor
$\caF^G_P \colon  \caM(G)\rightarrow \caM(P)$ and the functor $j_N \colon \caM(P) \rightarrow \caM(M)$.
We saw in \ref{Jacmod} that the latter is the left adjoint of $\check{}\, \caF^M_P \simeq \caF^M_P$.
Since $\caF_P^G$ is the left adjoint of $\Ind_P^G$, we deduce the first assertion. In summary:
\begin{align*}
  j_N \colon & \caM(G)  \stackrel{\caF_P^G}{\longrightarrow } \caM(P)  \stackrel{j_N}{\longrightarrow } \caM(M)\\
\Ind_P^G \colon & \caM(M)\stackrel{\caF_P^M=\check{}\,  \caF_P^M}{\longrightarrow } \caM(P) \stackrel{\Ind_P^G  }{\longrightarrow }\caM(G)
\end{align*}
are adjoints by composition of adjoint functors.

The exactness of $j_N \colon \caM(G) \rightarrow \caM(M)$ comes from the exactness of the forgetful
functor $\caF_P^G$ and the exactness of $j_N \colon \caM(P) \rightarrow  \caM(M)$ (Corollary
\ref{FoncJac1}). The exactness of $\Ind_P^G \colon \caM(M)\rightarrow \caM(G)$ comes from the exactness
of $\Ind_P^G \colon  \caM(P)\rightarrow \caM(G)$ (Proposition \ref{Ind}), and since the forgetful functor
$\caF_P^ M$ is isomorphic to the pseudo-forgetful functor $\check{}\, \caF_P^ M$, it admits a right and
left adjoint, and is therefore exact.
\end{proof}

From the properties of adjoint functors (\ref{adjfonct}), we deduce the:

\begin{cor}
The functor $j_N$ preserves colimits and the functor $\Ind_P^G$ preserves limits.
\end{cor}

\begin{rmq}
It turns out that $\Ind_P^G$ also admits a right adjoint. Indeed, since the variety $G/P$ is compact,
according to Lemma \ref{IndAdm},
\[
\Ind_P^G=\Ind_P^G \circ \caF^M_P = \ind_P^G  \circ \caF^M_P.
\]
But the forgetful functor $\caF^M_P$ admits a right adjoint, the functor $I^M_P$. On the other hand,
according to Theorem \ref{indcomp},
\[
\ind_P^G= P_P^G(\, \bullet \, \otimes \delta_{P\backslash G}^{-1}),
\]
and $P_P^G$ (resp. $\bullet \otimes \delta_{P\backslash G}^{-1}$) admits the pseudo-forgetful functor
$\check{}\, \caF_P^G$ (resp. $\bullet \otimes \delta_{P\backslash G}$, cf. Remark \ref{OHom}, 3) as a
right adjoint. Consequently, $\Ind_P^G$ preserves colimits.
\end{rmq}

\subsection{Normalization}\label{norm}
Let $P=MN$ be a parabolic subgroup of $G$. We normalize the parabolic induction functor as in Section
\ref{herm}, taking into account Lemma \ref{calcfonctmod}, to obtain a functor that commutes with the
functor $\pi\mapsto \tilde \pi$ of Section \ref{contrag}. We therefore set, for any smooth representation
$(\tau,E)$ of $M$,
\begin{equation}\label{eq:indmorm}
 i_P^G(\tau,E)= \Ind_P^G
  (\tau\otimes \delta_{P}^{-1/2},E)
\end{equation}
\index[not]{iPG@$i_P^G$} Similarly, we normalize the parabolic restriction functor $j_N$ to obtain the
left adjoint of $i_P^G$. For any smooth representation $(\pi,V)$ of $G$, we set
\begin{equation}
r_P^G(\pi,V)= (\pi_N\otimes \delta_{P}^{1/2},V_N ).
\end{equation}

\begin{rmq}
It seems that our conventions differ from those generally adopted, which replace $\delta_P$ by its
inverse in the formulas above. This comes from our initial choice of the modular function in
\ref{foncmod}. We note that the calculations of Lemma \ref{calcfonctmod} also reveal this problem, and
that while non-standard, our conventions are consistent.
\end{rmq}

\begin{prop}
\index[not]{rPG@$r_P^G$} The functors
\[
i_P^G \colon  \caM(M) \rightarrow \caM(G), \quad r_P^G \colon   \caM(G) \rightarrow \caM(M)
\]
are exact. The first is the right adjoint of the second.
\end{prop}

\begin{proof}
This follows without difficulty from Theorem \ref{frobPG}, since the effects of the normalizations cancel
each other out (cf. Remark \ref{OHom}, 3).
\end{proof}

The functors $i_P^G$ and $r_P^G$ will be respectively called parabolic induction and restriction
functors. \index[ter]{induction!parabolic} \index[ter]{restriction!parabolic}

\begin{cor}
The functor $i_P^G$ preserves limits. By adapting the remark of the previous paragraph, we easily see
that it also preserves colimits. The functor $r_P^G$ preserves colimits.
\end{cor}

\begin{rmqs}
Since $\delta_{P}$ takes values in $\bbR^\times_+$, its restriction to any compact subgroup of $P$ is
trivial. Since the unipotent radical $N$ is the union of its compact open subgroups, the restriction of
$\delta_{P}$ to $N$ is trivial. The representation $i_P^G (\tau,E)$ therefore has as its space the set of
functions $f:G\rightarrow E$ such that for all $m\in M$, $n\in N$, $g\in G$,
\begin{equation}\label{inductionnormalisee}
 f(mng)=\delta_{P}^{-1/2}(m)\tau(m)\cdot f(g).
\end{equation}
\end{rmqs}

\subsection{The functor $r_P^G$ preserves finite generation} \label{finiIJ}

\begin{prop}
Let $P=MN$ be a parabolic subgroup of $G$. Let $(\pi,V)$ be a finitely generated smooth representation of
$G$. Then $(\pi_N,V_N)$ is finitely generated. The same is true for $r_P^G(\pi,V)$.
\end{prop}

\begin{proof}
Let $\{v_1,\ldots ,v_l\}$ be a set of $l$ vectors generating $V$ as a representation of $G$. Let us
choose a compact open subgroup $K$ of $G$ fixing all the $v_i$, for $i=1,\ldots,l$. Since $P\backslash G$
is compact, $P\backslash G/K$ is finite. Let $\{g_1,\ldots, g_n\}$ be a system of representatives of
these double cosets. Then $V$ is generated as a $P$-module by the vectors of the form
$\pi(g_j)\cdot v_i$, $i=1,\ldots,l$, $j=1,\ldots,n$. Since $N$ acts trivially on $V_N$, $V_N$ is
generated as an $M$-module by the images of the $\pi(g_j)\cdot v_i$ in $V_N$. The normalization does not
affect this, and we deduce the last assertion.
\end{proof}

\subsection{Transitivity of $r_P^G$ and $i_P^G$. }\label{associa}
Let $P=MN$ and $Q=LU$ be two semi-standard parabolic subgroups of $G$ with $P \subset Q$. Then
$L\cap P=M(L\cap N)$ is a parabolic subgroup of $L$, $N=(N\cap L)U$ and $(L\cap P)U=MN=P$.

\begin{lemme}
Let $P=MN \subset Q=LU$ be two semi-standard parabolic subgroups of $G$. Then
$i_P^G=i_Q^G\circ  i_{P\cap L}^L$ and $r_P^G=r_{L\cap P}^L \circ r_Q^G$.
\end{lemme}
\begin{proof}
Note first that we are committing an abuse of language, what we write as equalities of functors are in
fact only natural isomorphisms. We have
\[
r_{L\cap P}^L \circ r_Q^G= (\, \cdot \,  \otimes \delta_{P\cap
L}^{1/2})\circ j_{N\cap L}\circ \caF_{P\cap L}^L \circ
(\, \cdot \,  \otimes \delta_{Q}^{1/2})\circ j_{U}\circ \caF_Q^G.
\]
Now, we easily see that
\begin{align*}
 &\caF_{P\cap L}^L \circ (\, \cdot \,  \otimes \delta_{Q}^{1/2})=(\,
 \cdot \,  \otimes \delta_{Q}^{1/2})\circ
\caF_{P\cap L}^L,\\
 &\caF_{P\cap L}^L \circ j_U=j_U
\circ \caF^Q_{(L\cap P)U},\\
&j_{N\cap L}\circ  (\, \cdot \,  \otimes \delta_{Q}^{1/2})= (\, \cdot
\,  \otimes
\delta_{Q}^{1/2}) \circ  j_{N\cap L}.
\end{align*}
This last equality comes from the fact that $\delta_Q$ is trivial on $N\cap L$. We therefore obtain
\begin{align*}
 r_{L\cap P}^L \circ r_Q^G&= (\, \cdot \,  \otimes \delta_{P\cap
  L}^{1/2})\circ (\, \cdot \,  \otimes \delta_{Q}^{1/2})\circ  j_{N\cap L}\circ j_U \circ \caF_{(P\cap L)U}^Q
\circ \caF_Q^G\\
&= (\, \cdot \,  \otimes \delta_{P\cap
  L}^{1/2}\delta_{Q}^{1/2})\circ  j_{N}\circ \caF_P^G\\
&= (\, \cdot \,   \otimes \delta_{P}^{1/2})\circ  j_{N}\circ
\caF_P^G\\
&=r_P^G
\end{align*}
We used the transitivity of the functors $j$ (Lemma \ref{FoncJac1}) and of the forgetful functors and the
equalities
\[
\delta_{P}= \delta_{N},\quad \delta_{Q}=\delta_U
\]
from Section \ref{calcfonctmod}. The expression of the latter in terms of roots shows the equality
$\delta_{N\cap L}\delta_U=\delta_N$ on $A_M$, and it extends to $M$ by the same argument as the one
already given in \ref{calcfonctmod}, namely the Cartan decomposition of $M$. We have therefore shown the
second assertion. The first follows by adjunction of the functors $r$ and $i$ and uniqueness of the
adjoint (cf. \ref{adjfonct}).
\end{proof}

\section{Supercuspidal representations and admissibility}

\subsection{Supercuspidal representations}\label{supercusp}

\begin{defi}
A smooth representation $(\pi,V)$ of $G$ is said to be supercuspidal \index[ter]{supercuspidal} if for
any proper parabolic subgroup $P=MN$ of $G$, $r_P^G(V)$ is zero.
\end{defi}

\begin{rmqs}
1. A smooth representation $(\pi,V)$ of $G$ is supercuspidal if for any proper standard parabolic
subgroup $P=MN$ of $G$, $r_P^G(V)$ is zero.

2. Since $r_P^G$ is exact, any subquotient of a supercuspidal representation is supercuspidal.
\end{rmqs}

The subgroup ${}^0G$ was defined in \ref{carnonram}.

\begin{thm}
Let $(\pi,V)$ be a smooth representation of $G$. The following conditions are equivalent:

$a)$ $(\pi,V)$ is supercuspidal,

$b)$ The matrix coefficients of $(\pi,V)$ are compactly supported modulo the center $Z(G)$,

$c)$ the restriction of $(\pi,V)$ to ${}^0G$ is compact.
\end{thm}

\begin{proof}
[$b) \Rightarrow c)$] Let $v \in V$, $\lambda \in \widetilde V$. We want to show that the restriction of
the matrix coefficient $\phi_{v,\lambda}$ to ${}^0G$ is compactly supported. Let $C$ denote the support
of $\phi_{v,\lambda}$ in $G$ and let ${}^0C= C\cap {}^0G$, $\overline{{}^0C}$ be the image of ${}^0C$ in
${}^0G/({}^0G\cap Z(G))$ and $\bar C$ the image of $C$ in $G/Z(G)$. The natural inclusion
${}^0G/({}^0G\cap Z(G)) \hookrightarrow G/Z(G)$ allows us to view ${}^0G/({}^0G\cap Z(G))$ as a closed
subgroup of $G/Z(G)$. Since by hypothesis $\bar C$ is compact,
\[
\overline{{}^0C}\subset \bar C \cap ({}^0G/({}^0G\cap Z(G)))
\]
is compact.

Theorem 2.6 of \cite{Hoch} asserts that if $H$ is a compact subgroup of a topological group $G$, and if
$G/H$ is compact, then $G$ is compact. Let $\bar K$ be a compact open subgroup of
${}^0G/({}^0G\cap Z(G))$, and $K$ its inverse image in ${}^0G$. By virtue of this result, and the fact
that ${}^0G\cap Z(G)$ is compact (\ref{Lambda}), $K$ is a compact open subgroup. We can cover
$\overline{{}^0C}$ by a finite number of translates of $\bar K$, and we can therefore cover ${}^0C$ by a
finite number of translates of $K$. This shows that ${}^0C$ is compact.

[$c) \Rightarrow a)$] Let $P=MN$ be a proper parabolic subgroup of $G$, with split component $A$. In
particular, $A_G$ is strictly included in $A$. Without loss of generality, we can assume that $P$ is
standard. Let $v \in V$. We want to show that $v\in V(N)$. Let $K$ be a compact open subgroup of $G$
admitting an Iwahori decomposition with respect to $P$,
\[
K=K_{\overline{N}} K_M K_N
\]
(cf. \ref{KN}) and such that $v\in V^K$.

Let $t \in A^{++}\cap {}^0G$ (such an element exists according to \ref{pasmieux}). For any root
$\alpha \in \Delta(P)$, $|\alpha(t)|_\bbF<1$, and therefore, for any compact subset $C$ of
$\bbR^\times_+$, there exists $m_0 \in \bbN$ such that if $|m|\geq m_0$, $|\alpha(t^m)|_\bbF$ is not in
$C$. It follows that for any compact subset of $A \cap {}^0G$, there exists $m_0 \in \bbN$ such that if
$|m|\geq m_0$, $t^m$ is outside this compact subset. Since the restriction of $(\pi,V)$ to ${}^0G$ is
compact by hypothesis, we have
\[
\pi(e_K)\pi(t^m)\cdot v=0
\]
for $m$ large enough. Therefore, for $m$ large enough,
\[
\pi(e_{t^{-m}Kt^m})\cdot
v=\pi(e_{t^{-m}K_{N}t^m})\pi(e_{t^{-m}K_{M}t^m})\pi(e_{t^{-m}K_{\overline{N}}t^m})\cdot v=0.
\]
Since $t\in A$, it is central in $M$ and $t^{-m}K_Mt^m=K_M \subset K$. Similarly,
$t^{-m}K_{\overline{N}}t^m \subset K_{\overline{N}}\subset K$, hence
\[
\pi(e_{t^{-m}K_{N}t^m}) \cdot v=0.
\]
Since $t^{-m}K_{N}t^m$ is a compact subgroup of $N$, we have $v \in V(N)$ according to Proposition
\ref{FoncJac1}.

[$a) \Rightarrow b)$] We now show that if $(\pi,V)$ is not compact modulo the center, then there exists a
parabolic $P=MN$ such that $V_N \neq 0$. Up to conjugation, and by transitivity of the restriction
functors, we may assume that such a parabolic is standard and maximal. Let $\phi_{v,\lambda}$ be a matrix
coefficient whose support is not compact modulo the center. The Cartan decomposition
$G=K_0  F_\emptyset^{+} C_\emptyset^+ K_0$ (cf. \ref{Cartan}, notice that $F_\emptyset$ normalizes
$C_\emptyset^+$) shows that there exists an $f \in F_\emptyset^{+}$ and a sequence of distinct elements
$t_n \in C_\emptyset^+$, $n\in \bbN$ such that
\[
\supp\phi_{v,\lambda}\cap K_0 ft_n K_0 \neq \emptyset,
\]
and such that the set of $t_n$ is not contained in any compact subset modulo the center of $G$. There
then exists a root $\alpha \in \Delta_\emptyset$ such that $\{ |\alpha(t_n)|_\bbF\}_{n\in\bbN}$ is not
contained in any compact subset of $\bbR^\times_+$, and therefore, by extracting a subsequence if
necessary, we may assume that
\begin{equation}\label{LIM0}
\lim_{n \rightarrow \infty}|\alpha(t_n)|_\bbF=0.
\end{equation}

Consider the maximal standard parabolic subgroup $P=MN=P_\alpha=M_\alpha N_\alpha$, associated to the
subset $\Delta_\emptyset \setminus \{\alpha\}$ of $\Delta_\emptyset$. Let us choose a compact open
subgroup $K$ of $G$, open and normal in $K_0$, admitting an Iwahori decomposition with respect to $P$,
$K=K_{\overline{N}}K_M K_N$, and such that $v \in V^K$, $\lambda \in \widetilde V^K$. Since $[K_0:K]$ is
finite, we can find $k_1,k_2$ in $K_0$ such that
\[
\supp\phi_{v,\lambda}\cap k_1K ft_nKk_2 \neq \emptyset
\]
for infinitely many $n$. For each of these $n$, let us choose $k_n'$ and $k''_n$ in $K$ such that
\[
\phi_{v,\lambda}(k_1k'_n ft_nk''_nk_2)\neq 0.
\]
We have:
\begin{align*}
 \phi_{v,\lambda}(k_1k'_n
 ft_nk''_nk_2)&=\lambda(\pi(k_1k'_nk_1^{-1})\pi(k_1ft_nk_2)\pi(k_2^{-1}k_n''k_2)\cdot
 v)\\
&= \lambda(\pi(k_1ft_nk_2)\cdot v)=\lambda_1(\pi(t_n)\cdot v_1)\\
&=(\tilde\pi(e_{K_N})\cdot\lambda_1)(\pi(t_n)\cdot v_1) \\
&=\lambda_1 ( \pi(e_{K_N}) \pi(t_n) \cdot v_1)\\
&= \lambda_1 (\pi(t_n)\pi(e_{t^{-1}_nK_Nt_n})\cdot v_1)\neq 0.
\end{align*}
The second equality comes from the fact that $K$ is normal in $K_0$ and that $v$ and $\lambda$ are fixed
by $K$. We then set $v_1= \pi(fk_2)\cdot v$ and $\lambda_1=\tilde \pi(k_1^{-1})\cdot \lambda$ (note that
$t_n$ and $f$ commute). Since $\lambda_1$ is fixed by $K$, and therefore by $K_N$, the rest of the
calculation follows easily. We obtain that $\pi(e_{t^{-1}_nK_N t_n})\cdot v_1\neq 0$. But (\ref{LIM0})
implies that $N=\bigcup_{n\in \bbN}  t^{-1}_nK_Nt_n$ and therefore $v_1 \notin V(N)$ by the
characterization of the Jacquet kernel $V(N)=\bigcup_{n\in \bbN} \ker \pi(e_{t^{-1}_nK_N t_n})$
established in Proposition \ref{FoncJac1}. This shows that $V_N$ is non-zero.
\end{proof}

\smallskip

\begin{rmqs}
1. If the center of $G$ is compact, a representation is compact if and only if it is supercuspidal.

2. Another equivalent condition for a representation $(\pi,V)$ to be supercuspidal is that the functions
$f_{K,v}$ defined in \ref{comprep} are compactly supported modulo $Z(G)$. This is easily seen by adapting
the proof of Theorem \ref{comprep}.
\end{rmqs}

\begin{cor}
Let $(\pi,V)$ be an irreducible smooth representation of $G$. Then there exists a parabolic subgroup
$P=MN$ of $G$ (which we can assume standard) such that $r_P^G(\pi,V)$ is a supercuspidal representation
of $M$. By Frobenius reciprocity, there exists an irreducible supercuspidal representation $(\tau,E)$ of
$M$ such that $(\pi,V)$ is a subrepresentation of $i_P^G(\tau,E)$.
\end{cor}
\begin{proof}
Let $P=MN$ be a standard parabolic subgroup of $G$, minimal with respect to the property
$r_P^G(V) \neq 0$ (this exists because for $P=G$, $r_P^G(V)=V \neq 0$). Recall (cf. \ref{semst}) that the
standard parabolic subgroups of $M$ are the traces of the standard parabolics $P'=M'N'$ of $G$ such that
$P'\subset P$. According to Lemma \ref{associa}, we have $r^M_{M\cap P'}\circ r_P^G(V)=r_{P'}^G(V)$. By
minimality of $P$, we have $r_{P'}^G (V)=0$ if $P'\neq P$, and therefore $r_P^G(\pi,V)$ is supercuspidal.
Moreover, Proposition \ref{finiIJ} asserts that $r_P^G(\pi,V)$ is finitely generated, hence admits an
irreducible quotient. Let $(\tau,E)$ be an irreducible quotient of $r_P^G(\pi,V) $. Since the functors
$r$ are exact, $(\tau,E)$ is supercuspidal and we have by Frobenius reciprocity:
\[
0\neq \Hom_M(r_P^G(\pi,V),(\tau,E))=\Hom_G((\pi,V), i_P^G  (\tau,E)).
\]
\end{proof}

We denote by $\caM(G)_{sc}$ \index[not]{M(G)sc@$\caM(G)_{sc}$} the full subcategory of $\caM(G)$ whose
objects are the supercuspidal representations of $G$. Since $\caM(G)_{sc}$ is stable under passing to
subquotients, it is an abelian category. We denote by $\mathbf{Irr}(G)_{sc}$
\index[not]{Irr(G)sc@$\mathbf{Irr}(G)_{sc}$} the set of isomorphism classes of irreducible supercuspidal
representations of $G$.

\subsection{Admissibility of irreducible representations}\label{admirr}

The relationship established between supercuspidal representations and compact representations and the
properties of parabolic induction allow us to show the following result.

\begin{thm}
Let $(\pi,V)$ be an irreducible smooth representation of $G$. Then $(\pi,V)$ is admissible.
\end{thm}

\begin{proof}
According to Proposition \ref{comprep}, any finitely generated compact representation is admissible.
Recall, according to Theorem \ref{supercusp}, that an irreducible supercuspidal representation has a
restriction to ${}^0G$ which is compact. Since the quotient $G/{}^0GZ(G)$ is finite (Proposition
\ref{Lambda}), and $Z(G)$ acts by scalars on an irreducible representation according to Schur's lemma
\ref{schur}, it is also finitely generated. This restriction is therefore admissible. Since any compact
subgroup of $G$ is contained in ${}^0G$, the admissibility of a representation depends only on its
restriction to ${}^0G$. It is therefore clear that an irreducible supercuspidal representation is
admissible. In the general case, by the corollary of the previous section, we can embed $(\pi,V)$ into a
representation of the form $i_P^G (\tau,E)$, with $(\tau,E)$ irreducible supercuspidal, and therefore
admissible. According to Lemma \ref{IndAdm}, $i_P^G (\tau,E)$ is admissible, and therefore $(\pi,V)$ is
too.
\end{proof}

\begin{cor}
Let $(\pi,V)$ be a smooth representation of $G$. Then $\pi$ is irreducible (resp. irreducible
supercuspidal) if and only if $\tilde \pi$ is irreducible (resp. irreducible supercuspidal).
\end{cor}

\begin{proof}
This follows from the theorem and Proposition \ref{Admis}.
\end{proof}

\subsection{Uniform admissibility theorem} \label{unifadm}
\index[ter]{uniform admissibility} The theorem of the previous paragraph shows that if $K$ is a compact
open subgroup of $G$, for any irreducible smooth representation $(\pi,V)$ of $G$, $\dim V^K$ is finite.
But a priori, this dimension can become arbitrarily large when $(\pi,V)$ varies. In fact, this is not the
case, as shown by the following:
\begin{thm}
Let $K$ be a compact open subgroup of $G$. There exists a positive constant $c=c(G,K)$ such that for any
irreducible smooth representation $(\pi,V)$ of $G$, $\dim V^K \leq c$.
\end{thm}

\begin{proof}
We can reformulate this result using Theorem \ref{piK}. It is equivalent to the fact that all simple
modules of the algebra $\caH(G,K)$ are of dimension bounded by a positive constant $c=c(G,K)$. The main
ingredient in the proof is the decomposition
\[
\caH(G,K)=  \sum_{j=1}^q\caH_0  \caA D_j \caH_0
\]
of Theorem \ref{HCDH}. Let $(\rho,W)$ be a simple $\caH(G,K)$-module. The admissibility result
\ref{admirr}, reinterpreted as above, tells us that such a module is finite-dimensional, say $\dim W=k$.
By Burnside's theorem (\cite{Lang}, Corollary XVII.3.3), $\rho \colon \caH(G,K)\rightarrow \End(W)$ is
surjective. Since $\rho(\caA)\subset \End(W)$ is commutative, since elements in $C_{A_G}$ are central and
act by a scalar on the simple module $W$ (via Schur's lemma) and since the cone
$S:=C_\emptyset^+/C_{A_G}$ is generated by the images of the $t_\alpha$, $\alpha\in\Delta_\emptyset$, and
the elements of the finite set $R$ (see \ref{SSM}), we see that $\rho(\caA)$ is generated by $l=d+|R|$
generators (representatives in $C_\emptyset^+$ of these elements) and the identity, and according to
Lemma \ref{dimborn}, we have $\dim \rho(\caA)\leq k^{2-2^{1-l}}$. Set $h=\dim \caH_0$. We then have
\[
k^2=\dim \End(W)=\dim \rho(\caH(G,K))\leq h^2 q \dim \rho(\caA)\leq
h^2 q k^{2-2^{1-l}},
\]
hence $k\leq (h^2q)^{2^{l-1}}$. We can therefore take $c=c(G,K)=(h^2 q)^{2^{l-1}}$.
\end{proof}

\begin{rmq}
An analogous result is true for irreducible smooth representations of the group ${}^0G$. We deduce it
from the result for $G$ as follows. First, recall that any compact subgroup $K$ of $G$ is included in
${}^0G$. Let $(\sigma,E)$ be an irreducible smooth representation of ${}^0G$. Let us form
$(\pi,V)=\ind_{{}^0G}^{\; \, G} (\sigma,E)$.

Frobenius reciprocity for compact induction then yields (cf. \ref{indres}):
\[
\Hom_G((\pi,V), (\tau,W))\simeq  \Hom_{{}^0G}((\sigma,E),
\res_{{}^0G}^{\; \,  G}(\tau,W)).
\]
Since $V=\ind_{{}^0G}^{\; \,  G} (E)=P_{{}^0G}^{\; \,  G} (E)= \caH(G)\otimes_{\caH({}^0G)}E$, it is
clear that $V$ is a finitely generated representation. Indeed, if $v \in E$, $v\neq 0$, and if $K$ is a
compact open subgroup of ${}^0G$ fixing $v$, then
\begin{align*}
 \caH(G) \otimes_{\caH({}^0G) } E&= \caH(G)\otimes_{\caH({}^0G)}
\caH({}^0G)\cdot v=  \caH(G)\otimes_{\caH({}^0G)}v\\&=
\caH(G)\otimes_{\caH({}^0G)}e_K\cdot  v=\caH(G)(e_K\otimes v).
\end{align*}
The representation $(\pi,V)$ therefore admits an irreducible quotient. Substituting this irreducible
quotient for $(\tau,W)$ in the Frobenius reciprocity formula above, we see that its restriction to
${}^0G$ contains $(\sigma,E)$, and therefore $\dim E^K \leq \dim W^K\leq c(G,K)$.
\end{rmq}

\subsection{Consequences of uniform admissibility} \label{finitude}
Let $K$ be a compact open subgroup of $G$. Given an irreducible supercuspidal representation $(\rho,W)$
of $G$ and a vector $v\in W$, the function
\[
g\mapsto f_{K,v}(g)=\rho(e_K)\rho(g^{-1})\cdot v
\]
has its support compact modulo $Z(G)$ in $G$. By the uniform admissibility theorem, we will establish the
following result.

\begin{prop}
Given a compact open subgroup $K$ of $G$, there exists an open subset $\Omega$ in $G$, compact modulo the
center, such that for any irreducible supercuspidal representation $(\rho,W)$ of $G$, and for any
$v\in W^K$, the support of $g\mapsto \rho(a_{g,K})\cdot v$ is included in $\Omega$.
\end{prop}

\begin{proof}
Since $v\in W^K$, we have
$\rho(a_{g,K})\cdot v =\rho(e_K)\rho(g)\rho(e_K)\cdot v= \rho(e_K)\rho(g)\cdot v=f_{K,v}(g^{-1})$. If
$K'\subset K$ is another compact open subgroup of $G$, since $e_Ke_{K'}=e_{K'}e_K=e_K$, we see that
$\supp f_{K,v}\subset (\supp f_{K',v})K$. Therefore it suffices to prove the result for certain
well-chosen compact open subgroups $K$. We therefore assume that $K$ is normal in $K_0$ and admits an
Iwahori decomposition along the standard parabolic subgroups. The uniform admissibility theorem shows
that there exists a constant $N\in \bbN$ such that $\dim W^K\leq N$. Let $z \in C_\emptyset^+$ such that
its class modulo the equivalence relation (\ref{eqC}) is non-trivial. Then the sequence of successive
powers $z,z^2,z^3,\ldots$ is not contained in any compact subset modulo $Z(G)$. Indeed, according to the
Cartan decomposition \ref{Cartan}, any compact subset of $G$ is contained in a finite disjoint union of
subsets of the form $K_0aK_0$ where the $a \in \widetilde C_\emptyset$. Therefore any compact subset
modulo $Z(G)$ is contained in a finite union of subsets of the form $K_0aK_0 Z(G)$. Since the abelian
group $Z(G)$ admits $Z(G)_{an}A_G$ as a finite-index subgroup, where $Z(G)_{an}$ is a compact abelian
group contained in $K_0$, we see that any compact subset modulo $Z(G)$ is contained in a finite union of
subsets of the form $K_0aA_GK_0$, $a \in Z(G)\widetilde C_\emptyset$. This is not the case for the
sequence $(z^n)$.

For all $w\in W$, let us show that $\rho(a_{z^N,K})\cdot w=\rho(e_K)\rho(z^N)\rho(e_K)\cdot w=0$. Let us
set $v'=\rho(e_K)\cdot w\in W^K$. If $v'=0$, the assertion is true. Suppose then that $v'\in W^K$ is
non-zero. Since $(\rho,W)$ is supercuspidal, the function
\[
g\mapsto  \rho(e_K) \rho(g)\cdot v' =f_{K,v'}(g^{-1})
\]
is compactly supported modulo $Z(G)$ in $G$ (cf. Remark \ref{supercusp}). There therefore exists an
integer $n_0$ such that $\rho(e_K)\rho(z^{n_0-1})\cdot v'\neq 0$ and $\rho(e_K)\rho(z^{n_0})\cdot v'=0$.
We want to show that $n_0\leq N$. We will show that the family
$\{\rho(e_K)\rho(z^{j-1})\cdot v' \}_{j=1, \ldots, n_0}$ is linearly independent in $W^K$, which is
obviously sufficient. Let $c_1,\ldots c_{n_0}$ be scalars such that
\[
\sum_{j=1}^{n_0}c_j \;  \rho(e_K)\rho(z^{j-1})\cdot v' =0.
\]
Since $\rho(e_K)\rho(z^{j-1})\cdot v'=\rho(a_{z^{j-1},K})\cdot v'$ and according to the second lemma in
\ref{HCDH},
\[
\rho(a_{z^{j'},K})\rho(a_{z^j,K})=\rho(a_{z^{j+j'},K}),
\]
we recursively obtain $c_{1}=0$, $c_{2}=0$, \ldots, $c_{n_0}=0$. The assertion is therefore proved.

According to \ref{SSM}, the class modulo $C_{A_G}$ of every $z\in C_\emptyset^+$ can be written
$r\prod_{\alpha\in\Delta_\emptyset}t_\alpha^{a_\alpha}$, with $r\in R$ and $a_\alpha\in\bbN$. Let $S_N$
be the set of $z \in C_\emptyset^+$ whose class modulo $C_{A_G}$ admits such an expression with
$a_\alpha\leq N$ for all $\alpha$. If $z \in C_\emptyset^+$ is not in $S_N$, there exists $\alpha$ such
that $a_\alpha>N$, and $z'=zt_\alpha^{-N}\in C_\emptyset^+$; by the second lemma of \ref{HCDH},
$a_{z,K}=a_{z',K}a_{t_\alpha^N,K}$, and therefore, by the above applied to $t_\alpha$,
$\rho(a_{z,K})\cdot w=0$ for all $w\in W$.

Now let $K''\subset K\cap\bigcap_{f\in F_\emptyset^+}f^{-1}Kf$ be a compact open subgroup, normal in
$K_0$ and admitting an Iwahori decomposition with respect to the standard parabolic subgroups, and let
$N''$ be the corresponding uniform admissibility constant. By the Cartan decomposition (\ref{Cartan}),
every $g\in G$ can be written $g=k_1fzk_2$, with $k_1,k_2\in K_0$, $f\in F_\emptyset^+$ and
$z\in C_\emptyset^+$. Since $K_0$ normalizes $K$,
$\rho(a_{g,K})\cdot v=\rho(k_1)\rho(a_{fz,K})\rho(k_2)\cdot v$, with $\rho(k_2)\cdot v\in W^K$. Now, for
$v'\in W^K\subset W^{K''}$, since $K''\subset f^{-1}Kf$, we have $e_K*\delta_f*e_{K''}=e_K*\delta_f$, and
therefore
\[
\rho(a_{fz,K})\cdot v'=\rho(e_K)\rho(f)\rho(e_{K''})\rho(z)\rho(e_{K''})\cdot v'=\rho(e_K)\rho(f)\rho(a_{z,K''})\cdot v',
\]
which is zero if $z\notin S_{N''}$, by the above applied to $K''$. The support of
$g\mapsto\rho(a_{g,K})\cdot v$ is therefore contained in the finite union of the $K_0fzC_{A_G}K_0$,
$f\in F_\emptyset^+$, $z$ running over a system of representatives of $S_{N''}$ modulo $C_{A_G}$, which
is a compact subset modulo $Z(G)$. This completes the proof of the proposition.
\end{proof}

\section{Decomposition of $\caM(G)$: the cuspidal part}

\subsection{Action of $G$ on $\mathbf{Irr}({}^0G)$}\label{Gacts0G}

We will use the results of paragraph \ref{decMG} concerning compact representations to obtain analogous
results on supercuspidal representations, by exploiting the relationship between supercuspidal
representations and compact representations established in \ref{supercusp}.

Since ${}^0G$ is a normal subgroup of $G$, $G$ acts on $\mathbf{Irr}({}^0G)$ by
$(g,\sigma)\mapsto \sigma^g$ (notation of Example \ref{FOA}).

\begin{lemme}
$(i)$ The orbits of the action of $G$ on $\mathbf{Irr}({}^0G)$ are of finite cardinality.

$(ii)$ Let $(\pi,V)$ be a smooth representation of $G$, $(\sigma,W)$ an irreducible smooth representation
of ${}^0G$, and $V(\sigma)$ the $\sigma$-isotypic component of the restriction of $(\pi,V)$ to ${}^0G$.
We then have, for all $g \in G$,
\[
\pi(g)\cdot V(\sigma)= V( \sigma^g).
\]
\end{lemme}

\begin{proof}
$(i)$ The action factors through the group $G/{}^0GZ(G)$, which is finite. The orbits are therefore of
finite cardinality.

$(ii)$ Recall that by definition $V(\sigma)$ is the image of
\[
\Hom_{{}^0G}((\sigma,W),(\pi,V))\otimes W \rightarrow V
\]
under the canonical map $f \otimes w\mapsto f(w)$. Similarly $V(\sigma^g)$ is the image of
\[
\Hom_{{}^0G}((\sigma^g,W),(\pi,V))\otimes W \rightarrow V.
\]
Let us define, for all $f \in \Hom_{{}^0G}((\sigma,W),(\pi,V))$, $\tilde f(w)=\pi(g)\cdot f(w)$,
$w\in W$. We check that $\tilde f \in \Hom_{{}^0G}((\sigma^g,W),(\pi,V))$:
\begin{align*}
 \tilde f((\sigma^g)(h)\cdot w)&=\tilde f(\sigma(g^{-1}hg)\cdot w)
 =\pi(g)\cdot  f(\sigma(g^{-1}hg)\cdot w)\\
&=  \pi(g)\pi(g^{-1}hg)\cdot f(w)= \pi(hg)\cdot f(w)\\
& = \pi(h)\cdot(\pi(g)\cdot f(w))=\pi(h)\cdot \tilde f(w)
\end{align*}

It is clear that $f \mapsto \tilde f$ realizes an isomorphism
\[
\Hom_{{}^0G}((\sigma,W),(\pi,V)) \simeq \Hom_{{}^0G}((\sigma ^g,W),(\pi,V)).
\]
The desired result follows immediately.
\end{proof}

\subsection{Restriction to ${}^0G$}\label{restra0G}
Recall that we denoted by $\caX(G)$ the set of unramified characters of $G$, i.e., the characters of $G$
trivial on ${}^0G$. We saw in \ref{carnonram} that $\caX(G)$ has the structure of a complex algebraic
torus. \enlargethispage{3\baselineskip} \needspace{165pt}
\begin{prop}
Let $(\pi,V)$ be an irreducible representation of $G$.

$(i)$ The elements of $\mathbf{Irr}({}^0G)$ appearing in the restriction of $(\pi,V)$ to ${}^0G$ form a
single $G$-orbit.

$(ii)$ The restriction of $(\pi,V)$ to ${}^0G$ is semisimple and of finite length.

$(iii)$ Let $(\pi_i,V_i)$, $i=1,2$, be irreducible representations of $G$. The following conditions are
equivalent:

$a)$ $\res_{{}^0G}^{\; \,  G}(\pi_1,V_1)= \res_{{}^0G}^{\; \,  G}(\pi_2,V_2)$.

$b)$ There exists $\chi \in \caX(G)$ such that $\pi_2=\pi_1 \otimes \chi$.

$c)$ $\Hom_{{}^0G}(\res_{{}^0G}^{\; \,  G} \pi_1,\res_{{}^0G}^{\; \,  G}  \pi_2)\neq 0$.
\end{prop}

\begin{proof}
Let us fix a set $\Gamma$ of representatives in $G$ for the cosets of $G/(Z(G){}^0G)$: this set is finite
(Proposition \ref{Lambda}). According to Schur's lemma, $Z(G)$ acts by scalars on $V$ (let $\chi_\pi$
denote the central character of $\pi$), and we deduce that $\res_{{}^0G}^{\; \,  G} (\pi,V)$ is finitely
generated and therefore possesses an irreducible quotient $(\sigma,W)$. Let
$f \in \Hom_{{}^0G}(\res_{{}^0G}^{\; \,  G} (\pi,V),(\sigma,W))$ be non-zero. Then, for all
$\gamma \in \Gamma$,
\[
f \circ \pi(\gamma) \in   \Hom_{{}^0G}(\res_{{}^0G}^{\; \,  G} (\pi,V),(\sigma^{\gamma^{-1}},W)).
\]
Consider
\[
\phi= \bigoplus_{\gamma \in \Gamma}  f \circ \pi(\gamma) \in   \Hom_{{}^0G}(\res_{{}^0G}^{\; \,  G} (\pi,V),
\oplus_{\gamma \in \Gamma}(\sigma^{\gamma^{-1}},W)).
\]
We verify that $\ker \phi$ is stable under the action of $G$. Let $v \in \ker \phi$, i.e.,
$f\circ \pi(\gamma)(v)=0$, for all $\gamma \in \Gamma$, and let $g \in G$, which we write $g=\gamma'yz$,
$\gamma' \in \Gamma$, $y \in {}^0G$, $z\in Z(G)$. For all $\gamma \in \Gamma$, let us set
$\gamma \gamma' =y''z'' \gamma''$, with $\gamma'' \in \Gamma$, $y'' \in {}^0G$, $z''\in Z(G)$. We then
have
\begin{align*}
\phi(\pi(g)\cdot v)&= \sum_{\gamma \in \Gamma}  f\circ \pi(\gamma)(\pi(\gamma')\pi(y)\pi(z)\cdot v)\\
&= \chi_\pi(z)\sum_{\gamma \in \Gamma}  f(\pi(\gamma \gamma'y \gamma'^{-1}\gamma^{-1} )\cdot (\pi(\gamma \gamma')\cdot v))\\
&= \chi_\pi(zz'')  \sum_{\gamma \in \Gamma}  \sigma(\gamma \gamma'y \gamma'^{-1}\gamma^{-1} )
\sigma(y'')\cdot f(\pi(\gamma'')\cdot v)\\
&=0.
\end{align*}
Since $\phi$ is non-zero, we conclude by the irreducibility of $V$ that the kernel of $\phi$ is zero.
Thus $\phi$ is injective and $\res_{{}^0G}^{\; \,  G} (\pi,V)$ is isomorphic to a subrepresentation of
$\oplus_{\gamma \in \Gamma}(\sigma^{\gamma^{-1}},W)$, which is of finite length and semisimple. We
deduce, for instance by using the semisimplicity criterion of Lemma \ref{semisimplicite} $(iii)$, that
the same is true for $\res_{{}^0G}^{\; \,  G} (\pi,V)$. We obtain $(ii)$.

Let us now write:
\[
\res_{{}^0G}^{\; \,  G} (\pi,V) = \bigoplus_{\sigma \in
\mathbf{Irr}({}^0G)} V(\sigma)
\]
where $V(\sigma)$ is the $\sigma$-isotypic component in $\res_{{}^0G}^{\; \,  G} (\pi,V)$. The identity
$\pi(g)\cdot V(\sigma)=V(\sigma^g)$ shows that the set of $\sigma$ such that $V(\sigma)\neq 0$ forms
exactly one $G$-orbit, because $(\pi,V)$ is irreducible. In particular, this set is finite, according to
Lemma \ref{Gacts0G}. This shows $(i)$.

Let us now prove $(iii)$. Since any unramified character of $G$ is by definition trivial on ${}^0G$, $b)$
implies $a)$. Since $a)$ trivially implies $c)$, it remains to show that $c)$ implies $b)$. Suppose that
$c)$ holds. According to $(ii)$ and Schur's lemma, $\Hom_{{}^0G}(V_1,V_2)$ is finite-dimensional and
non-zero. Let us equip this space with the action of $G$ given by
\[
g\cdot f:= \pi_2(g)\circ f \circ \pi_1(g)^{-1}, \quad (f \in
\Hom_{{}^0G}(V_1,V_2)), \; (g \in G).
\]
If $g \in {}^0G$, $g\cdot f=f$ for all $f \in \Hom_{{}^0G}(V_1,V_2)$, so this action of $G$ factors
through $G/{}^0G$, which is commutative because ${}^0G$ contains the derived subgroup of $G$. Since any
commutative family of linear operators on a non-zero finite-dimensional complex vector space admits a
common eigenvector, there exists a non-zero $f\in \Hom_{{}^0G}(V_1,V_2)$ and a character $\chi$ of
$G/{}^0G$ (i.e., an unramified character of $G$) such that $g\cdot f=\chi(g)\, f$, for all $g \in G$. By
construction, $f \in \Hom_{G}(\pi_1 \otimes \chi,\pi_2)$, and therefore by irreducibility,
$\pi_1 \otimes \chi=\pi_2$.
\end{proof}

\subsection{Inertial classes of supercuspidal representations}\label{clinertcusp}
Let $(\pi,V)$ be an irreducible supercuspidal representation of $G$. As above, let us introduce the
decomposition of $\res_{{}^0G}^{\; \,  G} (\pi,V)$ into isotypic components:
\[
V = \bigoplus_{i=1,\ldots ,m} V(\sigma_i),
\]
where the $\sigma_i$ form an orbit for the action of $G$ on $\mathbf{Irr}({}^0G)_c$. Note that this
decomposition actually depends only on the inertial class $[\pi]$. Let $(\rho,W)$ be a smooth
representation of $G$. As a representation of ${}^0G$ we have, according to Theorem \ref{cncomp},
\[
\res_{{}^0G}^{\; \,  G} (\rho,W)=\rho(e_{\sigma_1})\cdot W
\oplus\rho(e_{\sigma_2})\cdot W\oplus \cdots \oplus
\rho(e_{\sigma_m})\cdot W\oplus W'
\]
where $W'$ is the unique ${}^0 G$-invariant complement of
$\rho(e_{\sigma_1})\cdot W \oplus \cdots \oplus \rho(e_{\sigma_m})\cdot W$. Alternatively, $W'$ is
characterized as the unique subrepresentation of ${}^0 G$ none of whose composition factors is isomorphic
to one of the $\sigma_i$.

\begin{prop}
The subspace $W'$ of $W$ is $G$-invariant and for all $i=1,\ldots ,m$, for all $g \in G$,
\[
\rho(g)\circ \rho(e_{\sigma_i})= \rho(e_{g\cdot \sigma_i})\circ \rho(g).
\]
\end{prop}
\begin{proof}
Let $g \in G$. If $(\tau,E)$ is an irreducible subquotient of $(\rho,W')$, then $(\tau^g, E)$ is an
irreducible subquotient of $\rho(g)\cdot W'$, and since the $\sigma_i$ form a $G$-orbit in
$\mathbf{Irr}({}^0G)_c$, $\tau^g$ is not equivalent to any of the $\sigma_i$. It follows that
$\rho(e_{\sigma_i})\rho(g)\cdot W'=0$ for all $i$. By the uniqueness of $W'$, we have
$\rho(g)\cdot W'=W'$.

Let $w\in W$. Then we can write $w$ as
\[
w=w_1+\ldots + w_m+w', \quad \text{ with } w_i \in
\rho(e_{\sigma_i})\cdot W,\;  i=1,\ldots, m, \; w'\in W'.
\]
Let us fix $i$, let $g \in G$, and set $\sigma_j=\sigma_i^g$. Then, since
$$ \rho(e_{\sigma_i})(\rho(g)^{-1}\cdot w)=\rho(g)^{-1} \cdot w_j,$$ we have:
\[
\rho(g)    \rho(e_{\sigma_i})\rho(g)^{-1}\cdot w=w_j= \rho(e_{\sigma_j})\cdot w.
\]
This completes the proof of the proposition.
\end{proof}

\bigskip

Let us set $\rho(e_{[\pi]})= \sum_i  \rho(e_{\sigma_i})\colon W \rightarrow W$. It is an intertwining
operator from $(\rho,W)$ to itself. Since $\rho(e_{\sigma_i}) \rho(e_{\sigma_j})=0$ if $i\neq j$,
$\rho(e_{[\pi]})$ is an idempotent, and we therefore have a decomposition of $W$ into subrepresentations:
\[
W= \rho(e_{[\pi]})\cdot W \oplus (1-\rho(e_{[\pi]}))\cdot W.
\]

If we denote by $\caM(G)_{[\pi]}$ (resp. $(\caM(G)\setminus [\pi])$)
\index[not]{M(G)pi@$\caM(G)_{[\pi]}$} \index[not]{M(G)-pi@$\caM(G)\setminus [\pi]$} the full subcategory
of $\caM(G)$ whose objects are the representations of $G$ all of whose irreducible subquotients are in
$[\pi]$ (resp. none of whose irreducible subquotients are in $[\pi]$), we obtain a decomposition of
$\caM(G)$ into
\begin{equation}\label{dectaunotau}
\caM(G)=\caM(G)_{[\pi]} \times  (\caM(G)\setminus [\pi]).
\end{equation}

\begin{rmqs}
1. The inertial classes of supercuspidal representations are precisely the orbits of the action of
$\caX(G)$ on $\mathbf{Irr}(G)_{sc}$ defined in \ref{finXGpi}. We have seen that they are equipped
(non-canonically) with the structure of a complex algebraic torus. If $(\pi,V)$ is an irreducible
supercuspidal representation of $G$, the variety $[\pi]$ is isomorphic (non-canonically, it depends on
the choice of $\pi$ in its class $[\pi]$) to the complex torus $\caX(G)/\caX(G)(\pi)$, the quotient of
the torus $\caX(G)$ by its finite subgroup $\caX(G)(\pi)$, the stabilizer of $\pi$. The algebra of
functions on the affine variety $[\pi]$ is therefore the algebra of invariants
$\bbC[G/{}^0G]^{\caX(G)(\pi)}$.

2. The restriction of representations from $G$ to ${}^0G$ induces a bijection between the set
$[\mathbf{Irr}(G)_{sc}]$ of inertial classes of irreducible supercuspidal representations of $G$ and the
set of orbits under the action of $G$ on $\mathbf{Irr}({}^0G)_c$.
\end{rmqs}

\subsection{Decomposition of $\caM({}^0G)$}\label{DeCcusp}
Let us return to the decompositions of $\caM({}^0G)$ and $\caM({}^0G)_c$ obtained in paragraph
\ref{decMG}. We saw that given compact representations $\tau_1,\ldots,\tau_r$ of ${}^0G$, we have a
decomposition
\[
\caM({}^0G)= \prod_i \caM({}^0G)_{\tau_i} \times [\caM({}^0G)\setminus\tau_1,\ldots,\tau_r],
\]
and that on the other hand $\caM({}^0G)_{c}$ decomposes into
\[
\caM({}^0G)_{c}=\prod_{\tau \in \mathbf{Irr}({}^0G)_c} \caM({}^0G)_\tau.
\]
We will use the consequences of Theorem \ref{unifadm} to combine both into a finer decomposition of
$\caM({}^0G)$.
\begin{thm}
The category $\caM({}^0G)$ decomposes into
\[
\caM({}^0G)= \prod_{\tau \in \mathbf{Irr}({}^0G)_c }
\caM({}^0G)_{\tau}  \times \caM({}^0G)_{nc}
\]
where $\caM({}^0G)_{nc}$ denotes the full subcategory of smooth representations of ${}^0G$ none of whose
subquotients is a compact representation.
\end{thm}
\begin{proof}
According to Proposition \ref{decMG}, it suffices to verify condition $({\bf KF})$ given in that
proposition. Let $K$ be a compact open subgroup of $G$, and let $(\sigma,E)$ be an irreducible compact
representation of ${}^0G$. As in Remark \ref{unifadm}, let us choose an irreducible representation
$(\pi,V)$ of $G$ whose restriction to ${}^0G$ contains $(\sigma,E)$. It is then supercuspidal. Let
$v \in E^K\subset V^K$ and $\lambda \in \widetilde E^K \subset \widetilde V^K$ by means of which we form
the matrix coefficient
\[
g \mapsto \lambda(\pi(g)\cdot v)=\lambda(\pi(a_{g,K})\cdot v).
\]
According to Proposition \ref{finitude}, this matrix coefficient is supported in a compact subset
$\Omega$ modulo $Z(G)$. The argument of the proof of Theorem \ref{supercusp}, more precisely the
implication $[b)\Rightarrow c)]$, shows that the restriction of this matrix coefficient to ${}^0G$ is
compactly supported, say, on a compact set $\Omega_1$, and according to the choice of $v$ and $\lambda$,
it is equal to $g\mapsto \lambda(\sigma(a_{g,K})\cdot v)$. On the other hand, this matrix coefficient is
$K$-bi-invariant, hence in $\scrD(\Omega_1,K)$. But this algebra is finite-dimensional. By the
orthogonality of the matrix coefficients of irreducible compact representations (Corollary
\ref{orthschur}), we deduce that the number of irreducible compact representations of ${}^0G$ having
non-zero vectors fixed by $K$ is finite.
\end{proof}

\begin{rmq}
We note the following corollary of this proof and of Remark \ref{clinertcusp}: for any compact open
subgroup $K$ of $G$, the number of inertial classes of irreducible supercuspidal representations of $G$
having non-zero vectors fixed by $K$ is finite.
\end{rmq}

\subsection{Decomposition of $\caM(G)$}\label{decMG2}
We will now use the results of the previous paragraph to deduce an analogous decomposition of $\caM(G)$.

\needspace{168pt}
\begin{thm}
The category $\caM(G)$ decomposes into
\[
\caM(G)= \caM(G)_{sc}  \times \caM(G)_{ind} =\prod_{[\tau] \in [\mathbf{Irr}(G)_{sc}]}
\caM(G)_{[\tau]}  \times \caM(G)_{ind}
\]
where $\caM(G)_{sc}$ (resp. $\caM(G)_{ind}$) \index[not]{M(G)ind@$\caM(G)_{ind}$} denotes the full
subcategory of smooth representations of $G$ all of whose subquotients are supercuspidal representations
(resp. none of whose subquotients is supercuspidal).
\end{thm}
\begin{proof}
Let $(\pi,V)$ be a smooth representation of $G$. The results of the previous paragraph give us a
decomposition
\[
\res_{{}^0G}^{\; \, G} (V)= V_c\oplus V_{nc}
\]
By the uniqueness of this decomposition, the two factors on the right-hand side are stable under the
action of $G$. According to Theorem \ref{supercusp}, a representation of $G$ is supercuspidal if its
restriction to ${}^0G$ is compact, and therefore $V_c$ is in $\caM(G)_{sc}$ and $V_{nc}$ is in
$\caM(G)_{ind}$.

It remains to show that in turn $\caM(G)_{sc}$ decomposes as
\[
\caM(G)_{sc} =\prod_{[\tau] \in [\mathbf{Irr}(G)_{sc}]}\caM(G)_{[\tau]}.
\]
Let $(\pi, V) \in \caM(G)_{sc}$. Let us write
\[
\res_{{}^0G}^{\; \, G} (V)= \bigoplus_{(\sigma)  \in
\mathbf{Irr}({}^0G)_c/G}  V^{(\sigma)}
\]
where we sum over the orbits of the action of $G$ on $\mathbf{Irr}({}^0G)_c$ and $V^{(\sigma)}$ denotes
the direct factor of $\res_{{}^0G}^{\; \, G} (V)$ whose irreducible subquotients are in the orbit of
$\sigma$. Then it is clear that $V^{(\sigma)}$ is $G$-stable and supercuspidal, and that two irreducible
subquotients of $V^{(\sigma)}$ are in the same inertial class.
\end{proof}

\subsection{Injectivity and projectivity properties of supercuspidal representations}\label{pisc}

We saw in \ref{conscomp} that any compact representation is projective and injective in $\caM(G)$. If
${}^0G \neq G$, we would like an analogous result for supercuspidal representations in $\caM(G)$.
However, in general, a supercuspidal representation, even irreducible, is neither projective nor
injective in $\caM(G)$. We will nevertheless establish weaker results that will be useful to us later.

\begin{lemme}
Let $(\pi,V)$ be an admissible smooth representation of $G$ and let $(\tau,E)$ be an irreducible
supercuspidal representation of $G$. Suppose that $(\pi,V)$ admits a subquotient isomorphic to
$(\tau,E)$. Then $(\tau,E)$ is realized as a quotient and as a subrepresentation of $(\pi,V)$. In other
words, the spaces
\[
\Hom_G(\pi,\tau) \quad \text {and }\quad \Hom_G(\tau,\pi)
\]
are non-trivial.
\end{lemme}

\begin{proof}
According to the decomposition (\ref{decMG2}), we can assume that $(\pi,V)$ is in $\caM(G)_{[\tau]}$. All
irreducible subquotients of $(\pi,V)$ are then of the form $\tau\otimes \omega$, for some
$\omega \in \caX(G)$. Let $K$ be a compact open subgroup of $G$ such that $E^K\neq \{0\}$. Since
$K \subset {}^0G$, all non-trivial subquotients of $(\pi,V)$ have a non-trivial subspace fixed by $K$.
Since $(\pi,V)$ is admissible, we see that it must be of finite length. We want to show that
$\Hom_G(\tau,\pi) \neq \{0\}$. Let us consider the restriction to ${}^0G$ of these two representations.
Since these restrictions are both compact, hence semisimple, we have
\[
\Hom_{{}^0G}(\tau, \pi)\neq \{0\}.
\]
and since $(\pi,V)$ is of finite length, $\dim \Hom_{{}^0G}(\tau,\pi)< +\infty$.

Let $\Lambda= \Lambda(G)=G/{}^0G$. This group acts on $S=\Hom_{{}^0G}(\tau,\pi)$ by
\[
\bar g\cdot \phi = \pi(g)\circ  \phi \circ  \tau(g^{-1}),
\]
where $g$ is a representative in $G$ of the class $\bar g \in G/{}^0G$. We then have
$\Hom_G(\tau, \pi)=S^{\Lambda}$. Now, the hypothesis that $\tau$ is a composition factor of $\pi$
translates to the fact that there exist $V_1 \subset V_2$, subrepresentations of $V$, such that
$V_2/V_1 \simeq E$. Since $\tau_{|{}^0G}$ is projective in $\caM({}^0G)$, we have an exact sequence
\[
0\rightarrow \Hom_{{}^0G} (E,V_1)\rightarrow \Hom_{{}^0G}
(E,V_2)\rightarrow \Hom_{{}^0G} (E,V_2/V_1)\rightarrow 0.
\]
The group $\Lambda$ acts on each of the terms of this sequence, which becomes a sequence of
$\Lambda$-modules. On the other hand, the inclusion $V_2 \hookrightarrow V$ induces an injection of
$\Lambda$-modules
\[
\Hom_{{}^0G}(E,V_2) \hookrightarrow  \Hom_{{}^0G}(E,V).
\]
We therefore see that the space
\[
T= \Hom_{{}^0G} (E,V_2/V_1)
\]
is a subquotient of $S$ and that $T^\Lambda= \Hom_{G} (E,V_2/V_1) \neq \{ 0\}$. It therefore suffices to
prove that if $S$ is a finite-dimensional $\Lambda$-module (as a vector space over $\bbC$) having a
subquotient $T$ such that $T^\Lambda \neq \{ 0\}$, then $S^\Lambda \neq \{ 0\}$.

Recall that $\Lambda$ is a finitely generated free abelian group, hence isomorphic to $\bbZ^l$ for some
non-negative integer $l$. Defining a $\Lambda$-module structure on the vector space $S$ therefore amounts
to giving $l$ commuting endomorphisms $a_1,\ldots, a_l$ of $S$. The existence of a subquotient $T$ such
that $T^\Lambda \neq \{ 0\}$ means that there exist subspaces $S_1 \subset S_2$ of $S$ stable under
$a_1,\ldots, a_l$ such that the endomorphisms induced by $a_1,\ldots, a_l$ on $T=S_2/S_1$ are the
identity. We deduce that $1$ is a simultaneous eigenvalue of $a_1,\ldots, a_l$. In other words
$S^\Lambda \neq \{0\}$. The proof of the fact that $\Hom_G(\pi,\tau) \neq \{0\}$ is similar.
\end{proof}

\bigskip

This result admits a variant: let us fix a (smooth) central character
\[
\chi  \colon Z(G) \longrightarrow \bbC^\times
\]
and denote by $\caM(G)_\chi$ the full subcategory of $\caM(G)$ whose objects are the smooth
representations admitting the central character $\chi$.

\begin{prop}
Let $(\tau,W)$ be an irreducible supercuspidal representation of $G$ with central character $\chi$. Then
$(\tau,W)$ is projective and injective in $\caM(G)_\chi$.
\end{prop}

\begin{proof}
Let us prove projectivity, the proof of injectivity being similar. We must show that if
$p : (\pi,V) \rightarrow (\tau,W)$ is a non-zero $G$-morphism (hence surjective since $\tau$ is
irreducible), then there exists a $G$-equivariant section $s\colon (\tau,W)\rightarrow (\pi,V)$ such that
$p\circ s=\Id_W$. As in the proof of the lemma, we can assume $(\pi,V)$ is in $\caM(G)_{[\tau]}$. Since
the restrictions of $(\pi,V)$ and $(\tau,W)$ to ${}^0G$ are semisimple, there exists a
${}^0G$-equivariant section $s_0:(\tau,W)\rightarrow (\pi,V)$ such that $p\circ s_0=\Id_W$. Let $\caS$ be
the vector space (non-zero, as we have just seen) of these sections. Since the action of $Z(G)$ on $V$
and $W$ is given by $\chi$, $s_0$ is in fact ${}^0GZ(G)$-equivariant, and the finite abelian group
$\bar G= G/{}^0GZ(G)$ acts on $\caS$ by
\[
\bar g\cdot s = \pi(g)\circ  s \circ  \tau(g^{-1}).
\]
Let $s_0\in \caS$ be non-zero. Let us set $s=|\bar G|^{-1} \sum_{\bar g \in \bar G} g\cdot s_0$. It is
clear that $s$ is $G$-equivariant and that $p\circ s=\Id_W$.
\end{proof}

\section{The center of $\caM(G)_{[\pi]}$}\label{centreMpi}

We will now study more precisely the category $\caM(G)_{[\pi]}$, in particular its center. The definition
of the center of a category is found in \ref{centrecat}.

\subsection{Progenerator of $\caM(G)_{[\pi]}$} \label{progenpi}
Throughout this section, we fix an irreducible supercuspidal representation $(\pi,V)$ of $G$. Our goal is
to describe the category $\caM(G)_{[\pi]}$ as the category of unital right modules over a certain unital
ring. To do this, we will exhibit a small progenerator (in the sense of \ref{rmqfonct}) of
$\caM(G)_{[\pi]}$. Theorem \ref{rmqfonct} then ensures that such a description is possible. Since Theorem
\ref{rmqfonct} is stated without proof, and we need an explicit version tailored to our context, we will
provide complete proofs for all results.

We saw in Proposition \ref{restra0G} that the restriction $\res_{{}^0G}^{\; \,  G} (\pi)$ is a direct sum
of irreducible compact representations of ${}^0G$. Let $(\rho,W)$ be one of these representations of
${}^0G$, and let us set
\[
(\Pi,V_\Pi)=\ind_{{}^0G}^{\; \,  G} (\rho,W) .
\]
Note that the class of $\Pi$ does not depend on the choice of $(\rho,W)$. It actually depends only on
$[\pi]$. This follows from the fact established in Proposition \ref{restra0G} that all the irreducible
components of the restriction of $\pi$ to ${}^0G$ are conjugate.

\begin{prop}
The representation $\Pi$ is a small progenerator of $\caM(G)_{[\pi]}$.
\end{prop}

\begin{proof}
Since ${}^0G$ is open in $G$, and unimodular, the adjunction (\ref{indres}) gives us, for any smooth
representation $(\tau,E)$ in $\caM(G)_{[\pi]}$, a natural isomorphism:
\begin{equation}\label{Frt}
 \Hom_G(\Pi, \tau)\simeq    \Hom_{{}^0G}(\rho, \res_{{}^0G}^{\; \,  G}\tau).
\end{equation}
Since $\rho$ is projective in $\caM({}^0G)$, according to Proposition \ref{conscomp}, and since according
to Proposition \ref{restra0G} and Corollary \ref{conscomp}, the restriction to ${}^0G$ of any
representation in $\caM(G)_{[\pi]}$ is semisimple, it follows that $\Pi$ is projective. (This is a
general result: a left adjoint of an exact functor preserves projectives.) On the other hand, $\Pi$ is in
$\caM(G)_{[\pi]}$: its restriction to ${}^0G$ is a sum of conjugates of $\rho$, so that every irreducible
subquotient of $\Pi$ is in $[\pi]$ by Proposition \ref{restra0G}. Let us now show that
$\Hom_G(\Pi, \tau)$ is non-zero, for all $(\tau,E)\in \caM(G)_{[\pi]}$. Any irreducible component of the
restriction of $\tau$ to ${}^0G$ is isomorphic to an irreducible component of the restriction of $\pi$ to
${}^0G$. Since these are all conjugate under $G$, and since $\tau$ is $G$-stable, we see that the
restriction of $\tau$ to ${}^0G$ contains an irreducible component isomorphic to $\rho$. Therefore, since
the restriction to ${}^0G$ of any representation in $\caM(G)_{[\pi]}$ is semisimple,
$\Hom_{{}^0G}(\rho, \res_{{}^0G}^{\; \,  G}\tau)$ is non-zero. This shows that $\Pi$ is a generator
according to a remark made in \ref{rmqfonct}.

Moreover, $\Pi$ is finitely generated. Indeed, the space ${}^0G\backslash G$ is discrete since ${}^0G$ is
open in $G$. Let us choose a system of representatives $(g_i)_i$ in $G$ for ${}^0G\backslash G$. The
space $V_\Pi= \ind_{{}^0G}^{\; \,  G} (W)$ is then generated by the functions $(f_{i,w})_{i,w}$
characterized by
\begin{equation}\label{fiw}
 f_{i,w}(g_i)=w \quad \text{ and }   \supp
  f_{i,w}\subset {}^0Gg_i.
\end{equation}
Now we have $\Pi(g_{k}^{-1}g_i)\cdot f_{i,w}=f_{k,w}$. The space $\ind_{{}^0G}^{\; \,  G} (W)$ is then
generated as a representation of $G$ by the functions $f_{i_0,w}$, where we have arbitrarily fixed an
index $i_0$. It is convenient to take $i_0$ such that ${}^0Gg_{i_0}={}^0G$, and we may assume without
loss of generality that $g_{i_0}=\mathbf{1}_G$. Let us then set $f_{i_0,w}=f_w$. On the other hand, the
vector space generated by the $f_w$, $w\in W$ is stable under the action of ${}^0G$, and isomorphic to
$(\rho,W)$. Since $(\rho,W)$ is finitely generated, we conclude that $\Pi$ is finitely generated and this
finishes proving the proposition.
\end{proof}

\begin{rmq}
The progenerator $(\Pi,V_\Pi)$ is not unique; other choices are possible. We provide another one here
which will be useful later. Let us set
\[
(\Pi_1,V_{\Pi_1})= \ind_{{}^0 G}^{\; \,  G} (\res_{{}^0 G}^{\; \,  G} (\pi,V)).
\]
The proof that it is indeed a progenerator is similar to that given for $(\Pi,V_\Pi)$. To show that
$\Pi_1$ is finitely generated, we use the fact that $\res_{{}^0 G}^{\; \,  G} (\pi,V)$ is of finite
length. On the other hand, the second Mackey isomorphism \ref{Mackey}, whose hypotheses are satisfied by
the final remark of Section \ref{indcomp}, gives us
\[
\Pi_1 \simeq  (\ind_{{}^0 G}^{\; \,  G}  \mathrm{Triv})\otimes \pi,
\]
where $\mathrm{Triv}$ is the trivial representation of ${}^0G$. Recall that
$\ind_{{}^0 G}^{\; \,  G} \mathrm{Triv}$ is the group algebra $\bbC[\Lambda(G)]$, which is the algebra of
functions on the variety $\caX(G)$ of unramified characters of $G$ (see \ref{carnonram}).
\end{rmq}

Let us set $\caR=\End_G (\Pi)=\Hom_G(\Pi,\Pi)$. \index[not]{R@$\caR$} For any smooth representation
$(\tau,E)$ of $G$, $\Hom_G(\Pi,\tau)$ is a right $\caR$-module. Let us now consider the functor
\begin{equation}\label{EQC1}
 F_\Pi \colon \caM(G)_{[\pi]}  \longrightarrow  \caM(\caR)_d, \,
\quad (\tau,E)\mapsto \Hom_G(\Pi,\tau).
\end{equation}
Since $V_\Pi$ is naturally a left $\caR$-module, for any unital right $\caR$-module $M$, we can form
$M\otimes_\caR V_\Pi$. The group $G$ acts on $M\otimes_\caR V_\Pi$ by
$g\cdot (m\otimes f)=m\otimes \Pi(g)\cdot f$, where $g\in G$, $m\in M$, and $f \in V_\Pi$. Note that this
representation of $G$ is in $\caM(G)_{[\pi]}$, because $M\otimes_\caR V_\Pi$ is a quotient of a direct
sum of representations isomorphic to $\Pi$. This defines a functor
\begin{equation}\label{EQC2}
 G_\Pi \colon  \caM(\caR)_d  \longrightarrow  \caM(G)_{[\pi]}, \quad M
\mapsto   M\otimes_\caR V_\Pi.
\end{equation}

\begin{thm}
The functor $F_\Pi$ establishes an equivalence of categories between $\caM(G)_{[\pi]}$ and
$\caM(\caR)_d$, a quasi-inverse being given by $G_\Pi$.
\end{thm}
\begin{proof}
Let us first look at the composition
\[
\caF \colon \caM(G)_{[\pi]} \rightarrow \caM(G)_{[\pi]}, \quad (\tau,E)\mapsto \Hom_G(\Pi,\tau)\otimes_{\caR}V_\Pi
\]
Let us construct a natural transformation between $\caF$ and the identity functor of $\caM(G)_{[\pi]}$ by
setting:
\[
\mathrm{ev}_E \colon  \Hom_G(\Pi,\tau)\otimes_{\caR}V_\Pi \rightarrow E,\quad \alpha\otimes w \mapsto \alpha(w).
\]
We leave the (elementary but tedious) verification that this is indeed a natural transformation to the
reader. We show that this natural transformation is an isomorphism. Note for this that
\[
\caF(V_\Pi)= \Hom_G(\Pi,\Pi)\otimes_{\caR}V_\Pi=\caR\otimes_{\caR}V_\Pi
\]
and that $\mathrm{ev}_\Pi$ is the natural morphism $\caR\otimes_{\caR}V_\Pi \simeq V_\Pi$. On the other
hand, the functor $F_\Pi$ is exact ($(\Pi,V_\Pi)$ is projective) and commutes with direct sums
($(\Pi,V_\Pi)$ is finitely generated, hence small in the sense of category theory). The functor $G_\Pi$
is right exact and commutes with direct sums (it is a left adjoint, cf. \ref{OHom}). It follows that the
composite functor $\caF$ also possesses these properties. Let $(\tau,E)$ be an object of
$\caM(G)_{[\pi]}$. Since $\Pi$ is a progenerator of this category, we can realize $(\tau,E)$ as the
cokernel of a morphism between two modules which are direct sums (possibly infinite) of modules
isomorphic to $V_\Pi$ (Lemma \ref{rmqfonct}).

We therefore have an exact sequence of the form
\[
\oplus_{i\in I} V_\Pi  \rightarrow   \oplus_{i\in J} V_\Pi
\rightarrow  E \rightarrow 0 .
\]
Since $\caF$ commutes with direct sums and is right exact, we obtain a commutative diagram

\begin{equation}
\begin{CD}
\oplus_{i\in I}\caF(V_\Pi) @>>>   \oplus_{j\in J}\caF(V_\Pi) @>>>  \caF(E) @>>> 0  \\
@VVV              @VVV                          @VV{\mathrm{ev}_E}V\\
\oplus_{i\in I} V_\Pi  @>>> \oplus_{j\in J} V_\Pi @>>>  E @>>> 0
\end{CD}
\end{equation}
Since $\caF(V_\Pi)\simeq V_\Pi$, the two left vertical arrows are isomorphisms. It follows that
$\mathrm{ev}_E$ is an isomorphism, and therefore that $\mathrm{ev}$ is a natural isomorphism.

Now let $\caG$ be the functor
\[
\caG \colon \caM(\caR)_d \longrightarrow \caM(\caR)_d,\quad  M \mapsto \Hom_G(V_\Pi,M\otimes_{\caR}V_\Pi).
\]
To complete the proof, we must show that there exists a natural isomorphism from $\caG$ to the identity
of $\caM(\caR)_d$. For any module $M$ in $\caM(\caR)_d$, let us define a morphism of right $\caR$-modules
\[
\theta_M: M \rightarrow \Hom_G(V_\Pi,M\otimes_{\caR}V_\Pi).
\]
Let us set, for all $m\in M$, and for all $v$ in $V_\Pi$, $f_m(v)=m\otimes v$ and define $\theta_M$ by
$\theta_M(m)=f_m$. It is clear that the $\theta_M$ define a natural transformation from the identity of
$\caM(\caR)_d$ to $\caG$, and that moreover $\theta_{\caR}$ realizes the natural isomorphism between
$\caR$ and $\Hom_G(V_\Pi,\caR \otimes_{\caR}V_\Pi)$. As for $\caF$, we see that $\caG$ is right exact and
commutes with direct sums, and that for any module $M$, we can find an exact sequence of the form
\[
\oplus_{i\in I} \caR  \rightarrow   \oplus_{i\in J} \caR \rightarrow  M \rightarrow 0.
\]
We can then use the same argument as above to show that each $\theta_M$ is an isomorphism.
\end{proof}

\subsection{Another description of $\caM(G)_{[\pi]}$}
We continue with the notation of the previous section. In particular, $(\pi,V)$ is an irreducible
supercuspidal representation of $G$, $(\rho,W)$ is an irreducible (compact) component of the restriction
of $(\pi,V)$ to ${}^0G$ and $(\Pi,V_\Pi)=\ind_{{}^0G}^{\; \,  G} (\rho,W)$.

It is well known that for finite groups, the algebra of self-intertwining operators of an induced
representation is isomorphic to a certain convolution algebra. We will do the same here and establish an
isomorphism between $\caR$ and a convolution algebra. Let $\caH(G,\rho)$
\index[not]{H(G,rho)@$\caH(G,\rho)$} denote the set of functions $\phi \colon G\rightarrow \End_\bbC (W)$
satisfying:
\begin{align}
 &(i) \quad  \phi(hgh')=\rho(h)\phi(g)\rho(h'), \; (g\in G), \; (h,h' \in {}^0G)\\
\nonumber  (ii) \quad \supp(\phi)& \text{ is a finite union of left cosets modulo }
{}^0G
\end{align}

Let us define on $\caH(G,\rho)$ the convolution product given by
\[
\phi*\psi(x)= \sum_{\bar g\in G/{}^0G} \phi(g)\psi(g^{-1}x).
\]
(We sum over a system of representatives and verify that the result does not depend on the choice of
these. The associativity of the convolution product is verified by calculation). The identity element
$\phi_e$ of this algebra is characterized by $\phi_e(\mathbf{1}_G)=\Id_W$ and $\phi_e(g)=0$ if
$g\notin {}^0G$.

Recall that we defined in (\ref{fiw}), for all $w\in W$, a function $f_w$ in
$V_\Pi=\ind_{{}^0G}^{\; \,  G}(W)$ characterized by the properties $f_w(\mathbf{1}_G)=w$ and
$\supp f_w \subset {}^0G$.

\begin{prop}
The maps:
\begin{align*}
 \caR &\longrightarrow \caH(G,\rho),\qquad \alpha \mapsto \phi_\alpha \\
\text { where } \quad \phi_\alpha(g)(w)&= \alpha(f_w)(g),\quad (g\in G), \; (w\in W).
\end{align*}
and
\begin{align*}
\caH(G,\rho)&\longrightarrow\caR ,\qquad \phi \mapsto \alpha_\phi \\
\text { where } \quad   \alpha_\phi(f)(x)=\sum_{\bar g \in {}^0G \backslash G} &\phi(g)(f(g^{-1}x)), \quad (x\in G),
(f\in V_\Pi).
\end{align*}
are $\bbC$-algebra isomorphisms inverse to each other.
\end{prop}

\begin{proof}
It is straightforward to verify that $\alpha_{\phi_\alpha}=\alpha$ and $\phi_{\alpha_\phi}=\phi$, once we
note that any function $f \in V_\Pi$ decomposes as
\begin{equation}\label{decf}
 f=\sum_{\bar g \in {}^0G\backslash G} r(g^{-1})\cdot f_{f(g)},
\end{equation}
this sum being finitely supported according to the support property of $f$. We now show that
$\alpha \mapsto \phi_\alpha$ is an algebra morphism. Let $\alpha, \beta \in \caR$. We calculate, for all
$x\in G$ and all $w\in W$,
\begin{align}\label{EQ1}
&\phi_{\alpha\beta}(x)(w)=(\alpha\beta)(f_w)(x)=\alpha(\beta(f_w))(x)\\
\nonumber &= \alpha\left(\sum_{\bar g \in {}^0G  \backslash G} r(g^{-1})\cdot f_{\beta(f_w)(g)}\right) (x)
= \sum_{\bar g \in {}^0G  \backslash G    }  r(g^{-1})\cdot \alpha(f_{\beta(f_w)(g)})(x)\\
\nonumber &=  \sum_{\bar g \in {}^0G  \backslash G  }  \alpha(f_{\beta(f_w)(g)})(xg^{-1})
\end{align}
We used in this calculation, besides the finiteness of the sums, the decomposition (\ref{decf}) of
$\beta(f_w)$, and the fact that $\alpha$ commutes with the action of $G$ by right translation. On the
other hand we have
\begin{align}\label{EQ2}
&((\phi_\alpha*\phi_\beta)(x))(w)=  \sum_{\bar y \in G / {}^0G}\phi_\alpha(y)(\phi_\beta(y^{-1}x)(w))\\
\nonumber &=\sum_{\bar y\in G/{}^0G} \alpha(f_{\phi_\beta(y^{-1}x)(w)})(y)
=\sum_{\bar y\in G / {}^0G} \alpha(f_{\beta(f_w)(y^{-1}x)})(y)
\end{align}
A simple change of variable in this sum shows that we indeed have equality between (\ref{EQ1}) and
(\ref{EQ2}).
\end{proof}

Let us set $\caA=\caH(G,\rho)$ to lighten the notation. \index[not]{A@$\caA$} We can interpret the
equivalence of categories between $\caM(G)_{[\pi]}$ and $\caM(\caR)_d$ as an equivalence between
$\caM(G)_{[\pi]}$ and $\caM(\caA)_d$, by the isomorphism $\caR \simeq \caA$.

Let $(\tau,E)$ be a representation in $\caM(G)_{[\pi]}$. The algebra $\caR=\End_G (\Pi)$ acts (on the
right) on $\Hom_G(\Pi,\tau)$. We obtain a right $\caA$-module structure on $\Hom_G(\Pi,\tau)$ by
transport of structure, i.e.,
\[
s\cdot \phi  = s \circ  \alpha_\phi,\quad  (s\in \Hom_G(\Pi,\tau)),\,
(\phi \in \caA).
\]
Note that this depends on $\tau$, and not just on $\res_{{}^0G}^{\; \,  G} \tau$.

On the other hand, the algebra $\caA$ acts on the right on
$\Hom_{{}^0G}(\rho,\res_{{}^0G}^{\; \,  G} \tau)$ by
\[
(\psi\cdot \phi)(w)=  \sum_{\bar g \in G/{}^0G} \tau(g)\cdot \psi(\phi(g^{-1})\cdot w),
\]
for all $\psi \in \Hom_{{}^0G}(\rho,\res_{{}^0G}^{\; \,  G} \tau)$, for all $\phi \in \caA$, and for all
$w \in W$. The natural adjunction isomorphism (\ref{Frt})
\[
\Hom_G(\Pi, \tau)\simeq    \Hom_{{}^0G}(\rho, \res_{{}^0G}^{\; \, G}\tau)
\]
is then a morphism of right $\caA$-modules. This is easily verified using the explicit form of the
adjunction isomorphism. It is realized by the map
\begin{equation}\label{def1iso}
 \alpha \mapsto \psi_\alpha ,\quad \psi_\alpha(w)=\alpha(f_w), \quad
(w\in W)
\end{equation}
where $f_w$ is the function in $\ind_{{}^0G}^{\; G} W$ characterized by $\supp f_w={}^0G$ and
$f_w(\mathbf{1}_G)=w$. As shown by a straightforward calculation (using (\ref{decf})), its inverse is
\begin{equation}\label{def2iso}
 \psi \mapsto \alpha_\psi,\quad \alpha_\psi(f)= \sum_{\bar g \in
  {}^0G\backslash G} \tau(g)\cdot  \psi(f(g^{-1})), \quad (f \in
\ind_{{}^0G}^G W).
\end{equation}

It easily follows that the functor
\[
\caM(G)_{[\pi]} \rightarrow  \caM(\caA)_d, \quad \tau \mapsto
\Hom_{{}^0G}(\rho, \res_{{}^0G}^{\; \,  G}\tau).
\]
is an equivalence of categories whose inverse is given by
\[
\caM(\caA)_d \rightarrow \caM(G)_{[\pi]} ,\quad   M \mapsto M\otimes_\caA V_\Pi
\]
where we consider $V_\Pi$ as a left $\caA$-module by transport of structure.

On the other hand, in the case $\tau=\Pi$, the isomorphisms (\ref{def1iso}) and (\ref{def2iso}) give
\[
\caR= \Hom_G(\Pi, \Pi) \simeq    \Hom_{{}^0G}(\rho,\res_{{}^0G}^{\; \,  G}\Pi).
\]
Let us examine more closely the right-hand side. Recall that the space ${}^0G\backslash G$ is discrete,
and let us fix a system of representatives $\{g_i\}_i$. The space $V_\Pi$ is generated by the functions
$f_{i,w}$ defined in (\ref{fiw}). Let $(V_\Pi)_i$ denote the subspace of $V_\Pi$ of functions supported
in ${}^0Gg_i$. Then it is clear that
\[
V_\Pi= \bigoplus_i  (V_\Pi)_i,
\]
and moreover, since ${}^0G$ is normal in $G$, each $(V_\Pi)_i$ is stable under the action of ${}^0G$.
Indeed, for all $g_0 \in {}^0G$, $\Pi(g_0)\cdot f_{i,w}$ is still supported in ${}^0Gg_i$. We calculate
\begin{align*}
 \Pi(g_0)\cdot f_{i,w}(g_i)&=
f_{i,w}(g_ig_0)=f_{i,w}((g_ig_0g_i^{-1})g_i)=\rho(g_ig_0
g_i^{-1})\cdot f_{i,w}(g_i)\\
&= \rho^{g_i^{-1}}(g_0)\cdot w,
\end{align*}
which shows that
\[
W \rightarrow (V_\Pi)_i, \quad  w \mapsto f_{i,w}
\]
intertwines the representations $(\rho^{g_i^{-1}},W)$ and $(\res_{{}^0G}^{\; \,  G} \Pi,  (V_\Pi)_i)$.

As a vector space, we therefore see that
\[
\Hom_{{}^0G}(\rho, \res_{{}^0G}^{\; \,  G}\Pi)=
\Hom_{{}^0G}(\rho,\oplus_i(V_\Pi)_i )=\bigoplus_i\Hom_{{}^0G}(\rho,\rho^{g_i^{-1}}).
\]
If $\rho^{g_i^{-1}} \simeq \rho$ then $\Hom_{{}^0G}(\rho,\rho^{g_i^{-1}})$ is $1$-dimensional, and zero
otherwise.

Assume that $\rho^{g_i^{-1}} \simeq \rho$. An intertwining operator between $\rho$ and $\rho^{g_i^{-1}}$
is an endomorphism $A$ of $W$ such that, for all $g_0 \in {}^0G$
\begin{equation}\label{compat}
\rho(g_i g_0 g_i^{-1}) \circ A  = A \circ \rho(g_0)
\end{equation}
If $A$ is such an operator, let us define the function
\[
f_{i,A}\colon G \rightarrow \End_\bbC (W),
\]
supported in ${}^0G g_i$, and characterized by
\[
f_{i,A}(g_i)=A \quad \text{ and } f_{i,A}(g_0g_ig_1)= \rho(g_0)\circ
A \circ \rho(g_1)
\]
for all $g_0,g_1 \in {}^0G$. Relation (\ref{compat}) shows that $f_{i,A}$ is well-defined. Conversely,
such a function $f_{i,A}$ defines an intertwining operator $A$ between $\rho$ and $\rho^{g_i^{-1}}$. The
preceding discussion naturally reveals the convolution algebra $\caA=\caH(G,\rho)$ and the isomorphism
between $\Hom_{{}^0G}(\rho,\res_{{}^0G}^G\Pi)$ and $\caA=\caH(G,\rho)$.

\subsection{Structure of $\caA$}\label{strucMA}
We continue with the notation of the previous sections. For all $x \in G$, let $\caF_x(\rho)$ denote the
space of functions $\phi$ in $\caA$ whose support is contained in ${}^0Gx$. Let $\phi \in \caF_x(\rho)$.
We then have:
\[
\phi(x)\rho(x^{-1}gx)=  \phi(gx)=\rho(g)\phi(x),\quad (g\in{}^0G),
\]
and therefore we have for all $x \in G$ a linear isomorphism
\[
\phi \mapsto \phi(x),\quad \caF_x(\rho) \simeq \Hom_{{}^0G}( \rho^x,\rho).
\]
We deduce that $\dim \caF_x(\rho)\leq 1$ for all $x \in G$ and that $\dim \caF_x(\rho)=1$ if and only if
$x$ is in the stabilizer $N=N_G(\rho)$ of $\rho$ (i.e., $N$ is the set of $g\in G$ such that
$\rho^g\simeq \rho$, with the notation of \ref{entrel}). In particular, we have $\caA=\caH(N,\rho)$.

Let $\phi,\psi \in \caA$ with $\supp \phi \subset  {}^0Gx$ and $\supp \psi \subset {}^0Gy$, $x,y$ in $N$.
A direct calculation shows that
\begin{equation}\label{Gxy}
  \supp \phi *\psi \subset  {}^0Gxy \quad \text{ and }
\phi *\psi(xy)=\phi(x)\psi(y).
\end{equation}
We deduce that if $\phi$ is non-zero, it is invertible in $\caA$, its inverse being given by an element
whose support is in ${}^0Gx^{-1}$. On the other hand, since ${}^0G\backslash G$ is abelian
(\ref{Lambda}), ${}^0Gxy ={}^0Gyx$, and therefore if $\phi$ and $\psi$ are as above,
\begin{equation*}
\phi *\psi=c\;  \psi*\phi \; \text{ for some } c\in \bbC^\times.
\end{equation*}

Since $N$ contains ${}^0G Z(G)$ and the latter is of finite index in $G$ (Proposition \ref{Lambda}),
${}^0G\backslash N$ is a free abelian group of the same rank as $\Lambda(G)={}^0G\backslash G$. Let us
then choose a basis ${}^0Gx_1,\ldots, {}^0Gx_n$ of the lattice ${}^0G\backslash N$ and fix non-zero
elements $\phi_i$, $i=1,\ldots, n$ such that $\phi_i \in \caF_{x_i}(\rho)$. For each pair $(i,j)$ of
indices, we have a complex number $c_{ij}$ such that
\begin{equation}\label{laurentpol}
\phi_i*\phi_j=c_{ij}\; \phi_j * \phi_i.
\end{equation}
The algebra $\caA$ is generated by the $\phi_i$, $i=1,\ldots, n$, and the relations (\ref{laurentpol})
allow us to write any element of $\caA$ as a linear combination of monomials in the $\phi_i$ and their
inverses. In other words, as an algebra, $\caA$ is a ``twisted'' Laurent polynomial algebra in $n$
variables. It is then easy to see that the group of invertible elements $\caA^\times$ consists of the
monomials
\begin{equation}\label{invA}
  c \; \phi_1^{e_1}*\cdots*\phi_n^{e_n},\quad c\in \bbC^\times,
 e_1,\ldots ,e_n \in \bbZ.
\end{equation}

\begin{lemme}
For any right $\caA$-module $M$, let
\[
\caA \rightarrow \End_\bbC(M) ,\quad \phi \mapsto \phi_M
\]
denote the morphism giving the $\caA$-module structure. Let us fix a system of representatives of the
isomorphism classes of simple right $\caA$-modules, and consider the morphism
\[
\caA \rightarrow \prod  \End_\bbC(M) ,\quad (\phi \mapsto \phi_M)
\]
where the product is over this system of representatives. Then this morphism is injective.
\end{lemme}

\begin{proof}
The separation lemma \ref{nilp} shows that an element $\phi$ of the kernel of this morphism is
necessarily nilpotent. It immediately follows that $u=\mathbf{1_\caA}+\phi$ is invertible. Similarly,
$u'=\mathbf{1_\caA}-\phi^2$ is invertible and $u'=2u-u^2$. Recall that the invertible elements of $\caA$
are of the form given by (\ref{invA}). It follows easily that $u$ must be a scalar. Thus $\phi$ is a
scalar, and then necessarily $\phi=0$.
\end{proof}

\begin{cor}
$(i)$ The algebra $\caA$ is an integral domain, hence has no non-zero nilpotent elements.

$(ii)$ The algebra $\caA$ is commutative if and only if the restriction of $\pi$ to ${}^0G$ is
multiplicity-free.
\end{cor}

\begin{proof}
Point $(i)$ follows from the structure of $\caA$ described above: $\caA$ is graded by the finitely
generated free abelian group ${}^0G\backslash N$, the homogeneous component of degree ${}^0Gx$ consisting
of the elements supported in ${}^0Gx$; it is one-dimensional, and its non-zero elements are invertible
(\ref{Gxy}). Choose a total order on ${}^0G\backslash N$ compatible with the group law; the highest
degree component of a product $ab$ is the product of the highest degree components of $a$ and $b$, which
is non-zero. Hence $\caA$ is an integral domain. Let $m$ be the multiplicity of $\rho$ in
$\res_{{}^0G}^{\; \,  G} \pi$. The equivalence of categories $ \caM(G)_{[\pi]} \simeq \mathbf{mod-\caA}$
shows that the simple right $\caA$-modules are isomorphic to
\[
\Hom_{{}^0G}(\rho, \res_{{}^0G}^{\; \,  G} (\pi\otimes \omega))
\]
for some unramified character $\omega$ of $G$. Of course\footnote{Be careful, the following equality is
an equality of ${}^0G$-modules, but not of $\caA$-modules},
$\res_{{}^0G}^{\; G} (\pi\otimes \omega)=\res_{{}^0G}^{\; \,  G} (\pi) $ and therefore each simple right
$\caA$-module is of dimension $m$. If $m=1$, the lemma then shows that $\caA$ is commutative. Conversely,
if $\caA$ is commutative, each simple right $\caA$-module is of dimension $1$.
\end{proof}

\subsection{The center of $\caA$}\label{centreA}
\label{T0GIrr} We continue with the same notation as in the previous sections. We begin with some
preliminary remarks. We have already introduced the notation $N=N_G(\rho)$ for the stabilizer in $G$ of
$\rho$. The space $W$ of the representation $\rho$ is a ${}^0G$-stable subspace of $V$.

Let us set
\[
H=\{ g\in G\mid \pi(g)\cdot W =W  \}
\]
and let $\rho_H$ be the representation of $H$ on $W$ obtained by restriction of $\pi$. Since $\rho_H$
extends $\rho$, it is irreducible. On the other hand, $H$ contains ${}^0GZ(G)$, and the group
$G/{}^0GZ(G)$ is abelian, so $H$ is normal in $G$ and of finite index in $G$. For all $g \in G$, the
subspace $ \pi(g)\cdot W$ is then stable under the action of ${}^0G$. Let us denote by
$(\pi_{|\pi(g)\cdot W}, \pi(g)\cdot W)$ this representation of ${}^0G$. Thus for example, with
$g=\mathbf{1}_G$, $(\pi_{| W}, W)=(\rho,W)$. The map $\pi(g)$ intertwines $(\rho^{g},W)$ and
$(\pi_{|\pi(g)\cdot W}, \pi(g)\cdot W)$. It is then clear that $H\subset N$. Moreover, the spaces
$\pi(g)\cdot W$ are stable under $H$, and isomorphic to $\rho_H^{g}$ as a representation of $H$.

Moreover, since $V$ is irreducible, we have
\begin{equation}\label{sumdir}
V=\sum_{\bar g \in G/H} \pi(g)\cdot W.
\end{equation}

Let us introduce the group
\[
\bar H=\{ g \in G \mid \rho_H^g\simeq \rho_H \},
\]
the normalizer in $G$ of the representation $\rho_H$ of $H$. If $H=\bar H$, the representations
$\{ \pi(g)\cdot W \}$, where $g$ varies in a system of representatives of $G/H$, are pairwise
inequivalent. Indeed, if for $g_1,g_2 \in G$,
\[
\pi(g_1)\cdot W \simeq \pi(g_2)\cdot W
\]
as a representation of $H$, we have
\[
\pi(g_1^{-1}g_2)\cdot W \simeq  W,
\]
hence $g_1^{-1}g_2 \in \bar H=H$. We deduce that the components $\pi(g)\cdot W$ are linearly independent
and therefore that the sum (\ref{sumdir}) is direct.

If $H\neq \bar H$, let $H_1 \neq H$, $H_1 \subset \bar H$ such that $H_1/H$ is cyclic (recall that
$\bar H/H$ is finite abelian) and let $e$ be its order. Then $(\rho_H,W)$ admits an extension
$(\rho_{H_1},W)$ to $H_1$. Let $h_0 \in H_1$ be such that its image in $H_1/H$ is a generator of this
group, and set $h_0^{e}=k\in H$. Fix a non-zero intertwining operator
$\phi : (\rho_H,W) \rightarrow (\rho_H^{h_0}, W)$, bijective by irreducibility, and normalized (using
Schur's lemma) so that $\phi^{e}=\rho_H(k)$. If $h_1=h_0^j h$, $h_1 \in H_1$, $h \in H$, we then set
\[
\rho_{H_1}(h_1)= \phi^j \rho_H(h),
\]
which defines $\rho_{H_1}$.

Let us set $E=\sum_{\bar g \in H_1/H} \pi(g)\cdot W$. It is a representation of $H_1$, whose restriction
to $H$ is a direct sum of representations isomorphic to $\rho_H$. Let $(\rho_1,W_1)$ be an irreducible
representation of $H_1$, $W_1 \subset E$. We therefore have
\[
\Hom_H(W,W_1) \neq 0.
\]
Moreover, $H_1$ acts on this space by
\[
h\cdot \phi=\rho_1(h)\circ \phi \circ \rho_{H_1}(h)^{-1}, \quad h\in H_1, \;
\phi \in \Hom_H(W,W_1).
\]
This action factors through the cyclic group $H_1/H$, and therefore there exists a non-zero simultaneous
eigenvector $\phi$ for all $h \in H_1/H$, that is to say there exists $\chi \in \widehat{H_1/H}$ and
$\phi \in \Hom_H(W,W_1)$ such that $h\cdot \phi=\chi(h)\phi$ for all $h \in H_1/H$, which can be
translated as
\[
\phi \in \Hom_{H_1}(\rho_{H_1}, \rho_1 \otimes \chi^{-1}).
\]
Since $\rho_{H_1}$ and $\rho_1$ are irreducible, $\phi$ is an isomorphism, and we deduce that the
restriction of $\rho_1$ to $H$ is irreducible. Replacing $W$ with $W_1$ replaces $H$ with a group
containing $H_1$, and by iterating this process, we may eventually assume that we have chosen $W$ such
that $H=\bar H$ (which we will do in what follows), and thus
\begin{equation} \label{VgW}
V=\bigoplus_{\bar g \in G/H} \pi(g)\cdot W.
\end{equation}
We then have $\pi \simeq \ind_{H}^G (\rho_H)$. More precisely, the decomposition (\ref{VgW}) is a
decomposition into irreducible subrepresentations of $H$, and $\pi$ is isomorphic to the induced
representation of any of these components. It follows that $m$, the multiplicity of $\rho$ in
$\res_{{}^0G}^{\; \,  G} \pi$ is equal to $[N:H]$.

Let us now introduce another subgroup of $G$ containing ${}^0G$. To do this, let us recall some elements
of the duality between complex tori and lattices. If $\bbU$ is a complex torus, its group of algebraic
characters $X^*(\bbU)$ is a lattice (a free $\bbZ$-module of rank the dimension of $\bbU$) and
$\bbC[X^*(\bbU)]$ is the algebra of polynomial functions on $\bbU$. If $J$ is a finite subgroup of
$\bbU$, the quotient $\bbU/J$ is again a complex torus, and $X^*( \bbU/J)$ is identified with the
sublattice $L$ of $X^*(\bbU)$ of characters trivial on $J$. The algebra of polynomial functions on
$\bbU/J$ is $\bbC[L]\simeq \bbC[X^*(\bbU)]^J$.

Let us apply this to $\caX(G)$, which according to \ref{varXG} is a torus whose group of algebraic
characters is $\Lambda(G)=G/{}^0G$ and to its finite subgroup (see \ref{finXGpi})
\[
\caX(G)(\pi)=\{  \omega \in \caX(G)\mid \pi\otimes \omega \simeq \pi \}.
\]
The group of algebraic characters of the torus $\caX(G)/\caX(G)(\pi)$ is therefore a sublattice of
$\Lambda(G)=G/{}^0G$, which we can view as a subgroup of $G$ containing ${}^0G$. More explicitly, it is
the group:
\[
T= \bigcap_{ \omega \in \caX(G)(\pi) } \ker \omega .
\]
If $ \omega \in \caX(G)(\pi)$ then the restriction of $\omega$ to $Z(G)$ is trivial, hence
${}^0G Z(G)\subset T$ and therefore $T$ is of finite index in $G$.

\begin{rmq}
The duality between tori and lattices shows that $\caX(G)(\pi)$ consists exactly of the elements of
$\caX(G)$ whose restriction to $T$ is trivial. Moreover the algebra of polynomial functions on
$\caX(G)/\caX(G)(\pi)$ is then $\bbC[T/{}^0G]\simeq\bbC[\Lambda(G)]^{\caX(G)(\pi)}$. On the other hand,
recall (Remarks \ref{clinertcusp}) that the torus $\caX(G)/\caX(G)(\pi)$ is identified with the variety
$\mathbf{Irr}(G)_{[\pi]}$. The algebra of polynomial functions on $\mathbf{Irr}(G)_{[\pi]}$ is therefore
identified with $\bbC[T/{}^0G]$.
\end{rmq}

\begin{lemme}
We have $T \subset H \subset N$ with $[N:H]=[H:T]=m$. Let $\rho_T$ be the restriction of $\rho_H$ to $T$.
Then, any irreducible component of $\res_{T}^G (\pi)$ whose restriction to ${}^0G$ contains $\rho$ is
isomorphic to $\rho_T$. The restriction morphism
\begin{equation}\label{resT0G}
 \Hom_T(\rho_T,\res_T^G(\pi)) \rightarrow
\Hom_{{}^0G}(\rho,\res_{{}^0G}^{\; \,  G}(\pi))
\end{equation}
is an isomorphism.
\end{lemme}

\begin{proof}
We have already seen that $H\subset N$ with $[N:H]=m$. Let us show that $T \subset H$: let
$\omega \in \caX(G)$. We then have
\[
\pi \otimes \omega \simeq (\ind_H^G \rho_H)\otimes \omega \simeq
\ind_H^G (\rho_H\otimes  \omega_{|H}).
\]
Consequently, if $\omega_{|H}$ is trivial, $ \pi \otimes \omega \simeq \pi$, so $\omega \in \caX(G)(\pi)$
and $\omega_{|T}=1$. The duality between tori and lattices then shows that $T \subset H$.

Let $(\tau,E)$ be an irreducible component of the restriction of $\pi$ to $T$ whose restriction to
${}^0G$ has a non-trivial $\rho$-isotypic component. Since $\rho$ is also equal to the restriction of
$\rho_T$ to ${}^0G$, this means that we have
\[
\Hom_{{}^0G}(\res_{{}^0G}^{\; \,  T}\rho_T,\res_{{}^0G}^{\; \, T}\tau)\neq \{0\}.
\]
According to Proposition \ref{restra0G} $(iii)$ (more precisely, one must note that the proof of point
$c) \Rightarrow b)$ adapts without difficulty to our context), there exists a character $\eta$ of $T$,
trivial on ${}^0G$ such that
\[
\tau \simeq \rho_T\otimes \eta.
\]
The representation $\tau$ is contained in the restriction to $T$ of an irreducible component, say
$\tilde \tau$, of the restriction of $\pi$ to $H$. We have
\[
\tilde  \tau \simeq \rho_H \otimes \tilde \eta
\]
for some character $\tilde \eta$ of $H$ extending $\eta$ (we again apply point $c) \Rightarrow b)$ of
Proposition \ref{restra0G}, (iii)). We deduce that $\pi\simeq \ind_H^G (\rho_H\otimes \tilde \eta)$,
since $\rho_H\otimes \tilde \eta$ is isomorphic to one of the components of the decomposition
(\ref{VgW}). Let us further extend $\tilde \eta$ to a character $\hat \eta$ of $G$. Then
$\hat \eta_{|T}=\eta$ and
\[
\pi\simeq \ind_H^G (\rho_H\otimes \tilde \eta)\simeq (\ind_H^G \rho_H )
\otimes \hat \eta =\pi\otimes \hat \eta.
\]
This shows that $\hat \eta \in \caX(G)(\pi)$ so $\hat \eta_{|T}=1$. Consequently, $\eta=1$ and
$\tau\simeq \rho_T$.

We now use the adjunction isomorphism
\begin{equation}\label{eq:THG}
\Hom_T(\rho_T, \res_T^G (\pi))=\Hom_G(\ind_T^G (\rho_T),\pi).
\end{equation}
We will now show that the above implies that the left-hand side is of dimension $m$. Indeed, the
restriction to ${}^0G$ defines a natural map from $\Hom_T(\rho_T, \res_T^G (\pi))$ to
$\Hom_{{}^0G}(\rho, \res_{{}^0G}^{\; \,  G} \pi)$, since the spaces of $\rho$ and $\rho_T$ are the same.
Conversely, we have just shown that any morphism in $\Hom_{{}^0G}(\rho, \res_{{}^0G}^{\; \,  G} (\pi))$
extends uniquely to a $T$-equivariant morphism. The two spaces are therefore isomorphic (as vector
spaces), and this proves that (\ref{resT0G}) is an isomorphism. Let us examine the right-hand side of the
adjunction isomorphism. We have, since $H/T$ is finite abelian,
\begin{equation*}
 \ind_T^H (\rho_T)\simeq \bigoplus_{\xi \in
    \Hom_\bbZ(H/T, \bbC^\times)}\rho_H \otimes \xi.
\end{equation*}
Each $\xi \in \Hom_\bbZ(H/T,\bbC^\times)$ extends to an unramified character $\hat \xi$ of $G$ trivial on
$T$, and consequently in $\caX(G)(\pi)$. It follows that
\begin{align*}
  \ind_T^G (\rho_T)&\simeq  \ind _H^G (\ind_T^H (\rho_T))\simeq \ind_H^G
(  \oplus_{\xi} \rho_H\otimes \xi )\\
&\simeq  \bigoplus_{\xi}\ind_H^G
(\rho_H\otimes \xi)=   \bigoplus_{\xi} \ind_H^G (\rho_H)\otimes \hat \xi\\
&= [H:T] \; \pi.
\end{align*}
The right-hand side of (\ref{eq:THG}) is therefore of dimension $[H:T]$, and therefore $[H:T]=m$.
\end{proof}

We will now describe the center $\caZ$ of the algebra $\caA$.

\begin{thm}
The center of $\caA=\caH(N,\rho)$ is $\caZ=\caH(T,\rho)$. In particular, $\caA$ is a free $\caZ$-module
of rank $m^2$.
\end{thm}

\begin{proof}
Let $\phi \in \caA$. Let us write the support of $\phi$ as a disjoint union of right cosets ${}^0Gg_i$,
$i=1,\ldots, r$ and $\phi$ as a sum $\phi=\phi_1+\cdots +\phi_r$ where $\supp \phi_i={}^0Gg_i$. We
described in Section \ref{strucMA} the algebra $\caA$ as a twisted Laurent polynomial algebra. The
decomposition of $\phi$ above corresponds to a decomposition of a polynomial into monomials. It is then
clear that $\phi$ is in $\caZ$ if and only if each $\phi_i$ is in $\caZ$. We may therefore assume in what
follows that the support of $\phi$ consists of a single coset ${}^0Gg$. Suppose that $\phi \in \caZ$. Let
us show that $\phi \in \caH(T,\rho)$, i.e., that $g \in T$. For all $\omega \in \caX(G)$, $\phi$ acts on
the simple $\caA$-module $\Hom_{{}^0G}(\rho, \pi \otimes \omega)$ by a certain scalar, which we will
denote by $\lambda_{\omega}(\phi)$. The inclusion $\iota\colon W \rightarrow V$ defines an element of
$\Hom_{{}^0G}(\rho, \pi \otimes \omega)$. We then have, according to the definition of the action of
$\caA$ on $\Hom_{{}^0G}(\rho, \pi \otimes \omega)$,
\begin{equation}\label{lop}
 \lambda_{\omega}(\phi)w= (\iota \cdot \phi)(w)=
\pi(g^{-1})\omega(g^{-1})(\phi(g)w),\quad (w \in W).
\end{equation}
In particular, if $\omega$ is trivial, $\lambda_{triv}(\phi)\pi(g)=\phi(g)$. Since $\phi$ is invertible,
$\lambda_{triv}(\phi)\neq 0$ and therefore $g \in H$. Suppose now that $\pi\otimes \omega \simeq \pi$.
Then $\lambda_{\omega}(\phi)=\lambda_{triv}(\phi)$ and this implies that $\omega(g)=1$, for all
$\omega \in \caX(G)(\pi)$, hence $g\in T$.

Our goal is now to show the converse, i.e., that if $g\in T$, then $\phi \in \caZ$. Let $n$ be an element
of $N$. According to the previous lemma, $\rho_T^n\simeq\rho_T$. Let $T_n$ be an intertwining operator
between $\rho_T$ and $\rho_T^n$ and $F_n=\pi(n)T_n$. We then have for all $w \in W$ and all
$h \in {}^0G$:
\begin{align}
\nonumber
F_n\circ \rho(h)(w)&=\pi(n) \circ T_n \circ \rho(h)(w)=\pi(n) \circ
\rho(n^{-1}hn) \circ T_n(w)\\
&=\pi(hn)\circ T_n(w)= \pi(h)\circ F_n(w)
\end{align}
and therefore $F_n$ is in $\Hom_{{}^0G}(\rho,\pi)$. We can easily show using (\ref{VgW}) that as $n$
varies in a system of representatives $\caT$ of $N/H$, the set of $F_n$ thus obtained is a basis of
$\Hom_{{}^0G}(W,V)$. For $g \in T$, $\phi(g)$ is equal, up to a non-zero scalar factor, to $\rho_T(g)$.
We can renormalize $\phi$ to have $\phi(g)=\rho_T(g)$, which we assume in what follows. We have, for all
$w \in W$, for all $n \in \caT$, and all $\omega \in \caX(G)$, by definition of the right action of
$\caA$ on $\Hom_{{}^0G}(\rho,\pi\otimes \omega)$,
\begin{align}
(F_n\cdot \phi)(w)&=\omega(g^{-1}) \pi(g^{-1})\circ F_n(\phi(g)\cdot w)\\
\nonumber  &=\omega(g^{-1}) \pi(g^{-1}) \circ \pi(n)\circ T_n(\rho_T(g)\cdot w)\\
\nonumber &=\omega(g^{-1}) \pi(g^{-1})\circ  \pi(n)\circ \rho_T(n^{-1}gn)T_n(w)\\
\nonumber  &=\omega(g^{-1})\pi(n)\circ  T_n(w)=\omega(g^{-1})F_n(w)
\end{align}
and therefore $\phi$ acts on $\Hom_{{}^0G}(\rho,\pi \otimes \omega)$ by multiplication by
$\omega(g^{-1})$. In particular, $\phi$ acts by multiplication by a scalar on any simple right
$\caA$-module. We deduce that $\phi \in \caZ$. Indeed, the image of $\phi$ under the injection of Lemma
\ref{strucMA} is central. We know that $[N:T]=m^2$. Let us choose a system of representatives
$n_1,\ldots ,n_{m^2}$ for the right cosets of $T$ in $N$ and for each $i$, $1\leq i\leq m^2$, a non-zero
function $\phi_i$ in $\caA=\caH(N,\rho)$ with support ${}^0G n_i$. We easily verify using the support
properties (\ref{Gxy}) that $(\phi_i)_{1\leq i\leq m^2}$ is a basis of $\caA$ over its center
$\caH(T,\rho)$.
\end{proof}

We can deduce another important consequence from the above results:

\begin{cor}
The $\bbC$-algebra $\caA$ is Noetherian.
\end{cor}

\begin{proof}
Indeed, $\caZ$ is isomorphic to a Laurent polynomial algebra in several variables, hence Noetherian
(\cite{Douady}, 3.6.2). On the other hand, $\caA$ is a finitely generated module over $\caZ$. We conclude
by using (\cite{Douady}, 3.6.1). The category $\caM(G)_{[\pi]}$ is therefore equivalent to the category
of right modules over a Noetherian algebra.
\end{proof}

\needspace{126pt}
\begin{prop}
For all $g \in T$, let $\psi_g$ denote the element of $\caZ$ whose support is ${}^0Gg^{-1}$ normalized as
above, i.e., $\psi_g(g^{-1})=\rho_T(g^{-1})$. Then the map $g \mapsto \psi_g$ induces an isomorphism
$\bbC[T/{}^0G]\simeq \caZ$. This shows that $\caZ$ is isomorphic to the algebra of functions of the
variety $\mathbf{Irr}(G)_{[\pi]}$.
\end{prop}

\begin{proof}
The fact that $g \mapsto \psi_g$ induces an isomorphism $\bbC[T/{}^0G]\simeq \caZ$ follows easily from
(\ref{laurentpol}). We noted in \ref{T0GIrr} that the algebra of functions on $\mathbf{Irr}(G)_{[\pi]}$
is $\bbC[T/{}^0G]$. The last assertion follows from this.
\end{proof}

To conclude this section, we summarize the results obtained on the center of $\caM(G)_{[\pi]}$. The
equivalence of categories between $\caM(G)_{[\pi]}$ and $\caM(\caA)_d$ shows that the center of
$\caM(G)_{[\pi]}$ is isomorphic to $\caZ$, the center of $\caA$. Furthermore, $\caZ$ is isomorphic to the
algebra of polynomial functions on the variety $\mathbf{Irr}(G)_{[\pi]}$. We thus obtain a description of
the center of $\caM(G)_{[\pi]}$ as the algebra of polynomial functions on $\mathbf{Irr}(G)_{[\pi]}$. The
above isomorphisms are not canonical: $\mathbf{Irr}(G)_{[\pi]}$ is described as a homogeneous space for
$\caX(G)$, and this depends on the choice of a base point (here $\pi$) in the inertial class $[\pi]$. On
the other hand, we chose a particular progenerator of $\caM(G)_{[\pi]}$, on which the definition of the
algebra $\caR$ depends, and therefore that of $\caA$. In the next section, we will see how to make this
description independent of the choices.

\subsection{The center of $\caM(G)_{[\pi]}$: conclusion}\label{centreMpiconcl}

Let $M$ be a simple right $\caA$-module. The center $\caZ$ of $\caA$ acts by a character
$\lambda_M \colon \caZ \rightarrow \bbC$. Let us set $I_M=\ker \lambda_M$, so that the two-sided ideal
$I_M\caA$ is contained in the annihilator of $M$. The structure morphism
$$\phi_M \colon \caA \rightarrow \End_\bbC (M)$$ induces an algebra morphism, surjective according to
Burnside's theorem (\cite{Lang}, Corollary XVII.3.3):
\[
\tilde  \phi_M \colon  \caA / I_M\caA  \rightarrow \End_\bbC (M).
\]
Note that $\caA / I_M\caA$ is of dimension $m^2$ over $\bbC$ because $\caA$ is free over $\caZ$ of rank
$m^2$, and $\caZ / I_M \simeq \bbC$. Since $\End_\bbC (M)$ is also of dimension $m^2$, we deduce that
$\tilde  \phi_M$ is an algebra isomorphism. Moreover $\End_\bbC (M)$ is isomorphic as a right
$\End_\bbC (M)$-module to a sum of $m$ copies of $M$. We deduce that the same is true for
$\caA / I_M\caA$. It follows that $M$ is entirely determined by the action of $\caZ$, or more precisely
by the ideal $I_M$ of $\caZ$. Conversely, for any maximal ideal $I$ of $\caZ$, the finitely generated
right $\caA$-module $\caA / I\caA$ admits a simple quotient $M$, $I=I_M$ and $\caA / I\caA \simeq m M$.
In summary, if $M_I$ denotes the unique simple composition factor of the semisimple right $\caA$-module
$\caA / I\caA$, then
\[
I \mapsto M_I,\quad \mathrm{SpecMax}(\caZ) \rightarrow \mathbf{Irr}(\caA).
\]
defines a bijection between the set of maximal ideals of $\caZ$ and the set of isomorphism classes of
simple right $\caA$-modules. This equips $\mathbf{Irr}(\caA)$ with the structure of an affine algebraic
variety over $\bbC$, with function algebra $\caZ$.

\bigskip

The equivalence of categories between $\caM(\caA)_d$ and $\caM(G)_{[\pi]}$ therefore induces an algebra
isomorphism between $\caZ$ and the center $\caZ_{[\pi]}$ of the category $\caM(G)_{[\pi]}$, denoted by
$\phi \mapsto z_\phi$. This equivalence also induces a bijection between $\mathbf{Irr}(\caA)$ and
$\mathbf{Irr}(G)_{[\pi]}$.

For all $I \in \mathrm{SpecMax}(\caZ)$, the module $\caA / I\caA$ in $\caM(\caA)_d$ corresponds to the
module
\[
(\caA / I\caA )\otimes_\caA \Pi \simeq \Pi/I\Pi.
\]
There exists a unique irreducible representation (up to isomorphism) $\pi_I$ in $\caM(G)_{[\pi]}$ such
that $\Pi/I\Pi \simeq m \pi_I$. Let $\lambda_I$ denote the algebra morphism from $\caZ$ to $\bbC$ such
that $\ker\lambda_I=I$. Let $\phi \in \caZ$. Then $\phi$ acts on $M_I$ by multiplication by
$\lambda_I(\phi)$. We deduce that $z_\phi$ acts on $\pi_I$ by this same scalar $\lambda_I(\phi)$.

\begin{prop}
Let $z \in \caZ_{[\pi]}$ and let $\tau \in \mathbf{Irr} (G)_{[\pi]}$. Let us interpret $z$ as an element
of the algebra of functions of the affine algebraic variety $\mathbf{Irr}(G)_{[\pi]}$. Then the
evaluation of $z$ at the point $\tau$ is the scalar $z(\tau)$ by which $z$ acts on $\tau$.
\end{prop}

\begin{proof}
Let $\omega \in \caX(G)$ and $x \in T$. It suffices to prove the assertion in the case where
$z=z_{\psi_x}$, $\psi_x$ being the element of $\caZ$ defined in Proposition \ref{centreA}. According to
(\ref{lop}), the function $\psi_x$ acts on the simple right module
$\Hom_{{}^0G}(\rho,\pi \otimes \omega)$ by the scalar $\omega(x)$. Thus, the element $z_{\psi_x}$ acts by
the scalar $\omega(x)$ on $\pi \otimes \omega$.
\end{proof}

We note another consequence of all this:

\begin{cor}
The category $\caM(G)_{[\pi]}$ is indecomposable.
\end{cor}

\begin{proof}
When an abelian category splits into a direct sum of two categories, the projections onto each of the
factors are idempotents in the center of the category. Since $\caZ$ is an integral domain (because
$\caZ\subset\caA$ and $\caA$ is an integral domain, Corollary \ref{strucMA}), it contains no non-trivial
idempotent element.
\end{proof}

\section{Induced representations}

We fix in all that follows a minimal parabolic subgroup $P_\emptyset=M_\emptyset N_\emptyset$ of $G$. We
adopt the notation of Chapter \ref{chapstruct}, in particular
\[
W_G=W(A_\emptyset)=N_G(A_\emptyset)/Z_G(A_\emptyset)=N_G(A_\emptyset)/M_\emptyset.
\]

\subsection{The geometric lemma}\label{geomlemma}
\index[ter]{geometric lemma}

In this section, we show the following result:
\begin{thm}
Let $P=MN$ and $Q=LU$ be parabolic subgroups of $G$ and let $(\tau,E)$ be a smooth representation of $M$.
The representation
\[
r_Q^G \, i_P^G (\tau,E)
\]
of $L$ admits a filtration whose associated graded components are isomorphic to
\[
(i_{L\cap w^{-1}\cdot P}^{L} \circ \; w \circ  r_{w\cdot Q\cap M}^M )(\tau,E),
\]
where $w$ runs through the subset $\caW^{Q,P}$ of $G$ defined in \ref{PWQsemi}.
\end{thm}

\begin{rmqs}
1. The set $\caW^{Q,P}$ is a system of representatives in $G$ of the double cosets $P \backslash G/Q$.
According to the properties of $\caW^{Q,P}$ established in \ref{PWQsemi}, $L\cap w^{-1}\cdot P$ is a
parabolic subgroup of $L$ with Levi decomposition
\[
L\cap w^{-1}\cdot P=(L\cap w^{-1}\cdot M)(L\cap w^{-1}\cdot N),
\]
and $w\cdot Q\cap M$ is a parabolic subgroup of $M$ with Levi decomposition
\[
w\cdot Q\cap M=(w\cdot L\cap M)(w\cdot U\cap M) .
\]

\noindent --- 2. We have simply denoted by $w$ the forgetful functor
\[
\sigma \mapsto  {}^w \sigma,\quad      \caM(w\cdot L \
\cap M) \rightarrow \caM(L\cap w^{-1}\cdot M)
\]
associated to the isomorphism $L\cap w^{-1}\cdot M \rightarrow w\cdot L \cap M $ induced by conjugation
by $w$ ({\sl cf.} Example \ref{FOA}). We draw the reader's attention to the fact that many authors prefer
to denote it $w^{-1}$.
\end{rmqs}

We will first formulate this result more precisely. Consider the totally disconnected quotient space
$X=P\backslash G$ equipped with the right translation action of $G$:
\[
G\times X \rightarrow X,\quad (g,Ph)\mapsto Phg^{-1}.
\]
The parabolic subgroup $Q$ acts on $X$ with a finite number of orbits $Z_1,\ldots ,Z_k$ (Corollary
\ref{WMWGWM}). By Proposition \ref{actiontdprop}, these orbits are all locally closed in $X$. Moreover,
we can assume that this numbering of the orbits is such that the sets
\[
Y_1=Z_1, \; Y_2=Z_1 \cup Z_2,\; \ldots , Y_k=Z_1 \cup Z_2 \cup
\ldots  \cup Z_k=X
\]
are open subsets of $X$. We will assume this is the case in what follows.

Fix a $Q$-orbit $Z$ in $X$, and a point $Pz$ of this orbit, with $z \in \caW^{Q,P}$. For brevity, let
$\iota$ denote the automorphism $\Int(z)$ of $G$. Set
\[
M'=M \cap \iota(L), \quad L'=L \cap \iota^{-1}(M), \quad N'=M \cap
\iota(U),\quad  U'=L \cap \iota^{-1}(N).
\]
Then $P'=M'N'$ is a parabolic subgroup of $M$, $Q'=L'U'$ is a parabolic subgroup of $L$ and
$M' =\iota(L')$. Also set $U''=\iota^{-1}(N')$, so that $Q''=L'U''$ is a parabolic subgroup of
$\iota^{-1}(M)$. We have at our disposal the functors
\[
r_{P'}^M\colon \caM(M) \rightarrow \caM(M'), \quad i_{Q'}^L\colon \caM(L') \rightarrow \caM(L).
\]

Then define the functor $\Phi_Z$ by
\[
\Phi_Z \colon \caM(M) \rightarrow \caM(L), \quad \Phi_Z=i_{Q'}^L \circ \iota  \circ  r_{P'}^M.
\]
Here again we denote by $\iota$ the forgetful functor associated to the isomorphism
$\iota \colon L' \rightarrow M'$, i.e., the functor which to any representation $(\sigma,E)$ of
$\caM(M')$ associates the representation $({}^ \iota \sigma ,E)$ of $L'$ with the notation of Example
\ref{FOA}.

We can now reformulate the theorem in the following (more precise) way. The notation is the same as
above.
\begin{thm}
(reformulation) The functor
\[
F=r_Q^G \circ  i_P^G  \colon \caM(M)   \rightarrow \caM(L)
\]
admits a filtration by subfunctors
\[
0=F_0 \subset F_1 \subset \ldots \subset F_k= F
\]
such that $F_i/F_{i-1} \simeq \Phi_{Z_i}$, $i=1,\ldots ,k$.
\end{thm}

Recall that a subfunctor of a functor $F$ from a category $\scrC$ taking values in a category of modules
$\scrD$ is a functor $G$ from $\scrC$ to $\scrD$ such that for any object $X$ of $\scrC$, $G(X)$ is a
submodule of $F(X)$.

\begin{proof}
Let $Y$ be a $Q$-invariant open subset of $X$. For any smooth representation $(\sigma,E)$ of $M$, let us
realize as in \ref{Ind} the space of the induced representation $i_P^G(\sigma,E)$ as a space of functions
on $G$, and consider the subspace $i_Y(E)$ of $i_P^G(E)$ consisting of functions $f$ supported in $PY$
where
\[
PY=\{ g \in G \mid Pg \in Y  \}.
\]
Note in passing that $PY$ is an open subset of $G$, as the inverse image of an open subset by a
continuous map, in this case the canonical projection from $G$ to $X$. It is clear that $i_Y(E)$ is
stable under the action of $Q$ obtained on $i_P^G(E)$ by restriction of the action of $G$. We can
therefore apply to $i_Y(E)$ the functor $r_U=j_U \otimes \delta_{U}^{1/2}$, and form
\[
F_Y(E)=r_U(i_Y(E)).
\]
Since the functor $r_U$ is exact, $r_U(i_Y(E))$ is a subrepresentation of $F(E)=r_U (i_P^G(E))$. In other
words, $F_Y$ is a subfunctor of $F$.

\begin{prop}
Let $Y$ and $Y'$ be two $Q$-invariant open subsets of $X$. We then have, with the obvious notation
\[
F_{Y \cap Y'}=F_{Y}\cap F_{Y'},\quad  F_{Y \cup Y'}=F_{Y} +
F_{Y'},\quad  F_\emptyset =0 ,\quad F_X=F.
\]
\end{prop}
\begin{proof}
Since by Proposition \ref{norm}, the functor $r_U$ is exact, it suffices to show that
\[
i_{Y \cap Y'}=i_{Y}\cap i_{Y'},\quad  i_{Y \cup Y'}=i_{Y} +
i_{Y'},\quad  i_\emptyset =0 ,\quad i_X=i_P^G.
\]
The last two equalities are trivial. Moreover, we have
\[
P(Y\cup Y')=PY \cup PY', \quad  P(Y\cap Y')=PY \cap PY'.
\]

This shows that $i_{Y \cap Y'} = i_{Y}\cap i_{Y'}$ and $i_Y+i_{Y'}\subset i_{Y\cup Y'}$. It remains to
show that $  i_{Y\cup Y'}\subset  i_Y+i_{Y'}$.

Let $f$ be in $i_{Y \cup Y'}$. Set $\bar W=p(\mathrm{Supp}(f))\subset P\backslash G$, where
$\mathrm{Supp}(f)$ is of course the support of $f$. We thus have $\bar W \subset Y\cup Y'$. Now there
exist two functions $g_1$ and $g_2$ in $\scrC^{\infty}(Y)$ and $\scrC^{\infty}(Y')$ respectively such
that $g_1+g_2\equiv 1$ on $\bar W$. We then have $f=f(g_1 \circ p)+f(g_2 \circ p)$, with
$f(g_1 \circ p)\in i_Y$ and $f(g_2 \circ p)\in i_{Y'}$. Let us construct these functions $g_1$ and $g_2$.
Cover $\bar W$ by open sets $(U_x)_{x\in \bar W}$ such that $x \in U_x \subset Y$ if $x \in Y$ and
$x \in U_x \subset Y'$ if $x \in Y'$. Extract from this cover of $\bar W$ by the $(U_x)_{x\in \bar W}$ a
finite subordinate cover by disjoint compact open sets $(U_i)$ (see Lemma \ref{recouv}). Take for $g_1$
the sum of the characteristic functions of the open sets $U_i$ contained in $Y$, and for $g_2$ the sum of
the characteristic functions of the other $U_i$, which are contained in $Y'$. These functions then have
the required properties.
\end{proof}

\bigskip

By this proposition, we can extend the definitions of $i_Y$ and of the functor $F_Y$ to the case where
$Y$ is only a locally closed $Q$-invariant subset of $X$ as follows: write
\[
Y=T \cap V
\]
where $T$ is closed in $X$ and $V$ is open. By replacing $T$ with $\bigcap_{g \in Q}g\cdot T$ and $V$
with $\bigcup_{g \in Q} g\cdot V$, we may assume $T$ and $V$ invariant under the action of $Q$. Set
$Y_1=(X\setminus T) \cap V$. We have $Y\cap Y_1=\emptyset$, $Y\cup Y_1$ is open in $X$, and $Y_1$ is
$Q$-invariant. Then set $i_Y= i_{Y\cup Y_1}/i_{Y_1}$ and $F_Y= F_{Y\cup Y_1}/F_{Y_1}$. The proposition
shows that the functors thus constructed for different choices of $Y_1$ are naturally isomorphic.

On the other hand, the functor $F_Y$ can be described as the composition
\begin{equation}\label{FY}
 F_Y= r_U \circ  i_Y= r_U \circ  \Gamma_c   \circ \caR_{Y,Q}^{X,G}   \circ
 \caI_P^G\circ \delta_N^{-1/2}  \circ \caF_P^M
\end{equation}
where
\[
\caI_P^G \colon  \caM(P)\longrightarrow \scrC^\infty_{X,G}-\caM od,
\]
\[
\caR_{Y,Q}^{X,G}  \colon  \scrC^\infty_{X,G}-\caM od \longrightarrow \scrC^\infty_{Y,Q}-\caM od
\]
\[
\Gamma_c  \colon  \scrC^\infty_{Y,Q}-\caM od  \longrightarrow \caM(Q)
\]
are the functors defined in \ref{indetfaisc}. Indeed, that $F_Y$ is equal to this composition is clear if
$Y$ is a $Q$-invariant open subset of $X$, and the case of a locally closed subset then follows from the
above proposition and the second proposition of \ref{FCev}.

With the notation of the theorem, we now see that the functor $F$ admits a filtration
\[
0=F_0\subset F_{Y_1} \subset F_{Y_2} \subset \ldots \subset
F_{Y_k}=F
\]
and that for all $i=1,\ldots,k$, $F_i/F_{i-1} =F_{Z_i}$. It therefore suffices to prove that for any
$Q$-orbit $Z$ in $X$, we have a natural isomorphism $\Phi_{Z}\simeq F_{Z}$.

Set $\hat P=\iota^{-1}(P)$, $\hat M=\iota^{-1}(M)$, $\hat N=\iota^{-1}(N)$. Let us express the functors
$\Phi_Z$ and $F_Z$ in terms of similar functors defined by replacing $P,M,N$ by $\hat P$, $\hat M$,
$\hat N$. The automorphism $\iota$ of $G$ induces equivalences of categories (which we will all still
denote by $\iota$)
\[
\iota \colon \caM(G)\rightarrow \caM(G),\quad    \caM(P)\rightarrow
\caM(\hat P), \quad   \caM(M)\rightarrow
\caM(\hat M), \quad  \caM(M')\rightarrow \caM(L').
\]
It is clear that we have
\[
\iota ^{-1}\circ r_{Q''}^{\hat M} \circ \iota =  r_{P'}^M,
\]
whence,
\[
\Phi_Z=i_{Q'}^L \circ \iota \circ  r_{P'}^M =
i_{Q'}^L \circ \iota \circ  \iota ^{-1}\circ r_{Q''}^{\hat M} \circ \iota  = i_{Q'}^L \circ r_{Q''}^{\hat M} \circ \iota.
\]
Set $\hat Z= \hat P Q \in \hat P\backslash G /Q$ and define the functor
$\Phi_{\hat Z}= i_{Q'}^L \circ r_{Q''}^{\hat M} \colon \caM(M) \rightarrow \caM(L)$. We then have
\[
\Phi_Z = \Phi_{\hat Z}\circ  \iota.
\]
Similarly, $\iota$ induces a homeomorphism
\[
\iota \colon \hat P \backslash G \rightarrow P\backslash G ,\quad \hat P g
\mapsto P \cdot zg.
\]
Let us still call $\iota$ the inverse image functor
$$\iota \colon   \scrC^\infty_{P\backslash G}-\caM od \longrightarrow \scrC^\infty_{\hat P \backslash G}-\caM od$$
for this homeomorphism. Note that the orbit $PzQ$ has as its inverse image the subset $\hat PQ $ of
$\hat P \backslash G$.

It is clear that with this notation, we have $\iota\circ \caI_P^G \circ \iota ^{-1} = \caI_{\hat P}^G$,
whence
\[
F_Z= r_U \circ  \Gamma_c   \circ \caR_{Z,Q}^{X,G}    \circ  \iota^{-1}
\circ \caI_{\hat P}^G \circ \iota  \circ \delta_N^{-1/2}   \circ \caF_P^M.
\]
On the other hand, we also obviously have $\iota \circ \caF_P^M= \caF_{\hat P}^{\hat M}  \circ  \iota$
and $\iota\circ \delta_N^{-1/2}=\delta_{\hat N}^{-1/2}\circ \iota $ whence finally:
\[
F_Z= r_U \circ  \Gamma_c \circ  \caR_{Z,Q}^{X,G}    \circ  \iota^{-1}\circ \caI_{\hat P}^G
\circ  \delta_{\hat N}^{-1/2} \circ \caF_{\hat P}^{\hat M}  \circ  \iota.
\]
Denote
$F_{\hat Z}= r_U \circ  \Gamma_c \circ  \caR_{Z,Q}^{X,G}    \circ  \iota^{-1}\circ \caI_{\hat P}^G \circ   \delta_{\hat N}^{-1/2} \circ \caF_{\hat P}^{\hat M}$,
so that $F_Z=F_{\hat Z}\circ \iota$. Also set
$\caI_{\hat Z}=   \Gamma_c \circ  \caR_{Z,Q}^{X,G}    \circ  \iota^{-1}\circ \caI_{\hat P}^G $, so that
$F_{\hat Z}=  r_U \circ \caI_{\hat Z}  \circ  \delta_{\hat N}^{-1/2} \circ \caF_{\hat P}^{\hat M}$. We
have thus reduced the problem to showing that $F_{\hat Z}$ is naturally isomorphic to $\Phi_{\hat Z}$. To
do this, we will analyze the following diagram of functors.

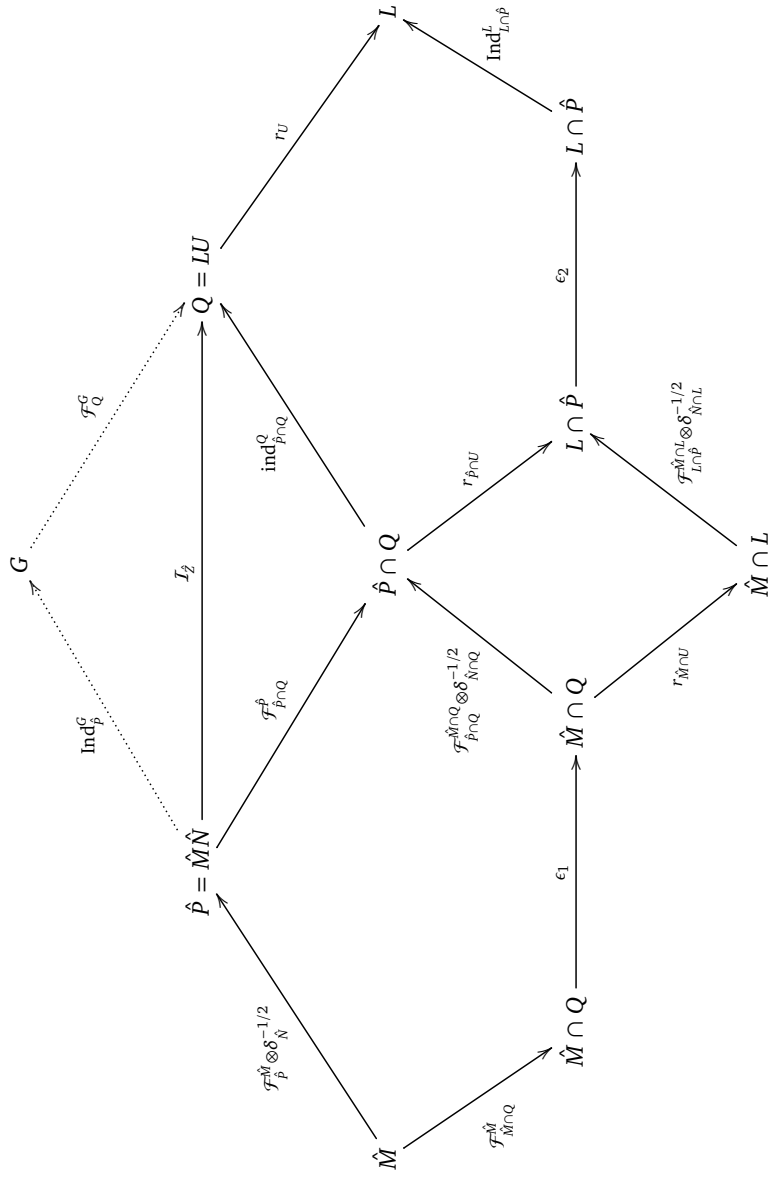
\vfill
\pagebreak
\begin{landscape}
\[
\xymatrix{
& &  & &  G    \ar@{.>}[ddrr]^{\caF_Q^G}  &  & & & \\
&& && && && \\
& &  \hat P=\hat M \hat N \ar@{.>}[uurr]^{\Ind_{\hat P}^G}
\ar[ddrr]^{\caF_{\hat P\cap
Q}^{\hat P}}\ar[rrrr]^{\caI_{\hat Z}}& & &   &   Q=LU
\ar[ddrr]^{r_U} & & \\
&& && && && \\
\hat M \ar[uurr]^{\caF_{\hat P}^{\hat M}\otimes \delta_{\hat N}^{-1/2}  }  \ar[ddr]_{\caF_{\hat M\cap Q}^{\hat
M}} &  & &   & \hat P\cap
Q\ar[uurr]^{\ind_{\hat P\cap Q}^Q} \ar[ddr]^{r_{\hat P \cap U}}& & & & L
\\
&& && && && \\
& \hat M\cap Q \ar[rr]^{\epsilon_1}& &\hat M\cap Q
\ar[uur]^{\caF_{\hat P\cap Q}^{\hat M\cap Q} \otimes \delta_{\hat
N\cap Q}^{-1/2}}
\ar[ddr]_{r_{\hat M\cap U}} &
&  L\cap \hat P \ar[rr]^{\epsilon_2}& &     L\cap \hat P
\ar[uur]_{\Ind_{L\cap \hat P}^L} & \\
&& && && &&\\
& & & & \hat M\cap L  \ar[uur]_{\caF_{L\cap \hat P}^{\hat M\cap L}
\otimes \delta_{\hat N \cap L}^{-1/2}} & &
& &\\
}
\]
\end{landscape}
\pagebreak

The source and target categories of the functors are represented by the vertices of this graph. The
vertex subscripted by a group $H$ represents the category of smooth representations $\caM(H)$. The arrows
represent the functors. The dotted arrows do not directly intervene in what follows and should be
considered decorative. The functor $F_{\hat Z}$ is obtained by following the ``top'' path leading from
$\hat M$ to $L$, and the functor $\Phi_{\hat Z}$ is obtained by following the ``bottom'' path, except
that we have introduced a twist by the characters
$\epsilon_1= \delta_{\hat N}^{-1/2}\delta^{1/2}_{\hat N \cap Q}$ of $\hat M \cap Q$ and
$\epsilon_2= \delta_{U}^{-1/2}\delta^{1/2}_{U \cap \hat P}$ of $\hat P \cap L$. To show that these
functors are naturally isomorphic, we will show that each small subdiagram with four (or three for the
top one) vertices is commutative up to natural isomorphism, then that the effects of $\epsilon_1$ and
$\epsilon_2$ cancel each other out. This is obviously sufficient. We will explain the notation not yet
introduced in this diagram during the proof.

\medskip

Consider the diagram
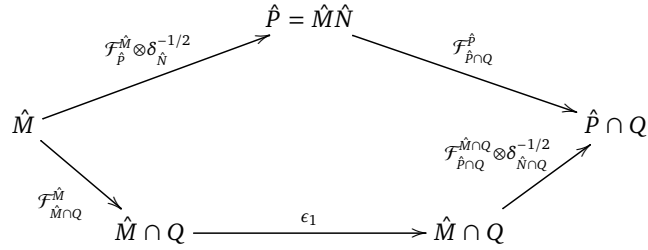
\[
\xymatrix{
& &  \hat P=\hat M \hat N \ar[drr]^{\caF_{\hat P \cap Q}^{\hat P}} & &
\\
\hat M \ar[urr]^{\caF_{\hat P}^{\hat M}\otimes \delta_{\hat N}^{-1/2}  }  \ar[dr]_{\caF_{\hat M\cap Q}^{\hat
M}} &  & &   & \hat P\cap Q \\
& \hat M\cap Q \ar[rr]^{\epsilon_1}& &\hat M\cap Q
\ar[ur]^{\caF_{\hat P\cap Q}^{\hat M\cap Q} \otimes \delta_{\hat
N\cap Q}^{-1/2}} &\\
}
\]
We have, using the properties of the chosen representatives ({\sl cf.} \ref{PWQsemi}),
\[
\hat P\cap Q=(\hat M \cap L)(\hat N \cap L)(\hat M \cap U)(\hat
N\cap U).
\]
Since the groups $\hat N$, $U$ are unions of their compact open subgroups, all the modular characters
intervening in this diagram are trivial on these groups. It therefore suffices to see how an element $l$
of $\hat M\cap L$ acts on the representations obtained from a smooth representation $(\sigma,E)$ of
$\hat M$, by following the top path versus the bottom path. For the first, $l$ acts by
$\delta_{\hat N}^{-1/2}(l) \sigma(l)$ and for the second by
$\epsilon_1(l)\delta_{\hat N \cap Q}^{-1/2}(l) \sigma(l)$. But
\[
\epsilon_1 \delta_{\hat N \cap Q}^{-1/2}=\delta_{\hat N}^{-1/2}
\delta_{\hat N \cap Q}^{1/2} \delta_{\hat N \cap Q}^{-1/2}= \delta_{\hat N}^{-1/2},
\]
which proves the commutativity of the diagram.

\medskip

Now consider the diagram
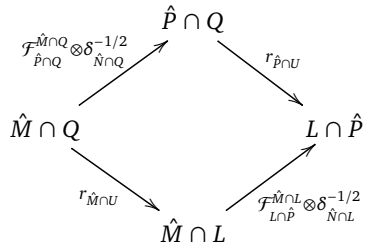
\[
\xymatrix{
&  \hat P\cap Q   \ar[dr]^{r_{\hat P\cap U}}&    \\
\hat M\cap Q  \ar[ur]^{\caF^{\hat M\cap Q}_{\hat P\cap Q} \otimes \delta_{\hat
N\cap Q}^{-1/2}}
\ar[dr]_{r_{\hat M\cap U}} & & L\cap \hat P  \\
& \hat M\cap L  \ar[ur]_{\caF_{L\cap \hat P}^{ \hat M\cap L} \otimes \delta_{\hat
N\cap L}^{-1/2}}  &
}
\]
The functor $r_{\hat P\cap U}$ is the functor $j_{\hat P \cap U}$ twisted by the character
$\delta^{1/2}_{\hat P \cap U}$ and the functor $r_{\hat M\cap U}$ is the functor $j_{\hat M \cap U}$
twisted by the character $\delta^{1/2}_{\hat M \cap U}$. Let $(\sigma,E)$ be a smooth representation of
$\hat M\cap Q$. Let us first verify that the representations $\pi_1$ and $\pi_2$ of $L\cap \hat P$
obtained respectively by following the top path and the bottom path in the diagram indeed act in the same
space. The space of $\pi_1$ is $E$ quotiented by the space generated by the vectors of the form
\[
\sigma(u) \delta_{\hat  N\cap Q}^{-1/2}(u)\cdot v-v, \, (v \in E),\, (u \in \hat P\cap U).
\]
But since $\hat P\cap U=(\hat M\cap U)(\hat N\cap U)$ and $\hat N\cap U$ acts trivially, this space is in
fact the space generated by the vectors of the form
\[
\sigma(u) \delta_{\hat  N\cap Q}^{-1/2}(u)\cdot v-v, \, (v \in E),\, (u \in \hat M\cap U).
\]
On the other hand, $\hat M \cap U$ being a union of its compact open subgroups,
$ \delta_{\hat  N\cap Q}^{-1/2}(u)=1$ for all $u \in \hat M\cap U$. Finally, the space of $\pi_1$ is the
quotient of $E$ by the space generated by the vectors of the form
\[
\sigma(u) \cdot v-v, \, (v \in E),\, (u \in \hat M\cap U),
\]
and this is also the space of $\pi_2$. We must now calculate the character twisting the action induced by
$\sigma$ on this quotient in both cases. For $\pi_1$, it is
$\delta_{\hat N \cap Q}^{-1/2} \delta_{\hat P \cap U}^{1/2} $ and for $\pi_2$ it is
$\delta_{\hat N \cap L}^{-1/2} \delta_{\hat M \cap U}^{1/2} $. Since
$\delta_{\hat N \cap Q}=\delta_{\hat N \cap U}\delta_{\hat N \cap L}$ and
$\delta_{\hat P \cap U}^{1/2} =   \delta_{\hat M \cap U}^{1/2} \delta_{\hat N \cap U}^{1/2}$, we see that
we indeed have $\pi_1=\pi_2$.

\medskip

We now turn to

\begin{equation}\label{diagrammedifficile}
\xymatrix{
&   &   Q=LU \ar[drr]^{r_U} & & \\
\hat P\cap Q\ar[urr]^{\ind_{\hat P\cap Q}^Q} \ar[dr]^{r_{\hat P \cap
    U}}& & & & L \\
  &  L\cap \hat P \ar[rr]^{\epsilon_2}& &     L\cap \hat P
 \ar[ur]_{\Ind_{L\cap \hat P}^L} & \\
}
\end{equation}

We will place ourselves in a slightly more general context, which by specialization will imply the
commutativity of a version without the twists of the diagram above. We locally adopt new notation:

Let $J$ be a totally disconnected group, $H$ a closed subgroup and $N$ a closed normal subgroup which is
an increasing directed union of its compact open subgroups (so in particular $N$ is unimodular) and
suppose that $HN$ is closed. Let $(\sigma,E)$ be a smooth representation of $H$. Consider the
representation
\[
(j_{H\cap N},E_{H \cap N}) \quad \text{of} \quad L'=H/H\cap N\simeq
HN/N\subset J/N .
\]
We claim that we then have a natural isomorphism
\begin{equation}\label{JHNL}
  j_N \circ \ind_H^J(\sigma,E)  \simeq
  \ind_{HN/N}^{J/N}(j_{H\cap N}(\sigma,E)\otimes \delta^{-1}),
\end{equation}
where $\delta$ is the modular character of $H$ on $H\cap N\backslash N$.

We first give a proof of (\ref{JHNL}) that does not yield the explicit form of the isomorphism. For
$h\in H$, denote by $\delta_J(h)$, $\delta_H(h)$, $\delta_N(h)$, $\delta_{H\cap N}(h)$, $\delta_{HN}(h)$
the moduli of the automorphism $\Int(h)$ of these groups (\ref{foncmod}), so that
$\delta=\delta_N/\delta_{H\cap N}$ and $\delta_{H\backslash J}=\delta_J/\delta_H$ (\ref{invmesquo}). By
the theorem of \ref{indcomp} and the proposition of \ref{Jacmod}, we have natural isomorphisms
\[
\ind_H^J(\sigma,E)\simeq\caH(J)\otimes_{\caH(H)}(E\otimes\delta_{H\backslash J}^{-1}),\qquad j_N(V)\simeq\caH(J/N)\otimes_{\caH(J)}V,
\]
hence, by associativity of the tensor product,
\[
j_N\circ \ind_H^J(\sigma,E)\simeq\caH(J/N)\otimes_{\caH(H)}(E\otimes\delta_{H\backslash J}^{-1}).
\]
On the other hand, since $HN$ is closed, $HN/N$ is a closed subgroup of $J/N$, isomorphic to $H/H\cap N$,
and the same results give
\begin{align*}
\ind_{HN/N}^{J/N}(j_{H\cap N}(E)\otimes\delta^{-1})&\simeq\caH(J/N)\otimes_{\caH(HN/N)}\big(\caH(HN/N)\otimes_{\caH(H)}(E\otimes\delta^{-1}\delta_{HN\backslash J}^{-1})\big)\\ &\simeq\caH(J/N)\otimes_{\caH(H)}(E\otimes\delta^{-1}\delta_{HN\backslash J}^{-1}),
\end{align*}
where we have used that the characters $\delta$ and $\delta_{HN\backslash J}$ are trivial on $H\cap N$,
which is a union of its compact open subgroups, and that $\delta_{(HN/N)\backslash(J/N)}$ identifies with
$\delta_{HN\backslash J}$, the moduli being multiplicative with respect to the normal subgroup $N$:
$\delta_J=\delta_{J/N}\delta_N$ and $\delta_{HN}=\delta_{HN/N}\delta_N$. For the same reason,
$\delta_{HN}=\delta_{H/H\cap N}\,\delta_N=\delta_H\delta_N/\delta_{H\cap N}$ on $H$, hence
\[
\delta_{H\backslash J}=\frac{\delta_J}{\delta_{HN}}\cdot\frac{\delta_{HN}}{\delta_H}=\delta_{HN\backslash J}\;\delta,
\]
which proves (\ref{JHNL}).

We will now prove (\ref{JHNL}) by giving the explicit form of the natural isomorphism. For all $v \in E$,
denote by $\bar v$ its image in $j_{H\cap N}(E)$. Since any function $f \in  \ind_H^J(E)$ is compactly
supported modulo $H$, for all $j \in J$, the function $ (r(j)\cdot f)_{|N}$ is compactly supported modulo
$H\cap N$. Moreover, for any $h \in H \cap N$, $n \in N$ and $j \in J$, we have
\[
\overline{f(hnj)}=\overline{\sigma(h)\cdot f(nj)}= \overline{f(nj)},
\]
and thus $\overline{(r(j)\cdot f)_{|N}}\in \scrD(N,H\cap N,\delta_{H \cap  N\backslash N}=1)$. Thus,
having fixed an invariant measure on $H\cap N\backslash N$ in the sense of Section \ref{invmesquo}, the
integral
\[
\bar f(j):= \int_{H\cap N\backslash N} \overline{f(nj)}\;  d\nu_{H\cap N\backslash N}(n)
\]
is well defined.

We have for all $u \in N$, for all $j \in J$, using the right invariance of the measure
$\nu_{H\cap    N\backslash N}$,
\begin{small}
\begin{equation*}
\bar f(uj)= \int_{H\cap N\backslash N} \overline{f(nuj)}\; d\nu_{H\cap N\backslash N}(n) =\int_{H\cap
N\backslash N} \overline{f(nj)} \; d\nu_{H\cap N\backslash N}(n)= \bar f(j) .
\end{equation*}
\end{small}
We can thus consider $\bar f$ as a function on $N\backslash J=J/N$ ($N$ is normal).

For all $h \in H$, for all $j \in J$,
\begin{small}
\begin{align*}
\bar f(hj)&= \int_{H\cap N\backslash N} \overline{f(nhj)}\; d\nu_{H\cap N\backslash N}(n)= \int_{H\cap
N\backslash N} \overline{f(h(h^{-1}nh)j)} \; d\nu_{H\cap N\backslash N}(n) \\
&= \int_{H\cap N\backslash N} \overline{ \delta^{-1}(h) \sigma(h) \cdot f(nj)} \; d\nu_{H\cap N\backslash
N}(n)= \delta^{-1}(h) \sigma_{H\cap N}(h) \cdot \bar f(j).
\end{align*}
\end{small}

The function $\bar f$ is clearly compactly supported modulo $HN$ in $J$. It follows from the above that
$f \mapsto \bar f$ defines a map from $\ind_H^J(\sigma,E)$ to
$\ind_{HN/N}^{J/N}(j_{H\cap N}(\sigma,E)\otimes \delta^{-1})$. It is clearly a $J$-morphism, and it
factors into a $J/N$-morphism
\begin{equation}\label{1P}
 j_N (\ind_H^J(\sigma,E)) \rightarrow  \ind_{HN/N}^{J/N}(j_{H\cap
   N}(\sigma,E)\otimes \delta^{-1}) \quad j_N(f) \mapsto \bar f.
\end{equation}

Let us now show the surjectivity of $f \mapsto \bar f$. Recall that
\[
\ind_{HN/N}^{J/N}(j_{H\cap N}(E)\otimes \delta^{-1})
\]
is the set of functions
\[
\phi\colon J/N \rightarrow j_{H\cap N}(E)
\]
smooth with compact support modulo $HN/N$ in $J/N$ such that, if $\bar x$ denotes the image of $x \in J$
in $J/N$,
\[
\phi(\bar h \bar n\bar
j)= \delta^{-1}(h)\sigma_{H\cap N}(\bar h \bar n)\phi(\bar j), \quad (\bar j \in
J/N), (\bar h\bar n \in HN/N)
\]
which we can therefore identify with the set of functions
\[
\phi\colon J \rightarrow j_{H\cap N}(E)
\]
smooth with compact support modulo $HN$ in $J$ such that
\[
\phi(hnj)= \delta^{-1}(h) \sigma_{H\cap N}(\bar h)\phi(j), \quad (j \in J), (h\in
H), (n \in N).
\]

Such functions are constructed as follows: fix an element $v$ of $E$ and an element $x$ of $J$. Then
there exists a compact open subgroup $K_0$ of $J$ such that the formula
\begin{equation}\label{P0}
\phi(hnxk)=\delta^{-1}(h)\sigma_{H\cap N}(h)\cdot \bar v, \quad (h\in H), (n\in
N), (k \in K_0)
\end{equation}
indeed defines a function $\phi$ of $\ind_{HN/N}^{J/N}(j_{H\cap N}(E)\otimes \delta^{-1})$ supported in
the open set $HNxK_0$. To do this, we must ensure that if
\begin{equation}\label{P1}
 h_1n_1xk_1=h_2n_2xk_2, \quad (h_1,h_2\in H), (n_1,n_2\in
N), (k_1,k_2 \in K_0)
\end{equation}
then
\begin{equation}\label{P2}
 \delta^{-1}(h_1)\sigma_{H\cap N}(h_1)\cdot \bar v=
 \delta^{-1}(h_2) \sigma_{H\cap N}(h_2)\cdot \bar v
\end{equation}
Now (\ref{P1}) is equivalent to
\begin{align}\label{P3}
\nonumber h_2^{-1}h_1&=n_2xk_2k_1^{-1}x^{-1}n_1^{-1}\\
&= xk_2k_1^{-1}x^{-1} ((xk_2k_1^{-1}x^{-1})^{-1} n_2
(xk_2k_1^{-1}x^{-1})) n_1^{-1}.
\end{align}
Thus $h=h_2^{-1}h_1$ is written $h=xkx^{-1}n$ for a certain $k \in K_0$ and a certain $n \in N$.

Since $v$ is fixed by a certain compact open subgroup of $H$, there exists an open subgroup $K$ of $J$
such that $xKx^{-1}\cap H$ fixes $v$. Let us project $xKx^{-1}\cap H$ into $H/H\cap N$. The inverse image
of this subset in $H$ is $(xKx^{-1}\cap H)(H\cap N)$. Then take $K_1$ $\subset K$ so that
\[
(xK_1x^{-1})N\cap H\subset
(xKx^{-1}\cap H)(H\cap N);
\]
this is possible because $(xKx^{-1}\cap H)(H\cap N)$ is an open neighborhood of $1$ in $H$, and $H$ is
closed in $J$, so that the $(xK_1x^{-1})N\cap H$, $K_1\subset K$, form a basis of neighborhoods of
$H\cap N$ in $H$, modulo $H\cap N$. Then according to (\ref{P3}),
\[
h=h_2^{-1}h_1=xkx^{-1}n \in (xKx^{-1}\cap H)(H\cap N)
\]
satisfies
\[
\sigma_{H\cap N}(h)\cdot \bar v= \bar v.
\]
On the other hand, the modular character $\delta$ is trivial on the compact open subgroups of $H$ and on
$H\cap N$, so $\delta(h)=1$. This shows (\ref{P2}) for such a choice of $K_0$.

The space $\ind_{HN/N}^{J/N}(j_{H\cap N}(E)\otimes \delta^{-1})$ is generated by the functions $\phi$
constructed in this way. It therefore suffices to find $f\in \ind_H^J(E)$ such that $\bar f=\phi$ for
$\phi$, $x$, $v$, $K_0$ as in (\ref{P0}). Take $f$ supported in $HxK_0$ such that
\[
f(hxk)=\sigma(h)\cdot v , \quad (h \in H), (k \in K_0).
\]
This is well defined, because if $h_1xk_1=h_2xk_2$,
\[
h_2^{-1}h_1= xk_2k_1^{-1}x^{-1} \in xK_0x^{-1}\cap H \subset  xKx^{-1}\cap H
\]
fixes $v$. It is then clear that
\begin{equation}\label{P4}
\bar f(j)= \int_{H\cap N\backslash N} \overline{f(nj)}\;  d\nu_{H\cap N\backslash N}(n)
\end{equation}
is non-zero only if $j \in NHxK_0=HNxK_0$. It therefore suffices to verify that
$\bar f(x)=\phi(x)=\bar v$. Now
\begin{equation}\label{P5}
\bar f(x)= \int_{H\cap N\backslash N} \overline{f(nx)}\;  d\nu_{H\cap N\backslash N}(n)
\end{equation}
and if $nx=hxk$, with $h \in H$, $k \in K_0$, we have
\[
h=nxk^{-1}x^{-1}  \in (xK_0x^{-1})N \cap H \subset  (xKx^{-1}\cap H )(N
\cap H)
\]
and thus $\overline {f(nx)}= \overline{f(hxk)}=\overline {\sigma(h)\cdot v}=\bar v$, therefore
\[
\bar f(x)= c \, \bar v,
\]
where $c$ is the measure of a certain compact open set in $H\cap N\backslash N$ for the measure
$\nu_{H\cap N\backslash N}$. We thus find that $\bar f$ is equal to a non-zero scalar multiple of $\phi$,
which is obviously sufficient.

\medskip

Let us now show the injectivity of (\ref{1P}). Let $f \in \ind_H^J(E)$ such that $\bar f = 0$. We must
show that $f$ has a zero image in $j_N(\ind_H^J(E))$, i.e., that there exists a compact open subgroup
$N_0 \subset N$ such that $\int_{N_0} f(jn) \, dn = 0$ for all $j \in J$.

Since $f$ is smooth and compactly supported modulo $H$, there exists a compact open subgroup
$K \subset J$ such that $f$ is right-invariant by $K$, and a compact subset $X \subset J$ such that
$\supp(f) \subset HX$. We can cover $X$ by a finite number of left cosets modulo $K$, say
$X \subset \bigcup_{i=1}^r x_i K$.

For each $i \in \{1, \dots, r\}$, define the function $\phi_i : N \rightarrow E$ by
$\phi_i(n) = f(n x_i)$. The function $\phi_i$ belongs to $\ind_{H \cap N}^N(E)$. The hypothesis
$\bar f(x_i) = 0$ means that the image of $\phi_i$ in the Jacquet module $E_{H \cap N}$ is zero. Let us
show that there exists a compact open subgroup $N_i\subset N$ such that $\int_M\phi_i(n)\,dn=0$ for every
compact open subgroup $M\subset N$ containing $N_i$. Set $w_M=\int_M\phi_i(n)\,dn\in E$. Since
$\phi_i(n_1n)=\sigma(n_1)\cdot\phi_i(n)$ for $n_1\in H\cap N$, and $\phi_i$ is compactly supported modulo
$H\cap N$, there exists $M_0$ such that for every $M\supset M_0$, $w_M$ is a non-zero multiple of
$e_{M\cap H}\cdot w_{M_0}$, the image of $w_{M_0}$ in $E_{H\cap N}$ being a multiple of $\bar f(x_i)=0$.
Thus $w_{M_0}\in E(H\cap N)$, and there exists a compact open subgroup $N'$ of $H\cap N$ such that
$e_{N'}\cdot w_{M_0}=0$. It suffices to take for $N_i$ a compact open subgroup of $N$ containing $M_0$
and $N'$: if $M\supset N_i$, then $e_{M\cap H}=e_{M\cap H}e_{N'}$, and therefore $w_M=0$.

We are looking for a single subgroup $N_0 \subset N$ which annihilates the integral for all $j \in J$.
Consider the set:
\[
C = \bigcup_{i=1}^r \bigcup_{k \in K} k^{-1} x_i^{-1} N_i x_i  k.
\]
Since $N$ is normal in $J$, $x_i^{-1} N_i x_i \subset N$. Since $K$ is compact, $C$ is a compact subset
of $N$. The group $N$ being the increasing directed union of its compact open subgroups, there exists a
compact open subgroup $N_0 \subset N$ containing $C$.

Set $F(j) = \int_{N_0} f(jn) \, dn$. By construction, the function $F$ is right-invariant by $N_0$. Let
us evaluate $F$ on $HX$. If $j \in HX$, we can write $j = h x_i k$ with $h \in H$,
$i \in \{1, \dots, r\}$ and $k \in K$. We then have:
\begin{align*}
F(h x_i k) &= \int_{N_0} f(h x_i k n) \, dn = \sigma(h) \int_{N_0} f(x_i k n k^{-1} k) \, dn \\
&= \sigma(h) \int_{k N_0 k^{-1}} f(x_i n') \, dn' \quad (\text{with } n' = k n k^{-1}) \\
&= \delta_N(x_i)\,\sigma(h) \int_{x_i k N_0 k^{-1} x_i^{-1}} f(n'' x_i) \, dn'' \quad (\text{with } n'' = x_i n' x_i^{-1}).
\end{align*}
Here $\delta_N(x_i)$ denotes the modulus of the automorphism $\Int(x_i)$ of $N$ (\ref{fmod}); the first
change of variables introduces no factor, since $k$ belongs to the compact subgroup $K$. Set
$M_{i,k} = x_i  k N_0 k^{-1} x_i^{-1}$. By definition of $N_0$, we have
$N_0 \supset k^{-1} x_i^{-1}  N_i x_ik$, which implies $M_{i,k} \supset N_i$. Since $M_{i,k}$ is a
compact open subgroup of $N$ containing $N_i$, the integral $\int_{M_{i,k}} \phi_i(n'') \, dn''$ is zero.
Thus, $F(j) = 0$ for all $j \in HX$.

For an arbitrary element $j \in J$, two cases arise: 1. If $j \in HXN_0$, we write $j = x n_0$ with
$x \in HX$ and $n_0 \in N_0$. The right-invariance of $F$ by $N_0$ gives $F(j) = F(x n_0) = F(x) = 0$. 2.
If $j \notin HXN_0$, then for all $n \in N_0$, $jn \notin HX$. Since $\supp(f) \subset HX$, we have
$f(jn) = 0$ for all $n \in N_0$, whence $F(j) = \int_{N_0} 0 \, dn = 0$.

In conclusion, $F(j) = 0$ for all $j \in J$, which proves that $f$ has a zero image in $j_N(\ind_H^J(E))$
and completes the proof of injectivity.

\bigskip

Let us return to diagram (\ref{diagrammedifficile}) by applying the above to $J=Q$, $H=\hat P\cap Q$,
$N=U$, $J/N=Q/U \simeq L$, $L'=H/H\cap N \simeq HN/N = \hat P\cap Q/\hat P \cap U\simeq \hat P \cap L$,
since $\hat P\cap Q=(\hat P\cap L)(\hat P\cap U)$ by the decomposition of $\hat P\cap Q$ recalled above.
Here $\delta=\delta_U\delta_{\hat P\cap U}^{-1}$, so that
$\delta^{-1}\delta_U^{1/2}=\delta_{\hat P\cap U}\delta_U^{-1/2}=\delta_{\hat P\cap U}^{1/2}\epsilon_2$:
the effects of the modular characters cancel out. The natural isomorphism between
$r_U\circ \ind_{\hat P \cap Q}^Q(E)$ and $\ind_{L \cap \hat P}^L (r_{\hat P\cap U}(E)\otimes \epsilon_2)$
is obtained, by passing to the quotient, from
\[
f \mapsto \bar f , \quad \ind_{\hat P \cap Q}^Q(E)
\rightarrow \ind_{L \cap \hat
P}^L (r_{\hat P\cap U}(E)\otimes \epsilon_2)
\]
\begin{equation}\label{forMule}
\bar f(l)= \delta_U(l)^{1/2}\int_{\hat P \cap U\backslash U} \overline{f(ul)} \;
d\nu_{\hat P \cap U\backslash U}(u).
\end{equation}

\bigskip

Let us now verify the commutativity of
\[
\xymatrix{
\hat P=\hat M \hat N  \ar[dr]^{\caF_{\hat P\cap Q}^{\hat P}}
\ar[rr]^{\caI_{\hat Z}}&    &Q=LU \\
&  \hat P\cap Q\ar[ur]^{\ind_{\hat P\cap Q}^Q} &
}
\]
Now $\caI_{\hat Z}=  \Gamma_c \circ  \caR_{Z,Q}^{X,G}   \circ  \iota^{-1}\circ \caI_{\hat P}^G $ and
$\ind_{\hat P\cap Q}^Q=   \Gamma_c \circ   \caI_{\hat P\cap Q}^Q$. Note that the stabilizer in $Q$ of the
orbit $Z=PzQ$ is none other than $\hat P \cap Q$, so that $(\hat P \cap Q)\backslash Q\simeq Z$. It
therefore suffices to prove that
\[
\xymatrix{
\caM(\hat P)\ar[r]^{\caI_{\hat P}^G} \ar[d]^{\caF_{\hat P\cap Q}^{\hat P}}  & \scrC^\infty_{\hat P\backslash G,G}
-\caM od \ar[r]^{\iota}& \scrC^\infty_ {P\backslash G,G}-\caM od
\ar[r]^{ \caR_{Z,Q}^{X,G}} & \scrC^\infty_{Z,Q}-{\caM od} \ar[d]\\
\caM(\hat P\cap Q)\ar[rrr]^{\caI_{\hat P\cap Q}^Q}& & & \scrC^\infty_{(\hat P\cap Q)\backslash Q,Q}-\caM od}
\]
commutes.

Take a representation $(\sigma,E)$ in $\caM(\hat P)$. The sheaf corresponding to it in
$\scrC^\infty_{\hat P\backslash G,G}-\caM od$ by the equivalence of categories of \ref{indetfaisc}
therefore has fiber $E$ over $\hat P\in \hat P\backslash G$. Let $X$ be the variety of parabolic
subgroups conjugate to $\hat P$ (or $P$). We have
\[
\hat P\backslash G \simeq X \simeq P\backslash G.
\]
These isomorphisms are obtained by choosing respectively $\hat P=Pz$ (careful, this is an equality of
points in $X$, the equality of parabolic subgroups of $G$ corresponding to it is $\hat P=z^{-1}Pz$) and
$P$ as base point in $X$, the isomorphism between the left term and the right one being realized by
$\iota$. The sheaf corresponding to $(\sigma,E)$ in $\scrC^\infty_{\hat P\backslash G,G}-\caM od$
therefore has fiber $E$ over $Pz$. Its restriction to $Z=PzQ$ likewise, so by identifying $Z$ and
$(\hat P\cap Q\backslash Q)$, we obtain a sheaf on $\scrC^\infty_{\hat P\cap Q\backslash Q, Q}-\caM od$,
with fiber $E$ over the base point $\hat P\cap Q$. The action of $Q$ being transitive, this fiber
characterizes the sheaf. On the other hand, it is now clear that this is also the fiber over this point
of the sheaf $\caI_{\hat P\cap Q}\circ \caF_{\hat P\cap Q}^{\hat P}(\sigma,E)$. This shows the
commutativity of the diagram.

\medskip

It remains to get rid of the characters $\epsilon_1$ and $\epsilon_2$. Since these characters are trivial
respectively on $\hat M \cap U$ and $\hat N \cap L$ we have
$\epsilon_1\circ r_{\hat M\cap U}=r_{\hat M\cap U}\circ \epsilon_1$ and
$\epsilon_2\circ \caF_{L \cap \hat P}^{L \cap \hat M } =\caF_{L \cap \hat P}^{L \cap \hat M}\circ \epsilon_2$.
It therefore suffices to see that the character
\[
\epsilon=\epsilon_1\epsilon_2= \delta_{\hat N}^{-1/2}\delta_{\hat
N\cap Q}^{1/2} \delta_{U}^{-1/2} \delta_{U\cap \hat P}^{1/2}
\]
is trivial on $L'=\hat M\cap L$. We can assume, by reducing to this case by conjugation, that the
parabolic subgroups $P$ and $Q$ are standard. Write the Cartan decomposition of $L'$ in the form
$L'=K_0'C_\emptyset^+F'K_0'$, where $F'$ is finite and some power of each element of $F'$ lies in
${}^0M_\emptyset A_\emptyset$ (see \ref{calcfonctmod}). Since the character $\epsilon$ with values in
$\bbR^\times_+$ is trivial on the compact subgroups of $L'$, hence on $K_0'$, it suffices to verify that
$\epsilon$ is trivial on $A_\emptyset$.

We retain the notation of \ref{Weylgroups}. For any unipotent subset $S$ of $\Sigma_{\emptyset}$, we have
\[
\delta_{U(S)}(a)=\prod_{\gamma \in S} |\gamma(a)^{m_\gamma+2m_{2\gamma}} |, \quad
(a \in A_\emptyset).
\]
Let us denote this simply $\delta_S$. Then we have
\[
\epsilon^2= \delta_{\hat \caN}^{-1} \delta_{\hat \caN \cap \caQ}  \delta_{\caU}^{-1}\delta_{ \caU \cap \hat \caP }
\]
Now,
\begin{align*}
&\hat \caN \setminus (\hat \caN \cap \caQ )=\hat \caN
\cap(\Sigma_{\emptyset}\setminus \caQ ) =\hat \caN \cap(-\caU) \\
&\caU \setminus ( \caU \cap \hat \caP )= \caU
\cap(\Sigma_{\emptyset} \setminus \hat  \caP )= \caU \cap(-\hat \caN),
\end{align*}
whence, since the subsets $\hat\caN\cap(-\caU)$ and $\caU\cap(-\hat\caN)$ are opposite to each other, and
$\delta_{-S}=\delta_S^{-1}$ by the above formula,
\[
\epsilon^2=\delta_{\hat \caN \cap(-\caU)}^{-1}\delta_{\caU \cap(-\hat \caN)}^{-1} =1.
\]
Since $\epsilon$ takes values in $\bbR^\times_+$, we deduce that $\epsilon=1$.
\end{proof}

It will be useful later to have a slightly more explicit form of the natural isomorphism of the first
theorem of \ref{geomlemma}. Let us resume the notation from the beginning of this section: $P=MN$ and
$Q=LU$ are two parabolic subgroups of $G$, $Z$ is an orbit of $Q$ in $P\backslash G$, and $Y$, $Y'$ are
two open subsets of $P\backslash G$, unions of $Q$-orbits, such that $Y$ is the disjoint union of $Y'$
and $Z$. Let $(\sigma,E)$ be a smooth representation of $M$. The functor $F_Z$ is a quotient of the
functor $i_Y$: we first quotient by the subfunctor $i_{Y'}$, then we apply the functor $r_U$ which is
itself also a passage to the quotient. The natural isomorphism between
\[
F_Z(\sigma,E)= r_U(i_Y(\sigma,E))/ r_U(i_{Y'}(\sigma,E))
\]
and
\[
\Phi_Z=i_{L\cap \iota^{-1}(P)}^L\circ \iota \circ r^M_{M\cap \iota(Q)}
\]
is induced by the map
\[
f \mapsto \bar f, \quad i_Y(E) \rightarrow  i_{L\cap \iota^{-1}(P)}^L\circ \iota \circ r^M_{M\cap \iota(Q)}(E)
\]
where
\begin{equation}\label{IndRes}
  \bar f(l)=\delta_U(l)^{1/2} \; \int_{U\cap
  \iota^{-1}(P)\backslash U} j_{M\cap \iota(U)}(f(\iota(ul)))\; d\nu_{U\cap
  \iota^{-1}(P)\backslash U}(u).
\end{equation}
This is obtained from (\ref{forMule}) and the other natural isomorphisms of the diagram.

\subsection{Cuspidal data}\label{donneescusp}
\index[ter]{cuspidal datum} From the geometric lemma \ref{geomlemma}, we deduce information about the
Jordan-Hölder series of representations induced from supercuspidal representations. More precisely, let
$P=MN$ be a parabolic subgroup of $G$, and $(\rho,W)$ an irreducible supercuspidal representation of $M$.
We want to study the induced representation $\pi=i_P^G\rho$.

\needspace{165pt}
\begin{defi}
A cuspidal datum is a pair $(M,(\rho,W))$ where $M$ is a Levi factor of a parabolic subgroup $P=MN$ of
$G$ and $(\rho,W)$ is an irreducible supercuspidal representation of $M$.

We say that two cuspidal data $(M_1,(\rho_1,W_1))$ and $(M_2,(\rho_2,W_2))$ are associated if there
exists $g\in G$ such that
\[
gM_1g^{-1} = M_2 \quad \text { and }  \quad \rho_2 \simeq \rho_1^g
\]

This defines an equivalence relation on the set of cuspidal data. We denote by $(M,(\rho,W))_G$
\index[not]{M(rho,W)_G@$(M,(\rho,W))_G$} the class of $(M,(\rho,W))$ and $\Omega(G)$
\index[not]{ZZOmegaG@$\Omega(G)$} the set of equivalence classes.
\end{defi}

Recall that according to Corollary \ref{supercusp}, for any irreducible smooth representation $(\pi,V)$
of $G$, one can find a parabolic subgroup $P=MN$ of $G$ and an irreducible supercuspidal representation
$(\tau,E)$ of $M$ such that
\begin{equation}\label{Homrmpi}
\Hom_M(r_P^G(\pi,V),(\tau,E))\neq \{0\}.
\end{equation}
Equivalently, by Frobenius reciprocity
\begin{equation} \label{Homrmpi2}
\Hom_G((\pi,V),i_P^G(\tau,E)) \neq \{0\}.
\end{equation}

The following result ``separates'' induced representations and supercuspidal representations:

\begin{lemme}
Let $P=MN$ be a proper parabolic subgroup of $G$ and let $(\rho,W)$ be a smooth representation of $M$.
Set $(\pi,V) =i_P^G(\rho,W)$. Then no irreducible subquotient of $(\pi,V)$ is supercuspidal.
\end{lemme}
\begin{proof}
Let $\pi_{sc}$ be the supercuspidal part of $\pi$ in the decomposition of Theorem \ref{decMG2}. We have
\[
\Hom_G(\pi_{sc},\pi)  =  \Hom_M(r_P^G( \pi_{sc}),\rho)=0
\]
and this implies that $\pi_{sc}$ is zero.
\end{proof}

\subsection{Consequences of the geometric lemma} \label{geolemmacusp}
The reader is referred to Section \ref{Weylgroups} for certain notation concerning Weyl groups. One
easily deduces from the geometric lemma \ref{geomlemma} the following results:

\begin{prop}
Let $P=MN$ and $Q=LU$ be parabolic subgroups of $G$, $(\rho,W)$ a supercuspidal representation of $M$ and
set $\tau= r_Q^Gi_P^G\rho$. Then

$(i)$ if $L$ admits no Levi subgroup conjugate to $M$ in $G$, then $\tau=0$,

$(ii)$ if $M$ is not conjugate to $L$ in $G$, then $\tau$ has no supercuspidal subquotient,

$(iii)$ suppose $M$ and $L$ are standard and conjugate, then $\tau$ admits a filtration by
subrepresentations whose successive quotients are of the form ${}^w \rho$, $w \in W(L,M)/W_L$. In
particular $\tau$ is supercuspidal.
\end{prop}

\begin{proof}
$(i)$ Since $\rho$ is supercuspidal, $r_{P'}^M(\rho)=0$ for any proper parabolic subgroup $P'=M'N'$ of
$M$. Thus, $\tau$ is zero unless $P'=M'=M$ with $M'=M\cap w\cdot L$, for a certain $w \in \caW^{Q,P}$.
The group $M$ is then a Levi subgroup of $w\cdot L$.

\noindent $(ii)$ This is an immediate consequence of $(i)$ and Lemma \ref{donneescusp}.

\noindent $(iii)$ This follows from Lemma \ref{WeylgroupsS}, $(iii)$.
\end{proof}

\subsection{Jordan-Hölder series of induced supercuspidal representations}\label{compindsup}
Let $M$ be a standard Levi subgroup of $G$. Recall that
\[
W(*,M)= \bigcup_L W(L,M),\quad l(M)=|W_M\backslash W(*,M)|
\]
where the sum is over the set of standard Levi subgroups of $G$ (see \ref{PGQ1} and \ref{WoM} for the
notation). If $M$ is a non-standard Levi subgroup, we set $l(M)=l(M')$ where $M'$ is a standard Levi
subgroup of $G$ conjugate to $M$.

\begin{thm}
Let $P=MN$ be a parabolic subgroup of $G$, $(\rho,W)$ an irreducible supercuspidal representation of $M$
and $(\pi,V)=i_P^G(\rho,W)$. Then the representation $\pi$ is of finite length, less than or equal to
$l(M)$. On the other hand, if $P'=M'N'$ is another parabolic subgroup of $G$, if $(\rho',W')$ is an
irreducible supercuspidal representation of $M'$, and if we set $\pi'=i_{P'}^G \rho'$ then the following
conditions are equivalent:

$(i)$ there exists $g \in G$ such that $gMg^{-1}=M'$ and $\rho^g=\rho'$, {\sl i.e.,} the cuspidal data
$(M,(\rho,W))$ and $(M',(\rho',W'))$ are associated,

$(ii)$ $\Hom_G(\pi,\pi')\neq 0$,

$(iii)$ the Jordan-Hölder series of $\pi'$ and $\pi$ are equivalent,

$(iv)$ the Jordan-Hölder series of $\pi'$ and $\pi$ have at least one element in common.

Finally, if we assume $P$ and $P'$ are standard, and if we set
\[
W(\rho,\rho') = \{w \in W_G \mid w\cdot M=M', {}^w\rho'=\rho    \},
\]
then $\dim \Hom_G(\pi,\pi')\leq |W_M \backslash W(\rho,\rho')|$.
\end{thm}
\begin{proof}
Up to conjugation, we may assume $P$ and $P'$ are standard, which we assume henceforth. Let $(\pi_0,V_0)$
be an irreducible subquotient of $(\pi,V)$. Corollary \ref{supercusp} asserts that there exists a
standard parabolic subgroup $Q=LU$ of $G$ such that $r_Q^G \pi_0$ is supercuspidal, non-zero. Proposition
\ref{geolemmacusp} $(ii)$ then shows that $M$ and $L$ are conjugate. Moreover, there exists an
irreducible supercuspidal subrepresentation $(\rho_0,W_0)$ of $r_Q^G (\pi_0,V_0)$ such that
$(L,(\rho_0,W_0))$ and $(M,(\rho,W))$ define associated cuspidal data. Since the functor $r_Q^G$ is
exact, we get that $(\rho_0,W_0)$ is a subquotient of $r_Q^G \pi=r_Q^G i_P^G \rho$ and the assertion then
follows from Proposition \ref{geolemmacusp}, $(iii)$.

Define a function from $\caM(G)$ taking values in $\bbN \cup \{+\infty\}$ by
\[
l'(\tau,E)= \sum_{L \sim M,\, Q=LU } l(r_Q^G \tau)
\]
where the sum is over the standard Levi subgroups of $G$ conjugate to $M$, $Q=LU$ is the standard
parabolic subgroup with Levi factor $L$, and $l$ is the function on $\caM(G)$ taking values in
$\bbN \cup \{+\infty\}$ giving the length of a representation. Since the functors $r_Q^G$ are exact, $l'$
is additive, i.e., if $(\tau_1,E_1)$ is a subrepresentation of $(\tau,E)$, then
$l'(E)= l'(E_1)+l'(E/E_1)$. From the above, we see that $1=l(\pi_0)\leq l'(\pi_0)$. By additivity, and
induction on the length, we see that $l(\pi,V)\leq l'(\pi,V)$. Now
$$l'(\pi,V)=\sum_{L \sim M,\, Q=LU }|W_M\backslash W(L,M)|=|W_M\backslash W(*,M)|=l(M)$$ according to
\ref{geolemmacusp}, $(iii)$. We deduce the first assertion of the theorem. This allows us to speak of a
Jordan-Hölder series ({\sl cf.} \ref{abcat}, 10) and in particular, ensures the existence of irreducible
subrepresentations.

It is clear that $(ii)\Rightarrow (iv)$ and that $(iii)\Rightarrow (iv)$. Let us now show that
$(iv)\Rightarrow (i)$. Let $(\pi_0,V_0)$ be an irreducible subquotient common to $(\pi,V)$ and
$(\pi',V')$. We have seen above how to attach to $(\pi_0,V_0)$ a cuspidal datum $(L,(\rho_0,W_0))$, which
is then associated to $(M,(\rho,W))$. The same reasoning applies to $(M',(\rho',W'))$. This shows that
the cuspidal data $(M,(\rho,W))$ and $(M',(\rho',W'))$ are associated.

We now show $(i) \Rightarrow (ii)$. By the adjunction of the functors $i_{P'}^G$ and $r_{P'}^G$, we have:
\begin{equation}\label{pipi}
\Hom_G(\pi, \pi')\simeq \Hom_{M'}(r_{P'}^G i_P^G \rho,\rho').
\end{equation}
Now according to Proposition \ref{geolemmacusp} $(iii)$, $\rho'$ is a composition factor of
$r_{P'}^G i_P^G \rho$. According to Lemma \ref{pisc}, $\rho'$ appears as a quotient of
$r_{P'}^G i_P^G \rho$ and thus (\ref{pipi}) is non-zero.

Moreover (\ref{pipi}) shows that we also have, still using Proposition \ref{geolemmacusp} $(iii)$
\[
\dim \Hom_G(\pi, \pi')=\dim  \Hom_{M'}(r_{P'}^G i_P^G \rho,\rho')\leq
| W_M\backslash  W(\rho,\rho ') |,
\]
which proves the last assertion of the theorem.

It remains to show $(i)\Rightarrow (iii)$. We first assume that $M$ is a maximal Levi subgroup, i.e.,
$l(M)=2$ (see \ref{levmax}). Since we have shown that $(i)\Rightarrow (ii)$, consider non-zero
intertwining operators $A: V \rightarrow V'$ and $A': V' \rightarrow V$. We then have
$l'(\pi)=l'(\pi')=2$ and thus $l(\pi)=1$ or $2$ (ditto for $l(\pi')$). Suppose $l(\pi)=1$, i.e., $\pi$ is
irreducible. Then $A$ is injective, and since
\[
l'(V'/A(V))= l'(V')-l'(A(V))=l'(V')-l'(V)=0
\]
we obtain $V'=A(V)$ and $\pi'$ is irreducible. This completes the proof in this case, and by symmetry, in
the case where $\pi'$ is irreducible. Suppose for the sake of contradiction that $l(\pi)=l(\pi')=2$ and
let $(\pi_0,V_0)$, $(\pi'_0,V'_0)$ be irreducible subrepresentations respectively of $(\pi,V)$ and
$(\pi',V')$. We then have $ l'(V/V_0)=l'(V'/V'_0)=l(V/V_0)=l(V'/V'_0)=1$. Since the case $M=M'$ and
$\rho \simeq \rho'$ is trivial, we assume this is not the case. Since
\[
\Hom_G(\pi_0,\pi)= \Hom_M(r_P^G \pi_0,\rho)\neq \{0 \}
\]
and $l'(\pi_0)=1$, we have
\[
\Hom_G(\pi_0,\pi')= \Hom_{M'}(r_{P'}^G \pi_0,\rho')= \{0 \}
\]
and thus $A(V_0)=0$. Similarly $A'(V'_0)=0$. It follows that $V_0$ is the unique irreducible
subrepresentation of $V$ (otherwise we would have $A=0$) and similarly $V'_0$ is the unique irreducible
subrepresentation of $V'$. It follows that $V/V_0\simeq V'_0$, $V'/V'_0\simeq V_0$, and the composition
factors of $V$ and $V'$ are therefore $V_0$ and $V'_0$. This completes the proof in the case $l(M)=2$.

Let $w \in W(L,M)$. We say that $w$ is elementary if there exists a standard parabolic subgroup $T=SV$ of
$G$ containing $P$ and $Q$, such that $w \in W_S$ and such that $l(M)=2$ as a standard Levi subgroup of
$S$, so we can apply the $l(M)=2$ case to $S$. Since $r_Q^G=r_{Q\cap S}^S \circ r_T^G$ and
$i_P^G= i_T^G \circ i_{P\cap S}^S$, we deduce $(iii)$ in the case where $w$ is elementary and
$\rho'={}^w \rho$.

We can conclude the proof using the following lemma (\cite{BeZe2}, Lemma 2.17).
\begin{lemme}
Let $M$ and $L$ be standard Levi subgroups of $G$ and $w \in W(L,M)$. Then there exists a sequence of
standard Levi subgroups of $G$
\[
L=L_0, L_1, \ldots, L_m =M
\]
and for each $i=1,\ldots, m$, $w_i\in W(L_{i-1},L_i)$ elementary, such that $w=w_m w_{m-1}\ldots w_1$.
\end{lemme}

This completes the proof of the theorem.
\end{proof}

\needspace{170pt}
\begin{cor}
Let $(\rho,W)$ be an irreducible supercuspidal representation of a Levi subgroup $M$ of $G$.

$(i)$ Let $P$ and $Q$ be two parabolic subgroups of $G$ with Levi factor $M$. Then the Jordan-Hölder
series of $i_P^G\rho$ and $i_Q^G\rho$ are equivalent.

$(ii)$ There always exists a non-zero intertwining operator between $i_P^G\rho$ and $i_Q^G\rho$.

$(iii)$ Let $P\in \caP(M)$ and $(\pi,V)$ an irreducible subquotient of $i_P^G \rho$. Then there exists
$Q \in \caP(M)$ such that $(\pi,V)$ is a subrepresentation of $i_Q^G(\rho,W)$.
\end{cor}

\begin{proof}
Only the third point deserves a brief explanation. We know from Corollary \ref{supercusp} that there
exists a parabolic subgroup $Q'=L'U'$ of $G$ and an irreducible supercuspidal representation $(\tau,E)$
of $L'$ such that $(\pi,V)$ is a subrepresentation of $i_{Q'}^G(\tau,E)$. According to the theorem, the
cuspidal data $(L',(\tau,E))$ and $(M,(\rho,W))$ are associated, i.e., there exists $g \in G$ such that
$g\cdot L'=M$ and $\tau^g=\rho$. We then take $Q=g \cdot Q'$.
\end{proof}

\section{Finiteness theorems} \label{thmsfinit}

In this section, we show that the functors $i_P^G$ and $r_P^G$ preserve certain properties of
representations. We already know from Lemma \ref{IndAdm} that $i_P^G$ preserves admissibility and from
Proposition \ref{finiIJ} that $r_P^G$ preserves finite generation. We show that $r_P^G$ preserves
admissibility and finite length and that $i_P^G$ preserves finite length. We will show later
(\ref{noether}) that $i_P^G$ preserves finite generation.

\subsection{Jacquet's lemma} \label{LemmeJacquet}
\index[ter]{Jacquet's lemma} The ``Jacquet lemma'' (here promoted to the rank of a theorem) is proved
under an admissibility hypothesis. This hypothesis will be removed later (see Section
\ref{LemmeJacquet2}), and the proof will be considerably more difficult. The weak version proved here is
sufficient to obtain as a corollary the fact that the functors $r_P^G$ preserve admissibility. The notion
of an Iwahori decomposition of a compact open subgroup and the related notation were introduced in
\ref{KN}.

\begin{thm}
Let $P=MN$ be a parabolic subgroup of $G$, with split component $A$, and let $(\pi,V)$ be a smooth
admissible representation of $G$. Let $K$ be a compact open subgroup of $G$ admitting an Iwahori
decomposition with respect to $P$. Then the projection map $V\rightarrow  V_N$ maps $V^K$ surjectively
onto $(V_N)^{K_M}$.
\end{thm}
\begin{proof}
Let $j:V\rightarrow  V_N$ denote the projection map for the $P$-module structure. If $v \in V^K$, then
for all $k\in K_M$ we have:
\[
\pi_N(k)j(v)=j(\pi(k)\cdot v)=j(v)
\]
and thus the image of $V^K$ under $j$ is indeed in $(V_N)^{K_M}$. We then choose $t \in C_{A}^{++}$. The
set $\{ t^{-m}K_{\overline{N}}t^m|\, m\in \bbN \}$ forms a basis of neighborhoods of the identity element
in $\overline{N}$ (Theorem \ref{KN}). Let $v \in  V^K$. Since
\[
t^{-1}K_{\overline{N}}t \subset K_{\overline{N}}, \quad
t^{-1}K_{M}t=K_M,
\]
we have:
\[
\pi(e_{K_{\overline{N}}})\pi(t)\cdot v= \pi(e_{K_{M}})\pi(t)\cdot v= \pi(t)\cdot v,
\]
and thus:
\begin{align} \label{calcul}
j(a_{t,K}\cdot v)&= j(e_{K}*\delta_t* e_K\cdot v) = j(\pi(e_{K})\pi(t)\cdot v)\\
\nonumber &= j(\pi(e_{K_N}*e_{K_{M}}*e_{K_{\overline{N}}}  )\pi(t)\cdot v)\\
\nonumber &=  j(\pi(e_{K_N}) \pi(t)\cdot v)= \pi_N(e_{K_N})\pi_N(t)\cdot j(v)\\
\nonumber  &= \pi_N(t)\cdot j(v).
\end{align}
The last equality comes from the fact that $N$, and therefore $K_N$, acts trivially on $V_N$. This shows
that $\pi_N(t)\cdot j(V^K) \subset j(V^K) $. In fact, equality holds, because $\pi_N(t)$ acts invertibly
on $V_N$ and $j(V^K)$ is finite-dimensional. It follows that $\pi_N(t^m)\cdot j(V^K) =j(V^K) $, for all
$m\in \bbZ$. Let $\bar v \in (V_N)^{K_M}$, let $v'$ be any lift of $\bar v$ in $V$, and let
$v=\pi(e_{K_M})\cdot v'$. We have $v \in V^{K_M}$ and $j(v)=\bar v$. Fix $m\in \bbN$ such that $v$ is
fixed by $t^{-m}K_{\overline{N}}t^m$, so that $K_{\overline{N}}$ fixes $\pi(t^m)\cdot v$. As before, we
have
\[
j(\pi(e_{K})\pi(t^m)\cdot v)=  \pi_N(e_{K_N})\pi_N(t^m)\cdot
j(v)=\pi_N(t^m)\cdot \bar v,
\]
and thus $\bar v \in \pi_N(t^{-m})j(V^K)=j(V^K)$.
\end{proof}

\begin{cor}
Let $P=MN$ be a parabolic subgroup of $G$, and let $(\pi,V)$ be a smooth admissible representation of
$G$. Then the representation $(\pi_N,V_N)$ of $M$ is admissible. The same is true for $r_P^G (\pi,V)$.
\end{cor}
\begin{proof}
This follows from the theorem and the fact that there exists a basis of neighborhoods of the identity in
$G$ consisting of compact open subgroups admitting an Iwahori decomposition with respect to $P$.
\end{proof}

We will now construct a section $s_P^K$ of the surjective morphism $j\colon V^K \rightarrow V_N^{K_M}$.
For all $v \in V^K \cap \ker j=V^K \cap V(N)$, there exists, by Proposition \ref{FoncJac1}, a compact
open subgroup $\Gamma$ of $N$ such that $e_\Gamma\cdot v=0$. Since $V^K$ is finite-dimensional, we can
even choose $\Gamma$ such that $e_\Gamma\cdot v=0$ for all $v \in V^K \cap \ker j$. According to Theorem
\ref{KN} $(iv)$, there exists $\epsilon >0$ such that for all $t \in C^+_A(\epsilon)$,
$\Gamma \subset t^{-1}K_Nt$. For such a $t$, and for all $v  \in V^K \cap \ker j$ we have:
\[
a_{t,K} \cdot v= e_K*\delta_t*e_K \cdot v=e_{K_{\overline{N}}}*
e_{K_M}*e_{K_{N}}*\delta_t \cdot v= e_{K_{\overline{N}}}*
e_{K_M}* \delta_t*e_{t^{-1}K_{N}t} \cdot v= 0.
\]
This shows that for such a $t$:
\begin{equation} \label{VKkerj}
 V^K \cap \ker j=V^K \cap V(N) \subset \ker \pi(a_{t,K}).
\end{equation}

We will now show that $j$ induces an isomorphism from $\im\pi(a_{t,K}) \subset V^K$ onto $V_N^{K_M}$. Let
us first show injectivity. Suppose that $v \in V^K$ is such that $j(a_{t,K}\cdot v)=0$. From
(\ref{calcul}), we have $\pi_N(t)\cdot j(v)=0$, whence $j(v)=0$ because $\pi_N(t)$ is invertible, and we
conclude from (\ref{VKkerj}) that $a_{t,K}\cdot v=0$.

We will now show surjectivity. Let $\bar v \in V_N^{K_M}$. Since $t$ is in the center of $M$, we have
$\pi_N(t^{-1})\cdot \bar v \in V_N^{K_M}$. By Jacquet's lemma, there exists $v \in V^K$ such that
$j(v)= \pi_N(t^{-1})\cdot \bar v$. Set $v'=a_{t,K}\cdot v \in \im\pi(a_{t,K})$. From (\ref{calcul}),
$j(v')= j(a_{t,K}\cdot v)= \pi_N(t)\cdot j(v)= \bar v$.

We deduce from the above that:
\[
\dim V^K= \dim \ker j_{|V^K} + \dim V_N^{K_M}= \dim \ker
\pi(a_{t,K})+ \dim \im \pi(a_{t,K}),
\]
and $\dim \im \pi(a_{t,K})= \dim V_N^{K_M}$, whence $\ker j_{|V^K}= \ker \pi(a_{t,K})$ and
$$V^K= \ker \pi(a_{t,K})\oplus  \im \pi(a_{t,K}).$$ Set $V_0^K=\ker \pi(a_{t,K}) $,
$V_*^K=\im \pi(a_{t,K})$.

\begin{prop}
The decomposition $V^K= V_0^K\oplus V^K_*$ obtained above is independent of the choice of
$t \in C_A^+(\epsilon)$. By inverting the restriction of $j$ to $V_*^K$, we obtain a section, depending
only on $P$ and $K$,
\[
s_P^K \colon  V_N^{K_M} \rightarrow V^K_* \subset V^K.
\]
\end{prop}

\begin{proof}
Let $t'$ be another element satisfying the same hypotheses as $t$. Then $tt'$ still satisfies these
hypotheses, and since by the second lemma of \ref{HCDH},
\[
a_{t,K}a_{t',K}=a_{tt',K}=a_{t',K} a_{t,K},
\]
we see that
\begin{align*}
\im a_{tt',K}&\subset \im a_{t,K},\\
 \im a_{tt',K}&\subset \im a_{t',K},\\
 \ker a_{t,K} &\subset \ker a_{tt',K},\\
 \ker a_{t',K} &\subset \ker a_{tt',K}
\end{align*}
By equality of dimensions, we see that all these inclusions are in fact equalities.
\end{proof}

\subsection{Induction and finite length}\label{longfin}

We will now show that the induction functor $i_P^G$ preserves modules of finite length. The proof
crucially uses Theorem \ref{compindsup}, which follows from the geometric lemma of induction-restriction.
This will be used in the proof of Proposition \ref{var0G}.

\begin{lemme}
Let $P=MN$ be a parabolic subgroup of $G$ and let $(\tau,E)$ be an irreducible representation of $M$.
Then the representation $i_P^G (\tau,E)  $ of $G$ is of finite length.
\end{lemme}

\begin{proof}
Suppose that $(\tau,E)$ is a supercuspidal representation of $M$. Then $i_P^G (\tau,E)  $ is of finite
length by Theorem \ref{compindsup}. The general case follows because $(\tau,E)$ is a subrepresentation of
a representation of the form $i_{P'}^M(\rho,W)$, where $(\rho,W)$ is an irreducible supercuspidal
representation of a Levi subgroup $M'$ of $M$ (Corollary \ref{supercusp}), the functor $i_P^G$ is exact
and thus $i_P^G (\tau,E)  $ is a subrepresentation of $i_{P'N}^G(\rho,W)=i_P^G i_{P'}^M(\rho,W)$ which is
of finite length by the above.
\end{proof}

\subsection{A theorem of Howe}\label{ThmHowe}
\index[ter]{Howe's theorem} We now want to show the following result:
\begin{thm}
Let $(\pi,V)$ be a smooth representation of $G$. Then $(\pi,V)$ is of finite length if and only if it is
finitely generated and admissible.
\end{thm}

To prove this result, we need some preparation.

\begin{lemme}
Let $P=MN$ be a standard parabolic subgroup of $G$ and $K$ a compact open subgroup of $G$, contained and
normal in $K_0$. If $(\pi,V)$ is generated by $V^K$, then $V_N$ is generated by $V_N^{K\cap M}$.
\end{lemme}
\begin{proof}
Since $G=K_0P$ and $K$ is normal in $K_0$, the space $V$ is clearly generated as a vector space by
vectors of the form $\pi(p)\cdot v$, where $v\in V^K$ and $p\in P$. From this, we deduce that $V_N$ is
generated by vectors of the form $j_N(\pi(p)\cdot v)$, $v\in V^K$ and $p\in P$, where $j_N$ is the
projection of $V$ onto $V_N$. Write $p=mn$, $m\in M$, $n\in N$. Then
\[
j_N(\pi(p)\cdot v)= \pi_N(m)\cdot j_N(v).
\]
Since $j_N(v)\in V_N^{K\cap M}$, the lemma is proved.
\end{proof}

\begin{prop}
Let $K$ be a compact open subgroup of $G$, contained and normal in $K_0$ and let $(\pi,V)$ be a smooth
representation of $G$. If $V^K$ generates $V$, then any non-zero supercuspidal subquotient of $V$
contains non-zero vectors fixed by $K$.
\end{prop}
\begin{proof}
Let $(\tau,E)$ be an irreducible supercuspidal representation of $G$. We use (\ref{dectaunotau}) to
decompose $V$ into a direct sum of subrepresentations $V=V(\tau)\oplus W$, where $V(\tau)$ is a
representation all of whose irreducible subquotients are in the inertial class $[\tau]$, and where no
irreducible subquotient of $W$ is in $[\tau]$. Since $V$ is generated by $V^K$, the same is true for
$V(\tau)$. Consider the restriction of $(\tau,E)$ to ${}^0G$. It is a finite direct sum of irreducible
compact representations of ${}^0G$. Let us call these $(\tau_i,E_i)$. The restriction of $V(\tau)$ to
${}^0G$ is therefore a direct sum of representations isomorphic to one of the $(\tau_i,E_i)$. Since
$V(\tau)$ is generated by $V(\tau)^K$, if $V(\tau) \neq \{0\}$, necessarily one of the $E_i$ contains a
non-zero vector fixed by $K$. Since $K \subset {}^0G$, this then shows that $E$ contains a non-zero
vector fixed by $K$.

Now let $(\rho,W)$ be a non-zero supercuspidal subquotient of $(\pi,V)$. Then any irreducible subquotient
$(\tau,E)$ of $(\rho,W)$ is also an irreducible subquotient of $(\pi,V)$, and from what has just been
established, it contains a non-zero vector fixed by $K$. Since the functor $V\mapsto V^K$ is exact
(Proposition \ref{piK}), we deduce that $(\rho,W)$ contains a non-zero vector fixed by $K$.
\end{proof}

\begin{cor}
Let $K$ be a compact open subgroup of $G$, contained and normal in $K_0$ and admitting an Iwahori
decomposition with respect to the standard parabolic subgroups. Let $(\pi,V)$ be a smooth admissible
representation of $G$ such that $V^K$ generates $V$. Then any subquotient of $V$ is generated by its
vectors fixed by $K$.
\end{cor}
\begin{proof}
Let $(\pi',V')$ be a subquotient of $(\pi,V)$. If the representation $(\pi',V')$ is non-zero and
supercuspidal, the proposition shows that $(V')^K \neq \{0\}$. Otherwise, we can find a standard
parabolic subgroup $P=MN$ such that $V'_N$ is supercuspidal (Corollary \ref{supercusp}, whose proof holds
for any non-zero representation). By the lemma, $V_N$ is generated by its vectors fixed under $K\cap M$.
Since $V_N'$ is a subquotient of $V_N$, we deduce from the proposition that $V'_N$ contains non-zero
fixed vectors under $K\cap M$. By Jacquet's lemma (\ref{LemmeJacquet}), the projection
$j_N \colon V' \rightarrow V'_N$ maps $(V')^K$ surjectively onto $(V'_N)^{K\cap M}$. It follows that
$(V')^K$ is non-zero. It remains to show that $(V')^K$ generates $V'$. Let $W$ be the submodule of $V'$
generated by $(V')^K$. We want to show that $W=V'$. Since $W^K=(V')^K$, and the functor $j_K$ is exact
(Proposition \ref{piK}), no subquotient of $V'/W$ contains a non-zero vector fixed by $K$. Since any
subquotient of $V'/W$ is also a subquotient of $V$, the above implies that $V'=W$.
\end{proof}


\begin{proof}[Proof of the theorem]
Suppose $V$ is generated by $v_1,\ldots,v_m$ and is admissible. There exists a compact open subgroup $K$
of $G$, contained and normal in $K_0$ and admitting an Iwahori decomposition with respect to the standard
parabolic subgroups such that $v_i\in V^K$ for all $i=1,\ldots,m$. Thus $V$ is generated by $V^K$. Since
$V$ is admissible, $\dim V^K$ is finite. If
\[
\{0\}\subsetneq  V_1\subsetneq  V_2 \subsetneq \cdots  \subsetneq V_s=V
\]
is a filtration by submodules, by the corollary above,
\[
\{0\}\subsetneq  V_1^K\subsetneq  V_2^K \subsetneq \cdots  \subsetneq V_s^K=V^K,
\]
and thus the length of $V$ is at most $\dim V^K$.

Conversely, if $(\pi,V)$ is of finite length, an argument by induction on the length easily shows that
$(\pi,V)$ is admissible. Indeed, if $(\pi,V)$ is irreducible, the assertion is Theorem \ref{admirr}. For
a length $n\geq 2$, we use the exactness of the functor $V\mapsto V^K$ and the induction hypothesis. On
the other hand, it is easy to see (for example by induction) that a representation of finite length is
finitely generated.
\end{proof}

\subsection{$r_P^G$ preserves representations of finite length}\label{rPGlongfin}

This is a consequence of Howe's theorem \ref{ThmHowe}, Proposition \ref{finiIJ} and Corollary
\ref{LemmeJacquet}.

\section{Bernstein decomposition}

\subsection{Variety structure of $\Omega(G)$}\label{var0G}

We begin by slightly reformulating a part of Theorem \ref{compindsup}.

\needspace{126pt}
\begin{lemme}
Let $P_i=M_iN_i$, $i=1,2$, be parabolic subgroups of $G$, $(M_i,(\rho_i,W_i))$ cuspidal data and
$(\pi,V)$ an irreducible smooth representation of $G$, such that $(\rho_i,W_i)$ is a composition factor
of $r_{P_i}^G(\pi,V)$. Then the cuspidal data $(M_1,(\rho_1,W_1))$ and $(M_2,(\rho_2,W_2))$ are
associated.
\end{lemme}
\begin{proof}
First of all, $ r_{P_i}^G(\pi,V)$ is of finite length by \ref{rPGlongfin} and admissible by Corollary
\ref{LemmeJacquet}. By Lemma \ref{pisc}, for each $i$, $(\rho_i,W_i)$ is a quotient of
$ r_{P_i}^G(\pi,V)$, and Frobenius reciprocity gives
\[
\{0\}\neq \Hom_{M_i} ( r_{P_i}^G(\pi,V), (\rho_i,W_i)) \simeq
\Hom_G((\pi,V),
i_{P_i}^G(\rho_i,W_i)).
\]
All the $i_{P_i}^G(\rho_i,W_i)$ therefore have a composition factor (a subrepresentation) in common, and
we can then conclude by Theorem \ref{compindsup}.
\end{proof}

\bigskip

This lemma allows us to define a map
\begin{align}\label{defSc}
 \mathbf{Sc} \colon &\mathbf{Irr}(G) \rightarrow \Omega(G)\\
\nonumber  (\pi,V)&\mapsto
(M,(\rho,W))_G.
\end{align}
where $(M,(\rho,W))$ is a cuspidal datum for which there exists $P\in \caP(M)$, such that $(\rho,W)$ is a
composition factor of $r_P^G(\pi,V)$ ($\mathbf{Sc}$ for ``cuspidal support'').

\begin{prop}
The map $\mathbf{Sc}$ is a surjection with finite fibers.
\end{prop}
\begin{proof}
The fact that $ \mathbf{Sc}$ is a surjection follows from Frobenius reciprocity, by considering for
$(M,(\rho,W))_G$ an irreducible subrepresentation of the induced representation.

Let $P \in \caP(M)$ and $(M,(\rho,W))$ be a cuspidal datum. We want to show that there exists only a
finite number of representations $(\pi,V) $ in $ \mathbf{Irr}(G)$ such that $(\rho,W)$ is a composition
factor of $r_P^G(\pi,V)$. Since $r_P^G(\pi,V)$ is admissible (Corollary \ref{LemmeJacquet}), we can, by
Lemma \ref{pisc}, realize $(\rho,W)$ as a quotient of $r_P^G(\pi,V)$, and thus
\[
\Hom_M(r_P^G(\pi,V), (\rho,W))=\Hom_G((\pi,V), i_P^G (\rho,W)) \neq \{0\}.
\]
Now $i_P^G (\rho,W)$ is of finite length, by Lemma \ref{longfin}.
\end{proof}

\bigskip

We define an equivalence relation weaker than being associated on the set of cuspidal data. We say that
two cuspidal data $(M_1,(\rho_1,W_1))$ and $(M_2,(\rho_2,W_2))$ define the same inertial support if there
exists \index[ter]{inertial support} $g\in G$ and an unramified character $\omega$ of $M_2$ such that
\[
gM_1g^{-1} = M_2 \quad \text { and }  \quad \rho_2 \simeq \rho_1^g
\otimes \omega
\]
We denote by $[M,(\rho,W)]_G$ \index[not]{M(rho,W)G@$[M,(\rho,W)]_G$} the class of $(M,(\rho,W))$ for
this equivalence relation, and $\caB(G)$ \index[not]{B(G)@$\caB(G)$} the set of equivalence classes.

\begin{exemple}
The inertial classes of supercuspidal representations of $G$ define inertial supports. With the notation
of \ref{clinertcusp}, if $(\pi,V)$ is a supercuspidal representation of $G$, we have
$[\pi]=[G,(\pi,V)]_G=\mathbf{Irr}(G)_{[\pi]}$. We have seen that this set is equipped with the structure
of an affine algebraic variety. We will extend this to the other inertial supports.
\end{exemple}

\begin{thm}
The set of cuspidal data $\Omega(G)$ is equipped with the structure of an algebraic variety, whose
connected components, indexed by $\caB(G)$, are quotients of complex tori by the action of finite groups.
\end{thm}
\begin{proof}
The first step of the proof is to note that for any Levi subgroup $M$ of $G$, the set
$\mathbf{Irr}(M)_{sc}$ is an algebraic variety. Indeed, it is the disjoint union of its inertial classes,
which we know admit an algebraic variety structure studied in detail in \ref{centreMpi}. Let us write
this
\[
\mathbf{Irr}(M)_{sc}=\coprod_i D_i, \quad D_i=[\rho_i],
\]
with algebra of functions $\prod_i \caZ_{[\rho_i]}$. The group $W(A)=N_G(M)/M$ acts naturally on this
variety.

The quotient of $\mathbf{Irr}(M)_{sc}$ by the action of $W(A)$ is obtained first by identifying two
inertial classes $D_1$ and $D_2$ (i.e., two connected components) of $\mathbf{Irr}(M)_{sc}$ when there
exists $w \in W(A)$ such that $w\cdot D_1=D_2$. Let $D$ be an inertial class of $\mathbf{Irr}(M)_{sc}$,
and let $W(D)$ be the stabilizer of $D$ in $W(A)$. Let $(\rho,W)$ be a supercuspidal representation of
$M$ in the inertial class $D$. Let us now describe the variety structure of the connected component of
$\Omega(G)$ containing $(M,(\rho,W))_G$: it is the quotient of the variety $D$ (a complex torus) by the
action of $W(D)$. It is therefore again an algebraic variety. Note that in this description, the choice
of $D$ is not unique. Now consider a system of representatives $M_1,\ldots ,M_r$ of the conjugacy classes
of Levi subgroups of $G$. Set $W_i = N_G(M_i)/M_i$, $i=1,\ldots, r$. It is obvious from the definitions
that
\[
\Omega(G)=\coprod_{i=1}^r   \mathbf{Irr}(M_i)_{sc}/W_i.
\]

We can parameterize the connected components of $\Omega(G)$ in another way. Consider the canonical
projection from $\Omega(G)$ to $\caB(G)$. It is then clear from the above that the connected components
of $\Omega(G)$ are the fibers of this projection.
\end{proof}

\bigskip

Let us define the map \index[not]{Si@$ \mathbf{Si}$}
\[
\mathbf{Si} \colon \mathbf{Irr}(G) \rightarrow \caB(G)
\]
obtained by composing $\mathbf{Sc}$ with the projection from $\Omega(G)$ to $\caB(G)$ ($\mathbf{Si}$ for
``inertial support''). The fibers of this map provide a partition of $\mathbf{Irr}(G)$. Let
$\mathbf{Irr}(G)_\frs$ \index[not]{Irr(G)_sc@$\mathbf{Irr}(G)_\frs $} denote the fiber over an element
$\frs$ of $\caB(G)$. If $\Omega$ is the connected component of $\Omega(G)$ corresponding to $\frs$, we
also denote $\mathbf{Irr}(G)_\frs=\mathbf{Irr}(G)_\Omega $.
\index[not]{Irr(G)_Omega@$\mathbf{Irr}(G)_\Omega$}

\subsection{The Bernstein decomposition theorem}\label{dcompthm}

\index[ter]{decomposition!Bernstein} \index[ter]{decomposition theorem} We will now show that the
decomposition $\mathbf{Irr}(G) =\coprod  \mathbf{Irr}(G)_\Omega$ of the previous paragraph induces a
decomposition of the category $\caM(G)$.

If $\Omega$ is a connected component of $\Omega(G)$ and if $(\pi,V)$ is a smooth representation, we
denote by $V(\Omega)$ \index[not]{V(Omega)@$V(\Omega)$} the maximal subrepresentation of $V$ such that
all irreducible subquotients of $V(\Omega)$ are in $\mathbf{Irr}(G)_\Omega$. We say that $(\pi,V)$ is
split \index[ter]{split!(representation)} with respect to $\Omega(G)$ if $V=\bigoplus_\Omega V(\Omega)$
(see \ref{deccat}).

\begin{thm}
Any representation $(\pi,V)$ in $\caM(G)$ is split with respect to $\Omega(G)$. The category $\caM(G)$
decomposes into a product of categories:
\[
\caM(G)= \prod_\Omega   \caM(G)_\Omega.
\]
where $\caM(G)_\Omega$ is the full subcategory of $\caM(G)$ of \index[not]{M(G)_Omega@$\caM(G)_\Omega$}
representations all of whose irreducible subquotients are in $\mathbf{Irr}(G)_\Omega$.
\end{thm}

\begin{proof}[Proof of the theorem]
The proof is done in several steps.

\begin{lemme}
If $V'$ is a subrepresentation of $V$, and $V$ is split with respect to $\Omega(G)$, then $V'$ is split
with respect to $\Omega(G)$.
\end{lemme}
\begin{proof}
This is a special case of Lemma \ref{deccat}.
\end{proof}

The proof of the theorem now consists of showing that any representation of $\caM(G)$ embeds into a
representation split with respect to $\Omega(G)$. We can then use the lemma to conclude.

Let us define $\mathrm{Cusp}(G)$ \index[not]{Cusp@$\mathrm{Cusp}(G)$} as the product over the standard
parabolic subgroups $P=MN$ of $G$ of the categories $\caM(M)_{sc}$. Then define the functors
\[
I \colon \mathrm{Cusp}(G) \rightarrow \caM(G),\quad (\rho_M,W_M)_P \mapsto \bigoplus_P \; i_P^G (\rho_M,W_M)
\]
and
\[
R \colon  \caM(G) \rightarrow   \mathrm{Cusp}(G),\quad  (\pi,V) \mapsto (r_P^G(\pi,V)_{sc})_P
\]
where by $r_P^G(\pi,V)_{sc}$ we mean the cuspidal part of the representation $r_P^G(\pi,V)$.

\begin{lemme}
$(i)$ The functor $R$ is the left adjoint of the functor $I$.

$(ii)$ $R$ is exact, faithful, and preserves finitely generated representations. The functor $I$ is exact
and faithful.

$(iii)$ For any representation $(\pi,V)$ in $\caM(G)$, the adjunction morphism $\alpha_V$ from $(\pi,V)$
to $IR (\pi,V)$ is an injection.
\end{lemme}

\begin{proof}
Let us calculate, for any object $(\rho_M,W_M)_P$ of $\mathrm{Cusp}(G)$, and any representation $(\pi,V)$
of $\caM(G)$,
\begin{align*}
 \Hom_{\mathrm{Cusp}(G)}(R(\pi,V), (\rho_M,W_M)_P)
\simeq & \prod_P \Hom_{M}( r_P^G(\pi,V)_{sc}, (\rho_M,W_M))\\
                                   \simeq  \prod_P \Hom_{M}( r_P^G(\pi,V), (\rho_M,W_M))
                                   \simeq &\prod_P \Hom_{G}( (\pi,V),i_P^G(\rho_M,W_M))\\
                                    \simeq \Hom_{G}( (\pi,V), \oplus_{P} i_P^G(\rho_M,W_M))
                                    \simeq &\Hom_{G}( (\pi,V), I((\rho_M,W_M)_P))
\end{align*}
We used the adjunction of the functors $i_P^G$ and $r_P^G$ for any parabolic subgroup $P$ of $G$, and the
fact that there are only finitely many standard parabolic subgroups. Since all these isomorphisms are
natural, this shows $(i)$. Similarly, the exactness of $R$ and $I$ follows from the exactness of the
functors $i_P^G$ and $r_P^G$, and the fact that the projection $\caM(M) \rightarrow  \caM(M)_{sc}$ also
defines an exact functor.

To show that $R$ is faithful, it suffices to show that if $(\pi,V)$ is non-zero in $\caM(G)$, then
$R(\pi,V) \neq 0$ (see Lemma \ref{exaCt}). Let $P$ be minimal among the standard parabolic subgroups of
$G$ such that $r_P^G(\pi,V)\neq 0$. Then $r_P^G(\pi,V)$ is cuspidal (the argument was used in the proof
of Corollary \ref{supercusp}), and thus $R(\pi,V) \neq 0$. The proof of the fact that $I$ is faithful is
similar. Indeed, each functor $i_P^G$ is faithful, since clearly, if $(\rho,W)$ is a non-zero
representation of $M$, $i_P^G(\rho,W)$ is non-zero. If $(\pi,V)$ is a finitely generated representation
of $\caM(G)$, we have seen that $r_P^G(\pi,V)$ is finitely generated (Proposition \ref{finiIJ}). The
assertion about $R$ follows.

We now show $(iii)$. Let $V'$ denote the kernel of $\alpha_V$ and $i$ the inclusion of $V'$ in $V$. Since
the adjunction morphism is natural, we have a commutative diagram:
\[
\begin{CD}
V' @>i>>  V \\
@V\alpha_{V'}VV     @V{\alpha_V}VV\\
IR(V') @>{IR(i)}>> IR(V)
\end{CD}
\]
Since $0=\alpha_V \circ i = IR(i)\circ \alpha_{V'}  $, and by the exactness of the functors $I$ and $R$,
$IR(i)$ is a monomorphism, we have $ \alpha_{V'}=0  $. Now the morphism $\alpha_{V'}$ corresponds by the
adjunction isomorphism to the identity of $R(V')$. This shows that $R(V')=0$, and thus that $V'=0$, the
functor $R$ being faithful.
\end{proof}

\bigskip

This lemma shows that for any representation $(\pi,V)$ in $\caM(G)$, we can embed $(\pi,V)$ into a
representation of the form
\[
\bigoplus_P \; i_P^G (\rho_M,W_M),
\]
where the sum is over the standard parabolic subgroups $P=MN$ of $G$ and $(\rho_M,W_M)$ is a
supercuspidal representation of the Levi subgroup $M$. To complete the proof of the theorem, it therefore
suffices to prove that the representations $i_P^G (\rho_M,W_M)$ are split with respect to $\Omega(G)$. We
will now work with only one index $M$, and we simplify the notation by setting $(\rho,W)=(\rho_M,W_M)$.
Since $(\rho,W)$ is supercuspidal, we can write it $W=\bigoplus_D W_D$, where $D$ runs through the
inertial classes of supercuspidal representations of $M$. We have,
\begin{equation}\label{DWD}
i_P^G(\oplus_D W_D)\simeq   \bigoplus_D \; i_P^G(W_D).
\end{equation}
Note that if we restrict the sum to a finite number of indices, the isomorphism is obvious. This shows in
particular that the right-hand side injects into the left-hand side. Let us realize the left-hand side in
the usual way as a function space on which $G$ acts by right translation. Let $f$ be a function in this
space. Then $f$ is fixed by a certain compact open subgroup $K$ of $G$, and thus $f$ is completely
determined by its values on a system of representatives of the double cosets $P\backslash G/K$. Since
this set is finite, it follows that the image of $f$ is contained in a finite sum of subspaces $W_D$ and
thus $f$ is in a finite sum of subspaces $ i_P^G(W_D)$. This shows that the injection of the right-hand
side into the left-hand side is also a surjection, hence an isomorphism.

It remains to show that $i_P^G (W_D)$ is split with respect to $\Omega(G)$, but this is immediate from
the definitions. Indeed, let $\Omega$ be the component of $\Omega(G)$ which is a quotient of the variety
$D$. If $(\sigma,E)$ is a composition factor of $i_P^G (W_D)$, then $r_P^G(\sigma,E)$ is a subquotient of
$r_P^Gi_P^G (W_D)$ by exactness of the functor $r_P^G$. By Proposition \ref{geolemmacusp} $(iii)$ and
(\ref{WMWA}), $r_P^Gi_P^G (W_D)$ admits a filtration whose successive quotients are of the form
${}^w(W_D)$, $w \in W(A)=N_G(M)/M$. Any irreducible subquotient of $r_P^G(\sigma,E)$ is therefore a
subquotient of a ${}^w(W_D)$, $w \in W(A)=N_G(M)/M$, which means that the cuspidal support of
$(\sigma,E)$ is in $\Omega$.

The connected components of $\Omega(G)$ being parameterized by the inertial classes of cuspidal data, we
can reformulate the first assertion of the theorem by saying that any representation $(\pi,V)$ of
$\caM(G)$ splits into a direct sum
\begin{equation}\label{VVs}
  V=\bigoplus_{\frs \in \caB(G)} V_\frs,
\end{equation}
where each $V_\frs$ is a smooth representation of $G$ all of whose irreducible subquotients map to $\frs$
by the map $\mathbf{Si}$. If $\frs$ and $\frt$ are distinct elements of $\caB(G)$, we deduce that for any
$(\pi,V)$ and $(\rho,W)$ in $\caM(G)$, we have $\Hom_G(V_\frs,W_\frt)=0$. This shows that any
$G$-morphism from $V$ to $W$ must map $V_\frs$ to $W_\frs$, and thus we have a natural isomorphism
\[
\Hom_G(V,W)=\bigoplus_{\frs \in \caB(G)} \Hom_G(V_\frs,W_\frs)  .
\]
The decomposition (\ref{VVs}) of the objects of $\caM(G)$ is preserved by morphisms, which shows that we
have a decomposition of categories
\begin{equation}\label{MGMGs}
\caM(G)=\prod_{\frs \in \caB(G)} \caM(G)_\frs.
\end{equation}
This completes the proof of the Bernstein decomposition theorem.
\end{proof}

\bigskip

Let $\Omega$ be a connected component of $\Omega(G)$. From the preceding proof, we can deduce other
characterizations of the representations of the subcategory $\caM(G)_\Omega$ of $\caM(G)$. Note in
passing that since $\caM(G)_\Omega$ is stable under passing to subquotients, it is an abelian category.
Fix a (standard) Levi subgroup $M$ of $G$ and $D=[\rho]$ an inertial class of irreducible supercuspidal
representations of $M$ such that $(M,\rho)_G \in \Omega$. Let $P$ be the standard parabolic subgroup of
$G$ containing $M$.

\begin{prop}
A representation $(\pi,V)$ is in $\caM(G)_\Omega$ if and only if one of the following equivalent
conditions is satisfied:

$(i)$ $(\pi,V)$ embeds into a sum of representations of the form $i_P^G (\tau,E)$, with
$(\tau,E)\in \caM(M)_D$.

$(ii)$ $(\pi,V)$ is a subquotient of a sum of representations as in $(i)$.

$(iii)$ For any standard parabolic subgroup $Q=LU$ of $G$ and inertial class $D'$ of irreducible
supercuspidal representations of $L$, such that $(M,D)$ and $(L,D')$ are not conjugate, the component in
$\caM(L)_{D'}$ of $r_Q^G(\pi,V)$ is zero.
\end{prop}

\begin{proof}
As previously stated, apart from $(ii)$, this is a reformulation of the proof of the theorem. Indeed, the
decomposition of a representation $(\pi,V)\in \caM(G)$ was obtained in the following way: we constructed
an embedding
\[
(\pi,V) \rightarrow IR(\pi,V),
\]
where $R(\pi,V)$ is the sum over the standard parabolic subgroups $Q=LU$ of the supercuspidal components
of $r_Q^G(\pi,V)$. Each of these components, in turn, decomposes into a direct sum of components
according to the inertial classes of supercuspidal representations $D'$ of $L$. We thus have $R(\pi,V)$
of the form
\[
R(\pi,V)=\bigoplus_Q \bigoplus_{D'} (\tau_{Q,D'},E_{Q,D'})
\]
where $(\tau_{Q,D'},E_{Q,D'}) \in \caM(L)_{D'}$. From this we deduce that $(\pi,V)$ is a
subrepresentation of
\[
\bigoplus_Q  i_Q^G(\oplus_{D'}(\tau_{Q,D'},E_{Q,D'}))=\bigoplus_Q
\bigoplus_{D'}i_Q^G (\tau_{Q,D'},E_{Q,D'})
\]
(the equality between the two sides was seen in \ref{DWD}). We then group the pairs $(Q,D')$ conjugate
under $G$, i.e., those defining the same connected component $\Omega$ of $\Omega(G)$. We obtain an
embedding
\[
(\pi,V) \hookrightarrow \bigoplus_\Omega \left(  \bigoplus_{(Q,D')|
[L,D']_G=\Omega } i_Q^G (\tau_{Q,D'},E_{Q,D'}) \right).
\]
The right-hand side is a representation split with respect to $\Omega(G)$, the decomposition appearing
explicitly. The decomposition of $(\pi,V)$ is obtained by taking the sum of the intersections of $V$ with
the factors of the right-hand side. A factor
\[
\bigoplus_{(Q,D') \mid   [L,D']_G=\Omega } i_Q^G (\tau_{Q,D'},E_{Q,D'})
\]
can be rewritten $i_P^G(\tau,E)$ with $(\tau,E)\in \caM(M)_D$ because the data $(M,D)$ and $(L,D')$ are
conjugate under $G$. This discussion shows that $(\pi,V)$ is in a single component $\caM(G)_\Omega$ if
and only if $(i)$ or $(iii)$ is satisfied.

We will now show the equivalence with $(ii)$. It is clear that $(i)\Rightarrow (ii)$. Suppose that
$(\pi,V)$ satisfies $(ii)$. Since the functor $r_Q^G$ is exact, for every standard parabolic subgroup
$Q$, $r_Q^G(\pi,V)$ is a subquotient of a representation of the form $r_Q^G (\oplus i_P^G (\tau,E))$,
where the $(\tau,E)$ are in $\caM(M)_D$. We then use Proposition \ref{geolemmacusp} to conclude that
$(\pi,V)$ satisfies $(iii)$.
\end{proof}

\subsection{Bernstein decomposition, induction and restriction} \label{BIR}

Let $P=MN$ be a parabolic subgroup of $G$. A cuspidal datum $(L,(\sigma,E))$ of $M$ is also a cuspidal
datum for $G$, since a Levi subgroup of $M$ is also a Levi subgroup of $G$, and we thus have a map:
\[
i_{MG}\colon \Omega(M)\rightarrow \Omega(G), \quad
(L,(\sigma,E))_M \mapsto  (L,(\sigma,E))_G
\]
This map induces another one, denoted in the same way, between sets of inertial supports:
\[
i_{MG}\colon \caB(M)\rightarrow \caB(G), \quad
[L,(\sigma,E)]_M \mapsto  [L,(\sigma,E)]_G
\]
\index[not]{i_MG@$ i_{MG}$}

\begin{prop}
Let $\frs \in \caB(M)$ and let $(\rho,W)$ be in $\caM(M)_\frs$. Then $i_P^G(\rho,W)$ is in
$\caM(G)_{i_{MG}(\frs)}$. Let $\frt \in \caB(G)$ and let $(\pi,V)$ be in $\caM(G)_\frt$. Then
$r_P^G(\pi,V)$ is in $\prod_{\frs \in i_{MG}^{-1}(\frt)} \caM(M)_\frs$.
\end{prop}

\begin{proof}
This follows directly from the characterizations of the representations of $\caM(M)_\frs$ and
$\caM(G)_{i_{MG}(\frs)}$ obtained in the previous section.
\end{proof}

\subsection{Central idempotents}\label{IdCent}

The decomposition theorem allows us to write the category $\caM(G)$ as a product of smaller categories,
parameterized by the connected components of the algebraic variety $\Omega(G)$. For any connected
component $\Omega$ of $\Omega(G)$, the projection
\[
e_\Omega \colon  \caM(G) \rightarrow \caM(G)_\Omega
\]
\index[not]{e_Omega@$e_\Omega$} defines an idempotent of the center of the category $\caM(G)$ and
\[
\Id_{\caM(G)}= \sum_{\Omega} e_\Omega.
\]
If the inertial support of $\Omega$ is $\frs \in \caB(G)$, we also denote
$e_\frs=e_\Omega$.\index[not]{e_s@$e_\frs$}

The Hecke algebra $\caH(G)$ therefore decomposes into
\[
\caH(G)=\bigoplus_{\frs \in \caB(G)} e_\frs*\caH(G)*e_\frs.
\]
Each $e_\frs*\caH(G)*e_\frs$ is a two-sided ideal of $\caH(G)$. The category $\caM(G)_\frs$ is equivalent
to the category $\caM(e_\frs*\caH(G)*e_\frs)$ of non-degenerate $e_\frs*\caH(G)*e_\frs$-modules.

\begin{rmq}
The idempotent $e_\frs$ is not an element of the Hecke algebra $\caH(G)$, but of its completion ({\sl
cf.} \ref{completeA}) $\overline{\caH(G)}$.
\end{rmq}

It remains to describe the center of each category $\caM(G)_\Omega$. When we take a connected component
$\Omega$ of $\Omega(G)$ given by an inertial class of supercuspidal representations of $G$, we were able
to describe the corresponding category as a category of unital modules over a unital algebra, and
describe the center of this category (i.e., the center of the unital algebra in question), as the algebra
of functions on $\Omega$. We must obtain similar results for an arbitrary connected component $\Omega$.
The first step consists of finding a progenerator for these categories. Suppose that $\Omega$ is given by
a Levi subgroup $M$ and an inertial class $D$ of supercuspidal representations of $M$. The naive idea is
to obtain this progenerator by inducing from $M$ to $G$ the progenerator of $\caM(M)_D$. The difficulty
that appears in this approach is to show that we indeed obtain in this way a finitely generated
projective object. The proof of this fact is delicate, and uses the generalization of Jacquet's lemma to
non-admissible representations, as well as the Noetherian properties of the category $\caM(G)$.

\subsection{Noetherianity of $\caM(G)$}\label{noether}

In this section, we show the following result which completes those obtained in \ref{thmsfinit}: let
$P=MN$ be a parabolic subgroup of $G$, then the functor $i_P^G$ preserves finitely generated
representations. This involves Noetherian properties of the category $\caM(G)$ coming from the
decomposition theorem \ref{dcompthm}.

In an abelian category, an object is Noetherian if any increasing sequence of subobjects is stationary
({\sl cf.} Appendix \ref{sousobj}). An abelian category such that any finitely generated object is
Noetherian will be called Noetherian. For example, the category $A-\mathbf{mod}$, where $A$ is a
Noetherian ring, is Noetherian. \index[ter]{Noetherian (category)}
\begin{prop}
The category $\caM(G)$ is Noetherian.
\end{prop}

\begin{proof}
Let us first consider a supercuspidal component $\caM(G)_\frs$ of the decomposition \ref{dcompthm}, say
$\frs=[G,(\rho,W)]_G$ where $(\rho,W)$ is thus an irreducible supercuspidal representation of $G$. We saw
in Corollary \ref{centreA} that the category $\caM(G)_\frs=\caM(G)_{[\rho]}$ is equivalent to the
category of right unital modules over a certain Noetherian $\bbC$-algebra. It is therefore a Noetherian
category.

Let us return to the functor $R\colon \caM(G)\rightarrow \mathrm{Cusp}(G)$ of Section \ref{dcompthm}.
Since each category $\caM(M)_{sc}$ which constitutes $\mathrm{Cusp}(G)$ decomposes (in the sense of
\ref{deccat}) according to \ref{decMG2} into a product $\caM(M)_{sc}=\prod_D \caM(M)_D$ where $D$ runs
through the inertial classes of supercuspidal representations of $M$, we can view $R$ as a functor
$R \colon  \caM(G)\rightarrow \prod_{(M,D)} \caM(M)_D$ where $M$ runs through the standard Levi subgroups
of $G$ and $D$ the inertial classes of supercuspidal representations of $M$. Since each category
$\caM(M)_D$ is Noetherian from the above, the same is true for their product: indeed, a finitely
generated object has only a finite number of non-zero components in $\prod_{(M,D)} \caM(M)_D$ since by
definition, this holds for a finite system of generators, and the assertion easily follows. Suppose that
$(\pi,V)$ is a finitely generated representation of $G$ and let
\[
V_1 \subset V_2 \subset \cdots \subset V_n \subset \cdots
\]
be a sequence of subrepresentations of $V$. Since $R$ is faithful and exact (\ref{dcompthm}), if this
sequence is not stationary,
\[
R(V_1) \subset R(V_2) \subset \cdots \subset R(V_n) \subset \cdots
\]
is not either. But since we know on the other hand that $R(V)$ is finitely generated, this cannot happen
since $\mathrm{Cusp}(G)$ is Noetherian.
\end{proof}

\begin{thm}
Let $P=MN$ be a parabolic subgroup of $G$. Then the functor $i_P^G$ preserves finitely generated
representations.
\end{thm}

\begin{proof}
We can of course without loss of generality assume that $P$ is a standard parabolic subgroup of $G$. Let
$(\rho,W)$ in $\caM(M)$ be a finitely generated representation. By the proposition above, $(\rho,W)$ is
Noetherian. We will show that the representation $i_P^G (\rho,W)$ is Noetherian. This implies of course
that $i_P^G (\rho,W)$ is finitely generated.

Suppose first that $(\rho,W)$ is supercuspidal. Let us apply the functor $R$ to $ i_P^G (\rho,W)$. By
Proposition \ref{geolemmacusp}, only the standard parabolic subgroups $Q=LU$ such that $L$ and $M$ are
conjugate contribute. For such a standard parabolic subgroup $Q$, the same proposition asserts that
$r_Q^G\circ i_P^G (\rho,W)$ has a filtration whose successive quotients are the $w\cdot (\rho,W)$,
$w \in W(L,M)/W_L$, which are finitely generated supercuspidal representations and therefore Noetherian.
It follows that $r_Q^G\circ i_P^G (\rho,W)$ is Noetherian, and thus the same holds for
$R( i_P^G (\rho,W))$. Since the functor $R$ is exact and faithful, it does not annihilate any of these
subquotients $w\cdot (\rho,W)$. We deduce by an argument given in the proof of the proposition above that
$i_P^G (\rho,W)$ is Noetherian.

In the general case, the second lemma of \ref{dcompthm} exhibits an embedding
\[
(\rho,W) \hookrightarrow IR(\rho,W)= \bigoplus_{P'} i_{P'}^M (\rho_{M'},W')
\]
where the sum is over the standard parabolic subgroups $P'=M'N'$ of $M$ and the $(\rho_{M'},W')$ are
supercuspidal representations. We obtain by the exactness of the functor $i_P^G$ and the transitivity
property \ref{associa} an embedding
\[
i_P^G (\rho,W) \hookrightarrow IR(\rho,W)= i_P^G(\oplus_{L}
i_{P'}^M \rho_{M'})=\bigoplus_{P'}  i_{P'}^G (\rho_{M'}).
\]
Each factor of the right-hand side is Noetherian, from the above, and we easily deduce that
$i_P^G (\rho,W)$ is as well.
\end{proof}

\section{Algebraic families of representations}

\subsection{$(G,B)$-modules} \label{GBmod}

Let $B$ be a reduced Noetherian finitely generated commutative $\bbC$-algebra. We refer the reader to
\cite{BourAC1} for the notion of a flat module.

\begin{defi}
A $(G,B)$-module is a $\bbC$-vector space $V$ equipped with a smooth representation $(\pi,V)$ of $G$ and
a flat $B$-module structure, where the two actions commute. Let $\caM(G,B)$ denote the category of
$(G,B)$-modules. If $K$ is a compact open subgroup of $G$, the space $V^K$ is again a $B$-module. A
$(G,B)$-module $V$ is admissible if for any compact open subgroup $K$ of $G$, $V^K$ is finitely generated
over $B$.
\end{defi}
\index[ter]{GBmodule@$(G,B)$-module} \index[ter]{admissible!($(G,B)$-module)}

\begin{prop}
Let $V$ be a $\bbC$-vector space equipped with a smooth representation $(\pi,V)$ of $G$ and a $B$-module
structure, where the two actions commute. Then $V$ is an admissible $(G,B)$-module if and only if for any
compact open subgroup $K$ of $G$, $V^K$ is a finitely generated projective $B$-module.
\end{prop}
\begin{proof}
According to \cite{BalX}, Chap. X, \S 1, no.~5, corollary of Proposition 8, and no.~4, Proposition 5
$(iii)$, there is an equivalence between finitely presented flat modules and finitely generated
projective modules. Since $B$ is Noetherian, a finitely generated module is finitely presented
(\cite{BalX}, Chap. X, \S 1, no.~4, Proposition 5 $(ii)$), so a finitely generated $B$-module is flat if
and only if it is projective. On the other hand, for any $B$-module $E$, the natural morphism
\[
E\otimes_B V^K \rightarrow (E \otimes_B V)^K
\]
is an isomorphism whose inverse is given by
\[
e\otimes v \mapsto e\otimes e_K\cdot v.
\]
Finally, the functor $j_K:V \mapsto V^K$ is exact, so if $V$ is a flat $B$-module, i.e., if the functor
$\bullet \otimes_B V$ is exact, then the functor $$j_K(\bullet \otimes_B V)\simeq \bullet \otimes_B V^K$$
is exact, which implies that $V^K$ is flat. This shows that if $V$ is an admissible $(G,B)$-module, then
$V^K$ is a finitely generated projective module. Conversely, if the functor $$\bullet \otimes_B V^K$$ is
exact for any compact open subgroup $K$ of $G$, then for any short exact sequence of $B$-modules
\[
0\rightarrow E_1\rightarrow E_2 \rightarrow E_3\rightarrow 0,
\]
we still have an exact sequence
\[
0\rightarrow E_1 \otimes_B V^K\rightarrow E_2 \otimes_B V^K
\rightarrow E_3 \otimes_B V^K \rightarrow 0.
\]
By taking $K$ smaller and smaller, we obtain, since the representations $E_i$ of $G$ are smooth
\[
0\rightarrow E_1 \otimes_B V\rightarrow E_2 \otimes_B V
\rightarrow E_3 \otimes_B V \rightarrow 0,
\]
which shows that $V$ is flat.
\end{proof}

\begin{rmq}
At the end of the proof, we used the following argument, which will be useful to us again: if $V$ is a
$\bbC$-vector space equipped with a smooth representation $(\pi,V)$ of $G$ and a $B$-module structure,
where the two actions commute, then $V$ is flat as a $B$-module if and only if all the $V^K$ are, as $K$
runs through the set of compact open subgroups of $G$.
\end{rmq}

A reduced Noetherian, finitely generated and commutative $\bbC$-algebra $B$ is the algebra of polynomial
functions of an affine algebraic variety $S$ ($S$ is identified with the maximal spectrum
$\mathrm{SpecMax} (B)$ of $B$). \index[ter]{maximal spectrum} \index[not]{SpecMax@$\mathrm{SpecMax} (B)$}
Recall that by definition the maximal spectrum of $B$ is the set of unital algebra morphisms
$\chi: B\rightarrow \bbC$. If $\chi \in S$, we denote by $\Psi_\chi \in \mathrm{SpecMax} (B)$ the
corresponding algebra morphism, $I_\chi$ its kernel (a maximal ideal of $B$), and we define
\[
V_\chi = V\otimes_B B/I_\chi
\]
the specialization \index[ter]{specialization} of $V$ at $\chi$. This remains a smooth representation of
the group $G$. We can view the $(G,B)$-module $(\pi,V)$ as an algebraic family (indexed by
$S\simeq \mathrm{SpecMax} (B)$) of smooth representations of $G$. We denote by $\mathbf{sp}_\chi$ the
canonical morphism from $V$ to $V_\chi$. Since $B/I_\chi \simeq \bbC$ as a vector space, we also denote
by $\bbC_{\Psi_\chi}$ the $B$-module $\bbC$ equipped with the action
\[
b\cdot z= \Psi_\chi(b) z , \quad (b\in B,\, z\in \bbC).
\]
Thus $V_\chi= V\otimes_B \bbC_{\Psi_\chi}$.\index[not]{V_chi@$V_\chi$}

\begin{lemme}
If $V$ is an admissible $(G,B)$-module, then for all $\chi \in S$, $V_\chi$ is an admissible
representation of $G$.
\end{lemme}
\begin{proof}
If $V^K$ is generated by $v_1,\ldots,v_r$ as a $B$-module, then
\[
(V\otimes_B \bbC_\chi)^K=V^K\otimes_B \bbC_\chi
\]
is generated by $v_1\otimes 1,\ldots,v_r\otimes 1$ as a $\bbC$-vector space.
\end{proof}

\begin{exemple}
Let $(\pi,V)$ be a smooth representation of $G$. Let $B=\bbC[\Lambda(G)]$ denote the algebra of
polynomials on the variety of unramified characters $\caX(G)$ of $G$. Recall that evaluation at the
points $g$ of $G$ defines a morphism
\[
\chi_{un}: G \rightarrow B^\times.
\]
(For all this see \ref{varXG}.)

Set, for all $v \in V$, $b \in B$, $g \in G$
\[
V_B= V\otimes_\bbC B, \quad \pi_B(g)(v \otimes b)= \pi(g)\cdot v \otimes
\chi_{un}(g)b.
\]
Let us verify that this is a $(G,B)$-module, with $B$ acting by multiplication on the second factor. On
the one hand, the action of $G$ and that of $B$ commute: for any $v \in V,\, b,b' \in B, \, g \in G$,
\[
b'\cdot (\pi_B(g)\cdot(v\otimes b))=b'\cdot (\pi(g)\cdot v\otimes
\chi_{un}(g)b)= \pi(g)\cdot v\otimes  \chi_{un}(g)b' b,
\]
and
\[
\pi_B(g) \cdot (b'\cdot(v\otimes b))= \pi_B(g)\cdot (v\otimes
b'b)= \pi(g)\cdot v\otimes \chi_{un}(g) b'b.
\]
On the other hand $V\otimes_\bbC B$ is a free $B$-module, hence flat.

If $g$ is in ${}^0G$, $\chi_{un}(g)=1$, so the action of $g$ on $B$ is trivial. This is true in
particular for any element of a compact subgroup $K$ of $G$. It follows that $V_B^K=V^K\otimes_\bbC B $,
and thus $V_B$ is an admissible $(G,B)$-module as soon as $V$ is admissible.

If $\chi \in \caX(G)$, the specialization of $\pi_B$ at $\chi$ is isomorphic to $\pi \otimes \chi$:
\[
V_B\otimes_B B/I_\chi= (V\otimes B)\otimes_B \bbC_{\Psi_\chi}=V\otimes \bbC_{\chi},
\]
where $\bbC_\chi$ is $\bbC$ viewed as the space of the representation $\chi$ of $G$. Indeed the action of
$G$ is given by
\begin{align*}
 &g \cdot (v\otimes b \otimes z)=\pi_B(g)\cdot(v\otimes b)  \otimes z=
\pi(g)\cdot v\otimes \chi_{un}(g) b \otimes z \\
&= \pi(g)\cdot v\otimes b
\otimes \Psi_\chi(\chi_{un}(g))z=\pi(g)\cdot v\otimes b
\otimes \chi(g)z\\
&=(\pi\otimes \chi)(g)\cdot v \otimes b
\otimes z.
\end{align*}
\end{exemple}
We used (\ref{Psichi}) for the penultimate equality.

\subsection{$(G,B)$-modules and induction}\label{GBind}

If $P=MN$ is a parabolic subgroup of $G$, and if $(\sigma,E)$ is an $(M,B)$-module, then $i_P^G (E)$ is
equipped with a $B$-module structure (we let $B$ act on the target space $E$ of the functions in
$i_P^G (E)$). It is obvious that the action of the group, by right translation of functions, commutes
with this action of $B$. To show that $i_P^G (\sigma,E)$ is a $(G,B)$-module, it remains to show that it
is a flat $B$-module. We use for this Remark \ref{GBmod}: it suffices to show that for any compact open
subgroup $K$ of $G$, $i_P^G (E)^K$ is a flat $B$-module. Now we described $i_P^G (E)^K$ in \ref{Ind}. Let
$\Omega$ be a system of representatives in $G$ of the double cosets $P\backslash G/K$, and for all
$ g \in \Omega$, let $K_{M}^g$ be the projection of $P\cap g Kg^{-1}$ onto $M\simeq P/N$. Recall that
since $P\backslash G$ is compact, $\Omega$ is a finite set. We then have
\begin{equation}\label{ipgek}
i_P^G (E)^K\simeq \bigoplus_{g \in \Omega} E^{K_M^g}
\end{equation}
and this isomorphism is an isomorphism of $B$-modules. Since $E$ is a flat $B$-module, all the
$E^{K_M^g}$ are flat, which shows that $i_P^G (E)^K$ is flat. This shows us that we have a functor
\[
i_P^G \colon  \caM(M,B) \rightarrow \caM(G,B).
\]
It is then easy to generalize Lemma \ref{IndAdm}:
\begin{lemme}
Let $P=MN$ be a parabolic subgroup of $G$. Then the functor
\[
i_P^G \colon \caM(M ,B) \rightarrow \caM(G,B)
\]
maps admissible modules to admissible modules.
\end{lemme}

Let $C$ be another reduced Noetherian commutative $\bbC$-algebra, and suppose that $C$ is a $B$-algebra.
The considerations of Section \ref{Oublietadjoints}, applied here to the case of $B$ and $C$ which are
commutative and unital, give us on the one hand a pseudo-forgetful functor (isomorphic to the forgetful
functor in this case)
\[
\caM(C) \longrightarrow  \caM(B), \quad E \mapsto \Hom_C(C,E)
\simeq E
\]
and its left adjoint
\[
\caM(B) \longrightarrow  \caM(C), \quad W \mapsto C\otimes_B W,
\]
which is called in this context the base change functor \index[ter]{base change} from $B$ to $C$. If $V$
is a flat $B$-module, then for any $C$-module $E$,
\[
E \otimes_C( C\otimes_B V)\simeq  ( E \otimes_C C) \otimes_B V.
\]
Since the forgetful functor from $C$ to $B$ is trivially exact, we see that $C\otimes_B V$ is a flat
$C$-module. This shows that we have a base change functor
\[
\caM(G,B) \longrightarrow  \caM(G,C), \quad  V \mapsto V \otimes_B C.
\]

\begin{prop}
Base change commutes with induction. More precisely, if $E$ is an $(M,B)$-module, then
$i_P^G(E)\otimes_B C$ is naturally isomorphic (as a $(G,C)$-module) to $i_P^G(E\otimes_B C)$.
\end{prop}
\begin{proof}
We first show that for any compact open subgroup $K$ of $G$,
\[
(i_P^G(E)\otimes_B C)^K \quad \text{ and } \quad i_P^G(E\otimes_B C)^K
\]
are naturally isomorphic as vector spaces. We use the description of the $K$-invariants of a parabolic
induction recalled above (\ref{ipgek}):
\[
(i_P^G(E)\otimes_B C)^K\simeq i_P^G(E)^K\otimes_B C \simeq \bigoplus_{g \in \Omega} E^{K_M^g} \otimes_B C
\]
and on the other hand
\[
i_P^G(E\otimes_B C)^K \simeq \bigoplus_{g \in \Omega} (E \otimes_B
C)^{K_M^g}\simeq \bigoplus_{g \in \Omega} E^{K_M^g} \otimes_B C.
\]
All the isomorphisms appearing in the two equations above are natural, and by passing to the limit as $K$
becomes smaller and smaller, we obtain $i_P^G(E)\otimes_B C \simeq i_P^G(E\otimes_B C)$. Let us now give
a formula for this isomorphism:
\begin{align} \label{ivbc}
 i_P^G(E)\otimes_B C \rightarrow i_P^G(E\otimes_B C), \quad f\otimes c
\mapsto \tilde f
\end{align}
where the function $\tilde f \colon G \rightarrow E\otimes_B C $ is given by
\[
\tilde f (g)= f(g)\otimes c.
\]
It is clear that (\ref{ivbc}) is an isomorphism of $(G,C)$-modules.
\end{proof}

\bigskip

We can apply this in the special case of specialization:
\begin{cor}
Specialization commutes with induction. More precisely, if $E$ is an $(M,B)$-module, then for all
$\chi \in S$, $\mathbf{sp}_\chi(i_P^G(E))$ is naturally isomorphic to
$i_P^G(\mathbf{sp}_\chi(E))= i_P^G(E_\chi)$, the isomorphism being induced by $f \mapsto \tilde f$,
$f \in i_P^G(E)$, where
\[
\tilde f(g)= \mathbf{sp}_\chi(f(g)).
\]
\end{cor}

\subsection{$(G,B)$-modules and restriction}\label{GBres}
If $P=MN$ is a parabolic subgroup of $G$, and if $(\pi,V)$ is a $(G,B)$-module, then $r_P^G (\pi,V)$ is
equipped with a $B$-module structure. It is obvious that the action of $M$ on $r_P^G(V)$ commutes with
this action of $B$. We will now show that $r_P^G (V)$ is a flat $B$-module. Let $E$ be any $B$-module.
Since, as a vector space
\[
r_P^G(V)=V_N=V/V(N)=V\otimes_{\caH(N)} \bbC
\]
where $\bbC$ is the trivial $\caH(N)$-module, we have
\[
E\otimes_B r_P^G(V)=E\otimes_B (V\otimes_{\caH(N)} \bbC) \simeq (E\otimes_B V)\otimes_{\caH(N)} \bbC.
\]
Now $V$ is a flat $B$-module, so the functor $ \bullet \otimes _B V$ is exact. On the other hand, the
restriction functor $\bullet \otimes_{\caH(N)} \bbC$ is also exact. This shows that $r_P^G(V)$ is a flat
$B$-module. Thus we obtain a functor
\[
r_P^G\colon \caM(G,B) \longrightarrow \caM(M,B)
\]

Proposition \ref{finiIJ} then generalizes to $(G,B)$-modules:
\begin{prop}
The functor $r_P^G\colon \caM(G,B) \longrightarrow \caM(M,B)$ preserves finitely generated $B$-modules.
\end{prop}
Here again, the crucial point is the compactness of $P\backslash G$, and thus the finiteness of
$P\backslash G/K$ for any compact open subgroup $K$ of $G$.

\subsection{Use of flatness} \label{utilplat}

As in the previous sections, $B$ is a reduced Noetherian, finitely generated and commutative
$\bbC$-algebra, and $S$ the affine algebraic variety corresponding to it. In applications, we will need
to reduce to the case where, when $(\pi,V)$ is an admissible $(G,B)$-module, and $K$ is any compact open
subgroup of $G$, the $B$-module $V^K$ is finitely generated free (which is of course a stronger condition
than finitely generated projective, any projective module being a direct summand of a free module). Now a
finitely generated projective module is locally free. \index[ter]{locally free} More precisely, if
$f\in B$ is non-zero, let $T_f$ denote the multiplicative subset of $B$ consisting of the powers of $f$,
$B_f$ the localization of $B$ at $T_f$, for any module $M$, $M_f$ the localization of $M$ (see the third
and fourth examples of \ref{propuniv}) and $S_f$ the Zariski open set of $S$ of maximal ideals not
containing $f$. \index[not]{B_f@$B_f$} \index[not]{T_f@$T_f$}
\index[not]{S_f@$S_f$}\index[not]{M_f@$M_f$}

\begin{thm}
A $B$-module $M$ is finitely generated projective if and only if there exist generators
$(f_i)_{i=1,\ldots ,r}$ of $B$, such that for all $i=1,\ldots ,r$ the localization $M_{f_i}$ is a
finitely generated free $B_{f_i}$-module.
\end{thm}

\begin{proof}
See \cite{BourAC1}, Chap. II, \S 5, no.~2, Theorem 1.
\end{proof}

We will now see applications of this result. Let $(\pi,V)$ be an admissible $(G,B)$-module, and let
$\phi$ be a $(G,B)$-module morphism from $V$ to itself, and for any compact open subgroup $K$ of $G$, let
\[
\phi_K\colon V^K \rightarrow V^K
\]
be the induced $B$-module morphism. Since $V^K$ is finitely generated projective, we find a finite family
$(f_i)_{i=1,\ldots ,r}$ of elements of $B$ generating $B$ over itself, such that for all $i=1,\ldots ,r$
the localization $V^K_{f_i}$ is a finitely generated free $B_{f_i}$-module. Let $f$ be one of the $f_i$,
and for any point $x$ in $S_f$, let
\begin{align*}
\mathbf{sp}_x(\phi) \colon  \mathbf{sp}_x(V)&\rightarrow \mathbf{sp}_x(V) \\
\mathbf{sp}_x(\phi_K) \colon   \mathbf{sp}_x(V^K)&\rightarrow   \mathbf{sp}_x(V^K)\\
\end{align*}
denote the morphisms obtained by specialization.

Note in passing that we of course have $ \mathbf{sp}_x(V_f^K)=\mathbf{sp}_x(V^K)$.

\begin{prop}
Let $P=MN$ be a parabolic subgroup of $G$, and let $(\rho,W)$ be an irreducible representation of $M$.
For all $\psi$ in $\caX(M)$, set $\pi_\psi= i_P^G (\rho\otimes \psi)$. If $\pi_{\psi_0}$ is irreducible
for a certain $\psi_0$ in $\caX(M)$, then $\pi_\psi$ is irreducible for all $\psi$ in a Zariski open set
of the algebraic variety $\caX(M)$.
\end{prop}

\begin{proof}
Set $B=\bbC[\Lambda(M)]$. This is the algebra of functions of the algebraic variety $\caX(M)$. Consider
the $(M,B)$-module
\[
(\rho_B,W_B)= (\rho \otimes \chi_{un},  W \otimes_\bbC  B)
\]
of Example \ref{GBmod}, and the induced $(G,B)$-module ({\sl cf.} \ref{GBind})
\[
(\pi_B,V_B)=i_P^G(\rho_B,W_B).
\]

Let $K$ be a compact open subgroup of $G$. The space $V_B^K$ is a $B$-module, but also has a module
structure for the Hecke algebra $\caH(G,K)$, inherited from the action of $\caH(G)$. The actions of $B$
and $\caH(G,K)$ on $V_B^K$ commute.

Let us now specialize at $\psi \in \caX(M)$. We have, by Corollary \ref{GBind}
\[
\mathbf{sp}_\psi (\pi_B,V_B)= \mathbf{sp}_\psi(i_P^G(\rho_B,W_B))=i_P^G(\rho\otimes \psi).
\]
Let $V_\psi$ denote the space of this representation. We then have
\[
V_\psi^K= i_P^G(\rho\otimes \psi)^K\simeq \mathbf{sp}_\psi(V_B)^K=\mathbf{sp}_\psi(V_B^K).
\]

For all $h \in \caH(G,K)$, let $\pi_B(h)$ denote the action of $h$ on $V_B^K$ and $\pi_\psi(h)$ the
action of $h$ on $\mathbf{sp}_\psi(V_B^K)$. This defines morphisms
\[
\pi_B^K\colon \caH(G,K) \rightarrow \End_B(V_B^K),
\]
\[
\pi_\psi^K\colon \caH(G,K) \rightarrow \End(V_\psi^K).
\]

We now wish to use Theorem \ref{piK} to show that $ i_P^G (\rho \otimes \psi)$ is generically irreducible
by reducing to showing that for a sufficiently small compact open subgroup $K$ of $G$, $V_\psi^K$ is an
irreducible $\caH(G,K)$-module. To this end, we must ensure that all irreducible subquotients of
$i_P^G(\rho\otimes \psi)$ have non-zero vectors fixed by $K$. For fixed $\psi$, this is obvious, but $K$
must satisfy this property for all $\psi$ to be able to apply the argument. We can of course assume that
$M$ is a standard Levi subgroup of $G$. We note that all irreducible subquotients of the
$i_P^G(\rho\otimes \psi)$ have the same inertial support (determined by that of $\rho$), say
$\Omega=[L,(\sigma,E)]_G$, with $L$ a standard Levi subgroup of $M$ (and thus of $G$). It therefore
suffices to show that there exists a compact open subgroup $K$ of $G$ such that all irreducible
representations $(\tau,W)$ with inertial support $\Omega$ satisfy $W^K\neq \{0\}$. Let $(\tau,W)$ be an
irreducible representation of $G$ with inertial support $\Omega$. Let $Q$ be the standard parabolic
subgroup of $G$ with Levi factor $L$. By Theorem \ref{dcompthm}, or rather its proof, $r_Q^G(\tau,W)$
admits an irreducible subquotient in an inertial class $[L, (\sigma',E')]_L$ with
$[L, (\sigma',E')]_G=[L, (\sigma,E)]_G$. Let $K$ be a compact open subgroup admitting an Iwahori
decomposition with respect to the standard parabolic subgroups and such that $(E')^{K_L}\neq \{0\}$ for
$(\sigma',E')$ running through a (finite) set of representatives of the inertial classes
$[L, (\sigma',E')]_L$ with $[L, (\sigma',E')]_G=[L, (\sigma,E)]_G$. Moreover $(E')^{K_L}$ is independent
of the choice of representative $\sigma'$ in its inertial class. Thus, for such a $K$,
$(r_Q^G(\tau,W))^{K_L}\neq \{0\}$. We deduce from Jacquet's lemma (Theorem \ref{LemmeJacquet}) that
$W^K\neq \{0\}$. We fix such a compact open subgroup $K$.

By Burnside's theorem (see \cite{GW} 3.1.3), $V_\psi^K$ is an irreducible $\caH(G,K)$-module exactly when
$\pi_\psi^K$ is surjective. Suppose that $\pi_{\psi_0}^K$ is surjective, so that $V_{\psi_0}^K$ is an
irreducible $\caH(G,K)$-module, and note that the description of the $K$-invariants (\ref{ipgek}) in the
induced representation $V_\psi$ shows that the fact that $V_{\psi}^K$ is non-zero is independent of
$\psi$ (what depends on $\psi$ of course is the $\caH(G,K)$-module structure of $V_{\psi}^K$). We must
show that $\pi_{\psi}^K$ is surjective for all $\psi$ in a Zariski open set of the algebraic variety
$\caX(M)$.

Let us resume the discussion preceding the statement of the proposition (so here $S=\caX(M)$). This
$\psi_0$ is in one of the principal open sets $S_f$, where the localization $V_{B,f}^K$ is free over
$B_f$. Let us still denote
\[
\pi_B^K\colon \caH(G,K) \rightarrow \End_{B_f}(V_{B,f}^K),
\]
and choose a basis $\caB$ of $\End_{B_f}(V_{B,f}^K)$ over $B_f$. Let us specialize at $\psi \in S_f$:
$\caB_\psi=\mathbf{sp}_\psi(\caB)$ is a basis of the $\bbC$-vector space
\[
\mathbf{sp}_\psi(\End_{B_f}(V_{B,f}^K))= \End(V_\psi^K).
\]
Choose a linearly independent family $\caF$ of $\caH(G,K)$ such that the image of $\caF$ by
$\pi_{\psi_0}$ is exactly the basis $\caB_{\psi_0}$. Let $F$ be the subspace of $\caH(G,K)$ generated by
$\caF$, and let $M=(m_{ij})_{ij}$ denote the matrix of the restriction of $\pi_B^K$ to $F$: it is a
matrix with coefficients in $B_f$, and $\Psi_\psi(m_{ij})$ is a matrix with complex coefficients, namely
the matrix of the restriction of $\pi_\psi^K$ to $F$ in the bases $\caF$ and $\caB_\psi$. The determinant
of this matrix is a polynomial in the $m_{ij}$, which is equal to $1$ at $\psi_0$. It is therefore
non-zero on a Zariski open set of $S_f$. On this open set, $\pi_\psi^K$ is therefore surjective.
\end{proof}

\subsection{Generic irreducibility of induced representations}\label{irrindgen}
\index[ter]{generic irreducibility}

Using the geometric lemma and Proposition \ref{utilplat}, we establish the following theorem.

\begin{thm}
Let $P=MN$ be a parabolic subgroup of $G$, and let $(\rho,W)$ be an irreducible supercuspidal
representation of $M$. For all $\psi$ in $\caX(M)$, set $\pi_\psi= i_P^G (\rho\otimes \psi)$. Then
$\pi_\psi$ is irreducible for all $\psi$ in a non-empty Zariski open set of the algebraic variety
$\caX(M)$.
\end{thm}
\begin{proof}
By Proposition \ref{utilplat}, it suffices to show that $\pi_{\psi_0}$ is irreducible for a certain
$\psi_0$ in $\caX(M)$.

Let $A$ be the split component of $M$. Let $\chi$ be the central character of $\rho$. The restriction of
the character $|\chi|$ of $Z(M)$ to $A$ takes values in $\bbR^\times_+$, and thus is trivial on ${}^0A$.
It therefore defines a character of $\Lambda(A)=A/{}^0A$. Since this lattice injects with a finite index
into $\Lambda(M)$ ({\sl cf.} \ref{Lambda}), we can extend $|\chi|$ to a character $\psi \in \caX(M)$. We
obtain $|\psi^{-1} \chi|=1$ on $A$, so $|\psi^{-1} \chi \vert $ factors into a character of the compact
group $Z(M)/A$ and since any character of a compact group is unitary, $|\psi^{-1} \chi|=1$ on $Z(M)$.
Since $\rho \otimes \psi^{-1}$ is supercuspidal with unitary central character, it is unitary. Indeed, a
supercuspidal representation is compact modulo the center by Theorem \ref{supercusp}, thus essentially
square-integrable modulo the center, and we can then conclude by Lemma \ref{unit2intlem}. By translating
by such a character $\psi^{-1}$, we may therefore assume, in the statement of the theorem, that $\rho$ is
unitary. Moreover, for any unitary $\psi \in \caX(M)$, $\rho \otimes \psi$ is then unitary, and thus also
$ i_P^G (\rho \otimes \psi,W)$. In particular, this latter representation is completely reducible.

The representation $ r_P^G i_P^G (\rho \otimes \psi,W)$ admits by \ref{geolemmacusp} $(iii)$ a filtration
whose quotients are of the form $w \cdot (\rho \otimes \psi)$ for $w \in W(M,M)/W_M$. We will show that
for all $\psi$ outside a Zariski closed set of the real variety $\im(\caX(M))$, the central characters of
$\rho\otimes \psi$ and of $w \cdot (\rho \otimes \psi)$ are different if $w$ is non-trivial. As above,
let $\chi$ denote the central character of $\rho$. The central character of $\rho\otimes \psi$ is
therefore $\chi \psi_{|Z(M)}$ and that of $w \cdot (\rho \otimes \psi)$ is $w \cdot (\chi \psi_{|Z(M)})$.
Since the lattice $\Lambda(A)$ injects with a finite index into $\Lambda(M)$, it suffices to show that
$\chi \alpha$ is different from $w \cdot (\chi\alpha)$ for all $\alpha$ outside a Zariski closed set of
$\caX(A)$. But this follows from the fact that $W(M,M)/W_M \simeq W(A)$ (\ref{WMWA}) which acts
faithfully on $\fra^*$ (\ref{actfid}) and from the description of $\im(\caX(M))$ in \ref{ReChi}. We
conclude by using the following lemma.

\begin{lemme}
Let $(\pi,V)$ be a representation of $G$ admitting a filtration
\[
\{0\}=V_0\subset V_1 \subset \cdots \subset V_n=V
\]
whose successive quotients $V'_i=V_i/V_{i-1}$ have distinct central characters. Then
\[
V\simeq \bigoplus_{i=1}^n V'_i.
\]
\end{lemme}
\begin{proof}
It is clear that the result is true for representations restricted to $Z(G)$. But since $G$ and $Z(G)$
commute, the action of $G$ preserves the decomposition into distinct characters of $Z(G)$.
\end{proof}

\bigskip

We deduce from the above that $ r_P^G i_P^G (\rho \otimes \psi,W)$ decomposes into a direct sum of
modules of the form $w \cdot (\rho \otimes \psi)$, $w \in W(M,M)/W_M$ whose central characters are
pairwise distinct, for all $\psi$ outside a Zariski closed set of $\im(\caX(M))$. We then have using
Frobenius reciprocity
\begin{align*}
 &\Hom_G(\pi_\psi,\pi_\psi)=  \Hom_M( r_P^G(\pi_\psi),\rho\psi)\\
&=\bigoplus_{w\in W(M,M)/W_M}  \Hom_M(w \cdot (\rho \psi)  ,\rho\psi)=
\Hom_M(\rho \psi  ,\rho\psi)
\end{align*}
which is clearly of dimension $1$ since $\rho\psi$ is irreducible. We deduce that $\pi_\psi$ is
irreducible. This completes the proof of the theorem.
\end{proof}

\subsection{Another application}\label{autreapp}

We return to the situation of \ref{utilplat}, preceding the statement of the proposition. In this
context, let $\Xi$ be a Zariski dense set of points in $S$ (since $B$ is reduced, the $x\in \Xi$ separate
the elements of $B$).

Fix a basis of $V^K_{f}$ and let $(m_{kl})_{k,l}$ be the matrix of the endomorphism
\[
\phi_{K,f}\colon V^K_{f} \rightarrow V^K_{f}
\]
in this basis. Since $\Xi$ is Zariski dense in $S$, its intersection with the open set $S_{f}$ is still
Zariski dense, and for all $x \in S_{f} \cap \Xi$, the matrix $(\Psi_x(m_{k,l}))_{k,l}$ is a matrix with
complex coefficients. It is the matrix of $\mathbf{sp}_x(\phi_{K})$.

\begin{prop}
Suppose that for all $x \in \Xi$, the specialization at $x$ of $V$ is irreducible. Then $\phi$ is the
action of an element $b \in B$.
\end{prop}
\begin{proof}
By Schur's lemma, for all $x \in \Xi$, $\mathbf{sp}_x(\phi)$ is a scalar operator, so the same is true
for $\mathbf{sp}_x(\phi_K)$. Thus the matrix $(\Psi_x(m_{k,l}))_{k,l}$ satisfies $\Psi_x(m_{k,l})=0$ if
$k \neq l$ and $\Psi_x(m_{k,k})=\Psi_x(m_{l,l})$ for all $k,l$. Since $\Xi$ is Zariski dense, we deduce
that $(m_{kl})_{k,l}$ is scalar. It follows that $\phi_{K,f}$ is given by the action of a $b \in B_f$.
This holding for any compact open subgroup $K$ of $G$ and for all the $f_i$, we deduce that $\phi$ is
given by the action of a $b \in B$ (the elements $b$ thus obtained glue together, since an element of
$B_f$ is determined by its action, $B$ being reduced and $\Xi$ Zariski dense).
\end{proof}

\begin{exemple}
Let $P=MN$ be a parabolic subgroup of $G$, and let $(\rho,W)$ be an irreducible supercuspidal
representation of $M$. We retain the notation of the proof of Proposition \ref{utilplat}. Let
\[
\phi\colon V_B \rightarrow V_B
\]
be a $(G,B)$-module morphism. By the generic irreducibility theorem, the hypotheses of the proposition
above are satisfied, and we conclude that there exists a $b \in B$ such that $\phi$ is the action of $b$.

Now consider the case $M=P=G$. Let $z$ be an element of the center of the category $\caM(G)_{[\rho]}$. It
is by definition a natural transformation of the identity functor which gives in particular a
$G$-equivariant endomorphism of $W_B=W\otimes B$. We obtain that the action of $z$ in $W_B=W\otimes B$ is
given by the action of an element $b \in B$ and we thus recover the results of \ref{centreMpi}.
\end{exemple}

\section{The second adjunction theorem}

\subsection{Generalized Jacquet lemma} \label{LemmeJacquet2}
\index[ter]{Jacquet's lemma!(generalized)} The statement of the following theorem differs from Theorem
\ref{LemmeJacquet} only in that we no longer make the admissibility hypothesis on the representation. The
proof is considerably more involved.

\begin{thm}
Let $P=MN$ be a parabolic subgroup of $G$, with split component $A$, and let $(\pi,V)$ be a smooth
representation of $G$. Let $K$ be a compact open subgroup of $G$ admitting an Iwahori decomposition with
respect to $P$. Then the projection map $j \colon V\rightarrow  V_N$ maps $V^K$ surjectively onto
$(V_N)^{K_M}$.
\end{thm}
\begin{proof}
We retain the steps of the proof of Theorem \ref{LemmeJacquet} which do not rely on the admissibility
hypothesis of $(\pi,V)$. First, let us verify that the image of $V^K$ under $j$ is indeed in
$(V_N)^{K_M}$. We then choose $t \in  C_A^{++}$. The set $ \{ t^{-m}K_{\overline{N}}t^m|\, m\in \bbN \}$
forms a basis of neighborhoods of the identity in $\overline{N}$, and the union of the
$ t^{-m}K_{N}t^m, \,  m\in \bbN$ is equal to $N$ (Theorem \ref{KN}).

Let $v \in V^K$, and let $m \in \bbN$. We still have, as in \ref{LemmeJacquet}
\begin{equation}\label{jvk}
 j( \pi (e_K) \pi (t^m) \cdot v) = \pi_N(t^m) \cdot j(v).
\end{equation}
This shows that
\begin{equation} \label{jVK}
 \pi_N(t^m) \cdot j(V^K) \subset  j(V^K)
\end{equation}
Although $\pi_N(t^m)$ acts invertibly on $j(V)=V_N$, the inclusion above may be strict without the
admissibility hypothesis.

Let $\bar v \in V_N^{K_M}$. Still without using the admissibility hypothesis, we have shown that there
exists $m\in \bbN$ such that
\begin{equation}\label{locenA}
 \pi_N(t^m)\cdot \bar v \in j(V^K),
\end{equation}
i.e., we have
\begin{equation}\label{locenB}
 \bigcup_{m \in \bbN}  \pi_N(t^{-m}) \cdot  j(V^K)=V_N^{K_M}.
\end{equation}

Consider the distribution $a_{t,K}=e_K*\delta_t*e_K  \in \caH(G,K)$, and let $A$ denote the endomorphism
of $V^K$ given by the action of $a_{t,K}$.
\begin{lemme}
We have
$$V(N)\cap V^K= \left(\bigcup_{i\in \bbN}  \ker \pi(a_{t^i,K})\right) \cap V^K= \bigcup_{i\in \bbN}  \ker A^i.$$
\end{lemme}

\begin{proof}
Note first that by the second lemma of \ref{HCDH}, we have $a_{t^i,K}= (a_{t,K})^i$. Set
$N_i= t^{-i}  K_N t^i$ and recall that $\bigcup_{i\in \bbN} N_i=N$. By Proposition \ref{FoncJac1}, we
then have $V(N)=\bigcup_{i\in \bbN} \ker e_{N_i}$, whence
\[
V(N)\cap V^K= \left( \bigcup_{i\in \bbN} \ker e_{N_i}\right) \cap
V^K=\bigcup_{i\in \bbN} \ker {e_{N_i}}_{|V^K} .
\]
Now,
\begin{align}
\nonumber a_{t^i,K} & = e_K * \delta_{t^i} * e_K= e_{K_N} * e_{K_M} * e_{K_{\bar  N}}
*\delta_{t^i} * e_K \\
\nonumber &=  e_{K_N} * e_{K_M} * \delta_{t^i} * e_{ t^{-i} K_{\overline{N}} t^i }*e_K=
e_{K_N} * \delta_{t^i} * e_{ t^{-i} K_M t^i }*e_K \\
\label{KNiK}&= \delta_{t^i} * e_{ t^{-i} K_N t^i }* e_K=  \delta_{t^i} *e_{N_i}*e_K.
\end{align}
We deduce that the restriction to $V^K $ of $a_{t^i,K}$ coincides with the restriction of
$\delta_{t^i} * e_{N_i}$. Since $\delta_{t^i}$ is invertible,
$\ker\pi(a_{t^i,K})_{|V^K}=\ker A^i= \ker {e_{N_i}}_{|V^K}$.
\end{proof}

\bigskip

We now need the notion of the localization of a vector space at an endomorphism of this space and of an
eventually stable endomorphism. The definitions and the results used are developed in the appendices, in
the fifth example of \ref{propuniv} and in \ref{stable}.
\begin{lemme}
The localization of $V^K$ at $A$ is naturally isomorphic to $$(V_N^{K_M},\pi_N(t)),$$ the canonical
morphism $\iota \colon  (V^K,A) \rightarrow (V_N^{K_M},\pi_N(t))$ being given by $j$.
\end{lemme}
\index[ter]{localization}
\begin{proof}
Recall the results of Lemma \ref{propuniv} which describe the localization. The localization of $V^K$ at
$A$ is, according to this lemma, (naturally) isomorphic to the localization of
$V^K/ (\cup_{i\in \bbN} \ker A^i)$ at $A'$, where $A'$ is the endomorphism induced by $A$. By the lemma
above
\[
V^K/ (\cup_{i\in \bbN} \ker A^i) = V^K /(V(N)\cap V^K)\simeq j(V^K)
\]
and the endomorphism $A'$ is given by $\pi_N(t)$ (this follows from (\ref{jvk})). It now remains to show
that $(V_N^{K_M},\pi_N(t))$ is indeed this localization. By Lemma \ref{propuniv}, it suffices to verify
that $\pi_N(t)$ is invertible on $(V_N)^{K_M}$ and that
\begin{equation} \label{jVKM}
 V_N^{K_M} = \bigcup_{i\in \bbN} \pi_N(t)^{-i}  j(V^K).
\end{equation}
Since
\[
\pi(e_{K_M})\pi(t)= \pi(t)\pi (e_{t^{-1} K_M t})= \pi(t)\pi (e_{K_M }),
\]
we have $\pi_N(e_{K_M})\pi_N(t)=\pi_N(t)\pi_N(e_{K_M})$ and thus $\pi_N(t)$ indeed defines an
endomorphism of $V_N^{K_M}$, and similarly for $\pi_N(t^{-1})$. Equality (\ref{jVKM}) is a reformulation
of (\ref{locenB}).
\end{proof}

\begin{prop}
The morphism $A$ is eventually stable. More precisely, there exists a constant $b=b(G,K)$ depending only
on $G$ and $K$ such that for any smooth representation $(\pi,V)$, $A^b$ is a stable endomorphism of
$V^K$. Set $V^K_0 =\ker A^b$, $V_*^K=\im A^b$. We then have
\[
V^K=V_0^K \oplus V_*^K.
\]
The localization of $V^K$ at $A$ is then also (naturally) isomorphic to $(V_*^K, A_{|V_*^K})$, the
canonical morphism $\iota$ being given by the projection onto $V_*^K$ parallel to $V_0^K$. Moreover, the
constant $b$ can be chosen to be less than or equal to the constant $c(G,K)$ of the uniform admissibility
theorem \ref{unifadm}.
\end{prop}

\index[not]{V_0K@$V_0^K$} \index[not]{V_*K@$V_*^K$}

Before moving on to the proof of this proposition, which will occupy several paragraphs, we now complete
the proof of the theorem. We deduce from the two lemmas and the proposition that $j$ realizes a (natural)
isomorphism from $V_*^K$ to $V_N^{K_M}$ and that the kernel of $j$ on $V^K$ is $V_0^K$. We have done more
than simply show surjectivity in the statement of the theorem, because the proof provides a section of
$j \colon V^K \rightarrow V_N^{K_M}$, namely the subspace $V_*^K$.
\end{proof}

\begin{rmqs}
1. The decomposition $V^K=V^K_0\oplus V^K_*$ is independent of the choice of $t$ in $C_A^{++}$. If
$(\pi,V)$ is admissible, we recover the decomposition of Proposition \ref{LemmeJacquet}.

--- 2. There exist compact open subgroups $C$ of $N$ and $\bar C$ of $\overline{N}$ such that for all
$(\pi,V)\in \caM(G)$,
\[
V_0^K= V^K \cap \ker e_C,\quad V_*^K= e_K*e_{\bar C}\cdot V
\]

--- 3. There exists $\epsilon>0$ such that for all $t \in A^+(\epsilon)$,
\[
V_0^K=\ker \pi(a_{t,K}), \quad   V_*^K=\im \pi(a_{t,K}).
\]
\end{rmqs}
\begin{proof}
The first point is the generalization of Proposition \ref{LemmeJacquet}. It is proved in the same way
using the fact that for any other element $s$ of $C_A^{++}$, by the second lemma of \ref{HCDH}, we have
for all $n$ in $\bbN$,
\[
a_{s^n,K}a_{t^n,K}= a_{(st)^n,K}.
\]
It also follows from the second point.

For the second point, it suffices to choose $C$ such that $ t^{-b} K_N t^{b} \subset C$ and
$ t^{b} K_{\overline{N}} t^{-b} \subset  \bar C$. We then have by (\ref{KNiK}), for all
$v \in V_0^K=\ker \pi(a_{t^b,K})$,
\[
0=\pi(a_{t^b,K})\cdot v=e_K*\delta_{t^b}*e_K\cdot v=
\delta_{t^b}*e_{N_b}\cdot v.
\]
Since the action of $\delta_{t^b}$ is invertible, we see that $v \in \ker e_{N_b}$, and since
$N_b \subset C$, we have $ e_C=e_C* e_{N_b}$, and thus $v \in \ker  e_C \cap V^K$. Conversely, if
$v \in \ker e_C$, then $v \in V(N)$ (Proposition \ref{FoncJac1}).

To show $V_*^K= e_K*e_{\bar C}\cdot V$, we use the equality
\[
a_{t^b,K}=e_K*e_{\overline{N}_b}*\delta_{t^b},
\]
where $\overline{N}_b=t^b K_{\overline{N}} t^{-b}$, which is shown like (\ref{KNiK}). Since
$\overline{N}_b \subset \bar C$ by hypothesis, we have $e_{\bar N_b}*e_{\bar C}=e_{\bar C}$, and we
obtain
\[
V_*^K= a_{t^b,K}\cdot V =  e_K*e_{\overline{N}_b}*\delta_{t^b}\cdot V
= e_K*e_{\overline{N}_b}\cdot V \supset  e_K*e_{\bar C}\cdot V.
\]
Conversely, since $A^b$ is stable, $V_*^K=\pi( a_{t^{bm},K})\cdot V^K$ for all $m\in \bbN^*$, and for $m$
large enough
\[
\bar C\subset \overline{N}_{bm}=t^{bm} K_{\overline{N}} t^{-bm},
\]
which gives the inclusion in the other direction.

The third point follows quite easily from the second. Indeed, we know (Theorem \ref{KN}) that there
exists $\epsilon>0$ such that, $C$ and $\bar C$ being as above, for all $t \in A^+(\epsilon)$,
\[
C\subset t^{-1}K_Nt , \quad \bar C \subset   tK_{\overline{N}}t^{-1}.
\]
We then have as above, for all $v \in V^K$,
\[
\pi(a_{t,K})\cdot v  =\delta_t*e_{t^{-1}K_Nt}\cdot v,
\]
and we deduce that
$\ker \pi(a_{t,K})_{|V^K}= \ker \pi(e_{t^{-1}K_Nt})_{|V^K}  \subset V^K \cap V(N)= V_0^K$ (see the proof
of the lemma above). Since $ C\subset t^{-1}K_Nt$, $e_{t^{-1}K_Nt}= e_{t^{-1}K_Nt}*e_C$, any
$v \in V^K_0=V^K \cap \ker e_C$ is in $\ker \pi(e_{t^{-1}K_Nt})_{|V^K}$, thus in
$\ker \pi(a_{t,K})_{|V^K}$. Similarly
\[
V_*^K=  e_K*e_{\bar C}\cdot V \supset   e_K*e_{\bar C}* e_{ tK_{\overline{N}}t^{-1}} \cdot V
=  e_K* e_{ tK_{\overline{N}}t^{-1}} \cdot V =\pi(a_{t,K})\cdot V.
\]
The converse comes from
\[
V_*^K= \pi(a_{t^b,K})\cdot V \subset \pi(a_{t,K})\cdot V.
\]
\end{proof}

The following paragraphs are devoted to the proof of the proposition.

\subsection{A special case}
Let $D$ be an inertial class of irreducible supercuspidal representations of a Levi subgroup $L$ of $G$,
and $(\rho,W)$ a supercuspidal representation in this inertial class. Recall that we introduced in Remark
\ref{progenpi} a progenerator
$(\Pi_1,V_{\Pi_1})= \ind_{{}^0L}^{\; \,  L}( \res_{{}^0L}^{\; \,  L} (\rho,W))$ of $\caM(L)_D$. On the
other hand, it was noted that
\[
V_{\Pi_1}  \simeq \ind_{{}^0L}^{\; \,  L} ( \mathrm{Triv})\otimes
W \simeq F \otimes W=W_F,
\]
where $F$ is the algebra of regular functions on $\caX(L)$. Indeed,
$\ind_{{}^0L}^{\; \,  L} ( \mathrm{Triv})$ is isomorphic to $\bbC[\Lambda(L)]$, the group algebra of
$\Lambda(L)=L/{}^0L$. Let us simply denote $(\Pi_1,V_{\Pi_1})= (\Pi_1,V_1)$ to lighten the notation.
Clearly, $V_1$ is an admissible $(L,F)$-module (a special case of Example \ref{GBmod}). By Lemma
\ref{GBind}, if $Q$ is a parabolic subgroup with Levi factor $L$, then
\[
(\Pi,V_\Pi)= i_Q^G(\Pi_1,V_1)
\]
is an admissible $(G,F)$-module.

Corollary \ref{GBind} asserts that the specialization of $V_{\Pi}$ at $\chi \in \caX(L)$ is
\begin{equation}\label{spind}
\mathbf{sp}_\chi(\Pi,V_\Pi)= \mathbf{sp}_\chi(i_Q^G (\Pi_1,V_1))=i_Q^G
(\mathbf{sp}_\chi(\Pi_1,V_1))=i_Q^G(\rho\otimes  \chi,W).
\end{equation}

We will now prove Jacquet's lemma in the special case where $(\Pi,V_\Pi)=i_Q^G(\Pi_1,V_1)$. Let us fix a
parabolic subgroup $P=MN$ and a compact open subgroup $K$ as in the statement of Jacquet's lemma of the
previous section, whose notation we also resume. Since $\Pi_1$ is supercuspidal, by Proposition
\ref{geolemmacusp}, $r_P^G i_Q^G (\Pi_1,V_1)$ admits a filtration whose quotients are of the form
$i_{w\cdot L}^M({}^w \Pi_1)$ if $L$ is conjugate to a Levi subgroup of $M$, and zero otherwise. It
follows that $r_P^G i_Q^G (\Pi_1,V_1)$ is an admissible $(M,F)$-module. By the second lemma of Section
\ref{LemmeJacquet2}, we see that the endomorphism $A$ of $V_\Pi^K=(i_Q^G V_1)^K$ satisfies the hypotheses
of Proposition \ref{stable}. We deduce that it is eventually stable, i.e., that there exists $b\in \bbN$
such that
\begin{equation}\label{VPiK}
V_\Pi^K = (V_\Pi)_0^K\oplus (V_\Pi)_*^K
\end{equation}
with $A^b( (V_\Pi)_0^K)=0$ and $A^b$ invertible on $(V_\Pi)_*^K$. We now want to show that we can take
$b \leq c(G,K)$, the uniform admissibility constant of Theorem \ref{unifadm}. To do this, by Lemma
\ref{stable}, $(iv)$, it suffices to see that $A^c$ annihilates $(V_\Pi)^K_0$.

By Theorem \ref{irrindgen}, the representation $\mathbf{sp}_\chi(\Pi,V_\Pi)= i_Q^G(\rho \otimes \chi) $,
$\chi \in \caX(L)$, is irreducible for all $\chi$ outside a Zariski closed set of $\caX(L)$. In
particular, for all these $\chi$,
\begin{equation}\label{spcGK}
 \dim \mathbf{sp}_\chi(V_\Pi)^K= \dim (i_Q^G(W\otimes \bbC_\chi))^K\leq c= c(G,K),
\end{equation}
where $c(G,K)$ is the uniform admissibility constant \ref{unifadm}.

The $F$-module $V_\Pi^K$ is also an $\caH(G,K)$-module. The actions of $F$ and $\caH(G,K)$ commute, and
the decomposition (\ref{VPiK}) is a direct sum of $F$-modules. For all $h\in \caH(G,K)$, we denote by
$\Pi(h)$ the action of $h$ on $V_\Pi^K$ and $\Pi_\chi(h)$ the action of $h$ on
$\mathbf{sp}_\chi(V_\Pi)^K$. The situation is similar to that encountered in \ref{utilplat}. The image
under $\mathbf{sp}_\chi$ of $(V_\Pi)_0^K \subset (V_\Pi)^K$ is annihilated by the operator
$\Pi_\chi(A^c)$, for all $\chi$ such that (\ref{spcGK}) is satisfied.

We use the same technique as in \ref{utilplat}. Consider one of the principal open sets $S_f$ of
$\caX(L)$, where the localization $V_{\Pi,f}^K$ is free over $F_f$. Let us still denote by $\Pi(h)$ the
action of an element $h \in \caH(G,K)$ on $V_{\Pi,f}^K$. We deduce from (\ref{VPiK}) that
\[
V_{\Pi,f}^K= (V_{\Pi,f})_0^K\oplus (V_{\Pi,f})_*^K
\]
and choose a basis $\caB$ of $V_{\Pi,f}^K$ over $F_f$, the union of the bases $\caB_0$ and $\caB_1$ of
$(V_{\Pi,f})_0^K$ and $ (V_{\Pi,f})_*^K$ respectively (by localizing again if necessary, so that the
direct factors of the right-hand side are also free). Let us specialize at $\chi \in S_f$:
$\caB_\chi=\mathbf{sp}_\chi(\caB)$ is a basis of the $\bbC$-vector space $i_Q^G(W\otimes \bbC_\chi)^K$.
The matrix of $\Pi(A^c)$ in the basis $\caB$ has coefficients in $F_f$, and its image under $\Psi_\chi$
is the matrix (with complex coefficients) of $\Pi_\chi(A^c)$. Since $\Pi_\chi(A^c)$ annihilates $\caB_0$
for all $\chi$ in a Zariski open set of $S_f$, and $F$ is reduced, we deduce (by an argument previously
employed in \ref{autreapp}) that $\Pi(A^c)$ annihilates $\caB_0$ and thus $(V_{\Pi,f})_0^K$. Since this
holds for all the principal open sets $S_f$ of Theorem \ref{utilplat}, we obtain that $\Pi(A^c)$
annihilates $(V_{\Pi})_0^K$, and thus that $b\leq c(G,K)$.

This completes the proof of Jacquet's lemma in the special case of the representation $i_Q^G(\Pi_1,V_1)$.

\subsection{End of the proof of Jacquet's lemma}

We will now extend the result of Proposition \ref{LemmeJacquet2} to all smooth representations of $G$. We
begin with induced representations. Let $Q=LU$ be a parabolic subgroup of $G$, $D$ an inertial class of
irreducible supercuspidal representations of $L$, and $(\rho,W)$ a representation of $\caM(L)_D$. We want
to show that $i_Q^G(\rho,W)$ satisfies Proposition \ref{LemmeJacquet2}. In the previous section, we
established this for the representation $i_Q^G(\Pi_1,V_1)$, where $(\Pi_1,V_1)$ is a small progenerator
of $\caM(L)_D$. By Lemma \ref{rmqfonct}, the representation $(\rho,W)$ is the cokernel of a morphism
between two representations which are direct sums of representations isomorphic to $(\Pi_1,V_1)$, i.e.,
we have an exact sequence of the form
\[
\bigoplus_\alpha V_1 \rightarrow \bigoplus_\beta V_1 \rightarrow   W
\rightarrow 0.
\]
Since the functors $i_Q^G$ and $j_K$ are exact and preserve direct sums (Corollary \ref{norm}), we obtain
an exact sequence
\[
\bigoplus_\alpha (i_Q^G V_1)^K \rightarrow \bigoplus_\beta (i_Q^G V_1)^K \rightarrow  (i_Q^G W)^K
\rightarrow 0.
\]
Since $A^b$ is stable on $(i_Q^G V_1)^K$, the endomorphism induced by $a_{t^b,K}$ on
$(\bigoplus_\alpha i_Q^G V_1)^K$ is stable (ditto for the second sum), for the same constant $b\in \bbN$
as before. The result for $(i_Q^G W)^K$ follows immediately (Lemma \ref{stable}, $(ii)$). This
establishes Jacquet's lemma for $i_Q^G(\rho, W)$.

Let us now move on to an arbitrary representation $(\pi,V)$ of $\caM(G)$. The decomposition theorem
\ref{dcompthm} allows us to write
\[
V=\bigoplus_{\frs \in \caB(G)}  V_\frs,
\]
and we are thus reduced to showing the proposition for each component $V_\frs$. For the supercuspidal
components, this is trivial by definition. There remain the induced components. By the second lemma of
\ref{dcompthm}, $(iii)$, each $V_\frs$ embeds into a finite direct sum of induced representations of the
form
\[
\bigoplus_{i\in I} i_{Q_i}^G (\rho_i,W_i).
\]
Consider the cokernel $C$ of this injection. In the same way, $C$ embeds into a sum of the same form. We
deduce an exact sequence of the form
\[
0\rightarrow V_\frs \rightarrow \bigoplus_{i\in I} i_{Q_i}^G (\rho_i,W_i)
\rightarrow \bigoplus_{j\in J}  i_{Q_j}^G (\rho_j,W_j)
\]
which realizes $V_\frs$ as the kernel of a $G$-morphism between two representations for which the
proposition is established. We then use the exactness of the functor $j_K$, and the fact that it
preserves finite direct sums to obtain an exact sequence
\[
0\rightarrow V_\frs^K \rightarrow \bigoplus_{i\in I} i_{Q_i}^G (\rho_i,W_i)^K
\rightarrow \bigoplus_{j\in J}  i_{Q_j}^G (\rho_j,W_j)^K.
\]
Lemma \ref{stable} $(iii)$ allows us to conclude.
\qed

\subsection{A consequence} \label{consLJG}

A first consequence of the generalized Jacquet lemma is that Corollary \ref{ThmHowe} is now valid without
the admissibility hypothesis. Let us restate it:

\begin{cor}
Let $K$ be a compact open subgroup of $G$, contained and normal in $K_0$ and admitting an Iwahori
decomposition with respect to the standard parabolic subgroups. Let $(\pi,V)$ be a smooth representation
of $G$ such that $V^K$ generates $V$. Then any subquotient of $V$ is generated by its vectors fixed by
$K$.
\end{cor}

We deduce from this, by virtue of Theorem \ref{35}:

\begin{prop}
Let $K$ be a compact open subgroup of $G$, contained and normal in $K_0$ and admitting an Iwahori
decomposition with respect to the standard parabolic subgroups. Then the full subcategory of $\caM(G)$ of
representations generated by their $K$-invariant vectors is equivalent to the category of unital modules
over the Hecke algebra $\caH(G,K)$.
\end{prop}

\subsection{Another consequence} \label{cons2}
Let us place ourselves in the setting of Section \ref{GBres}. The restriction functor
\[
r_P^G:  \caM(G,B) \longrightarrow  \caM(M,B)
\]
preserves $B$-admissibility. The argument is the same as that of Corollary \ref{LemmeJacquet}.

\subsection{Bernstein's second adjunction theorem}\label{secondeadjonction}
\index[ter]{second adjunction} Let $P=MN$ be a parabolic subgroup of $G$. The functor $r_P^G$ is the left
adjoint of the functor $i_P^G$, but the latter also admits a right adjoint.

This fact is not particularly deep: the functor $i_P^G$ is the composition of the forgetful functor from
$M$ to $P$ (which admits a right adjoint), a normalization functor (which is inconsequential here), and
the parabolic induction functor $\Ind_P^G$. Since $G/P$ is compact, it coincides with the compact
induction functor and thus admits a right adjoint (cf. \ref{indcomp}).

The explicit determination of this right adjoint, however, is much deeper. Let $\bar P=M\overline{N}$ be
the parabolic subgroup opposite to $P$.
\begin{thm}
The functor $r_{\bar P}^G$ is the right adjoint of the functor $i_P^G$. For any representation $(\pi,V)$
in $\caM(G)$ and for any representation $(\tau,E)$ in $\caM(M)$, we therefore have a natural isomorphism
\begin{equation}\label{2AD}
  \Hom_G(i_P^G(\tau,E),(\pi,V))\simeq
\Hom_M((\tau,E),r_{\bar P}^G(\pi,V)).
\end{equation}
\end{thm}
This result is known as Bernstein's second adjunction theorem. We will see that it follows from Jacquet's
lemma. We first show that it is equivalent to the following statement:
\begin{thm}
For any smooth representation $(\pi,V)$ of $G$,
\begin{equation} \label{2ad2}
r_{\bar P}^G \tilde \pi= (r_P^G \pi)^\sim .
\end{equation}
\end{thm}
\begin{proof}
We have, for any smooth representation $(\tau,E)$ of $M$ and any smooth representation $(\pi,V)$ of $G$,
\begin{align}\label{rpdual}
\nonumber  \Hom_G(i_P^G \tau,\tilde \pi) &\simeq  \Hom_G(\pi, (i_P^G\, \tau)^\sim) \\
 \nonumber          &\simeq  \Hom_G(\pi ,i_P^G \, \tilde \tau)\\
  \nonumber          & \simeq  \Hom_M(r_P^G(\pi), \tilde     \tau)\\
      &\simeq  \Hom_M(\tau,r_P^G(\pi)^\sim ),
\end{align}
The first equality is Lemma \ref{VWWV}. The second is the fact that normalized induction commutes with
duality. The third is the adjunction of the functors $i_P^G$ and $r_P^G$ and the last is again Lemma
\ref{VWWV}. If we assume that $r_{\bar P}^G \tilde \pi= (r_P^G \pi)^\sim$, we then obtain (\ref{2AD})
with $\tilde \pi$ instead of $\pi$. To obtain (\ref{2AD}) with $\pi$ which is not of the form
$\tilde \sigma$, we proceed as follows. The representation $\pi$ injects into $\tilde{\tilde \pi}$, and
the cokernel $\pi_1$ of this injection injects into $\tilde{\tilde \pi}_1$. We obtain an exact sequence
\[
0 \hookrightarrow \pi \hookrightarrow \tilde{\tilde \pi}
\rightarrow \tilde{ \tilde \pi}_1.
\]

By the left exactness of the functors $\Hom_G(X,\bullet)$ and $r_{\bar P}^G$, this yields a commutative
diagram with exact rows and whose vertical arrows are the isomorphisms (\ref{2AD}) established above for
$\tilde {\tilde \pi}$ and $\tilde{\tilde \pi}_1$:
\begin{equation*}
\begin{CD}
0  @>>>   \Hom_G(i_P^G\tau, \pi) @>>>   \Hom_G(i_P^G\tau,
\tilde{\tilde \pi})    @>>>   \Hom_G(i_P^G\tau, \tilde{\tilde \pi}_1)    \\
& &            & &                          @VVV   @VVV  \\
0  @>>>  \Hom_M(\tau, r_{\bar P}^G \pi)  @>>> \Hom_M(\tau, r_{\bar
  P}^G \tilde{\tilde \pi})  @>>> \Hom_M(\tau, r_{\bar P}^G
\tilde{\tilde \pi}_1).
\end{CD}
\end{equation*}
We easily deduce that $ \Hom_G(i_P^G\tau, \pi)$ and $ \Hom_M(\tau, r_{\bar P}^G \pi)$ are isomorphic and
that this isomorphism is natural.

Conversely, if we assume that the second adjunction theorem is established, then we obtain
\[
\Hom_M(\tau,r_{\bar P}^G( \tilde \pi))\simeq
\Hom_G(i_P^G(\tau),\tilde \pi) \simeq  \Hom_M (\tau,(r_P^G(\pi))\,
\tilde{}\, ),
\]
the second isomorphism being (\ref{rpdual}). Since this is true for all $\tau \in \caM(M)$, Principle
\ref{egalitecat} gives us the existence of a natural isomorphism
$r_{\bar P}^G \tilde \pi\simeq   (r_P^G \pi)^\sim$.
\end{proof}
\bigskip

Let us now prove this theorem. We reformulate it in the following form.

\begin{prop}
There exists an $M$-equivariant duality
\[
\bilo_P \colon r_P^G (V) \times  r_{\bar P}^G (\widetilde V) \rightarrow \bbC
\]
satisfying, for all $v \in V$, for all $\lambda \in \widetilde V$, for all $t \in C_A^{++}$, for all
$m \in \bbN$ large enough,
\begin{equation} \label{casselmanpairing}
\bil{r_P^G(\pi)(t^m).\bar v}{\bar \lambda}_P=
\delta_P^{1/2}(t^m)\lambda(\pi(t^m)\cdot v),
\end{equation}
where $\bar v=j_N(v)$, $\bar \lambda=j_{\overline{N}}(\lambda)$.

This duality induces an isomorphism $r_{\bar P}^G(\widetilde V)\simeq r_{P}^G(V)^\sim$.

For all $v \in V$, for all $\lambda \in \widetilde V$, there exists $\epsilon>0$ such that for all
$t \in  C_A^+(\epsilon)$,
\begin{equation} \label{casselmanpairingadm}
\bil{r_P^G(\pi)(t).\bar v}{\bar \lambda}_P=
\delta_P^{1/2}(t)\lambda(\pi(t)\cdot v).
\end{equation}
\end{prop}

\begin{proof}
First, we must find a duality between $r_P^G \pi$ and $r_{\bar P}^G \tilde \pi$. The respective spaces of
these representations are $V_N$ and $ \widetilde V_{\bar N}$. Let $K$ be a compact open subgroup of $G$
admitting an Iwahori decomposition with respect to $P$. We have, in the proof of Jacquet's lemma,
established a decomposition $V^K=V^K_0\oplus V^K_*$ such that $V^K_*\simeq V_N^{K_M}$, and in the same
way, we have a decomposition $\widetilde V^K=\widetilde V^K_0\oplus \widetilde V^K_*$ with
$\widetilde V^K_*\simeq \widetilde V_{\bar N}^{K_M}$. Note that $V_*^K$ is orthogonal to
$\widetilde V_0^K$ and that $V_0^K$ is orthogonal to $\widetilde V_*^K$. Indeed, let $v\in V_*^K$,
$\lambda \in \widetilde V_0^K$, and write $v=a_{t^m,K}\cdot v'$, with the notation of Section
\ref{LemmeJacquet2}, where $m$ is a sufficiently large integer. We then have
\[
\lambda(v)=  \lambda(a_{t^m,K}\cdot v') =((a_{t^m,K})\, \check{}\,
\cdot  \lambda)(v')= (a_{t^{-m},K}\cdot \lambda)(v').
\]
Here, the element $t$ is in $C_A^{++}$, a subset defined relative to the parabolic subgroup $P$, and more
precisely to its unipotent radical $N$. The subset which corresponds to the parabolic subgroup
$\bar P=M\overline{N}$ is
\[
C_A^{--}:= \{ a^{-1}|\, a \in C_A^{++} \}.
\]
The subspace $\widetilde  V_0^K$ of $\widetilde  V^K$ is therefore the kernel of $a_{t^{-m},K}$, for $m$
sufficiently large. We thus have $a_{t^{-m},K}\cdot \lambda=0$, which shows that $V_*^K$ is orthogonal to
$\widetilde V_0^K$. The other assertion is proved in the same way.

The restriction of $j$ to $V_*^K$ is an isomorphism onto $V_N^{K_M}$. Let $s_P^K$ denote its inverse.
Similarly, let $s_{\bar P}^K$ denote the isomorphism from $\widetilde V_{\overline{N}}^{K_M}$ onto
$\widetilde V_*^K$ inverse to $j$. Let $\bar v\in V_N^{K_M}$ and
$\bar \lambda \in \widetilde  V_{\overline{N}}^{K_M}$, and denote $v=s_P^K (\bar v)$ and
$\lambda=s_{\bar P}^K (\bar \lambda)$. Set
\[
\bil{\bar v}{\bar \lambda}= \bil{v}{\lambda},
\]
where the second duality bracket is the natural duality between $V$ and $\widetilde V$. It is clear that
this indeed defines a duality between $V_N^{K_M}$ and $ \widetilde V_{\bar N}^{K_M}$. Indeed, by Lemma
\ref{contrag} we have $(\widetilde V)^K=(V^*)^K=(V^K)^*$ and thus this follows from the fact that $V_0^K$
is orthogonal to $\widetilde V_*^K$ and that $V_*^K$ is orthogonal to $\widetilde V_0^K$. This also shows
that $(\widetilde V)^{K_M}_{\overline{N}}$ is indeed the entire dual of $V_N^{K_M}$, i.e., that
$(V_N^{K_M})^*= \widetilde V_{\bar  N}^{K_M}$.

As $K$ runs through the family of compact open subgroups of $G$ admitting an Iwahori decomposition with
respect to $P$, the family of $K_M$ runs through a basis of neighborhoods of the identity in $M$. We thus
have $V_N=\bigcup_{K} V_N^{K_M}$ and
$\widetilde V_{\overline{N}}= \bigcup_{K} \widetilde V_{\bar N}^{K_M}$. We must verify that the dualities
defined above are compatible, and can thus extend to a duality between $V_N$ and
$(\widetilde V)_{\overline{N}}$. To do this, it suffices to look at what happens for two compact open
subgroups $K'\subset K$ of $G$. We of course have in this case $V^K \subset V^{K'}$. Let
$v \in V_0^{K'}\cap V^K$. We thus have for $m$ sufficiently large,
\[
a_{t^m,K}\cdot v= e_K*\delta_{t^m}*e_K\cdot v=e_K*e_{K'}
*\delta_{t^m}*e_{K'}*e_K\cdot v= e_K*a_{t^m,K'}\cdot v=0.
\]
This shows that $V_0^{K'}\cap  V^{K} \subset V_0^{K} $. Similarly, if $v\in V_*^K$, let us write
$v=a_{t^m,K}\cdot v'$ for a certain $v'\in V^K$. We then have
\[
v= e_K*\delta_{t^m}*e_K\cdot v'= e_K*e_{K'}*\delta_{t^m}*e_{K'}*e_K\cdot v'=  e_K*a_{t^m,K'}\cdot v',
\]
whence $V^K_*\subset e_K\cdot V^{K'}_*$. We want to show that these inclusions are equalities. Let $n$ be
an integer greater than $b$, the constant of Proposition \ref{LemmeJacquet2}. Assuming that $n$ is large
enough, so that $t^{-n}K_{\overline{N}}t^n \subset K'$ and $t^{n}K_{N}t^{-n} \subset K'$ and calculating
as in the proof of the second lemma of \ref{HCDH}, we obtain
\[
a_{t^n,K}*  a_{t^n,K'}= e_K* a_{t^{2n},K'}, \quad  a_{t^n,K'}*  a_{t^n,K}=  a_{t^{2n},K'}*e_K
\]
which easily shows that we have
\[
V_0^{K'}\cap  V^{K} =V_0^{K}, \quad  V^K_*=  e_K\cdot V^{K'}_*.
\]
We will now show that on $V_N^{K_M}$,
\begin{equation} \label{spk}
 e_K \circ s_P^{K'}\circ i=s_P^K,
\end{equation}
where $i$ is the inclusion of $V_N^{K_M}$ into $V_N^{K'_M}$. Let then $\bar v \in V_N^{K_M}$, and write
\[
\bar v =j_N(\pi(a_{t^n,K})\cdot v)
\]
for a certain $v \in V^K$, so that
\[
s_P^K(\bar v)=\pi(a_{t^n,K})\cdot v.
\]
We have by (\ref{jvk}),
\[
j_N(\pi(a_{t^n,K})\cdot v) = \pi_N(t^n)j_N(v) =   j_N(\pi(a_{t^n,K'})\cdot v)
\]
whence
\begin{align*}
   &e_K \circ s_P^{K'}(\bar v) =   e_K \circ s_P^{K'}(j_N(\pi(a_{t^n,K'})\cdot v))=
 e_K \cdot (\pi(a_{t^n,K'})\cdot v)\\
&= \pi(e_K*a_{t^n,K'}*e_K)\cdot v=\pi(a_{t^n,K})\cdot v= s_P^K(\bar v).
\end{align*}
This shows (\ref{spk}).

We will also use that for all $\bar v \in V_N^{K_M}$
\begin{equation} \label{spk2}
 (1-e_K)s_P^{K'}(\bar v) \in V_0^{K'}=  V^{K'} \cap \ker j_N .
\end{equation}
Indeed, by (\ref{spk}),
\[
(1-e_K)s_P^{K'}(\bar v) =s_P^{K'}(\bar v)-e_K s_P^{K'}(\bar v)=
s_P^{K'}(\bar v)-s_P^{K}(\bar v).
\]
Now,
\[
j_N(s_P^{K'}(\bar v)-s_P^{K}(\bar v))=\bar v-\bar v=0,
\]
which proves the assertion.

We can now calculate:
\begin{align*}
& \bil{ s_P^K(\bar v) } { s_{\bar P}^K (\bar \lambda) }-  \bil{
  s_P^{K'}(\bar v) }{s_{\bar P}^{K'}(\bar \lambda) }     \\
=&\bil {e_K \cdot s_P^{K'} (\bar v) }{ e_K \cdot s_{\bar P}^{K'} (\bar
  \lambda) } -   \bil{ s_P^{K'}(\bar v) }{ s_{\bar P}^{K'}(\bar
  \lambda) }\\
=& \bil{ s_P^{K'} (\bar v) }{ e_K \cdot  s_{\bar P}^{K'}(\bar \lambda)-
  s_{\bar P}^{K'}(\bar \lambda) } +
  \bil{ e_K\cdot s_P^{K'}(\bar v)-  s_P^{K'}(\bar v) }{ e_K \cdot
    s_{\bar P}^{K'}(\bar \lambda)}\\
=& \bil{ s_P^{K'} (\bar v) }{ (e_K-1) \cdot  s_{\bar P}^{K'}(\bar \lambda)} +
  \bil{ (e_K-1)\cdot s_P^{K'}(\bar v)}{ e_K \cdot
    s_{\bar P}^{K'}(\bar \lambda)}\\
\end{align*}
Now
$ \bil{ (e_K-1)\cdot s_P^{K'}(\bar v)}{ e_K \cdot s_{\bar P}^{K'}(\bar \lambda)}= \bil{ e_K*(e_K-1)\cdot s_P^{K'}(\bar v)}{ s_{\bar P}^{K'}(\bar \lambda)}=0$
and $\bil{ s_P^{K'} (\bar v) }{ (e_K-1) \cdot  s_{\bar P}^{K'}(\bar \lambda)} =0$ by (\ref{spk2}) and the
orthogonality of $\widetilde V_0^{K'}$ and $ V_*^{K'}$.

This shows that the duality $\bilo_P$ is in fact defined independently of the choice of $K$.

Let us now verify that this duality is $M$-equivariant. Let $\bar v\in V_N$ and
$\bar \lambda \in \widetilde V_{\overline{N}}$, and $m \in M$. We will now show that
\[
\bil{\pi_N(m)\cdot \bar v}{\tilde \pi_{\overline{N}}(m)\cdot \bar \lambda}_P=
\bil{\bar v}{ \bar \lambda}_P.
\]
Let us choose a compact open subgroup as above such that $K_M$ fixes $\bar v$ and $\bar \lambda$ and let
$v\in V_*^K$, $\lambda \in \widetilde V^K_*$ lifting respectively $\bar v\in V_N$ and $\bar \lambda$, so
that $\bil{\bar v} { \bar \lambda}_P=\bil{v} { \lambda}$. Since
$\pi_N(m)\cdot \bar v= \overline {\pi(m)\cdot v }$,
$\tilde \pi_{\overline{N}}(m)\cdot \bar \lambda = \overline {\tilde \pi(m)\cdot \lambda }$, we see that
we will have the desired equality if $\pi(m)\cdot v \in  V_*^{K'}$, and
$\tilde \pi(m)\cdot \lambda \in \widetilde  V_*^{K'}$, for another compact open subgroup $K'$ of $G$
satisfying the same conditions as $K$, since we have shown that the duality does not depend on the choice
of $K$. Indeed, we then have
\[
\bil{\pi_N(m)\cdot \bar v}{\tilde \pi_{\overline{N}}(m)\cdot \bar \lambda}_P=
\bil{\pi(m)\cdot v}{\tilde \pi(m)\cdot \lambda}= \bil{ v}{\lambda}_P.
\]
Now we also have, by taking $K'=mKm^{-1}$, and using the fact that $t$ is central in $M$,
\begin{align*}
\delta_m
*a_{t,K}&=\delta_m*e_K*\delta_t*e_K=e_{K'}*\delta_m*\delta_t*e_K
=e_{K'}*\delta_t*\delta_m*e_K\\
&=e_{K'}*\delta_t*e_{K'}*\delta_m= a_{t,K'}*\delta_m.
\end{align*}
This easily implies that $\pi(m)\cdot v \in  V_*^{K'}$, and
$\tilde \pi(m)\cdot \lambda \in \widetilde  V_*^{K'}$. We immediately deduce that
\[
\bil{r_P^G(\pi)(m)\cdot \bar v}{r_{\bar P}^G(\tilde \pi)(m)\cdot \bar \lambda}_P=
\bil{\bar v}{ \bar \lambda}_P.
\]
Indeed, the normalizations by $\delta_P^{1/2}(m)$ and $\delta_{\bar P}^{1/2}(m)$ cancel each other out
(Lemma \ref{calcfonctmod}).

It remains to verify (\ref{casselmanpairing}) and (\ref{casselmanpairingadm}). Let us fix a compact open
subgroup $K$ of $G$ admitting an Iwahori decomposition with respect to the standard parabolic subgroups
and fixing $v$ and $\lambda$. We then have
\begin{align*}
\delta_P^{1/2}(t^m)\lambda(\pi(t^m)\cdot v)
&=\delta_P^{1/2}(t^m)(e_K \cdot \lambda) (\pi(t^m)\cdot e_K\cdot  v)\\
&=   \delta_P^{1/2}(t^m)\lambda(a_{t^m,K}\cdot v)\\
&=  \delta_P^{1/2}(t^m)\bil{a_{t^m,K}\cdot v}{\lambda}\\
&=  \delta_P^{1/2}(t^m) \bil{a_{t^m,K}\cdot v}{ s_{\bar  P}^K (j_{\overline{N}}(\lambda))}\\
&=\delta_P^{1/2}(t^m)\bil{ s_{ P}^K (j_{N}(a_{t^m,K}\cdot v))}{
  s_{\bar P}^K (j_{\overline{N}}(\lambda))}
\end{align*}
For the last two equalities, we use the fact that for $m$ large enough, $a_{t^m,K}\cdot v \in V_*^K$, the
orthogonality of $\widetilde V_0^K$ and $V_*^K$ (resp. of $V_0^K$ and $\widetilde V_*^K$) and the fact
that $s_P^K$ is a section of $j_N$ on $V^K$ (resp. $s_{\bar P}^K$ a section of $j_{\overline{N}}$ on
$\widetilde V^K$). By the definition of the duality $\bilo_P$, we obtain
\begin{align*}
&\delta_P^{1/2}(t^m)\lambda(\pi(t^m)\cdot v)
=\delta_P^{1/2}(t^m)\bil{j_{N}(a_{t^m,K}\cdot v)}{j_{\bar
    N}(\lambda)}_P\\
&= \bil{  \delta_P^{1/2}(t^m)  \pi_N(t^m)\cdot \bar  v}{\bar \lambda}_P\\
&= \bil{  r_P^G(\pi)(t^m)\cdot \bar  v}{\bar \lambda}_P
\end{align*}
We used (\ref{jvk}), an equality proved in \ref{LemmeJacquet}. The proof of (\ref{casselmanpairingadm})
is identical, using the fact (third remark of \ref{LemmeJacquet2}) that for $t \in C_A^+(\epsilon)$,
\[
V^K=V_0^K\oplus V_*^K= \ker a_{t,K} \oplus \im a_{t,K}.
\]
This completes the proof of the proposition.
\end{proof}

\subsection{Second adjunction and completion}\label{secadjcomp}

We provide an interpretation of the second adjunction theorem in terms of module completion ({\sl cf.}
Section \ref{compmod}).

\begin{thm}
Let $P=MN$ be a parabolic subgroup of $G$ and let $\bar P=M\overline{N}$ be the opposite parabolic
subgroup. Then, for any representation $(\pi,V)$ of $\caM(G)$, we have a natural isomorphism
\[
\phi:\bar V^N \simeq \overline{V_{\overline{N}}}.
\]
For all $v \in \bar V^N$, $\phi(v)$ is characterized by the following property:

--- for any compact open subgroup $K$ of $G$ admitting an Iwahori decomposition with respect to $P$,
$\phi(e_K \cdot v)=e_{K_M}\cdot \phi(v)$.
\end{thm}

\begin{proof}
By definition, the module $\overline{V_{\overline{N}}}$ is the projective limit of the
$e\cdot V_{\overline{N}}$, where $e$ runs through the set of idempotents of $\caH(M)$. Since the
$e_{K_M}$ form a directed system of idempotents of $\caH(M)$ as $K$ runs through the set of compact open
subgroups of $G$ admitting an Iwahori decomposition with respect to $P$, we have
\[
\overline{V_{\overline{N}}}= \varprojlim_{K} e_{K_M}\cdot V_{\overline{N}} =
\varprojlim_{K} V_{\overline{N}}^{K_M}.
\]
In the previous section (the roles of $\bar P$ and $P$ must be reversed), we exhibited a subspace $V^K_*$
of $V^K$ such that $j_{\overline{N}}$ realizes an isomorphism $V_{\overline{N}}^{K_M}\simeq V^K_*$. We
thus obtain
\[
\overline{V_{\overline{N}}} = \varprojlim_{K} V_{\overline{N}}^{K_M} \simeq
\varprojlim_{K} V^K_* \subset \varprojlim_{K} V^K= \bar V.
\]
This identifies the space $\overline{V_{\overline{N}}}$ with the space:
\[
\overline{V_*}:=\{ v \in \bar V\, | \, e_K\cdot v \in   V_*^K, (\forall K)       \}.
\]

According to the second remark of \ref{LemmeJacquet2}, for any sufficiently large compact open subgroup
$C$ of $N$, we have $V_*^K=e_Ke_C\cdot V$. We deduce that $v\in  \overline{V_*}$ if and only if for all
$K$ as above and for any sufficiently large compact open subgroup $C$ of $N$ (depending explicitly on
$K$, see Remark \ref{LemmeJacquet2}), $e_K\cdot v \in e_Ke_C\cdot V$. This is further equivalent to:

- For any compact open subgroup $C$ of $N$, for any sufficiently small $K$ as above (depending explicitly
on $C$ (this again provides a directed system of idempotents)),
$e_K\cdot v \in e_Ke_C\cdot V= e_Ce_K\cdot V$.

Thus: $v\in  \overline{V_*}$ if and only if for any compact open subgroup $C$ of $N$,
$v \in e_C\cdot \bar V$, and since $N$ is the union of its compact subgroups, we obtain
$ \overline{V_{\overline{N}}}=\overline{V_*}\simeq  \bar V^N$.
\end{proof}

We obtained this theorem as a consequence of Jacquet's lemma \ref{LemmeJacquet2}. Let
$(\pi, V)\in \caM(G)$. By applying the result to $(\tilde \pi,\widetilde V)$, we obtain an isomorphism
\[
\phi \colon \overline{\widetilde V}^N \rightarrow \overline{\widetilde V_{\overline{N}}}.
\]
Now $\overline{\widetilde V}^N \simeq (V^*)^N \simeq (V_N)^*$ ({\sl cf.} (\ref{dualetcompeq})), and thus
$\widetilde V_{\overline{N}}$ is identified with the smooth part of $(V_N)^*$ ({\sl cf.} \ref{lisscomp}),
which by definition is $\widetilde{V_N}$. We thus recover the canonical isomorphism:
\[
r_{\bar P}^G \tilde \pi \simeq ( r_P^G\pi)^\sim.
\]
The proof given here is of course essentially equivalent to the one given in \ref{secondeadjonction}, but
the formalism of module completions allows for a more streamlined exposition. It also allows us to obtain
the second adjunction theorem in a slightly more elegant way. Indeed, for any representation $(\tau,E)$
in $\caM(M)$ and for any representation $(\pi,V)$ in $\caM(G)$, we have, writing the equality sign for
natural isomorphisms:
\begin{align*}
 &\Hom_G(i_P^G(E),V)= \Hom_G(\ind_P^G(\caF_P^M( \delta_P^{-1/2} \otimes E)),V)\\
&=\Hom_G( P_P^G( \caF_P^M(\delta_P^{1/2} \otimes E)),V)\\
&=\Hom_P(\caF_P^M(\delta_P^{1/2} \otimes E), \check{}\, \caF_P^G(V))\\
&=\Hom_P(\caF_P^M(\delta_P^{1/2} \otimes E), \Hom_G(\caH(G),V)_{\caH(P)})\\
&=\Hom_P(\caF_P^M(\delta_P^{1/2} \otimes E), \bar V_{\caH(P)})=\Hom_P(\caF_P^M(\delta_P^{1/2} \otimes E), \bar V) \\
&=\Hom_M(\delta_P^{1/2} \otimes E, \bar V^N )=\Hom_M( E,
\delta_P^{-1/2} \otimes\bar V^N )\\
&= \Hom_M( E, \delta_P^{-1/2} \otimes \overline{V_{\overline{N}}} )=\Hom_M(
E, \delta_{\bar P}^{1/2} \otimes V_{\overline{N}})\\
&= \Hom_M( E, r_{\bar P}^G V)
\end{align*}

We have successively used the identity between the compact induction functor and the functor $P$
(\ref{indcomp}), the adjunction of $P$ and the pseudo-forgetful functor and the expression of the
pseudo-forgetful functor in terms of the completion (\ref{Oublietadjoints}), the fact that a morphism
from a non-degenerate module to an arbitrary module always takes values in the non-degenerate part, the
adjunction between the forgetful functor $\caF_P^M$ and the functor consisting of taking the
$N$-invariants (which is easily verified), and finally the isomorphism
$\bar V^N \simeq \overline{V_{\overline{N}}}$ proved above.

\section{The center of $\caM(G)_\Omega$}

\subsection{A progenerator of $\caM(G)_\Omega$} \label{progenMOm}

Let $\Omega$ be a connected component of $\Omega(G)$ and let $(M,(\rho,W))$ be a cuspidal datum whose
class is in $\Omega$. We may assume that $M$ is a standard Levi subgroup of $G$, a Levi factor of the
standard parabolic subgroup $P=MN$. Let $D$ be the inertial class in $M$ of $(\rho,W)$, i.e.,
$D=[M,(\rho,W)]_M$. Our goal is to find a small progenerator of the category $\caM(G)_\Omega$. Let
$\Pi_D$ be the progenerator of the category $\caM(M)_D$ constructed in \ref{progenpi}. A natural
candidate \index[not]{ZZPi_Omega@$\Pi_\Omega$} is $\Pi_\Omega= i_P^G(\Pi_D)$. The first step is to show
that $\Pi_\Omega$ is independent of the choices made. To do this, we will use a result which will only be
proved in \ref{corang1} (where the reader is also invited to refer for notation).

\begin{lemme}
Assume $M$ is maximal (see \ref{corang1}), and let $w$ be the non-trivial element of $\caW(M,*)$. Set
$M'=w\cdot M$ and $\rho'=w\cdot \rho$, $D'=[M',\rho']_{M'}$. Let $P'$ denote the standard parabolic
subgroup with Levi factor $M'$. Then
\[
i_P^G (\Pi_D) \simeq i_{P'}^G (\Pi_{D'}).
\]
\end{lemme}

It is important to note that the isomorphism above is not canonical. On the other hand, the assertion
would be obvious with $w\cdot P$ instead of $P'$. But here, $w\cdot P=\bar P'$.

\begin{proof}
If $w\cdot(M,D)=(M,D)$, this is trivial. Otherwise, we are under the hypotheses of Corollary
\ref{corang1}. The equivalence of categories of this corollary shows that we must verify that
$r_Di_{P'}^G (\Pi_{D'})=\Pi_{D}$, but this is immediate from the geometric lemma \ref{geolemmacusp}.
\end{proof}

\bigskip

We can now dispense with the maximality hypothesis on $M$. \needspace{163pt}
\begin{prop}
Let $(M,D)$ be as above, and let $w \in \caW(M,*)$. Set $M'=w\cdot M$, $D'=w\cdot D$ and let $P'$ be the
standard parabolic subgroup with Levi factor $M'$. Then
\[
i_P^G (\Pi_D) \simeq i_{P'}^G (\Pi_{D'}).
\]
\end{prop}

\begin{proof}
We use the fact that we can write $w$ as a product of elementary transformations (Lemma
\ref{compindsup}). This reduces us to the case where $w$ is itself elementary. There then exists a Levi
subgroup $L$ of $G$ containing $M$ as a maximal Levi subgroup. Let $Q=LU$ be the standard parabolic
subgroup with Levi factor $L$. By the lemma above
\[
i_{P\cap L}^L (\Pi_D) \simeq i_{P'\cap L}^L (\Pi_{D'}).
\]
We deduce
\[
i_P^G(\Pi_D)=i_Q^G  i_{P\cap L}^L (\Pi_D) \simeq i_Q^G i_{P'\cap
L}^L (\Pi_{D'})= i_{P'}^G (\Pi_{D'}).
\]
\end{proof}

\begin{cor}
Let $(M,D)$ be as above, and let $P' \in \caP(M)$. Then $i_P^G(\Pi_D)\simeq i_{P'}^G(\Pi_D)$.
\end{cor}

\begin{proof}
Since $P'$ contains $A_\emptyset$, there exists $w\in N_G(A_\emptyset)$ such that $P''=w\cdot P'$ is
standard; its standard Levi factor is $M''=w\cdot M$, so that $w$ defines an element of $\caW(M,*)$. We
then have
\[
i_{P'}^G(\Pi_D)= i_{w^{-1}\cdot P''}^G(\Pi_D) \simeq  i_{P''}^G(\Pi_D^{w}) =  i_{P''}^G(\Pi_{D''})
\]
where $D''=w\cdot D$, and the previous proposition gives $i_{P''}^G(\Pi_{D''})\simeq i_P^G(\Pi_D)$.
\end{proof}

The consequence of all this is that the (isomorphism class of the) representation $\Pi_\Omega$ is in fact
independent of the datum $(M,D)$ chosen to construct it. We can now prove the desired result:
\begin{thm}
The representation $\Pi_\Omega$ is a finitely generated progenerator of $\caM(G)_\Omega$.
\end{thm}
\begin{proof}
Bernstein's second adjunction theorem asserts that $i_P^G$ is the left adjoint of the functor
$r_{\bar P}^G$. The latter being exact, we immediately deduce that $i_P^G$ preserves projectives. This
shows that $\Pi_\Omega$ is projective. We also showed in \ref{noether} that $i_P^G$ preserves finitely
generated representations. It follows that $\Pi_\Omega$ is finitely generated. Now a finitely generated
projective object is small in the sense of category theory (\ref{rmqfonct}). It remains to show that
$\Pi_\Omega$ is a generator of $\caM(G)_\Omega$, and to do this, it suffices to see that any irreducible
representation $(\pi,V)$ of $\caM(G)_\Omega$ is a quotient of $\Pi_\Omega$. Indeed, we want to show that
for any representation $(\pi',V')$ in $\caM(G)_\Omega$, $\Hom_G(\Pi_\Omega,V')$ is not zero. But we know
that $(\pi',V')$ admits an irreducible subquotient $(\pi,V)$, say that $(\pi,V)$ is a quotient of a
subrepresentation $(\pi_1,V_1)$. If we know that $\Hom_G(\Pi_\Omega,\pi)$ is not zero, by the
projectivity of $\Pi_\Omega$, $\Hom_G(\Pi_\Omega,\pi_1)$ is not zero, and a fortiori
$\Hom_G(\Pi_\Omega,\pi')$. Suppose therefore that $(\pi,V)$ is irreducible in $\caM(G)_\Omega$. Let
$P'=M'N'$ be a standard parabolic subgroup of $G$ such that $r_{\bar P'}^G(\pi)$ is supercuspidal. The
functor $i_{P'}^G$ is the left adjoint of $r_{\bar P'}^G$, and we have the adjunction morphism associated
to the identity of $r_{\bar P'}^G(\pi)$:
\[
\beta_\pi \colon   i_{P'}^G r_{\bar P'}^G(\pi) \rightarrow \pi.
\]
Since $\pi$ is irreducible, this morphism is surjective (because it is non-zero) and by the exactness of
the functor $i_{P'}^G$ there exists a composition factor $\rho'$ of $r_{\bar P'}^G(\pi)$ such that
$ i_{P'}^G (\rho')$ maps surjectively onto $\pi$. Let $D'$ be the inertial class of $\rho'$. By
definition, $\Omega$ is also the connected component of $\Omega(G)$ associated to $(M',D')$. Since
$\rho' \in D'$, and $\Pi_{D'}$ is a progenerator of $\caM(M')_{D'}$, there exists a surjective morphism
$\Pi_{D'}  \rightarrow \rho'$. The functor $ i_{P'}^G$ being exact, the image morphism
$i_{P'}^G (\Pi_{D'})  \rightarrow i_{P'}^G \rho'$ remains surjective, and by composition, we obtain the
desired surjective morphism $i_{P'}^G (\Pi_{D'})  \rightarrow \pi$. The proposition above allows us to
conclude.
\end{proof}

\subsection{The equivalence $\caM(G)_\Omega \simeq \caM(\caR_\Omega)_d$} \label{eQuiv}
We retain the notation of the previous section. Set
\[
\caR_\Omega=\End_G(\Pi_\Omega).
\]
\index[not]{R_Omega@$ \caR_\Omega$} Exactly as in \ref{progenpi}, we introduce the functors
\[
F_{\Pi_\Omega}  \colon \caM(G)_\Omega \rightarrow \caM(\caR_\Omega)_d,  \quad (\pi,V) \mapsto \Hom_G(\Pi_\Omega,\pi),
\]
\[
G_{\Pi_\Omega}  \colon   \caM(\caR_\Omega)_d\rightarrow \caM(G)_\Omega, \quad M \mapsto M\otimes_{ \caR_{\Omega}} \Pi_\Omega.
\]

The proof of Theorem \ref{progenpi} can be copied word for word and we obtain:
\begin{thm}
The functors are quasi-inverses and realize an equivalence of categories between $\caM(G)_\Omega$ and
$ \caM(\caR_\Omega)_d$.
\end{thm}

\subsection{An estimate of the length of induced supercuspidal representations}\label{longindsc}

Let $(M,(\rho,W))$ be a cuspidal datum of $G$, where $M$ is a standard Levi subgroup of $G$. Let $P$ be
the standard parabolic subgroup of $G$ with Levi factor $M$ and as in the previous section, set
\[
[M,(\rho,W)]_M=D, \quad   [M,(\rho,W)]_G=\frs.
\]
Let $\Pi_D$ be the progenerator of $\caM(M)_D$ constructed in \ref{progenpi}, and let
$\Pi_\Omega=i_P^G \Pi_D$, which we have just shown to be a progenerator of the category $\caM(G)_\frs$.
Let $r_D$ denote the composition of the functor $r_ P^G$ from $\caM(G)$ to $\caM(M)$ and the projection
onto the component $\caM(M)_D$ of $\caM(M)$.

\begin{prop}
The functor $r_D : \caM(G)_\frs \longrightarrow \caM(M)_D$ maps non-zero objects to non-zero objects.
\end{prop}

\begin{proof}
Let $\pi$ be a (non-zero) smooth representation of $\caM(G)_\frs$. Since $i_{\bar P}^G\Pi_D$ is a
progenerator of $\caM(G)_\frs$ (Corollary \ref{progenMOm}), we have by the second adjunction theorem,
\[
\{0\}\neq \Hom_G(i_{\bar P}^G\Pi_D, \pi)\simeq \Hom_M(\Pi_D, r_P^G\pi)=  \Hom_M(\Pi_D, r_D \pi)
\]
and thus $r_D\pi \neq 0$.
\end{proof}

\begin{cor}
The length of $i_P^G(\rho)$ is at most $|W(D)|$, where
\[
W(D)=\{ w \in W(M,M), w\cdot D =  D    \}/W_M.
\]
In particular, if $ |W(D)|=1$, $i_P^G(\rho)$ is irreducible.
\end{cor}

\begin{proof}
The proposition and the exactness of the functor $r_D$ show that the length of $i_P^G(\rho)$ is less than
or equal to the length of $r_Di_P^G(\rho)$. But the geometric lemma shows that this latter representation
is of length at most $|W(D)|$.
\end{proof}

\subsection{The center of $\caM(G)_\Omega$} \label{centreMGOm}

The notation is the same as in the previous sections. Set $\caR_D= \End_M(\Pi_D)$. For all
$\phi \in \caR_D$, $i_P^G(\phi)$ is in
\[
\End_G(i_P^G(\Pi_D))= \End_G(\Pi_\Omega)=\caR_\Omega.
\]
This defines by functoriality an algebra morphism
\[
i_P^G\colon \caR_D  \rightarrow   \caR_\Omega.
\]
This algebra morphism is injective, since the functor $i_P^G$ is faithful. We can therefore identify
$\caR_D$ with a subalgebra of $\caR_\Omega$. In particular $\caR_\Omega$ is an $\caR_D$-bimodule.

The $\caR_D$-bimodule structure of $\caR_\Omega$ gives an induction functor:
\begin{equation}\label{iDuc}
  \caM(\caR_D)_d\longrightarrow   \caM(\caR_\Omega)_d, \quad M
\mapsto M \otimes_{\caR_D} \caR_\Omega.
\end{equation}
We will now show that this functor, through the equivalences of categories
\begin{equation}\label{eQui}
 \caM(\caR_D)_d\simeq \caM(M)_D \quad \text
  {  and }\quad   \caM(\caR_\Omega)_d  \simeq \caM(G)_\Omega,
\end{equation}
is isomorphic to the induction functor $i_P^G$.

Let $(\pi,V)$ be in $\caM(G)_\Omega$. By Bernstein's second adjunction theorem, we have a natural
isomorphism
\[
\Hom_G(\Pi_\Omega,\pi)=\Hom_G(i_P^G(\Pi_D),\pi)\simeq \Hom_M(\Pi_D,
r_{\bar P}^G (\pi))
\]
The left-hand side is a unital right $\caR_D$-module via the injection of $\caR_D$ into $\caR_\Omega$ and
the right-hand side is as well. The naturality of the adjunction isomorphism implies that this
isomorphism is an isomorphism of right $\caR_D$-modules. In other words, the forgetful functor
$$  \caM(\caR_\Omega)_d \rightarrow \caM(\caR_D)_d $$ induced by $\caR_D \hookrightarrow \caR_\Omega$
corresponds via the equivalences of categories (\ref{eQui}) to the functor $r_{\bar P}^G$. Now the
induction functor (\ref{iDuc}) is the left adjoint of the forgetful functor, and $i_P^G$ is the left
adjoint of $r_{\bar P}^G$. By the uniqueness up to isomorphism of the adjoint, this proves the assertion.

\bigskip

Let $W(D)$ \index[not]{W(D)@$W(D)$} denote the subgroup of $W(A_M)=N_G(M)/M$ stabilizing the
supercuspidal inertial class $D$.

\begin{lemme}
The right $\caR_D$-module $\caR_\Omega$ is free of rank $|W(D)|$. It is therefore free of finite rank
over the center $\frZ_D$ of $\caR_D$.
\end{lemme}

\begin{proof}
By Bernstein's second adjunction theorem, we have a natural isomorphism
\[
\caR_\Omega= \Hom_G(i_P^G(\Pi_D), i_P^G(\Pi_D))\simeq
\Hom_M(\Pi_D,r_{\bar P}^G  i_P^G (\Pi_D)).
\]
Denoting (in this paragraph) by $r_D$ \index[not]{r_D@$r_D$} the composition of the functor
$r_{\bar P}^G$ (and not $r_P^G$ as previously) with the projection onto the component $\caM(M)_D$ of
$\caM(M)$, we obtain a natural isomorphism
\begin{equation}\label{ROm}
\caR_\Omega \simeq \Hom_M(\Pi_D, r_D i_P^G(\Pi_D)).
\end{equation}
We have seen above that (\ref{ROm}) is an isomorphism of right $\caR_D$-modules.

By equation (\ref{WMWA}), Proposition \ref{geolemmacusp} $(iii)$, and the definition of $W(D)$, there
exists a filtration of $r_D  i_P^G(\Pi_D) $ whose successive quotients are the $w \cdot \Pi_D$,
$w \in W(D)$. Since all these $w \cdot \Pi_D$ are projective, because $\Pi_D$ is, this filtration splits
into a direct sum
\[
r_D  i_P^G(\Pi_D) = \bigoplus_{w \in W(D)}  w \cdot \Pi_D.
\]
Now, if $w \in W(D)$, $ w \cdot \Pi_D \simeq \Pi_D$.

We obtain
\[
\caR_\Omega \simeq  \Hom_M(\Pi_D, \bigoplus_{w \in W(D)}
\Pi_D)=\prod_{w \in W(D)} \Hom_M(\Pi_D,\Pi_D)=\caR_D^{|W(D)|}.
\]
Note that this isomorphism is an isomorphism of right $\caR_D$-modules, but that the left structure is
lost. The second assertion follows from Theorem \ref{centreA}.
\end{proof}

\bigskip

Let $\frZ_D$ and $\frZ_\Omega$ denote respectively the \index[not]{Z_Omega@$\frZ_\Omega$} center of
$\caR_D$ and $\caR_\Omega$. Let us recall the results of Sections \ref{centreA} and \ref{centreMpi}. For
any maximal ideal $I$ of $\frZ_D$, $\Pi_D/I\Pi_D$ is isomorphic to $m\sigma_I$, for a unique (up to
isomorphism) irreducible object $\sigma_I \in \caM(M)_D$ (the integer $m$ is a certain multiplicity
defined in \ref{centreMpi}). On the other hand $I \mapsto \sigma_I$ is a bijection between
$\mathrm{SpecMax}(\frZ_D)$ and $\mathbf{Irr}(M)_D$, thus equipping $\mathbf{Irr}(M)_D$ with the structure
of a complex affine variety whose polynomial algebra is $\frZ_D$. The equivalence of categories
\ref{eQui} maps $\caR_D/I\caR_D$ to
\[
\caR_D/I\caR_D \otimes_{\caR_D} \Pi_D \simeq     \Pi_D/I\Pi_D
\]
and thus, since $i_P^G$ corresponds to the induction functor (\ref{iDuc}), the representation
$i_P^G(\Pi_D/I\Pi_D)$ corresponds to the right $\caR_\Omega$-module
\[
\caR_D/I\caR_D \otimes_{\caR_D} \caR_\Omega \simeq \caR_\Omega/I\caR_\Omega.
\]

Let $z \in \frZ_\Omega$. Suppose that $i_P^G(\sigma_I)$ is irreducible, so that $z$ acts on it by a
certain scalar. Since
\[
i_P^G(\Pi_D/I\Pi_D) \simeq  i_P^G(m\sigma_I)\simeq m \,  i_P^G(\sigma_I),
\]
we see that $z$ acts by this same scalar on $ i_P^G(\Pi_D/I\Pi_D)$. From the above, $z$ acts by this same
scalar on $\caR_\Omega/I\caR_\Omega$. By the lemma above, $\caR_\Omega$ is a finitely generated free
module over $\frZ_D$. Let us fix a basis of $\caR_\Omega$ over $\frZ_D$. Then the action of $z$ on
$\caR_\Omega$ is expressed in this basis by a certain matrix $(z_{ij})$. The chosen basis gives us for
any maximal ideal $I$ of $\frZ_D$, a basis over $\bbC$ of $\caR_\Omega/I\caR_\Omega$ and the action of
$z$ on $\caR_\Omega/I\caR_\Omega$ is expressed in this basis by a certain matrix $(\bar z_{ij})$, where
$\bar z_{ij}=\Psi_I(z_{ij})$, $\Psi_I$ being the unique element of $\Hom_{\bbC-alg}(\frZ_D,\bbC)$ with
kernel $I$.

Recall that by Theorem \ref{irrindgen}, $i_P^G (\sigma_I)$ is irreducible for all $I$ in a dense Zariski
open set of $\mathrm{SpecMax}(\frZ_D)$. For these $I$, the matrix $(\bar z_{ij})$ is scalar, and it
follows that the matrix $(z_{ij})$ is scalar. This shows that $z \in \frZ_D$. We have thus obtained the
following lemma.

\begin{lemme}
The center $\frZ_\Omega$ of $\caR_\Omega$ is contained in the center $\frZ_D$ of $\caR_D$.
\end{lemme}

In fact, we easily obtain a slightly stronger result. The group $W(D)$ acts on the variety
$\mathbf{Irr}(M)_D$ and thus on its polynomial algebra $\frZ_D$. Indeed, if $\sigma \in D$, we saw above
that $z \in  \frZ_\Omega$ acts on $i_P^G \sigma$ by exactly the same scalar as its action on $\sigma$.
Corollary \ref{compindsup} asserts that the composition factors of $i_P^G \sigma$ are the same as those
of $i_{P'}^G \sigma$ for any other parabolic subgroup $P'$ of $G$ admitting $M$ as a Levi factor. This is
true in particular for $P'=w\cdot P$, $w \in W(D)$. It follows that $z$ acts on $i_P^G(w\cdot \sigma)$ by
the same scalar as the one by which it acts on $i_P^G(\sigma)$. We thus have
\[
\frZ_\Omega \subset\frZ_D^{W(D)}.
\]
In particular
\[
\frZ_\Omega= \scrC_{\frZ_D}(\caR_\Omega)
\]
where $ \scrC_{\frZ_D}(\caR_\Omega) =\{r\in \frZ_D \mid   ar=ra \; (\forall a\in \caR_\Omega) \}$ is
defined by the $\frZ_D$-bimodule structure of $\caR_\Omega$.

We must now see that the inclusion above is in fact an equality. Since $\caR_\Omega$ is an
$\caR_D$-bimodule, the center $\frZ_\Omega$ of $\caR_\Omega$ is characterized using only this structure:
\[
\frZ_\Omega \simeq \scrC_{\caR_D}(\caR_\Omega)=\{ r \in \caR_D \, |\,
ar=ra, \;  (\forall a \in  \caR_\Omega) \}.
\]

Let us examine more closely the left $\caR_D$-module structure of
$\caR_\Omega \simeq \Hom_M(\Pi_D, r_D i_P^G(\Pi_D))$. The action on $\Hom_M(\Pi_D, r_D i_P^G(\Pi_D))$
comes from that on $r_D i_P^G(\Pi_D)$ obtained by functoriality.

Let us resume the filtration of $r_D  i_P^G\Pi_D $ whose successive quotients are the $w \cdot \Pi_D$,
$w \in W(D)$. Since all these $w \cdot \Pi_D$ are projective, this filtration induces one on the space
\[
\Hom_M(\Pi_D, r_D i_P^G(\Pi_D)),
\]
compatible with the $\caR_D$-bimodule structure, with quotients isomorphic to
$\Hom_M(\Pi_D, w\cdot \Pi_D)$. This space admits a natural left $\caR_D$-module structure, obtained from
the action, twisted by $w$, of $\caR_D$ on $ w\cdot \Pi_D$.

We will use the following lemma:

\begin{lemme}
Let $A$ be a commutative unital ring, and $\scrC$ the category of (unital) $A$-bimodules. For all
$M \in \scrC$, we denote
\[
\scrC_{A}(M)=\{ a \in A \mid ma=am, \;  (\forall m \in  M) \}.
\]
Suppose that $M \in \scrC$ admits a filtration
\[
0 =M_{n+1} \subset    M_n \subset M_{n-1} \subset \cdots \subset   M_1
\subset M_0=M
\]
whose successive quotients $Q_i=M_{i}/M_{i+1}$, $i=0,\ldots, n$, satisfy
\[
\Hom_\scrC(Q_i,Q_j)=0  \quad (i\neq j).
\]
Then $\scrC_A(M)=\bigcap_{i=0}^n \scrC_A(Q_i)$.
\end{lemme}
\begin{proof}
We proceed by induction on $n$, the case $n=0$ being trivial. We may therefore assume that
$\scrC_A(M_1)=\bigcap_{i=1}^n \scrC_A(Q_i)$. From the short exact sequence
\[
0 \rightarrow M_1 \rightarrow M  \rightarrow Q_0 \rightarrow 0
\]
we easily deduce that $\scrC_A(M) \subset \scrC_A(M_1) \cap \scrC_A(Q_0)$. Conversely, if
$\beta \in \scrC_A(M_1) \cap \scrC_A(Q_0)$, then we can define
\[
Q_0=M/M_1 \rightarrow M_1, \quad m \mapsto \beta m - m \beta,
\]
an element of $\Hom_\scrC(Q_0,M_1)$. Now our hypothesis implies that this morphism is zero.
\end{proof}

We must now apply this to the category $\scrC$ of unital $\frZ_D$-bimodules, which we do by observing
that if $w \neq w'$ in $W(D)$, then
\[
\Hom_{\scrC}(\Hom_M(\Pi_D, w\cdot \Pi_D),
\Hom_M(\Pi_D, w'\cdot \Pi_D))=\{0\},
\]
which is clear: $W(D)$ acts faithfully on $\frZ_D$ (an element acting trivially on the torus $D$ would
act trivially on $\caX(M)$, which is excluded by \ref{varXG} and \ref{ReChi}), so that for $w\neq w'$ the
two $\frZ_D$-module structures on $\Hom_M(\Pi_D,w\cdot\Pi_D)$ and $\Hom_M(\Pi_D,w'\cdot\Pi_D)$ differ. We
then obtain
\[
\frZ_\Omega= \bigcap_{w \in W(D)} \scrC_{\frZ_D}( \Hom_M(\Pi_D, w\cdot \Pi_D))=\frZ_D^{W(D)}.
\]

From this, it is now easy to deduce the following theorem.
\begin{thm}
Let $\Omega$ be a connected component of the variety $\Omega(G)$. Then the center of the category
$\caM(G)_\Omega$ is isomorphic to the algebra of polynomial functions on the variety $\Omega$. Similarly,
the center $\frZ(G)$ of the category $\caM(G)$, the product of the $\frZ_\Omega$, is identified with the
algebra of polynomial functions on $\Omega(G)$.
\end{thm}
\begin{proof}
Let $\caO(\Omega)$ denote the algebra of polynomial functions on $\Omega$. With the preceding notation,
the variety $\Omega$ is the quotient of the variety $D$ under the action of the group $W(D)$ (see
\ref{var0G}). Now the algebra of polynomial functions on $D$ is isomorphic to $\frZ_D$, and thus
$\caO(\Omega) \simeq \frZ_D ^{W(D)}$. Now $ \frZ_D ^{W(D)} \simeq \frZ_\Omega$. The assertion about the
center of $\caM(G)$ follows immediately.
\end{proof}

\begin{cor}
The center $\frZ_\Omega$ of the category $\caM(G)_\Omega$ is a Noetherian $\bbC$-algebra (isomorphic to
the algebra of regular functions of the quotient of a complex algebraic torus by the action of a finite
group).
\end{cor}

\begin{rmqs}
1. It may be useful to give explicitly what appears only implicitly above: the value at a point
$(M,(\rho,W))_G$ of $\Omega(G)$ of the function $z\in \frZ(G)$ is the scalar by which $z$ acts on
$i_P^G (\rho,W)$, where $P$ is a parabolic subgroup of $G$ with Levi factor $M$.

--- 2. The central idempotents $e_\Omega$ of \ref{IdCent} are the characteristic functions of the
connected components $\Omega$ of $\Omega(G)$.

--- 3. Corollary \ref{centreMpiconcl} generalizes: the categories $\caM(G)_\Omega$ are indecomposable.
\end{rmqs}

\subsection{Harish-Chandra homomorphisms}
\index[ter]{Harish-Chandra!(homomorphism)}

The terminology used comes from the following analogy with the theory of real reductive groups \cite{HC}:

- center of $\caM(G)$ $\leftrightarrow$ center of the enveloping algebra,

- cuspidal support of representations $\leftrightarrow$ infinitesimal character.

Let $P=MN$ be a parabolic subgroup of $G$. Recall the map defined in \ref{BIR}, which we now easily see
is a finite morphism of algebraic varieties:
\[
i_{MG}\colon \Omega(M)\rightarrow \Omega(G)
\]
This induces a morphism between algebras of polynomial functions:
\[
i_{MG}^*\colon \frZ(G) \rightarrow \frZ(M)
\]
called the Harish-Chandra homomorphism. Since the morphism $i_{MG}$ is finite, $ \frZ(M)$ is a finitely
generated $ \frZ(G)$-module. The results of \ref{BIR} and the first remark of the previous section show
that:

\begin{prop}
Let $z \in \frZ(G)$ and set $z_M= i_{MG}^*(z)$. Then, for any representation $(\sigma,E)$ of $\caM(M)$,
$i_P^G(z_M)$ is a $G$-equivariant endomorphism of $i_P^G(\sigma,E)$ which coincides with the endomorphism
defined by $z$, viewed as an element of the center of the category. For any representation $(\pi,V)$ of
$\caM(G)$, $r_P^G(z)$ is an $M$-equivariant endomorphism of $r_P^G(\pi,V)$ which coincides with the one
defined by $z_M$.
\end{prop}
\index[not]{i_MG*@$i_{MG}^*$}

\subsection{Bernstein center and compact open subgroups} \label{KetB}

The following result is implicit in the proof of the decomposition theorem \ref{dcompthm}.

\begin{thm}
Let $(\pi,V)$ be a representation in $\caM(G)$, and let
\[
V=\bigoplus_{\frs \in \caB(G)} V_\frs
\]
be its Bernstein decomposition. Let $K$ be a compact open subgroup of $G$. Then $V_\frs^K=0$ except for a
finite number of components $V_\frs$.

More precisely, for such a $K$, let $\mathbf{Irr}(G)_K$ \index[not]{Irr(G)_K@$\mathbf{Irr}(G)_K$} denote
the set of equivalence classes of irreducible smooth representations of $G$ admitting non-zero vectors
fixed by $K$, and $\Omega(G)_K$ \index[not]{ZZOmega(G)_K@$\Omega(G)_K$} the set of their cuspidal
supports. Then $\Omega(G)_K$ is a finite set of connected components of $\Omega(G)$ (and we denote by
$\caB(G)_K$ their inertial support). \index[not]{B(G)_K@$\caB(G)_K$}
\end{thm}

\begin{proof}
First, recall that the number of inertial classes of irreducible supercuspidal representations of $G$
admitting non-zero vectors fixed by $K$ is finite (Remark \ref{DeCcusp}) and that this fact is crucial in
the proof of the decomposition theorem. Now let $\frs=[M,(\rho,W)]_G \in \caB(G)$ and let $D$ be the
inertial class of $(\rho,W)$ in $M$. For all $\psi \in \caX(M)$, we can identify $i_P^G(\rho\psi,W)^K$
with a space not depending on $\psi$. Indeed, by Lemma \ref{Ind},
\[
i_P^G(\rho\psi,W)^K \simeq  \bigoplus_{\bar g\in  P  \backslash G /K} W^{p(K_g)},
\]
where $K_g=P \cap gKg^{-1}$ and $p$ denotes the projection of $P$ onto $P/N \simeq M$. Recall that
$ P\backslash G /K$ is finite (Lemma \ref{IndAdm}). For any cuspidal support
$\theta=(M,(\rho \psi,W))_G$, the fiber $\mathbf{Sc}^{-1}(\theta)\subset \mathbf{Irr}(G)$ is the set of
composition factors of $i_P^G(\rho\psi,W)$, and thus $\theta \in \Omega(G)_K$ if and only if
$\bigoplus_{\bar g \in  P \backslash G /K} W^{p(K_g)} \neq \{0\}$. Thus, either $\frs$ is contained in
$\Omega(G)_K$, or $\frs$ does not intersect $\Omega(G)_K$, i.e., $\Omega(G)_K$ is a union of connected
components of $\Omega(G)$. Moreover $\frs \in \caB(G)_K$ if and only if $D \in \caB(M)_{p(K_g)}$ for a
certain $\bar g \in  P \backslash G /K$. This reduces us to the case of supercuspidal representations
treated in Remark \ref{DeCcusp}.
\end{proof}

\begin{prop}
Let $K$ be a compact open subgroup of $G$, contained and normal in $K_0$ and admitting an Iwahori
decomposition with respect to the standard parabolic subgroups.

$(i)$ The subcategory $\caM(G)_K$ of $\caM(G)$ of representations generated by their $K$-invariant
vectors is stable under passing to subquotients and
\[
\caM(G)_K \simeq \prod_{\frs \in \caB(G)_K} \caM(G)_\frs.
\]
$(ii)$ The functor $V \mapsto V^K$ realizes an equivalence of categories from $\caM(G)_K$ to
$\caM(\caH(G,K))$, whose inverse is
\[
M \mapsto \caH(G)\otimes_{\caH(G,K)} M.
\]
$(iii)$ For all $(\pi,V)\in \caM(G)$, $(\pi,V) \in \caM(G)_K$ if and only if
$(\tilde \pi,\widetilde V)\in\caM(G)_K $.

$(iv)$ The algebra $\caH(G,K)$ is Noetherian.
\end{prop}

\begin{proof}
The first assertion of $(i)$ is established in Proposition \ref{consLJG}. From the above, we deduce the
decomposition of $\caM(G)_K$.

$(ii)$ follows from Theorem \ref{35}.

$(iii)$ Let $(\pi,V) \in \caM(G)_K$. We decompose $\widetilde V$ into $\widetilde V=W_1\oplus W_2$, where
$W_1$ is the subrepresentation generated by $\widetilde V^K$, thus in
\[
\caM(G)_K \simeq \prod_{\frs \in \caB(G)_K} \caM(G)_\frs
\]
and $W_2$ is in $\prod_{\frs \notin \caB(G)_K} \caM(G)_\frs$. We then have $W_2^K=\tilde \pi(e_K)W_2=0$.
We want to show that $W_2$ is zero. Suppose that $w \in W_2$ is non-zero. Then there exists $v\in V$ such
that $w(v)\neq 0$, and since $V$ is generated by $V^K$, there exist $g_1,\ldots, g_l \in G$,
$v_1,\ldots,v_l \in V^K$ such that $v=\sum_i \pi(g_i)\cdot v_i$. We then have
\[
0 \neq  w(v)=\sum_i (\tilde \pi(g_i^{-1})\cdot w)(v_i)
\]
and thus one of the factors $(\tilde \pi(g_i^{-1})\cdot w)(v_i)$ is non-zero. Now
\[
0 \neq (\tilde \pi(g_i^{-1})\cdot w)(v_i)= (\tilde
\pi(g_i^{-1})\cdot w)( \pi(e_K)\cdot v_i)= (\tilde \pi(e_K)\tilde \pi(g_i^{-1})\cdot w)(v_i),
\]
and we obtain a contradiction since $(\tilde \pi(e_K)\tilde \pi(g_i^{-1})\cdot w)\in W_2^K=0$.
Conversely, if $(\tilde \pi,\widetilde V)\in  \caM(G)_K$, $V\hookrightarrow \widetilde{\widetilde V}$ is
in $\caM(G)_K$.

$(iv)$ We have seen that each component $\caM(G)_\frs$ is a Noetherian category, so the same is true for
$\caM(G)_K$, and thus for $\caM(\caH(G,K))$. Now the algebra $\caH(G,K)$ is finitely generated over
itself, hence Noetherian.
\end{proof}

\medskip

The inverse problem, which consists of realizing, for any $\frs \in \caB(G)$, the category $\caM(G)_\frs$
as a category of the form $\caM(\caH)$, where $\caH$ is a Hecke algebra of the form $\caH=\caH(G,K)$ (in
fact slightly more general Hecke algebras are allowed) is called the theory of types. We refer the
interested reader to \cite{BuKu}.

\bigskip

Let $K$ be a compact open subgroup of $G$, contained and normal in $K_0$ and admitting an Iwahori
decomposition with respect to the standard parabolic subgroups and let $P=MN$ be a standard parabolic
subgroup. Then $K_M=K\cap M$ has the same properties as $K$ relative to the group $M$. Let us resume the
notation of Section \ref{BIR}.

\begin{cor}
The subset $\caB(M)_{K_M}$ of $\caB(M)$ is the inverse image under $i_{MG}$ of the subset $\caB(G)_{K}$.
\end{cor}
\begin{proof}
Let $\frs=[M,(\sigma,E)]_G$ be an element of $\caB(G)$. We show that $\frs \in  \caB(G)_K$ if and only if
$\frt=[M,\sigma]_M \in \caB(M)_{K_M}$. Suppose $\frt=[M,\sigma]_M \in \caB(M)_{K_M}$, i.e., that $E$
admits non-zero vectors fixed by $K_M$. Consider the representation $\pi=i_P^G(\sigma)$, with space $V$.
Since by Jacquet's lemma, $V^K$ surjects onto $V_N^{K_M}$ which contains a subquotient isomorphic to
$E^{K_M}$ by the geometric lemma, we have $V^K \neq 0$. Thus there exists an irreducible representation
in $\caM(G)_\frs$ which is in $\caM(G)_K$, which shows from the above that
$\caM(G)_\frs \subset \caM(G)_K$, and thus $\frs \in \caB(G)_K$. Conversely, if $\frs \in \caB(G)_K$, by
the proof of the theorem, $E^{p(K_g)} \neq 0$ for a certain $g\in G$. Write $g=pk$, with $k \in K_0$ and
$p \in P$. Then, since $K$ is assumed normal in $K_0$,
\[
K_g=gKg^{-1}\cap P=  pkKk^{-1}p^{-1}\cap P=p(K\cap P)p^{-1}
\]
and all the $p(K_g)$ are conjugate (to $K_M$) in $M$. Thus $E^{K_M} \neq 0$. Let us now prove the
corollary. Let $\frt \in \caB(M)$; we will now show that $\frt \in \caB(M)_{K_M}$ if and only if
$i_{MG}(\frt) \in \caB(G)_K$. We have $\frt=[L,\sigma]_M$ for a certain standard Levi subgroup $L$
contained in $M$. If $\frt \in \caB(M)_{K_M}$, then $\sigma$ admits a non-zero vector fixed by
$K_L=K\cap L$, and $i_{MG}(\frt) =[L,\sigma]_G$ admits a non-zero vector fixed by $K$, thus
$i_{MG}(\frt) \in \caB(G)_K$. Conversely, if $i_{MG}(\frt) \in \caB(G)_K$, then $\sigma$ admits a
non-zero vector fixed by $K_L=K\cap L$ and thus $\frt \in \caB(M)_{K_M}$.
\end{proof}

\subsection{Finitely generated representations}\label{reptf}

For any representation $(\pi,V)$ of $\caM(G)$, let $V_\frs$ denote the component of $V$ in the category
$\caM(G)_\frs$, $\frs \in \caB(G)$ ({\sl cf.} \ref{VVs}).

\begin{lemme}
Let $(\tau,E)$ be a finitely generated representation in $\caM(G)$. Then $E_\frs$ is zero except for a
finite number of $\frs \in \caB(G)$.
\end{lemme}
\begin{proof}
Each generator of $(\tau,E)$ has components in only a finite number of $E_\frs$.
\end{proof}

\begin{thm}
Let $(\tau,E)$ be a finitely generated smooth representation of $G$. Then $ (\tau,E)$ is
$\frZ(G)$-admissible, i.e., for any compact open subgroup $K$ of $G$, $E^K$ is a finitely generated
$\frZ(G)$-module.
\end{thm}

Note that this is a generalization of the admissibility of irreducible representations: if $(\tau,E)$ is
irreducible, $\frZ(G)$ acts by scalars on $E^K$, and thus if $E^K$ is $\frZ(G)$-admissible, this space is
finite-dimensional.

\begin{proof}
Using the lemma, we may assume that $(\tau,E)$ is in a category $\caM(G)_\frs$, for a certain
$\frs \in \caB(G)$. Since $\frZ(G)$ then acts on $(\tau,E)$ via its quotient $\frZ_\frs$, we must show
that $E^K$ is a finitely generated $\frZ_\frs$-module.

Let us start with the case where $\frs \in \caB(G)$ is the inertial support of an irreducible
supercuspidal representation $(\pi,V)$ of $G$. We then resume the notation of Sections \ref{progenpi} and
\ref{centreMpi}. Let $(\Pi,V_\Pi)$ be a progenerator of the category $\caM(G)_{[\pi]}$, for example the
one constructed in Remark \ref{progenpi} (denoted there by $(\Pi_1,V_{\Pi_1})$). We know that
$(\tau,E) \in\caM(G)_{[\pi]} $ is isomorphic to a quotient of a sum of copies of $(\Pi,V_\Pi)$, more
precisely $(\tau,E) $ fits into an exact sequence of $G$-modules of the form
\[
\bigoplus_{i\in I}V_\Pi \rightarrow  \bigoplus_{j\in J} V_\Pi
\rightarrow E\rightarrow 0
\]
If $(\tau,E)$ is finitely generated, we can then take $J$ to be finite. It therefore suffices to show
that $V_\Pi^K$ is a finitely generated $\caZ_{[\pi]}$-module. Recall that
\[
(\Pi,V_\Pi)=\ind_{{}^0G}^{\; \,   G} (\res_{{}^0G}^{\; \,   G} (\pi,V)) \simeq V\otimes
F,
\]
where $F=\bbC[\Lambda(G)]$ is the algebra of polynomial functions on the variety $\caX(G)$. Since $\pi$
is irreducible supercuspidal, it is admissible. It therefore remains to show that $F$ is a finitely
generated $\caZ_{[\pi]}$-module. But this is obvious, given the description of $\caZ_{[\pi]}$ as the
algebra of polynomial functions on
\[
\mathbf{Irr}(G)_{[\pi]} \simeq \caX(G)/\caX(G)(\pi).
\]

Let us now look at the general case. Let $(M,(\rho,W))$ be a cuspidal datum such that
$\frs=[M,(\rho,W)]_G$, where $M$ is a Levi subgroup of $G$. Let $P$ be a parabolic subgroup of $G$ with
Levi factor $M$. Let $(\Pi_\frs,V_\frs)$ be the progenerator of the category $\caM(G)_\frs$ constructed
in \ref{progenMOm} (denoted there by $(\Pi_\Omega,V_\Omega)$). The same reasoning as above shows that it
suffices to see that $V_\frs^K$ is a finitely generated $\frZ_\frs$-module. Now
\[
(\Pi_\frs,V_\frs)=i_P^G (\Pi,V_\Pi)
\]
where $(\Pi,V_\Pi)$ is the progenerator of the category $\caM(M)_{[\rho]}$ constructed as above. From the
above, $(\Pi,V_\Pi)$ is $\frZ_{[\rho]}$-admissible, where $\frZ_{[\rho]}$ is the center of the category
$\caM(M)_{[\rho]}$. As in Section \ref{GBmod}, it is easy to see that parabolic induction preserves
$\frZ_{[\rho]}$-admissibility. Moreover, $\frZ_{[\rho]}$ is a finitely generated $\frZ_\frs$-module (see
\ref{centreMGOm}), and thus $V_\frs$ is $\frZ_\frs$-admissible.
\end{proof}

\begin{cor}
Let $K$ be a compact open subgroup of $G$. Then the Hecke algebra $\caH(G,K)$ is a finitely generated
module over $\frZ(G)$. In particular, $\caH(G,K)$ is a finitely generated module over its center, which
contains $e_K*\frZ(G)*e_K$.
\end{cor}

\begin{proof}
Let us take as a representation of $G$ the space $\scrD(G/K)$ of locally constant compactly supported
functions on $G/K$ (which can also be viewed as the space of functions on $G$ right-invariant by $K$ and
compactly supported). It is a cyclic representation generated by the characteristic function of $K$. The
group $G$ acts on the left, and $\scrD(G/K)^K\simeq \scrD(G,K)\simeq \caH(G,K)$ as a $\frZ(G)$-module.
\end{proof}

\begin{prop}
For a representation $(\tau,E)$ of $G$ to be finitely generated, it is necessary and sufficient that it
be $\frZ(G)$-admissible and have non-zero components in only a finite number of $\caM(G)_\frs$.
\end{prop}
\begin{proof}
The condition is necessary by the lemma and the theorem above. It is sufficient, because if
\[
E=\bigoplus_{\frs \in F} E_\frs
\]
where $F$ is a finite subset of $\caB(G)$, we then choose a compact open subgroup $K$ of $G$ satisfying
the conditions of Corollary \ref{KetB} such that
\[
F \subset \caB(G)_K
\]
and thus $(\tau,E)$ is generated by $E^K$ where $E^K$ is finitely generated over $\frZ(G)$, and thus a
fortiori finitely generated over $\caH(G,K)$, by the corollary.
\end{proof}

\begin{variante}
We can reformulate and prove the theorem and the proposition above for $(G,B)$-modules: a $(G,B)$-module
$(\tau,E)$ is finitely generated if and only if it is $B\otimes_\bbC \frZ(G)$-admissible and has non-zero
components in only a finite number of $\caM(G)_\frs$. We deduce that the induction functors $i_P^G$
preserve finitely generated $B$-modules.
\end{variante}

\subsection{Return to Howe's theorem}\label{Ho2}

Let $(\pi,V)$ be a smooth representation of $G$. We say that $(\pi,V)$ is $\frZ(G)$-finite if the
annihilator of $(\pi,V)$ in $\frZ(G)$ is an ideal of finite codimension.

We conclude our study of the relations between admissibility, finite generation and the action of
$\frZ(G)$ with the following theorem.

\needspace{169pt}
\begin{thm}
Let $(\pi,V)$ be a smooth representation of $G$ and consider the following possible properties of
$(\pi,V)$:

$(1)$ $(\pi,V)$ is finitely generated,

$(2)$ $(\pi,V)$ is admissible,

$(3)$ $(\pi,V)$ is $\frZ(G)$-finite.

Then any two of these properties imply the third, and the representation $(\pi,V)$ is then of finite
length.
\end{thm}

\begin{proof}
Howe's theorem \ref{ThmHowe} asserts that $(1)$ and $(2)$ are equivalent to $(\pi,V)$ being of finite
length. It is easy to see by induction on the length that a representation of finite length is
$\frZ(G)$-finite, starting with the case where $(\pi,V)$ is irreducible (this is then Schur's lemma).

Suppose that $(\pi,V)$ satisfies $(1)$ and $(3)$. We saw in \ref{reptf} that $(1)$ implies that $(\pi,V)$
is $\frZ(G)$-admissible, i.e., that for any compact open subgroup $K$ of $G$, $V^K$ is a finitely
generated $\frZ(G)$-module. But $(3)$ tells us that $V^K$ is annihilated by an ideal of finite
codimension of $\frZ(G)$, thus $V^K$ is finite-dimensional. This shows that $(\pi,V)$ is admissible.

Suppose that $(\pi,V)$ satisfies $(2)$ and $(3)$. By the characterization of finitely generated
representations of Proposition \ref{reptf}, we must show that $(\pi,V)$ is $\frZ(G)$-admissible, and has
only a finite number of non-trivial Bernstein components $V_\frs$. The first point is clear, since for
any compact open subgroup $K$ of $G$, $V^K$ is finite-dimensional, thus finitely generated over
$\frZ(G)$. The second is as well, since an ideal of finite codimension of $\frZ(G)$ will contain the
components $\frZ_\frs$ except for a finite number of $\frs \in \caB(G)$, and on the other hand, the
identity of $\frZ_ \frs$ acts by the identity on $V_\frs$.
\end{proof}

\section{Notes on Chapter VI}

Generalities on induction and restriction functors can be found in many places in the literature. The
same is true for the main results on supercuspidal representations (\cite{BeZe1}, \cite{BeZe2},
\cite{Ca}). The attribution of results is sometimes delicate. Thus, the various characterizations of
supercuspidal representations (Theorem \ref{supercusp}) seem to be due to H. Jacquet for
$G=\mathrm{GL}(n,\bbF)$ and to Harish-Chandra in general, and the name Harish-Chandra's theorem seems to
have prevailed (perhaps also because other important results bear Jacquet's name). I also used \cite{DeB}
for the writing of the proof of this theorem. The admissibility and uniform admissibility theorems, as
well as Proposition \ref{finitude} are due to Bernstein (\cite{Be74}, \cite{BeZe1}). Here again, the
notes of S. DeBacker \cite{DeB} were useful to me. The study of the category $\caM(G)_{[\pi]}$ and the
determination of its center is based on the notes \cite{Ro}. The geometric lemma is due independently to
Bernstein-Zelevinskii \cite{BeZe2} and Casselman \cite{Ca}; the proof largely follows \cite{BeZe2}, but
also \cite{Ca} at a crucial point, which allows for an explicit integral formula (\ref{IndRes}). The
consequences of the geometric lemma on the composition series of parabolic inductions of supercuspidal
representations \ref{geolemmacusp} and \ref{compindsup} are obtained in \cite{BeZe2}. The generic
irreducibility of induced representations \ref{irrindgen} is taken from \cite{Be1}. The proof of
Jacquet's lemma \ref{LemmeJacquet} written here borrows from several sources, in particular \cite{BeZe1},
\cite{Be1}, \cite{DeB}. Howe's theorem is taken from \cite{BeZe1}. The proof of the Bernstein
decomposition theorem essentially follows the notes \cite{Ro}, notes themselves inspired by \cite{Be1} to
which I referred from time to time. The second adjunction theorem is established, for the main lines of
the proof, in \cite{Be1} (see also \cite{Be2} and \cite{Dat}). I tried to make it more accessible by
adding some details. It nevertheless remains a difficult and subtle proof, the crucial point of which is
the stabilization lemma \ref{stable} used in the proof of the generalized Jacquet lemma
\ref{LemmeJacquet2}. The second adjunction theorem is equivalent to a strong form of Casselman's duality,
valid without the admissibility hypothesis (\cite{Be1} and \cite{Be2}).

The results on the Bernstein center are presented here following the ideas of \cite{Be1} and \cite{Be2},
taken up and detailed in the notes \cite{Ro}. The approach, as we have explained, is somewhat different
from the one followed in the earlier draft \cite{Del}, but the important corollaries of the main result
are taken from \cite{Del}. The statement of Theorem \ref{Ho2} is inspired by the analogous result for
real groups, which can be found scattered in \cite{KV}.





\chapter{Langlands classification}

The goal of this chapter is to state and prove the Langlands classification theorem, which reduces the
problem of determining $\mathbf{Irr}(G)$, i.e., the isomorphism classes of irreducible smooth
representations of a $p$-adic reductive group $G$, to the problem of determining the classes of
irreducible {\sl tempered} representations. The definition of a tempered representation given here relies
solely on normalized exponents, defined in the first section. Casselman's criterion provides a necessary
and sufficient condition on the normalized exponents of a representation for it to be square-integrable
modulo the center: they must lie in certain open cones. A tempered representation is then defined as one
whose normalized exponents all lie in the closures of these open cones. A tempered representation is
therefore close to being a representation square-integrable modulo the center, and moreover, we show that
an irreducible tempered representation appears as a subrepresentation of a parabolic induction of a
representation square-integrable modulo the center. In particular, since it embeds into the parabolic
induction of a unitary representation, it is unitary. We show on the other hand that parabolic induction
preserves the class of tempered representations. This is not the case for parabolic restriction, but one
can define a notion of tempered parabolic restriction, which preserves the class of tempered
representations, and even obtain a tempered version of Frobenius reciprocity and the geometric lemma.

We then study the intertwining operators between representations of the form $i_P^G(\rho,W)$ and
$i_Q^G(\rho,W)$, where $P$ and $Q$ are two parabolic subgroups with the same Levi factor $M$ and
$(\rho,W)$ is a smooth representation of $M$. When $(\rho,W)$ satisfies a certain criterion
($PQ$-regularity), we exhibit a canonical intertwining operator between these representations. The
transitivity properties of these intertwining operators are studied in detail. A Langlands triplet
consists of a parabolic subgroup $P=MN$ of $G$, an irreducible tempered representation $(\sigma,E)$ of
$M$ and an unramified character $\psi$ of $M$ satisfying a certain positivity condition. We then show
that the representation $\sigma \otimes \psi$ is $P\bar P$-regular, and thus that there exists an
intertwining operator between $i_P^G(\sigma\otimes\psi)$ and $i_{\bar P}^G(\sigma\otimes\psi)$. A more
refined study shows that the space of such operators is $1$-dimensional, and that
$i_P^G(\sigma\otimes\psi)$ (resp. $i_{\bar P}^G(\sigma\otimes\psi)$) admits a unique irreducible quotient
(resp. subrepresentation). We call this quotient the Langlands quotient of the triplet. The
classification theorem asserts that any irreducible representation of $G$ appears as a Langlands
quotient, in an essentially unique way.

\vfill

\section{Casselman's criterion and applications}

\subsection{Exponents}\label{exposants}
Let $P=MN$ be a parabolic subgroup of $G$ with split component $A$. Let $(\pi,V)$ be a smooth admissible
representation of $G$. The representation $r_P^G (\pi,V)$ is then admissible by Corollary
\ref{LemmeJacquet}. For any compact open subgroup $K_M$ of $M$, and for any $t \in A$, since
$tK_Mt^{-1}=K_M$, we can write
\[
(r_P^G V)^{K_M}= \bigoplus_{\chi} (r_P^G V)^{K_M}_\chi
\]
where $\chi$ runs through the set of smooth characters of $A$ and $(r_P^G V)^{K_M}_\chi$ is the
characteristic subspace of $ (r_P^G V)^{K_M}$ for the character $\chi$, i.e.,
\[
(r_P^G V)^{K_M}_\chi= \{v \in (r_P^G V)^{K_M} \, | \, \exists d \in \bbN, \,
\forall t\in A,
(r_P^G(\pi)(t)-\chi(t) \Id_V)^d\cdot v=0 \}.
\]

If $K'_M$ is another compact open subgroup of $M$, $K'_M \subset K_M$, then
$(r_P^G V)_\chi^{K_M}\subset (r_P^G V)_\chi^{K'_M}$. We can therefore set
\[
(r_P^G V)_\chi= \bigcup_{K_M}(r_P^G V)_\chi^{K_M}
\]
where $K_M$ runs through the compact open subgroups of $M$. If this space is non-trivial, we say that
$\chi$ is a normalized exponent \index[ter]{normalized exponent} of $\pi$ for the parabolic subgroup
$P=MN$. We denote by $\mathrm{Exp}(A,r_P^GV)$ the set of these exponents.
\index[not]{Exp(A,rPG)@$\mathrm{Exp}(A,r_P^GV)$} We then obtain
\begin{equation}\label{secar}
(r_P^G V) = \bigoplus_{\chi \in \mathrm{Exp}(A,r_P^GV) } (r_P^G V)_\chi.
\end{equation}
If $(\pi,V)$ is of finite length, the same is true for $r_P^G(\pi,V)$ ({\sl cf.} \ref{rPGlongfin}), and
the set $\mathrm{Exp}(A,r_P^GV)$ is finite, by Schur's lemma.

\begin{rmq}
If we take $P=G$, we obtain a decomposition:
\[
V = \bigoplus_{\chi}  V_\chi
\]
which is a slightly weaker version of the decomposition of Proposition \ref{caraccentral} (since
$A_G\subset Z(G)$).
\end{rmq}

Recall that parabolic induction and restriction preserve admissibility. We will compare the
decompositions of a representation and that of its induced or restricted representation.
\needspace{284pt}
\begin{lemme}
$(i)$ Let $P=MN$ be a parabolic subgroup of $G$ and $(\rho,W)$ an admissible representation of $M$. Let
\[
W=\bigoplus_{\psi \in  \mathrm{Exp}(A_M,W)}  W_\psi, \qquad i_P^G W=
\bigoplus_{\chi \in \mathrm{Exp}(A_G,i_P^G W)}( i_P^G W)_\chi
\]
be the respective decompositions into characteristic subspaces of $W$ and of $ i_P^G W$ for the actions
of $A_M$ and of $A_G$ respectively. If $$( i_P^G W)_\chi\neq \{0\},$$ then $\chi$ is the restriction to
$A_G$ of a character $\psi$ of $A_M$ such that $ W_\psi\neq \{0\}$.

--- $(ii)$ Let $P=MN$ be a parabolic subgroup of $G$ and $(\pi,V)$ an admissible representation of $G$.
Let
\[
V=\bigoplus_{\psi \in \mathrm{Exp}(A_G,V) }  V_\psi, \quad r_P^G V=
\bigoplus_{\chi\in \mathrm{Exp}(A_M,r_P^G V) } ( r_P^G V)_\chi
\]
be the respective decompositions into characteristic subspaces of $V$ and of $ r_P^G V$ for the actions
of $A_G$ and of $A_M$ respectively. Then the restriction of a character $\chi$ of $A_M$ appearing in the
second decomposition to $A_G$ is a character $\psi$ appearing in the first decomposition.
\end{lemme}

\begin{proof}
$(i)$ For all $f \in i_P^G W$, let us simply denote by $g\cdot f$ the action (by right translation) of an
element $g$ of $G$ on $f$. Suppose that $f$ is a non-zero eigenvector for the action of all $a \in A_G$
for the eigenvalue $\chi(a)$. We then have, for all $x \in G$,
\[
\chi(a)f(x)=(a\cdot f)(x)=f(xa)=f(ax)=\rho(a)\cdot f(x).
\]
Take $x\in G$ such that $w=f(x)\neq 0$. Then $w$ is an eigenvector of $\rho(a)$ for the eigenvalue
$\chi(a)$. This shows that $\chi$ is the restriction to $A_G$ of a character $\psi$ of $A_M$ intervening
non-trivially in the decomposition of $W$.

$(ii)$ This is obvious from the definition of $r_P^G V$, once we note that for all $a \in A_G$,
$\delta_P(a)=1$.
\end{proof}

\subsection{Casselman's criterion}\label{CasCrit}
\index[ter]{Casselman's criterion} In this section, we establish a criterion showing that a
representation is square-integrable modulo the center.
\begin{thm}
Let $(\pi,V)$ be an admissible smooth representation of $G$ admitting a unitary central character. Then
$\pi$ is square-integrable modulo the center if and only if, for every parabolic subgroup $P=MN$ of $G$,
and for every normalized exponent $\chi$ of $\pi$ for the parabolic subgroup $P=MN$,
$\Re (\chi) \in {}^+[\fra_M^*]_P^G$.
\end{thm}

\begin{proof}
We refer the reader to Section \ref{carreint} for the definitions and notation concerning
square-integrable representations modulo the center and to \ref{APG} and \ref{ReChi} for the definitions
of $\Re (\chi)$ and ${}^+[\fra_M^*]_P^G$. We will slightly adapt the definitions of square-integrable
representations to the context of $p$-adic reductive groups. Indeed, the quotient $Z(G)/A_G$ is compact,
and since $A_G=C_{A_G}{}^0A_G$, with ${}^0A_G$ compact, even $Z(G)/C_{A_G}$ is compact (see
\ref{castordep}). We therefore fix a unitary character $\chi$ of $C_{A_G}$ and a $G$-invariant measure
$dg^*$ on $G/C_{A_G}$, and we define the space $L^2(G,\chi,dg^*)$ as in \ref{L2chi} and the
square-integrable representations modulo the center as in \ref{unit2int}, $C_{A_G}$ replacing $Z(G)$.
Since $Z(G)/C_{A_G}$ is compact, it is easy to see that this new definition is equivalent to the old one.

Let $(\pi,V)$ be as in the statement of the theorem. We seek necessary and sufficient conditions for
\begin{equation} \label{I1}
\int_{G/C_{A_G}} |\phi_{v,\lambda}(g)|^2 \; dg^*
\end{equation}
to converge for all $v \in  V$, $\lambda \in \widetilde V$. Let $K\subset K_0$ be a compact open subgroup
of $G$, normal in $K_0$, admitting an Iwahori decomposition with respect to the standard parabolic
subgroups of $G$ and fixing $v$ and $\lambda$. Let $k_1,\ldots ,k_r$ be a system of representatives of
the cosets of $K$ in $K_0$. The Cartan decomposition \ref{Cartan} of $G$ gives us:
\[
G=  \bigcup_{i,j=1}^r  Kk_i F_\emptyset^+ C_\emptyset^+ k_j K.
\]
Let us now also fix a system of representatives $S$ of the equivalence classes in $C_\emptyset^+$ for the
relation (\ref{eqC}), so that we can make the Cartan decomposition even finer,
\[
G=  \bigcup_{i,j=1}^r  \bigcup_{z\in S}    Kk_i F_\emptyset^+ z
C_{A_G}  k_j K,
\]
and deduce that
\[
G/C_{A_G}=  \bigcup_{i,j=1}^r  \bigcup_{z\in S}    Kk_i F_\emptyset^+
z  k_j K.
\]

These unions are not disjoint, but a given double coset under $K$ appears in them at most
$r^2|F_\emptyset^+|$ times: for fixed $i,j,f$, $z$ is determined by the Cartan decomposition. Up to
constants, this allows us to bound (\ref{I1}) from above and below by
\begin{align*}
& \sum_{i,j=1}^r   \sum_{f \in F_\emptyset^+} \sum_{z \in S}
|\phi_{v,\lambda}(k_i f z k_j)|^2  \mu_G(KfzK)\\
&= \sum_{i,j=1}^r   \sum_{f \in F_\emptyset^+} \sum_{z\in S}
|\phi_{k_j\cdot v,k_i^{-1}\cdot \lambda}(f z)|^2  \mu_G(KfzK).
\end{align*}

We thus obtain a finite sum of series with non-negative terms, all of the form (replacing $k_j\cdot v$
and $k_i^{-1}\cdot \lambda$ simply by $v$ and $\lambda$, and noting that $k_j\cdot v$ and
$k_i^{-1}\cdot \lambda$ are still fixed by $K$)
\begin{equation}\label{S1}
   \sum_{f \in F_\emptyset^+} \sum_{z\in S}
|\phi_{v,\lambda}(f z)|^2  \mu_G(KfzK).
\end{equation}
Let us look more closely at the general term of this series. Let us take again the elements $t_\alpha$,
$\alpha \in\Delta_\emptyset$, the submonoid $S_0\simeq\bbN^d$ of $S$ that they generate and the finite
set $R$ of \ref{SSM}, so that $S=\bigcup_{r\in R}rS_0$; we have $|\alpha(t_\alpha)|_\bbF<1$ and
$|\beta(t_\alpha)|_\bbF=1$ if $\beta \neq  \alpha$.

Fix a positive integer $n$, and, for every standard parabolic subgroup $P=MN$, set
\[
S_M(n)= \{ z=\prod_{\alpha \in \Delta_\emptyset} t_\alpha^{i_\alpha}\mid \forall \alpha \in \Delta_\emptyset \setminus
\Delta_\emptyset^M, \,    i_\alpha>n, \;   \forall \alpha \in  \Delta_\emptyset^M, \,  0\leq i_\alpha\leq n\}.
\]
It is then clear that
\[
S_0=\coprod_{P=MN \text{ standard }} S_M(n),\qquad S=\bigcup_{r\in R}\ \bigcup_{P=MN \text{ standard }} rS_M(n).
\]
Note that every element $z$ of $rS_M(n)$, $r\in R$, can be written $z=z_1t'$ with $z_1$ in the finite set
$S_M(n)_1$ of the $r\prod_{\alpha \in \Delta_\emptyset^M}  t_\alpha^{i_\alpha}$, $r\in R$,
$0\leq i_\alpha\leq n$, and $t'$ in the set
\[
S_M'(n)=   \{z'= \prod_{\alpha \in \Delta_\emptyset\setminus  \Delta_\emptyset^M } t_\alpha^{i_\alpha},
\, i_\alpha >n\}.
\]

Now fix a standard parabolic subgroup $P=MN$, and decompose $K$ with respect to the parabolic subgroup
$P=MN$, with the usual notation:
\[
K=K_{N}K_MK_{\overline{N}}.
\]
Let $t'\in S'_M(n)$. For $\beta\in\Delta_\emptyset^M$, we have $\beta(t')\in\varpi^\bbZ$ and
$|\beta(t')|_\bbF=1$, hence $\beta(t')=1$; thus $t'$, viewed as an element of $X_*(A_\emptyset)$, lies in
$X_*(A_\emptyset)\cap\fra_M=X_*(A_M)$, that is, $t'\in C_{A_M}^+$ (see (\ref{coneA})). In particular,
$t'$ is central in $M$, and by \ref{KN}, we have
\[
t'^{-1}K_{\overline{N}}t' \subset K_{\overline{N}}\subset K,\quad   t'^{-1}K_{M}t' =
K_{M} \subset K,
\]
and
\[
t'^{-1}K_{N}t'=  K_{t'^{-1}Nt'}
\supset K_N.
\]
Hence
\begin{align*}
 Kt'K&= K_{ N}K_MK_{\overline{N}} t'K= t' (t'^{-1}K_Nt') ( t'^{-1}K_Mt')(t'^{-1}K_{\overline{N}}t') K  \\
&= t' K_{t'^{-1}Nt'}K .
\end{align*}
This implies that
\begin{small}
\[
\mu_G(Kt'K)=\mu_G( t' K_{t'^{-1}Nt'}K)=\mu_G( K_{t'^{-1}Nt'}K)= \mu_G(
K_{t'^{-1}Nt'}K_MK_{\overline{N}})
\]
\[
= [K_{t'^{-1}Nt'}: K_N] \mu_G(K_NK_MK_{\overline{N}})=[K_{t'^{-1}Nt'}: K_N] \mu_G(K) .
\]
\end{small}

The computations of modular functions carried out in \ref{calcfonctmod} for the elements of $A_M$ give
$[K_{t'^{-1}Nt'}: K_N]=\delta_P(t')$, so that we finally have
\[
\mu_G(Kt'K)=\delta_P(t') \mu_G(K).
\]
Now let $f\in F_\emptyset^+$, $r\in R$ and $z\in rS_M(n)$. Write $z=z_1t'$ as above, with
$z_1\in S_M(n)_1$ and $t'\in S'_M(n)$, and set $m=fz_1$. The element $m$ runs over the finite set
$F_\emptyset^+S_M(n)_1\subset M_\emptyset$, and it commutes with $t'$. For $x,y\in G$,
$KxyK\subset KxK\,yK$ is a union of at most $\mu_G(KxK)/\mu_G(K)$ left translates of $KyK$, hence
$\mu_G(KxyK)\le \mu_G(KxK)\mu_G(KyK)/\mu_G(K)$. Applying this to $(x,y)=(m,t')$ and to
$(x,y)=(m^{-1},mt')$, we obtain constants $c,C>0$, depending only on the finite set of the $m$, such that
\[
c\,\delta_P(t')\mu_G(K)\le \mu_G(KfzK)=\mu_G(Kmt'K) \le C\,\delta_P(t')\mu_G(K).
\]
Moreover, $\phi_{v,\lambda}(fz)=\phi_{v,\lambda}(t'm)=\phi_{m\cdot v,\lambda}(t')$. Every element of $S$
belongs to one of the sets $rS_M(n)$, and to at most $|R|$ of them, and the map $(f,z)\mapsto(m,t')$ has
fibers of bounded cardinality. Since the number of standard parabolic subgroups of $G$ is finite, as is
$F_\emptyset^+S_M(n)_1$, the problem of the convergence of the series (\ref{S1}), for all $\lambda,v$,
therefore reduces to the problem of the convergence of the series of the form:
\[
\sum_{z\in S_M'(n)} |\phi_{v,\lambda}(z)|^2   \delta_P(z)    .
\]

Since $(\pi,V)$ is admissible, we can choose $\epsilon>0$ (depending on $v$ and $\lambda$) such that
(\ref{casselmanpairingadm}) holds for all $z \in C_{A}^+(\epsilon)$. We then choose $n$ large enough, so
that $ S_M'(n) \subset  A_M^+(\epsilon)$. We then have
\begin{align*}
 &\sum_{z\in S_M'(n)} |\phi_{v,\lambda}(z)|^2
 \delta_P(z)    \\
&= \sum_{z\in S_M '(n)}
|\delta_P(z)^{1/2}\lambda(\pi(z)\cdot v)  |^2 \\
 &=     \sum_{z\in S_M'(n)}
| \langle  r_P^G(\pi)(z)\cdot j_N(v), j_{\overline{N}}(\lambda) \rangle_P  |^2,
\end{align*}
and we are reduced to the study of the series
\[
\sum_{z\in S_M'(n)}| \langle  r_P^G(\pi)(z)\cdot j_N(v), j_{\overline{N}}(\lambda) \rangle_P  |^2.
\]

Identify the subgroup $S_M$ generated by $S_M'(n)$ in $A_M$ with $\bbZ^d$, where
$d=|\Delta_\emptyset\setminus \Delta_\emptyset^M   |$. The normalized exponents of $(\pi,V)$ for the
parabolic subgroup $P$, restricted to this subgroup, therefore identify with characters of $\bbZ^d$. Seen
in this way, such a normalized exponent $\chi$ is given by a $d$-tuple $(z_1,\ldots,z_d)$ of non-zero
complex numbers (see \ref{repsZd}). More explicitly, a basis of $S_M$ is given by the $t_{\alpha_i}$,
$i=1,\ldots, d$, where
\[
\{ \alpha_i\}_{i=1,\ldots, d}= \Delta_\emptyset\setminus \Delta_\emptyset^M,
\]
and $z_i=\chi(t_{\alpha_i})$.

Let $K\subset K_0$ be a compact open subgroup of $G$, normal in $K_0$, admitting an Iwahori decomposition
with respect to the standard parabolic subgroups of $G$ and fixing $v$ and $\lambda$. We then have
$j_N(v)\in V_N^{K_M}$ and $j_{\overline{N}}(\lambda) \in \widetilde V_{\overline{N}}^{K_M}$. Proposition
\ref{repsZd} gives us the existence of polynomials $Q_\chi$ (viewed as functions on $S_M$) such that
\[
\langle  r_P^G(\pi)(z)\cdot j_N(v), j_{\overline{N}}(\lambda) \rangle_P=
\sum_\chi  Q_\chi(z)\chi(z),
\]
the (finite) sum running over the normalized exponents of $(\pi,V)$ for the parabolic subgroup $P$.

It is then easy to see that the series
\[
\sum_{z\in S_M'(n)} | \langle  r_P^G(\pi)(z)\cdot j_N(v), j_{\bar
N}(\lambda) \rangle_P   |^2 =
\sum_{z\in S_M'(n)}|\sum_\chi  Q_\chi(z)\chi(z)|^2
\]
converges if for every normalized exponent $\chi$, $|\chi(t_{\alpha_i})|<1$, for all $i=1,\ldots, d$.
Conversely, if one of the normalized exponents does not satisfy the condition, say
$\chi \in \mathrm{Exp}(A,r_P^GV)$, for some standard parabolic subgroup $P$, then, taking $v$ such that
$j_N(v)$ lies in the corresponding eigenspace, and $\lambda \in \widetilde V$ such that
$\langle  j_N(v), j_{\overline{N}}(\lambda)  \rangle_P\neq 0$, we obtain
\[
\sum_{z\in S_M'(n)}
| \langle  r_P^G(\pi)(z)\cdot j_N(v), j_{\overline{N}}(\lambda) \rangle_P  |^2 =
\sum_{z\in S_M'(n)} c |\chi(z)|^2=+\infty,
\]
$c$ being a positive constant. Let us translate this condition. For $\beta\in\Delta_\emptyset$, we have
$|\beta(t_{\alpha_i})|_\bbF=1$ if $\beta\neq\alpha_i$ and $|\alpha_i(t_{\alpha_i})|_\bbF<1$, that is,
$\langle\beta,H_M(t_{\alpha_i})\rangle=0$ and $\langle\alpha_i,H_M(t_{\alpha_i})\rangle>0$:
$H_M(t_{\alpha_i})\in\fra_M$ is, modulo $\fra_G$, a positive multiple of the fundamental coweight
associated with $\alpha_i$. Write $\Re(\chi)=\sum_{j=1}^d c_j\,{\alpha_j}_{|\fra_M}+\mu_G$, with
$\mu_G\in\fra_G^*$. Since the central character of $\pi$ is unitary, and $\chi_{|A_G}$ is this central
character restricted to $A_G$, we have $\mu_G=0$. Since
$|\chi(t_{\alpha_i})|=q^{-\langle\Re(\chi),H_M(t_{\alpha_i})\rangle}$ and
$\langle\Re(\chi),H_M(t_{\alpha_i})\rangle=c_i\langle\alpha_i,H_M(t_{\alpha_i})\rangle$, the condition
$|\chi(t_{\alpha_i})|<1$ for $i=1,\ldots,d$ is equivalent to $c_i>0$ for all $i$, that is, by \ref{bof},
to $\Re(\chi)\in{}^+[\fra_M^*]_P^G$. This is the condition of the statement. This completes the proof of
the theorem.
\end{proof}

\subsection{An application}\label{corang1}
Recall that a Levi subgroup $M$ of $G$ is said to be maximal if there exists a maximal proper parabolic
subgroup $P$ of $G$ of which $M$ is a Levi factor. Suppose that $M$ is a maximal standard Levi subgroup
of $G$ (this condition is equivalent to $l(M)=2$, {\sl cf.} \ref{levmax}). The set
\index[not]{W(M,*)@$ \caW(M,*) $}
\[
\caW(M,*):=\{ g \in G \mid gMg^{-1} \text{ standard Levi subgroup of } G   \}/M
\]
then has two elements, say $\caW(M,*)=\{ 1, w\}$. Note that with the notation of \ref{WoM}, we have
$\caW(M,*)\simeq W(M,*)/W_M$.

If we set $M'=w\cdot M$, and if $P'$ is the standard parabolic subgroup of $G$ admitting $M'$ as a Levi
factor, then $w\cdot P=\overline{P'}$, the parabolic subgroup opposite to $P'$. We provide necessary
conditions for the reducibility of supercuspidal induced representations from $P$ to $G$.

\begin{thm}
Let $P=MN$ be as above and let $(\rho,W)$ be an irreducible supercuspidal representation of $M$. Set
$\rho'= \rho^w \in \caM(M')_{sc}$. Suppose that $(\pi,V)=i_P^G (\rho,W)$ is reducible. Then $M=M'$ and
$\rho'=\rho \psi$ for a certain unramified character $\psi$ of $M$ with values in $\bbR_+^*$.
\end{thm}

\begin{proof}
We can assume that $\pi$ has a unitary central character. Indeed, if this is not the case, an argument
from the proof of Theorem \ref{irrindgen} tells us that there exists an unramified character $\chi$ of
$G$ such that the central character of $\pi_\chi=\pi\otimes \chi$ is unitary. By Mackey's isomorphism
\ref{Mackey}, we have $\pi_\chi=i_P^G(\rho)\otimes \chi\simeq i_P^G(\rho\otimes \chi_{|M})$. The
representation $\pi$ is therefore reducible if and only if $\pi_\chi$ is, and the conclusion about $\rho$
follows from that about $\rho\otimes \chi_{|M}$.

We first assume that $M \neq M'$. Set $r(\pi)=\{ r_Q^G(\pi) \}$, where $Q$ runs through the proper
standard parabolic subgroups of $G$. By Proposition \ref{geolemmacusp}, the only contributions come from
the parabolic subgroups $P$ and $P'$, more precisely
\[
r_P^G(\pi)=r_P^G i_P^G(\rho)=\rho, \quad  r_{P'}^G(\pi)=r_{P'}^G i_P^G(\rho)=\rho'.
\]
Theorem \ref{compindsup} asserts that $\pi$ is of length at most 2, thus exactly 2, since $\pi$ is
assumed reducible. Set $(\pi',V')=i_{P'}^G (\rho',W')$. A similar analysis applies to $\pi'$. By
adjunction we have
\begin{equation}\label{aDJon}
 \Hom_G(\pi,\pi') \simeq \Hom_{M'}(r_{P'}^G
  i_P^G (\rho), \rho')\simeq \Hom_{M'}(\rho',\rho'),
\end{equation}
and the dimension of this space is 1 because $\rho'$ is irreducible. Similarly $\Hom_G(\pi',\pi)$ is of
dimension 1. We have seen that $l(\pi)=2$, so there exist irreducible representations $\pi_1$ and $\pi_2$
of $G$ and an exact sequence
\begin{equation}\label{pi12}
 0\rightarrow \pi_1 \rightarrow \pi \rightarrow \pi_2 \rightarrow
0.
\end{equation}
The composition factors of $\pi'$ and of $\pi$ are the same (Theorem \ref{compindsup}), but since
$\Hom_G(\pi',\pi)$ is of dimension 1, $\pi$ and $\pi'$ are not completely reducible.

Since the functor $r_P^G$ is exact, we obtain an exact sequence
\[
0\rightarrow r_P^G(\pi_1) \rightarrow r_P^G(\pi)=\rho  \rightarrow r_P^G(\pi_2) \rightarrow 0.
\]
Since $\rho$ is irreducible, either $r_P^G(\pi_1)=\rho$ and $r_P^G(\pi_2)=0$, or vice versa. But since
$\Hom_G(\pi_1,\pi)\neq 0$, we have by adjunction $$\Hom_{M}(r_P^G  \pi_1,\rho)\neq 0$$ and thus
$r_P^G(\pi_1)\neq 0$, whence $r_P^G(\pi_1)=\rho$ and $r_P^G(\pi_2)=0$.

A similar argument using the functor $r_{P'}^G$ shows that either $r_{P'}^G(\pi_1)=\rho'$ and
$r_{P'}^G(\pi_2)=0$, or vice versa. Since the functor $R$ of the second lemma of \ref{dcompthm} is
faithful, we cannot simultaneously have $ r_P^G(\pi_1)=r_{P'}^G(\pi_1)=0$ and similarly for $\pi_2$.
Thus, $r_P^G(\pi_1)=\rho$, $r_P^G(\pi_2)=0$, $r_{P'}^G(\pi_1)=0$ and $r_{P'}^G(\pi_2)=\rho'$.

Since $M$ is maximal, $A_M/A_G$ is of dimension $1$. The restriction to $A_M$ of the central character
$\chi$ of $\rho$ therefore satisfies one of these three mutually exclusive possibilities:

$\bullet$ $|\chi(a)|=1$ for all $a \in A_M$,

$\bullet$ $|\chi(a)|>1$ for all $a \in A_M^{++}$,

$\bullet$ $|\chi(a)|<1$ for all $a \in A_M^{++}$.

In the first case, the central character of $\rho$ is unitary, and since $\rho$ is supercuspidal, $\rho$
is then unitary. In this case, $\pi=i_P^G(\rho)$ is again unitary, thus completely reducible, which is
excluded.

In the second case, replacing $\rho$ by its contragredient reduces the problem to the third case, treated
below. The conclusion of the theorem for $\tilde \rho$ obviously implies the result for $\rho$.

In the third case, Casselman's criterion applies to the representation $\pi_1$, since $r_Q^G\pi_1=0$ for
any proper standard parabolic subgroup $Q \neq P$, and $r_P^G\pi_1=\rho$ and moreover, we ensured that
the central character of $\pi$, thus that of $\pi_1$, is unitary. The representation $\pi_1$ is
square-integrable modulo the center, thus in particular unitary. If we denote by $\bar \pi_1$ the
conjugate representation of $\pi_1$, we then have $\bar \pi_1=\tilde \pi_1$, whence, using (\ref{2ad2}):
\[
\overline{r_{\bar P}^G (\pi_1)}=  r_{\bar P}^G (\bar \pi_1)= r_{\bar P}^G (\tilde \pi_1) = r_{P}^G ( \pi_1) ^\sim=\tilde \rho.
\]
We saw above that $r_{P'}^G\pi_1=0$. Now $r_{\bar P}^G$ is obtained from $r_{P'}^G$ by composition with
$w$. We deduce that $r_{\bar P}^G\pi_1=0$, and thus that $\rho=0$, which is absurd.

We have thus shown that the hypothesis $M\neq M'$ always leads to a contradiction. We therefore have
$M=M'$, and in this case, $P=P'$, $w\cdot P=\bar P$, and $r_P^G i_P^G \rho$ has composition factors
$\rho$ and $\rho'$ (this is still Proposition \ref{geolemmacusp}). We assume $\rho \neq \rho'$, as
otherwise the conclusion of the theorem is trivial. Since $\rho$ and $\rho'$ are supercuspidal, by Lemma
\ref{pisc}, we even have
\[
r_P^G i_P^G \rho= \rho\oplus \rho',
\]
and thus
\[
\Hom_G(\pi,\pi')\simeq\Hom_{M'}(r_{P}^G
i_P^G (\rho), \rho')\simeq\Hom_{M'}(\rho \oplus \rho',\rho'),
\]
which is of dimension 1. We repeat the previous reasoning. The representations $\pi$ and $\pi'$ are of
length 2 and we write the exact sequence as in (\ref{pi12}). The exactness of the Jacquet functor gives
us as before $r_P^G \pi_1=\rho$, $r_P^G \pi_2=\rho'$.

As above, we reduce to the case where Casselman's criterion implies that $\pi_1$ is square-integrable
modulo the center. We thus have $\pi_1^h=\pi_1$, or equivalently $\tilde \pi_1=\bar \pi_1$.

Since $w\cdot P=\bar P$,
\[
{}^w(r_{\bar P}^G(\pi_1))= r_{P}^G(\pi_1)=\rho,
\]
and we obtain, using (\ref{2ad2}):
\[
\rho'=\rho^w= r_{\bar P}^G(\pi_1)= r_{\bar P}^G(\tilde{\tilde
{\pi_1}})=
(r_{P}^G(\tilde  \pi_1))^\sim=  (r_{P}^G(\bar
\pi_1))^\sim=    \   \tilde{\bar \rho}=\rho^h.
\]
In the same way as in the proof of Theorem \ref{irrindgen}, we can write $\rho$ in the form
$\rho=\psi \rho_0$ where $\rho_0$ is unitary supercuspidal, and $\psi$ is an unramified character of $M$
with values in $\bbR_+^\times$. We thus have
\[
\rho'=\rho^h=(\psi \rho_0)^h=\psi^{-1}\rho_0,
\]
whence $\rho'=\psi^{-2}\rho$.
\end{proof}

\begin{cor}
Let $P=MN$, with $M$ a maximal Levi subgroup of $G$, and let $w$ be the non-trivial element of
$\caW(M,*)$. Let $D$ be an inertial class of supercuspidal representations of $M$ such that
$w\cdot D \neq D$ (this is the case for example if $w\cdot M \neq M$). Let $\Omega$ be the connected
component of $\Omega(G)$ associated to $(M,D)$. Then

$(i)$ $i_P^G$ maps the representations of $D$ (i.e., the irreducible representations of $\caM(M)_D$) to
the irreducible representations of $\caM(G)_\Omega$.

$(ii)$ Let $r_D$ be the composition of the functor $r_P^G$ and the projection onto $\caM(M)_D$. Then
$r_D$ is a faithful functor from $\caM(G)_\Omega$ to $\caM(M)_D$.

$(iii)$ For any representation $(\pi,V)$ in $\caM(G)_\Omega$, the natural morphism
$\pi \rightarrow i_P^Gr_D(\pi)$ is an isomorphism. The functor
\[
r_D \colon \caM(G)_\Omega \rightarrow \caM(M)_D
\]
is an equivalence of categories whose inverse is $i_P^G$.
\end{cor}
\begin{proof}
The first point follows immediately from the theorem. Set $w\cdot M=M'$, $w\cdot D =D'$ and let $P'$ be
the standard parabolic subgroup with Levi factor $M'$. We have under these conditions $w\cdot P=\bar P'$.

We will now show that the restriction to $\caM(G)_\Omega$ of the functor $r_D$ is faithful. By Lemma
\ref{exaCt}, it suffices to show that if $\pi\in \caM(G)_\Omega$, then $r_D(\pi)\neq 0$. We can assume
$\pi$ is irreducible by the exactness of $r_D$. Then by $(i)$, $\pi$ is written $\pi=i_P^G(\rho)$ with
$\rho \in D$, or else $\pi=i_{P'}^G(\rho')$ with $\rho'\in D'$. In both cases, Proposition
\ref{geolemmacusp} gives $r_D(\pi)$ non-zero.

The functor $i_P^G \colon  \caM(M)_D \rightarrow  \caM(G)_\Omega$ is the right adjoint of
$r_D \colon \caM(G)_\Omega  \rightarrow \caM(M)_D$. The adjunction morphism (see \ref{adjfonct})
$$\alpha_\pi\colon \pi  \rightarrow i_P^Gr_D(\pi)$$ is injective: since
$\beta_{r_D(\pi)}\circ r_D(\alpha_\pi)$ is the identity, $r_D(\alpha_\pi)$ is injective, where $\beta$ is
the adjunction morphism below; since $r_D$ is exact, we therefore have $r_D(\ker\alpha_\pi)=0$, hence
$\ker\alpha_\pi=0$ since $r_D$ is faithful. This is the argument of the second lemma of \ref{dcompthm},
$(iii)$, applied to $r_D$. On the other hand, the geometric lemma \ref{geolemmacusp} shows that the
adjunction morphism $$\beta_\tau \colon   r_Di_P^G \tau \rightarrow \tau,$$ $\tau \in \caM(M)_D$, is an
isomorphism. Applying the exact functor $r_D$ to $\alpha_\pi$, we obtain an injection:
\[
r_D(\alpha_\pi) \colon  r_D(\pi)  \rightarrow r_Di_P^Gr_D(\pi).
\]
Composing this with $\beta_{r_D(\pi)}$ we obtain
\[
r_D(\pi)\stackrel{ r_D(\alpha_\pi)}  {\longrightarrow} r_Di_P^Gr_D(\pi)
\stackrel{ \beta_{r_D(\pi)}}  {\longrightarrow} r_D(\pi)
\]
which is the identity of $r_D(\pi)$ using the properties of adjunction morphisms ({\sl cf.} Appendix
\ref{adjfonct}). This implies that $r_D(\alpha_\pi)$ is an isomorphism. We deduce, since $r_D$ is exact
and faithful, that $\alpha_\pi$ is a (natural) isomorphism. The adjunction morphisms $\alpha$ and $\beta$
are therefore respectively natural isomorphisms between the identity of $\caM(G)_\Omega$ and
$i_P^G \circ r_D$ and between $r_D \circ i_P^G$ and the identity of $\caM(M)_D$. This shows the
equivalence of categories.
\end{proof}

\section{Tempered representations}\label{sec:temprep}
In this section, we develop the algebraic theory of tempered representations. For the analytic aspects
(asymptotic behavior of the character and matrix coefficients), we refer the reader to \cite{Wald}.

\subsection{Definition by normalized exponents} \label{temp1}

\enlargethispage{2\baselineskip}
\begin{defi}
A representation $(\pi,V)$ of $\caM(G)$ is said to be tempered if it is admissible, and if for any
parabolic subgroup $P=MN$, any normalized exponent $\chi$ of $\pi$ with respect to $P$ satisfies
$\Re (\chi) \in {}^+\overline{[\fra_M^*]}_P^G$.
\end{defi}

\begin{rmqs}
1. If $(\pi,V)$ is tempered, any normalized exponent $\chi$ of $\pi$ with respect to $G$ is unitary ({\sl
i.e.,} $\Re (\chi)=0$). In particular, if $(\pi,V)$ is tempered and admits a central character, it is
unitary.

--- 2. If $(\pi,V)$ is irreducible and square-integrable modulo the center, then by Casselman's criterion
\ref{CasCrit}, $\pi$ is tempered.

--- 3. Since any parabolic subgroup of $G$ is conjugate to a standard parabolic subgroup, for $(\pi,V)$
to be tempered, it suffices that the condition $\Re (\chi) \in {}^+\overline{[\fra_M^*]}_P^G$ be
satisfied for any normalized exponent $\chi$ of $\pi$ with respect to any standard parabolic subgroup $P$
of $G$.
\end{rmqs}

Our immediate goal is to determine how tempered representations behave with respect to the induction
functors $i_P^G$ and restriction functors $r_P^G$.

\subsection{Tempered representations and induction} \label{temp2}

Parabolic induction functors preserve tempered representations:
\begin{lemme}
Let $P=MN$ be a parabolic subgroup of $G$ and $(\sigma,E)$ a tempered representation of $M$. Then
$i_P^G (\sigma,E)$ is tempered.
\end{lemme}

\begin{proof}
It suffices to verify this for standard parabolic subgroups $P$. Let $Q=LU$ be a standard parabolic
subgroup of $G$. Let $\chi$ be a normalized exponent of $i_P^G(\sigma,E)$ with respect to $Q$. We want to
show that $\Re(\chi) \in {}^+\overline{[\fra_L^*]}_Q^G$.

By the geometric lemma \ref{geomlemma}, $(r_Q^G  i_P^G (\sigma,E))_\chi$ admits a filtration whose
associated graded module has as components the representations
\[
(i_{L\cap w^{-1}\cdot P }^L \circ w \circ r_{w\cdot Q \cap M}^M (\sigma,E))_\chi,
\]
where $w$ runs through the set $W^{L,M}$ defined in (\ref{Wml}). Since $\sigma$ is tempered, any
normalized exponent $\psi$ with respect to $w \cdot Q \cap M$ satisfies
\[
\Re (\psi) \in   {}^+\overline{[\fra_{w \cdot L\cap M}^*]}_{w\cdot
Q\cap M}^M,
\]
i.e., we can write
\[
\Re (\psi) =\sum_{\alpha \in \Delta^M (w \cdot Q \cap M)} c_\alpha
\alpha
\]
where the $c_\alpha$ are $\geq 0$.

We have
\[
\fra_\emptyset^*= \fra_{w \cdot L \cap M}^* \oplus
(\fra_\emptyset^{w \cdot L \cap M}  )^*
\]
and the $\alpha \in  \Delta^M (w \cdot Q \cap M)$ are elements of $\fra_{w \cdot L \cap M}^*$ which can
therefore be viewed as elements of $\fra_\emptyset^*$. Such an $\alpha \in  \Delta^M (w \cdot Q \cap M)$
can then be written
\[
\alpha= \beta+\mu
\]
where $\beta \in  \Delta^M_\emptyset \setminus \Delta^{w \cdot L \cap M}_\emptyset$ and
$\mu \in (\fra_\emptyset^{w \cdot L \cap M}  )^*   $.

This allows us to write $\Re (\psi) $ in the form
\[
\Re (\psi) =\sum_{\alpha \in  \Delta^M_\emptyset \setminus \Delta^{w \cdot L \cap
M}_\emptyset } c_\alpha
\alpha  \, +\mu
\]
with $\mu \in (\fra_\emptyset^{w \cdot L \cap M}  )^*   $.

Then
\[
\Re (w^{-1}  \cdot \psi) =\sum_{\alpha \in  \Delta^M_\emptyset \setminus \Delta^{ M \cap
w\cdot  L}_\emptyset } c_\alpha \;  w^{-1}\cdot \alpha  \, +
w^{-1}\cdot \mu
\]
with $w^{-1}\cdot \mu \in (\fra_\emptyset^{L \cap w^{-1}\cdot M}  )^*   $. Since
\[
\fra_\emptyset^*= \fra_{L \cap w^{-1}  \cdot  M}^* \oplus
(\fra_\emptyset^{ L \cap w^{-1}  \cdot M}  )^* =   \fra_{L}^* \oplus
(\fra^L_{ L \cap w^{-1}  \cdot M}  )^*  \oplus
(\fra_\emptyset^{ L \cap w^{-1}  \cdot M}  )^*,
\]
we have $(w^{-1}\cdot \mu)_{|\fra_L}=0$.

The choice of the system of representatives $W^{L,M}$ is now crucial, because by (\ref{Wml}), if
$\alpha \in \Delta^M_\emptyset$, then $w^{-1} \cdot \alpha \in \Sigma_\emptyset^+$. This shows that
$$ \Re (w^{-1}  \cdot \psi) - w^{-1}\cdot \mu   \in  {}^+\overline{[\fra_\emptyset^*]}^G_{P_\emptyset} .$$

By Lemma \ref{exposants}, any character $\chi_1$ of $A_L$ such that
\[
(i_{L\cap w^{-1}\cdot P }^L \circ w \circ r_{w\cdot Q \cap M}^M
\sigma)_{\chi_1} \neq 0
\]
is the restriction to $A_L$ of a character $\psi_1$ of $A_{L\cap w^{-1}\cdot M}$ such that
$$(w \circ  r_{w\cdot Q \cap M}^M \sigma)_{\psi_1}\neq 0.$$ Such a character $\psi_1$ is of the form
$w^{-1}  \cdot \psi$ with $\psi$ as above. This shows that $\Re(\psi_1)$ is the restriction to $\fra_L$
of an element of ${}^+\overline{[\fra_\emptyset^*]}^G_{P_\emptyset} $. Such an element is in
${}^+\overline{[\fra_L^*]}^G_{Q} $.
\end{proof}

\subsection{Tempered representations and restriction} \label{temp3}

We now turn to the study of the behavior of tempered representations under the action of the functor
$r_P^G$.

\begin{lemme}
Let $(\pi,V)$ be a tempered representation of $G$. Let $P=MN$ be a standard parabolic subgroup of $G$.
Decompose the Jacquet module $r_P^G V$ into a direct sum of representations of $M$,
\[
r_P^G V=  (r_P^G V)^0  \oplus (r_P^G V)^+
\]
where
\[
(r_P^G V)^0  = \bigoplus_{\chi \in \mathrm{Exp}(A_M,r_P^G V)|
\Re(\chi)=0}  (r_P^G V)_\chi,
\]
\[
(r_P^G V)^+  =\bigoplus_{\chi \in \mathrm{Exp}(A_M,r_P^G V)|
\Re(\chi)\neq 0} (r_P^G V)_\chi
\]
The representation $(r_P^G \pi)^0$ is then a tempered representation of $M$.
\end{lemme}

\begin{proof}
Let us fix $\chi\in \mathrm{Exp}(A_M,r_P^G V)$ such that $\Re(\chi)=0 $. A standard parabolic subgroup of
$M$ is the intersection with $M$ of a standard parabolic subgroup of $G$, say $Q=LU$, contained in $P$.
The transitivity of the restriction functors (see \ref{associa}) shows us that any normalized exponent
$\chi_Q$ with respect to $M\cap Q$ of $(r_P^G \pi)_\chi$ is a normalized exponent of $\pi$ with respect
to $Q$. Since $\pi$ is tempered, we have,
\begin{equation} \label{chiQ}
\Re(\chi_Q) \in {}^+\overline{[\fra_L^*]}_Q^G.
\end{equation}

Recall the decomposition \ref{traceP}:
\[
\fra_L^*=  \fra_M^* \oplus (\fra_L^M)^*
\]
By Lemma \ref{exposants} $(ii)$, the restriction of $\chi_Q$ to $A_M$ is equal to $\chi$. The projection
of $\Re (\chi_Q) $ into $\fra^*_M$ is therefore equal to $\Re (\chi)=0$. We deduce that
$ \Re (\chi_Q) \in  (\fra_L^M)^*$, whence
\[
\Re (\chi_Q) \in  (\fra_L^M)^* \cap  {}^+\overline{[\fra_L^*]}_Q^G
\subset   {}^+\overline{[(\fra_L^M)^*]}_{Q\cap M}^M .
\]
\end{proof}

\subsection{Frobenius reciprocity for tempered representations} \label{caractemp}

We obtain the following version of Frobenius reciprocity for tempered representations:
\begin{prop}
Let $P=MN$ be a parabolic subgroup of $G$. For any tempered representation $(\sigma,E)$ of $M$ and for
any tempered representation $(\pi,V)$ of $G$, we have a natural isomorphism
\[
\Hom_G(V,i_P^G E) \simeq \Hom_M(r_P^G (V)^0, E)
\]
\end{prop}
\begin{proof}
This result is the conjunction of Frobenius reciprocity \ref{frobPG} and Lemma \ref{caraccentral}.
\end{proof}

\begin{lemme}
Let $(\pi,V)$ be a smooth admissible tempered representation of $G$ admitting a central character. Then
$\pi$ is square-integrable modulo the center if and only if for any proper standard parabolic subgroup
$P=MN$, $(r_P^G \pi)^0=\{0\}$.
\end{lemme}

\begin{proof}
This follows from the definitions, the fact that the central character is unitary (first remark of
\ref{temp1}) and Casselman's criterion.
\end{proof}

\subsection{Geometric lemma for tempered representations}\label{geolemtemp}

We will now establish a version of the geometric lemma for tempered representations.

\begin{prop}
Let $P=MN$ and $Q=LU$ be two parabolic subgroups of $G$ and $(\rho,W)$ a tempered representation of $M$.
Then $r_Q^G(i_P^G(\rho,W))^0$ admits a filtration whose associated graded components are the
\[
i_{L\cap w^{-1}\cdot P}^L \circ w \circ (r_{w\cdot Q \cap M}^M
W)^0,
\]
where $w$ runs through the set $\caW^{Q,P}$ defined in \ref{PWQsemi} parameterizing the orbits of $Q$ on
$X=P\backslash G$.
\end{prop}

\begin{proof}
The only difference with respect to the first theorem of \ref{geomlemma} is that we take the tempered
part $(r_{w\cdot Q \cap M}^M W)^0$ of $r_{w\cdot Q \cap M}^M W$ in the formula. We can assume $P$ and $Q$
are standard, and in this case $\caW^{Q,P}$ is the set denoted $W^{L,M}$ in \ref{WeylgroupsS}. It is
clear that if $w \in W^{L,M}$ and if
\[
\chi \in \mathrm{Exp}(A_L, i_{L\cap w^{-1}\cdot P}^L \circ w \circ (r_{w\cdot Q \cap M}^M
W)^0)
\]
then $\chi$ is unitary (first remark of \ref{temp1}), since the representation
$i_{L\cap w^{-1}\cdot P}^L \circ w \circ (r_{w\cdot Q \cap M}^M W)^0$ is tempered, thus $\Re(\chi)=0$.
Conversely, we will now show that if $w \in W^{L,M}$ and if
\[
\chi \in \mathrm{Exp}(A_L,  i_{L\cap w^{-1}\cdot P}^L \circ w \circ (r_{w\cdot Q \cap M}^M
W)^+),
\]
we have $\Re(\chi)\neq 0$. Let us fix such $w,\chi$. As in the proof of Lemma \ref{temp2}, there exists
$\chi' \in \mathrm{Exp}(A_{w\cdot L\cap M}, r_{w\cdot Q \cap M}^M W)$ with $\Re(\chi')\neq 0$ such that
$\Re(\chi)$ is the restriction of $\Re(w^{-1}  \cdot \chi')$ to $\fra_L$. We write, as in the proof of
Lemma \ref{temp2},
\[
\Re (\chi') =\sum_{\alpha \in  \Delta^M_\emptyset \setminus \Delta^{w \cdot L \cap
M}_\emptyset } c_\alpha \alpha  \, +\mu
\]
with $\mu \in (\fra_\emptyset^{w \cdot L \cap M}  )^*   $. Since $\Re(\chi')\neq 0$, one of the
$c_\alpha$ is $>0$.

We then have
\[
\Re (w^{-1}  \cdot \chi') =\sum_{\alpha \in  \Delta^M_\emptyset
\setminus \Delta^{ w\cdot L  \cap  M}_\emptyset } c_\alpha \;  w^{-1}\cdot \alpha  \, +
w^{-1}\cdot \mu
\]
with $w^{-1}\cdot \mu \in (\fra_\emptyset^{L \cap w^{-1}\cdot M}  )^*   $ and thus
$(w^{-1} \cdot \mu)_{|\fra_L}=0$.

It therefore suffices to prove that $\sum_\alpha c_\alpha(w^{-1}\cdot \alpha)_{|\fra_L} \neq 0$. Let us
first show that if $\alpha \in  \Delta^M_\emptyset \setminus \Delta^{w \cdot L \cap M}_\emptyset $, then
$(w^{-1}\cdot \alpha)_{|\fra_L} \neq 0$. Now, by definition, there exists a non-zero subspace
$\frg_\alpha$ of the Lie algebra of $w \cdot U \cap  M$ in which $A_\emptyset$ acts by the root $\alpha$.
Then $A_\emptyset$ acts in $w^{-1}\cdot \frg_\alpha$, a subspace of the Lie algebra of
$U \cap w^{-1}\cdot M$, by $w^{-1}\cdot \alpha$. Since $0$ is not an eigenvalue of the action of $A_L$ in
the Lie algebra of $U$, we indeed have $(w^{-1}\cdot \alpha)_{|\fra_L} \neq 0$, and it is a root of $A_L$
in $U$, that is, a linear combination with non-negative coefficients, not all zero, of the elements of
$\Delta(Q)$. Since the $c_\alpha$ are $\geq 0$, and one of them is $>0$, there can be no cancellation in
the sum, which is therefore non-zero.
\end{proof}

\subsection{Irreducible tempered representations}

The following result shows that tempered representations are obtained as subrepresentations of induced
representations of discrete series.

\needspace{151pt}
\begin{thm}
Let $(\pi,V)$ be an irreducible smooth representation of $G$. Then $(\pi,V)$ is tempered if and only if
there exists a standard parabolic subgroup $P=MN$ of $G$ and an irreducible square-integrable
representation modulo the center $(\sigma,E)$ of $M$ such that $(\pi,V)$ is a subrepresentation of
$i_P^G (\sigma,E)$.

If $(P_1=M_1N_1, (\sigma_1,E_1))$ and $(P_2=M_2N_2, (\sigma_2,E_2))$ are two pairs satisfying this
property, then there exists $g \in G$ such that $gM_1g^{-1}=M_2$ and $\sigma_1^g \simeq \sigma_2$.
\end{thm}
\begin{proof}
The second remark of \ref{temp1} asserts that $(\sigma,E)$ is a tempered representation of $M$, and we
have seen that normalized induction preserves tempered representations. This shows that the condition is
sufficient. To prove necessity, let us choose a standard parabolic subgroup $P=MN$ of $G$, minimal with
respect to the property that $(r_P^G V)^0 \neq 0$. Since this property is satisfied for $P=G$, the
existence of such a parabolic subgroup is ensured. Since $\pi$ is irreducible, $r_P^G (\pi,V)^0$ is
finitely generated (Proposition \ref{finiIJ}) and thus admits an irreducible quotient $(\sigma,E)$, by
\ref{JoHo}. By Lemma \ref{caractemp} and the minimality of $P$, we see that $(\sigma,E)$ is
square-integrable modulo the center. By Frobenius reciprocity, $\Hom_G(V,i_P^G E)\neq \{0\}$. Since
$(\pi,V)$ is irreducible, it appears as a subrepresentation of $i_P^G (\sigma,E)$. This completes the
proof of the first point.

For the second point, note first that a representation square-integrable modulo the center $(\sigma,E)$
is unitary, and thus that its normalized induced representation $i_P^G (\sigma,E)$ is as well. If
moreover $(\sigma,E)$ is irreducible, $i_P^G (\sigma,E)$ is of finite length, thus is semisimple, and any
subquotient is also a subrepresentation and a quotient.

We have $\Hom_G(V,i_{P_2}^G E_2)\neq 0$ and thus by Frobenius reciprocity for tempered representations,
$$\Hom_{M_2}((r_{P_2}^G V)^0,E_2)\neq 0.$$ Since $V$ is a quotient of $i_{P_1}^G E_1$, by the exactness
of the functor $r_{P_2}^G$, $(r_{P_2}^G V)^0$ is a quotient of
\[
(r_{P_2}^G i_{P_1}^G E_1)^0
\]
thus $\Hom_{M_2} ((r_{P_2}^G i_{P_1}^G E_1)^0, E_2)\neq 0$. By Proposition \ref{geolemtemp}, there exists
$w \in W^{M_2,M_1}$ such that
\[
\Hom_{M_2}(i_{M_2\cap w^{-1}P_1}^{M_2} \circ w \circ (r_{M_1 \cap
w\cdot P_2}^{M_1} E_1)^0,E_2)  \neq 0.
\]
Let us fix such a $w$. In particular $(r_{M_1 \cap w\cdot P_2}^{M_1} E_1)^0\neq 0$ and thus by Lemma
\ref{caractemp}, $M_1 \cap  w\cdot P_2=M_1$, i.e., $M_1 \subset w\cdot P_2$. We then have
\[
\Hom_{M_2}(i_{M_2\cap w^{-1}P_1}^{M_2} \circ w \circ E_1,E_2)  \neq 0,
\]
which gives, since the representations appearing are semisimple,
\[
\Hom_{M_2}(E_2, i_{M_2\cap w^{-1}P_1}^{M_2} \circ w \circ E_1)  \neq 0.
\]
By Frobenius reciprocity, we obtain
\[
\Hom_{M_2 \cap w^{-1}\cdot M_1} ((r_{M_2\cap w^{-1}\cdot P_1}^{M_2} \sigma_2)^0,  \sigma_1^w)\neq 0.
\]
In particular $(r_{M_2\cap w^{-1}\cdot P_1}^{M_2} E_2)^0\neq 0$, which shows that
$M_2\cap w^{-1}\cdot P_1 = M_2$, i.e., $M_2 \subset w^{-1}\cdot P_1$. We deduce $M_2=w^{-1}\cdot M_1$ and
$$\Hom_{M_2}(\sigma_2, \sigma_1^{w})\neq 0,$$ whence $\sigma_2\simeq\sigma_1^w$.
\end{proof}

\begin{cor}
An irreducible tempered representation is unitary.
\end{cor}
\begin{proof}
This follows from the theorem because an irreducible square-integrable representation is unitary
(\ref{unit2int}) and parabolic induction preserves unitarity.
\end{proof}

\begin{rmq}
Without the irreducibility hypothesis in the preceding corollary, the conclusion of the latter can be
false.
\end{rmq}

\section[Self-intertwining of induced representations]{Intertwining operators of induced representations}\label{opentre1}

\subsection{Bruhat order}\label{opentre}

Let $P$ and $Q$ be two parabolic subgroups of $G$ with the same Levi factor $M$. Let $(\rho,W)$ be a
smooth representation of $M$. This section aims to construct a canonical intertwining operator between
$i_P^G\rho$ and $i_Q^G\rho$, i.e., an element of $\Hom_G(i_P^G\rho, i_Q^G\rho)$. By Frobenius reciprocity
\[
\Hom_G(i_P^G\rho, i_Q^G\rho) \simeq \Hom_M(r_Q^Gi_P^G\rho,\rho).
\]
Thus we seek to exhibit an element of $ \Hom_M(r_Q^Gi_P^G\rho,\rho)$.

We assume that the Levi factor $M$ is standard, and thus that the parabolic subgroups $P$ and $Q$ are
semi-standard. This restriction is not essential, but of pure convenience, allowing us to use notation
already introduced. The reader can easily adapt the proofs to the general case, or deduce the results
from the special case, by conjugation.

Recall that the geometric lemma \ref{geomlemma} gives a filtration of the functor $r_Q^G\circ i_P^G$:
\[
0=F_0\subset F_1\ldots \subset F_k=r_Q^G\circ i_P^G\colon \caM(M) \rightarrow \caM(M)
\]
such that:
\[
F_i/F_{i-1} \simeq i_{M\cap w_i^{-1}\cdot P}^M \circ w_i
\circ r_{M \cap w_i \cdot Q}^M,
\]
where the $w_i$ are representatives of the double cosets $\bar w_i \in P\backslash G /Q$, numbered
according to a total order specified in \ref{geomlemma}. In this context, it is convenient to
parameterize $P\backslash G /Q$ not by the system of representatives $\caW^{Q,P}$ of \ref{PWQsemi}, which
depends on $P$ and $Q$, but simply by the set of double cosets $W_M\backslash W_G /W_M$ (see
\ref{WMWGWM}) which is independent of them. On the other hand, we equip this set with an order which
depends on $P$ and $Q$:

\begin{defi}
(Bruhat order) We equip $W_M\backslash W_G /W_M$ with the partial order
\[
\bar w \leq_{PQ}  \bar w'  \text{ if } P\bar w Q\subset \overline{ P\bar w' Q}.
\]
When $P$ and $Q$ are clearly indicated by the context, we write $\leq$ instead of $\leq_{PQ}$.
\end{defi}

We denote by $\bar 1 \in  W_M\backslash W_G /W_M$ the double coset of the identity element of $W_G$. We
then choose a total order $\preceq$ on $W_M\backslash W_G /W_M$ which refines the Bruhat order defined
above. Thus, if $P=Q$, $\bar 1$ is the smallest element for the order $\preceq$. On the contrary, if
$Q=\bar P$, it is the largest. Recall that in \ref{geomlemma}, we had equipped the set of orbits of $Q$
in $P\backslash G$ with a total order
\[
Z_1<Z_2 \cdots < Z_m
\]
so that for all $i=1,\ldots,m$, $\bigcup_{j\leq i}Z_j$ is open in $P\backslash G$. This order can be
chosen such that it is the inverse of the order $\preceq$ on
$W_M\backslash W_G /W_M\simeq P\backslash G/Q$.

\subsection{Filtrations}

Set, for any smooth representation $(\rho,W)$ of $M$, and any element
$\bar w \in  W_M\backslash W_G /W_M$,
\[
\widetilde F^{\preceq \bar w}_{PQ}(W):= \{f \in i_P^G(W), \;
\mathrm{Supp} (f) \cap \left( \bigcup_{\bar w' \preceq \bar w}
P\bar w' Q      \right)=\emptyset     \}
\]
and similarly let us define $\widetilde F^{\prec \bar w}_{PQ}(W)$ with $\prec$ instead of $\preceq$. Also
set
\[
\widetilde F^{< \bar w}_{PQ}(W):= \{f \in i_P^G(W), \;
\mathrm{Supp} (f) \cap \overline{P\bar w Q}\subset P\bar w Q    \}.
\]
and
\[
\widetilde F^{\leq   \bar w}_{PQ}(W):= \{f \in i_P^G(W), \;
\mathrm{Supp} (f) \cap \overline{P\bar w Q}=\emptyset   \}.
\]

Note that we then always have:
\[
\widetilde F^{\preceq \bar w}_{PQ}(W) \subset \widetilde F^{\prec
\bar w}_{PQ}(W),
\quad \widetilde F^{\prec \bar w}_{PQ}(W) \subset \widetilde F^{< \bar w}_{PQ}(W),
\quad \widetilde F^{\leq \bar w}_{PQ}(W) \subset \widetilde F^{< \bar w}_{PQ}(W).
\]

All the subspaces of $i_P^G(W)$ defined above are stable under the action of $Q$. Let
\[
F^{\preceq \bar w}_{PQ}(W), F^{\prec \bar w}_{PQ}(W), F^{\leq \bar w}_{PQ}(W),
F^{< \bar w}_{PQ}(W),
\]
denote the respective images of
\[
\widetilde F^{\preceq \bar w}_{PQ}(W), \widetilde F^{\prec
\bar w}_{PQ}(W),  \widetilde F^{\leq \bar w}_{PQ}(W), \widetilde
F^{< \bar w}_{PQ}(W)
\]
in $r_Q^Gi_P^G(W)$.

The functors $F^{\preceq \bar w}_{PQ}$, $F^{\prec \bar w}_{PQ}$, $F^{\leq\bar w}_{PQ}$,
$F^{< \bar w}_{PQ}$ are subfunctors of $r_Q^G\circ i_P^G$. It is clear that
\[
F^{ \prec \bar w}_{PQ}/F^{\preceq \bar w}_{PQ}\simeq   F^{< \bar w}_{PQ}/ F^{\leq \bar w}_{PQ},
\]
and the geometric lemma is reformulated as
\[
F^{ \prec \bar w}_{PQ}/ F^{\preceq \bar
w}_{PQ}\simeq   i_{M\cap w^{-1}\cdot P}^M \circ w \circ r_{w\cdot Q \cap M}^M.
\]

For the element $\bar 1$, this gives
\[
F^{ \prec \bar 1 }_{PQ}/ F^{\preceq \bar
1}_{PQ}\simeq F^{< \bar 1}_{PQ}/  F^{\leq \bar
1}_{PQ}\simeq  \Id_{\caM(M)},
\]
the right-hand side being the identity functor of the category $\caM(M)$.

In other words, if $(\rho,W)$ is a smooth representation of $M$, $\rho$ appears as a subquotient in
$r_Q^G\circ i_P^G(\rho)$. We show that under certain conditions, this subquotient appears as a quotient
of $r_Q^G\circ i_P^G(\rho)$, thereby providing a non-zero element of
\[
\Hom_{M}(r_Q^G\circ i_P^G(\rho), \rho).
\]

\subsection{$PQ$-regularity}\label{PQreg}

Let $(\sigma,E)$ be a smooth representation of $M$. Let us restrict this representation to the split
component $A_M$ of $M$, and then extend it to the group algebra $\bbC[A_M]$. This gives us an algebra
morphism
\[
\sigma_{|\bbC[A_M]} \colon  \bbC[A_M] \rightarrow    \End_\bbC(E)
\]
Let us then define the following ideals of $\bbC[A_M]$:

\begin{align*}
 I_\sigma&= \ker  \sigma_{|\bbC[A_M]} \\
 I_\sigma^{PQ} &= \ker \left[ r_Q^G\circ i_P^G/F_{PQ}^{< \bar 1}(E)  \right]_{|\bbC[A_M]}      \\
I_\sigma^{Q,\bar w}&= \ker \left[w \circ r_{M\cap w \cdot Q}^M(E)
\right]_{|\bbC[A_M]}  , \quad \bar w \in W_M\backslash W_G/W_M
\end{align*}

Similarly, when $\chi$ is a character of $A_M$, we still denote by $\chi$ its extension to an algebra
morphism from $\bbC[A_M]$ to $\bbC$.

\begin{lemme}
The following properties are equivalent:

$(i)$ $I_\sigma + I_{\sigma}^{PQ}=\bbC[A_M]$

$(ii)$ $I_\sigma + \bigcap_{\bar w < \bar 1 }  I_\sigma^{Q,\bar w}=\bbC[A_M]$

Moreover, $(i)$ or $(ii)$ implies:

$(iii)$
$\forall \bar w < \bar 1, \quad \mathrm{Exp}(A_M,\sigma) \cap \mathrm{Exp}(A_M,  w \circ r_{M\cap w \cdot Q}^M(\sigma))= \emptyset.$
Conversely, $(iii)$ implies $(i)$ and $(ii)$ if $(\sigma,E)$ is of finite length.

When one of the first two conditions is satisfied, we say that $\sigma$ is $PQ$-regular.
\index[ter]{PQreg@$PQ$-regular}
\end{lemme}

\begin{proof}
We choose the total order $\preceq$ on $W_M\backslash W_G/W_M$ finer than the Bruhat order $\leq_{PQ}$ so
that $\{\bar w <\bar 1\}=\{ \bar w \prec \bar 1\}$. The representation
$(r_Q^G\circ i_P^G/F_{PQ}^{< \bar 1})(E)$ admits a filtration whose subquotients are isomorphic to the
\[
i_{M\cap w^{-1}\cdot P}^M \circ w
\circ r_{M \cap w \cdot Q}^M  (\sigma),  \quad (\bar w < \bar 1).
\]
Since the parabolic induction functors are faithful, we have
\[
\ker \left[  i_{M\cap w^{-1}\cdot P}^M \circ w
\circ r_{M \cap w\cdot Q}^M (\sigma)  \right]_{|\bbC[A_M]}=
\ker \left[   w \circ r_{M \cap w \cdot Q}^M (\sigma)  \right]_{|\bbC[A_M]}.
\]
Since the number of $\bar w < \bar 1$ is bounded by $|W_G|$, we have
\[
\left(  \bigcap_{\bar w < \bar 1 } I_\sigma^{Q,\bar w} \right)^{|W_G|}  \subset
I_{\sigma}^{PQ} \subset     \bigcap_{\bar w < \bar 1 } I_\sigma^{Q,\bar w} .
\]
The second inclusion shows that $(i)$ implies $(ii)$. For the converse, it suffices to see that if $I$
and $J$ are two ideals of a commutative unital ring $A$ such that $I+J=A$, then $I^n+J^m=A$, for all
strictly positive integers $n$ and $m$. We can write $\una=i+j$, with $i\in I$ and $j \in J$, and thus
$\una=(\una)^{n+m}=(i+j)^{n+m}$ which is in $I^n+J^m$ by the binomial formula, which shows the assertion.

Suppose now that $(\sigma,E)$ is of finite length. By the definition of normalized exponents, and the
fact that they appear in finite number, there exists an integer $d$ such that:
\[
\left( \bigcap_{\chi \in \mathrm{Exp}(A_M,E)} \ker \chi \right)^d
\subset I_\sigma \subset \bigcap_{\chi \in \mathrm{Exp}(A_M,E)} \ker \chi
\]
and
\[
\left( \bigcap_{\chi \in  \cup_{\bar w < \bar 1}
\mathrm{Exp}(A_M,w \circ r_{M\cap w \cdot Q }^M (E))} \ker \chi \right)^d
\subset I_\sigma^{PQ} \subset \bigcap_{\chi \in  \cup_{\bar w < \bar 1}
\mathrm{Exp}(A_M,w \circ r_{M\cap w \cdot Q }^M (E))} \ker \chi
\]
Note that the two right inclusions do not assume that $(\sigma,E)$ is of finite length.

If $(iii)$ fails due to some $\chi_0$, then
\[
I_\sigma \subset \bigcap_{\chi \in \mathrm{Exp}(A_M,E)} \ker
\chi \subset \ker \chi_0,
\]
\[
I_\sigma^{PQ}  \subset   \bigcap_{\chi \in  \cup_{\bar w < \bar 1}
\mathrm{Exp}(A_M,w \circ r_{M\cap w \cdot Q }^M (E))}
\ker \chi  \subset \ker \chi_0
\]
and thus
\[
I_\sigma +  I_\sigma^{PQ}  \subset \ker \chi_0 \neq \bbC[A_M],
\]
which shows that $(i)$ is not satisfied.

Conversely, suppose $(iii)$, and suppose that
\[
\left( \bigcap_{\chi \in \mathrm{Exp}(A_M,E)} \ker \chi
\right)  + \left( \bigcap_{\chi \in  \cup_{\bar w < \bar 1}
\mathrm{Exp}(A_M,w \circ r_{M\cap w \cdot Q}^M (E))} \ker \chi \right) \neq  \bbC[A_M] .
\]
Then there exists a (proper) maximal ideal $\frM$ of $\bbC[A_M]$ containing the left-hand side, i.e.,
that
\[
\left( \bigcap_{\chi \in \mathrm{Exp}(A_M,E)} \ker \chi
\right) \subset \frM , \quad   \left( \bigcap_{\chi \in  \cup_{\bar w < \bar 1}
\mathrm{Exp}(A_M,w \circ r_{M\cap w \cdot Q }^M (E))} \ker \chi \right) \subset \frM  .
\]
Since the sets of exponents are finite, and each $\ker \chi$ is a maximal ideal of $\bbC[A_M]$, we see
that $\frM$ is equal to $\ker \chi_0$ for a certain $\chi_0$ both in $\mathrm{Exp}(A_M,E)$ and in
$\cup_{\bar w < \bar 1} \mathrm{Exp}(A_M,w \circ r_{M\cap w \cdot Q }^M (E))$, which is impossible by
$(iii)$. We therefore have
\[
\left( \bigcap_{\chi \in \mathrm{Exp}(A_M,E)} \ker \chi
\right)  + \left( \bigcap_{\chi \in  \cup_{\bar w < \bar 1}
\mathrm{Exp}(A_M,w \circ r_{M\cap w \cdot Q}^M (E))} \ker \chi \right) =  \bbC[A_M] ,
\]
whence by an argument given above,
\[
\left( \bigcap_{\chi \in \mathrm{Exp}(A_M,E)} \ker \chi
\right)^d + \left( \bigcap_{\chi \in  \cup_{\bar w < \bar 1}
\mathrm{Exp}(A_M,w \circ r_{M\cap w \cdot Q }^M (E))} \ker \chi \right)^d= \bbC[A_M]
\]
and thus $I_\sigma +  I_\sigma^{PQ}  =\bbC[A_M]$.
\end{proof}

\subsection{Intertwining operators}\label{INTOP}
Consider the inclusion
\[
E \simeq F_{PQ}^{< \bar 1}/F_{PQ}^{\leq \bar 1}(E)
\stackrel{i_\sigma}{\hookrightarrow}
(r_Q^G \circ i_P^G/ F_{PQ}^{\leq \bar 1})(E).
\]
If the representation $(\sigma,E)$ is $PQ$-regular, since all the other subquotients (other than
$\sigma$) of the right-hand side have by hypothesis exponents which are not those of $\sigma$, the
representation theory of the group $A_M$ shows that this inclusion $i_\sigma$ admits a retraction
$r_\sigma$. More explicitly, this retraction is given by the action of an element
$\epsilon_\sigma \in I_\sigma^{PQ}$ such that $\mathbf{1}_{\bbC[A_M]}\in \epsilon_\sigma + I_\sigma$. In
particular, it commutes with the action of $M$. The composition
\[
K_{Q|P}(\sigma)   \colon r_Q^G \circ i_P^G (E)
\rightarrow   (r_Q^G \circ i_P^G/ F_{PQ}^{\leq \bar 1})(E)\stackrel{r_\sigma}{\longrightarrow} E
\]
is a non-zero element in
\[
\Hom_M(r_Q^G\circ i_P^G(\sigma),\sigma)
\]
which by Frobenius reciprocity gives a non-zero intertwining operator
\[
J_{Q|P}(\sigma) \in \Hom_G(i_P^G(\sigma),i_Q^G(\sigma)).
\]

It is useful to provide a more explicit form of the operator $K_{Q|P}$. Recall that the natural
isomorphism
\[
F_{PQ}^{< \bar 1}/F_{PQ}^{\leq \bar 1}(E)\simeq E
\]
is obtained by passing to the quotient from the morphism
\begin{align*}
    \widetilde F_{PQ}^{< \bar 1} &\rightarrow E   \\
              f &\mapsto \int_{U\cap N \backslash U} f(u) \, d\nu_{U\cap N
              \backslash U}(u)=\int_{U\cap \overline{N} } f(u)  \, du
\end{align*}
This is indeed a special case of formula (\ref{IndRes}), and the equality between the two integrals comes
from the decomposition $$U=(U\cap N)(U\cap \overline{N}).$$ We can see directly from this formula that
these integrals are convergent, since by definition of $\widetilde F_{PQ}^{< \bar 1}$,
$\mathrm{Supp} (f)\cap PQ$ is compact modulo $P$, thus $\mathrm{Supp} (f_{|U})$ is compact modulo
$P\cap U=N\cap U$.

\begin{prop}
Let $P=MN$ and $Q=MU$ be two semi-standard parabolic subgroups of $G$ and $(\sigma,E)$ a smooth
representation of $M$ which is assumed to be $PQ$-regular. Then the intertwining operator
$J_{Q|P}(\sigma)$ is the unique element of $\Hom_G(i_P^G(\sigma), i_Q^G(\sigma))$ satisfying
\[
(\forall f \in \widetilde F_{PQ}^{<\bar 1}(E)), \quad
J_{Q|P}(\sigma)(f)(\mathbf{1}_G)=
\int_{\overline{N} \cap U} f(u) \; du       .
\]
\end{prop}

\begin{proof}
The remark preceding the proposition shows that this integral is convergent. Frobenius reciprocity
associates to any operator $J \in \Hom_G(i_P^G(\sigma),i_Q^G(\sigma))$ the unique operator
$K \in \Hom_M(r_Q^G\circ i_P^G(\sigma),\sigma)$ such that
\[
(\forall f \in i_P^G(E)), \quad      J(f) (\mathbf{1}_G)=K(\bar f)
\]
where $\bar f$ is the image of $f$ in $r_Q^G\circ i_P^G(E)$ (see the proof of the Frobenius reciprocity
theorem in \ref{Frob}). It therefore suffices to prove that $K_{Q|P} (\sigma)$ is the unique element of
$\Hom_M(r_Q^G\circ i_P^G(\sigma),\sigma)$ such that
\[
(\forall f \in  \widetilde F_{PQ}^{<\bar 1}(E)  ), \quad K_{Q|P}(\sigma)(\bar f)=
\int_{\overline{N} \cap U} f(u) \; du,
\]
Now an element of $\Hom_M(r_Q^G\circ i_P^G(\sigma),\sigma)$ is determined by its restriction to
$F_{PQ}^{<\bar 1}(E)$, since $\Hom_M\big((r_Q^G\circ i_P^G/F_{PQ}^{<\bar 1})(E),\sigma\big)=0$ by
$PQ$-regularity: the subquotients of this quotient have exponents which are not those of $\sigma$.
Uniqueness follows, and the above formula is indeed the one which gives $K_{Q|P}(\sigma)$.
\end{proof}

\subsection{Properties of intertwining operators}

For any semi-standard parabolic subgroups $P$ and $Q$ of $G$ having the same Levi factor $M$, set:
\[
d(P,Q)=| \Sigma(P)\cap \Sigma(\bar Q)   |
\]
This is the number of reduced roots of $A_M$ that are positive for $P$ and negative for $Q$. When
$M=M_\emptyset$, it is the length of the element $w$ of the Weyl group $W(A_\emptyset)$ such that
$w\cdot P=Q$; in general, such an element of $W(A_M)$ does not exist (for example for $G=GL_3$ and
$M=GL_2\times GL_1$).

\begin{prop}
Let $M$ be a semi-standard Levi subgroup of $G$ and let $(\sigma,E)$ be a smooth representation of $M$.
Let $P_1,P_2,P_3$ be three semi-standard parabolic subgroups of $G$ with Levi factor $M$.

$(i)$ Suppose that
\[
d(P_1,P_2)+d(P_2,P_3)=d(P_1,P_3), \quad d(P_1,P_2)=1.
\]
If $\sigma$ is $P_1P_3$-regular, then it is $P_2P_3$-regular and $P_1P_2$-regular. Moreover
\[
J_{P_3|P_2} (\sigma)\circ J_{P_2|P_1}(\sigma)= J_{P_3|P_1}(\sigma).
\]

$(ii)$ If $Q=LU$ is a semi-standard parabolic subgroup of $G$ containing $P_1$ and $P_2$, and if $\sigma$
is $P_1P_2$-regular, then with the natural identifications and the obvious notation $\sigma$ is
$(P_1\cap L)(P_2 \cap L)$-regular and
\[
J_{P_2|P_1} (\sigma)=i_Q^G(J_{(P_2\cap L)|(P_1 \cap L)}(\sigma)).
\]

$(iii)$ If $L$ is a standard Levi subgroup of $G$ containing $M$ such that $P_1\cap L=P_2\cap L$ and such
that $Q_1=P_1L=LU_1$, $Q_2=P_2L=LU_2$ are parabolic subgroups (with unipotent radical respectively $U_1$
and $U_2$) of $G$ and if $i_{P_1\cap L}^L(\sigma)$ is $(P_1L)(P_2L)$-regular, then $\sigma$ is
$P_1P_2$-regular and moreover, with the natural identifications and the obvious notation
\[
J_{P_2|P_1}(\sigma)= J_{Q_2|Q_1}( i_{P_1\cap L}^L(\sigma)).
\]
\end{prop}

\begin{proof}
We begin by proving the assertions concerning the regularity of the representations. Our hypothesis in
$(i)$ is that
\[
I_\sigma + \bigcap_{\bar w <_{P_1P_3} \bar 1 }  I_\sigma^{P_3,\bar
w}=\bbC[A_M].
\]

We want to show that
\begin{equation}\label{Dat1}
  I_\sigma + \bigcap_{\bar w <_{P_1P_2} \bar 1 }  I_\sigma^{P_2,\bar
   w}=\bbC[A_M] \text{ and }  I_\sigma + \bigcap_{\bar w <_{P_2P_3} \bar
   1 }   I_\sigma^{P_3,\bar w}=\bbC[A_M].
\end{equation}
By definition, $\bar w <_{P_2P_3}\bar 1$ implies\footnote{The reader should be careful not to confuse the
notation for closure and that for opposite parabolics. In the formulas that follow, it refers to the
former.}: $P_2wP_3 \subset \overline{ P_2P_3}$, thus
\[
P_1 w P_3 \subset   P_1 P_2wP_3 \subset P_1\overline{ P_2P_3} \subset \overline{P_1
P_2P_3}=\overline{ P_1P_3},
\]
the last equality being a consequence of $d(P_1,P_2)+d(P_2,P_3)=d(P_1,P_3)$. This shows that
$\bar w <_{P_2P_3}\bar 1$ implies $\bar w <_{P_1P_3}\bar 1$, which suffices to show the second assertion
of (\ref{Dat1}).

As above, we see that $\bar w <_{P_1P_2}\bar 1$ implies $\bar w <_{P_1P_3}\bar 1$. It remains to compare
the $I_\sigma^{P_3,\bar w}$ and the $I_\sigma^{P_2,\bar w}$, and to do this, we must compare the
representations $r_{M\cap w\cdot P_3}^{M}(\sigma)$ and $r_{M\cap w\cdot P_2}^{M} (\sigma)$ for an element
$\bar w <_{P_1P_2} \bar 1$ of $W_M\backslash W_G/W_M$. Since we assumed that $d(P_1,P_2)=1$,
$\Sigma(P_1)\cap \Sigma(\bar P_2)$ is a singleton, say a root $\alpha$. Let $N_\alpha$ be the
corresponding root subgroup, and let $U$ be the unipotent subgroup generated by all root subgroups
$U_\beta$, $\beta \in \Sigma(P_1)\cap \Sigma(P_2)$. Let $L$ be the Levi subgroup of $G$ generated by $M$,
$N_\alpha$ and $N_{-\alpha}$. We then have
\[
N_1=UN_\alpha, \quad N_2=UN_{-\alpha}, \quad L\cap P_1= MN_\alpha, \quad L\cap P_2= MN_{-\alpha}.
\]
Moreover, $Q=LU$ is a semi-standard parabolic subgroup of $G$ containing $P_1$ and $P_2$, and since
$d(P_1,P_2)+d(P_2,P_3)=d(P_1,P_3)$, we have $P_3 \cap L=P_2 \cap L$. From this, we easily deduce that
\footnote{The reader should be careful not to confuse the notation for closure and that for opposite
parabolics. In the formulas that follow, it refers to the former.}:
\[
\overline{P_1P_2}=   \overline{MN_\alpha N_{-\alpha}U}= \overline{L U}=Q.
\]
It follows that if $\bar w <_{P_1P_2} \bar 1$ then
\[
P_1 w P_2 \subset  \overline{P_1P_2}= Q=LU,
\]
whence
\[
MN_\alpha w MN_{-\alpha} = (P_1\cap L) w (P_2\cap L)\subset L,
\]
which entails, by the Bruhat decomposition of $L$, that $\bar w \in W_M\backslash W_L/W_M$, and we deduce
\[
M\cap w \cdot P_3=  M\cap w \cdot (P_3\cap L)=  M\cap
w \cdot (P_2\cap L)= M\cap w \cdot P_2.
\]
Then $I_\sigma^{P_3,\bar w}=I_\sigma^{P_2,\bar w}$, which implies the first equality in (\ref{Dat1}).

In $(ii)$, the hypothesis is that
\[
I_\sigma + \bigcap_{\bar w <_{P_1P_2} \bar 1 }  I_\sigma^{P_2,\bar
w}=\bbC[A_M]
\]
and we want to show that
\[
I_\sigma + \bigcap_{\bar w <_{(P_1\cap L)(P_2\cap L)} \bar 1 }
I_\sigma^{(P_2\cap L),\bar
w}=\bbC[A_M].
\]
We consider $W_M\backslash W_L/W_M$ as a subset of $W_M\backslash W_G/W_M$. For any
$\bar w \in W_M\backslash W_L/W_M$, we note that
\[
M\cap w \cdot P_2= M\cap w \cdot( P_2\cap L)
\]
so that $I_\sigma^{P_2,\bar w}=I_\sigma^{(P_2\cap L),\bar w}$. On the other hand, if
\[
(P_1 \cap L)\bar w (P_2 \cap L) \subset
\overline{(P_1 \cap L)\ (P_2 \cap L)},
\]
then
\[
P_1 \bar w P_2  =  P_1(P_1 \cap L) \bar w (P_2 \cap L)P_2
\subset P_1\overline{(P_1 \cap L)\ (P_2 \cap L)}P_2 \subset  \overline{P_1P_2},
\]
and thus $\bar w <_{(P_1\cap L)(P_2\cap L)} \bar 1$ implies that $\bar w <_{P_1P_2} \bar 1$. The desired
conclusion follows from these two facts.

In $(iii)$, the hypothesis is now that
\[
I_{i_{P_1\cap L}^L(\sigma)} + \bigcap_{\bar w <_{(P_1L)(P_2L)}
\bar 1 }   I_{i_{P_1\cap L}^L(\sigma)}^{(P_2L),\bar w}=\bbC[A_L]
\]
and we want to show that
\[
I_\sigma + \bigcap_{\bar w <_{P_1P_2} \bar 1 }  I_\sigma^{P_2,\bar
w}=\bbC[A_M].
\]
Note that if $\bar w \in W_M\backslash W_G/W_M$, with $\bar w <_{P_1P_2} \bar 1$, then
\[
Q_1\bar w Q_2= (P_1L)\bar w (P_2L)=LP_1 \bar w P_2 L\subset L \overline{P_1P_2}L
\subset \overline{(LP_1)(LP_2)}=\overline{Q_1Q_2}
\]
and thus $\bar w <_{(P_1L)(P_2L)} \bar 1$.

We need to compare the $I_{i_{P_1\cap L}^L(\sigma)}^{(P_2L),\bar w}$ and the $I_\sigma^{P_2,\bar  w}$.
Fix a $\bar w \in W_M\backslash W_G/W_M$, and apply the geometric lemma to the representation
$w \circ (r_{L \cap w \cdot(P_2 L)}^L\circ i_{P_1 \cap L}^L(\sigma)   )$. This representation admits a
filtration whose subquotients are of the form
\[
w\circ (i_{L \cap w \cdot L \cap v^{-1} \cdot P_1 }^{L \cap w \cdot L}
\circ v \circ r_{M \cap (vw)\cdot (P_2  L)}^{M}(\sigma)),
\]
where $\bar v$ runs through the set $W_M\backslash W_L/W_{L \cap w\cdot L}$. The annihilator of such a
representation in $\bbC[A_M]$ is the same as that of the representation
$$w\circ v\circ ( r_{M \cap (vw) \cdot (P_2  L)}^{M}(\sigma))=(vw) \circ ( r_{M \cap (vw) \cdot (P_2  L)}^{M}(\sigma))  $$
because the parabolic induction functor is faithful. This annihilator is contained in that of
$(vw)\circ ( r_{M \cap (vw) \cdot P_2 }^{M}(\sigma))$, equal by definition to
$I_\sigma^{P_2,\overline{vw}} \cap \bbC[A_L]$. This shows that
\[
I_{i_{P_1\cap L}^L(\sigma)}^{(P_2  L),\bar w} \subset
\bigcap_{\bar w' \mapsto \bar w}  (I_\sigma^{P_2,\overline{w}'}
\cap \bbC[A_L])
\]
where the $\overline{w}'$ describe the fiber above $\bar w$ of the projection
\[
W_M\backslash W_G/W_M  \rightarrow  W_L\backslash W_G/W_L.
\]
We deduce that
\[
I_\sigma \cap \bbC[A_L]+  \bigcap_{\bar w <_{P_1P_2} \bar 1 }  (I_\sigma^{P_2,\bar
w}\cap \bbC[A_L])=\bbC[A_L],
\]
which implies the desired equality.

We show the equalities between intertwining operators of the proposition. Let us resume point $(i)$. The
hypothesis that $d(P_1,P_3)=d(P_1,P_2)+d(P_2,P_3)$ can be translated as\footnote{The reader should be
careful not to confuse the notation for closure and that for opposite parabolics. In the formulas that
follow, it refers to the latter.}:
\[
\overline{N}_1 \cap N_3= (\overline{N}_1 \cap N_2)(\overline{N}_2 \cap N_3).
\]

On the other hand, we have the following characterization of $\widetilde F_{P_1P_3}^{<\bar 1}(E)$:
\[
\widetilde F_{P_1P_3}^{< \bar 1}(E) = \{f \in i_{P_1}^G(E) \, | \, \mathrm{Supp}
f  \cap (\overline{N}_1 \cap N_3) \text{ is compact }  \},
\]
because the natural map
\[
(\overline{N}_1 \cap N_3)\rightarrow  P_1\backslash P_1P_3
\]
is a homeomorphism. It follows that if $f \in \widetilde F_{P_1P_3}^{<\bar 1}(E)$, we have

$\bullet$ For all $u \in \overline{N}_2 \cap N_3$, the function $r(u)\cdot f \in i_{P_1}^G(E)$ has
compact support on $\overline{N}_1 \cap N_2$, and thus is in $\widetilde F_{P_1P_2}^{<\bar 1}(E)$.

$\bullet$ For all $u \in \overline{N}_2 \cap N_3$, we have
\[
J_{P_2|P_1}(\sigma)(f)(u)=   J_{P_2|P_1}(\sigma)(r(u)\cdot
f)(\mathbf{1}_G)= \int_{\overline{N}_1 \cap N_2} f(vu)\, dv.
\]
This shows that the function $J_{P_2|P_1}(\sigma)(f) \in i_{P_2}^G(E)$ is compactly supported on
$\overline{N}_2\cap N_3$, and thus is in $\widetilde F_{P_2P_3}^{< \bar 1}(E)$.

It follows from these two points that
\begin{align*}
   J_{P_3|P_2}(\sigma)\circ J_{P_2|P_1}(\sigma)(f)(\mathbf{1}_G)&=\int_{\overline{N}_2 \cap N_3}
\int_{\overline{N}_1 \cap N_2}f(vu) \, dv  \, du\\
&= \int_{\overline{N}_1\cap N_3} f(u)\, du = J_{P_3|P_1}(\sigma)(f)(\mathbf{1}_G)
\end{align*}

According to Proposition \ref{INTOP}, this equality is in fact valid for any $f \in i_{P_1}^G(E)$, and
this therefore proves that
\[
J_{P_3|P_2}(\sigma)\circ J_{P_2|P_1}(\sigma)=J_{P_3|P_1}(\sigma).
\]

We now turn to point $(ii)$. First fix, for $i=1,2$, a natural isomorphism between $i_{P_i}^G(E)$ and
$ i_Q^G \circ i_{P_i\cap L}^L(E)$:
\[
f \mapsto \tilde f
\]
where
\[
\tilde f(g)= f_g\,\colon l\in L \mapsto \delta_U(l)^{1/2}f(lg).
\]
We easily verify that for all $g \in G$, $f_g \in i_{P_i \cap L}^L (E)$ (note that $U \subset N_i$), that
$\tilde f \in  i_Q^G \circ i_{P_i\cap L}^L(E)$ and that the inverse isomorphism is given by
\[
\tilde f \mapsto [ g \mapsto \tilde f(g)(\mathbf{1}_{L})].
\]

If $f \in  i_{P_1}^G(E)$, we have for all $g \in G$, by definition of the induced operator
$i_Q^G ( J_{P_2\cap L|P_1\cap  L}(\sigma))$,
\[
i_Q^G ( J_{P_2\cap L|P_1\cap L}(\sigma))(\tilde f)= J_{P_2\cap L|P_1\cap
L}(\sigma)\circ \tilde f.
\]
Whence
\begin{equation}\label{alacon}
 \left( \left( i_Q^G ( J_{P_2\cap L|P_1\cap L}(\sigma))(\tilde f) \right)(g) \right) (\mathbf{1}_{L}) =
\left(J_{P_2\cap L|P_1\cap  L}(\sigma)(\tilde f(g))\right) (\mathbf{1}_{L}).
\end{equation}
Let $f \in \widetilde F_{P_1P_2}^{<\bar 1}(E)$. Since $\mathrm{Supp}(f)\cap (\overline{N}_1 \cap N_2)$ is
compact, it is also the case for
\[
\mathrm{Supp}(\tilde f(\mathbf{1}_G)) \cap ((\overline{N}_1\cap L) \cap (N_2\cap  L)),
\]
since $((\overline{N}_1\cap L) \cap (N_2 \cap L))= \overline{N}_1 \cap N_2$. The characterization of the
spaces $\widetilde F^{<\bar 1}$ used above then shows that
$\tilde f(\mathbf{1}_G) \in \widetilde F_{(P_1\cap L)(P_2\cap L)}^{<\bar 1}(E)$ and we can then evaluate
the expression (\ref{alacon}) at $g=\mathbf{1}_G$.
\[
\left(  J_{P_2\cap L|P_1\cap L}(\sigma)(\tilde f (\mathbf{1}_G)) \right) (\mathbf{1}_L)
=\int_{\overline{N}_1 \cap N_2} f(u)\, du= J_{P_2|P_1}(\sigma)(f)  (\mathbf{1}_G).
\]
This suffices to show that $i_Q^G(J_{P_2\cap L|P_1\cap L}(\sigma))=J_{P_2|P_1}(\sigma)$, according to
Proposition \ref{INTOP}.

Finally, let us finish proving $(iii)$. We have $N_i=(N_i \cap L)(N_i \cap U_i)$, $i=1,2$, which implies
that the natural inclusion of $\bar U_1 \cap U_2$ into $\overline{N}_1 \cap N_2$ is a homeomorphism.
Identify $i_{P_1}^G(E)$ and $i_{P_1L}^G \circ i_{P_1\cap L}^L(E)$ as in the proof of $(ii)$.

If $\tilde f \in i_{P_1L}^G \circ i_{P_1\cap L}^L(E)$, the function $f \in i_{P_1}^G(E) $ which
corresponds to it satisfies
\[
\delta_{U_1}^{1/2}(l)f(lg)=\tilde f(g)(l), \quad (l \in L), \, (g \in G).
\]
In particular, for all $l \in L$, $\mathrm{Supp}(r(l)\cdot f)\subset \mathrm{Supp}(\tilde f)$ and thus
for all $l \in L$, $r(l)\cdot f \in \widetilde F_{P_1P_2}^{<\bar 1}(E)$ as soon as
$\tilde f \in \widetilde F_{(P_1L)(P_2L)}^{<\bar 1}(i_{P_1\cap L}^L(E))$. We can then calculate, noting
at the last step that the modular characters cancel out,
\begin{align*}
&\widetilde{ J_{P_2|P_1}(\sigma)(f)}(\mathbf{1}_G)(l)= \delta_{U_2}^{1/2}(l) J_{P_2|P_1}(\sigma)(f)(l)\\
&=\delta_{U_2}^{1/2}(l) J_{P_2|P_1}(\sigma)(r(l)\cdot f)(\mathbf{1}_G)= \delta_{U_2}^{1/2}(l) \int_{\overline{N}_1 \cap N_2} f(ul)\, du \\
&=\delta_{U_2}^{1/2}(l) \int_{\bar U_1 \cap U_2} f(ul)\, du = \delta_{U_2}^{1/2}(l)\delta_{\bar U_1\cap U_2}(l)^{-1}  \int_{\bar U_1 \cap U_2} f(lu)\, du \\
&=\delta_{U_2}^{1/2}(l) \delta_{\bar U_1\cap U_2}(l)^{-1}\delta_{U_1}^{-1/2}(l) \int_{\bar U_1 \cap U_2} \tilde f(u)(l)\, du \\
&= J_{P_2L|P_1L}(i_{P_1 \cap L}^L(\sigma))(\tilde f)(\mathbf{1}_G)(l).
\end{align*}
Once again, the conclusion follows from Proposition \ref{INTOP}.
\end{proof}

\section{Langlands classification}

\subsection{Langlands triplets}\label{langtrip}

Langlands triplets give a parameterization of $\mathbf{Irr}(G)$.

\begin{defi}
A Langlands triplet consists of:

- a standard parabolic subgroup $P=MN$,

- an irreducible tempered representation $(\sigma, E)$ of $M$,

- an unramified character $\psi \in \caX(M)$ such that $\Re(\psi) \in {}^G_P[\fra_M^*]^+$
\end{defi}

\enlargethispage{2\baselineskip}
\begin{lemme}
If $(P=MN,(\sigma, E),\psi)$ is a Langlands triplet, then the representation $\sigma \otimes \psi$ is
$P\bar P$-regular. In particular, there exists a non-zero intertwining operator, unique up to a scalar
factor:
\[
J_{\bar P|P}(\sigma\psi)\colon i_P^G(\sigma \otimes \psi)\rightarrow
i_{\bar P}^G(\sigma \otimes \psi).
\]
\end{lemme}

\begin{proof}
By Lemma \ref{PQreg}, we must show that if $\bar w <_{P\bar P}\bar 1$,
$\bar w \in W_M\backslash W_G/W_M$, then
\[
\mathrm{Exp} (A_M,\sigma \otimes \psi)  \cap  \mathrm{Exp}
(A_M,w\circ r_{M\cap w\cdot \bar P}^M(\sigma \otimes \psi))=\emptyset.
\]

Now, the elements of $\mathrm{Exp} (A_M,w\circ r_{M\cap w\cdot \bar P}^M(\sigma \otimes \psi))$ are of
the form
\[
(w^{-1}\cdot \chi)_{|A_M} (w^{-1}\cdot \psi)_{|A_M}.
\]

Let us look more closely at these two terms. On the one hand
$$\chi \in  \mathrm{Exp} (A_{M\cap w\cdot M},  r_{M\cap w\cdot \bar P}^M(\sigma)).$$ Since $\sigma$ is
tempered, we have
\[
\Re (\chi)\in {}^+\overline{[\fra_{M\cap w\cdot M}^*]}_{M\cap w\cdot \bar P}^M.
\]
The proof of Lemma \ref{temp2} then shows that
\begin{equation}\label{L0}
 \Re (w^{-1}\cdot \chi_{|A_M}) \in
{}^+\overline{[\fra_{M}^*]}_{\bar P}^G= {}^-\overline{[\fra_{M}^*]}_{P}^G  .
\end{equation}

On the other hand, $\psi \in \caX(M)$, thus $w^{-1}\cdot \psi \in \caX(w^{-1}\cdot M)$. Since
\[
A_{M} \subset A_{M\cap  w^{-1}\cdot M} \subset
A_\emptyset \subset M\cap w^{-1} \cdot M \subset w^{-1}\cdot M,
\]
$(w^{-1}\cdot \psi)_{|A_M} \in \caX(A_M)$, and $\Re((w^{-1}\cdot \psi)_{|A_M})$ is obtained as follows:
from $\mu:=\Re(\psi)\in \fra_M^*$: we consider $\mu$ as an element of
$\fra_{M\cap w \cdot M}^*=\fra_M^*\oplus (\fra_{M\cap w \cdot M}^M)^*$, and thus
\[
w^{-1}\cdot \mu =  \Re(w^{-1}\cdot \psi) \in   \fra_{M\cap w^{-1}
\cdot M}^*=\fra_M^*\oplus (\fra_{M\cap w^{-1}   \cdot M}^M)^*.
\]
We then project $w^{-1}\cdot \mu$ onto $\fra_M^*$, i.e.,
\[
\Re((w^{-1}\cdot \psi)_{|A_M})= (w^{-1}\cdot \mu)_{|\fra_M}.
\]

If $\bar w \neq \bar 1$, we now show that we have
\begin{equation}\label{L1}
\Re (w^{-1}\cdot \psi_{|A_M}) \in \Re (\psi)+
({}^-\overline{[\fra_{M}^*]}_{P}^G \setminus \{0\}) .
\end{equation}
To do this, we first show that
\begin{equation}\label{L2}
(\forall \mu \in  {}^G_{P_\emptyset}\overline{[\fra_\emptyset^*]}^+), \; (\forall
w \in W_G), \quad  \mu- w^{-1}\cdot \mu \in  {}^+\overline{[\fra_{\emptyset}^*]}_{P_\emptyset}^G.
\end{equation}
We proceed by induction on the length of $w$. Let $\alpha_1,\ldots, \alpha_l$ denote the elements of
$\Delta_\emptyset$, and $\beta_1,\ldots ,\beta_l$ the dual basis. We want to establish that for all
$i=1, \ldots ,l$,
\[
\langle  \mu- w^{-1} \cdot \mu , \beta_i \rangle \geq 0.
\]
Now,
\[
\langle  \mu- w^{-1}\cdot \mu , \beta_i \rangle = \langle  \mu, \beta_i-
w \cdot \beta_i \rangle .
\]
If $l(w)=1$, $w=s_{\alpha_j}$ for a certain root $\alpha_j \in \Delta_\emptyset$ and in this case,
\[
\langle  \mu, \beta_i-
w \cdot \beta_i \rangle  =  \langle  \mu, \beta_i-
s_{\alpha_j} \cdot \beta_i \rangle= \begin{cases}   \langle  \mu,
\check\alpha_i \rangle \; \text {if } i=j\\ 0  \; \text {if } i\neq j \end{cases}.
\]
(we have $s_{\alpha_j}\cdot\beta_i=\beta_i-\langle\alpha_j,\beta_i\rangle\check\alpha_j$). This quantity
is therefore non-negative, since $\mu$ is in the closed positive chamber.

If $l(w)>1$, we write $w=w's_{\alpha_j}$ with $l(w')=l(w)-1$, and we then have
\begin{align*}
&\langle  \mu, \beta_i-w \cdot \beta_i \rangle= \langle  \mu, \beta_i- w' s_{\alpha_j} \cdot \beta_i \rangle\\
=& \langle  \mu, \beta_i- w'\cdot \beta_i \rangle +   \langle \mu,   w'\cdot \beta_i -w' s_{\alpha_j} \cdot \beta_i
\rangle\\
=& \langle  \mu, \beta_i- w'\cdot \beta_i \rangle + \langle {w'}^{-1}\cdot \mu,  \beta_i -  s_{\alpha_j} \cdot \beta_i \rangle\\
=& \langle  \mu- {w'}^{-1}\cdot \mu , \beta_i \rangle +\langle {w'}^{-1}\cdot \mu,  \beta_i -  s_{\alpha_j} \cdot \beta_i \rangle.
\end{align*}
By the induction hypothesis, the first term is $\geq 0$. For the second, if $i\neq j$, we obtain $0$, and
if $i=j$, we obtain
\[
\langle {w'}^{-1}\cdot \mu, \check\alpha_i \rangle=  \langle \mu, w' \cdot \check\alpha_i \rangle.
\]
But this is $\geq 0$, since $\mu$ is in the closed positive chamber and $w'\cdot\check\alpha_i$ is a
positive coroot, because $l(w')=l(w's_{\alpha_i} )-1$ is equivalent to
\[
w'\cdot \alpha_i \in  {}^+\overline{[\fra_{\emptyset}^*]}_{P_\emptyset}^G.
\]
This completes the proof of (\ref{L2}). Let us deduce (\ref{L1}) from it, by setting $\Re(\psi)=\mu$. By
hypothesis, $\mu=\Re(\psi) \in {}^G_P[\fra_M^*]^+$, and thus
$\mu\in {}^G_{P_\emptyset}\overline{[\fra_\emptyset^*]}^+$, and from the above
\[
\mu- w^{-1}\cdot \mu \in
{}^+\overline{[\fra_{\emptyset}^*]}_{P_\emptyset}^G .
\]
We deduce that
\[
\mu- (w^{-1}\cdot \mu)_{|\fra_M^*} \in
{}^+\overline{[\fra_M^*]}_{P}^G .
\]
because the projection of $\fra_\emptyset^*$ onto $\fra_M^*$ maps
${}^+\overline{[\fra_{\emptyset}^*]}_{P_\emptyset}^G$ to ${}^+\overline{[\fra_M^*]}_{P}^G$. It remains to
show that $\mu- w^{-1}\cdot \mu \neq 0$ for $\bar w \neq \bar 1 \in W_M\backslash W_G/W_M$. Suppose the
contrary: $w^{-1}\cdot \mu=\mu$. Let us choose the set $W^{M,M}$ defined in \ref{WeylgroupsS} as a system
of representatives of the double cosets $W_M\backslash W_G/W_M$. Then
$w(M \cap P_\emptyset)\subset P_\emptyset$ and thus if
$\alpha \in \Delta^M(M\cap P_\emptyset)=\Delta_\emptyset^M$, we have
$w\cdot \alpha \in {}^+\overline{[\fra_{\emptyset}^*]}_{P_\emptyset}^G$. If
$\alpha \in \Delta_\emptyset^G \setminus \Delta^M_\emptyset$, we have
\[
( \mu, w\cdot \alpha )=( w^{-1}\cdot \mu,  \alpha)= (\mu,  \alpha )>0.
\]
By Lemma \ref{fait24}, we then obtain
\[
w \cdot \alpha \in {}^+\overline{[\fra_{\emptyset}^*]}_{P_\emptyset}^G.
\]
Finally, we have shown that for all $\alpha \in \Delta^G_\emptyset$,
$ w\cdot \alpha \in {}^+\overline{[\fra_{\emptyset}^*]}_{P_\emptyset}^G$, which implies that $l(w)=0$.
This contradicts the hypothesis $w\neq 1$. Finally, $(w^{-1}\cdot\mu)_{|\fra_M}\neq\mu$: otherwise
$\nu=\mu-w^{-1}\cdot\mu$ would be orthogonal to $\fra_M^*$, hence $(\mu,\nu)=0$, and since
$|w^{-1}\cdot\mu|=|\mu|$, we would have $|\nu|^2=2(\mu,\mu)-2(\mu,w^{-1}\cdot\mu)=2(\mu,\nu)=0$, which
contradicts the above. Thus
$\mu-(w^{-1}\cdot\mu)_{|\fra_M}\in{}^+\overline{[\fra_M^*]}_{P}^G\setminus\{0\}$, which establishes
(\ref{L1}).

We can now complete the proof of the lemma. From (\ref{L0}) and (\ref{L1}),
\[
\Re(w^{-1}\cdot \chi_{|A_M})+\Re(w^{-1}\cdot \psi_{|A_M}) \in \Re (\psi)
+({}^-\overline{[\fra_{M}^*]}_{P}^G \setminus \{0\}).
\]
Now an element of $\mathrm{Exp}(A_M,\sigma\otimes \psi)$ is written $\chi_1\psi$, with
$\chi_1 \in \mathrm{Exp}(A_M,\sigma)$. The representation $\sigma$ being tempered, $\Re(\chi_1)=0$, thus
\[
\Re(\chi_1\psi)=\Re(\psi) \in {}^G_P[\fra_M^*]^+.
\]
This shows that if $\bar w \neq \bar 1$,
\[
\mathrm{Exp} (A_M,\sigma \otimes \psi)  \cap  \mathrm{Exp}
(A_M,w\circ r_{M\cap w\cdot \bar P}^M(\sigma \otimes \psi))=\emptyset.
\]

Since $\bar 1$ is the maximal element for the Bruhat order $\leq_{\bar P P}$, we then see that by
Frobenius reciprocity $\Hom_G(i_P^G(\sigma\otimes \psi), i_{\bar P}^G(\sigma\otimes \psi))$ is
$1$-dimensional, generated by $J_{\bar P|P }(\sigma \psi)$.
\end{proof}

\subsection{The classification theorem}\label{classthm}

We can now prove the main theorem of this chapter.

\enlargethispage{2\baselineskip}
\begin{thm}
$(i)$ Let $(P=MN,(\sigma,E),\psi)$ be a Langlands triplet. Then the representation
$i_P^G(\sigma\otimes \psi,E)$ admits a unique irreducible quotient. We denote it by $J(P,\sigma,\psi)$.

$(ii)$ If $J(P_1,\sigma_1,\psi_1)= J(P_2,\sigma_2,\psi_2) $, then $M_1=M_2$, $\Re(\psi_2\psi_1^{-1})=0$
and $\sigma_1\psi_1 \simeq \sigma_2 \psi_2$.

$(iii)$ Let $(\pi,V)\in \mathbf{Irr}(G)$. Then there exists a Langlands triplet $(P=MN,(\sigma,E),\psi)$
such that $(\pi,V)\simeq J(P,\sigma,\psi)$.
\end{thm}

\begin{proof}
$(i)$ We have seen that $\sigma \otimes \psi$ is $P\bar P$-regular and that
$$\Hom_G(i_P^G(\sigma\otimes \psi), i_{\bar P}^G(\sigma\otimes \psi))$$ is $1$-dimensional, generated by
$J_{\bar P|P }(\sigma \psi)$. Let $(\pi,V)$ be a quotient of $i_P^G(\sigma\otimes \psi)$. By the second
adjunction theorem,
\[
\Hom_G(i_P^G(\sigma\otimes \psi), \pi)\simeq \Hom_M(\sigma\otimes
\psi,r_{\bar P}^G(\pi))
\]
and we obtain a non-zero morphism
\[
\sigma \otimes \psi  \rightarrow r_{\bar P}^G(\pi)
\]
which is injective since $\sigma \otimes \psi$ is irreducible. Since $\sigma \otimes \psi$ is
$P\bar P$-regular, this embedding admits a retraction
\[
r_{\sigma\psi}\colon r_{\bar P}^G(\pi)  \rightarrow  \sigma \otimes \psi
\]
and by Frobenius reciprocity, we obtain a non-zero morphism in
\[
\Hom_G(\pi, i_{\bar P}^G(\sigma\otimes \psi)).
\]
If $(\pi,V)$ is irreducible, this morphism is an injection. The composition
\[
i_P^G(\sigma\otimes \psi) \rightarrow \pi  \rightarrow i_{\bar
P}^G(\sigma\otimes \psi)
\]
is a non-zero element of $\Hom_G(i_P^G(\sigma\otimes \psi), i_{\bar P}^G(\sigma\otimes \psi))$, and is
therefore equal to $J_{\bar P|P }(\sigma \psi)$ up to a scalar factor. This shows that $(\pi,V)$ is the
unique irreducible quotient of $i_P^G(\sigma\otimes \psi)$.

Moreover, we see from the above that $(\pi,V)$ is characterized as an irreducible subquotient of
$ i_P^G(\sigma\otimes \psi)$ by the property
\[
\Re(\psi) \in \Re( \mathrm{Exp}(A_M, r_{\bar P}^G(\pi))) .
\]

Let us now prove $(iii)$. Let $(\pi,V)\in \mathbf{Irr}(G)$. Consider the set
\[
\Re(\mathrm{Exp}(V)):=  \bigcup_M \Re(\mathrm{Exp}(A_M,  r_{\bar P}^G(\pi))),
\]
where $M$ runs through the standard Levi subgroups of $G$.

Let $\mu \in   \Re(\mathrm{Exp}(V))$, say $\mu \in \Re(\mathrm{Exp}(A_M, r_{\bar P}^G(\pi)))$ for a
certain standard parabolic subgroup $P=MN$ of $G$, and let
\[
\mu =\mu^++\mu^-+\mu_G
\]
the decomposition of $\mu$ given by the Langlands combinatorial lemma \ref{comblangl}, with
$\mu^+ \in {}^G_Q [(\fra_L^G)^*]^+$ and $\mu^- \in  p_M^L ({}^- \overline{[\fra_M^*]_P^G}  )$ where
$Q=LU$ is a certain standard parabolic subgroup of $G$ containing $M$.

Choose $\mu$ such that the norm of $\mu^+$ is maximal. Since $\mu^-_{|A_L}$ is trivial, we have by Lemma
\ref{exposants} and the transitivity of parabolic restriction functors,
\[
\mu_{|A_L}=\mu^+ +\mu_G \in \Re(\mathrm{Exp}(A_L, r_{\bar Q}^G (V))).
\]
We can therefore find an irreducible subquotient $(\rho,W)$ of $r_{\bar Q}^G (\pi,V)$ whose central
character $\chi_\rho$ satisfies
\[
\Re(\chi_\rho)= \mu^+ +\mu_G .
\]
By decomposing $ r_{\bar Q}^G (V)$ according to the central character, we can even assume that $(\rho,W)$
is a subrepresentation of $r_{\bar Q}^G (\pi,V)$.

Let $\psi \in \caX(L)$ such that $\Re(\psi) =  \mu^+ +\mu_G$. Then set $\sigma=\rho \otimes \psi^{-1}$.
We will now show that $\sigma$ is tempered. To do this, consider a standard parabolic subgroup
$P_1=M_1N_1$ of $G$ contained in $Q$ and let us examine more closely the exponents of
\[
r_{\bar P_1 \cap L}^L(\sigma)= r_{\bar P_1 \cap L}^L(\rho\otimes \psi^{-1}).
\]
Let $\lambda$ be the real part of such an exponent. We see that $\lambda \in (\fra_{M_1}^L)^*$ because
indeed $\Re(\chi_\rho\psi^{-1})=0$ and thus by $(ii)$ of Lemma \ref{exposants}, $\lambda_{|\fra_L}=0$. We
deduce that $\lambda$ and $\mu^+$ are orthogonal, and thus
\[
|\lambda|^2+|\mu^{+}|^2=|\lambda+ \mu^+|^2=|(\lambda+ \mu^+)^+|^2+|(\lambda+ \mu^+)^-|^2 .
\]
Since $\rho \hookrightarrow r_{\bar Q}^G(\pi)$, we have
\[
r_{\bar P_1\cap L }^L (\rho) \hookrightarrow  r_{\bar P_1\cap L }^L
\circ
r_{\bar Q}^G(\pi)=  r_{\bar P_1}^G(\pi),
\]
and we then see that
$\lambda+\Re(\psi)=\lambda+ \mu^++\mu_G \in \Re(\mathrm{Exp}(A_{M_1}, r_{\bar P_1}^G(V)))$. The choice of
$\mu$ implies that
\[
|\mu^{+}|\geq |(\lambda+ \mu^+)^+|
\]
whence
\[
|(\lambda+ \mu^+)^-|\geq |\lambda|.
\]
On the other hand, since
\[
{}^G_{P_\emptyset}\overline{ [\fra_\emptyset^*] }^+ \subset
-\mu^++ {}^G_{P_\emptyset}\overline{  [\fra_\emptyset^*] }^+ ,
\]
we have
\[
|  (\lambda+ \mu^+)^- |= \mathrm{Dist}( \lambda+ \mu^+,
{}^G_{P_\emptyset}\overline{  [\fra_\emptyset^*] }^+  )\leq \mathrm{Dist}( \lambda,
{}^G_{P_\emptyset}\overline{  [\fra_\emptyset^*] }^+  )  =|\lambda^-|.
\]
Finally, we obtain
\[
|\lambda^-|\leq |\lambda|\leq |(\lambda+ \mu^+)^-|\leq  |\lambda^-|,
\]
whence
$\lambda=\lambda^-\in {}^-\overline{[\fra_{M_1}^*]}_{P_1}^G = {}^+\overline{[\fra_{M_1}^*]}_{\bar P_1}^G$.
By changing $\Delta_\emptyset^L$ to its opposite, so that the $\bar P_1 \cap L$ are the standard
parabolic subgroups of $L$, we see that $\sigma$ satisfies the desired conditions to be tempered.
Moreover, $\Re(\psi)= \mu^+ +\mu_G \in {}_Q^G[\fra_L^*]^+$.

We thus have $\rho=\sigma \otimes \psi \hookrightarrow r_{\bar Q}^G(\pi)$, whence we obtain by the second
adjunction theorem a non-zero morphism in $\Hom_G(i_Q^G(\sigma\otimes \psi),\pi)$, which shows that $\pi$
is a quotient of $i_Q^G(\sigma\otimes \psi)$. Thus $J(Q,\sigma,\psi)=\pi$.

It remains to prove $(ii)$, i.e., uniqueness. Suppose that two Langlands triplets
\[
(P_i=M_iN_i,  (\sigma_i,E_i), \psi_i),\quad i=1,2,
\]
give the same Langlands quotient $(\pi,V)$:

\begin{equation}
\begin{CD}
i_{P_1}^G(\sigma_1\psi_1) @>>>  \pi   @>>>   i_{\bar P_1}^G(\sigma_1\psi_1)    \\
&  &                        @VVV   &   \\
i_{P_2}^G(\sigma_2\psi_2) @>>>  \pi   @>>>   i_{\bar P_2}^G(\sigma_2\psi_2).    \\
\end{CD}
\end{equation}
The vertical arrow is the identity of $(\pi,V)$. We deduce the existence of a non-zero morphism in
$\Hom_G(i_{P_1}^G(\sigma_1\psi_1),  i_{\bar P_2}^G(\sigma_2\psi_2))$, and by Frobenius reciprocity, of a
non-zero morphism in $\Hom_{M_2}( r_{\bar   P_2}^G \circ i_{P_1}^G(\sigma_1\psi_1),   \sigma_2\psi_2)$.

By the geometric lemma, there exists $\bar w \in W_{M_1}\backslash W_G /W_{M_2}$ such that
\[
\Hom_{M_2}( i_{M_2 \cap w^{-1} \cdot P_1 }^{M_2} \circ  w \circ   r_{M_1
\cap w \cdot   \bar P_2}^{M_1}(\sigma_1\psi_1),   \sigma_2\psi_2)
\neq \{ 0 \},
\]
whence, by the second adjunction theorem,
\[
\Hom_{M_2\cap w^{-1}\cdot M_1}(  w \circ   r_{M_1
\cap w \cdot   \bar P_2}^{M_1}(\sigma_1\psi_1),   r_{M_2 \cap
w^{-1} \cdot  \bar P_1}^{M_2} (\sigma_2\psi_2))\neq \{ 0 \}.
\]

Set $\mu_i=\Re(\psi_i) \in {}_{P_i}^G[\fra_{M_i}^*]^+$. The above shows that there exist
\[
\lambda_1 \in \mathrm{Exp}(A_{M_1   \cap   w \cdot   M_2},   r_{M_1   \cap w \cdot   \bar  P_2}^{M_1}(\sigma_1))
\]
and
$\lambda_2 \in \mathrm{Exp}(A_{M_2   \cap  w^{-1} \cdot   M_1},   r_{M_2  \cap w^{-1} \cdot \bar  P_1}^{M_2}(\sigma_2))$
such that
\[
w^{-1}\cdot (\lambda_1+\mu_1)=\lambda_2+\mu_2
\]
Since $\sigma_1$ and $\sigma_2$ are tempered,
\[
\lambda_1 \in {}^+\overline{[\fra_{M_1\cap w\cdot M_2}^*]}_{M_1   \cap w \cdot   \bar  P_2}^{M_1}
\subset {}^+\overline{[\fra_{M_1\cap w\cdot     M_2}^*]}_{\bar R_1}^{G}= {}^-\overline{[\fra_{M_1\cap w\cdot M_2}^*]}_{R_1}^{G}
\]
\[
\lambda_2 \in {}^+\overline{[\fra_{M_2\cap w^{-1}\cdot M_1}^*]}_{M_2   \cap w^{-1} \cdot   \bar  P_1}^{M_2}
\subset {}^+\overline{[\fra_{M_2\cap w^{-1}\cdot   M_1}^*]}_{\bar R_2}^{G}= {}^-\overline{[\fra_{M_2\cap w^{-1}\cdot    M_1}^*]}_{R_2}^{G},
\]
where $R_1$ (resp. $R_2$) is the standard parabolic subgroup of $G$ with Levi factor $M_1\cap w\cdot M_2$
(resp. $M_2\cap w^{-1}\cdot M_1$). Recall (Lemma \ref{WeylgroupsS}) that $M_1   \cap w \cdot P_2$ (resp.
$M_2   \cap w^{-1} \cdot  P_1$) is the standard parabolic subgroup of $M_1$ (resp. of $M_2$) with Levi
factor $M_1\cap w\cdot M_2$ (resp. $M_2\cap w^{-1}\cdot M_1$), so that
\[
M_1   \cap w \cdot P_2 \subset R_1, \quad M_2   \cap w^{-1} \cdot P_1 \subset R_2.
\]

We also have
\[
w^{-1}\cdot \lambda_1 \in   {}^+\overline{[\fra_{M_2\cap
w^{-1}\cdot M_1}^*]}_{w^{-1}\cdot M_1 \cap   \bar
P_2}^{w^{-1}\cdot M_1} \subset {}^+\overline{[\fra_{M_2\cap w^{-1}\cdot
M_1}^*]}_{\bar R_2}^{G}= {}^-\overline{[\fra_{M_2\cap w^{-1}\cdot
M_1}^*]}_{R_2}^{G}.
\]

Note that
\[
\mu_1 \in  {}_{P_1}^G[\fra_{M_1}^*]^+ \subset
{}_{R_1}^G  \overline{ [\fra_{M_1\cap w\cdot M_2}^*]}^+, \quad \mu_2 \in  {}_{P_2}^G [\fra_{M_2}^*]^+ \subset
{}_{R_2}^G \overline{ [\fra_{M_2\cap w^{-1}\cdot M_1}^*]}^+.
\]
And finally
\[
w^{-1}\cdot   \mu_1 \in  {}_{w^{-1}\cdot P_1}^G[\fra_{w^{-1}\cdot
M_1}^*]^+ \subset {}_{R_2}^G \overline{ [\fra_{M_2\cap w^{-1}\cdot M_1}^*]}^+.
\]

Set
\[
\mu= \lambda_2+\mu_2= w^{-1}\cdot \lambda_1+ w^{-1}\cdot
\mu_1.
\]

The existence and uniqueness in the statement of the Langlands combinatorial lemma \ref{comblangl} then
give us
\[
\mu^-= \lambda_2 =w^{-1}\cdot \lambda_1,
\]
\[
\mu^+=\mu_2^+= (w^{-1}\cdot \mu_1)^+
\]
\[
\mu_G=(\mu_2)_G = (w^{-1}\cdot \mu_1)_G.
\]

Let us now show that this implies that $w=1$, by an argument already used in the proof of Lemma
\ref{langtrip}. Indeed, the choice of the system of representatives $W^{M_2,M_1}$ gives that
$w(M_2 \cap P_\emptyset)\subset P_\emptyset$ and thus if
$\alpha \in \Delta^{M_2}(M_2\cap P_\emptyset)=\Delta_\emptyset^{M_2}$, we have
$w\cdot \alpha \in {}^+\overline{[\fra_{\emptyset}^*]}_{P_\emptyset}^G$. If
$\alpha \in \Delta_\emptyset^G \setminus \Delta^{M_2}(M_2 \cap P_\emptyset)$, we have
\[
\langle \mu_1, w\cdot \alpha \rangle=\langle w^{-1}\cdot \mu_1,  \alpha
\rangle = \langle \mu_2,  \alpha \rangle>0.
\]
By Lemma \ref{fait24}, we then obtain
\[
w \cdot \alpha \in {}^+\overline{[\fra_{\emptyset}^*]}_{P_\emptyset}^G.
\]
Finally, we have shown that for all $\alpha \in \Delta^G_\emptyset$,
$ w\cdot \alpha \in {}^+\overline{[\fra_{\emptyset}^*]}_{P_\emptyset}^G$, which implies that $l(w)=0$.

Thus $\mu_1=\mu_2$. Since $\mu_1 \in   {}_{P_1}^G[\fra_{M_1}^*]^+$, the set of
$\alpha\in\Delta_\emptyset$ such that $\langle\mu_1,\check\alpha\rangle>0$ is
$\Delta_\emptyset\setminus\Delta_\emptyset^{M_1}$ (for $\alpha\in\Delta^{M_1}_\emptyset$, we have
$\langle\mu_1,\check\alpha\rangle=0$ since $\mu_1\in\fra_{M_1}^*$), and similarly for $\mu_2$. This shows
that $M_1=M_2$. The non-zero morphism in
\[
\Hom_{M_2\cap w^{-1}\cdot M_1}(  w \circ   r_{M_1   \cap w \cdot   \bar P_2}^{M_1}(\sigma_1\psi_1),
r_{M_2 \cap   w^{-1} \cdot  \bar P_1}^{M_2} (\sigma_2\psi_2))
\]
is in fact a non-zero morphism in
\[
\Hom_{ M_1}(  \sigma_1\psi_1, \sigma_2\psi_2)
\]
which is necessarily an isomorphism since $ \sigma_1\psi_1$ and $\sigma_2\psi_2$ are irreducible.
\end{proof}

\section{Notes on Chapter VII}

The proof of Casselman's criterion written here borrows from \cite{Be1}, \cite{Ca} and especially
\cite{DeB}. The application \ref{corang1} is taken from \cite{Be1}; recall that this is a crucial point
to show that the small progenerators of the components of the Bernstein decomposition exhibited in
\ref{progenMOm} are indeed generators. The section on tempered representations is taken from \cite{Wald},
the one on intertwining operators from \cite{Dat}. The exposition of the Langlands classification also
follows that of \cite{Dat}. The reader can also refer to \cite{Si2} for another exposition. For other
write-ups of Casselman's criterion, one may consult \cite{Ca}, Theorem 4.4.6, and \cite{Wald}, III.1.1.



\appendix
\renewcommand{\thesection}{\thechapter.\Roman{section}}
\renewcommand{\thesubsection}{\thesection.\arabic{subsection}}
\renewcommand{\theHequation}{\thesubsection.\arabic{equation}}

\chapter{Elements of category theory}

\section{Categories and functors}\label{catfonc}

We assume that the reader possesses some rudiments of category theory, i.e., at least the basic
definitions and terminology (object, morphisms, functors, subcategory...), as can be found in any good
book on the subject, for example \cite{Freyd} or \cite{Par}. For a concise and rigorous introduction (in
French) to this subject, covering all the aspects we need in this book, we refer to Chapter II of
\cite{Douady}. In addition to general statements which, for the convenience of the reader, we recall ---
and sometimes prove ---, we develop a certain number of heterogeneous examples, because they appeared in
the main body of the book. We will not delve into the logical foundations of category theory, which we
use from a naive point of view, hoping that everything we say can be made rigorous by standard procedures
(Grothendieck universes...).

Let $\scrC$ be a category. We simply write $X\in \scrC$ to say that $X$ is an object of $\scrC$ and
$f:X\rightarrow Y$ (or $X \stackrel{f}{\rightarrow}  Y$), to say that $f$ is a morphism from $X$ to $Y$.
We denote by $\Hom_\scrC(X,Y)$ the set of morphisms from $X$ to $Y$ and $\Id_X$ the identity of $X$.

The notation for functors is standard. If $F$ is a functor from the category $\caA$ to the category
$\caB$, if $X$, $Y$ are objects of $\caA$ and if $f$ is a morphism from $X$ to $Y$, we denote
respectively by $F(X)$, $F(Y)$ and $F(f)$ their images under the functor $F$. Often, in the mathematical
literature, functors are explicitly defined only on objects. Their effect on morphisms is generally
obvious, and it is part of the reader's task to make it explicit.

\subsection{Notable morphisms}

We begin by introducing some definitions concerning morphisms. A morphism
$X \stackrel{f}{\rightarrow}  Y$ is an {\sl isomorphism}\index[ter]{isomorphism}, if there exists a
morphism $Y \stackrel{g}{\rightarrow}  X$ such that $gf=\Id_X$ and $fg=\Id_Y$. It is a {\sl monomorphism}
\index[ter]{monomorphism} if the only pairs of morphisms $Z  \stackrel{h_1}{\rightarrow}  X$ and
$Z \stackrel{h_2}{\rightarrow}  X$ such that
\[
Z  \stackrel{h_1}{\rightarrow}  X  \stackrel{f}{\rightarrow}  Y= Z
\stackrel{h_2} {\rightarrow}  X  \stackrel{f}{\rightarrow}  Y
\]
are those such that $h_1=h_2$. It is an {\sl epimorphism} \index[ter]{epimorphism} if the only pairs of
morphisms $Y  \stackrel{h_1}{\rightarrow}  Z$ and $Y \stackrel{h_2}{\rightarrow}  Z$ such that
\[
X  \stackrel{f}{\rightarrow}  Y   \stackrel{h_1}{\rightarrow}  Z= X
\stackrel{f} {\rightarrow}  Y  \stackrel{h_2}{\rightarrow}  Z
\]
are those such that $h_1=h_2$.

An isomorphism is both an epimorphism and a monomorphism, but the converse is false in general. It is
true in abelian categories.

\subsection{Subobjects, quotient objects}\label{sousobj}

Let us now introduce the notions of {\sl subobjects}\index[ter]{subobject} and {\sl quotient objects}
\index[ter]{quotient object} of an object $X$ in a category $\scrC$. Consider the preorder $\leq$ on the
monomorphisms $Y \stackrel{f}{\rightarrow}  X$:
\[
Y_1 \stackrel{f_1}{\rightarrow}  X  \leq Y_2 \stackrel{f_2}{\rightarrow}  X
\]
if there exists a morphism $Y_1  \stackrel{i}{\rightarrow} Y_2$ such that
\[
Y_1 \stackrel{i} {\rightarrow} Y_2  \stackrel{f_2}{\rightarrow}  X = Y_1 \stackrel{f_1}{\rightarrow}  X.
\]
Note then that $i$ is a monomorphism. We say that $f_1$ and $f_2$ are equivalent if $f_1 \leq f_2$ and
$f_2\leq f_1$. A subobject of $X$ is an equivalence class of monomorphisms into $X$ for this equivalence
relation. We denote by \index[not]{S_X@$S_X$} $S_X$ the collection of subobjects of $X$. It is clear that
$\leq$ induces an order relation on $S_X$. Dually, we define the quotient objects of $X$: we equip the
collection of epimorphisms $X \stackrel{f}{\rightarrow}  Y$ with the preorder $\leq$ given by
\[
X \stackrel{f_1}{\rightarrow}  Y_1  \leq X  \stackrel{f_2}{\rightarrow}  Y_2
\]
if there exists a morphism $Y_2  \rightarrow Y_1$ such that
\[
X \stackrel{f_2} {\rightarrow} Y_2 \rightarrow   Y_1 = X
\rightarrow Y_1.
\]
We say that $f_1$ and $f_2$ are equivalent if $f_1 \leq f_2$ and $f_2\leq f_1$. A quotient of $X$ is an
equivalence class of epimorphisms and $\leq$ induces an order on the collection $Q_X$
\index[not]{Q_X@$Q_X$} of quotients of $X$.

If $C$ and $D$ are two subobjects of $X$, we call the {\sl intersection} \indexter{intersection} of $C$
and $D$, and we denote by $C \cap D$, a greatest lower bound of $C$ and $D$ in $S_X$. Similarly, we call
the {\sl union} \indexter{union} of $C$ and $D$, and we denote by $C \cup D$, a least upper bound of $C$
and $D$ in $S_X$. Of course, such objects do not necessarily exist but if they exist, they are unique. We
can also obviously extend these notions by defining the intersection and the union of an arbitrary family
of subobjects of $X$. We will say that the category $\scrC$ admits arbitrary (resp. finite) intersections
(resp. unions), if for any object $X$ of $\scrC$ and any family (resp. finite family) of subobjects of
$X$, the intersection (resp. the union) of this family exists.

Very often, instead of saying: (the class of) $Y  \stackrel{f}{\rightarrow}  X$ is a subobject of $X$, we
just say: $Y$ is a subobject of $X$, and we write $Y \subseteq  X$.

\subsection{Noetherian and finitely generated objects}

Let $X$ be an object of $\scrC$. Let $U$ be a subset of $S_X$ and $A$ a subobject of $X$ in $U$. We say
that $A$ is {\sl maximal} in $U$, if for any subobject $B$ of $X$ in $U$ such that $ A \subseteq B$, we
have $B=A$. We say that $X$ is {\sl Noetherian} if for any non-empty subset $U$ of $S_X$, there exists a
maximal element in $U$. We can also define the notion of a Noetherian object by a {\sl maximal chain
condition}. A chain in $S_X$ is a totally ordered subset of $S_X$. It is then easy to verify that the
object $X$ is Noetherian if and only if any chain in $S_X$ contains a maximal element (\cite{Douady},
Proposition 1.6.1). In practice, it suffices to consider well-ordered chains, or even increasing
sequences of subobjects. \index[ter]{Noetherian (object)}

Another very closely related notion is that of a {\sl finitely generated} object
\index[ter]{finitely generated (object)} in a category. We assume here that the category $\scrC$ admits
arbitrary unions. We say that an object $X$ of $\scrC$ is finitely generated if for any chain of proper
subobjects $(A_i)$ of $X$, $\bigcup_i A_i$ is still a proper subobject of $X$ (proper means not
isomorphic to $X$). We can characterize finitely generated objects by a compactness property. The object
$X$ of $\scrC$ is finitely generated if for any family of subobjects $(A_i)$ of $X$ such that
$\bigcup_i A_i=X$, there exists a finite number of subobjects $A_{i_1}, \ldots ,A_{i_n}$ such that
$\bigcup_{j=1}^n A_{i_j}=X$ (see \cite{Par} 4.10, Theorem 1).

It is clear from the definitions that a subobject of a Noetherian object is finitely generated. Indeed,
consider a chain of proper subobjects $(A_i)_{i}$ of a subobject $Y$ of a Noetherian object $X$. Then by
the Noetherian property, this chain admits a maximal element, say $B$, and it is clear that this maximal
object is a least upper bound of the family $(A_i)_{i}$, i.e., $B=\bigcup_i A_i\subset Y$. It follows
that $B=\bigcup_i A_i$ is a proper subobject of $Y$, which shows that $Y$ is finitely generated.
Conversely, if every subobject of $X$ is finitely generated, then $X$ is Noetherian. Indeed, let $B_i$ be
an increasing sequence of subobjects and set $B=\bigcup_i B_i$. Since $B$ is finitely generated, there
exist $i_1 \leq \ldots \leq  i_r$ such that $B=\bigcup_{j=1}^r B_{i_j}=B_{i_r}$, which shows that the
increasing sequence of the $B_i$ is stationary.

A category where every finitely generated object is Noetherian will be called {\sl Noetherian}.
\indexter{Noetherian (category)}

\begin{exemple}
A module $M$ over a unital ring $A$ is a finitely generated object in the category $\caM(A)$ if and only
if $M$ is finitely generated over $A$ in the algebraic sense (\cite{Douady}, Lemma 3.6.1). The category
$\caM(A)$ is Noetherian if and only if the ring $A$ is (left) Noetherian.
\end{exemple}

\subsection{Initial, final, zero object}\label{obnul}

\indexter{initial (object)}\indexter{final (object)}\indexter{zero (object)} Let us now define the
notions of {\sl initial}, {\sl final} or {\sl zero} object. An object $X$ in a category $\scrC$ is said
to be initial (resp. final) if for any object $Y$ of $\scrC$, $\Hom_\scrC(X,Y)$ (resp. $\Hom_\scrC(Y,X)$)
is a singleton. A zero object is an object that is both initial and final. Initial, final, or zero
objects are unique up to isomorphism. If $\scrC$ admits a zero object $0_\scrC$, all the sets
$\Hom_\scrC (X,Y)$ possess a distinguished element, namely the composition
\[
X \rightarrow  0_\scrC  \rightarrow Y
\]
which we will simply denote by $0$. Under these conditions, we can define the {\sl
kernel}\indexter{kernel} of a morphism $X \stackrel{f} {\rightarrow} Y $. It is a morphism
$K \rightarrow X$ such that
\[
K \rightarrow X    \stackrel{f} {\rightarrow} Y= K   \stackrel{0} {\rightarrow} Y
\]
and such that for any morphism $L \rightarrow X$ having the same property, namely
\[
L \rightarrow X    \stackrel{f} {\rightarrow} Y= L   \stackrel{0}
{\rightarrow} Y,
\]
there exists a unique morphism $L \rightarrow K$ such that
\[
L \rightarrow K  \rightarrow X= L \rightarrow X.
\]
The notion of {\sl cokernel}\indexter{cokernel} is dual: the cokernel of a morphism
$X \stackrel{f} {\rightarrow} Y $ is a morphism $Y \rightarrow K$ such that
\[
X    \stackrel{f} {\rightarrow} Y \rightarrow K = Y  \stackrel{0} {\rightarrow} K
\]
and such that for any morphism $Y \rightarrow L$ having the same property, namely
\[
X    \stackrel{f} {\rightarrow} Y \rightarrow L = Y   \stackrel{0}
{\rightarrow} L,
\]
there exists a unique morphism $K \rightarrow L$ such that
\[
Y \rightarrow K  \rightarrow L= Y \rightarrow L.
\]
We denote by $\ker f$ and $\coker f$ respectively the kernel and the cokernel of $f$ (when they exist).

\subsection{Product category, opposite category}

Let $\caA$ and $\caB$ be two categories. We form the category $\caA\times \caB$ whose objects are the
pairs $(X,Y)$, $X\in \caA$, $Y\in \caB$, and the morphisms between $(X,Y)$ and $(X',Y')$ are the pairs of
morphisms $(f,g)\in \Hom_\caA(X,X')\times \Hom_\caB (Y,Y')$. We also form the category $\caA^{op}$ whose
objects are the same as $\caA$, and for any pair $(X,X')$ of objects,
$\Hom_{\caA^{op}}(X,X')=\Hom_\caA(X',X)$. The use of the opposite category allows us to dispense with the
notion of a contravariant functor; consequently, all functors considered will always be covariant unless
explicitly stated otherwise. Another conceptual contribution of the introduction of opposite categories
is the notion of duality. As we have seen in the various definitions given above, these often come in
pairs, one being obtained from the other ``by reversing the direction of the arrows'', i.e., by working
in the opposite categories. Thus, the pairs monomorphism/epimorphism, subobject/quotient object, initial
object/final object, kernel/cokernel are dual. We will see many other instances of this principle in what
follows. In practice, any statement in category theory admits a dual statement, whose proof is obtained
by reversing the direction of the arrows.

\section[Equivalence of categories]{Natural transformation, equivalence of categories}

\subsection{Natural transformations}
\indexter{natural transformation} Let $\caA$ and $\caB$ be two categories, and $F,G$ two functors from
$\caA$ to $\caB$. A natural transformation $\theta$ from $F$ to $G$ is the data for any object $X$ of
$\caA$ of a morphism
\[
\theta_X\colon F(X)\rightarrow G(X)
\]
such that for any morphism $f:X\rightarrow Y$ in $\caA$, the diagram
\begin{equation*}
\begin{CD}
F(X) @>{F(f)}>> F(Y)\\
       @VV{\theta_X}V              @VV{\theta_Y}V\\
G(X)@>{G(f)}>>G(Y)
\end{CD}
\end{equation*}
commutes.

For any pair of categories $(\caA,\caB)$, we can form a category $\caF_{\caA,\caB}$
\index[not]{F_AB@$\caF_{\caA,\caB}$} whose objects are the functors from $\caA$ to $\caB$, and the
morphisms are the natural transformations between these functors. We therefore denote by $\Hom(F,G)$ the
set of natural transformations from $F$ to $G$ and $\Id_F$ the trivial natural transformation from the
functor $F$ to itself. A natural isomorphism between $F$ and $G$ is a natural transformation
$\theta \in  \Hom(F,G)$ admitting an inverse, i.e., a natural transformation $\phi \in \Hom(G,F)$ such
that $\theta\phi=\Id_G$, $\phi\theta=\Id_F$.

Any category $\caA$ is equipped with an identity functor, which we denote by $\caI d_\caA$:
\[
\caI d_\caA: \caA \rightarrow \caA,\quad X\mapsto X, \quad f\mapsto f,
\]
for any object $X$ and any morphism $f$ of $\caA$.

\subsection{Equivalences of categories}
\indexter{equivalence of categories} One of the important ideas at the origin of category theory is that,
in many cases, one is more interested in the isomorphism classes of objects of a given category than in
the objects themselves. This is particularly true in representation theory for example, where the
problems of classifying irreducible representations are always ``up to isomorphism''. Heuristically, if
we start from a category $\scrC$, imagine that we form the subcategory $\bar \scrC$ whose objects are
given by a choice of representatives of the isomorphism classes of objects of $\scrC$, the morphisms
between two representatives being all the morphisms between these objects in $\scrC$ (in other words, it
is a full subcategory). We would like a definition of equivalence of categories that makes $\scrC$ and
$\bar \scrC$ equivalent. Requiring the existence of functors $F: \scrC \rightarrow \bar \scrC$ and
$G: \bar \scrC \rightarrow  \scrC$ inverse to each other, i.e., such that $FG=\caI d_{\bar \scrC}$,
$GF=\caI d_{ \scrC}$ as a definition of equivalence is too strong a constraint to obtain this. The
appropriate definition is the following.

\begin{defi}
Let $F:\caA\rightarrow\caB$ be a functor between two categories $\caA$ and $\caB$. We say that $F$ is an
equivalence of categories if there exists a functor $G:\caB\rightarrow\caA$ such that $FG$ is naturally
isomorphic to $\caI d_{\caB}$, and $GF$ naturally isomorphic to $\caI d_{ \caA}$. The categories $\caA$
and $\caB$ are then said to be equivalent and the functors $F$ and $G$ are said to be quasi-inverses.
\end{defi}

\begin{exemple}
Let $\mathbf{Vect_k^n}$ be the category of vector spaces of dimension $n$ over the field $k$ whose
morphisms are linear maps. Let $V_k^n$ be the category whose unique object is the vector space $k^n$ and
whose morphisms are the endomorphisms of $k^n$. As functor $F$, we take the inclusion of $V_k^n$ into
$\mathbf{Vect_k^n}$. It is an equivalence of categories.

This example well illustrates the idea expressed above, since all the objects of $\mathbf{Vect_k^n}$ are
isomorphic to $k^n$. Note that the construction of a quasi-inverse $G$ requires the choice of such an
isomorphism for any object (in other words the choice of a basis). There is therefore no uniqueness of
the quasi-inverse, and its construction uses the axiom of choice.
\end{exemple}

The functor $F:\caA\rightarrow\caB$ is said to be {\sl full} \indexter{full (functor)} (resp. {\sl
faithful}) \indexter{faithful (functor)} if for any pair of objects $(X,Y)$ in $\caA$, $F$ realizes a
surjection (resp. an injection) between $\Hom_\caA(X,Y)$ and $\Hom_\caB(F(X),F(Y))$. A full and faithful
functor is said to be fully faithful. A subcategory $\caA$ of a category $\caB$ is said to be full if the
inclusion functor $\caA \rightarrow \caB$ is fully faithful.

We often use the following criterion to show that a functor is an equivalence of categories
(\cite{GelfMan}, Theorem 1.13 or \cite{Douady}, Theorem 2.3.5).

\begin{thm}
A functor $F:\caA\rightarrow\caB$ is an equivalence of categories if it is fully faithful and if any
object of $\caB$ is isomorphic to an object of the form $F(X)$, where $X \in \caA$.
\end{thm}

\section{Universal problems and representable functors}

We now introduce the notion of a universal problem \indexter{universal problem} and of a solution to such
a problem. Certain examples chosen to illustrate this notion have appeared in the text, such as for
example the localization of a vector space with respect to an endomorphism. We take this opportunity to
establish a lemma and its corollary, used in \ref{LemmeJacquet2}.

\subsection{Universal problems}\label{propuniv}
\begin{defi}
Let $F:\caA\rightarrow\caB$ be a functor and let $Y$ be an object of $\caB$. We say that the object $X$
of $\caA$ and the morphism $i:Y \rightarrow F(X)$ are solutions of the right universal problem
\indexter{universal problem} posed by $F$ and $Y$ if for any object $Z$ of $\caA$ and any morphism
$f: Y\rightarrow F(Z)$, there exists a unique morphism $g:X\rightarrow Z$ satisfying $F(g)\circ i=f$.
Similarly we say that $X \in \caA$ and $j: F(X) \rightarrow Y$ are solutions of the left universal
problem posed by $F$ and $Y$ if for any object $Z$ of $\caA$ and any morphism $f: F(Z) \rightarrow Y$,
there exists a unique morphism $g:Z\rightarrow X$ satisfying $j\circ F(g)=f$.
\end{defi}

\begin{rmq}
If two pairs $(X,i)$ and $(X',i')$ are solutions of the (right) universal problem posed by $F$ and $Y$,
then $X$ and $X'$ are isomorphic, up to a unique isomorphism, as can be seen by taking as morphism $f$
successively the morphisms $i$ and $i'$. We indeed obtain two morphisms $g_1: X\rightarrow X'$ and
$g_2: X'\rightarrow X$ satisfying $F(g_1)\circ i=i'$ and $F(g_2)\circ i'=i$, whence
$F(g_2 \circ g_1)\circ i=i$ and $F(g_1 \circ g_2)\circ i'=i'$. Since the morphisms $\Id_X$ and $\Id_{X'}$
also satisfy $F(\Id_X)\circ i=i$ and $F(\Id_{X'})\circ i'=i'$, by uniqueness, $g_2 \circ g_1=\Id_X$ and
$g_1 \circ g_2= \Id_{X'}$. The same holds for left universal problems.
\end{rmq}

\begin{exemples}\label{complete}

1. {\sl Coproducts, products.} Let $\caA$ be a category. Consider the ``diagonal'' functor
\[
\Delta: \caA \rightarrow \caA\times \caA,\quad X \mapsto (X,X),
\quad f \mapsto (f,f)
\]
and let us be given an element $(X,Y)$ of $\caA\times \caA$. A solution of the left universal problem
posed by $\Delta$ and $(X,Y)$ is called the product of $X$ and $Y$ and denoted $X\times Y$. Similarly the
coproduct (or ``direct sum''), denoted $X \amalg Y$ (or $X\oplus Y$), is the solution of the right
universal problem posed by $\Delta$ and $(X,Y)$.

--- 2. {\sl Completion of a topological space. } Let $\mathbf{EvTop}$ and $\mathbf{\overline{EvTop}}$ be
the categories of topological vector spaces and complete topological vector spaces respectively, the
morphisms being continuous linear maps. The functor we consider is the (fully faithful) inclusion of
$\mathbf{\overline{EvTop}}$ into $\mathbf{EvTop}$. We are given $Y \in \mathbf{EvTop}$. A solution of the
right universal problem posed by the inclusion and $Y$ is the completion of $Y$, and the continuous
inclusion $Y \hookrightarrow  \bar Y$. More generally, one can define the completion of a topological
space admitting a uniform structure, the uniform structure being necessary to define the notion of a
Cauchy filter (see \cite{BourTG}).

--- 3. {\sl Localization of a ring. } Let $\caA$ be the category whose objects are pairs $(A,S)$
consisting of a commutative unital ring $A$ and a multiplicative subset $S$ of $A$, the morphisms between
two objects $(A,S)$ and $(A',S')$ being unital ring morphisms $\phi\colon A \rightarrow A'$ such that
$\phi(S) \subset S'$. We consider the full subcategory $\caA'$ of $\caA$ whose objects are pairs $(A,S)$,
with $S$ contained in the invertible elements of $A$. The solution of the right universal problem posed
by an element $(A,S)$ of $\caA$ and the inclusion of $\caA'$ into $\caA$ is the localization $S^{-1}A$ of
the ring $A$ at $S$ and the natural morphism $A \rightarrow S^{-1}A$. Note that if $0\in S$, we have
$S^{-1}A=\{0\}$.

--- 4. {\sl Localization of a module. } Let $A$ be a commutative unital ring, and $S$ a multiplicative
subset of $A$. Consider the full subcategory $\caM(A)_S$ of $\caM(A)$ whose objects are modules $M$ such
that the action of any element $s \in S$ on $M$ is invertible (i.e., for all $s \in S$, for all
$m \in M$, there exists $t\in A$ such that $ts\cdot m=st\cdot m=m$). The right universal problem posed by
a module $M$ of $\caM(A)$ and by the forgetful functor from $\caM(A)_S$ to $\caM(A)$ admits as a solution
the localization $S^{-1}M$ and the natural morphism $M \rightarrow S^{-1}M$ (viewed simply as an
$A$-module, and not as an $S^{-1}A$-module). Here again, if $0\in S$, we have $S^{-1}M=\{0\}$.

--- 5. {\sl Localization of a vector space with respect to an endomorphism. } \indexter{localization} Let
$k$ be a field, $k[X]$ the ring of polynomials with coefficients in $k$, and $\caM(k[X])$ the category of
left unital $k[X]$-modules. A $k[X]$-module can be viewed as the data of a vector space over $k$ and an
endomorphism of this vector space. The category $\caM(k[X])$ is therefore equivalent to the one whose
objects are pairs $(L,a)$ where $L$ is a vector space over $k$ and $a \in \End_k (L)$, the morphisms
between two objects $(L,a)$ and $(L',a')$ being linear maps between $L$ and $L'$ intertwining $a$ and
$a'$. Consider in $k[X]$ the multiplicative subset $\{X^k\}_{k\in \bbN}$ of powers of $X$. The previous
example asserts the existence, for any $k[X]$-module $L$, of a localization with respect to the subset
$\{X^k\}_{k\in \bbN}$. Let us translate this in terms of vector spaces equipped with endomorphisms by the
equivalence of categories described above: given a pair $(L,a)$ where $L$ is a vector space over $k$ and
$a \in \End_k (L)$, there exists a pair $(L_a,\tilde a)$, where $L_a$ is a vector space over $k$ and
$\tilde a \in \End_k (L_a)$ is invertible, and a morphism $\iota\colon (L,a) \rightarrow (L_a,\tilde a)$
such that for any morphism $f: (L,a) \rightarrow (M,b)$, where $b \in \End_k (M)$ is invertible, there
exists a unique morphism $\tilde f: (L_a,\tilde a) \rightarrow (M,b)$ such that $f=\tilde f \circ \iota$.
We call $(L_a,\tilde a)$ the localization of $L$ at $a$. If $B$ is a $k$-algebra, we can make the same
construction by replacing $k$-vector spaces by $B$-modules, and endomorphisms of vector spaces by
endomorphisms of $B$-modules.
\end{exemples}

We provide a slightly more explicit description of the localization of $L$ at $a$ in the fifth example
above. We state the necessary results in the following lemma:

\begin{lemme}
Let $L$ be a vector space over a field $k$ and $a$ an endomorphism of $L$.

$(i)$ If $a$ is injective, the canonical morphism $\iota \colon (L,a) \rightarrow (L_a,\tilde a )$ is
also injective.

$(ii)$ Let $K:=\bigcup_n \ker a^n$. Then $K$ is stable under $a$. Let $a'$ be the endomorphism of
$L'=L/K$ induced by $a$. Then the localization $(L'_{a'}, \tilde a')$ of $L'$ at $a'$ is isomorphic to
the localization of $L$ at $a$. Moreover $a'$ is injective, as well as the canonical morphism $\iota'$
from $(L',a')$ to $(L'_{a'}, \tilde a')$. The space $L'_{a'}$ is therefore an extension of $L'$, such
that $\tilde a'$ coincides with $a'$ on $L'$, and $L'_{a'}= \bigcup_{n \in \bbN}  (\tilde a')^{-n} (L')$.
\end{lemme}

\begin{proof}
$(i)$ It suffices to construct an injective morphism $$f \colon  (L,a)\rightarrow (M,b),$$ with $b$
invertible. Since $f$ factors as $f=\tilde f \circ \iota$, it is clear that $\iota$ must then be
injective. Let us construct such a morphism. Let $R_0$ be a complement of $\im a$ in $L$. Let us
construct a sequence of vector spaces $R_i$ and isomorphisms $f_i:R_{i+1} \rightarrow R_i$. Set
$M=L\oplus(\bigoplus_{i=1}^\infty R_i)$ and $b\colon M  \rightarrow M$, $b=a+ \sum_{i=0}^\infty f_i$. By
construction, $b$ is injective and surjective, so it is an isomorphism. Then take as morphism $f$ the
inclusion of $L$ into $M$.

$(ii)$ Note that an element $l\in L$ such that $a^n(l)=0$, for a certain $n\in \bbN^*$, is in the kernel
of $\iota$, since $\iota \circ a^n=\tilde a^n \circ\iota$ and $\tilde a^n$ is invertible. This implies
that
\[
K:=\bigcup_n \ker a^n \subset \ker \iota.
\]
We can therefore factor $\iota$ through $L'=L/K$. Let $a'$ be the morphism induced by $a$ on $L'$, and
let $(L'_{a'},\tilde a')$ be the localization of $L'$ at $a'$. Using the universal properties satisfied
by $(L'_{a'},\tilde a')$ and $(L_{a},\tilde a)$, we see that these are isomorphic. On the other hand,
$a'$ is injective. By $(i)$, the canonical morphism $\iota'$ from $(L',a')$ to $(L'_{a'},\tilde a')$ is
also injective. We then identify $L'$ with a subspace of $L'_{a'}$ by this injection $\iota'$.

Set $M=\bigcup_{n \in \bbN}  (\tilde a')^{-n} (L')$. We have $L'\subset M \subset L'_{a'}$ and $M$ is
stable under $\tilde a'$. It easily follows that $M=L'_{a'}$. Conversely, if $(M,\beta)$ is an extension
of $L'$ satisfying $\beta_{|L'}=a'$, $\beta$ invertible and $M= \bigcup_{n \in \bbN} \beta^{-n} (L')$,
$(M,\beta)$ satisfies the universal property characterizing the localization $(L'_{a'}, \tilde a')$.
\end{proof}

\begin{cor}
Suppose that the canonical morphism $\iota \colon  L \rightarrow L_a$ is surjective. Let
$K=\bigcup_{n \in \bbN} \ker a^n$, and $L'=L/K$. Then $(L_{a},\tilde a)$ is isomorphic to $(L',a')$.
\end{cor}

\subsection{Representable functors} \label{representable}
Let $\scrC$ be a category. Consider the category $\caF_{\scrC^{op},\mathbf{Ens}}$ of contravariant
functors from $\scrC$ to $\mathbf{Ens}$. For each object $X$ of $\scrC$, we can construct a functor $h_X$
in $\caF_{\scrC^{op},\mathbf{Ens}}$ by setting for all $Y \in \scrC$,
\[
h_X(Y)=\Hom_\scrC(Y,X),
\]
and for any morphism $f:Y \rightarrow Z$ in $\scrC^{op}$ (i.e., $f$ is a morphism in $\Hom_\scrC(Z,Y)$),
\[
h_X(f): h_X(Y)\rightarrow h_X(Z), \quad \phi \mapsto \phi \circ f.
\]

Let $\phi: X_1 \rightarrow X_2$ be a morphism in $\scrC$. We associate to $\phi$ a natural transformation
$h_\phi$ from $h_{X_1}$ to $h_{X_2}$ by setting, for any object $Y$ of $\scrC$,
\[
h_{\phi,Y} : \Hom_\scrC(Y,X_1) \rightarrow  \Hom_\scrC(Y,X_2),\quad  g \mapsto
\phi \circ g.
\]
It is clear that $h_\psi h_\phi=h_{\psi\phi}$ when $\psi$ and $\phi$ are composed in $\scrC$. This
therefore defines a functor
\[
h : \scrC \rightarrow \caF_{\scrC^{op},\mathbf{Ens}}
\]

\begin{defi}
A functor $F$ in $\caF_{\scrC^{op},\mathbf{Ens}}$ is said to be representable if there exists a natural
isomorphism $\varpi\colon h_X \rightarrow F$, for a certain $X$ in $\scrC$. We then say that $(X,\varpi)$
represents the functor $F$.
\end{defi}

\begin{thm}[Yoneda]
The functor $h$ realizes an equivalence of categories between $\scrC$ and the full subcategory of
$\caF_{\scrC^{op},\mathbf{Ens}}$ whose objects are the representable functors. In particular, $h$ is
fully faithful, i.e., $\phi\mapsto h_\phi$ realizes a bijection
\[
\Hom_\scrC(X,Y) \simeq \Hom_{\caF_{\scrC^{op},\mathbf{Ens}}}(h_X,h_Y)
\]
\end{thm}

\begin{proof}
Let $F \in \caF_{\scrC^{op},\mathbf{Ens}}$ and $X$ be an object of $\scrC$. We will now show that the set
of natural transformations from $h_X$ to $F$ is in bijection with the set $F(X)$. Let $x \in F(X)$. For
all $Y \in \scrC$ set
\[
\theta_{x,Y}\colon h_X(Y)=\Hom_\scrC(Y,X)\rightarrow F(Y), \quad f \mapsto
F(f)(x).
\]
Then $\theta_x$ is a natural transformation from $h_X$ to $F$. Indeed, for any morphism
$g: Y\rightarrow Z$ in $\scrC^{op}$ (i.e., $g: Z\rightarrow Y$ in $\scrC$), the diagram
\[
\xymatrix{
h_X(Y)\ar[r]^{\theta_{x,Y}}\ar[d]_{h_X(g)}
&  F(Y) \ar[d]^{F(g)} \\
h_X(Z) \ar[r]_{\theta_{x,Z}}& F(Z)\\ }
\]
commutes.

In the other direction, suppose that $\theta$ is a natural transformation from $h_X$ to $F$. Set
$\eta(\theta)=\theta_X(\Id_X)\in F(X)$. We will now show that the maps $x \mapsto \theta_x$ and
$\theta\mapsto \eta(\theta)$ are inverse to each other. We have on the one hand
\[
\eta(\theta_x)=\theta_{x,X}(\Id_X)=\Id_{F(X)}(x)=x,
\]
and on the other hand for all $Y \in \scrC$ and for all $f\in h_X(Y)$
\[
\theta_{\eta(\theta),Y}(f)=F(f)(\eta(\theta))=F(f)( \theta_X(\Id_X))=\theta_Y(h_X(f)(\Id_X))=\theta_Y(f).
\]
This proves the assertion.

Let us apply this with $F=h_Y$. We obtain
\[
\Hom_\scrC(X,Y)= h_Y(X) \simeq
\Hom_{\caF_{\scrC^{op},\mathbf{Ens}}}(h_X,h_Y),\quad  \phi \mapsto \theta_
\phi.
\]
Now it is immediate to verify that $\theta_\phi=h_\phi$. This shows that $h$ is fully faithful.
\end{proof}

\begin{rmq}
The functor $F$ is represented by $(X,\varpi)$ if and only if $(X,\varpi)$ is a solution of the (left)
universal problem posed by $h$ and $F$. In the other direction, if we consider the left universal problem
posed by $F: \caA \rightarrow \caB$ and $Y \in \caB$, we obtain a solution $(X,j)$ as soon as the functor
\[
Z \mapsto \Hom_\caB(F(Z),Y),\quad \caA^{op} \rightarrow \mathbf{Ens}
\]
is representable by $(X,\varpi)$, and $j=\varpi_X(\Id_X)$.
\end{rmq}

\begin{proof}
Suppose that $F$ is represented by $(X,\varpi)$. We verify the desired universal property: let
$\theta \colon  h_Z\rightarrow F$ be a natural transformation. We want to show that there exists a unique
morphism $\phi\colon Z \rightarrow X$ such that $\varpi\circ h_\phi=\theta$. As in the proof of the
theorem, $\theta$ corresponds to an element $\eta(\theta)$ of $F(Z)$, and thus
$\varpi^{-1}(\eta(\theta))$ is an element of $h_X(Z)=\Hom_\scrC(Z,X)$. Take
$\phi=\varpi^{-1}_Z(\eta(\theta))$. We then have for all $Y\in \scrC$, and all
$f \in h_Z(Y)=\Hom_\scrC(Y,Z)$
\begin{align*}
& (\varpi\circ h_{\varpi^{-1}_Z(\eta(\theta))})_Y(f)
  =\varpi_Y(\varpi^{-1}_Z(\eta(\theta))\circ f)\\
=&\varpi_Y(\varpi^{-1}_Z(\theta_Z(\Id_Z))\circ f)=\varpi_Y(\varpi^{-1}_Y(\theta_Y(f)))=\theta_Y(f).
\end{align*}
Uniqueness is immediate, because $\varpi$ being an isomorphism, $\theta$ determines $h_\phi$, and thus
$\phi$ by the theorem.

Conversely, if $(X,\varpi)$ is a solution of the (left) universal problem posed by $h$ and $F$, we must
verify that $\varpi$ is an isomorphism, i.e., that for all $Y \in \scrC$,
\[
\varpi_Y: h_X(Y)\rightarrow F(Y)
\]
is an isomorphism. Let us construct its inverse: let $y \in F(Y)$. To it corresponds, by Yoneda's
theorem, a natural transformation $\theta_y: h_Y \rightarrow F $, which by the universal property, we can
factor through $\varpi$: $\theta_y=\varpi \circ h_\phi$, for a certain $\phi \in \Hom_\scrC(Y,X)=h_X(Y)$.
It is now easy to verify that the map $y \mapsto \phi$ thus defined is the inverse of $\varpi_Y$.

We leave the proof of the second assertion to the reader.
\end{proof}

\begin{cor}
If $F$ is a representable functor, the object representing $F$ is determined up to a unique isomorphism.
\end{cor}

\subsection{An isomorphism principle}\label{egalitecat}
We now reformulate Yoneda's lemma in a form useful for showing that two objects of a category are
isomorphic, or that two functors are naturally isomorphic.

\begin{prop}
Let $X,Y$ be two objects of a category $\caM$ such that for any object $Z$ of $\caM$, we have an
isomorphism, natural in $Z$:
\begin{equation}\label{hyphop}
\Hom_\caM(Z,X) \simeq \Hom_\caM(Z,Y).
\end{equation}
Then $X\simeq Y$. If $F,G: \scrC \rightarrow  \caM$ are two functors from a category $\scrC$ to $\caM$
such that we have a natural isomorphism in $X$ and $Z$:
\begin{equation}\label{hyphop2}
\Hom_\caM(Z,F(X)) \simeq \Hom_\caM(Z,G(X)).
\end{equation}
Then the functors $F$ and $G$ are naturally isomorphic.
\end{prop}

\section{Limits and colimits}\label{racines}

\subsection{Limits}\label{limites}
Let $\scrC$ and $\caI$ be two categories, the second, for reasons that will soon become apparent, being
called an ``index category'', or ``diagram scheme''. For any object $X$ of $\scrC$, let us define the
constant functor
\[
\Delta_X \colon \caI \rightarrow \scrC
\]
where $\Delta_X(i)=X$ for any object $i$ of $\caI$ and $\Delta_X(i\rightarrow j)=\Id_X$ for any morphism
$i\rightarrow j$ in $\caI$. It is clear that if $\caI$ is not empty, the constant functors from $\caI$ to
$\scrC$ correspond bijectively to the objects of $\scrC$, and that the natural transformations between
two constant functors $\Delta_X$ and $\Delta_Y$ correspond bijectively to the morphisms from $X$ to $Y$.
Indeed, if $X \stackrel{f}{\rightarrow} Y$ is a morphism in $\scrC$, we define the natural transformation
\[
\Delta_f \colon  \Delta_X \rightarrow \Delta_Y
\]
by setting
\[
\Delta_f(i)=f
\]
for any object $i$ of $\caI$. We obtain a functor
\[
\Delta \colon  \scrC \rightarrow \caF_{\caI,\scrC}.
\]

Let $D  \colon  \caI \rightarrow \scrC$ be a functor (such a functor is called a diagram in $\scrC$, of
scheme $\caI$). A {\sl limit} \indexter{limit} of $D$ (if it exists) consists of a constant functor
$\Delta_X \colon  \caI \rightarrow \scrC $ (which can be viewed simply as the data of the object $X$) and
a natural transformation $\eta$ from $\Delta_X$ to $D$, such that for any constant functor
$\Delta_Y\colon \caI \rightarrow \scrC $ and any natural transformation from $\Delta_Y$ to $D$, there
exists a unique natural transformation $\Delta_Y \rightarrow \Delta_X$ (i.e., a unique morphism
$Y\rightarrow X$) such that
\[
\Delta_Y \rightarrow \Delta_X \stackrel{\eta}{\rightarrow} D= \Delta_Y
\rightarrow D.
\]
In a way, the natural transformation $\Delta_X \rightarrow D$ is the ``best approximation'' of $D$ by a
constant functor.

\begin{rmq}
We see that $(\Delta_X,\eta)$ is a solution of the left universal problem posed by $\Delta$ and $D$.
There is therefore uniqueness of the limit up to a unique isomorphism, and we denote it by
$\varprojlim D$. On the other hand, using Remark \ref{representable} linking universal problem and
representable functor, we obtain for all $Y \in \scrC$ a natural isomorphism
\begin{equation}\label{limHom}
 \Hom_\scrC(Y,\varprojlim D)\simeq
\Hom_{\caF_{\caI,\scrC}}(\Delta_{Y}, D).
\end{equation}
\end{rmq}

More concretely, a limit of $D$ is therefore the data of an object $X$ of $\scrC$, and for any object
$i $ of $\caI$ of a morphism from $X$ to $D(i)$, such that for any morphism $i \rightarrow j$ in $\caI$,
\begin{align}\label{compracines}
 [X \rightarrow D(i) \rightarrow
  D(j)]= [X \rightarrow D(j)].
\end{align}
The fact that $\Delta_X \rightarrow D$ is the best approximation of $D$ by a constant functor translates
to: if there exists an object $Y $ of $\scrC$ and for any object $i$ of $\caI$ morphisms
$Y \rightarrow D(i)$ satisfying compatibility relations analogous to (\ref{compracines}), then there
exists a unique morphism $Y \rightarrow X$ such that for all $i\in \caI$
\[
Y  \rightarrow X   \rightarrow D(i)=  Y  \rightarrow D(i).
\]

\begin{exemples}

\noindent 1. { \sl Final objects}. A final object in a category $\scrC$ is a limit. Indeed, let us take
for $\caI$ the empty category (no objects, no morphisms). A limit of the empty diagram consists of an
object $T$ of $\scrC$ such that for every object $Y$ of $\scrC$, there exists a unique morphism
$Y \rightarrow T$: it is a final object.

\medskip

--- 2. {\sl Kernels.} Let $f\colon X \rightarrow Y$ be a morphism in a category $\scrC$ equipped with a
zero object. Let $\caI$ be the category having two objects $i$ and $j$ and two morphisms besides $\Id_i$
and $\Id_j$, namely $i \stackrel{1}{\rightarrow} j$ and $i \stackrel{0}{\rightarrow} j$. Let $D$ be the
functor such that $D(i)=X$, $D(j)=Y$, $D(1)=f$ and $D(0)=0$. Then the limit of the functor $D$ is the
kernel of $X \rightarrow Y$.

\medskip

--- 3. {\sl Products.} Let $X$ and $Y$ be two objects of a category $\scrC$. Consider the discrete
category $\caI$ formed by two objects $i$ and $j$ (discrete means that the only morphisms are the
identities of the objects). Let $F\colon \caI \rightarrow \scrC$ be the functor defined by $F(i)=X$ and
$F(j)=Y$. Then the limit of $F$ is the product of $X$ and $Y$, denoted by $X\times Y$.

More generally, for any family $(X_i)_{i\in I}$ of objects of $\scrC$, we form the discrete category
$\caI$ whose objects are the $i\in I$, and we define a functor $F :\caI \rightarrow \scrC$ by $F(i)=X_i$.
The limit of $F$ is then the product of the $X_i$, denoted $\prod_{i\in I} X_i$.

\medskip

--- 4. {\sl Fiber product.} Let $\scrC$ be a category and in $\scrC$ consider the morphisms
\begin{equation}\label{pullback}
 \xymatrix{  & X\ar[d]^f   \\
Y\ar[r]^g & Z    }
\end{equation}
Let $\caI$ be the category formed by three objects $i,j,p$, whose morphisms, besides the identities of
the objects, are $i\rightarrow p$ and $j\rightarrow p$, and let $D$ be the functor from $\caI$ to $\scrC$
such that $D(i)=X$, $D(j)=Y$, $D(p)=Z$, $D(i\rightarrow p)=f$ and $D(j\rightarrow p)=g$. The limit of $D$
is the fiber product, denoted $X\times_Z Y$, of the diagram (\ref{pullback}).

\medskip

--- 5. {\sl Projective or inverse limits.} Let $\caI$ be an {\sl upward directed poset} (i.e., for any
pair $(i, j)$ of elements of $\caI$, there exists an element $k$ of $\caI$ such that $i\leq k$,
$j\leq k$). We consider $\caI$ as a category, where the morphisms, besides the identities of the objects,
are given by the pairs $(i, j)$ with $i\leq j$, the corresponding morphism being $j\rightarrow i$. Let
$D:\caI \rightarrow \scrC$ be a functor, and set, for all $i$, $j$ in $\caI$, $A_i=D(i)$,
$D(j \rightarrow i)=f_{ji}$. We say that $((A_i)_{i\in I}, (f_{ji})_{i\leq j})$ is an {\sl inverse
system}. \indexter{inverse system} Then the projective limit (or inverse limit) of the system of
$(A_i,f_{ji})$ is the limit of $D$.

--- 6. {\sl Intersection of subobjects.} We leave it to the reader to interpret an intersection of
subobjects of an object of a category as a limit in an appropriate category.

\end{exemples}

\subsection{Colimits} \label{colimites}
We define dually the notion of colimit \indexter{colimit} of a functor $D \colon  \caI \rightarrow \scrC$
as in the previous paragraph: it consists of a constant functor
$\Delta_X \colon  \caI \rightarrow \scrC $ and a natural transformation from $D$ to $\Delta_X$, such that
for any constant functor $\Delta_Y \colon  \caI \rightarrow \scrC $ and any natural transformation from
$D$ to $\Delta_Y$, there exists a unique natural transformation $\Delta_X \rightarrow \Delta_Y$ (i.e., a
unique morphism $X\rightarrow Y$) such that
\[
D \rightarrow \Delta_X \rightarrow \Delta_Y = D \rightarrow \Delta_Y.
\]
A colimit is a solution of the right universal problem posed by $\Delta$ and $D$.

As examples of colimits, we have the dual notions of those of the examples of the previous paragraph. We
thus obtain respectively initial objects, cokernels, coproducts, pushouts, inductive (or direct) limits,
unions of subobjects. The notation used for colimits is $\varinjlim$. \index[not]{limi@$\varinjlim$}

Let us make the example of inductive limits a little more explicit. Let $I$ be an upward directed poset.
We consider $\caI$ as a category, where the morphisms, besides the identities of the objects, are given
by the pairs $(i, j)$ with $i\leq j$, the corresponding morphism being $i \rightarrow j$. Let
$D:\caI \rightarrow \scrC$ be a functor, and set, for all $i$, $j$ in $\caI$, $A_i=D(i)$,
$D(i \rightarrow j)=f_{ij}$. We say that $((A_i)_{i\in I}, (f_{ij})_{i\leq j})$ is a {\sl direct system}
\indexter{direct system}. Then the inductive limit (or direct limit) of the system of $(A_i,f_{ij})$ is
the colimit of $D$.

\begin{exemple}
Let $\caA$ be the category of abelian groups. In this category, any family $(G_i)_{i\in I}$ admits a
product, given by the Cartesian product of sets $\prod_{i \in I}G_i$, equipped with the product group law
of that of the $G_i$, and canonical projections $p_i: \prod_{i \in I}G_i \rightarrow G_i$. Any family
$(G_i)_{i\in I}$ also admits a coproduct, called the {\sl direct sum} of the $G_i$ and denoted
$\bigoplus_{i \in I}G_i$. If $I$ is finite, $\bigoplus_{i \in I}G_i$ is the product $\prod_{i \in I}G_i$
equipped with the canonical injections $i_i: G_i \rightarrow \prod_{i \in I}G_i$. If $I$ is infinite,
$\bigoplus_{i \in I}G_i$ is the subgroup formed by the elements of $\prod_{i \in I}G_i$ having only a
finite number of non-zero coordinates, equipped with these same injections.

Now let $((G_i)_{i\in I}, (f_{ij})_{i\leq j})$ be a direct system of abelian groups over an upward
directed poset $I$ (the group laws are denoted additively). Let us form the coproduct of the $G_i$ in the
category of sets, i.e., their disjoint union. Let us equip it with the following equivalence relation
(which is transitive precisely because $I$ is directed): $g_i \in G_i$ is equivalent to $g_j\in G_j$ if
there exists $k\geq i,j$ such that $f_{ik}(g_i)=f_{jk}(g_j)$ and let us call $\underline{G}$ the set of
equivalence classes. Then $\underline{G}$ is naturally equipped with a group structure, and we have
canonical projections $f_i \colon G_i \rightarrow \underline{G}$.

On the other hand, $\underline{G}$ satisfies the universal property of the colimit: for any system of
morphisms $\phi_i \colon  G_i \rightarrow Y$ compatible with the $f_{ij}$, the $\phi_i$ factor through
the $f_i$. We thus have $\varinjlim_{i \in I}G_i=\underline{G}$.

We immediately deduce the following fact, which is fundamental in sheaf theory: if $g_i \in G_i$ is such
that $f_i(g_i)=0$, then there exists $j \geq i$ such that $f_{ij}(g_i)=0$.

\end{exemple}

\section{Adjoint functors}\label{adjfonct}

There are two (equivalent) ways of viewing the notion of adjoint functors. Depending on the application,
one perspective may be more suitable than the other. \indexter{adjoint (functor)}

\subsection{Definition by $\Hom$ functors}

Let $\caA$ and $\caB$ be two categories, and $F \colon \caA \rightarrow \caB$,
$G \colon \caB \rightarrow \caA$ two functors. Let us define the functors:
\[
R_F \colon  \caA \times \caB^{op} \rightarrow \mathbf{Ens}, \quad (X,Y)\mapsto \Hom_\caB(Y,F(X)),
\]
and
\[
L_G \colon   \caA \times \caB^{op} \rightarrow \mathbf{Ens}, \quad (X,Y)\mapsto \Hom_\caA(G(Y),X).
\]

If $ (f,g)\in \Hom_\caA(X,X')\times \Hom_\caB(Y',Y)$ the morphism $R_F(f,g)$ is given by:
\[
R_F(f,g) \colon  \Hom_\caB(Y,F(X)) \rightarrow \Hom_\caB(Y',F(X')), \quad
\psi \mapsto F(f)\circ \psi\circ g,
\]
and the morphism $L_G(f,g)$ by:
\[
L_G(f,g) \colon   \Hom_\caA(G(Y),X) \rightarrow \Hom_\caA(G(Y'),X'), \quad \psi \mapsto f\circ \psi\circ G(g).
\]

\begin{defi}
We say that $F$ is the right adjoint of $G$ (or that $G$ is the left adjoint of $F$) if the functors
$R_F$ and $L_G$ are naturally isomorphic. Let $\theta$ denote this natural isomorphism. We then have,
\[
\theta_{Y,X}\colon \Hom_\caB(Y,F(X)) \simeq \Hom_\caA(G(Y),X),\quad (X\in \caA),\,
(Y\in \caB).
\]
\end{defi}

The following theorem links the notions of adjoint functors and representable functors.
\begin{thm}
Let $G:\caB \rightarrow \caA$ be a functor. Then $G$ admits a right adjoint if and only if for any object
$X$ of $\caA$, the contravariant functor
\[
Y \mapsto \Hom_\caA(G(Y),X)
\]
from $\caB$ to $\mathbf{Ens}$ is representable.
\end{thm}

\begin{proof}
The condition is obviously necessary, since if $G$ admits $F$ as a right adjoint, the functor
$Y \mapsto \Hom_\caA(G(Y),X)$ is represented by $(F(X),\theta)$. We will now show that it is sufficient.
Recall that by Yoneda's lemma, the functor
\[
h\colon \caB \rightarrow \caF_{\caB^{op}, \mathbf{Ens}}
\]
is fully faithful. It therefore realizes an equivalence of categories onto the full subcategory of
representable functors of $\caF_{\caB^{op}, \mathbf{Ens}}$. Let $\Lambda$ be a quasi-inverse of $h$ and
let $H\colon \caA \rightarrow  \caF_{\caB^{op}, \mathbf{Ens}} $ be the functor defined by
$H(X)(Y)=\Hom(G(Y),X)$, which by hypothesis takes values in the full subcategory of representable
functors of $\caF_{\caB^{op}, \mathbf{Ens}}$. Then $F=\Lambda \circ H$ is the desired adjoint.
\end{proof}

\begin{cor}
There is uniqueness up to isomorphism of the right adjoint of a functor $G \colon \caB \rightarrow \caA$.
The same is true for left adjoints.
\end{cor}

\subsection{Adjunctions}\label{moradj}

Let $F$ and $G$ be two adjoint functors as above. For any object $Y \in \caB$, we have the morphism
$\eta_Y \in \Hom_\caB(Y,FG(Y))$ corresponding to $\Id_{G(Y)}$. Similarly, for any $X$ in $\caA$, we have
a morphism $\epsilon_X  \in \Hom_\caA(GF(X),X)$ corresponding to $\Id_{F(X)}$. We obtain natural
transformations $\eta$ and $\epsilon$ respectively from $\Id_\caB$ to $FG$ and from $GF$ to $\Id_\caA$.
We call $\eta$ and $\epsilon$ the {\sl adjunction morphisms}\indexter{adjunction morphisms}, the first
sometimes being called the {\sl unit} and the second the {\sl counit}. These morphisms moreover satisfy
that the compositions
\begin{align}\label{moad1}
 F(X) \stackrel{ \eta_{F(X)} }{\longrightarrow} FGF(X)
\stackrel{F(\epsilon_X)} {\longrightarrow}F(X)
\end{align}
\begin{align}\label{moad2}
 G(Y) \stackrel { G(\eta_{Y})
  }{\longrightarrow} GFG(Y)   \stackrel{\epsilon_{G(Y)}}{\longrightarrow}  G(Y)
\end{align}
are equal to the respective identities of $F(X)$ and $G(Y)$.

Conversely, suppose that $\eta$ and $\epsilon$ are natural transformations respectively from $\Id_\caB$
to $FG$ and from $GF$ to $\Id_\caA$. Then we obtain for all $X \in \caA$ and $Y \in \caB$,
\begin{align*}
 \theta_{Y,X}\colon \Hom_\caA(G(Y),X) &\rightarrow \Hom_\caB(Y,F(X))\\
f&\mapsto [Y \stackrel{\eta_Y}{\rightarrow}FG(Y)\stackrel{F(f)}{\rightarrow} F(X)  ]
\end{align*}
and in the opposite direction
\begin{align*}
 \sigma_{Y,X}\colon \Hom_\caB(Y,F(X))&\rightarrow
  \Hom_\caA(G(Y),X)\\
g& \mapsto [G(Y)
\stackrel{G(g)}{\rightarrow}GF(X)\stackrel{\epsilon_X}{\rightarrow} X.  ]
\end{align*}
We easily deduce from the naturality of $\eta$ and $\epsilon$ that $\theta$ and $\sigma$ are natural
transformations respectively from $L_G$ to $R_F$ and from $R_F$ to $L_G$. On the other hand, if $\eta$
and $\epsilon$ satisfy (\ref{moad1}) and (\ref{moad2}), then $\theta$ and $\sigma$ are mutually inverse.

We call an {\sl adjunction} \indexter{adjunction} a quadruple $(F,G, \epsilon,\eta)$ where

- $F \colon \caA \rightarrow \caB$, $G \colon \caB \rightarrow \caA$ are two functors.

- $\eta \colon   \Id_\caB \rightarrow FG$ and $\epsilon\colon GF \rightarrow \Id_\caA$ are two natural
transformations satisfying (\ref{moad1}) and (\ref{moad2}).

\bigskip

The following result is used very often.

\begin{thm}
Let $F \colon  \caA \rightarrow \caB$, $G \colon \caB \rightarrow \caA$ be two functors, $G$ being the
left adjoint of $F$. Then $F$ preserves limits and $G$ colimits.
\end{thm}

\begin{proof}
We will now show the first assertion, the second being established in the same way. Let
$D \colon  \caI \rightarrow \caA$ be as in \ref{limites}, and suppose that $\varprojlim D$ exists. We
want to show that $\varprojlim (F\circ D)$ exists and is equal to $F(\varprojlim D)$. We use Proposition
\ref{egalitecat}. For all $Y \in \caB$, we have
\begin{align*}
& \Hom_\caB(Y,F(\varprojlim D))=\Hom_\caA(G(Y),\varprojlim D)\simeq
\Hom_{\caF_{\caI,\caA}}(\Delta_{G(Y)}, D) \\
\simeq &\Hom_{\caF_{\caI,\caB}}(\Delta_{Y}, F(D))=\Hom_{\caB}(Y,\varprojlim F\circ D),
\end{align*}
We have successively used the adjunction of $F$ and $G$, the natural isomorphism (\ref{limHom}), the fact
that any natural transformation $\theta \colon  \Delta_{G(Y)}\rightarrow D$ consists, for all
$i \in \caI$, of a morphism $G(Y)\rightarrow D(i)$ (satisfying certain compatibility conditions), and
that by the adjunction of $F$ and $G$, this is equivalent to the data of morphisms
$Y\rightarrow F(D(i))$, defining a natural transformation from $\Delta_Y$ to $F\circ D$. Finally, we use
the universal property of the limit. This shows that $F(\varprojlim D)$ is indeed the limit of
$F\circ D$.
\end{proof}

\begin{exemple}
To illustrate the notion of adjoint functor, we choose a simple example which is a special case of the
constructions of Section \ref{AoubliB}. Let $A$ be a unital $\bbC$-algebra and let $B$ be a subalgebra of
$A$. Any left unital $A$-module $M$ is also a $B$-module, and we obtain a forgetful functor:
\[
\caF \colon  \caM(A) \rightarrow \caM(B).
\]
A left adjoint of the forgetful functor, the base change functor, is given by
\[
P_B^A \colon  N \mapsto A\otimes_B N.
\]

To show this assertion, we must exhibit a natural isomorphism:
\[
\theta_{X,Y}\colon \Hom_A(A\otimes_B X,Y) \rightarrow  \Hom_B(X,\caF(Y)).
\]
If $f \colon  A\otimes_B X\rightarrow Y$ is an $A$-morphism, set
\[
\theta_{X,Y}(f)(x)=f(1\otimes x),\quad (x \in X).
\]
This indeed defines a $B$-morphism because:
\begin{align*}
\theta_{X,Y}(f)(b\cdot x)&=f(1\otimes b\cdot x)=f(b\otimes
 x)=f(b\cdot (1\otimes x))\\
=&b\cdot f(1\otimes x)=b\cdot(\theta_{X,Y}(f)(x))
\end{align*}

The inverse of $\theta$ is given by $\sigma$:
\[
\sigma_{X,Y} \colon \Hom_B(X,\caF(Y))\rightarrow  \Hom_A(A\otimes_B X,Y)
\]
defined for all $g \in \Hom_B(X,\caF(Y))$ by
\[
\sigma_{X,Y}(g)(a\otimes x)=a\cdot g(x).
\]
This is indeed an $A$-morphism because:
\[
\sigma_{X,Y}(g)(a'\cdot(a\otimes x))=\sigma_{X,Y}(g)(a'a\otimes x)=
a'a\cdot g(x)=a'\cdot(\sigma_{X,Y}(g)(a\otimes x)).
\]
We leave it to the reader to verify the naturality of $\theta$ and $\sigma$.

We check that $\theta$ and $\sigma$ are inverse to each other:
\[
\sigma_{X,Y}(\theta_{X,Y}(f))(a\otimes x)=a\cdot
(\theta_{X,Y}(f)(x))=a\cdot (f(1\otimes x))=f(a\otimes x).
\]
and
\[
\theta_{X,Y}( \sigma_{X,Y}(g))(x)=  (\sigma_{X,Y}(g))(1\otimes x)=g(x).
\]

We have therefore proved the adjunction of the forgetful and base change functors using the first
definition of adjoint functors. When we look more closely at this proof, we realize that the second point
of view is underlying. Indeed, the key point is the existence of a $B$-morphism
\[
\eta_X\colon x \mapsto 1\otimes x, \quad X \rightarrow A\otimes_B X
\]
and of an $A$-morphism
\[
\epsilon_Y\colon a\otimes y \mapsto a\cdot y, \quad   A\otimes_B Y \rightarrow Y
\]
These are adjunction morphisms and the fact that $\theta$ and $\sigma$ are inverse to each other comes
from the fact that $\eta$ and $\epsilon$ satisfy (\ref{moad1}) and (\ref{moad2}), i.e., in this case:
\[
x \mapsto 1\otimes x \mapsto 1\cdot x=x
\]
and
\[
a\otimes x \mapsto a\otimes 1\otimes  x \mapsto a\cdot (1\otimes x)=  a\otimes x.
\]
\end{exemple}

\section{Abelian categories} \label{catab}

In this section, we recall the fundamental results concerning abelian categories. Our references on the
subject will be the books \cite{Freyd} and \cite{Par}, to which the reader is invited to refer for
further details.

\subsection{Axiomatics of abelian categories}

The notion of an abelian category axiomatizes the fundamental properties of the category $A-\mathbf{mod}$
of left modules over a ring $A$. Any full subcategory of a category $A-\mathbf{mod}$ stable under passing
to submodules and quotients is again an abelian category, for example:

- the category $\caM(A)$ of left unital modules over a unital ring $A$,

- the category of sheaves (or presheaves) of abelian groups on a topological space is also abelian, but
it is not an example of the previous type (the quotient presheaf of a sheaf is not a sheaf in general).

Let us now give the axioms defining abelian categories. Those given here are not the most economical. Let
$\caA$ be a category. We say that $\caA$ is an abelian category if it satisfies the following
axioms\footnote{The numbering of the axioms adopted here, from $\mathbf{AB0}$ to $\mathbf{AB6}$, does not
coincide with that of Grothendieck in ``Tohoku''; for example, our axiom $\mathbf{AB6}$ expresses that
the category is locally finitely generated.}.

\noindent $\mathbf{AB0.}$ There exists a zero object in $\caA$, denoted $0$.

\medskip

\noindent $\mathbf{AB1.}$ The category $\caA$ admits finite products.

\noindent $\mathbf{AB1^*.}$ The category $\caA$ admits finite coproducts.

\medskip

\noindent $\mathbf{AB2.}$ Every morphism in $\caA$ admits a kernel.

\noindent $\mathbf{AB2^*.}$ Every morphism in $\caA$ admits a cokernel.

\medskip

\noindent $\mathbf{AB3.}$ Every monomorphism in $\caA$ is a kernel.

\noindent $\mathbf{AB3^*.}$ Every epimorphism in $\caA$ is a cokernel.

\medskip

The custom in abelian categories is to use the terminology ``direct sum'' and the notation $\bigoplus$
rather than ``coproduct'' and $\coprod$.

\subsection{Properties of abelian categories} \label{abcat}

We deduce many consequences from these axioms. Here are some of the most important ones.

--- 1. A morphism in $\caA$ is an isomorphism if and only if it is both a monomorphism and an
epimorphism.

--- 2. Let $A$ be an object of $\caA$, and let $S_A$ be the set of its subobjects, $Q_A$ the set of its
quotients. Then the maps
\[
\coker \colon  S_A \rightarrow Q_A,\quad \ker \colon  Q_A \rightarrow S_A
\]
which, respectively, associate to any subobject $B\rightarrow A$ its cokernel, and to any quotient
$A \rightarrow C$ its kernel, are inverse to each other. Moreover, they reverse the order on $S_A$ and
$Q_A$. We denote by $A/B$ the quotient $\coker(B\rightarrow A)$ for any subobject $B\rightarrow A$ of
$A$.

--- 3. Any two subobjects of an object $A$ admit an intersection and a union. As for coproducts, it is
often more expressive to change the terminology for unions in abelian categories by speaking instead of
sums. We then use the notation $+$ or $\sum$.

\bigskip

We define the {\sl image} of a morphism $f \colon A \rightarrow B$ in $\caA$ as the smallest subobject
$C \rightarrow B$ of $B$ such that $f$ factors through $C \rightarrow B$. Dually, we define the {\sl
coimage} of $f$ as the smallest quotient $A \rightarrow D$ of $A$ which factors $f$.

--- 4. Every morphism $f \colon  A \rightarrow B$ admits an image (resp. a coimage) which is identified
with $\ker(\coker f)$ (resp. with $\coker(\ker f)$).

--- 5. (Factorization of morphisms) Every morphism $A \rightarrow B$ in $\caA$ factors uniquely up to
isomorphism into $A \rightarrow I \rightarrow B$, where $A \rightarrow I $ is an epimorphism, and
$I \rightarrow B$ is a monomorphism. Moreover $A \rightarrow I $ is the coimage of $A \rightarrow B$, and
$I \rightarrow B$ is its image.

The existence of kernel and image for all morphisms in $\caA$ allows us to define the notion of an {\sl
exact sequence} of morphisms in $\caA$. We assume the reader is sufficiently familiar with this notion so
as not to have to redefine it.

--- 6. For any pair $(A,B)$ of objects of $\caA$, the set $\Hom_\caA(A,B)$ is equipped with an abelian
group structure, whose neutral element is the distinguished element $0$ whose existence is ensured by
$\mathbf{AB0}$. The operations of composition of morphisms are distributive with respect to addition on
the $\Hom_\caA(A,B)$.

--- 7. For any pair $(A,B)$ of objects of $\caA$, the product $A\times B $ and the sum $A \oplus B$ whose
existence is ensured by the axioms $\mathbf{AB1}$, $\mathbf{AB1^*}$ are naturally isomorphic.

--- 8. {\sl $3\times 3$ Lemma.} Consider the following diagram in $\caA$
\[
\xymatrix{
& 0\ar[d]  & 0\ar[d] &  0\ar[d] &   \\
0\ar[r] & A \ar[d]\ar[r]&  A'\ar[d]\ar[r]&  A''\ar[d]\ar[r]& 0 \\
0\ar[r]& B\ar[r]\ar[d]& B'\ar[r]\ar[d] & B''\ar[r]\ar[d] &  0\\
& C \ar[d]&  C' \ar[d]&  C'' \ar[d]&  \\
&  0&  0& 0 &
}
\]
where the columns and rows are exact. Then we can complete this diagram with a monomorphism
$C\rightarrow C'$ and an epimorphism $C' \rightarrow C''$ so that the enriched diagram is still
commutative and the sequence $0\rightarrow C\rightarrow C' \rightarrow C''\rightarrow 0$ is exact.

--- 9. {\sl Noether isomorphisms.} From the $3\times 3$ lemma, we can easily deduce the two Noether
isomorphism theorems:

$(i)$ given an object $C$ in $\caA$ and two subobjects $A$ and $B$ of $C$ with $A \subseteq B$, then
$B/A$ is a subobject of $C/A$ and $(C/A)/ (B/A)$ is isomorphic to $C/B$,

$(ii)$ given an object $C$ in $\caA$ and two subobjects $A$ and $B$ of $C$, we have
\[
B/(A \cap B) \simeq (A + B)/A.
\]

Assuming moreover that $A \cap B =0$, and $A +  B=C$, $(ii)$ gives $C \simeq A \oplus B$.

\bigskip

The existence of limits is of course an important property of categories. This leads us to consider the
following axioms:

\bigskip

\noindent $\mathbf{AB4.}$ Every family of objects admits a product in $\caA$.

\noindent $\mathbf{AB4^*.}$ Every family of objects admits a direct sum in $\caA$.
\bigskip

An abelian category satisfying $\mathbf{AB4.}$ is said to be complete, \index[ter]{complete} while an
abelian category satisfying $\mathbf{AB4^*.}$ is said to be cocomplete. \index[ter]{cocomplete}

\noindent $\mathbf{AB5.}$ $\caA$ is cocomplete, and for any increasing directed system $M_i$ of
subobjects of $M$, for any subobject $N$ of $M$, we have
\[
N\cap \left( \sum_i M_i  \right)=  \sum_i (N\cap M_i).
\]
Note that the sums exist by $\mathbf{AB4^*}$.

A Grothendieck category \indexter{Grothendieck category} is a cocomplete abelian category satisfying
$\mathbf{AB5}$ (we do not require here the existence of a generator).

Let us continue our list of properties of abelian categories.

--- 10. Suppose that $\caA$ is complete (resp. cocomplete). Then any diagram
$D\colon \caI \rightarrow \caA$ admits a limit (resp. colimit). In particular, $\caA$ admits inverse
(resp. direct) limits.

\subsection{Jordan-Hölder series and composition factors} \label{JoHo}
\index[ter]{Jordan-Hölder}

Let $\caA$ be an abelian category. A non-zero object $M$ of $\caA$ is said to be {\sl simple} if every
subobject of $M$ \index[ter]{simple} is either the zero object $0$, or $M$. We denote by
$\mathbf{Irr}(\caA)$ \index[not]{Irr(A)@$\mathbf{Irr}(\caA)$} the set of equivalence classes of simple
objects of $\caA$. Suppose that we have an increasing sequence of subobjects of $M$, all distinct
\[
0=B_0 \subset B_1 \subset \cdots \subset B_n =M.
\]
Such a sequence is called a {\sl composition series}, or {\sl Jordan-Hölder series} if the successive
quotients $B_i/B_{i-1}$, $i=1,\ldots, n$ are simple. The integer $n$ is called the {\sl length} of this
composition series and the quotients $B_i/B_{i-1}$ are called {\sl composition factors} of $M$. The {\sl
Jordan-Hölder theorem} then asserts that if the object $M$ admits a composition series, then all
composition series are of the same length and their factors are isomorphic, up to permutation. An object
admitting a composition series is said to be of finite length, its length being that of one of its
composition series. The Jordan-Hölder theorem is proved from the Noether isomorphisms \ref{abcat}, 9.

To obtain the results on the series and composition factors of an object generalizing those of the
category of $A$-modules, we must make some additional hypotheses on $\caA$. We therefore assume that
$\caA$ is a Grothendieck category. It therefore satisfies $\mathbf{AB4^*}$ and $\mathbf{AB5}$. We assume
moreover:

\noindent $\mathbf{AB6}$: every object of $\caA$ is a sum of finitely generated subobjects.

\begin{exemples}
Note that the categories $A-\mathbf{mod}$ (respectively $\caM(A)$), $A$ a ring (respectively $A$ a unital
ring) satisfy all these axioms.
\end{exemples}

\bigskip

The first result we will establish uses Zorn's lemma. Zorn's lemma can be considered as an axiom of set
theory. It is equivalent to the axiom of choice. We say that an ordered set $S$ is inductively ordered if
any non-empty totally ordered subset of $S$ admits an upper bound. We refer the reader to \cite{Lang} or
\cite{Douady} for the definitions of these terms. Zorn's lemma asserts that a non-empty inductively
ordered set admits a maximal element.

\begin{prop}
$(i)$ Let $M$ be a non-zero finitely generated object of $\caA$. Then $M$ admits a simple quotient.

$(ii)$ If $M$ is not finitely generated, there exist subobjects $N_1 \subset N_2$ of $M$ such that the
quotient $N_2/N_1$ is simple.
\end{prop}

\begin{proof}
It is clear that $(ii)$ follows from $(i)$ by taking for $N_2$ a non-zero finitely generated subobject in
$M$, whose existence is ensured by $\mathbf{AB6}$. To prove $(i)$, let us introduce the (ordered) set $S$
of proper subobjects of $M$. To be able to apply Zorn's lemma, we now show that $S$ is inductive, i.e.,
that if $(N_i)_i$ is a chain in $S$ (a totally ordered subset) then $N:=\bigcup_i N_i $ is in $S$. We
must show that $N$ is a proper subobject. Suppose the contrary, i.e., $M=\bigcup_i N_i $. Since $M$ is
finitely generated, there exists a certain index $i_0$ such that $M=N_{i_0}$, which is absurd since
$N_{i_0}$ is a proper subobject. The assertion is therefore proved. By Zorn's lemma, $S$ admits a maximal
element $M_0$, and by maximality, the quotient $M/M_0$ is simple.
\end{proof}

\bigskip

Let $\mathbf{JH}(M)$ denote the set of simple subquotients of an object $M$. We then have:

\begin{lemme}
$(i)$ If $M' \subset M$ is a subobject, then $$\mathbf{JH}(M) = \mathbf{JH}(M') \cup \mathbf{JH}(M/M').$$

$(ii)$ $\mathbf{JH}(M)=\emptyset$ if and only if $M=0$.

$(iii)$ If $(M_i)_{i\in I}$ is a family of subobjects of $M$,
\[
\mathbf{JH}(\sum_{i \in I} M_i)= \bigcup_{i \in I}  \mathbf{JH}(M_i) .
\]
\end{lemme}

\begin{proof}
$(i)$ is clear. $(ii)$ follows from the previous proposition. For $(iii)$, the assertion follows from
$(i)$ by induction if $I$ is finite. By replacing the $M_i$ by a family consisting of sums of $M_i$, we
reduce to the case where $(M_i)_{i\in I}$ is an increasing directed system. Let $Q=M'/M''$ be a simple
subquotient of $\sum_{i \in I} M_i$. Suppose that for all $i \in I$, $Q\notin \mathbf{JH}(M_i)$. Then for
all $i \in I$,
\[
M' \cap (M''+M_i)=M''.
\]
Thus, using $\mathbf{AB5}$, we obtain
\[
M'=M' \cap (\sum_i  (M''+M_i) )= \sum_i  M' \cap (M''+M_i)=M'',
\]
which contradicts the fact that $Q$ is simple.
\end{proof}

\section{Semisimplicity}\label{semisimplicite}
\label{AVII}

\begin{lemme}\label{lemmesemsimP}
Let $\caA$ be an abelian category satisfying the hypotheses of the previous section, and let $E$ be an
object of $\caA$. The following conditions are equivalent:

$(i)$ $E$ is isomorphic to a direct sum of simple objects (i.e., $E$ is {\sl completely
reducible}).\indexter{completely reducible}

$(ii)$ $E$ is the sum of its simple subobjects.

$(iii)$ For any subobject $E'$ of $E$, there exists a subobject $E''$ of $E$ such that $E=E' \oplus E''$.
\end{lemme}
\begin{proof}
It is clear that $(i)$ implies $(ii)$. Suppose $(ii)$; we now show $(iii)$. Let $E'$ be a subobject of
$E$. We can, by hypothesis, write $E$ as a sum (but not direct) of simple subobjects,
\[
E  = \sum_{i \in I} E_i
\]
for a certain index set $I$. Let us take a subobject of $E$ of the form $\sum_{j\in J} E_j$,
$J \subset I$, maximal with respect to inclusion and such that the sum $F=E' + \sum_{j\in J} E_j$ is
direct (the existence of such a $J$ is ensured by Zorn's lemma). Then this sum is equal to $E$. Indeed,
it suffices to see that each $E_i$, $i\in I$, is in $E' \oplus (\bigoplus_{j\in J} E_j)$. Since the
intersection of $F$ with $E_i$ is a subobject of $E_i$, this intersection is $E_i$ or $0$ since $E_i$ is
simple. If this intersection were zero, we could adjoin $i$ to $J$ which contradicts the maximality of
$J$. Thus $E_i \subset F$. Let us now show that $(iii)$ implies $(i)$. We start by showing that if
$E\neq 0$, $E$ then admits a simple subobject. Let $E'$ be a non-zero finitely generated subobject of
$E$, whose existence is ensured by the hypothesis $\mathbf{AB6}$. Consider the set of proper subobjects
of $E'$, ordered by inclusion. As in the proof of Proposition \ref{JoHo}, we see that this set is
inductively ordered, and thus that it admits by Zorn's lemma a maximal element $M$. This $M$ is a
subobject of $E$. By hypothesis, we can write $E=M\oplus H$, for a certain subobject $H$. We then have
\[
E' =M\oplus (E'\cap H).
\]
Since $M$ is maximal in $E'$, $E' \cap H$ is simple. This shows that $E$ admits a simple subobject. Now
let $E_0$ be a maximal direct sum of simple subobjects of $E$. If $E\neq E_0$, we write $E=E_0\oplus F$.
We now apply the previous remark to $F$, which still satisfies $(iii)$ (if $F'\subset F$ and
$E=F'\oplus G$, then $F=F'\oplus(G\cap F)$): there exists a simple subobject of $F$. This contradicts the
definition of $E_0$.
\end{proof}

\begin{cor}
A subobject of a semisimple object is semisimple.
\end{cor}

\begin{proof}
Suppose that $F$ is a subobject of $E$ with $E$ semisimple, and let $F'$ be a subobject of $F$. It is
also a subobject of $E$, and thus there exists a subobject $G$ of $E$ such that $E=F'\oplus G$. Set
$G'=G \cap F$. We then have $F=F'\oplus G'$.
\end{proof}

\section{Notable functors}

By notable functors, we mean functors preserving properties of the objects or morphisms to which they
apply. For example, we defined above faithful functors and full functors. We give other examples.

\subsection{Exact functors}\label{exaCt}
\index[ter]{exact (functor)}

Let $F \colon \caA \rightarrow \caB$ be a functor between two abelian categories. We say that $F$ is {\sl
additive} if it respects the group structure of the $\Hom$, i.e., if for any pair $(X,Y)$ of objects of
$\caA$, $$F \colon  \Hom_\caA(X,Y)\rightarrow  \Hom_\caB(F(X),F(Y))$$ is a morphism of abelian groups.
Functors between abelian categories will implicitly always be assumed to be additive. In what follows,
the categories are abelian and the functors are additive.

We say that it is {\sl left exact} (resp. {\sl right exact}) if for any exact sequence in $\caA$ of the
form
\[
0\rightarrow X \stackrel{f}{\rightarrow} Y\stackrel{g}{\rightarrow} Z,
\]
the sequence
\[
0\rightarrow F(X)  \stackrel{F(f)}{ \rightarrow} F(Y)  \stackrel{F(g)}{ \rightarrow} F(Z),
\]
is exact (resp. if for any exact sequence of the form
\[
X \stackrel{f}{\rightarrow} Y\stackrel{g}{\rightarrow} Z
\rightarrow 0,
\]
the sequence
\[
F(X)  \stackrel{F(f)}{ \rightarrow} F(Y)  \stackrel{F(g)}{ \rightarrow} F(Z)\rightarrow 0,
\]
is exact). If $F$ is right and left exact, we simply say that it is {\sl exact}. For a contravariant
functor $F \colon \caA \rightarrow \caB $, the convention for left or right exactness is the following.
We consider $F$ as a covariant functor $ \caA^{op} \rightarrow \caB $ and we transfer the definitions
above.

One can show that a left exact functor $F  \colon \caA \rightarrow \caB$ preserves finite limits (i.e.,
for any functor $D \colon  \caI \rightarrow \caA$, where the index category $\caI$ is finite, admitting a
limit $X$, then $F(X)$ is the limit of the functor $F\circ D$). Similarly a right exact functor preserves
finite colimits.

Recall that a functor $F \colon  \scrC \rightarrow \scrD$ is faithful if for any pair of objects $(X,Y)$
of $\scrC$, the map
\[
\Hom_\scrC(X,Y)\stackrel{F}{\longrightarrow} \Hom_\scrD(F(X),F(Y))
\]
is injective. If $\scrC$ and $\scrD$ are abelian categories and the functor $F$ is additive, $F$ is
faithful if for all $X\stackrel{f}{\rightarrow} Y$, $F(f) = 0$ implies that $f=0$. In an abelian category
every morphism $X\stackrel{f}{\rightarrow} Y$ factors uniquely into
\[
X\rightarrow I \rightarrow Y,
\]
where $X\rightarrow I$ is an epimorphism and $I \rightarrow Y$ a monomorphism (see property 5 of abelian
categories in Appendix \ref{abcat}). On the other hand, if $f\neq 0$, then $I\neq 0$. If $F$ is exact, it
preserves epimorphisms and monomorphisms, and we deduce that
\[
F(X)\rightarrow F(I) \rightarrow F(Y)
\]
is the factorization of $F(f)$. If $F(f)=0$, we then have $F(I)=0$. We deduce from this brief discussion
the following result.

\begin{lemme}
If $F \colon \scrC \rightarrow \scrD$ is an exact functor, and $F(I)=0$ implies that $I=0$, then $F$ is
faithful.
\end{lemme}

\subsection{Progenerators}\label{rmqfonct}

Let $X$ be an object of the abelian category $\caA$. We define the functors:
\begin{align*}
 \Hom_\caA(X,\bullet):\qquad  \caA &\rightarrow \bbZ-\mathbf{mod},\\
                                 Y&\mapsto \Hom_\caA(X,Y),\\
 [Y\stackrel{f}{\rightarrow}Z]&\mapsto
[\Hom_\caA(X,Y) \stackrel {f \circ \bullet}{\longrightarrow} \Hom_\caA(X,Z) ]
\end{align*}

\begin{align*}
 \Hom_\caA(\bullet ,X):\qquad  \caA &\rightarrow \bbZ- \mathbf{mod},\\
                                 Y&\mapsto \Hom_\caA(Y,X),\\
 [Y\stackrel{f}{\rightarrow}Z]&\mapsto
[\Hom_\caA(Z,X) \stackrel {\bullet \circ f }{\longrightarrow} \Hom_\caA(Y,X) ]
\end{align*}

The first of these functors is covariant, the second contravariant. They are both left exact. An object
$X$ of $\caA$ is said to be {\sl projective} \indexter{projective} if $\Hom_\caA(X,\bullet)$ is exact,
{\sl injective} \indexter{injective} if $\Hom_\caA(\bullet,X)$ is exact. It is a {\sl
generator}\indexter{generator} of $\caA$ (resp. a {\sl cogenerator}\indexter{cogenerator}) if
$\Hom_\caA(X,\bullet)$ (resp. $\Hom_\caA(\bullet,X)$) is faithful. A {\sl
progenerator}\indexter{progenerator} is a projective generator. A simple criterion is that a projective
object $X$ is a progenerator in $\caA$ if and only if for any non-zero object $Y$ in $\caA$,
$\Hom_\caA(X,Y)$ is non-trivial. This condition is sufficient by virtue of Lemma \ref{exaCt}.

Let $\caA$ be a cocomplete abelian category. A projective object $P$ in $\caA$ is said to be {\sl small}
\indexter{small (object)} if the functor $\Hom_\caA(P,\bullet)$ preserves direct sums. There exist
several equivalent characterizations of this notion. In particular, $P$ is a small projective if for any
direct sum $M=\bigoplus_{i\in I}M_i$, and for any morphism $\phi \colon  P \rightarrow M$, there exists a
finite set $J\subset I$ such that $\phi$ factors into
$P \rightarrow \bigoplus_{j\in J}M_j \rightarrow M$. In a Grothendieck category, an object is a small
projective if and only if it is finitely generated projective ({\sl cf.} \cite{Par}, 4.11, Lemma 1). The
following result gives us a criterion for a cocomplete abelian category to be equivalent to the category
of right modules over a certain ring $\caR$ (cf. \cite{Par}, 4.11, Theorem 1).

\begin{thm}
Let $\caA$ be an abelian category, cocomplete, and admitting a small progenerator $P$. Then $\caA$ is
naturally equivalent to the category of unital right modules over the ring $\caR_P=\End_\caA(P)$, the
equivalence being given by the functor $\Hom_\caA(P,\bullet)$.
\end{thm}

\begin{rmq}
We of course have the converse, since a unital ring $A$ is a small progenerator of $\caM(A)$.
\end{rmq}

We also note the following result.
\begin{lemme}
Let $\caA$ be a cocomplete abelian category, admitting a small progenerator $P$. Then any object $M$ of
$\caA$ is the cokernel of a morphism between two objects isomorphic to direct sums of objects isomorphic
to $P$.
\end{lemme}
\begin{proof}
Consider the morphism
\[
\bigoplus_{f\in \Hom_\caA(P,M)}P \rightarrow M
\]
whose restriction to the factor of index $f$ is $f$. It is an epimorphism: otherwise, its cokernel $Q$
would be non-zero, and since $P$ is a generator, there would exist a non-zero morphism
$g\colon P\rightarrow Q$; by projectivity of $P$, $g$ lifts to a morphism $f\colon P\rightarrow M$, whose
composite with $M\rightarrow Q$ is zero by construction, which is absurd. Consider the kernel $C$ of this
morphism. The same argument gives us the existence of an epimorphism of the form
$\bigoplus_{j\in J}P\rightarrow C$, from which we deduce an exact sequence
\[
\bigoplus_{j\in J}P  \rightarrow   \bigoplus_{f\in \Hom_\caA(P,M)}P \rightarrow M
\rightarrow 0.
\]
\end{proof}

\section{Decompositions of categories}\label{deccat}

Let $\caM$ be a Grothendieck category ($\mathbf{AB4^*}$ and $\mathbf{AB5}$) also satisfying
$\mathbf{AB6}$ and let $(\caM_i)_{i\in I}$ be a family of full subcategories.

We write $\mathcal{M} = \prod_{i \in I} \mathcal{M}_i$ if the full subcategories $\mathcal{M}_i$ are
mutually orthogonal (i.e., for $i \neq j$, $\Hom_{\mathcal{M}}(X_i, X_j) = \{0\}$ for any
$X_i \in \mathcal{M}_i$ and $X_j \in \mathcal{M}_j$) and if every object $M \in \mathcal{M}$ can be
uniquely written as a direct sum $M = \bigoplus_{i \in I} M_i$ with $M_i \in \mathcal{M}_i$. Under these
conditions, the functor $(M_i)_{i \in I} \mapsto \bigoplus_{i \in I} M_i$ realizes an equivalence of
categories between the categorical product $\prod_{i \in I} \mathcal{M}_i$ and $\mathcal{M}$.

Suppose that $\caM =\prod_{i\in I}\caM_i $. Then $\mathbf{Irr}(\caM)$ is a disjoint union
\[
\mathbf{Irr}(\caM)=\coprod_{i \in I} \mathbf{Irr}(\caM_i).
\]

In the other direction, suppose now that $S$ is a subset of $ \mathbf{Irr}(\caM)$, and let $\caM_S$
denote the full subcategory of $\caM$ whose objects are those all of whose irreducible subquotients are
in $S$. Then $\caM_S$ is an abelian subcategory, stable under passing to subquotients, extensions and
colimits. For any object $E$ of $\caM$, we denote by $E_S$ the sum of all the subobjects of $E$ which are
in $\caM_S$. Then it is clear that $E_S \in \caM_S$. Let $S'$ be another subset of $ \mathbf{Irr}(\caM)$,
such that $S\cap S'=\emptyset$. If $E \in \caM_S$ and $E' \in \caM_{S'}$, we have
$\Hom_\caM(E,E')=\{0\}$. Moreover, for all $E \in \caM$, $E_S\cap E_{S'}=\{0\}$, and thus
\[
E_S \oplus E_{S'} \subset E.
\]

\begin{defi}
We say that the subset $S \subset \mathbf{Irr}(\caM)$ splits an object $E$ in $\caM$, if
$E=E_S\oplus E_{\bar S}$, where $\bar S= \mathbf{Irr}(\caM) \setminus S$. We say that $S$ splits the
category $\caM$ if $S$ splits every object of $\caM$. More generally, suppose that we have a
decomposition into a disjoint union:
\[
\mathbf{Irr}(\caM) =\coprod_{\alpha} S_\alpha.
\]
We say that this decomposition $\{S_\alpha\}$ splits an object $E$ in $\caM$, if
$E=\bigoplus_\alpha  E_{S_\alpha}$ and that $\{S_\alpha\}$ splits the category $\caM$ if $\{S_\alpha\}$
splits every object of $\caM$.
\end{defi}

\begin{lemme}
Suppose that $\mathbf{Irr}(\caM) =\coprod_{\alpha} S_\alpha$ splits an object $E$ in $\caM$. Then any
subquotient of $E$ is also split according to this decomposition.
\end{lemme}

\begin{proof}
Let us write $E=\bigoplus_{\alpha} E_\alpha$. It is clear that it suffices to show that for any subobject
$L$ of $E$, we have
\[
L=\bigoplus_{\alpha} (E_\alpha \cap L).
\]
Set $C=L/\bigoplus_{\alpha} (E_\alpha \cap L)$. Then for all $\alpha$,
\[
\mathrm{JH}(C)\subset\mathrm{JH}(L/(E_\alpha \cap L))\subset
\mathrm{JH}(E/E_\alpha)\subset \bigcup_{\beta \neq
\alpha}\mathrm{JH}(E_\beta)\subset \mathbf{Irr}(\caM)\setminus S_\alpha.
\]
This shows that $\mathrm{JH}(C)\subset \bigcap_\alpha (\mathbf{Irr}(\caM)\setminus S_\alpha)=\emptyset$,
and thus $C=\{0\}$.
\end{proof}

\section{Center of a category}\label{centrecat}

Let $\scrC$ be a category. The center \index[ter]{center!of a category} $\frz(\scrC)$ of $\scrC$ is the
set of endomorphisms (i.e., natural transformations to itself) of the identity functor. In other words,
an element $z$ of the center of $\scrC$ is the data, for any object $X$ of $\scrC$, of a morphism
$z_X:X\rightarrow X$, so that given two objects $X$ and $Y$ of $\scrC$, and a morphism
$f: X\rightarrow Y$, we have $f\circ z_X=z_Y \circ f$.

\begin{exemple}
If $\scrC=\caM(A)$, the category of unital $A$-modules over the unital ring $A$, then $\frz(\scrC)$ is
naturally identified with the center of the ring $A$.
\end{exemple}


\chapter{Amitsur's theorem and corollaries}

\section{Amitsur's theorem} \label{amitsur}
\indexter{Amitsur}

The main result of this section is due to Amitsur. It can be viewed as a non-commutative version of
Hilbert's Nullstellensatz for algebras over the field $\bbC$.

\begin{defi}
Let $A$ be a non-zero $\bbC$-algebra with identity $\una$. The spectrum of an element $a$ of $A$ is the
set
\[
\mathrm{Spec}(a) =\{ \lambda \in \bbC\, | \; a-\lambda \una \text{
is not invertible }    \}
\]
\end{defi}

We identify the subalgebra $\bbC \una$ of $A$ with $\bbC$.

\begin{thm}
Let $A$ be a non-zero $\bbC$-algebra with identity $\una$. Suppose that $A$ is of at most countable
dimension over $\bbC$. Then

$a)$ if $A$ is a division algebra, $A=\bbC$,

$b)$ for all $a \in A $, $\mathrm{Spec} (a)$ is non-empty. Moreover $\mathrm{Spec} (a) =\{ 0 \}$ if and
only if $a$ is nilpotent.
\end{thm}

\begin{proof}
We will show that $b)$ implies $a)$. Suppose there exists $a\in A$, $a \notin \bbC$. Then for all
$\lambda \in \bbC$, $a-\lambda\una$ is non-zero, hence invertible if we assume that $A$ is a division
algebra. This gives us $\mathrm{Spec} (a) =\emptyset$, which contradicts $b)$.

Let us now prove $b)$. Suppose $\mathrm{Spec} (a) =\emptyset$, i.e., that for all $\lambda \in \bbC$,
$a-\lambda \una$ is invertible. The uncountable family
\[
\{  (a-\lambda \una)^{-1}\, |  \lambda \in \bbC      \}
\]
of elements of $A$ is therefore linearly dependent, since the dimension of $A$ is at most countable. This
shows that there exist a finite number of elements $\mu_i \in \bbC$, which we can assume to be non-zero,
and the same number of elements $\lambda_i$ in $\bbC$ such that
\[
\sum_i  \mu_i \;    (a-\lambda_i \una)^{-1} =0.
\]

Multiplying by $\prod_i (a-\lambda_i \una)$, we obtain a polynomial $P\in \bbC[t]$ of strictly positive
degree such that $P(a)=0$ (assuming the $\lambda_i$ distinct, $P=\sum_i\mu_i\prod_{j\neq i}(t-\lambda_j)$
is non-zero, since $P(\lambda_i)\neq 0$; it is not constant, since $P(a)=0$). This polynomial $P$ admits
a factorization over $\bbC$ of the form
\[
P(t)=c(t-\alpha_1)\ldots(t-\alpha_r),
\]
$c\in \bbC^\times$, $\alpha_1, \ldots ,\alpha_r \in \bbC$. This implies that
\[
(a-\alpha_1 \una )\ldots (a-\alpha_r \una)=0
\]
which contradicts the fact that all the $a-\alpha_j \una$ are invertible. This proves that
$\mathrm{Spec} (a) $ is non-empty.

Let us now prove the assertion concerning nilpotent elements. Let $a \in A$ and $n \in \bbN$ such that
$a^n=0$. Then $a$ cannot be invertible, which shows that $0  \in \mathrm{Spec} (a) $. On the other hand,
if $\lambda \in \bbC^\times$, $a-\lambda \una$ is invertible, with inverse
\[
(a-\lambda \una)^{-1}= -\lambda^{-1}(\una- \lambda^{-1}a)^{-1}=
-\lambda^{-1}\sum_{i=0}^n  (\lambda^{-1}a)^i.
\]
This shows that $\mathrm{Spec} (a) =\{ 0 \}$. Conversely, if $\mathrm{Spec} (a)=\{ 0 \}$, since
$\bbC \setminus \{0\}$ is still an uncountable family, the above argument shows that there exists a
non-constant polynomial $P \in \bbC[t]$ such that $P(a)=0$. This polynomial $P$ admits a factorization
over $\bbC$ of the form
\[
P(t)=c t^n(t-\alpha_1)\ldots(t-\alpha_r),
\]
$c,\alpha_1, \ldots ,\alpha_r \in \bbC^\times$. Since each $a-\alpha_j \una$ is invertible, we deduce
that $n>0$ and $a^n=0$.
\end{proof}

\section{Schur's lemma} \label{ASchur}

Let $A$ be a $\bbC$-algebra (which we do not necessarily assume to be unital). We will use the result of
the previous section to show the following theorem.

\begin{thm}
Let $M$ be a simple $A$-module and a $\bbC$-vector space. We then have

$a)$ $\End_A (M)$ is a division algebra.

$b)$ If the dimension of $M$ over $\bbC$ is at most countable, we have
\[
\End_A (M) \simeq \bbC.
\]
To do this, it suffices that the dimension of $A$ be at most countable.
\end{thm}

\begin{proof}
Note first that by replacing $A$ by $\widetilde A=A\oplus \bbC \una$ where as in Section \ref{AMODUn} we
add an identity to $A$, we can assume that $A$ is a unital algebra and $M$ a unital $A$-module.

Let $0 \neq f \in \End_A (M)$. Since $\ker f$ and $\im f$ are submodules of $M$, and $M$ is simple, we
have $\ker f=0$ and $\im f =M$. This shows that $f$ is invertible.

Let us now show that if the dimension of $A$ is at most countable, the same is true for that of $M$.
Consider $m\in M$ such that $A\cdot m \neq 0$ and
\[
A \rightarrow M, \quad a \mapsto a\cdot m.
\]
Since $M$ is simple, $A\cdot m=M$ and thus $M$ is of at most countable dimension over $\bbC$.

Similarly, let us prove $b)$: consider
\[
\End_A (M)  \rightarrow M, \quad f \mapsto f(m).
\]
Since $f(a \cdot m)=a \cdot f(m)$ and $A \cdot m=M$, we see that $f$ is determined by $f(m)$, and thus
that the morphism above is injective. We deduce that the dimension of $\End_A (M) $ is at most countable.
Amitsur's theorem then shows that $\End_A( M)\simeq \bbC $.
\end{proof}

\section{Separation lemma}\label{nilp}

Let us prove the following separation result:
\begin{lemme}
Let $A$ be an algebra over $\bbC$ of countable dimension. Let $a$ be a non-nilpotent element of $A$. Then
there exists a simple left module $M$ such that $a\cdot M \neq 0$.
\end{lemme}

\begin{proof}
By replacing $A$ with $\widetilde A=A \oplus \bbC\una$, we may assume that $A$ has a unit. Since $a$ is
not nilpotent, by Amitsur's theorem \ref{amitsur} and its proof, there exists $\lambda \in \bbC^\times$
such that $a-\lambda\una$ is not left invertible (the proof of $b)$ applies verbatim to left inverses, by
multiplying on the right by $\prod_i (a-\lambda_i \una)$). Then
\[
A \rightarrow  A,\quad b \mapsto  b(a-\lambda \una)
\]
is not surjective ($\una$ is not reached) and the module $N=A/A(a-\lambda \una)$ is therefore non-zero.
Since $N$ is finitely generated, by Proposition \ref{JoHo}, $N$ admits a simple quotient $M$. We then
have $a\cdot M\neq 0$. Indeed, if $M=N/I$, where $I$ is a submodule of $N$, and if $\widetilde I$ is the
inverse image of $I$ in $A$, then $\widetilde I$ is a left ideal of $A$ which contains
$A(a-\lambda \una)$. If moreover we assume that $a\cdot M=0$, then $\widetilde I$ contains $a \cdot A$.
But then $\widetilde I$ contains $a-(a-\lambda \una)=\lambda \una$, hence contains $\una$, and therefore
$\widetilde I=A$. Since $M \simeq A/\widetilde I$, we reach a contradiction.
\end{proof}

\chapter{Linear algebra}

\numberwithin{equation}{section}

\section{Commutative subalgebras of $\End(V)$}\label{dimborn}
We prove the following result on commutative subalgebras of $\End(V)$:

\begin{lemme}
Let $V$ be a vector space over $\bbC$ of dimension $k$ and let $R$ be a commutative subalgebra of
$\End (V)$ generated by $l$ elements $a_1,\ldots,a_l$ and the identity. Then
$\dim R \leq f_l(k)=k^{2-2^{1-l}}$.
\end{lemme}
\begin{proof}
Note that the case $l=1$ follows from the Cayley-Hamilton theorem. Consider the characteristic subspaces
of the endomorphism $a_1$. These are subspaces stable under the action of $R$. We decompose each of them
into characteristic subspaces for the action of $a_2$. We obtain a finer decomposition into subspaces
stable under the action of $R$. Let us reiterate the operation up to $a_l$. We thus decompose $V$ into a
direct sum of subspaces stable under the action of $R$, say $m$ subspaces of dimension $k_1,\ldots,k_m$,
where each $a_i$ is the sum of a scalar operator and a nilpotent operator. Since $2-2^{1-l}\geq 1$, $f_l$
is convex and $f_l(0)=0$, we have
\[
f_l(k)\geq f_l(k_1)+ \cdots +f_l(k_m)
\]
and it suffices to prove the lemma for each subspace. We can therefore assume that the $a_i$ are sums of
a scalar operator and a nilpotent operator. It is clear that we can also replace the $a_i$ by their
nilpotent parts, and thus ultimately assume them to be nilpotent.

Let $\phi_l(k)$ denote the largest possible dimension of a commutative subalgebra $R$ of $\End(V)$
generated by $l$ {\sl nilpotent} elements $a_1,\ldots,a_l$ and the identity.

Let $I$ be the ideal of $R$ generated by the $a_i$ and set $V_i=I^i\cdot V$. We thus have a sequence
\[
\{0\}=V_k\subset V_{k-1} \subset \cdots \subset V_1 \subset V_0=V
\]
of subspaces of $V$. Let $W$ be a complement of $V_1$ in $V$ and let $s$ be its dimension. Since
$I^j\cdot W$ generates $V_j \; \mathrm{mod} \, V_{j+1}$, it is clear that $R\cdot W=V$. Thus each
$a\in R$ is determined by its restriction to $W$. We deduce that $\dim R\leq sk$, whence
$\phi_l(k)\leq sk$. Let $R'$ be the subalgebra of $R$ generated by $a_2,\ldots,a_l$ and the identity, and
$R''$ the ideal generated by $a_1$, so that $R=R'+R''$. We have $\dim R' \leq \phi_{l-1}(k)$, and on the
other hand, since $a_1\cdot V \subset V_1$, $\dim R''$ is less than or equal to the dimension of the
restriction of $R$ to $V_1$ and thus less than or equal to $\phi_l(k-s)$. Since $\phi_l$ is increasing,
we have $\phi_l(k-s)\leq \phi_l(\lfloor k-\phi_l(k)/k \rfloor)$ ($\lfloor x \rfloor$ denotes the integer
part of the real number $x$) whence
\[
\phi_l(k)\leq \phi_l(\lfloor k-\phi_l(k)/k\rfloor)+\phi_{l-1}(k).
\]

We will now conclude by induction on $l$ and $k$, noting that for all $l$, the case $k=1$ is trivial
($\phi_l(1)=1$ for all $l$). Suppose therefore by contradiction that $\phi_l(k) > f_l(k)$ for a certain
$l\geq 1$ and a certain $k \geq 2$, assuming the inequality of the lemma established for any pair
$(l',k')$ with $l'<l$, or $l'=l$ and $k'<k$. We then have, by the inequality proved above
\[
f_l(k)< \phi_l(k) \leq \phi_l(\lfloor k-\phi_l(k)/k\rfloor)+\phi_{l-1}(k),
\]
and by the induction hypothesis
\[
f_l(k)< f_l(\lfloor k-\phi_l(k)/k\rfloor)+f_{l-1}(k).
\]
Since we assumed $\phi_l(k)> f_l(k)$, $ k-\phi_l(k)/k <  k-f_l(k)/k  $, and the function $f_l$ being
increasing, we obtain
\[
f_l(k)< f_l( k-\phi_l(k)/k)+f_{l-1}(k) \leq f_l( k-f_l(k)/k)+f_{l-1}(k).
\]
We will now show that this leads to a contradiction. Set $\epsilon=1-l$. The inequality above is
\[
k^{2-2^\epsilon}< \left (k-\frac{k^{2-2^\epsilon}}{k}\right)^{2-2^\epsilon}+k^{2-2^{\epsilon+1}},
\]
or equivalently
\[
1< \left (1- k ^{-2^{\epsilon}} \right)^{2-2^\epsilon}+ k^{2^\epsilon-2^{\epsilon+1} }=
\left (1- k ^{-2^{\epsilon}} \right)^{2-2^\epsilon} + k^{-2^{\epsilon} } .
\]
We obtain
\[
1 - k^{-2^{\epsilon} }  <  \left(1- k ^{-2^{\epsilon}} \right)^{2-2^\epsilon},
\]
which is impossible because $0<1 - k^{-2^{\epsilon} }<1$ and $2-2^\epsilon \geq 1$.
\end{proof}

\section{Stable endomorphisms} \label{stable}

\begin{defi}
Let $V$ be a vector space over $\bbC$ and $t \in \End_\bbC (V)$. We say that $t$ is stable if
$\ker t= \ker t^2$ and $\im t= \im t^2$. We say that $t$ is eventually stable if there exists
$n \in \bbN^*$ such that $t^n$ is stable. We use the same terminology for endomorphisms of $B$-modules,
where $B$ is a $\bbC$-algebra.
\end{defi}

\begin{lemme}
$(i)$ Let $t \in \End_\bbC (V)$. Then $t$ is stable if and only if
\[
V=\ker t  \oplus \im t.
\]

$(ii)$ Consider the commutative diagram
\[
\xymatrix{  A \ar[r]^{\phi}   \ar[d]_{f_A}  &B \ar[r]^{\psi}\ar[d]_{f_B} &C \ar[r]
\ar[d]_{f_C} &0  \\
A \ar[r]_{\phi}& B \ar[r]_{\psi} &C \ar[r] & 0}
\]
where the rows are exact and the endomorphisms $f_A$ and $f_B$ are stable. Then $f_C$ is stable.

$(iii)$ Consider the commutative diagram
\[
\xymatrix{  0 \ar[r] & A \ar[r]^{\phi}   \ar[d]_{f_A}  &B \ar[r]^{\psi}\ar[d]_{f_B} &C
\ar[d]_{f_C}   \\
0 \ar[r] & A \ar[r]_{\phi}& B \ar[r]_{\psi} &C }
\]
where the rows are exact and the endomorphisms $f_B$ and $f_C$ are stable. Then $f_A$ is stable.

$(iv)$ Suppose $t \in \End_\bbC(V)$ is eventually stable. Then the sequences $N_k=\ker t^k$ and
$I_k=\im t^k$ are respectively increasing and decreasing, both stationary starting from the first index
$k_0$ such that $N_{k_0}=N_{k_0+1}$ or $I_{k_0}=I_{k_0+1}$.
\end{lemme}

\begin{proof}
$(i)$ If $t$ is stable, and if $v \in \ker t \cap \im t$, let us write $v=t(w)$, for a certain $w \in V$.
We then have $ t(v)=t^2(w)=0 $, thus $w \in \ker t^2 =\ker t$, i.e., $v=t(w)=0$. This shows that
$\ker t \cap \im t =\{0\}$. Let us now show that any $v \in V$ is in $\ker t + \im t$. We have
$t(v) \in \im t=\im t^2$, thus $t(v)=t^2(w)$ for a certain $w \in V$. We then have $t(v-t(w))=0$, i.e.,
$v-t(w)\in \ker t$, which establishes the assertion. Conversely, if $V=\ker t \oplus \im t$, then for all
$v \in \ker t^2$, $t(v)\in \ker t \cap  \im t=\{0\}$, whence $v \in \ker t$, and thus
$\ker t^2 \subset \ker t$. The inclusion in the other direction being trivial, we obtain
$\ker t^2 =\ker t$. Similarly if $v=t(w)\in \im t$, we write $w=w_0+w_*$ with $w_0\in \ker t$ and
$w_* \in \im t$, and we obtain $v=t(w_0+w_*)=t(w_*)\in \im t^2$. Thus $\im t \subset \im t^2$, and thus
$\im t^2 = \im t$. The endomorphism $t$ is therefore stable.

$(ii)$ Let $c \in C$. Since $\psi$ is surjective, there exists $b \in B$ such that $\psi(b)=c$. Since
$f_B$ is stable, by $(i)$, we can write $b=b_0+b_*$ with $b_0 \in \ker f_B$ and $b_* \in \im f_B$. We
then have $c =\psi(b_0)+\psi(b_*)$ with $f_C(\psi(b_0))=\psi(f_B(b_0))=0$, thus $\psi(b_0)\in \ker f_C$,
and $\psi(b_*)\in \psi(f_B(B))=f_C(\psi(B))\subset \im f_C$. This shows that $C=\ker f_C + \im f_C$.

Suppose $c \in  \ker f_C \cap  \im f_C$. Write $c=f_C(c')$, with $c'=\psi(b'_0)+\psi(b'_*)$,
$b'_0 \in \ker f_B$ and $b_*' \in \im f_B$. We obtain
\[
c=f_C(c')=f_C(\psi(b'_0))+f_C(\psi(b'_*))=\psi(f_B(b'_0))+\psi(f_B(b'_*))=\psi(f_B(b'_*)).
\]
Moreover, $f_C(c)=0$ gives
\[
0=f_C(c) =f_C(\psi(f_B(b'_*)))=\psi(f_B^2(b'_*)).
\]
Thus $f_B^2(b'_*) \in \ker \psi=\im \phi$. There therefore exists $a \in A$ such that
$f_B^2(b'_*)=\phi(a)$. Writing $a=a_0+a_*$, with $a_0 \in \ker f_A$ and $a_* \in \im f_A$, we obtain
\[
f_B^2(b'_*)=\phi(a_0)+\phi(a_*).
\]
Now $f_B(\phi(a_0))=\phi(f_A(a_0))=0$, thus $\phi(a_0)\in \ker f_B$, and
$$\phi(a_*)\in \phi(f_A(A))=f_B(\phi(A))\subset \im f_B.$$ Since $B=\im  f_B \oplus \ker f_B$, we deduce
that $\phi(a_0)=0$ and
\[
f_B^2(b'_*)=\phi(a_*)=\phi(f_A(a'))=f_B(\phi(a'))
\]
for a certain $a'=f_A(a'') \in \im f_A$. This gives
\[
f_B(b'_*)-\phi(f_A(a''))= f_B(b'_*-\phi(a''))\in \ker f_B\cap \im f_B=\{0\}.
\]
We can conclude that
\[
c=\psi(f_B(b'_*))=\psi(\phi(a'))=0,
\]
and thus $\ker f_C \cap  \im f_C=\{0\}$.

$(iii)$ Suppose that $a\in \ker f_A \cap \im f_A$. Then $\phi(a)\in \ker f_B$ because
$f_B(\phi(a))=\phi(f_A(a))=0$ and $\phi(a)\in \im f_B$ because $\phi(a) \in \phi(f_A(A))=f_B(\phi(A))$.
Since $f_B$ is stable, by $(i)$, we then have $\phi(a)=0$, whence $a=0$ because $ \phi$ is injective.
This shows that $\ker f_A \cap \im f_A=\{0\}$.

We now want to show that any $a \in A$ is in $\ker f_A + \im f_A$. Write $\phi(a)=b_0+f_B(b_*)$ with
$b_0 \in \ker f_B$, and $b_* \in \im f_B$, which is permissible since $f_B$ is stable. Suppose that
\begin{equation}\label{fAfB}
b_0 \in \im \phi=\ker \psi, \quad b_*\in \im \phi=\ker \psi,
\end{equation}
and let $a_0,a_* \in A$ such that $b_0=\phi(a_0)$, $b_*=\phi(a_*)$. We then have
\[
\phi(a)=\phi(a_0)+ f_B(\phi(a_*))=\phi(a_0+f_A(a_*)),
\]
and since $\phi$ is injective, $a=a_0+f_A(a_*)$. On the other hand
\[
0=f_B(b_0)= f_B(\phi(a_0))=\phi(f_A(a_0)),
\]
whence since $\phi$ is injective, $f_A(a_0)=0$, and thus $a_0 \in \ker f_A$. Thus
$a \in \ker f_A + \im f_A$. It remains to prove (\ref{fAfB}). We have
\[
\psi(\phi(a))=0=\psi(b_0+f_B(b_*))=\psi(b_0)+\psi(f_B(b_*)).
\]
Now $f_C(\psi(b_0))=\psi(f_B(b_0))=0$, whence $\psi(b_0) \in \ker f_C$ and
$\psi(f_B(b_*))=f_C(\psi(b_*))\in \im (f_C)$. Since $f_C$ is stable, we deduce $\psi(b_0)=0$ and
$\psi(f_B(b_*))=f_C(\psi(b_*))=0$. This shows that $b_0,f_B(b_*)\in \ker \psi$. Set $\psi(b_*)=c_0+c_*$,
with $c_0 \in \ker f_C$ and $c_*\in \im f_C$. We then have
\[
f_C(\psi(b_*))=\psi(f_B(b_*))=0=f_C(c_0+c_*)=f_C(c_*).
\]
We then see that $c_*  \in \im f_C \cap \ker f_C=\{0\}$. We obtain
$\psi(b_*)=c_0\in \ker f_C \cap \im f_C=\{0\}$, and thus $b_*\in \ker \psi$.

$(iv)$ It is clear that the sequences $N_k$ and $I_k$ are respectively increasing and decreasing, and by
hypothesis, they both become stationary starting from a certain rank. Set $E_k=N_{k+1}/N_{k}$ and
$F_k=I_k/I_{k+1}$. The restriction of $t^k$ to $N_{k+1}$ takes values in $I_k$, is zero on $N_k$, and
thus induces a morphism $\bar t^k\colon E_k \rightarrow F_k$. The kernel of this morphism is
$N_k+(\im t\cap N_{k+1})/N_k$ and its image is $I_{k+1}+t^k(N_{k+1})/I_{k+1}$.

We also have a morphism $\bar t\colon E_{k+1} \rightarrow E_{k}$ induced by
\[
N_{k+2} \stackrel{t}{\rightarrow}  N_{k+1} \rightarrow   N_{k+1} /N_k=E_k
\]
whose kernel is exactly $N_{k+1}$. The morphism $\bar t$ is therefore injective, with image
\[
N_k+t(N_{k+2})/N_k=N_k+(\im t\cap N_{k+1})/N_k.
\]

Finally, we have a morphism $\bar t\colon F_{k} \rightarrow F_{k+1}$, induced by
\[
I_{k} \stackrel{t}{\rightarrow}  I_{k+1} \rightarrow   I_{k+1} /I_{k+2}=F_{k+1}
\]
whose kernel contains $I_{k+1}$. The morphism $\bar t$ is surjective, with kernel
$I_{k+1}+t^k(N_{k+1})/I_{k+1}$.

Thus, we have constructed an exact sequence
\[
0   \rightarrow  E_{k+1} \stackrel{\bar t}{\rightarrow} E_{k}  \stackrel{\bar t^k}{\rightarrow}
F_k  \stackrel{\bar t}{\rightarrow}  F_{k+1} \rightarrow 0
\]

We immediately deduce that the sequence $N_k$ becomes stationary from the first index $k_0$ such that
$N_{k_0}=N_{k_0+1}$ and that the same is true for the sequence $I_k$. Suppose that the sequence $N_k$
becomes stationary starting from the index $k_0$. Then $E_k=0$ for all $k\geq k_0$, and we deduce from
the exact sequence above that $F_k\simeq F_{k+1}$ for all $k\geq k_0$. Since the sequence $I_k$ becomes
stationary starting from a certain rank, for $k$ large enough $F_k=0$, and this shows that $F_k=0$ for
all $k\geq k_0$. The sequence $I_k$ therefore becomes stationary no later than the sequence $N_k$, but
since the same argument applies in the other direction, we deduce the final assertion of $(iv)$.
\end{proof}

\begin{prop}
Let $B$ be a Noetherian commutative $\bbC$-algebra, and $L$ a finitely generated $B$-module. Let
$a \in \End_B (L)$ and suppose that the localization $(L_a,\tilde a)$ of $L$ at $a$ (see the fifth
example of \ref{propuniv}) is a finitely generated $B$-module. Then $a$ is eventually stable. Let
$n \in \bbN^*$ such that $a^n$ is stable, and let $K=\ker a^n$, $I =\im a^n$. Then $L=K\oplus I$ and the
localization $(L_a,\tilde a)$ is isomorphic to $(I,a_{|I})$.
\end{prop}

\begin{proof}
We can equip $\End_B(L_a)$ with a $B$-module structure via
\[
(b\cdot \phi)(v)= b\cdot (\phi(v))=\phi(b\cdot v), \quad (\phi \in
\End_B(L_a), \, b \in B, v\,  \in L_a).
\]
Since $L_a$ is a finitely generated $B$-module, it is easy to see that $\End_B(L_a)$ is as well, and thus
we can find elements $b_i\in B$, $i=1,\ldots,m$ such that
\[
\tilde a^{-m}+b_1 \tilde a^{-m+1}+ \cdots +b_m=0
\]
whence $ \tilde a^{-1}=-(b_1+b_2 \tilde a+\cdots +b_m \tilde a^{m-1})$, and thus
$ \tilde a^{-1}\in B[\tilde a]$. We will now show that the morphism $\iota \colon L \rightarrow L_a$
(fifth example of \ref{propuniv}) is surjective. It is clear that $ \tilde a$ stabilizes $\iota(L)$, and
since $ \tilde a^{-1}\in B[\tilde a]$, $ \tilde a^{-1}$ stabilizes $\iota(L)$ and $\tilde a_{|\iota(L)}$
is therefore invertible. Consider the pair $(\iota(L),\tilde a)$. It is also manifestly a solution of the
same universal problem as $(L_a,\tilde a)$. We deduce the surjectivity of $\iota$.

We now use the lemma and Corollary \ref{propuniv}, which show that $(L_{a},\tilde a)$ is isomorphic to
$(L',a')$, where $L'=L/K$, with $K=\cup_{n\in \bbN}\ker a^n$. Since $L_a$ and $L$ are finitely generated
$B$-modules, $K$ is as well. Consequently, there exists a non-zero integer $n$ such that $K=\ker a^n$.
Set $I =\im a^n$. We then have $L=K\oplus I$: indeed, $K\cap I=0$, since if $a^ny\in K=\ker a^n$, then
$y\in\ker a^{2n}=K$, and $a^ny=0$; on the other hand, $\iota$ induces an isomorphism of $L/K$ onto $L_a$,
on which $\tilde a$ is invertible, hence $L=a^n(L)+K$. Since $\iota$ is surjective, it induces an
isomorphism $(I,a_{|I})\simeq (L_a,\tilde a)$. Now $\tilde a$ is invertible, and thus $a_{|I}$ is
invertible.
\end{proof}

\section{Representations of $\bbZ^d$}\label{repsZd}

Let $V$ be a finite-dimensional vector space over $\bbC$, and $u_1,\ldots, u_d$ invertible endomorphisms
of $V$ commuting pairwise. This is equivalent to the data of a representation $U$ of the abelian group
$\bbZ^d$ in $V$,
\[
U\colon \underline{n}=(n_1,\ldots ,n_d) \mapsto u_1^{n_1}\circ \cdots
\circ u_d^{n_d}.
\]

Let us fix $v \in V$ and $\lambda \in V^*$ and form the matrix coefficient
\[
F\colon \bbZ^d \rightarrow \bbC, \quad  F(\underline{n})=\lambda(
U(\underline{n})\cdot v).
\]

\begin{prop}
There exist characters $\chi_1, \ldots, \chi_l$ of $\bbZ^d$ and polynomials
\[
Q_1,\ldots ,Q_l\in \bbC[X_1,\ldots, X_d]
\]
such that
\[
F(\underline{n})=\sum_{j=1}^l \chi_j(\underline{n})Q_j(\underline{n}).
\]
\end{prop}

\begin{proof}
Since the $u_i$ commute pairwise, they admit common eigenspaces, and if $z_1, \ldots ,z_d$ are the
eigenvalues respectively of $u_1,\ldots, u_d$ on one of these eigenspaces, then $U(\underline{n})$ acts
on this eigenspace as the scalar
\[
\chi(\underline{n})=z_1^{n_1} \ldots z_d^{n_d}, \quad
(\underline{n}=(n_1,\ldots, n_d)\in \bbZ^d).
\]
Let $\chi_1, \ldots, \chi_l$ be all the pairwise distinct characters of $\bbZ^d$ thus obtained, the
character $\chi_j$ being given by a $d$-tuple $(z_1^{(j)}, \ldots, z_d^{(j)})$ of non-zero complex
numbers.

Let us choose a common triangularization basis of the $u_1,\ldots, u_d$ such that for all
$ \underline{n} \in \bbZ^d$, the matrix $A(\underline{n})$ of $U(\underline{n})$ in this basis is block
diagonal, with blocks $A_j(\underline{n})$, $A_j(\underline{n})$ being an upper triangular matrix whose
diagonal elements are of the form $\chi_j(\underline{n})$, $j=1,\ldots ,l$. It suffices to show that the
coefficients of the matrix $A(\underline{n})$ are of the desired form, and for this, it suffices to do it
for each block $A_j$. We therefore fix $j$. Let $(\underline{e_1},\ldots ,\underline{e_d})$ denote the
canonical basis of $\bbZ^d$. We then have
\[
A_j(\underline{e_k})= z_k^{(j)}(\Id+N_k^{(j)}), \quad (k=1, \ldots ,d)
\]
where the $N_k^{(j)}$ are strictly upper triangular matrices (thus nilpotent), commuting pairwise. This
gives us
\begin{align*}
 A_j(\underline{n})&= \prod_{i=1}^d  \left( z_i^{(j)}(\Id+N_i^{(j)})\right)^{n_i}
= \chi_j(\underline{n}) \prod_{i=1}^d (\Id+N_i^{(j)})^{n_i}\\
&=\chi_j(\underline{n}) \prod_{i=1}^d \left( \Id+ \binom{n_i}{1} N_i^{(j)}+  \binom{n_i}{2} (N_i^{(j)})^2 +  \binom{n_i}{3 } (N_i^{(j)})^3 +\cdots \right)
\end{align*}
the sum being finite, since the $N_i^{(j)}$ are nilpotent. The binomial coefficients being polynomials in
the $n_i$, this completes the proof.
\end{proof}

\cleardoublepage
\addcontentsline{toc}{chapter}{\bibname}
\bibliographystyle{plain}
\bibliography{English-26-padic}

\cleardoublepage
\addcontentsline{toc}{chapter}{List of notation}
\printindex[not]

\cleardoublepage
\addcontentsline{toc}{chapter}{Subject Index}
\printindex[ter]

\vfill
\pagebreak
\clearpage
\ifodd\value{page}\thispagestyle{empty}\null
\clearpage
\fi
\thispagestyle{empty}

\begin{tikzpicture}[remember picture, overlay]
    \fill[MainColor] (current page.north west) rectangle ([yshift=-2cm]current page.north east);

    \fill[MainColor] (current page.south west) rectangle ([yshift=2cm]current page.south east);
    \node[text=white, font=\large\scshape, anchor=west]
        at ([xshift=3cm, yshift=1cm]current page.south west)
        {};
\end{tikzpicture}

\vspace*{2cm}

\begin{center}
{\LARGE \bfseries \textcolor{MainColor}{Representations of $p$-adic Reductive Groups}}\\[0.5cm] {\Large
David Renard}
\end{center}

\vspace{1.5cm}

\noindent This book presents a part of the (complex) representation theory of $p$-adic reductive groups.
Starting from basics accessible to graduate students, it culminates with the theory of the ``Bernstein
center'' and the Langlands classification of irreducible smooth representations.

\vspace{0.5cm}

\noindent The book consists of seven chapters, Chapters VI and VII constituting the core of the book.
Chapter VI is devoted to the study of the category of smooth representations of a $p$-adic reductive
group, establishing among other things the Bernstein decomposition theorem and the description of its
center. Chapter VII deals with square-integrable and tempered representations, and the Langlands
classification theorem is proved there.

\vspace{0.5cm}

\noindent The first four chapters are placed in a more general framework and tackle, in order, the study
of algebras with idempotents, totally disconnected locally compact spaces and groups, smooth
representations of the latter, and particular classes of representations (compact, unitary,
square-integrable). Chapter V is a review of the structure results of $p$-adic reductive groups. The
appendices give the elements of category theory and some results in algebra used in the text.

\vfill
\
\end{document}